\documentclass[10pt]{article}
\usepackage{amsmath,amssymb,amscd}
\usepackage{makeidx}
\makeindex
\newcommand{\nbigl}{\mathcal{L}}
\newcommand{\nbigy}{\mathcal{Y}}
\newcommand{\nbigx}{\mathcal{X}}
\newcommand{\nbige}{\mathcal{E}}
\newcommand{\nbigo}{\mathcal{O}}
\newcommand{\nbigz}{\mathcal{Z}}
\newcommand{\nbigc}{\mathcal{C}}
\newcommand{\nbign}{\mathcal{N}}
\newcommand{\nbigm}{\mathcal{M}}
\newcommand{\nbigh}{\mathcal{H}}
\newcommand{\nbigf}{\mathcal{F}}
\newcommand{\nbiga}{\mathcal{A}}
\newcommand{\nbigb}{\mathcal{B}}
\newcommand{\nbigu}{\mathcal{U}}
\newcommand{\nbigd}{\mathcal{D}}
\newcommand{\nbigi}{\mathcal{I}}
\newcommand{\nbigp}{\mathcal{P}}
\newcommand{\nbigk}{\mathcal{K}}
\newcommand{\nbigs}{\mathcal{S}}
\newcommand{\nbigt}{\mathcal{T}}
\newcommand{\nbigq}{\mathcal{Q}}
\newcommand{\nbigr}{\mathcal{R}}
\newcommand{\nbigw}{\mathcal{W}}
\newcommand{\nbigv}{\mathcal{V}}
\newcommand{\nbigg}{\mathcal{G}}

\newcommand{\nbigss}{\mathcal{SS}}

\newcommand{\proj}{\mathbb{P}}

\newcommand{\seisuu}{\mathbb{Z}}
\newcommand{\rnum}{\mathbb{Q}}

\newcommand{\cnum}{\mathbb{C}}
\newcommand{\real}{\mathbb{R}}
\newcommand{\hyperh}{\mathbb{H}}

\newcommand{\gbigy}{\mathfrak Y}

\newcommand{\gbigr}{\mathfrak R}

\newcommand{\gbign}{\mathfrak N}

\newcommand{\gbigb}{\mathfrak B}
\newcommand{\gbigc}{\mathfrak C}

\newcommand{\gbigv}{\mathfrak V}
\newcommand{\gbigl}{\mathfrak L}
\newcommand{\gbigt}{\mathfrak T}
\newcommand{\gbigd}{\mathfrak D}

\newcommand{\gbigq}{\mathfrak Q}
\newcommand{\gbigg}{\mathfrak G}
\newcommand{\gbigf}{\mathfrak F}
\newcommand{\gbigi}{\mathfrak I}

\newcommand{\gminig}{\mathfrak g}

\newcommand{\gminih}{\mathfrak h}
\newcommand{\gminik}{\mathfrak k}
\newcommand{\gminip}{\mathfrak p}
\newcommand{\gminiw}{\mathfrak w}
\newcommand{\gminif}{\mathfrak f}

\newcommand{\gminir}{\mathfrak r}
\newcommand{\gminis}{\mathfrak s}

\newcommand{\gminiy}{\mathfrak y}
\newcommand{\gminic}{\mathfrak c}
\newcommand{\gminiq}{\mathfrak q}

\newcommand{\vexi}{\boldsymbol \xi }

\newcommand{\vecphi}{\boldsymbol \phi}
\newcommand{\vecdelta}{\boldsymbol \delta}

\newcommand{\vecL}{\boldsymbol L}

\newcommand{\vecf}{\boldsymbol f}

\newcommand{\vecv}{\boldsymbol v}
\newcommand{\vecu}{\boldsymbol u}
\newcommand{\vecy}{\boldsymbol y}

\newcommand{\vecN}{\boldsymbol N}
\newcommand{\vecr}{\boldsymbol r}
\newcommand{\veck}{\boldsymbol k}
\newcommand{\vecI}{\boldsymbol I}
\newcommand{\vecn}{\boldsymbol n}
\newcommand{\vecY}{\boldsymbol Y}
\newcommand{\veca}{\boldsymbol a}
\newcommand{\vecV}{\boldsymbol V}
\newcommand{\vecU}{\boldsymbol U}
\newcommand{\vecE}{\boldsymbol E}

\newcommand{\llarr}{\longleftarrow}

\newcommand{\rarr}{\rightarrow}
\newcommand{\lrarr}{\longrightarrow}

\newcommand{\pf}{{\bf Proof}\hspace{.1in}}
\newcommand{\qed}{\mbox{\rule{1.2mm}{3mm}}}

\def\Hom{\mathop{\rm Hom}\nolimits}
\def\Sym{\mathop{\rm Sym}\nolimits}
\def\Ext{\mathop{\rm Ext}\nolimits}

\def\Cok{\mathop{\rm Cok}\nolimits}
\def\Image{\mathop{\rm Im}\nolimits}
\def\CH{\mathop{\rm CH}\nolimits}
\def\Spec{\mathop{\rm Spec}\nolimits}
\def\type{\mathop{\rm type}\nolimits}
\def\pardeg{\mathop{\rm par\textrm{-}deg}\nolimits}
\def\Gr{\mathop{\rm Gr}\nolimits}

\def\Cone{\mathop{\rm Cone}\nolimits}
\def\rank{\mathop{\rm rank}\nolimits}

\def\Ker{\mathop{\rm Ker}\nolimits}

\def\ch{\mathop{\rm ch}\nolimits}
\def\Pic{\mathop{\rm Pic}\nolimits}
\def\id{\mathop{\rm id}\nolimits}
\def\br{\mathop{\rm br}\nolimits}
\def\top{\mathop{\rm top}\nolimits}
\def\depth{\mathop{\rm depth}\nolimits}
\def\Cr{\mathop{\rm Cr}\nolimits}
\def\Td{\mathop{\rm Td}\nolimits}
\def\YOK{\mathop{\rm YOK}\nolimits}
\def\Flag{\mathop{\rm Flag}\nolimits}
\def\TH{\mathop{\rm TH}\nolimits}
\def\SL{\mathop{\rm SL}\nolimits}
\def\GL{\mathop{\rm GL}\nolimits}
\def\SGL{\mathop{\rm SGL}\nolimits}
\def\PGL{\mathop{\rm PGL}\nolimits}
\def\sign{\mathop{\rm sign}\nolimits}
\def\spl{\mathop{\rm split}\nolimits}
\def\rel{\mathop{\rm rel}\nolimits}
\def\qcoh{\mathop{\rm qcoh}\nolimits}
\def\End{\mathop{\rm End}\nolimits}
\def\Ob{\mathop{\rm Ob}\nolimits}
\def\ob{\mathop{\rm ob}\nolimits}
\def\quo{\mathop{\rm quo}\nolimits}
\def\Quot{\mathop{\rm Quot}\nolimits}
\def\ori{\mathop{\rm or}\nolimits}
\def\ev{\mathop{\rm ev}\nolimits}
\def\tr{\mathop{\rm tr}\nolimits}
\def\Eu{\mathop{\rm Eu}\nolimits}
\def\inv{\mathop{\rm inv}\nolimits}
\def\mov{\mathop{\rm mov}\nolimits}

\def\Res{\mathop{\rm Res}\nolimits}
\def\SW{\mathop{\rm SW}\nolimits}
\def\Dec{\mathop{\rm Dec}\nolimits}
\def\Map{\mathop{\rm Map}\nolimits}
\def\Ind{\mathop{\rm Ind}\nolimits}
\def\cycl{\mathop{\rm cycl}\nolimits}
\def\pt{\mathop{\rm pt}\nolimits}

\newcommand{\res}{\Res}

\newcommand{\YOKtilde}{\widetilde{\YOK}}
\newcommand{\nYOK}{\mathcal{Y}\mathcal{O}\mathcal{K}}
\newcommand{\nYOKtilde}{\widetilde{\nYOK}}

\newcommand{\THhat}{\widehat{\TH}}
\newcommand{\Mhat}{\widehat{M}}
\newcommand{\Ehat}{\widehat{E}}
\newcommand{\phihat}{\widehat{\phi}}
\newcommand{\Ohat}{\widehat{O}}

\newcommand{\del}{\partial}

\newcommand{\nhom}{{\mathcal Hom}}
\newcommand{\nHom}{\nhom}
\newcommand{\nEnd}{{\mathcal End}}
\newcommand{\nrhom}{{\mathcal R \mathcal Hom}}

\newcommand{\poin}{{\mathcal Poin}}
\newcommand{\Poin}{{\mathcal Poin}}
\newcommand{\Or}{{\mathcal Or}}
\newcommand{\yw}{Y(W_{\cdot})}

\newcommand{\ywquo}{ Y_{\quo}(W_{\cdot}) }

\newcommand{\phibar}{\overline{\phi}}

\newcommand{\phitilde}{\widetilde{\phi}}

\newcommand{\wx}[1]{W_{#1\,X} }

\newcommand{\nbigahat}{\widehat{\nbiga}}

\newcommand{\vecyhat}{\widehat{\vecy}}

\newcommand{\vmx}{V_{m,X}}

\newcommand{\qmy}{Q(m,y)}

\newcommand{\Type}{{\mathcal Type}}
\newcommand{\Typetilde}{\widetilde{\Type}}

\newcommand{\yhat}{\widehat{y}}
\newcommand{\Psihat}{\widehat{\Psi}}
\newcommand{\Psitilde}{\widetilde{\Psi}}
\newcommand{\Amss}{\nbiga_{m}^{ss} \bigl(
  P_{\vecy}^{\alpha_{\ast},\delta}(m),\epsilon_{\ast}, \delta(m) \bigr)}
\newcommand{\Ams}{\nbiga_{m}^{s} \bigl(
  P_{\vecy}^{\alpha_{\ast},\delta}(m),\epsilon_{\ast}, \delta(m)\bigr)}
\newcommand{\Qmss}{Q^{ss}(m,\vecy,[L],\alpha_{\ast},\delta)}

\newcommand{\Nbar}{\underline{N}}
\newcommand{\ellbar}{\underline{\ell}}

\newcommand{\Qtilde}{\widetilde{Q}}
\newcommand{\nbigatilde}{\widetilde{\nbiga}}
\newcommand{\nbigqtilde}{\widetilde{\nbigq}}
\newcommand{\nbigltilde}{\widetilde{\nbigl}}
\newcommand{\nbigbtilde}{\widetilde{\nbigb}}
\newcommand{\rhobar}{\overline{\rho}}
\newcommand{\nbigmtilde}{\widetilde{\nbigm}}
\newcommand{\Gtilde}{\widetilde{G}}

\newcommand{\Ttilde}{\widetilde{T}}
\newcommand{\Tbar}{\overline{T}}
\newcommand{\tbar}{\overline{t}}

\newcommand{\nbigmbar}{\overline{\nbigm}}
\newcommand{\rhotilde}{\widetilde{\rho}}
\newcommand{\pitilde}{\widetilde{\pi}}
\newcommand{\ptilde}{\widetilde{p}}

\newcommand{\gtilde}{\widetilde{g}}
\newcommand{\Ftilde}{\widetilde{F}}
\newcommand{\Fbar}{\overline{F}}
\newcommand{\Ytilde}{\widetilde{Y}}

\newcommand{\cone}{\Cone}

\newcommand{\iotatilde}{\widetilde{\iota}}

\newcommand{\Bhat}{\widehat{B}}
\newcommand{\Btilde}{\widetilde{B}}
\newcommand{\Bbar}{\overline{B}}

\newcommand{\nbigmhat}{\widehat{\nbigm}}
\newcommand{\Wtilde}{\widetilde{W}}
\newcommand{\Obtilde}{\widetilde{\Ob}}
\newcommand{\Obcheck}{\overline{\Ob}}
\newcommand{\Obbar}{\overline{\Ob}}
\newcommand{\obtilde}{\widetilde{\ob}}
\newcommand{\Mtilde}{\widetilde{M}}
\newcommand{\Ctilde}{\widetilde{C}}
\newcommand{\Ubar}{\overline{U}}
\newcommand{\Zhat}{\widehat{Z}}
\newcommand{\psitilde}{\widetilde{\psi}}
\newcommand{\nbigehat}{\widehat{\nbige}}
\newcommand{\nbigvhat}{\widehat{\nbigv}}
\newcommand{\Zbar}{\overline{Z}}

\newcommand{\Phitilde}{\widetilde{\Phi}}
\newcommand{\nbiglhat}{\widehat{\nbigl}}
\newcommand{\nbigktilde}{\widetilde{\nbigk}}
\newcommand{\gbigbtilde}{\widetilde{\gbigb}}
\newcommand{\gbigvtilde}{\widetilde{\gbigv}}
\newcommand{\Vtilde}{\widetilde{V}}
\newcommand{\Qtildehat}{\widehat{\widetilde{Q}}}
\newcommand{\nbigbtildehat}{\widehat{\widetilde{\nbigb}}}
\newcommand{\Qhat}{\widehat{Q}}
\newcommand{\nbigatildehat}{\widehat{\widetilde{\nbiga}}}
\newcommand{\gminifhat}{\widehat{\gminif}}
\newcommand{\Ybar}{\overline{Y}}
\newcommand{\Phibar}{\overline{\Phi}}
\newcommand{\gminigbar}{\overline{\gminig}}
\newcommand{\Cbar}{\overline{C}}
\newcommand{\nSW}{\mathcal S \mathcal W}
\newcommand{\nbigsbar}{\overline{\nbigs}}
\newcommand{\Sbar}{\overline{S}}
\newcommand{\gbigihat}{\widehat{\gbigi}}
\newcommand{\vecYhat}{\widehat{\vecY}}
\newcommand{\gminiyhat}{\widehat{\gminiy}}
\newcommand{\vecxi}{\vexi}

\newcommand{\nbigvtilde}{\widetilde{\nbigv}}
\newcommand{\dettilde}{\widetilde{\det}}
\newcommand{\gminigtilde}{\widetilde{\gminig}}

\newcommand{\bdmath}{\begin{displaymath}}
\newcommand{\edmath}{\end{displaymath}}

\newcommand{\beqn}{\begin{equation}}
\newcommand{\eeqn}{\end{equation}}

\newcommand{\beqnarray}{\begin{eqnarray}}
\newcommand{\eeqnarray}{\end{eqnarray}}

\newcommand{\bitemize}{\begin{itemize}}
\newcommand{\eitemize}{\end{itemize}}

\newcommand{\benumerate}{\begin{enumerate}}
\newcommand{\eenumerate}{\end{enumerate}}

\newcommand{\bdescriprion}{\begin{description}}
\newcommand{\edescriprion}{\end{description}}
\newtheorem{thm}{Theorem}[section]
\newtheorem{cor}[thm]{Corollary}

\newtheorem{rem}[thm]{Remark}
\newtheorem{lem}[thm]{Lemma}
\newtheorem{prop}[thm]{Proposition}
\newtheorem{df}[thm]{Definition}

\newtheorem{problem}[thm]{Problem}
\newtheorem{example}[thm]{Example}
\newtheorem{notation}[thm]{Notation}
\newtheorem{condition}[thm]{Condition}

\begin{document}

\title{A Theory of the Invariants Obtained from
 the Moduli Stacks of \\
 Stable Objects on a Smooth Polarized Surface
}
\author{Takuro Mochizuki}
\date{}

\maketitle

\begin{abstract}

Let $X$ be a smooth polarized algebraic surface
over the compex number field.
We discuss the invariants obtained from
the moduli stacks of semistable sheaves of arbitrary ranks on $X$.
For that purpose,
we construct the virtual fundamental classes
of some moduli stacks,
and we show the transition formula
of the integrals over
the moduli stacks of the $\delta$-stable Bradlow pairs
for the variation of the parameter $\delta$.

Then, we study the relation among the invariants.
In the case $p_g=\dim H^2(X,\nbigo_X)>0$,
we show that the invariants are independent
of the choice of a polarization of $X$.
We also show that the invariants 
can be reduced to the invariants obtained from
the moduli of abelian pairs and the Hilbert schemes.
In the case $p_g=0$,
we obtain the weak wall crossing formula
and the weak intersection rounding formula,
which describes the dependence of
the invariants on the polarization.

\vspace{.1in}

\vspace{.1in}

\end{abstract}

\tableofcontents

\section{Introduction}

\subsection{Problems}
\subsubsection{Construction of the invariants}

Let $X$ be a smooth projective surface
with an ample line bundle $\nbigo_X(1)$
over the complex number field $\cnum$.
We assume that $X$ is simply connected in {\em Introduction},
for simplicity.
Let $y$ be an element of
the cohomology group $H^{\ast}(X,\rnum)$,
which is the Chern character of
a torsion-free coherent sheaf $E$ on $X$.
Let $a$ be the first Chern class determined by $y$.
We take a line bundle $\nbigl_a$ on $X$
such that $c_1(\nbigl_{a})=a$.
An oriented torsion-free sheaf on $X$ of type $y$
is defined to be a coherent torsion-free sheaf
$E$ such that $\ch(E)=y$
with an isomorphism $\rho:\det(E)\simeq \nbigl_a$.
Let $\nbigm^{ss}(\yhat)$ (resp. $\nbigm^{s}(\yhat)$)
denote the moduli stack of 
the semistable (resp. stable)
oriented torsion-free sheaves of type $y$.
One of our main problems is the following.
\begin{problem}
Construct an invariant from the moduli
stack $\nbigm^{ss}(\yhat)$.
\hfill\qed
\end{problem}

Let $\Ehat^u$ denote the universal sheaf
over $\nbigm^{ss}(\yhat)\times X$.
Let $\Phi=P(\Ehat^u)$ be a polynomial
of the slant products
$\ch_i(\Ehat^u)/b$ for elements $b\in H^{\ast}(X)$ 
and $i\in\seisuu_{\geq\,0}$.
Naively speaking, we would like to obtain the number:
\begin{equation}
\label{eq;06.6.8.1}
 \Phi(\yhat):=\int_{\nbigm^{ss}(\yhat)}\Phi
=\deg\bigl( \Phi\cap[\nbigm^{ss}(\yhat)]\bigr)
\end{equation}
Namely,
we would like to obtain the $0$-cycle
by taking the cap product
of $\Phi$ and the fundamental class
$[\nbigm(\yhat)]$,
and we would like to obtain the number
by taking the degree of the $0$-cycle.

There are two main problems
to make (\ref{eq;06.6.8.1}) well-defined.
\begin{description}
\item[(A)]
 The moduli stack $\nbigm^{ss}(\yhat)$
 is not smooth, 
 and hence it does not have the natural 
 fundamental class, in general.
\item[(B)]
 Even if $\nbigm^{ss}(\yhat)$ is smooth,
 the moduli stack $\nbigm^{ss}(\yhat)$
 is not Deligne-Mumford, in general.
 In such a situation, there are no known satisfactory definition
 of the degree of $0$-cycles,
 or in other words, 
 the push forward of cycles.
\end{description}

\begin{rem}
We do not give detail in Introduction
about how we consider the cohomology
and the evaluation on the Deligne-Mumford stacks.
See the subsection {\rm \ref{subsection;06.6.19.110}}.
\hfill\qed
\end{rem}

\begin{rem}
The construction of such invariants
was also discussed in {\rm \cite{Kronheimer}}
for real $4$-dimensional manifolds
from the differential geometric view point.
\hfill\qed
\end{rem}

\subsubsection{Virtual fundamental classes}

The problem like (A) was well discussed
and established in the course of their study of
Gromov-Witten invariants
(\cite{bf}, \cite{l-t}, \cite{f-o}).
In this paper,
we follow the method of Behrend and Fantechi
(\cite{bf}).
Namely, we will show that
some interested moduli stacks
are naturally provided with the perfect obstruction
theories in the sense of \cite{bf}.
(We will review the obstruction theory 
in the subsubsection \ref{subsubsection;06.6.8.3}.)
And, we will obtain the virtual fundamental classes.
Such an obstruction theory may be well known,
perhaps.
For example, 
it is a standard fact that
the first cohomology group $H^1(X,End(E))$
gives the space of the infinitesimal deformations,
and that the second cohomology group
$H^2(X,End(E))$ gives the space of the obstruction,
for any vector bundle $E$ on $X$.
However, the author does not know
an appropriate reference
to deal with the obstruction theory
in the sense of \cite{bf}
for the moduli stacks of reduced oriented $L$-Bradlow pairs
and the master spaces,
which is available for
our arguments using the localization.
Thus, we give a detailed argument
in the section \ref{section;06.6.4.250}.

\begin{rem}
\label{rem;06.6.30.1}
Recently,
the theory of ``derived stacks'' has been developed,
which seems to provide us
a general and powerful tool
to construct obstruction theories for some stacks.
(See {\rm \cite{toen}} for an overview of the theory.)
For example,
the results in {\rm\cite{toen-vaquie}} implies
the construction of the obstruction theory of the moduli stack
of semistable sheaves.
It is not clear to the author, at the present moment,
whether we can directly apply their results
to the moduli stacks of semistable oriented
reduced $L$-Bradlow pairs and the master spaces.
But, it would be quite hopeful to construct
the obstruction theory of such stacks 
and to redo the argument in this paper,
from that point of view.
The author would like to come back to this problem
in future.

However, the author also expects that
we will obtain the same ``invariants'',
even if we adopt the other way of
the construction of the obstruction theories.
(See the subsubsection {\rm\ref{subsubsection;06.7.2.1}}.)
\hfill\qed
\end{rem}

\subsubsection{A naive idea for the construction of the invariants}

To discuss the problem (B),
we recall $L$-Bradlow pair and reduced $L$-Bradlow pair.
Let $L$ be a line bundle on $X$.
Let $U$ be a scheme,
and let $E$ be a torsion-free sheaf on $X$,
which is flat over $U$.
A morphism $\phi:p_X^{\ast}L\lrarr E$
is called an $L$-section.
Such a pair $(E,\phi)$ is called an $L$-Bradlow pair.

Let $M$ be a line bundle on $U$,
and let $[\phi]$ be a morphism
$p_X^{\ast}L\otimes p_U^{\ast}M\lrarr E$
such that $[\phi]_{|\{u\}\times X}\neq 0$
for any $u\in U$.
Such a pair $(M,[\phi])$ is called a reduced $L$-section,
and a pair $\bigl(E,(M,[\phi])\bigr)$ is called a reduced
$L$-Bradlow pair.
We will often omit to denote $M$, for simplicity.

Let $\nbigp^{\br}$ denote the set of
polynomials $\delta$ such that $\deg(\delta)\leq 1$
and $\delta(t)>0$ for any sufficiently large $t$.
Let $\delta$ be any element of $\nbigp$.
Recall that the $\delta$-semistability 
and $\delta$-stability conditions
are defined for $L$-Bradlow pairs
and reduced $L$-Bradlow pairs.
(See the subsubsection \ref{subsubsection;06.6.8.2}.)
We use the notation $\nbigm^{ss}(y,L,\delta)$
(resp. $\nbigm^s(y,L,\delta)$)
to denote the moduli stack
of the semistable (resp. stable)
$L$-Bradlow pairs of type $y$.
We also use the notation $\nbigm^{ss}(\yhat,[L],\delta)$
(resp. $\nbigm^s(\yhat,[L],\delta)$)
to denote the moduli stack of
the $\delta$-semistable (resp. $\delta$-stable)
oriented reduced $L$-Bradlow pairs of type $y$.

We say that $\delta$ is critical,
if $\nbigm^{ss}(\yhat,[L],\delta)=
 \nbigm^{s}(\yhat,[L],\delta)$ does not hold.
It can be shown that 
there are only finitely many critical parameters.
For a non-critical parameter $\delta$,
the moduli stack $\nbigm^{ss}(\yhat,[L],\delta)$
is Deligne-Mumford.
Moreover, it is naturally provided with the perfect
obstruction theory.
Therefore, we have the integrals of cohomology classes
on $\nbigm^{ss}(\yhat,[L],\delta)$.

Let us consider the case $L=\nbigo(-m)$.
If $\delta$ is sufficiently small, it is non-critical,
and we have the naturally defined morphism
$\pi:\nbigm^{ss}(\yhat,[\nbigo(-m)],\delta)\lrarr \nbigm^{ss}(\yhat)$.
It is easy to observe that the morphism $\pi$ 
is smooth, if $m$ is sufficiently large.
The relative tangent bundle is denoted by $T_{\rel}$.
We put $H_y(m):=\int_X\Td(X)\cdot y\cdot \ch\bigl(\nbigo(m)\bigr)$,
which is same as $\dim H^0\bigl(X,E(m)\bigr)$
for any $(E,\rho)\in \nbigm^{ss}(\yhat)$.
We regard $\nbigm^{ss}\bigl(\yhat,[\nbigo(-m)],\delta\bigr)$
as a good approximation of $\nbigm^{ss}(\yhat)$,
and we would like to put as follows:
\begin{equation}
 \label{eq;06.6.8.5}
 \Phi(\yhat):=\int_{\nbigm^{ss}(\yhat,[\nbigo(-m)],\delta)}
 \Phi\cdot \frac{\Eu(T_{\rel})}{H_y(m)}
\end{equation}
Needless to say, we should ask the following:
\begin{description}
\item[(C)]
Is (\ref{eq;06.6.8.5}) independent of the choice of $m$?
\end{description}

In the case $\nbigm^{ss}(\yhat)=\nbigm^s(\yhat)$,
the morphism $\pi$ is $\proj^{H_y(m)-1}$-bundle.
Hence (\ref{eq;06.6.8.5}) is independent of the choice of $m$,
and it is compatible with the ordinary definition.
We will obtain the affirmative answer of $(C)$,
in general.

\subsubsection{Motivation of the study}

The Donaldson invariants for smooth projective surfaces
can be obtained from the moduli spaces of semistable
sheaves of rank $2$.
It is natural to ask what invariants are obtained
from the moduli stacks of semistable sheaves of higher ranks.
That is one of our motivations for the study.

Donaldson invariant has been studied intensively since the 1980s,
motivated by the application to the topology of
real four dimensional manifolds.
Nowadays, it is believed that
the topological information contained in Donaldson invariant
can be obtained from only Seiberg-Witten invariant,
essentially.
Similarly, 
even if we obtain the invariants
from the moduli stacks of the objects with higher ranks,
it is not so reasonable to expect
a new exciting application to topology.

However, it seems interesting to investigate the relation 
among the invariants.
For example, we can ask the following two problems,
which are related with the Kotschick-Morgan conjecture
and the Witten conjecture of Donaldson invariant.

\begin{problem}
\label{problem;06.6.8.10}
Clarify the dependence of $\Phi(\yhat)$
on the polarization $\nbigo_X(1)$.
\hfill\qed
\end{problem}

\begin{problem}
\label{problem;06.6.8.11}
Reduce $\Phi(\yhat)$ to the sum
of the integrals over the products of the moduli spaces 
of the objects with rank one.
\hfill\qed
\end{problem}

As for Problem \ref{problem;06.6.8.10},
we will show that $\Phi(\yhat)$ is independent 
of the choice of $\nbigo_X(1)$ in the case $p_g>0$,
and we will obtain the weak wall crossing formula
in the case $p_g=0$.
As for Problem \ref{problem;06.6.8.11},
we will give such a formula in the case $p_g>0$.

In principle, we can show the existence of the relations
of the invariants as above,
which are universal in some sense.
Moreover, we {\em expect} that
they can be described in terms of good functions
such as modular forms.
The author hopes that this work would,
at least tentatively, provide
a part of the foundation for such an interesting study.

We will explain our main results
in the next subsections.

\subsection{Transition Formulas in the Simple Cases}
\subsubsection{The case where the $2$-stability
 condition is satisfied}
\label{subsubsection;06.6.9.5}

Let $\nbigm(\yhat,[L])$ denote the moduli stack of
the oriented reduced $L$-Bradlow pairs of type $y$.
We have the relative tautological line bundle
 $\nbigo_{\rel}(1)$ on $\nbigm(\yhat,[L])$.
(See the subsubsection \ref{subsubsection;06.4.11.75}
for the definition.)
The restriction to 
$\nbigm^{ss}(\yhat,[L],\delta)$
is also denoted by $\nbigo_{\rel}(1)$.
Let $\omega$ denote the first Chern class of $\nbigo_{\rel}(1)$.
We consider the cohomology classes
which is described as a sum of cohomology classes
of the following form:
\[
 \Phi=P(\Ehat^u)\cdot \omega^k.
\]
If $\delta$ is a non-critical parameter,
i.e. $\nbigm^{ss}(\yhat,[L],\delta)=\nbigm^{s}(\yhat,[L],\delta)$,
we put as follows:
\[
 \Phi(\yhat,[L],\delta):=
 \int_{\nbigm^{s}(\yhat,[L],\delta)} \Phi
\]
Let $\delta$ be critical.
We take parameters $\delta_-<\delta<\delta_+$
such that $|\delta_{\kappa}-\delta|$ are sufficiently small.
We would like to describe the transition
$\Phi(\yhat,[L],\delta_+)-\Phi(\yhat,[L],\delta_-)$
as the sum of the integrals over
the products of the moduli stacks
of the objects with lower ranks.
Such a description is called the transition formula.

If the following condition is satisfied,
the problem is rather simple.
\begin{description}
\item[(2-stability condition)]
 We say that the $2$-stability condition
holds for $(y,L,\delta)$,
if the automorphism group of
$(E,\phi)\in\nbigm^{ss}(y,L,\delta)$
is $\{1\}$ or $G_m$.
\end{description}

To state the theorem,
we need some preparation.
Let $\Type$ denote the set of
cohomology classes of $X$
obtained as the Chern character
of some torsion-free sheaves on $X$.
For any $y\in\Type$,
the $H^0(X)$-part is denoted by $\rank(y)$.
We have the Hilbert polynomial $H_y(t)$ of $y$
which satisfies the following for any integer $m$:
\[
 H_y(m):=\int_X\Td(X)\cdot y\cdot \ch\bigl(\nbigo(m)\bigr).
\]
The reduced Hilbert polynomial $H_y/\rank(y)$
is denoted by $P_y$.
When a parameter $\delta$ is given,
we put $P^{\delta}_y:=(H_y+\delta)/\rank (y)$.
We put as follows:
\[
 S(y,\delta):=\bigl\{
(y_1,y_2)\in\Type^2\,\big|\,
 y_1+y_2=y,\,\,\,
P^{\delta}_y=
 P^{\delta}_{y_1}=P_{y_2}
 \bigr\}
\]

For a given $(y_1,y_2)\in \Type^2$,
we put $r_i=\rank y_i$.
We put as follows:
\[
 \nbigm(y_1,\yhat_2,L,\delta):=
 \nbigm^{ss}(y_1,L,\delta)\times
 \nbigm^{ss}(\yhat_2)
\]
On $\nbigm(y_1,\yhat_2,L,\delta)\times X$,
we have the sheaf $E^u_1$
which is obtained from the universal sheaf
on $\nbigm^{ss}(y_1,L,\delta)\times X$
via the natural projection.
We also have the sheaf $\Ehat_2^u$
which is obtained from the universal sheaf
on $\nbigm^{ss}(\yhat_2)\times X$
via the natural projection.

Let $G_m$ denote the one dimensional torus.
Let $e^{w\cdot t}$ denote the trivial line bundle
with the $G_m$-action of weight $w$.
We have the following element
of the $K$-group of the $G_m$-equivariant
coherent sheaves on
$\nbigm(y_1,\yhat_2,L,\delta)$:
\begin{multline}
 \gbign_0(y_1,y_2)=
-Rp_{X\,\ast}\Bigl(
 \nrhom\bigl(E_1^u\!\cdot\! e^{-t},\,
 \Ehat_2^u\!\cdot\! e^{r_1(t-\omega_1)/r_2}\bigr)
\Bigr)
-Rp_{X\,\ast}\Bigl(
 \nrhom\bigl(\Ehat_2^u\!\cdot\! e^{r_1(t-\omega_1)/r_2},\,
 E^u_1\!\cdot\! e^{-t}\bigr)\Bigr) \\
+Rp_{X\,\ast}\Bigl(
 \nhom\bigl(L\!\cdot\! e^{-t},\,
 \Ehat_2^u\!\cdot\! e^{r_1(t-\omega_1)/r_2}\bigr)
\Bigr)
\end{multline}
Here, we put $\omega_1:=c_1(\Or(E_1^u))/r_1$,
and $e^{w\cdot\omega_1}$ denotes $\Or(E_1^u)^{w/r_1}$
formally.

As a special case of Theorem \ref{thm;06.5.31.100},
we obtain the following theorem.
\begin{thm}
\label{thm;06.6.8.101}
Assume that the $2$-stability condition holds
for $(y,L,\delta)$.
Then, we have the following equality:
\begin{equation}
\label{eq;06.6.8.121}
 \Phi(\yhat,[L],\delta_+)
-\Phi(\yhat,[L],\delta_-)
=\sum_{(y_1,y_2)\in S(y,\delta)}
 \int_{\nbigm(y_1,\yhat_2,L,\delta)}
\underset{t=0}{\res}
\left(
 \frac{P\bigl(E_1^u\cdot e^{-t}\oplus
 \Ehat_2^u\cdot e^{r_1(t-\omega_1)/r_2}\bigr)
 \cdot t^k}
 {\Eu\bigl(\gbign_0(y_1,y_2)\bigr)}
\right)
\end{equation}
In the case $p_g>0$ and $\rank(y_1)>1$,
the contributions from $(y_1,y_2)$ are $0$.
\hfill\qed
\end{thm}

The following condition is called
the $i$-vanishing condition for $(y,L,\delta)$:
\begin{description}
\item[($i$-vanishing condition)]
We have
 $H^j(X,L^{-1}\otimes E)=0$
 for any $j\geq i$ and 
 for any
 $(E,\phi)\in \nbigm^{ss}(y,L,\delta)$.
\end{description}

The problem of the transition is comparatively simple
as in the following proposition.
\begin{prop}
Assume that the $2$-stability condition
and the $2$-vanishing condition
are satisfied for $(y,L,\delta)$.
Then, the transition at $\delta$ is trivial in the case $p_g>0$,
i.e.,
the equality
$\Phi(\yhat,\delta_+)=\Phi(\yhat,\delta_-)$ holds.
\hfill\qed
\end{prop}

If the $1$-vanishing condition holds for $(y,L,\delta)$,
we have the smooth morphism
of $\nbigm^{ss}(\yhat,[L],\delta)$
to the moduli stack $\nbigm(\yhat)$
of the oriented torsion-free sheaves of type $y$.
The relative tangent bundle is denoted by $T_{\rel}$.
We put $N_L(y):=\int_X\Td(X)\cdot y\cdot \ch(L^{-1})$,
which is same as $\rank T_{\rel}+1$.
We will be interested in the following integral:
\[
\Phi_1(\yhat,\delta):=
 \int_{\nbigm^{ss}(\yhat,[L],\delta)}
 \Phi_1,\quad
 \Phi_1:=P(\Ehat^u)\cdot \frac{\Eu(T_{\rel})}{N_L(y)}
\]
For $(y_1,y_2)\in S(y,\delta)$,
we put as follows:
\[
 \nbigm(\yhat_1,\yhat_2,[L],\delta):=
 \nbigm^{ss}(\yhat_1,[L],\delta)
\times
 \nbigm^{ss}(\yhat_2).
\]
On $\nbigm(\yhat_1,\yhat_2,[L],\delta)$,
we have the sheaf $\Ehat^u_1$
which is the pull back of the universal sheaf
on $\nbigm^{ss}(\yhat_1,[L],\delta)\times X$
via the natural projection.
Similarly, we have the sheaf $\Ehat^u_2$
which is the pull back of the universal sheaf
on $\nbigm^{ss}(\yhat_2)\times X$
via the natural projection.
Let $e^{w\cdot s}$ denote the trivial line bundle
with the $G_m$-action of weight $w$.
Let $e^{w\cdot s}$ denote the trivial
line bundle on $\nbigm(\yhat_1,\yhat_2,[L],\delta)$
with the $G_m$-action of weight $w$.
We have the following element of
the $K$-group of the $G_m$-equivariant
coherent sheaves
on $\nbigm(\yhat_1,\yhat_2,[L],\delta)$:
\[
 -Rp_{X\,\ast}\Bigl(
 \nrhom\bigl(\Ehat_1^u\!\cdot\! e^{-s/r_1},\,
 \Ehat_2^u\!\cdot\! e^{s/r_2}\bigr)
\Bigr)
-Rp_{X\,\ast}\Bigl(
 \nrhom\bigl(\Ehat_2^u\!\cdot\! e^{s/r_2},\,
 \Ehat^u_1\!\cdot\! e^{-s/r_1}\bigr)\Bigr) \\
\]
Let $Q\bigl(\Ehat_1\cdot e^{-s/r_1},\Ehat_2\cdot e^{s/r_2}\bigr)$
denote the equivariant Euler class.
As a special case of Theorem \ref{thm;06.6.4.10},
we obtain the following.
\begin{thm}
Assume that the $2$-stability condition
and the $1$-vanishing condition hold
for $(y,L,\delta)$.
In the case $p_g>0$,
we have $\Phi(\yhat,[L],\delta_+)=\Phi(\yhat,[L],\delta_-)$.
In the case $p_g=0$,
we have the following equality:
\[
 \Phi(\yhat,[L],\delta_+)
-\Phi(\yhat,[L],\delta_-)
=\sum_{(y_1,y_2)\in S(y,\delta)} 
 \int_{\nbigm(\yhat_1,\yhat_2,[L],\delta)}
 \Psi(y_1,y_2)
\]
The cohomology classes $\Psi(y_1,y_2)$ are given as follows:
\begin{equation}
 \Psi(y_1,y_2)=
 \frac{N_{L}(y_1)}{N_L(y)}\cdot
\underset{s=0}\Res\left(
 \frac{P\bigl(\Ehat_1^u\cdot e^{-s/r_1}\oplus \Ehat_2^u\cdot e^{s/r_2}\bigr)}
 {Q\bigl(\Ehat_1^u\cdot e^{-s/r_1},\,\Ehat_2^u\cdot e^{s/r_2} \bigr)}
 \right)
\cdot \frac{\Eu(T_{1,\rel})}{N_L(y_1)}
\end{equation}
Here, $T_{1,\rel}$ denote the relative tangent bundle
of the smooth morphism
$\nbigm^{ss}(\yhat_1,[L],\delta)
\lrarr \nbigm(\yhat_1)$.
\hfill\qed
\end{thm}

\subsubsection{Reduced $\vecL$-Bradlow pair}

Let $\vecL=(L_1,L_2)$ be a pair of line bundles on $X$.
Let $U$ be a scheme,
and let $E$ be a $U$-torsion-free sheaf on $U\times X$
which is flat over $U$.
Let $[\vecphi]$ be a pair $([\phi_1],[\phi_2])$
of reduced $L_i$-sections $[\phi_i]$ of $E$
such that $[\phi_i]_{|\{u\}\times X}\neq 0$
for any $u\in U$.
Such a pair $(E,[\vecphi])$ is called
a reduced $[\vecL]$-Bradlow pair.
Let $\vecdelta=(\delta_1,\delta_2)$
be an element of $\nbigp^{\br\,2}$.
We naturally have 
the notion of $\vecdelta$-semistability and $\vecdelta$-stability
for reduced $\vecL$-Bradlow pairs.
Let $\nbigm^{ss}(\yhat,[\vecL],\vecdelta)$
(resp. $\nbigm^s(\yhat,[\vecL],\vecdelta)$)
denote the moduli stack of 
$\vecdelta$-semistable (resp. $\vecdelta$-stable)
$\vecL$-Bradlow pairs of type $y$.
When $\delta_i$ are sufficiently small,
we have the morphism
$\nbigm^{ss}(\yhat,[\vecL],\vecdelta)
\lrarr \nbigm^{ss}(\yhat)$.

Let us move $\delta_1$ by fixing $\delta_2$.
We say that $\delta_1$ is critical,
if $\nbigm^{ss}(\yhat,[\vecL],\vecdelta)\neq 
 \nbigm^{s}(\yhat,[\vecL],\vecdelta)$ holds.
As in the subsubsection \ref{subsubsection;06.6.9.5},
we have the notion of $2$-stability condition
for $(y,\vecL,\vecdelta)$ (Definition \ref{df;06.6.8.15}).
When $\delta_i$ are sufficiently small,
it can be shown that the $2$-stability condition
is always satisfied even if $\delta_1$ is critical,
and we have a good transition formula.

Assume that $\delta_1$ is critical.
We take elements $\delta_-,\delta_+\in\nbigp^{\br}$
such that $\delta_-<\delta_1<\delta_+$
and that $|\delta_{\kappa}-\delta_1|$ $(\kappa=\pm)$
are sufficiently small.
We put $\vecdelta_{\kappa}=(\delta_{\kappa},\delta_2)$
for $\kappa=\pm$.
Let $T_{\rel}^{(1)}$ denote the relative tangent 
bundle of the smooth morphism
of $\nbigm^{ss}(\yhat,[\vecL],\vecdelta_{\kappa})$
to the moduli stack $\nbigm(\yhat,[L_2])$
of reduced $L_2$-Bradlow pairs.
We consider the following cohomology class:
\[
 \Phi=\frac{\Eu(T_{\rel}^{(1)})}{N_{L_1}(y)}
\cdot P(\Ehat^u)(\omega^{(2)})^k
\]
Here, $\omega^{(2)}$ denote the first Chern class
of the line bundle which is the pull back of
the relative tautological line bundle
of $\nbigm(\yhat,[L_2])$.
We put as follows, for $\kappa=\pm$:
\[
 \Phi(\yhat,[\vecL],\vecdelta_{\kappa}):=
\int_{\nbigm^{ss}(\yhat,[\vecL],\vecdelta_{\kappa})}
 \Phi
\]

We would like to discuss the transition formula
for $\Phi(\yhat,[\vecL],\vecdelta_+)-\Phi(\yhat,[\vecL],\vecdelta_-)$.
We put as follows:
\[
 S(y,\vecdelta):=
 \bigl\{
 (y_1,y_2)\in\Type^2\,\big|\,
 P_{y_1}=P_{y_2},\,\,\,
 \delta_1/\rank(y_1)=\delta_2/\rank(y_2)
 \bigr\}
\]
For any $(y_1,y_2)\in S(y,\vecdelta)$,
we put as follows:
\[
 \nbigm(\yhat_1,\yhat_2,[\vecL],\vecdelta)
:=\nbigm^{ss}(\yhat_1,[L_1],\delta_1)
\times
 \nbigm^{ss}(\yhat_2,[L_2],\delta_2)
\]
Let $\Ehat_i$ denote the sheaf on 
$\nbigm(\yhat_1,\yhat_2,[\vecL],\vecdelta)$
which is the pull back of the universal sheaf
over $\nbigm^{ss}(\yhat_i,[L_i],\delta_i)\times X$.
Let $\nbigo_{i,\rel}(1)$ denote the pull back
of the relative
tautological line bundles on
$\nbigm^{ss}(\yhat_i,[L_i],\delta_i)$.
We put $\omega_i:=c_1\bigl(\nbigo_{i,\rel}(1)\bigr)$.
Let $T_{1,\rel}$ denote the bundle
on $\nbigm(\yhat_1,\yhat_2,[\vecL],\vecdelta)$
induced by the relative tangent bundle
of the smooth morphism
$\nbigm^{ss}(\yhat_i,[L_i],\delta_i)
\lrarr \nbigm(\yhat_i)$.

Let $e^{w\cdot s}$ denote the trivial line bundle
with the $G_m$-action of weight $w$
as above.
We have the following element of the $K$-group
$K^{G_m}\bigl(\nbigm(\yhat_1,\yhat_2,[\vecL],
 \vecdelta)\bigr)$
of $G_m$-equivariant coherent sheaves
on $\nbigm(\yhat_1,\yhat_2,[\vecL],\vecdelta)$:
\[
  -Rp_{X\,\ast}\Bigl(
 \nrhom\bigl(\Ehat_1^u\!\cdot\! e^{-s/r_1},\,
 \Ehat_2^u\!\cdot\! e^{s/r_2}\bigr)
\Bigr)
-Rp_{X\,\ast}\Bigl(
 \nrhom\bigl(\Ehat_2^u\!\cdot\! e^{s/r_2},\,
 \Ehat^u_1\!\cdot\! e^{-s/r_1}\bigr)\Bigr) \\
\]
The equivariant Euler class is denoted by
$Q(\Ehat_1^u\cdot e^{s/r_1},\Ehat_2^u\cdot e^{s/r_2})$.
We also have the following element
of $K^{G_m}(\nbigm(\yhat_1,\yhat_2,[\vecL],\delta))$:
\[
 Rp_{X\,\ast}\nhom(L_2,\Ehat_1^u)\cdot e^{-s/r_1-s/r_2+\omega_2}
\]
The equivariant Euler class
is denoted by
$R(L_2\cdot e^{-\omega_2+s/r_2},\,
 \Ehat_1^u\cdot e^{-s/r_1})$.
We obtain the following proposition
as the special case of Proposition \ref{prop;06.6.4.15}
and Proposition \ref{prop;06.5.31.30}.
\begin{prop}
\label{prop;06.6.8.40}
In the case $p_g>0$,
we have the equality
$ \Phi(\yhat,[\vecL],\vecdelta_+)=\Phi(\yhat,[\vecL],\vecdelta_-)$.
In the case $p_g=0$,
the following equality holds:
\begin{equation}
\label{eq;06.6.8.35}
 \Phi(\yhat,[\vecL],\vecdelta_+)
-\Phi(\yhat,[\vecL],\vecdelta_-)
=\sum_{(y_1,y_2)\in S(y,\vecdelta)}
 \int_{\nbigm(\yhat_1,\yhat_2,[\vecL],\delta)}
 \Psi(y_1,y_2)
\end{equation}
Here, $\Psi(y_1,y_2)$ are given as follows:
\begin{equation}
 \Psi(y_1,y_2):=
 \frac{N_{L_1}(y_1)}{N_{L_1}(y)}\cdot
 \underset{s=0}\Res
\left(
 \frac{P(\Ehat_1^u\!\cdot\! e^{-s/r_1}\oplus 
 \Ehat_2^u\!\cdot\! e^{s/r_2})
 \cdot(\omega_2-s/r_2)^k}
 {Q(\Ehat_1^u\cdot e^{-s/r_1},\Ehat_2^u\cdot e^{s/r_2})
\cdot R(L_2\!\cdot\! e^{-\omega_2+s/r_2},
 \Ehat_1^u\!\cdot\! e^{-s/r_1})}
\right)
 \cdot \frac{\Eu(T_{1,\rel})}{N_{L_1}(y_1)} 
\end{equation}
\hfill\qed
\end{prop}

As a special case,
let us consider the integral of
the following cohomology class,
assuming that the $1$-vanishing condition
holds for $(y,L_2)$:
\begin{equation}
\label{eq;06.6.8.25}
 \Phi=\frac{\Eu(T^{(1)}_{\rel})}{N_{L_1}(y)}
 \cdot \frac{\Eu(T^{(2)}_{\rel})}{N_{L_2}(y)}
 \cdot P(\Ehat^u)
\end{equation}

\begin{lem}
\label{lem;06.6.24.1}
Assume that the $1$-vanishing condition holds
for $(y,L_2)$.
Let $\Phi$ be as in {\rm(\ref{eq;06.6.8.25})}.
Then, the cohomology class $\Psi(y_1,y_2)$
in {\rm(\ref{eq;06.6.8.35})} is given as follows:
\[
 \underset{s=0}\Res
 \left(
  \frac{P(\Ehat_1^u\cdot e^{-s/r_1}\oplus \Ehat_2^u\cdot e^{s/r_2})}
 {Q(\Ehat_1^u\cdot e^{-s/r_1},\Ehat_2^u\cdot e^{s/r_2})}
\right)
 \cdot \frac{\Eu(T_{1,\rel})}{N_{L_1}(y_1)}
 \cdot \frac{\Eu(T_{2,\rel})}{N_{L_2}(y_2)}
\]
\hfill\qed
\end{lem}

\subsubsection{Well-definedness of (\ref{eq;06.6.8.5})}

In the case $p_g>0$,
we can show the well-definedness of (\ref{eq;06.6.8.5})
by using Proposition \ref{prop;06.6.8.40}.
Assume that $L_1^{-1}$ is ample,
and that the $1$-vanishing condition holds
for $(y,L_2)$.
If we take a sufficiently large integer $m$,
the $1$-vanishing condition holds
for $(y,L_1^{m})$.
We put $\vecL^{(m)}:=(L_1^{m},L_2)$.
Let $T_{\rel}^{(1)}$ denote the relative tangent bundle
of the smooth map
$\nbigm^{ss}(\yhat,[\vecL^{(m)}],\vecdelta)
\lrarr \nbigm(\yhat,[L_2])$.
We use the notation $T_{\rel}^{(2)}$ in a similar meaning.
We consider the following number:
\[
 g(L_1^{m},L_2,\delta_1,\delta_2):=
\int_{\nbigm^{ss}(\yhat,[\vecL^{(m)}],\vecdelta)}
 P(\Ehat^u)\cdot
 \frac{\Eu(T^{(1)}_{\rel})}{N_{L_1^m}(y)}
 \cdot
 \frac{\Eu(T^{(2)}_{\rel})}{N_{L_2}(y)}
\]
We assume that both of $\delta_i$ are sufficiently small.
When $\delta_1$ is sufficiently smaller than $\delta_2$,
we have the following:
\[
 g(L_1^{m},L_2,\delta_1,\delta_2)=
 \int_{\nbigm^{ss}(\yhat,[L_2],\delta_2)}
 \Phi\cdot \frac{\Eu(T^{(2)}_{\rel})}{N_{L_2}(y)}
\]
When $\delta_2$ is sufficiently smaller than $\delta_1$,
we have the following:
\[
  g(L_1^{m},L_2,\delta_1,\delta_2)=
 \int_{\nbigm^{ss}(\yhat,[L_1^m],\delta_1)}
 \Phi\cdot \frac{\Eu(T^{(1)}_{\rel})}{N_{L_1^m}(y)}
\]
When we move $\delta_1$,
the transitions are trivial,
due to Proposition \ref{prop;06.6.8.40}.
Therefore, we obtain the following:
\[
 \int_{\nbigm^{ss}(\yhat,[L_1^m],\delta_1)}
 \Phi\cdot \frac{\Eu(T^{(1)}_{\rel})}{N_{L_1^m}(y)}
=
 \int_{\nbigm^{ss}(\yhat,[L_2],\delta_2)}
 \Phi\cdot \frac{\Eu(T^{(2)}_{\rel})}{N_{L_2}(y)}
\]
In particular, we obtain that (\ref{eq;06.6.8.5})
is independent of the choice of $m$.
Moreover, we can show the following equality,
assuming the $2$-vanishing condition for $(y,L)$
$d:=\int_X\Td(X)\cdot y\cdot\ch(L^{-1})-1\geq 0$:
\begin{equation}
 \Phi(\yhat)=\int_{\nbigm^{ss}(\yhat,[L],\delta)} \Phi\cdot \omega^d 
\end{equation}

In the case $p_g=0$,
the problem is more subtle.
We can derive it from Proposition \ref{prop;06.6.8.40}
and Lemma \ref{lem;06.6.24.1}
that (\ref{eq;06.6.8.5}) is independent
of the choice of $m$ in this case, too.
However, we need some more additional argument.

\begin{rem}
Although we do not discuss the parabolic structure
in this section,
we will consider the invariants 
obtained from the moduli stacks
of the oriented parabolic torsion-free sheaves
and 
the oriented parabolic reduced $L$-Bradlow pairs.
And, it is not clear whether {\rm (\ref{eq;06.6.8.5})}
is independent of the choice of $m$, in general.
Instead, we can show the existence of the following limit,
for a line bundle $L$ such that $L^{-1}$ is ample:
\begin{equation}
 \label{eq;06.6.8.100}
 \lim_{m\to\infty}
 \int_{\nbigm^{ss}(\vecyhat,[L^m],\alpha_{\ast},\delta)}
 \Phi\cdot \frac{\Eu(T_{\rel})}{N_{L^m}(y)}
\end{equation}
Here, $\vecy$ denotes a type of parabolic sheaves,
and $\alpha_{\ast}$ denotes a system of weight.
Moreover, the limit is independent of the choice of $L$.
Hence, we may adopt {\rm(\ref{eq;06.6.8.100})}
as the definition of  $\Phi(\vecyhat,\alpha_{\ast})$.
\hfill\qed
\end{rem}

\subsection{Rank $2$ Case}
\subsubsection{Dependence on the polarizations}
\label{subsubsection;06.6.8.151}

In the case $\rank(y)=2$,
the $2$-stability condition is always satisfied.
Hence, Theorem \ref{thm;06.6.8.101} provides
us a tool to discuss the problems \ref{problem;06.6.8.10}
and \ref{problem;06.6.8.11}.
We explain our result for Problem \ref{problem;06.6.8.10}
in this subsubsection.

To distinguish the dependence on the polarization $H$,
we use the notation
$\nbigm_H(\yhat)$ to denote the moduli stack
of torsion-free sheaves of type $y$
which are semistable with respect to $H$.
We also assume $a^2-4n\leq \xi^2<0$.
We put $W^{\xi}:=\bigl\{c\in
NS(X)\otimes\real\,\big|\,(c,\xi)=0\bigr\}$,
which is called the wall determined by $\xi$.
It is well known that
the ample cone is divided into the chambers
by such walls,
and the moduli $\nbigm_H(\yhat)$ depends
only on the chambers to which $H$ belongs.
We put as follows:
\[
 \Phi_H(\yhat):=\int_{\nbigm_H(\yhat)}\Phi
\]
We would like to discuss how $\Phi_H(\yhat)$ changes
when the polarizations vary across the wall $W^{\xi}$.

Let $a$ and $n$ be the first and second Chern classes of $y$.
We put as follows:
\[
 S_0(y,\xi):=\bigl\{
 (a_0,a_1)\in NS(X)^2\,\big|\,
 a_0-a_1=m\cdot\xi\,\,(m>0)
 \bigr\}
\]
For any $(a_0,a_1)\in S_0(y,\xi)$,
we put as follows:
\[
 X(a_0,a_1):=\coprod_{n_0+n_1=n-a_0\cdot a_1}X^{[n_0]}\times X^{[n_1]}
\]
On $X^{[n_0]}\times X^{[n_1]}\times X$,
we have the sheaf $\nbigi_i^u$, which is the pull back
of the universal ideal sheaf over $X^{[n_i]}\times X$
via the natural projection.
Let $e^{a_i}$ denote the holomorphic line bundle on $X$
whose first Chern class is $a_i$.
It is uniquely determined up to isomorphism,
since we have assumed that $X$ is simply connected
in {\em Introduction}.
Let $Q\bigl(\nbigi_0^ue^{a_0-s},\nbigi_1^ue^{a_1+s}\bigr)$
be the equivariant Euler class of the following
element of the $K$-group of
the $G_m$-equivariant coherent sheaves on
$X^{[n_0]}\times X^{[n_1]}$:
\[
-Rp_{X\,\ast}\left(
 \nrhom\bigl(\nbigi_0^u\cdot e^{a_0-s},\,\,
 \nbigi_1^u\cdot e^{a_1+s}\bigr)
 \right)
-Rp_{X\,\ast}\left(
 \nrhom\bigl(\nbigi_1^u\cdot e^{a_1+s},\,\,
 \nbigi_0^u\cdot e^{a_0-s}\bigr)
 \right)
\]

\begin{thm}
Let $C_+$ and $C_-$ be chambers
which are divided by the wall $W^{\xi}$.
Let $H_+$ and $H_-$ be ample line bundles
contained in $C_+$ and $C_-$, respectively.
We assume $(H_-,\xi)<0<(H_+,\xi)$.
\begin{itemize}
\item
 In the case $p_g>0$,
 we have $\Phi_{H_+}(\yhat)=\Phi_{H_-}(\yhat)$.
 Namely, the invariant does not depend
 on the choice of generic polarizations.
\item
 In the case $p_g=0$,
 we have the following equality:
\begin{equation}
 \label{eq;06.6.8.111}
 \Phi_{H_+}(\yhat)-\Phi_{H_-}(\yhat)
=\sum_{(a_0,a_1)\in S_0(y,\xi)}
 \int_{X(a_0,a_1)}
\underset{s=0}\Res
 \left(
 \frac{P\bigl(\nbigi_0^u\cdot e^{a_0-s}\oplus
 \nbigi_1^u\cdot e^{a_1+s}\bigr)}
 {Q\bigl(\nbigi_0^u\cdot e^{a_0-s},\,\nbigi^u_1\cdot e^{a_1+s}\bigr)}
 \right)
\end{equation}
We call {\rm(\ref{eq;06.6.8.111})} the weak wall crossing formula.
\hfill\qed
\end{itemize}
\end{thm}

\begin{rem}
Under the assumption that the wall
$W^{\xi}$ is good,
the weak wall crossing formula was proved
for the Donaldson invariant 
in {\rm\cite{eg}} and {\rm\cite{fq}},
which was refined in {\rm\cite{gny}}.
\hfill\qed
\end{rem}

\begin{rem}
Remarkably, 
L. G\"ottsche-H. Nakajima-K. Yoshioka
established the way
to derive the wall crossing formula
from the weak wall crossing formula.
(See {\rm\cite{gny}}.)
\hfill\qed
\end{rem}

\begin{rem}
K. Yamada proved the independence of the invariants
from the polarizations in some cases.
\hfill\qed
\end{rem}

\subsubsection{Reduction to the integrals over Hilbert schemes}
\label{subsubsection;06.6.8.201}

Let us discuss the problem \ref{problem;06.6.8.11}
in the case $\rank(y)=2$.
Let $a$ and $n$ denote the first and second Chern classes  of $y$.
We also assume $p_g>0$.
For any element
$a_1$ of the Neron-Severi group $NS^1(X)$,
we put $a_2:=a-a_1$.
Let $e^{a_i}$ denote the holomorphic line bundle
whose first Chern class is $a_i$.
Let $\nbigi_i^u$ denote the universal
ideal sheaves over $X^{[n_i]}\times X$.
Let $\nbigz_i$ denote the universal $0$-scheme
over $X^{[n_i]}\times X$.
Let $\Xi_i$ denote $p_{X\,\ast}\bigl(
\nbigo_{\nbigz_i}\otimes e^{a_i}\bigr)$.
We use the same notation to denote
the pull back of them via  appropriate morphisms.
In the case 
$\bigl(c_1(\nbigo(1)),a_1\bigr)<
\bigl(c_1(\nbigo(1)),a_2\bigr)$,
we put as follows:
\[
 \nbiga(a_1,y):=\sum_{n_1+n_2=n-a_1\cdot a_2}
\int_{X^{[n_1]}\times X^{[n_2]}}
 \underset{s=0}\Res
\left(
 \frac{P\bigl(\nbigi_1^u\cdot e^{a_1-s}\oplus 
 \nbigi_2^u\cdot e^{a_2+s}\bigr)}
 {Q\bigl(\nbigi_1^u\cdot e^{a_1-s},
 \nbigi_2^u\cdot e^{a_2+s}\bigr)}
 \cdot \frac{\Eu(\Xi_1)\cdot \Eu(\Xi_2\cdot e^{2s})}
 {(2s)^{n_1+n_2-p_g}}
\right)
\]
In the case
$\bigl(c_1(\nbigo(1)),a_1\bigr)
=\bigl(c_1(\nbigo(1)),a_2\bigr)$,
we put as follows:
\[
 \nbiga(a_1,y):=\sum_{\substack{n_1+n_2=n-a_1\cdot a_2\\ n_1>n_2}}
\int_{X^{[n_1]}\times X^{[n_2]}}
 \underset{s=0}\Res
\left(
 \frac{P\bigl(\nbigi_1^u\cdot e^{a_1-s}\oplus 
 \nbigi_2^u\cdot e^{a_2+s}\bigr)}
 {Q\bigl(\nbigi_1^u\cdot e^{a_1-s},
 \nbigi_2^u\cdot e^{a_2+s}\bigr)}
 \cdot \frac{\Eu(\Xi_1)\cdot \Eu(\Xi_2\cdot e^{2s})}
 {(2s)^{n_1+n_2-p_g}}
\right)
\]

Recall that an abelian pair is defined to be a pair
of a holomorphic line bundle $\nbigl$ and a section
$\phi:\nbigo\lrarr \nbigl$.
Let $M(c)$ denote the moduli of
abelian pairs $(\nbigl,\phi)$ such that
$c_1(\nbigl)=c$.
We can show the following proposition.
\begin{prop}
[Proposition \ref{prop;06.5.13.70}]
Assume $H^1(X,\nbigo_X)=0$
and $p_g=\dim H^2(X,\nbigo_X)> 0$.
Moreover,
we assume that the virtual fundamental class
of $M(c)$ is not $0$.
Then, the expected dimension of $M(c)$ is $0$.

We can regard the virtual fundamental class 
$[M(c)]$ as the number,
and then it is same as the following:
\[
 \SW(c):=\frac{\prod_{i=1}^{d(c)}(i-p_g)}{d(c)!}
\]
Here we put $d(c):=\dim M(c)=\dim H^0(X,e^c)-1$.
\hfill\qed
\end{prop}

Then, we put as follows:
\[
 \nSW(X,y):=\bigl\{
 a_1\in NS(X)
 \,\big|\,
 \SW(a_1)\neq 0 ,\,\,
 \bigl(a_1,c_1(\nbigo_X(1))\bigr)\leq
 \bigl(a,c_1(\nbigo_X(1))\bigr)/2
 \bigr\}
\]

\begin{thm}
[Theorem \ref{thm;06.6.6.10}]
\label{thm;06.6.8.200}
Assume $p_g>0$ and $H^1(X,\nbigo)=0$.
Assume $P_y>P_K$ and $\chi(y)-1\geq 0$,
where $K$ denotes the canonical line bundle of $X$,
and we put $\chi(y)=\int\Td(X)\cdot y$
for the Todd class $\Td(X)$.
Then, we have the following equality:
\[
 \int_{\nbigm(\yhat)}P(\Ehat^u)
+\sum_{a_1\in \nSW(X,y)}
\SW(a_1)\cdot 2^{1-\chi(y)}\cdot \nbiga(a_1,y)=0.
\]
\hfill\qed
\end{thm}

\subsection{Higher Rank Case}
\subsubsection{The case $p_g>0$}

If $p_g>0$ is satisfied,
the results in the rank $2$ case
can be rather easily generalized in the higher rank case.
Actually,
we have the formally same formula as (\ref{eq;06.6.8.121}).
\begin{thm}
 [Theorem \ref{thm;06.6.6.20}]
Assume $p_g>0$.
Then the following equality holds:
\begin{equation}
  \label{eq;06.6.8.130}
 \Phi(\yhat,[L],\delta_+)
-\Phi(\yhat,[L],\delta_-)
=\sum_{(y_1,y_2)\in S_1(y,\delta)}
 \int_{\nbigm(y_1,\yhat_2,L)}
\underset{t=0}{\res}
\left(
 \frac{P\bigl(E_1^u\cdot e^{-t}\oplus
 \Ehat_2^u\cdot e^{r_1(t-\omega_1)/r_2}\bigr)
 \cdot t^k}
 {\Eu\bigl(\gbign_0(y_1,y_2)\bigr)}
\right)
\end{equation}
Here, we put
$S_1(y,\delta):=
 \bigl\{(y_1,y_2)\in S(y,\delta)\,\big|\,\rank (y_1)=1\bigr\}$.
We use the notation $\nbigm(y_1,\yhat_2,L)$
instead of $\nbigm(y_1,\yhat_2,L,\delta)$,
because $\delta$-semistability condition is trivial
in the case $\rank(y_1)=1$.
\hfill\qed
\end{thm}

As for Problem \ref{problem;06.6.8.10},
we can show the independence of the invariant
from the polarization,
by using the formula (\ref{eq;06.6.8.130}).
Let $\Phi_H(\yhat)$ be as in the subsubsection
\ref{subsubsection;06.6.8.151}.

\begin{thm}
[Theorem \ref{thm;06.6.8.150}]
The invariant $\Phi_H(\yhat)$
is independent of the choice of a generic polarization
in the case $p_g>0$.
\hfill\qed
\end{thm}

As for Problem \ref{problem;06.6.8.11},
we obtain an immediate generalization of
Theorem \ref{thm;06.6.8.200}.
We do not reproduce it here.
See the subsubsection \ref{subsubsection;06.6.30.2}
for detail.

\subsubsection{Transition formula in the case $p_g=0$}

The transition formula is comparatively complicated
in the case $p_g=0$.
We restrict ourselves to the case
where the $1$-vanishing condition holds
for $(y,L,\delta)$,
and we discuss the integral of the cohomology class
of the following form:
\[
 \Phi= P(\Ehat^u)\cdot\frac{\Eu(T_{\rel})}{N_L(y)}
\]
We put as follows:
\[
 \Phi(\delta):=\int_{\nbigm^{ss}(\yhat,[L],\delta)}\Phi.
\]

For each positive integer $k$
we put as follows:
\[
 S_k(y,\delta):=
 \bigl\{
 \vecY=(y_1,\ldots,y_k)\in\Type^k\,\big|\,
 P_{y_i}
=P^{\delta}_{y}
 \bigr\}
\]
For each element
$\vecY=(y_1,\ldots,y_k)\in S_k(y,\delta)$,
we put $|\vecY|=\sum_{i=1}^ky_i$.
We also put as follows:
\[
 W(\vecY):=\prod_{i=1}^k
 \frac{\rank(y_i)}
 {\sum_{1\leq j\leq i}\rank(y_j)}
\]
We put as follows:
\[
 \Sbar(y,\delta):=
 \coprod_{k}\Sbar_k(y,\delta),
\quad\quad
 \Sbar_k(y,\delta):=
 \bigl\{
 (y_0,\vecY)\in\Type\times S_k(y,\delta)
 \,\big|\,y_0+|\vecY|=y
 \bigr\}
\]
For any $(y_0,\vecY)\in \Sbar(y,\delta)$,
we put as follows:
\[
 \nbigm(\yhat_0,\vecYhat,[L]):=
 \nbigm^{ss}\bigl(\yhat_0,[L],\delta_-\bigr)
\times
 \prod_{i=1}^k\nbigm^{ss}(\yhat_i)
\]
Let $\Ehat_0^u$ denote the sheaf
over $\nbigm(\yhat_0,\vecYhat,[L])\times X$
which is obtained as the pull back of
the universal sheaf over
$\nbigm(\yhat_0,[L],\delta_-)\times X$
via the natural projection.
We use the notation $\Ehat_i^u$ in a similar meaning.

When $(y_0,\vecY)\in \Sbar(y,\delta)$ is given,
we put as follows:
\[
 T_0=-\sum_{j>0}\frac{t_j}{\sum_{0\leq h< j}\rank(y_h)},
\quad
 T_i=-\sum_{j>i}\frac{t_j}{\sum_{0\leq h<j}\rank(y_h)}
+\frac{t_i}{\rank(y_i)}
\]
Here, $t_1,\ldots,t_k$ are the variables.
Let $G$ denote the $k$-dimensional torus
$\Spec k[\tau_1,\tau_1^{-1},\ldots,\tau_k,\tau_k^{-1}]$.
Let $e^{w\cdot t_i}$ denote the trivial line bundle
with $G$-action which is induced by the action of
$\Spec k[\tau_i,\tau_i^{-1}]$ of weight $w$.
We have the following element of 
the $K$-group
of the $G$-equivariant coherent sheaves
over $\nbigm(\yhat_0,\vecYhat,[L])$:
\[
 -Rp_{X\,\ast}\nrhom\bigl(\Ehat_i\cdot e^{T_i},\,\Ehat_j\cdot e^{T_j}\bigr)
-Rp_{X\,\ast}\nrhom\bigl(\Ehat_j\cdot e^{T_j},\,\Ehat_i\cdot e^{T_i}\bigr)
\]
The equivariant Euler class is denoted by
$Q\bigl(\Ehat_i\cdot e^{T_i},\Ehat_j\cdot e^{T_j}\bigr)$.
We put as follows:
\[
 Q\bigl(\Ehat_0\cdot e^{T_0},\Ehat_1\cdot e^{T_1},\ldots,
 \Ehat_k\cdot e^{T_k}\bigr):=
 \prod_{i<j}
 Q\bigl(\Ehat_i\cdot e^{T_i},\Ehat_j\cdot e^{T_j}\bigr)
\]
Let $T_{0\,\rel}$ denote the vector bundle
over $\nbigm(\yhat_0,\vecYhat,[L])$
obtained from the relative tangent bundle
of the smooth map 
$\nbigm(\yhat_0,[L],\delta)\lrarr 
 \nbigm(\yhat_0,[L])$.
Then, we put as follows:
\[
 \Psi(y_0,\vecY):=
 \underset{t_1=0}\Res\cdots
 \underset{t_k=0}\Res
\left(
\frac{P\bigl(\bigoplus_{i=0}^k\Ehat_i^u\cdot e^{T_i}\bigr)}
 {Q\bigl(\Ehat_0\cdot e^{T_0},\ldots,\Ehat_k\cdot e^{T_k}\bigr)}
\right)
 \cdot \frac{\Eu(T_{0,\rel})}{N_L(y_0)}
\]

In a sense,
much part of this paper is devoted to
the proof of the following theorem.

\begin{thm}
\label{thm;06.6.9.100}
We have the following formula:
\begin{equation}
\label{eq;06.6.8.250}
 \Phi(\delta_+)-\Phi(\delta_-)
=\sum_{(y_0,\vecY)\in \Sbar(y,\delta)}
 \frac{N_L(y_0)}{N_L(y)}\cdot
 W(\vecY)\cdot
 \int_{\nbigm(\yhat_0,\vecYhat,[L])}
 \Psi(y_0,\vecY)
\end{equation}
\hfill\qed
\end{thm}

\subsubsection{Weak Intersection rounding formula}
As for the problem \ref{problem;06.6.8.10},
we easily obtain a generalization of the weak wall crossing
formula
by using (\ref{eq;06.6.8.250}).
(See Theorem \ref{thm;06.6.7.100}.)
We do not reproduce it here.
The formula itself is not so easy to deal with,
partially because it contains the integrals
over the moduli stacks of semistable sheaves with higher ranks.
To derive a more accessible quantity from our invariants,
we consider the ``intersection rounding formula''.
The general case will be discussed in the subsection
\ref{subsection;06.6.7.400}.
In this subsubsection, we reproduce the result in the rank $3$ case.
See the subsubsection \ref{subsubsection;06.6.8.300}
for more detail.

We take an element $\vecxi=(\xi_1,\xi_2)\in NS^1(X)^2$
such that $\xi_1$ and $\xi_2$ are linearly independent,
and we put $W^{\vecxi}:=W^{\xi_1}\cap W^{\xi_2}$.
A connected component $T$ of
$ W^{\vecxi}\setminus
 \bigcup_{W\neq W^{\xi_i}}W$
is called a tile.
For each tile $T$,
there exist four chambers 
$C_{++}$, $C_{+-}$, $C_{--}$ and $C_{-+}$
with the following properties:
\begin{itemize}
\item
 The closure of $C_{\kappa_1,\kappa_2}$ contains $T$.
\item
 Take an ample line bundle $H_{\kappa_1,\kappa_2}\in 
 C_{\kappa_1,\kappa_2}$.
 Then the signature of 
 the pairing $(H_{\kappa_1,\kappa_2},\xi_i)$ is $\kappa_i$
\end{itemize}
Now, we put as follows:
\[
 \nbigd_{\vecxi}^T\Phi(\yhat):=
 \Phi_{H_{++}}(\yhat)
-\Phi_{H_{+-}}(\yhat)
-\Phi_{H_{-+}}(\yhat)
+\Phi_{H_{--}}(\yhat).
\]
We would like to express
$\nbigd_{\vecxi}^T\Phi(\yhat)$
as the sum of the integrals
over the products of Hilbert schemes.

Let $S(2,1)$ be the set of
$\veca=(a_0,a_1,a_2)\in NS^1(X)^3$
with the following property:
\begin{itemize}
\item
 $a_0+a_1-2a_2=A_1\cdot \xi_1$
 and $a_0-a_1=A_2\cdot \xi_2$
 for some $A_i>0$.
\end{itemize}
Let $S(1,2)$ be the set of
$\veca=(a_0,a_1,a_2)\in NS^1(X)^3$
with the following property:
\begin{itemize}
\item
 $2a_0-(a_1+a_2)=A_1\cdot \xi_1$
 and $a_1-a_2=A_2\cdot \xi_2$
 for some $A_i>0$.
\end{itemize}
For each $\veca$,
we put as follows:
\[
 X(y,\veca):=
 \coprod_{n_0+n_1+n_2=N(y,\veca)}
 \prod_{i=0}^2 X^{[n_i]},
\quad
 N(y,\veca)=n+\frac{a_0^2+a_1^2+a_2^2-a^2}{2}
\]
Here, $a$ and $n$ denote the first Chern class
and the second Chern class of $y$,
respectively.

\begin{prop}
\mbox{{}}
\begin{itemize}
\item
 $\nbigd_{\vecxi}^T\Phi(\yhat)$ is independent of the choice of 
 a tile $T$.
 Therefore, we can omit to denote $T$.
\item
 The following equality holds:
\begin{multline}
\label{eq;06.6.8.400}
 \nbigd_{\vecxi}\Phi(\yhat)=
 \sum_{\veca\in S(1,2)}
 \int_{X(y,\veca)}
 \underset{t_2}\Res\,
 \underset{t_1}\Res\left(
 \frac{P\bigl(\nbigi^u_0 e^{a_0-t_1}\oplus
 \nbigi^u_1e^{a_1+t_1/2-t_2}\oplus
 \nbigi^u_2e^{a_2+t_1/2+t_2}\bigr)}
 {Q\bigl(
 \nbigi^u_0 e^{a_0-t_1},\,\,
 \nbigi^u_1e^{a_1+t_1/2-t_2},\,\,
 \nbigi^u_2e^{a_2+t_1/2+t_2}
 \bigr) }
 \right)
 \\
+\sum_{\veca\in S(2,1)}
 \int_{X(y,\veca)}
 \underset{t_2}\Res\,
 \underset{t_1}\Res\left(
 \frac{P\bigl(\nbigi_0^ue^{a_0-t_1/2-t_2}\oplus
 \nbigi^u_1e^{a_1-t_1/2+t_2}\oplus
 \nbigi^u_2e^{a_2+t_1} \bigr)}
 {Q\bigl(\nbigi_0^ue^{a_0-t_1/2-t_2},\,\,
 \nbigi^u_1e^{a_1-t_1/2+t_2},\,\,\,
 \nbigi^u_2e^{a_2+t_1}\bigr)}
 \right)
\end{multline}
\end{itemize}
The formula {\rm(\ref{eq;06.6.8.400})}
is called
the weak intersection rounding formula
in the rank $3$ case.
See Theorem {\rm\ref{thm;06.6.7.200}}
for the general case.
\hfill\qed
\end{prop}

\begin{rem}
The author expects that $\nbigd_{\vecxi}\Phi(\yhat)$ can be described
more beautifully like the wall crossing formula
in the rank $2$ case
({\rm\cite{gottsche}}, {\rm\cite{gny}}).
It may be interesting to see the equality
obtained by exchanging the roles of $\xi_1$ and $\xi_2$,
\hfill\qed
\end{rem}

\subsection{Master Space}
\subsubsection{The master space of Thaddeus}
\label{subsubsection;06.6.9.3}

In this subsubsection,
we recall the picture of the master space given by M. Thaddeus,
which is one of the most important tools in this study.
Let $G$ be a linear reductive group.
Let $U$ be a projective variety with a $G$-action.
Let $\nbigl_i$ $(i=1,2)$ be $G$-polarizations of $U$.
Then, we have the open subset $U^{ss}(\nbigl_i)$
of semistable points of $U$ with respect to $\nbigl_i$.
It is interesting to compare the stacks
$\nbigm_i:=U^{ss}(\nbigl_i)/G$ $(i=1,2)$.

For that purpose,
Thaddeus introduced the idea of the master space.
Let us consider the $G$-variety
$\TH=\proj(\nbigl_1^{-1}\oplus \nbigl_2^{-1})$
on $U$.
We have the canonical polarization $\nbigo_{\proj}(1)$
on $\TH$.
We have the canonically defined $G$-action on $\TH$,
and $\nbigo_{\proj}(1)$ gives the $G$-polarization.
The set of the semistable points is denoted by $\TH^{ss}$.
Then we obtain the stack $M=\TH^{ss}/G$.

We have the $G_m$-action on $\TH$ given by $\rho(t)([x:y])=[t\cdot x:y]$,
where $[x:y]$ denotes the homogeneous coordinate of $\TH$
along the fiber direction.
The action $\rho$ commutes with the action of $G$.
Thus we have the $G_m$-action on $M$,
which we denote by $\rhobar$.

We have the natural inclusion
$\TH_i=\proj(\nbigl_i^{-1})\lrarr \TH$.
Due to $\nbigo_{\proj}(1)_{|\TH_i}=\nbigl_i$,
it is easy to observe $\TH_i^{ss}=U(\nbigl_i)^{ss}$.
Thus we have the inclusion $\nbigm_i\lrarr M$.
Since $\TH_i$ is a component of the fixed point set
with respect to the action $\rho$,
the stacks $\nbigm_i$ are the components
of the fixed point set of the action $\rhobar$.

We may have the fixed points of the action $\rhobar$,
which are not contained in $\nbigm_1\cup \nbigm_2$.
Let $x$ be a point of $\TH^{ss}$,
which is not a fixed point of the action $\rho$.
Assume $G_m\cdot x\subset G\cdot x$.
Then the point $\pi(x)$ of $M$ is a fixed point
of the action $\rhobar$,
where $\pi$ denotes the natural projection $\TH^{ss}\lrarr M$.
The component of such fixed points
is called the exceptional fixed point set.
In a sense,
the information on the difference of 
the quotient stacks $\nbigm_i$ $(i=1,2)$
are concentrated at the exceptional fixed point set.
We will later explain how to derive the information
by using $G_m$-localization in our case.

\begin{rem}
When we consider the categorical quotients
$M_i:=U^{ss}(\nbigl_i)/\!/G$
and $Y:=\TH^{ss}(\nbigl_i)/\!/(G\times G_m)$,
we obtain the morphisms $M_1\lrarr Y\llarr M_2$.
The diagram is called Thaddeus Flip.
It also provides us a significant tool
to compare the spaces $M_i$.
In particular,
if $M_i$ are smooth
and each morphism $M_i\lrarr Y$ is
the blow up along the smooth center,
the flip is quite useful.
\hfill\qed
\end{rem}

\subsubsection{A construction of the master space 
when the $2$-stability condition is satisfied}

\label{subsubsection;06.6.9.6}

We explain the way of the construction
of the master space for our moduli stacks,
when the $2$-stability condition is satisfied
(the subsubsection \ref{subsubsection;06.6.9.5}).
To begin with,
we give a remark.
It is known that the coarse moduli scheme of
semistable torsion-free sheaves
is obtained as the categorical quotient of 
the set of the semistable points of some projective variety
provided with a reductive group action.
But, we will discuss the moduli stacks of
semistable {\em parabolic} sheaves
or semistable {\em parabolic} $L$-Bradlow pairs,
although we omit to mention parabolic structure
in this {\em Introduction}.
For such parabolic objects,
the author does not know the reference to deal with
the problem whether the moduli spaces have such descriptions.
Hence, we need a modification of the construction
of the master space.

Let $y$ be an element of $\Type$.
We put $y(m):=\ch\bigl(\nbigo(m)\bigr)$.
Let $H_y$ be the Hilbert polynomial associated to $y$.
Let $V_m$ be an $H_y(m)$-dimensional vector space.
We put $V_{m,X}:=V_m\otimes\nbigo_X$.
The projectivization of $V_m$ is denoted by $\proj_m$.
Let $Q(m,y)$ denote the quot scheme.
It is the moduli scheme of the quotient sheaves of $V_{m,X}$
whose Chern character is $y(m)$.
It is easy to construct the moduli scheme
$Q(m,\yhat)$ of the quotient sheaves
of $V_{m,X}$ with {\em orientations}
whose Chern character is $y(m)$.
We fix an inclusion $\iota:\nbigo(-m)\lrarr L$.
Let $\delta$ be an element of $\nbigp^{\br}$.
Then, we obtain the subscheme
$Q^{ss}(m,y,[L],\delta)$ of $Q^{\circ}(m,y)\times\proj_m$,
which is the moduli scheme of the quotient sheaves
$(q:V_{m,X}\rarr \nbige)$
with non-trivial reduced $[L]$-sections
$\phi:L\lrarr \nbige(-m)$
with the following property:
\begin{itemize}
\item The Chern character of $\nbige$ is $y(m)$,
and the pair $(\nbige(-m),\phi)$ is $\delta$-semistable.
\end{itemize}
We put
$Q^{ss}(m,\yhat,[L],\delta):=
 Q^{ss}(m,y,[L],\delta)\times_{Q(m,y)}Q(m,\yhat)$.

On the other hand,
let $\det(y(m))$ denote the $H^2(X)$-part of $y(m)$.
We take a line bundle $\gbigl_a$ such that $c_1(\gbigl_a)=y(m)$.
We may assume
$H^i(X,\gbigl_a)=0$ $(i>0)$.
We denote by $Z_m$ the projectivization of $H^0(X,\gbigl_a)$,
which is called the Gieseker space.
We put $\nbiga_m:=Z_m\times\proj_m$

When a parameter $\delta\in\nbigp^{\br}$ is given,
the ample $\rnum$-line bundle $\nbigl$ 
is given on $\nbiga_m$ as follows:
\[
 \nbigl:=\nbigo_{Z_m}\bigl(P^{\delta}_{y}(m)\bigr)
\otimes \nbigo_{\proj_m}\bigl(\delta(m)\bigr)
\]
We have the naturally defined $\SL(V_m)$-action
on $\nbiga_m$,
and $\nbigl$ gives the $\SL(V_m)$-polarization.
Let $\nbiga_m(\nbigl)$ denote the open subset of
the semistable points with respect to $\nbigl$.
It is rather standard to show the following proposition.
\begin{prop}
[Proposition \ref{prop;06.4.11.50}]
\label{prop;06.6.9.1}
We have the $\SL(V_m)$-equivariant
closed immersion
$Q^{ss}(m,y,[L],\delta)\lrarr 
 \nbiga_m^{ss}(\nbigl)$.
\hfill\qed
\end{prop}

We take a large integer $k$
such that $\nbigl^{\otimes\,k}$ is actually a line bundle.
We take a rational number $\gamma$
whose absolute value is sufficiently small.
We put
$\nbigl_{\gamma}:=
 \nbigl^{\otimes\,k}\otimes\nbigo_{\proj_m}(\gamma)$.
Let $\nbiga_m^{ss}(\nbigl_{\gamma})$ denote 
the open subset of the semistable points
with respect to $\nbigl_{\gamma}$.
On the other hand,
we take $\delta_-<\delta<\delta_+$ sufficiently closely.
It is not difficult to observe the following
by using Proposition \ref{prop;06.6.9.1}:
\begin{lem}
\label{lem;06.6.9.2}
In the case $\gamma<0$,
 we have the closed immersion
 $Q^{ss}(m,y,[L],\delta_{-})\lrarr \nbiga^{ss}(\nbigl_{\gamma})$.
 In the case $\gamma>0$,
  we have the closed immersion
 $Q^{ss}(m,y,[L],\delta_+)\lrarr \nbiga^{ss}(\nbigl_{\gamma})$.
\hfill\qed
\end{lem}

We take rational numbers $\gamma_2<0<\gamma_1$
such that $|\gamma_i|$ are sufficiently small.
We take a large number $k'$ such that
$k'\cdot(\gamma_1-\gamma_2)=1$.
Let $\pi:\nbiga_m\lrarr \proj_m$ denote the projection.
We put 
$\nbigb:=
 \proj\bigl(\pi^{\ast}\nbigo_{\proj_m}(0)\oplus \nbigo_{\proj_m}(1)
 \bigr)$,
which is $\proj^1$-bundle over $\nbiga$.
The tautological line bundle is denoted by $\nbigo_{\proj}(1)$.
We put
$\nbigo_{\nbigb}(1):=
 \nbigo_{\proj}(1)\otimes \nbigl_{\gamma_1}^{\otimes\,k'}$.
It gives the $\SL(V_m)$-polarization 
of $\nbigb$.
Let $\nbigb^{ss}$ denote the open subset of
the semistable points of $\nbigb$ with respect to
$\nbigo_{\nbigb}(1)$.

We have the natural inclusion
$\nbigb_1:=\proj\bigl(\pi^{\ast}\nbigo_{\proj_m}(0)\bigr)
\subset \nbigb$
and $\nbigb_2:=\proj\bigl(\pi^{\ast}\nbigo_{\proj_m}(1)\bigr)\subset
 \nbigb$.
We remark $\nbigo_{\nbigb}(1)_{|\nbigb_i}=\nbigl_{\gamma_i}^{\otimes k'}$.
Let $\nbigb^{ss}_i$ denote the semistable points
of $\nbigb_i$ with respect to $\nbigo_{\nbigb}(1)$.
Then, we put as follows:
\[
 \THhat^{ss}:=
 Q^{ss}(m,\yhat,[L],\delta)\times_{\nbiga}\nbigb^{ss},
\quad
\THhat_i^{ss}:=
 Q^{ss}(m,\yhat,[L],\delta)\times_{\nbiga}\nbigb_i^{ss}.
\]
We put $\Mhat:=\THhat^{ss}/\GL(V_m)$,
which is the master space in this case.
We also put $\Mhat_i:=\THhat^{ss}/\GL(V_m)$.
Due to Lemma \ref{lem;06.6.9.2},
we have $\Mhat_1\simeq\nbigm^{s}(\yhat,[L],\delta_+)$
and $\Mhat_2\simeq\nbigm^s(\yhat,[L],\delta_-)$.

When the $2$-stability condition is satisfied,
it can be easily shown that $\Mhat$
is Deligne-Mumford and proper.
We have the naturally defined $G_m$-action on $\Mhat$,
as explained in the subsubsection \ref{subsubsection;06.6.9.3}.
The fixed point set is as follows:
\[
 \Mhat_1\sqcup \Mhat_2\sqcup
 \coprod_{(y_1,y_2)\in S(y,\delta)}
 \Mhat^{G_m}(y_1,y_2)
\]
Here, $\Mhat^{G_m}$ is isomorphic to the moduli stack
of the objects $(E_1,\phi,E_2;\rho)$ 
with the following properties:
\begin{itemize}
\item $(E_1,\phi)$ is $\delta$-stable $L$-Bradlow pair.
\item $E_2$ is semistable torsion-free sheaf.
\item $\rho$ is an orientation of $E_1\oplus E_2$.
\end{itemize}
It is easy to observe that
$\Mhat^{G_m}$ is isomorphic to
$\nbigm^{ss}(y_1,L,\delta)\times\nbigm^{ss}(\yhat_2)$
up to etale finite morphisms.

\subsubsection{The $G_m$-localization method in the case where the
   $2$-stability condition is satisfied}

\label{subsubsection;06.6.9.15}

We explain how to use the master space $\Mhat$
to obtain the transition formula (Theorem \ref{thm;06.6.8.101}),
when the $2$-stability condition is satisfied for $(y,L,\delta)$.
We use the notation in the subsubsection 
\ref{subsubsection;06.6.9.5}.
Let $\varphi:\Mhat\lrarr \nbigm(\yhat,[L])$ be the naturally defined
morphism.
Let $\nbigt(1)$ denote the trivial line bundle on $\nbigm(\yhat,[L])$
with the $G_m$-action of weight $1$.
We have the naturally defined $G_m$-action
on $\varphi^{\ast}\Ehat^u$,
$\varphi^{\ast}\nbigo_{\rel}(1)$
and $\varphi^{\ast}\nbigt(1)$.
Therefore, we obtain the following
$G_m$-equivariant cohomology classes on $\Mhat$:
\[
 \Phi_t:=P\bigl(\varphi^{\ast}\Ehat^u\bigr)\cdot 
 c_1\bigl(\varphi^{\ast}\nbigo_{\rel}(1)\bigr)^k,
\quad
 \Phibar_t:=\Phi_t\cdot c_1\bigl(\varphi^{\ast}\nbigt(1)\bigr)
\]
The master space is naturally provided with 
the perfect obstruction theory,
and hence we have the virtual fundamental class.
Therefore, we can consider the polynomial
$\int_{\Mhat}\Phibar_t\in\rnum[t]$.
Since $\varphi^{\ast}\nbigt(1)$ is the trivial line bundle,
the specialization at $t=0$ is trivial.
We use the localization theory of the virtual fundamental classes
due to Graber-Pandharipande (\cite{gp}).
We have the following equality in $\rnum[t,t^{-1}]$:
\[
 \int_{\Mhat}\Phibar_t
=\sum_{i=1,2}\int_{\Mhat_i}\frac{t\cdot \Phi_t}{\Eu(\gbign(\Mhat_i))}
+\sum_{(y_1,y_2)\in S(y,\delta)}
 \int_{\Mhat^{G_m}(y_1,y_2)}
 \frac{t\cdot \Phi_t}{\Eu\bigl(\gbign(\Mhat^{G_m}(y_1,y_2))\bigr)}
\]
Here, $\gbign(\Mhat_i)$ and $\gbign(\Mhat^{G_m}(y_1,y_2))$
denote the virtual normal bundles.
Therefore, we obtain the following equality:
\begin{equation}
\label{eq;06.6.9.5}
 \sum_{i=1,2}\int_{\Mhat_i}
 \underset{t=0}\Res
 \left(\frac{\Phi_t}{\Eu\bigl(\gbign(\Mhat_i)\bigr)}\right)
+\sum_{(y_1,y_2)\in S(y,\delta)}
 \int_{\Mhat^{G_m}(y_1,y_2)}
 \underset{t=0}\Res
 \left(\frac{\Phi_t}{\Eu\bigl(\gbign(\Mhat^{G_m}(y_1,y_2))\bigr)}
 \right)=0.
\end{equation}
It is easy to observe that the first term
of the left hand side of (\ref{eq;06.6.9.5}) can be rewritten
as follows:
\[
 -\int_{\Mhat_1}\Phi+\int_{\Mhat_2}\Phi
=-\Phi(\yhat,[L],\delta_+)
+\Phi(\yhat,[L],\delta_-)
\]
Note that $\Mhat^{G_m}(y_1,y_2)$
is isomorphic to $\nbigm^{ss}(y_1,\yhat_2,L,\delta)$
up to etale finite morphisms.
By a formal calculation,
it can be shown that the second term 
is same as the right hand side of (\ref{eq;06.6.8.121}).
Thus we obtain the transition formula
in the case where the $2$-stability condition is satisfied
for $(y,L,\delta)$.

In the case $p_g>0$,
the transition problem is rather easy,
because of the following proposition.
\begin{prop}
[Proposition \ref{prop;06.5.31.1} and 
Proposition \ref{prop;06.5.11.500}]
\label{prop;06.6.9.34}
\mbox{{}}
\begin{itemize}
\item
 In the case $\rank(y)>1$ and $p_g>0$,
 the virtual fundamental class of 
 $\nbigm^{ss}(y,L,\delta)$ is trivial.
 (We remark that the virtual fundamental class of
 $\nbigm^{ss}(\yhat,[L],\delta)$ is not necessarily
 trivial.)
\item
Even in the case $\rank(y)=1$ and $p_g>0$,
we have the vanishing $[\nbigm^{ss}(y,L,\delta)]=0$,
if the $2$-vanishing condition is satisfied
for $(y,L,\delta)$.
\item
On the other hand,
 we have
 $\kappa^{\ast}\bigl([\nbigm^{ss}(y,L,\delta)]\bigr)
=[\nbigm^{ss}(\yhat,[L],\delta)]$,
where $\kappa$ denotes the naturally defined 
 etale finite morphism
 $\nbigm^{ss}(\yhat,[L],\delta)\lrarr\nbigm^{ss}(y,L,\delta)$
 in the case $p_g=0$.
\hfill\qed
\end{itemize}
\end{prop}

Due to Proposition \ref{prop;06.6.9.34},
we obtain the vanishing of the contributions
in (\ref{eq;06.6.9.5})
from $(y_1,y_2)$ such that $\rank(y_1)>1$,
when $p_g>0$.
Even in the case $\rank(y_1)=1$,
the contributions vanish if $L$ is sufficiently ``negative''.

\subsubsection{The enhanced master space}
\label{subsubsection;06.6.9.50}

When the $2$-stability condition is not satisfied,
$\Mhat$ in the subsubsection \ref{subsubsection;06.6.9.6}
is not Deligne-Mumford.
So we need some modification,
which we explain in this subsubsection.
We use the notation in the subsubsection \ref{subsubsection;06.6.9.6}.
We put as follows:
\[
 \nbigbtilde:=\nbigb\times\Flag(V_m,\Nbar),
\quad
\nbigbtilde_i:=\nbigb_i\times\Flag(V_m,\Nbar).
\]
Here, $\Flag(V_m,\Nbar)$ denotes the full flag variety of $V_m$:
\[
\Flag(V_m,\Nbar)=\Bigl\{\nbigf_1\subset\nbigf_2\subset\cdots\subset
 \nbigf_N\,\Big|\,
 \dim\nbigf_i/\nbigf_{i-1}=1
 \Bigr\}
\]

Let $G_l(V_m)$ denote the Grassmannian variety of 
the $l$-dimensional subspaces of $V_m$.
We have the canonical polarization $\nbigo_{G_l}(1)$.
We have the morphism
$\rho_l:\Flag(V_m,\Nbar)\lrarr G_l(V_m)$
given by $\rho_l(\nbigf_{\ast})=\nbigf_l$.
Take small positive numbers $n_i$
$(i=1,2,\ldots,N)$,
then we obtain the $\SL(V_m)$-polarization of $\nbigbtilde$:
\[
 \nbigo_{\nbigbtilde}(1):=
 \nbigo_{\nbigb}(1)\otimes
 \bigotimes_{l=1}^N\rho_l^{\ast}\nbigo_{G_l(V_m)}(n_i)
\]
Let $\nbigbtilde^{ss}$ and $\nbigbtilde^{ss}_i$
denote the set of semistable points
with respect to $\nbigo_{\nbigbtilde}(1)$.
In this case, we put as follows:
\[
 \THhat^{ss}:=
  Q^{ss}(m,\yhat,[L],\delta)\times_{\nbiga}\nbigbtilde^{ss},
\quad
 \THhat_i^{ss}:=
 Q^{ss}(m,\yhat,[L],\delta)\times_{\nbiga}\nbigbtilde_i^{ss}.
\]
We use the notation $\Mhat$ to denote
$\THhat^{ss}/\GL(V_m)$,
and we call it the enhanced master space.
We also put $\Mhat_i:=\THhat^{ss}_i/\GL(V_m)$.
Due to Lemma \ref{lem;06.6.9.2},
it can be shown that 
$\Mhat_1$ is the full flag bundle
over $\nbigm^{s}(\yhat,[L],\delta_+)$
associated to the vector bundle $p_{X\,\ast}\Ehat^u(m)$.
Similarly, $\Mhat_2$ is the full flag bundle over
$\nbigm^{s}(\yhat,[L],\delta_-)$.

We can show the following.
\begin{prop}
[Proposition \ref{prop;06.5.16.50}
and Proposition \ref{prop;06.6.9.10}]
$\Mhat$ is Deligne-Mumford
and proper.
\hfill\qed
\end{prop}

We have the naturally defined $G_m$-action on $\Mhat$,
as explained in the subsubsection \ref{subsubsection;06.6.9.3}.
To describe the fixed point set,
we need some preparation.

\begin{df}
A decomposition type is defined to be
a datum $\gbigi:=(y_1,y_2,I_1,I_2)$ as follows:
\begin{itemize}
\item
 $y=y_1+y_2$ in $\Type$
 such that
 $P_{y_1}^{\delta}=P_{y}^{\delta}$.
\item
 $\Nbar=I_1\sqcup I_2$
 such that  $|I_i|=H_{y_i}(m)$,
 here $H_{y_i}$ denote the Hilbert polynomials
 associated to $y_i$.
\end{itemize}
The set of the decomposition types 
is denoted by $\Dec(m,y,\delta)$.
For $\gbigi=(y_1,y_2,I_1,I_2)\in\Dec(m,y,\delta)$,
we put $\gminik(\gbigi):=\min(I_2)-1$.
\hfill\qed
\end{df}

We introduce the notion of $(\delta,\ell)$-semistability.
\begin{df}
Let $(E,[\phi])$ be a reduced $L$-Bradlow pair on $X$,
and let $\nbigf$ be a full flag of $H^0\bigl(X,E(m)\bigr)$.
Let $\ell$ be any positive integer.
We say that $(E,[\phi],\nbigf)$ is $(\delta,\ell)$-semistable,
if the following conditions are satisfied:
\begin{itemize}
\item
 $(E,[\phi])$ is $\delta$-semistable.
\item
 Take any Jordan-H\"older filtration
 of $(E_{\ast},[\phi])$ with respect to $\delta$-semistability:
\[
 E^{(1)}\subset
 E^{(2)}\subset\cdots
 \subset E^{(i-1)}
 \subset (E^{(i)},\phi)
\subset\cdots\subset (E^{(k)},\phi)
\]
Then we have
$\nbigf_{\ell}\cap H^0\bigl(X,E^{(i-1)}(m)\bigr)=\{0\}$
and 
 $\nbigf_{\ell}\not\subset H^0\bigl(X,E^{(j)}(m)\bigr)$
 for $j<k$.
\end{itemize}

We denote by
$\nbigmtilde^{ss}\bigl(y,[L],(\delta,\ell)\bigr)$,
the moduli stack of such tuples
$(E,[\phi],\nbigf)$.
In the oriented case,
we use the notation
$\nbigmtilde^{ss}\bigl(\yhat,[L],\alpha_{\ast},(\delta,\ell)\bigr)$
as usual.

Similarly,
we have the notion of $(\delta,\ell)$-semistability 
for a $L$-Bradlow pairs $(E,\phi)$ 
with a full flag $\nbigf$ of $H^0(X,E(m))$
such that $\phi\neq 0$.
The moduli stack is denoted by
$\nbigmtilde^{ss}\bigl(y,L,(\delta,\ell)\bigr)$.
\hfill\qed
\end{df}

Then, the fixed point set is as follows:
\[
 \Mhat_1\sqcup \Mhat_2\sqcup
 \coprod_{\gbigi\in \Dec(m,y,\delta)}
 \Mhat^{G_m}(\gbigi)
\]
Here, $\Mhat^{G_m}(\gbigi)$ is isomorphic
to the moduli stack of the objects
$(E_1,\phi,E_2,\rho,\nbigf^{(1)},\nbigf^{(2)})$
with the following properties:
\begin{itemize}
 \item
 $(E_1,\phi)$ is a $\delta$-semistable $L$-Bradlow pair.
\item $\nbigf^{(1)}$ is a full flag of $H^0(X,E_1(m))$
 such that
 $(E_1,\phi,\nbigf^{(1)})$ is
 $\bigl(\delta,\gminik(\gbigi)\bigr)$-semistable.
\item
 $E_2$ is a semistable torsion-free sheaf.
\item
 $\nbigf^{(2)}$ is a full flag of $H^0(X,E_2(m))$
 such that
 $\bigl(E_2, \nbigf_{\min(I_2)}^{(2)}\bigr)$ 
 is $\epsilon$-semistable reduced $\nbigo(-m)$-Bradlow pair,
 where $\epsilon$ denotes any sufficiently small positive number.
\item
 $\rho$ is an orientation of $E_1\oplus E_2$.
\end{itemize}
Therefore, $\Mhat^{G_m}(\gbigi)$
is isomorphic to
$\nbigmtilde^{ss}\bigl(y_1,L,(\delta,\gminik(\gbigi))\bigr)
 \times\nbigmtilde^{ss}(\yhat_2,+)$
up to etale finite morphisms,
where $\nbigmtilde^{ss}(\yhat_2,+)$
denotes oriented semistable sheaf $(E_2,\rho_2)$
with a full flag $\nbigf^{(2)}$ as in the condition above.

\subsubsection{The $G_m$-localization method}

We explain how to use the enhanced master space $\Mhat$
to obtain the transition formulas.
The argument is essentially same as that
in the subsubsection \ref{subsubsection;06.6.9.15}.
Let $\nbigm(m,\yhat,[L])$ denote the open subset of
$\nbigm(\yhat,[L])$ determined by the condition $O_m$.
(See the subsubsection \ref{subsubsection;06.6.9.25}.)
We have the vector bundle
$p_{X\,\ast}\Ehat^u(m)$
on $\nbigm(m,\yhat,[L])$.
The associated full flag bundle is denoted by
$\nbigmtilde(m,\yhat,[L])$.
Let $\Ttilde_{\rel}$ denote the relative tangent bundle
of the smooth morphism
$\nbigmtilde(m,\yhat,[L])\lrarr \nbigm(m,\yhat,[L])$.
We have the naturally defined morphism
$\varphi:\Mhat\lrarr \nbigmtilde(m,\yhat,[L])$.
Let $\nbigt(1)$ denote the trivial line bundle on 
$\nbigmtilde(m,\yhat,[L])$
with the $G_m$-action of weight $1$.
We have the naturally defined $G_m$-action
on $\varphi^{\ast}\Ehat^u$,
$\varphi^{\ast}\nbigo_{\rel}(1)$,
$\varphi^{\ast}\Ttilde_{\rel}$
and $\varphi^{\ast}\nbigt(1)$.
We consider the following cohomology classes on $\Mhat$:
\[
 \Phitilde_t:=P(\varphi^{\ast}\Ehat^u)\cdot
 c_1\bigl(\varphi^{\ast}\nbigo_{\rel}(1)\bigr)^k
\cdot
 \frac{\Eu(\Ttilde_{\rel})}{H_y(m)!},
\quad
 \Phibar_t:=\Phitilde_t\cdot c_1\bigl(\varphi^{\ast}\nbigt(1)\bigr)
\]
By the argument in the subsubsection
\ref{subsubsection;06.6.9.15},
we obtain the following equality:
\begin{equation}
\label{eq;06.6.9.40}
 \sum_{i=1,2}\int_{\Mhat_i}
 \underset{t=0}\Res
 \left(\frac{\Phitilde_t}{\Eu\bigl(\gbign(\Mhat_i)\bigr)}\right)
+\sum_{\gbigi\in\Dec(m,y,\delta)}
 \int_{\Mhat^{G_m}(\gbigi)}
 \underset{t=0}\Res
 \left(\frac{\Phitilde_t}{\Eu(\gbign(\Mhat^{G_m}(\gbigi)))} 
 \right)=0.
\end{equation}
Here $\gbign(\Mhat_i)$ and $\gbign(\Mhat^{G_m}(\gbigi))$
denote the virtual normal bundles.
Since $\Mhat_i$ are the full flag bundles
over $\nbigm^{ss}(\yhat,[L],\delta_{\kappa})$ $(\kappa=\pm)$,
it is easy to observe that the first term
of the left hand side of (\ref{eq;06.6.9.5}) can be rewritten
as follows:
\[
 -\int_{\Mhat_1}\Phitilde+\int_{\Mhat_2}\Phitilde
=-\Phi(\yhat,[L],\delta_+)
+\Phi(\yhat,[L],\delta_-)
\]
Let us see the second term in the case $p_g>0$.
Recall that $\Mhat^{G_m}(\gbigi)$
is isomorphic to
$\nbigmtilde^{ss}\bigl(\yhat_1,L,(\delta,\gminik(\gbigi))\bigr)
\times
 \nbigmtilde^{ss}\bigl(\yhat_2,+\bigr)$
up to etale finite morphisms.
Due to the vanishing similar to Proposition
\ref{prop;06.6.9.34},
we obtain the vanishing of the contributions
in (\ref{eq;06.6.9.40}) from $\gbigi=(y_1,y_2,I_1,I_2)$
such that $\rank(y_1)>1$.
In the case $\rank(y_1)=1$,
the $\bigl(\delta,\gminik(\gbigi)\bigr)$-semistability condition
is trivial.
Therefore,
$\nbigmtilde^{ss}\bigl(y_1,L,(\delta,\gminik(\gbigi))\bigr)$
is the full flag bundle over $\nbigm(y_1,L)$
associated to the vector bundle $p_{X\,\ast}E_1^u(m)$.
On the other hand,
$\nbigmtilde^{ss}(\yhat_2,+)$
is the flag variety bundle over
$\nbigm^{ss}\bigl(\yhat_2,[\nbigo(-m)],\epsilon\bigr)$
for any sufficiently small positive number $\epsilon$.
(See the subsubsection \ref{subsubsection;06.5.17.20}
 for more detail.)
Therefore, we can easily observe that
the second term of the left hand side of
(\ref{eq;06.6.9.40}) is same as
the right hand side of (\ref{eq;06.6.8.130}).
Thus, we obtain the transition formula
in the case $p_g>0$.

On the other hand,
we need some more additional argument in the case $p_g=0$.
Due to the equality (\ref{eq;06.6.9.40}),
we obtain the expression
$\Phi(\yhat,[L],\delta_+)-\Phi(\yhat,[L],\delta_-)$
as the summation of the terms of the following form:
\[
 \int_{\nbigmtilde^{ss}(\yhat_1,[L],\delta,\ell)\times\nbigm^{ss}(\yhat_2)}
 \!\!\!\!\!\Psi(y_1,y_2)
\]
Then, we take the enhanced master space connecting
$\nbigmtilde^{ss}\bigl(\yhat_1,[L],(\delta,\ell)\bigr)$
and $\nbigmtilde^{ss}\bigl(\yhat_1,[L],\delta_-\bigr)$.
Such a space can be constructed by the method
in the subsubsection \ref{subsubsection;06.6.9.50}.
We have only to choose appropriate numbers
$n_i$ $(i=1,\ldots,N)$.
Let $S(\ell)$ denote the set of decomposition types
$\gbigi^{(2)}=(y_1^{(2)},y_2^{(2)},I_1^{(2)},I_2^{(2)})$
with the following property:
\begin{itemize}
\item
 $y_1^{(2)}+y_2^{(2)}=y_1$.
\item
 $I_1^{(2)}\sqcup I_2^{(2)}=\bigl\{1,\ldots,H_{y_1}(m)\bigr\}$
 and $|I_i^{(2)}|=H_{y_i^{(2)}}(m)$.
\item
 $\{1,\ldots,\ell\}\subset I_1^{(2)}$.
\end{itemize}
Then, applying the localization method again,
we obtain the following:
\[
  \int_{\nbigmtilde^{ss}(\yhat_1,[L],\delta,\ell)\times\nbigm^{ss}(\yhat_2)}
\!\!\! \Psi(y_1,y_2)=
\int_{\nbigmtilde^{ss}(\yhat_1,[L],\delta_-)\times\nbigm^{ss}(\yhat_2)}
\!\!\!\Psi(y_1,y_2)
+\sum_{\gbigi^{(2)}\in S(\ell)}
 \int_{\Mhat^{G_m}(\gbigi^{(2)})}\!\!\!\Psi^{(2)}(y_1^{(2)},y_2^{(2)},y_2)
\]
We can apply this procedure inductively.
Note $\rank(y_1^{(2)})<\rank(y_1)<\rank(y)$.
Hence, the inductive process will stop.
Thus, we can obtain the general transition formula,
in principle,
by using a rather easy combinatorial argument.
However, the general formula would be comparatively complicated,
and it is less interesting for the author at this moment.
Hence, we restrict ourselves to the transition formula
in the special case, as in Theorem \ref{thm;06.6.9.100}.

\subsection{Outline of the Paper}
\subsubsection{Section \ref{section;06.7.3.1}}

In the subsection \ref{subsection;06.7.3.2},
we prepare some notation.
In the subsection \ref{subsection;06.7.3.3},
we review basic results from the geometric invariant theory.
In particular, we recall a sufficient condition
for a quotient stack to be Deligne-Mumford and proper.
We also recall the numerical criterion and 
calculate some easy examples.
The results will be used in the section \ref{section;06.7.3.10}.

In the subsection \ref{subsection;06.7.3.4},
we review the basis of cotangent complex.
Then, we recall how to calculate some cotangent 
complexes of quotient stacks
in the subsubsection \ref{subsubsection;06.4.29.15},
which will be used in the section \ref{section;06.6.4.250}
in many times.
We also recall some more examples in 
the subsubsection \ref{subsubsection;06.5.25.2},
which will be used in the subsections \ref{subsection;06.6.1.10}
and \ref{subsection;06.6.1.11}.

In the subsection \ref{subsection;06.7.3.5},
we review the obstruction theory of
Behrend-Fantechi \cite{bf}.
Then, we explain a naive strategy
to construct the obstruction theory 
in the subsubsection \ref{subsubsection;06.6.22.1}.
We recall the obstruction theory of
the locally free subsheaves in the subsubsection
\ref{subsubsection;06.5.4.5}.
It gives the obstruction theory of the moduli
of torsion-free quotient sheaves over a smooth projective surface.
The result will be used in the subsection \ref{subsection;06.6.13.5}.
We also obtain the smoothness of the moduli
of quotient torsion-sheaves over a smooth projective curve,
although we do not use it later.
In the subsubsection \ref{subsubsection;06.6.24.5},
we recall the obstruction theory
of the filtration of a vector bundle on a smooth projective curve.
It will be used to see the obstruction theory of
the parabolic structure.

In the subsection \ref{subsection;06.7.3.6},
we recall some easy results for
equivariant complexes on Deligne-Mumford stacks
with GIT construction,
which will be used in the subsection \ref{subsection;06.6.17.2}.
In the subsection \ref{subsection;06.7.3.7},
we give some elementary remarks
on extreme sets,
which are used in the subsections 
\ref{subsection;06.5.21.110}--\ref{subsection;06.6.6.2}.
In the subsection \ref{subsection;06.5.17.5},
we give easy remarks on the twist of line bundles.

\subsubsection{Section \ref{section;06.5.17.20}}

In the subsection \ref{subsection;06.5.17.21},
we review the basic notion.
In the subsubsections
\ref{subsubsection;06.7.3.211}--\ref{subsubsection;06.7.3.212},
we recall the definition of some structure on torsion-free sheaves
such as orientation, parabolic structure,
$L$-section, and reduced $L$-section.
In the subsubsection \ref{subsubsection;06.4.24.10},
we prepare the notion of type
and the notation of some moduli stacks.
In the subsubsection \ref{subsubsection;06.4.11.75},
we recall the notion of
the relative tautological line bundle 
of the moduli stacks of oriented reduced $L$-Bradlow pair.
We also see the relation of the moduli stack
of oriented reduced $L$-Bradlow pair
and the moduli stack of unoriented unreduced $L$-Bradlow pair.

In the subsection \ref{subsection;06.5.17.22},
we recall the Hilbert polynomials.
Then, we have the naturally defined semistability conditions,
which is discussed in the subsection \ref{subsection;06.7.3.23}.
We recall the notion of Harder-Narasimhan filtration
and partial Jordan-H\"older filtration
in the subsubsection \ref{subsubsection;06.7.3.231}.
Then, we introduce the notion of $(\delta,\ell)$-semistability
in the subsubsection \ref{subsubsection;06.4.25.10},
which is useful to control the transition.

In the subsection \ref{subsection;06.7.3.24},
we review boundedness of some families.
We recall foundational theorems in the subsubsection 
\ref{subsubsection;06.7.3.241},
Then, we recall the boundedness of
semistable $L$-Bradlow pairs
when the parameter is varied
in the subsubsection \ref{subsubsection;06.7.3.242}.
The important observation is due to M. Thaddeus.
In the subsubsection \ref{subsubsection;06.7.3.243},
we see boundedness of Yokogawa family,
which will be used to show properness
of some morphisms in the section \ref{section;06.7.3.10}.

In the subsection \ref{subsection;06.7.3.25},
we recall the $1$-stability and $2$-stability conditions.
In the subsection \ref{subsection;06.5.15.15},
we recall some moduli schemes of quotient sheaves
with some structure.

\subsubsection{Section \ref{section;06.7.3.10}}

In the subsection \ref{subsection;06.7.3.11},
we review a basic result of the geometric invariant
theory for the construction of the moduli stack
of $\delta$-semistable parabolic $L$-Bradlow pairs.
In the subsection \ref{subsection;06.5.15.100},
we consider the perturbation of $\delta$-semistability
condition.
Namely, we multiply the full flag bundle
to the quot schemes,
and we discuss what is obtained
for small perturbation of the semistability conditions.

The results in the subsections
\ref{subsection;06.5.21.110}--\ref{subsection;06.6.6.2}
are the core of this paper,
which are significant to discuss the transition formula.
In the subsection \ref{subsection;06.5.21.110},
we construct the enhanced master space,
and we show that it is Deligne-Mumford and proper.
In the subsection \ref{subsection;06.6.6.2},
we see the fixed point set with respect to the torus action.

In the subsection \ref{subsection;06.5.11.100},
we construct the enhanced master space in the oriented case,
and we give a description of
the stack theoretic fixed point set with respect
to the natural torus action.
They are essentially reformulation of the results in the previous subsections.
We give a more convenient description
of the fixed point set in the subsection \ref{subsection;06.6.21.3},
i.e.,
we observe that they are isomorphic to the moduli stacks
of objects with the lower ranks,
up to etale finite morphisms.

In some cases, the problem is much simpler.
The statements for such cases are given
in the subsection \ref{subsection;06.7.3.15}.

\subsubsection{Section \ref{section;06.6.4.250}}

In the subsection \ref{subsection;06.7.3.31},
we discuss the deformation theory 
associated to torsion-free sheaf $E$ on $U\times X$,
where $X$ is a smooth projective surface.
We put 
$\gminig(V_{\cdot}):=\nhom(V_{\cdot},V_{\cdot})^{\lor}[-1]$
and $\Ob(V_{\cdot}):=
 Rp_{X\,\ast}\bigl(\gminig(V_{\cdot})\otimes\omega_X\bigr)$
for a resolution $V_{\cdot}$ of $E$.
In the subsubsection \ref{subsubsection;06.5.2.3},
we explain how we obtain the morphism
$\gminig(V_{\cdot})\lrarr L_{U\times X/X}$,
and hence $\Ob(V_{\cdot})\lrarr L_U$.
In the subsubsection \ref{subsubsection;06.5.9.20},
we see that $\gminig(V_{\cdot})$ is decomposed
into the trace-free part and the diagonal part,
and that the diagonal part is related to the determinant bundle.
In the subsubsection \ref{subsubsection;06.5.4.100},
we give some factorization which will be useful 
in the construction of the obstruction theory of the master space.
In the subsubsection \ref{subsubsection;06.7.3.311},
we give the obstruction theory of some open subset 
of the moduli stack of torsion-free sheaves,
by directly applying the construction in the subsubsection
\ref{subsubsection;06.5.2.3}.
In the subsubsection \ref{subsubsection;06.5.11.1},
we discuss the case of the moduli of line bundles,
which will be used in the construction
of the relative obstruction theory for orientation
in the subsection \ref{subsection;06.7.3.32}.

In the subsection \ref{subsection;06.7.3.33},
we discuss the relative obstruction theory for $L$-section. 
In the subsubsections
\ref{subsubsection;06.5.4.1}--\ref{subsubsection;06.5.4.50},
we give the construction and show the relative obstruction property.
In the subsubsection \ref{subsubsection;06.5.21.250},
we give some factorization which will be useful
to discuss the obstruction theory of the master space.
In the subsection \ref{subsection;06.7.3.34},
we discuss the relative obstruction theory for
reduced $L$-section. 
We need to make some modification
to the construction in the subsection \ref{subsection;06.7.3.33}.
In the subsection \ref{subsection;06.7.3.35},
we discuss the relative obstruction theory 
for parabolic structure.
By pulling them together,
we can easily construct the obstruction theory
of the moduli stacks of parabolic $L$-Bradlow pairs
and some related objects,
which is explained in the subsection \ref{subsection;06.6.13.5}.

Then, we discuss the obstruction theory of the master space
in the subsection \ref{subsection;06.5.23.1}.
Once we have the factorizations as in the subsubsections
\ref{subsubsection;06.5.4.100},
\ref{subsubsection;06.5.21.250} and
\ref{subsubsection;06.5.4.150},
the construction is easy.
We also obtain the obstruction theories for some related stacks.
In the subsubsections
\ref{subsubsection;06.5.21.600}--\ref{subsubsection;06.6.13.16},
we give only the statements in some easy cases
for explanation.

In the subsection \ref{subsection;06.6.21.2},
we discuss the obstruction theory of the fixed point set.
In the subsection \ref{subsection;06.6.17.2},
we discuss the equivariant obstruction theory of the master
space and the induced obstruction theory of the fixed point set.
We give the statements for the simple cases
in the subsubsections
\ref{subsubsection;06.7.3.36.1}--\ref{subsubsection;06.6.12.1}
and \ref{subsubsection;06.7.3.37.1}--\ref{subsubsection;06.7.3.37.2}.

\subsubsection{Section \ref{section;06.7.3.40}}

By showing the perfectness of the obstruction theories,
we obtain the virtual fundamental classes for some stacks,
which is discussed in the subsection \ref{subsection;06.7.3.41}.
We compare the virtual fundamental classes
of the moduli stack of
the $\delta$-stable oriented reduced $L$-Bradlow pairs
and the moduli stack of $\delta$-stable $L$-Bradlow pairs
in the subsection \ref{subsection;06.7.3.42}.
Although the moduli stacks are isomorphic up to etale finite morphisms,
the obstruction theories are not same in general,
and we obtain the vanishing of the virtual fundamental class of
the moduli of $\delta$-stable $L$-Bradlow pairs
in the case $p_g>0$.

In the subsection \ref{subsection;06.6.1.10},
we discuss the virtual fundamental classes
of the moduli stack of the objects with rank one.
In the subsubsection \ref{subsubsection;06.6.30.3},
we see the moduli of $L$-abelian pairs.
In particular,
we give a detailed description of the virtual fundamental class
in the case where $H^2(X,\nbigo)\neq 0$ and $H^1(X,\nbigo)=0$
are satisfied.
In the subsubsection \ref{subsubsection;06.7.3.431},
we discuss the obstruction theory of the parabolic Hilbert schemes.
In the rest of the subsection,
we show the splitting given in Proposition \ref{prop;06.5.13.35}.

In the subsections \ref{subsection;06.7.3.44}--\ref{subsection;06.6.1.11},
we discuss the relations of the obstruction theories
of some moduli stacks.

\subsubsection{Section \ref{section;06.6.30.5}}

In the subsection \ref{subsection;06.6.19.110},
we explain how we think the cohomology and the evaluation
for our theory.
In the subsection \ref{subsection;06.7.3.52},
we show the transition formulas in the simple cases.
They are sufficiently useful
for the construction of the invariants,
which is discussed in the subsection \ref{subsection;06.7.3.53}.
They also provide the sufficient tool
to discuss the transition problem in the rank $2$ case,
which is done in the subsection \ref{subsection;06.7.3.54}.

In the subsection \ref{subsection;06.7.3.55},
we discuss the transition formula
for the case $p_g>0$.
It is rather easy to show, and the formula is formally same
as that in the simple case.

Then, we discuss the transition formula
for the case $p_g=0$,
in the subsection \ref{subsection;06.6.21.30}.
By using it, we obtain the weak wall crossing formula
in the subsection \ref{subsection;06.7.3.57}.
We write down the weak wall crossing formula 
and the weak intersection rounding formula
for the rank $3$ case in the subsubsection
\ref{subsubsection;06.6.7.1}--\ref{subsubsection;06.6.8.300}.
We also give a transition formula
for a critical parabolic weight in the subsubsection
\ref{subsubsection;06.7.3.575}.

In the subsection \ref{subsection;06.6.7.400},
we derive the weak intersection rounding formula
from the weak wall crossing formula.

\subsection{Some Remarks}
\label{subsection;06.6.30.20}
\subsubsection{Further study}
\label{subsubsection;06.7.2.1}

This note is tentative, partially because of
the recent intensive development of the theory of stacks.
We use only the rather old results in this paper.
The author believes that it should be one of the main themes
of the study of the stacks to make it easy to deal with 
the subject and the formalism in this work.
Hence, he hopes that our arguments
would be replaced with the new ways.

For example,
as already mentioned in Remark \ref{rem;06.6.30.1},
the powerful theory of ``derived stack'' has been
developed (see \cite{toen} for overview).
It seems to be applicable to a wide range of similar problems,
contrast to our method in this paper.
The author thinks that our method is more elementary.
But the theory of ``derived stack'' will surely be standard
in algebraic geometry in very near future
(or, perhaps already?),
and hence he hopes that
our construction and argument
would be refreshed
from that point of view.

Then, it would be desired to compare the obstruction theories.
It is not clear for the author whether such comparison is easy or not.
However, he expects that
the comparison of the invariants would be easily done,
at least in the case $H^1(X,\nbigo)=0$,
if the new construction would be done appropriately.
We should have the transition formulas as given
in the section \ref{section;06.6.30.5},
due to which it can be reduced to the comparison of the invariants
obtained from the moduli stacks of the objects of rank one.
We should have the splitting as in Proposition \ref{prop;06.5.13.35},
and hence it can be reduced to the comparison
of the invariants obtained from the abelian pairs.
(See the subsubsection \ref{subsubsection;06.6.30.2}
for such a reduction in the case $H^2(X,\nbigo)\neq 0$.
We may also obtain such a reduction in the case $H^2(X,\nbigo)=0$,
although the formula would be more complicated.)
In the case $H^1(X,\nbigo)=0$,
the moduli of abelian pairs is smooth,
and the obstruction theory is given by the obstruction bundle
which should be as in the subsubsection \ref{subsubsection;06.6.30.3}.
Thus, the comparison of the invariants could be done.
The author hopes to clarify the points somewhere.

\vspace{.1in}
One of the important missing for our theory
is the blow up formula,
i.e., comparison of the invariants for $X$ and a blow up of $X$.
Originally, the author intended to develop the theory
``without blow up''.
But it seems to contain an interesting problem
for such a comparison,
and the author would like to come back to this problem,
if possible.

\subsubsection{Difference from the previous version}

In this second version,
we obtained the comparatively satisfactory transition formula
in the higher rank case,
which made us to obtain the weak wall crossing formula.
For that purpose, we introduced the notion of
$(\delta,\ell)$-semistability.
The other essential ideas are not changed.
The author also hopes that the readability is improved.

\subsection{Acknowledgement}
The author thanks H. Nakajima for his stimulating questions.
He also thanks F. Kato for some information.
The author thanks A. Ishii and Y. Tsuchimoto 
for their constant encouragement.
He is grateful to the colleagues of Department of Mathematics
at Kyoto University for their cooperations.
The author wrote the current version of this paper
during his stay at Max-Planck Institute for Mathematics.
He wrote the previous version of this paper
during his stay at Institute for Advanced Study.
He acknowledges the institutions
for their excellent hospitality and support.
He also wishes to acknowledge National Scientific Foundation
for a grant DMS 9729992,
although any opinions, findings and conclusions
or recommendations expressed in this material
do not necessarily reflect the views of
the National Science Foundation.
He thanks the supports
by Japan Society for the Promotion of Science
and the Sumitomo Foundation in some period.

\section{Preliminary}
\label{section;06.7.3.1}

\subsection{Notation}
\label{subsection;06.7.3.2}
\subsubsection{Vector bundles}

Let $Y$ be a variety.
Let $g:T\lrarr U$ be a morphism of stacks.
Then the naturally induced morphism $T\times Y\lrarr U\times Y$
is denoted by $g_Y$ or simply by $g$.

Let $V$ be a vector bundle on $Y$.
The sheaf of local sections of $V$
is also denoted by the same notation $V$,
if there are no risk of confusion.
But, we use some particular notation in the following case:
Let $V_1,V_2$ be vector bundles.
We have the sheaf $\nhom(V_1,V_2)$\index{$\nhom(V_1,V_2)$}
of the morphisms from $V_1$ to $V_2$.
The corresponding vector bundle
is denoted by $N(V_1,V_2)$.\index{$N(V_1,V_2)$}

Let $F$ be a vector bundle on $Y$.
The complement of the image of the $0$-section in $F$
is denoted by $F^{\ast}$,
i.e.,
$F^{\ast}:=F-Y$,
and the dual bundle of $F$ is denoted by
$F^{\lor}$.\index{$F^{\ast}$, $F^{\lor}$}
The projectivization of $F$
is denoted by $\proj(F^{\lor})$ or $\proj_F$.

\subsubsection{Coherent sheaves on a product}

Let $X$ be a scheme over $k$,
and let $U$ be a stack over $k$.
We denote by $p_X$\index{$p_X$}
the projection forgetting the $X$-component:
\[
 p_X:U\times X\lrarr U.
\]
Similarly,
$p_U$ denotes the projection $U\times X\lrarr X$.
We often denote $\{u\}\times X$ by $X_u$ for any point $u\in U$.

A coherent sheaf $E$ over $U\times X$
is called a $U$-coherent sheaf, if it is flat over $U$.
\index{$U$-coherent sheaf}
A $U$-coherent sheaf $E$ is called $U$-torsion free sheaf,
\index{$U$-torsion-free sheaf}
if $E_{|\{u\}\times X}$ is torsion-free
for each $u\in U$.
We will often omit to denote ``$U$-'',
if there are no risk of confusion.

For any coherent sheaf $E$ on $U\times X$,
we use the notation
$E(m)$ to denote
$E\otimes p_U^{\ast}\nbigo_X(m)$
for a given polarization $\nbigo_X(1)$ of $X$.

\subsubsection{Signature}

\label{subsubsection;06.4.29.30}

We follow the signature convention
in \cite{ks2}.
We recall some of them for later use
in our situation.

Let $X$ be an algebraic stack.
For two $\nbigo_X$-complexes 
$C^{\cdot}$ and $D^{\cdot}$,
let $\nhom(C^{\cdot},D^{\cdot})$
denote the complex whose $i$-th term is
$\bigoplus_{k-j=i}\nhom\bigl(C^{j},D^k\bigr)$
and whose differentials are given as follows:
\[
  \nhom\bigl(C^j,D^k\bigr)\lrarr
 \nhom\bigl(C^j,D^{k+1}\bigr)
\oplus
 \nhom\bigl(C^{j-1},D^{k}\bigr),
\quad
 a\longmapsto
 \bigl(
 d_D\circ a,\,\,
 (-1)^{k-j+1}a\circ d_C
 \bigr).
\]

Let us see some examples.
For a complex $C^{\cdot}$,
we denote the dual complex
$\nhom(C^{\cdot},\nbigo_X)$
by $C^{\cdot\,\lor}$.
The differential is as follows:
\[
 \nhom(C^{n},\nbigo_X)
\lrarr 
 \nhom(C^{n-1},\nbigo_X),
\quad
 a\longmapsto (-1)^{n+1}\cdot a\circ d_X
\]

For two term complexes
$C^{\cdot}=(C^{-1}\rarr C^0)$
and $D^{\cdot}=(D^{-1}\rarr D^0)$,
the complex $\nhom\bigl(C^{\cdot},D^{\cdot}\bigr)$
is given as follows:
\[
\begin{array}{ll}
 \nhom\bigl(C^{0},D^{-1}\bigr)
\lrarr
 \nhom\bigl(C^{0},D^{0}\bigr)
\oplus
 \nhom\bigl(C^{-1},D^{-1}\bigr),
 &
 a\longmapsto \bigl(d_D\circ a,\,a\circ d_C\bigr)
 \\
 \mbox{{}}\\
\nhom\bigl(C^{0},D^{0}\bigr)
\oplus
 \nhom\bigl(C^{-1},D^{-1}\bigr)
\lrarr
 \nhom\bigl(C^{-1},D^0\bigr),
 &
 (b_1,b_2)\longmapsto
 -b_1\circ d_C+d_D\circ b_2
\end{array}
\]
We will often use the dual
$\nhom(C^{\cdot},D^{\cdot})^{\lor}$,
which is given as follows:
\[
\begin{array}{ll}
 \nhom\bigl(D^{0},C^{-1}\bigr)
\lrarr 
 \nhom(D^0,C^0)
\oplus
 \nhom(D^{-1},C^{-1}),
 &
 a\longmapsto
 (-d_C\circ a,a\circ d_D)\\
 \mbox{{}}\\
 \nhom(D^0,C^0)
\oplus
 \nhom(D^{-1},C^{-1})
\lrarr
 \nhom\bigl(D^{-1},C^{0}\bigr),
 &
 (b_1,b_2)\longmapsto
 -b_1\circ d_D-d_C\circ b_2
\end{array}
\]

\subsubsection{Compatible diagrams}
\label{subsubsection;06.5.22.15}
\index{compatible diagrams}
Let $A_{i,j}$ $(i=1,2)$ $(j=1,2,3,4)$ be objects
in some category.
Assume that we are given morphisms
$\varphi_{j}:A_{1,j}\lrarr A_{2,j}$.
We also assume that we are given commutative
diagrams $(CD)_i$:
\[
 \begin{CD}
 A_{i,1} @>{a_i}>> A_{i,2}\\
 @V{b_i}VV @V{c_i}VV \\
 A_{i,3} @>{d_i}>> A_{i,4}
 \end{CD}
\]
We say that $(CD)_1$ and $(CD)_2$ are compatible
with respect to the morphisms $\varphi_{j}$ $(j=1,2,3,4)$,
if every face of the naturally obtained cube is commutative.
It is equivalent to the commutativity of the following diagrams:
\[
\begin{CD}
 A_{1,1} @>>> A_{1,2}\\
 @VVV @VVV \\
 A_{2,1} @>>> A_{2,2}
\end{CD}
\quad
\begin{CD}
 A_{1,1}@>>> A_{1,3}\\
 @VVV @VVV \\ 
 A_{2,1} @>>> A_{2,3}
\end{CD}
\quad
\begin{CD}
 A_{1,2}@>>> A_{1,4}\\
 @VVV @VVV \\
 A_{2,2} @>>> A_{2,4}
\end{CD}
\quad
\begin{CD}
 A_{1,3} @>>> A_{1,4}\\
 @VVV @VVV \\
 A_{2,3} @>>> A_{2,4}
\end{CD}
\]

\subsubsection{Filtrations and complexes on a curve}
\label{subsubsection;06.5.23.10}

Let $\nbigd$ be a smooth projective curve 
over a stack $S$.
Let $E_a$ $(a=1,2)$
be coherent $\nbigo_{\nbigd}$-modules
which are flat over $S$.
Assume that we are given a decreasing filtration 
$F(E_a)=(F_i(E_a)\,|\,i=1,\ldots,l)$
of $E_a$ such that $\Cok_i(E_a)=E_a/F_{i+1}(E_a)$
are flat over $S$.

Let $V_{a,\cdot}=(V_{a,-1}\rarr V_{a,0})$ be
a locally free resolution of $E_a$.
We put $V_a^{(1)}:=V_{a,0}$, $V_a^{(l+1)}=V_{a,-1}$
and $V_a^{(i)}:=\Ker\bigl(V_{a,0}\lrarr \Cok_i(E_a)\bigr)$
$(i=2,\ldots,l)$.
Let $f_i:V_D^{(i+1)}\lrarr V_D^{(i)}$,
$t_i:V_D^{(i)}\lrarr V_D^{(1)}$
and $s_i:V_D^{(l+1)}\lrarr V_D^{(i)}$
denote the inclusions.

Let us consider the complex $C_1(V_1^{\ast},V_2^{\ast})$
given as follows:\index{$C_1(V_1^{\ast},V_2^{\ast})$}
\[
\begin{CD}
 \nhom\bigl(V_1^{(1)},V_2^{(l+1)}\bigr)
@>{d^{-1}}>>
\bigoplus_{i=1}^{l+1}
 \nhom\bigl(
 V_1^{(i)},V_2^{(i)}
 \bigr)
@>d^0>>
 \bigoplus_{i=1}^l
 \nhom\bigl(V_1^{(i+1)},V_2^{(i)}\bigr)
\end{CD}
\]
Here the first term stands in the degree $-1$,
The differentials $d_{i}$ are given as follows:
\begin{equation}
 d^{-1}(a)=\bigl(
 s_i\circ a\circ t_i\,\big|\,i=1,\ldots,l+1
 \bigr),
\end{equation}
\begin{equation}
 \label{eq;06.6.14.2}
  d^0(b_1,\ldots,b_l)
=\bigl(-f_1\circ b_1+b_2\circ f_1,
 -f_2\circ b_2+b_3\circ f_2,\ldots,
 -f_l\circ b_l+b_{l+1}\circ f_l\bigr)
\end{equation}
We have the naturally defined morphism:
\begin{equation}
\label{eq;06.6.14.1}
\varphi=(\varphi_i):C_1(V_1^{\ast},V_2^{\ast})
\lrarr \nhom\bigl(V_{1,\cdot},V_{2,\cdot}\bigr)
\end{equation}
The morphism $\varphi_1$ 
is given by
$ \varphi_1(a_i)=\sum s_{i+1}\circ a_i\cdot t_i$.
The morphism $\varphi_0$ is the projection
by the identification
$V_0=V^{(1)}$ and $V_{-1}=V^{(l+1)}$.
The morphism $\varphi_2$ is the identity.
It can be directly checked that $\varphi$ is the morphism
of complexes.
We put
$C_2(V_1^{\ast},V_{2}^{\ast}):=
 \cone(\varphi)[-1]$.
\index{$C_2(V_1^{\ast}V_2^{\ast})$}
The following lemma is easy to check.
\begin{lem}
The complexes $C_i(V_1^{\ast},V_2^{\ast})$ 
and the morphism 
$\varphi:C_1(V_1^{\ast},V_2^{\ast})\lrarr 
 \nhom(V_{\cdot},V_{\cdot})$
depend only on  $(E_1,F)$ and $(E_2,F)$
in the derived category $D(\nbigd)$.
\hfill\qed
\end{lem}

\begin{notation}
We denote $C_i(V_1^{\ast},V_2^{\ast})$
by $\nrhom'_i(E_{1\ast},E_{2\ast})$.
\index{$\nrhom'_1(E_{1\ast},E_{2\ast}),
 \nrhom'_2(E_{1\ast},E_{2\ast})$}
\hfill\qed
\end{notation}

When $E_a$ and $E_a/F_j(E_a)$ are locally free sheaves,
then we have
$\nbigh^i\bigl(
\nhom_1'(E_1,E_2)\bigr)=0$ $(i\neq 0)$,
and $\nbigh^0\bigl(\nhom_1'(E_1,E_2)\bigr)$
is isomorphic to the sheaf of homomorphisms
of $E_1$ to $E_2$ which preserve the filtrations.

\subsection{Geometric Invariant Theory}
\label{subsection;06.7.3.3}
\subsubsection{GIT quotient and algebraic stacks}

Let $k$ be an algebraically closed field
with characteristic $0$.
Let $G$ be a linear reductive group over $k$.
Let $Y$ be a projective variety over $k$,
provided with a $G$-action $\rho$.
Let $L$ be an ample line bundle on $Y$,
such that $\rho$ can be lifted to the action on $L$.

We recall some basic definitions.
A point $y\in Y$ is a semistable point with respect to $L$,
if there exists a $G$-invariant section
$s$ of $L^{\otimes\,n}$ for some $n>0$
such that $s(y)\neq 0$.
A point  $y\in Y$ is defined to be a stable point
with respect to $L$,
if there exists a $G$-invariant section of $L^{\otimes\,n}$
for some $n>0$
such that $s(y)\neq 0$
and that any orbits of $G$ contained
in $Y-s^{-1}(0)$ are closed.
Let $Y^s(L)$ (resp. $Y^{ss}(L)$)
denote the set of the stable (resp. semistable)
points with respect to $L$.
The foundational theorem of D. Mumford
is the following.

\begin{prop}
[\cite{GIT}]
\label{prop;06.6.21.10}
\mbox{{}}
There exists the uniform categorical quotient 
$\pi:Y\lrarr Y^{ss}/\!/G$.
Moreover the following holds:
\begin{itemize}
\item
The map $\pi$ is affine and
 universally submersive.
\item
$Y^{ss}/\!/G$ is a projective variety.
\item
There exists the open subset $Y^{s}/\!/G$
of $Y^{ss}/\!/G$,
such that $\pi^{-1}(Y^s/\!/G)=Y^s$
and that $\pi:Y^{s}\lrarr Y^{s}/\!/G$
is the universally geometric quotient of $Y^{s}$.
\end{itemize}
\end{prop}
\pf
See Proposition 1.9, Theorem 1.10
and the page 40 in \cite{GIT}.
\hfill\qed

\vspace{.1in}

We combine it with some results 
of A. Vistoli in \cite{vistoli}.
Let $Y^{sf}$ denote the set of
stable points of $Y$ whose stabilizers are 
finite.
In this situation,
we obtain the quotient stack $Y^{sf}/G$,
which is Deligne-Mumford.
See the sections 2 and 7 of the paper of Vistoli
\cite{vistoli} for more detail about such a quotient stack.
We recall one of his results here.

\begin{prop}
[\cite{vistoli}]
 \label{prop;06.4.25.1}
The naturally induced morphism 
$Y^{sf}/G\lrarr Y^{sf}/\!/G$ is proper.
\end{prop}
\pf
The map $Y^{sf}\lrarr Y^{sf}/\!/G$
is universally geometric quotient.
In particular, it is universally submersive,
and the geometric fibers are precisely
the orbits of the geometric points of $X$.
Therefore,
$Y^{sf}/\!/G$ is a quotient of $Y^{sf}$ by $G$
in the sense of Vistoli
(see the page 630 of \cite{vistoli}).
Applying Proposition 2.11 of \cite{vistoli},
we can conclude that the map
$Y^{sf}/G\lrarr Y^{sf}/\!/G$  is proper.
\hfill\qed

\begin{cor}
\label{cor;06.4.25.100}
Let $Z$ be a variety over $k$ with a $G$-action.
Let $\Phi:Z\lrarr Y$ be a $G$-equivariant immersion
with the following property:
\begin{itemize}
\item
The stabilizer group of any point of $Z$
is finite.
\item
The image $\Phi(Z)$ is contained in $Y^{s}(L)$.
\item
$\Phi:Z\lrarr Y^{ss}(L)$ is proper.
\end{itemize}
Then the quotient stack $Z/G$ is proper and
Deligne-Mumford.
\end{cor}
\pf
We have the closed substack
$Z/G$ of $Y^{sf}/G$.
We also have the closed subscheme
$Z/\!/G\lrarr Y^{sf}/\!/G\subset Y^{ss}/\!/G$.
Since $Y^{ss}/\!/G$ is projective,
$Z/\!/G$ is also projective.
From the previous lemma,
we obtain the properness of the morphism
$Z/G\lrarr Z/\!/G$.
Therefore,
we obtain the properness of $Z/G$.
\hfill\qed

\subsubsection{Mumford-Hilbert criterion
 and some elementary examples}

Let $Y$ and $G$ be as above.
Let $\lambda:G_m\lrarr G$ 
be a one-parameter subgroup.
We put $P(\lambda):=\lim_{t\to 0}\lambda(t)\cdot P$.
Then $\lambda$-acts on
the fiber $L_{|P(\lambda)}$.
The weight is denoted by
$\mu_{\lambda}(P,L)$.
The criterion says
that the point $P$ is semistable (resp. stable)
with respect to $L$,
if and only if $\mu_{\lambda}(P,L)\geq 0$
($\mu_{\lambda}(P,L)>0$)
for any one-parameter subgroup $\lambda$.

\begin{rem}
We use the convention to identify
a vector bundle and 
the sheaf of its sections.
Hence the above definition of $\mu$ 
is same as that given in {\rm\cite{GIT}}.
\hfill\qed
\end{rem}

For later use,
we recall some elementary examples.
Let $V$ be a vector space over 
an algebraically closed field $k$,
with a base $u_1,\ldots,u_N$.
Let $\lambda$ be the one-parameter subgroup
of $\SL(V)$ given by
$\lambda(t)\cdot u_i=t^{w_i}\cdot u_i$,
where $\sum w_i=0$ and $w_{i}\leq w_{i+1}$.
Let $V^{(i)}$ denote the subspace generated
by $u_1,\ldots,u_i$.
Let $V=\bigoplus V_w$ denote the weight decomposition
of $\lambda$,
i.e.,
$\lambda$ preserves the decomposition,
and the action on $V_w$ is the multiplication of $t^w$.
We put $\nbigg_j:=\bigoplus_{w\leq j}V_w$.

We denote a point of $\proj(V^{\lor})$ by $[v]$
by using a representative $v\in V-\{0\}$.
Let us consider the right $\SL(V)$-action
on $\proj(V^{\lor})$,
given by $g\cdot [v]:=[g^{-1}(v)]$,
which can be lifted to the action on
$\nbigo_{\proj(V^{\lor})}(1)$.

\begin{lem}
[\cite{GIT}]
\label{lem;06.5.15.50}
$\mu_{\lambda}\bigl([v],\nbigo_{\proj(V^{\lor})}(1)\bigr)$
is same as
$\min\bigl\{i\,\big|\,v_i\in \nbigg_i\bigr\}$.
It can be reworded as follows:
\[
  \mu_{\lambda}\bigl([v],\nbigo(1)\bigr)=
 \sum_{i}w_i\cdot
 \bigl(\dim V^{(i)}\cap\langle v\rangle
-\dim V^{(i-1)}\cap \langle v\rangle\bigr)
=\sum_j j\cdot \bigl(
 \dim \nbigg_j\cap\langle v\rangle
-\dim \nbigg_{j-1}\cap \langle v\rangle
 \bigr).
\]
Here $\langle v\rangle$ denotes the subspace
generated by $v$.
\end{lem}
\pf
According to the weight decomposition $V=\bigoplus V_i$,
we have the decomposition $v=\sum v_i$.
In $\proj(V^{\lor})$,
we have the following:
\[
 \lambda(t)[v]=\bigl[\lambda(t)^{-1}v\bigr]
=\left[
 \sum t^{-i}\cdot v_i
 \right].
\]
We put 
$i_0:=\max\bigl\{i\,|\,v_i\neq 0\bigr\}
=\min\bigl\{i\,\big|\,v\in\nbigg_i\bigr\}$.
Then it is easy to see
$\lim_{t\to 0}\lambda(t)[v]=[v_{i_0}]$.
The weight of $\lambda$ on
$\nbigo_{\proj(V^{\lor})}(1)_{|[v_{i_0}]}$ is $i_0$.
Thus the first claim is obtained.
The second claim follows from the first one.
\hfill\qed

\vspace{.1in}

Let $G_l(V)$ denote the Grassmann variety
of the $l$-dimensional subspaces of $V$:
\[
G_l(V):=\bigl\{
 \iota:W\subset V\,\big|\,\dim W=l
 \bigr\}.
\]
We have the Pl\"ucker embedding
$G_l(V)\lrarr \proj\bigl(\bigwedge^{l}V^{\lor}\bigr)$
given by 
$W\longmapsto
 \bigwedge^lW\subset \bigwedge^lV$.
It gives the polarization $\nbigo_{G_l(V)}(1)$
of $G_l(V)$.
The $SL(V)$ has the right action on $G_l(V)$
given by $\iota\longmapsto g^{-1}\circ \iota$,
which can be lifted 
to the action on $\nbigo_{G_l(V)}(1)$.

\begin{lem}
 \label{lem;06.4.22.15}
For  any point $W$ of $G_l(V)$,
we have the following:
\[
 \mu_{\lambda}\bigl(W,\nbigo_{G_l(V)}(1)\bigr)
=\sum_{i=1}^N
 \bigl(
 \rank W\cap V^{(i)}-\rank W\cap V^{(i-1)}
 \bigr)\cdot w_i
=\sum_{j\in\seisuu}
j\cdot \dim \frac{W\cap \nbigg_j}{W\cap \nbigg_{j-1}}.
\]
\end{lem}
\pf
For any $J=(j_1<j_2<\cdots<j_l)$,
we put $u_J:=u_{j_1}\wedge \cdots \wedge u_{j_l}$,
and $w_J:=\sum_{i=1}^l w_{j_i}$.
Collection of such $u_J$ gives the base of $\bigwedge^l V$.
We have the naturally induced $\SL(V)$-action
on $\bigwedge^l V$.
Let $\widetilde{\lambda}$ denote the one-parameter
subgroup  of $\SL\bigl(\Lambda^lV\bigr)$ 
induced by $\lambda$.
We have $\widetilde{\lambda}(t)(u_J)=t^{w_J}\cdot u_J$.

Let us take a base $v_1,\ldots,v_l$ of $W$
such that
$v_h=u_{i_h}+\sum_{j<i_h}a_{h,j}\cdot u_{j}$.
Then $z:=v_1\wedge\cdots\wedge v_l$
is expressed as the sum
$\sum a_J \cdot u_{J}$,
where $a_J=1$ if $J=I=(i_1<\cdots <i_l)$
and $a_J=0$ if $w_J>w_I$.
We have
$\mu_{\lambda}\bigl(W,\nbigo_{G_l}(1)\bigr)
=\mu_{\widetilde{\lambda}}\bigl(
 z,\nbigo_{\proj(\bigwedge^l V^{\lor})}(1)
 \bigr)=w_I$
due to Lemma \ref{lem;06.5.15.50}.
Then, it is easy to derive the claim of the lemma.
\hfill\qed

\vspace{.1in}
We also have the Grassmann variety
$G'_l$ of quotients of $l$-dimensions:
\[
 G'_l:=\bigl\{ q:V\lrarr Q\,\big|\,\dim Q=l\bigr\}.
\]
We have the Plucker embedding
$G'_l\lrarr \proj\bigl(\bigwedge^{N-l}V \bigr)$
given by $q\longmapsto
 \bigwedge^{N-l}q:\bigwedge^{N-l}V\lrarr \bigwedge^{N-l}Q$.
It gives the polarization $\nbigo_{G_l'}(1)$.

\begin{lem}
 [\cite{GIT}, \cite{my}]
Let $q:V\lrarr Q$ be a point of $G_l'$.
We put $W:=\Ker(q)$.
Then we have the following:
\[
 \mu_{\lambda}\bigl(q,\nbigo_{G_l'}(1)\bigr)
=\sum_{i=1}^N
 w_i\cdot\bigl(\dim V^{(i)}\cap W-\dim V^{(i-1)}\cap W
 -1\bigr)
=\sum_{j=1}^N
 j\cdot \left(
 \dim \frac{W\cap \nbigg_j}{W\cap\nbigg_{j-1}}
-\dim\frac{\nbigg_j}{\nbigg_{j-1}}
 \right).
\]
\end{lem}
\pf
We put $W^{(i)}:=\nbigg_i\cap W\big/\nbigg_{i-1}\cap W$.
Since we have the natural isomorphism
$\nbigg_i/\nbigg_{i-1}\simeq V_i$,
we regard $W^{(i)}$ as the subspace of $V_i$.
It is easy to see that
the limit $\lim_{t\to 0}\lambda(t)\cdot q$ is given by
the quotient $\widehat{q}:V\lrarr \bigoplus V_i/W^{(i)}$.
The weight of $\lambda$
on $\nbigo_{G_l'}(1)_{|\widehat{q}}$ is
$-i\cdot \dim\bigl(V_i/W^{(i)}\bigr)$.
Then it is easy to derive the claim.
\hfill\qed

\begin{rem}
We have the obvious isomorphism
$G_l(V)\simeq G_{N-l}'(V)$.
It does not preserve the semistability conditions
on the varieties 
induced by the Pl\"ucker embeddings.
\hfill\qed
\end{rem}

\subsection{Cotangent Complex}
\label{subsection;06.7.3.4}
\subsubsection{Basis}

Recall some fundamental properties of the cotangent complexes
from \cite{ill}, \cite{lau} and \cite{olsson}.
Let $\nbigx$ and $\nbigy$ be Deligne-Mumford stacks
with the etale site.
For any morphism $f:\nbigx\lrarr \nbigy$
of Delinge-Mumford stacks,
the cotangent complex 
was constructed by L. Illusie \cite{ill}
as a complex of $\nbigo_{\nbigx}$-modules,
which is denoted by 
$L_{\nbigx/\nbigy}$ or $L_f$.
Recall that the cotangent complex controls
the deformation theory (Section 3 \cite{ill})
in the following sense.
Let $T$ be a scheme over $\nbigy$,
and let $h:T\lrarr \nbigx$ be a $\nbigy$-morphism.
Let $\Tbar$ be a $\nbigy$-scheme such that
$T$ is a closed $\nbigy$-subscheme of $\Tbar$
and the corresponding ideal $J$ is square zero,
i.e., $J^2=\{f\cdot g\,\big|\,f,g\in J\}=0$.

\begin{prop}
 [Illusie, \cite{ill}]
 \label{prop;06.4.29.5}
We have the obstruction class $o(h)$
in $Ext^1(h^{\ast}L_{\nbigx/\nbigy},J)$
with the following property:
\begin{itemize}
\item
The morphism $h$ can be extended over
$\Tbar$, if and only if $o(h)$ vanishes.
\end{itemize}
When $o(h)=0$,
the set of the extension classes is the torsor over
the group $Ext^0\bigl(h^{\ast}L_{\nbigx/\nbigy},J\bigr)$.
\hfill\qed
\end{prop}

The cotangent complex has a nice functorial property.
For example,
if $g:\nbigy\lrarr \nbigz$ be a morphism,
then we have the distinguished triangle,
$f^{\ast}L_{\nbigy/\nbigz}\lrarr
 L_{\nbigx/\nbigz}\lrarr L_{\nbigx/\nbigy}
 \lrarr f^{\ast}L_{\nbigy/\nbigz}[1]$
in the derived category $D(\nbigx)$.

\vspace{.1in}

As for general Artin stacks with the lisse-\'{e}tale site,
the cotangent complex with some good functorial property
was obtained by G. Laumon, L. Moret-Bailly and M. Olsson
(Section 17 of \cite{lau} and Section 8 of \cite{olsson}).
For any Artin stack $\nbigx$,
Olsson introduced the category $D'_{\qcoh}(\nbigx)$
of the projective systems
$K=(\cdots \rarr K_{\geq -n-1}\rarr K_{\geq -n}\rarr
 \cdots \rarr K_{\geq\,0})$ in 
$D^{+}(\nbigx)$
such that the morphisms
$K_{\geq\,-n}\lrarr \tau_{\geq-n}K_{\geq-n}$
and 
$\tau_{\geq\,-n}K_{\geq\,-n-1}
\lrarr \tau_{\geq\,-n}K_{\geq\,-n}$
are isomorphisms.
Here $\tau_{\geq\,-n}$ denotes the canonical
$n$-th truncation functor.
See \cite{olsson} for the functorial property
of $D'_{\qcoh}(\nbigx)$.
Let $f:\nbigx\lrarr \nbigy$ be a quasi-compact
and quasi-separated morphism of Artin stacks.
Then, we can associate
$L_{\nbigx/\nbigy}=L_f=
 \bigl(\cdots\rarr L_{\nbigx/\nbigy}^{\geq\,-n-1}
 \rarr L_{\nbigx/\nbigy}^{\geq\,-n}\rarr
 \cdots L_{\nbigx/\nbigy}^{\geq\,0}\bigr)
\in D'_{\qcoh}(\nbigx)$ to $f$
with the following property
 (Theorem 8.1 \cite{olsson}):
\begin{itemize}
\item
If $\nbigx$ and $\nbigy$ are algebraic spaces,
$L_{\nbigx/\nbigy}^{\geq\,-n}$ is isomorphic to
$\tau_{\geq\,-n}L_{\nbigx/\nbigy}$ in
$D^{+}_{\qcoh}(\nbigx)$.
Here the latter $L_{\nbigx/\nbigy}$
denotes the usual cotangent complex
defined by Illusie. 
\item
When we are given a $2$-commutative diagram 
of Artin stacks,
\[
 \begin{CD}
 \nbigx'@>{f}>>\nbigx \\
 @V{g}VV @VVV \\
 \nbigy'@>{h}>>  \nbigy\\
 \end{CD}
\]
we have the functorial morphism
$Lf^{\ast}L_{\nbigx/\nbigy}\lrarr 
 L_{\nbigx'/\nbigy'}$.
If the diagram is 2-Cartesian, and if one of $g$ or $h$ is flat,
then the morphism 
$Lf^{\ast}L_{\nbigx/\nbigy}\lrarr L_{\nbigx'/\nbigy'}$
is isomorphic.
\item
 Let $f:\nbigx\lrarr \nbigy$ be a morphism of Artin stacks.
 Let $g:\nbigy\lrarr \nbigz$ be another morphism.
 Then we have the distinguished triangle
 $Lf^{\ast}L_{\nbigy/\nbigz}\lrarr
 L_{\nbigx/\nbigz}\lrarr L_{\nbigx/\nbigy}
 \lrarr Lf^{\ast}L_{\nbigy/\nbigz}[1]$
 in $D'_{\qcoh}(\nbigx)$.
\end{itemize}
The following properties can be derived 
directly from the construction.
(See Section 8 of \cite{olsson}
for the construction of $L_{\nbigx/\nbigy}$.)
\begin{itemize}
\item
Each $L_{\nbigx/\nbigy}^{\geq\,-n}$
is an object in $D^{[-n,1]}_{\qcoh}(\nbigx)$.
\item
If $f$ is smooth and representable,
then $L_{\nbigx/\nbigy}$ is quasi isomorphic to the $0$-th cohomology sheaf,
which is isomorphic to the locally free sheaves of Kahler differentials
$\Omega_{\nbigx/\nbigy}$.
 In general, if $f$ is smooth,
 any $L_{\nbigx/\nbigy}^{\geq\,-n}$ is
 of perfect amplitude contained in $[0,1]$.
In particular, they are isomorphic to
 $L_{\nbigx/\nbigy}^{\geq \,0}$.
\end{itemize}

\begin{rem}
M. Aoki {\rm(\cite{aoki})}
generalized the deformation theory of Illusie.
He showed that a generalization of
Proposition {\rm{\ref{prop;06.4.29.5}}}
holds for the Artin stacks.
\hfill\qed
\end{rem}

\subsubsection{Quotient stacks}
\label{subsubsection;06.4.29.15}

Let $G$ be a smooth group $S$-scheme.
Let $Y$ be a smooth $S$-scheme with a $G$-action.
The quotient stack is denoted by $Y_G$.
Let $f:Y\lrarr Y_G$ be a morphism.
We have the corresponding $G$-torsor $P(f)$
over $Y$.
Since $f$ is smooth and representable,
the cotangent complex $L_f$ is isomorphic
to the sheaf $\Omega_f$ of the relative Kahler differentials.
\begin{lem}
 \label{lem;06.4.30.2}
$\Omega_f$ is isomorphic to
the sheaf of 
the $G$-invariant sections of $\Omega_{P(f)/Y}$.
\end{lem}
\pf
Let $\pi:Y\lrarr Y_G$ denote the morphism
corresponding to the trivial torsor.
We have the following diagram:
\begin{equation}
 \begin{CD}
 P(f)\times_YP(f) @>{\ptilde_1,\ptilde_2}>> P(f) @>{\pitilde}>> Y\\
 @VV{f_2}V @VV{f_1}V @VV{f}V \\
 Y\times_{Y_G}Y @>{p_1,p_2}>> Y @>{\pi}>> Y_G
 \end{CD}
\end{equation}
Then the sheaf $\Omega_{f}$ is the descent
of $\Omega_{f_1}$ by the cocycle condition
$\ptilde_1^{\ast}\Omega_{f_1}
\simeq
 \Omega_{f_2}
\simeq
 \ptilde_2^{\ast}\Omega_{f_2}$,
where the isomorphisms are given
by the naturally defined differentials.
We have $P(f)\times_{Y}P(f)\simeq P(f)\times G$
for which $p_1$ and $p_2$ correspond
to the natural projection
and the $G$-action respectively.
Then the claim is obvious.
\hfill\qed

\vspace{.1in}

Let $Z$ be an Artin stack over $S$
with a morphism $F:Z\lrarr Y_G$.
We have the corresponding $G$-torsor
$P(F)$ over $Z$
and the $G$-equivariant map $\widetilde{F}:P(F)\lrarr Y$.
\[
 \begin{CD}
 P(F) @>{\widetilde{F}}>> Y\\
 @VVV@V{\pi}VV \\
 Z @>{F}>> Y_G
 \end{CD}
\]
Let us describe the pull back of the cotangent
complex $F^{\ast}L_{Y_G/S}$
on $Z$.
We have the map
$\alpha:\widetilde{F}^{\ast}\Omega_{Y/S}
 \lrarr 
 \Omega_{P(F)/Z}$
on $P(F)$,
which is the composite of the differential
$\Ftilde^{\ast}\Omega_{Y/S}
\lrarr \Omega_{P(F)/S}$
and the natural projection
$\Omega_{P(F)/S}\lrarr \Omega_{P(F)/Z}$.

\begin{prop}
 \label{prop;06.4.30.1}
$F^{\ast}L_{Y_G/S}$ is represented 
by the decent of $\Cone(-\alpha)[-1]$
with respect to the natural $G$-action.
\end{prop}
\pf
We recall the construction of $L_{Y_G/S}$
in this case.
We put $Y^{(m)}:=
 \overbrace{Y\times_{Y_G}\cdots\times_{Y_G}Y}^{m+1}$.
We have the natural morphisms
$Y^{(m)}\lrarr Y_G\lrarr S$.
We have the complexes
$C^{(m)}:=\bigl(
\Omega_{Y^{(m)}/S}\lrarr \Omega_{Y^{(m)}/Y_G}
\bigr)$
on $Y^{(m)}$,
where $\Omega_{Y^{(m)}/S}$ stands in the degree $0$.
We have the strictly simplicial structure
given by the naturally defined quasi isomorphisms
$\pi_i:\pi^{\ast}C^{(m-1)}\lrarr C^{(m)}$
 $(i=0,1,\ldots,m)$.
Then $L_{Y_G/S}\in D'_{\qcoh}(Y_G)$ is obtained
as the gluing of
$\bigl(C^{(m)}\,\big|\,m=0,1,\ldots\bigr)$ simplicially.

We put $P(F)^{(m)}:=
 \overbrace{P(F)\times_Z\cdots\times _ZP(F)}^{m+1}$.
We have the naturally defined morphisms
$F^{(m)}:P(F)^{(m)}\lrarr Y^{(m)}$.
Then $F^{\ast}L_{Y_G/S}$ is obtained as the gluing
of $\bigl(
F^{(m)\ast}C^{(m)}\,\big|\,m=0,1,\ldots\bigr)$.
We have the following commutative diagram:
\[
 \begin{CD}
 F^{(m)\,\ast}\Omega_{Y^{(m)}/S}
 @>>> F^{(m)\,\ast}\Omega_{Y^{(m)}/Y_G}\\
 @V{=}VV @V{\simeq}VV \\
 F^{(m)\,\ast}\Omega_{Y^{(m)}/S}
 @>>> \Omega_{P(F)^{(m)}/Z}
 \end{CD}
\]
Here the bottom morphism is 
same as the composite
$F^{(m)\,\ast}\Omega_{Y^{(m)}/S}
\lrarr
 \Omega_{P(F)^{(m)}/S}
\lrarr \Omega_{P(F)^{(m)}/Z}$,
where the first one is the differential
and the second one is the natural projection.

Let $q_i:Y\times G^m\lrarr Y$ $(i=0,1,\ldots,m)$
be a morphism given by
$q_i(y,g_1,\ldots,g_m)=y\cdot g_1\cdot\cdots\cdot g_{i}$.
They induce the isomorphism
$Y\times G^m\lrarr Y^{(m)}$.
Under the identification,
$q_i$ is the projection onto the $i$-th component.
Similarly,
we have the identification
$P(F)\times G^m\simeq P(F)^{(m)}$,
under which $F^{(m)}$ is given by
$F^{(m)}(y,g_1,\ldots,g_m)=
 \bigl(\widetilde{F}(y),g_1,\ldots,g_m\bigr)$.
Let $\rho_m$ denote the projection of
$P(F)\times G^m$ onto $G^m$.
We have the subcomplex
$\bigl(\rho_m^{\ast}\Omega_{G^m}
\stackrel{\id}{\lrarr}
 \rho_m^{\ast}\Omega_{G^m}\bigr)$
of $F^{\ast}C^{(m)}$.
It is compatible with the simplicial structure.
The quotients are denoted by
$\widehat{C}^{(m)}$,
and then the gluing of 
$\bigl(\widehat{C}^{(m)}\,\big|\,m=0,1,\ldots\bigr)$
also gives
$F^{\ast}L_{Y_G/S}$ in $D(Z)$.
Then, it follows that
$F^{\ast}L_{Y_G/S}$ is given as the descent of
$\widehat{C}^{(0)}= 
\bigl(\widetilde{F}^{\ast}\Omega_{Y/S}
\stackrel{\alpha}{\lrarr} \Omega_{P(F)/Z}\bigr)$
with respect to the natural $G$-action.
\hfill\qed

\vspace{.1in}
Let $H$ denote the composite of $F$
and the canonical map $Y_G\lrarr S_G$.
Let $P(H)$ denote the $G$-torsor over $Z$
corresponding to $h$.
Since we have the natural isomorphism
$P(H)\simeq P(F)$,
we do not distinguish them.
Let $\widetilde{H}:P(F)\lrarr S$ be the lift of $H$.
Let $\pi$ denote the projection $P(F)\lrarr Z$.
We have the canonical isomorphism
$\pi^{\ast}H^{\ast}L_{S_G/S}[1]\simeq
 \widetilde{H}^{\ast}\Omega_{S/S_G}
 \simeq \Omega_{P(F)/Z}$.
We also have the canonical isomorphism
$\pi^{\ast}F^{\ast}L_{Y_G/S_G}
\simeq
 \widetilde{F}^{\ast}\Omega_{Y/S}$.
We obtain the following corollary.

\begin{cor}
\label{cor;06.5.9.2}
The morphism
$F^{\ast}L_{Y_G/S_G}\lrarr H^{\ast}L_{S_G/S}[1]$
on $Z$ is obtained as the descent of
$\alpha:\widetilde{F}^{\ast}\Omega_{Y/S}
\lrarr \Omega_{P(F)/Z}$.
\end{cor}
\pf
We have the distinguished triangle
$H^{\ast}L_{S_G/S}\lrarr 
 F^{\ast}L_{Y_G/S}\lrarr F^{\ast}L_{Y_G/S_G}
\lrarr H^{\ast}L_{S_G/S}[1]$.
Due to Proposition \ref{prop;06.4.30.1},
we know the morphism
$H^{\ast}L_{S_G/S}\lrarr F^{\ast}L_{Y_G/S}$.
Then,
we know the morphism
$F^{\ast}L_{Y_G/S}\lrarr H^{\ast}L_{S_G/S}[1]$.
\hfill\qed

\begin{example}
Let $E$ be a vector bundle on 
a variety $X$ of rank $R$,
and let $f:X\lrarr k_{\GL(R)}$
be the corresponding map.
Then we have
$f^{\ast}L_{k_{GL(R)}/k}\simeq \End(E)[-1]$.
\hfill\qed
\end{example}

\begin{rem}
The expression in Proposition {\rm\ref{prop;06.4.30.1}}
is natural, in the following sense.
Let $Y_i$ $(i=1,2)$ be $S$-schemes
with $G$-actions,
and let $g:Y_1\lrarr Y_2$
be a $G$-equivariant morphism.
Let $g_G:Y_{1\,G}\lrarr Y_{2\,G}$ denote 
the induced morphism.
Let $h_1:Z\lrarr Y_{1\,G}$ be a morphism.
The composite $g_G\circ h_1$ is denoted by $h_2$.
We have the corresponding torsor $P$ over $Z$
and the $G$-equivariant morphisms
$\widetilde{h}_i:P\lrarr Y_i$.
We have the natural commutative diagram
of $G$-equivariant sheaves on $P$:
\[
\begin{CD}
 \widetilde{h}_2^{\ast}\Omega_{Y_2/S}
 @>{\alpha_2}>>
 \Omega_{P/Z}\\
 @VVV @VVV \\
 \widetilde{h}_1^{\ast}\Omega_{Y_1/S}
 @>{\alpha_1}>>
 \Omega_{P/Z}
\end{CD}
\]
Then the morphism
 $h_2^{\ast}L_{Y_{2\,G}/S}\lrarr 
 h_1^{\ast}L_{Y_{1\,G}/S}$ is 
the descent of the induced morphism
$\Cone(-\alpha_2)[-1]\lrarr \Cone(-\alpha_2)[-1]$.
\hfill\qed
\end{rem}

\begin{rem}
\label{rem;06.5.22.500}
Let $G_1$ be a smooth group scheme over $S$.
Assume that $Y$ is provided with the $G_1$-action,
which commutes with the $G$-action.
It induces the $G_1$-action on $Y_G$.
Moreover, assume that $Z$ is also provided
with the $G_1$-action
such that $F$ is $G_1$-equivariant.
Then, we have the naturally induced $G_1$-actions
on the complex $\cone(-\alpha)[-1]$,
which commutes with the $G$-action.
It induces the $G_1$-action
on the descent of $\cone(-\alpha)[-1]$ on $Z$.
In particular,
we obtain the $G_1$-equivariant representative
of $F^{\ast}L_{Y_G/S}$.
\hfill\qed
\end{rem}

Let $\pi:Y\lrarr Y_G$ denote the canonical projection.
Due to Proposition \ref{prop;06.4.30.1},
$L_{Y_G/S}$ on $Y_G$ is the descent of
$\bigl(\Omega_{Y/S}\stackrel{\alpha}{\lrarr}
 \Omega_{Y/Y_G}\bigr)$
given on $Y$ with the natural $G$-action,
where $\Omega_{Y/S}$ stands in the degree $0$.
For simplicity, we consider the case $S=\Spec(k)$.
Due to Lemma \ref{lem;06.4.30.2},
we have $\Omega_{Y/Y_G}\simeq 
\gminig^{\lor}\otimes\nbigo_Y$,
where $\gminig$ denotes the tangent space
of $G$ at the unit element,
or equivalently the vector space of the right invariant
vector fields, and $\gminig^{\lor}$ denotes the dual.
Let $\Theta_{Y/S}$ denote the relative tangent sheaf
of $Y/S$. The $G$-action on $Y$
induces the map
$A:\gminig\otimes\nbigo_Y\lrarr  \Theta_{Y/S}$.

\begin{lem}
 \label{lem;06.4.30.3}
The map
$\alpha:\Omega_{Y/S}\lrarr \gminig^{\lor}\otimes\nbigo_Y$
is given by the dual of $-A$.
Namely,
we have
$\pi^{\ast}L_{Y_G/S}
\simeq
 \Cone(A)[-1]$.
\end{lem}
\pf
Let $p_i:Y\times_{Y_G}Y\lrarr Y$ denote 
the projection on the $i$-th component.
We have the following factorization
of $p_1^{\ast}\alpha$:
\[
 p_1^{\ast}\Omega_{Y/S}
\lrarr
 \Omega_{Y\times_{Y_G}Y/S}
\lrarr
 \Omega_{p_2}
\simeq
 p_1^{\ast}\Omega_{Y/Y_G}
\]
Each morphisms are induced by the natural differential.
Let us take the identification
$Y\times_{Y_G}Y\simeq Y\times G$,
for which $p_1$ and $p_2$
correspond to the natural projection onto $Y$
and the $G$-action,
respectively.

Let $y$ be any closed point of $Y$,
and let $e$ be the identity element of $G$.
We have $p_1(y,e)=p_2(y,e)=y$.
We denote the differential of $p_i$
at $(y,e)$ by $T_{(y,e)}p_i$.
Let us consider the specialization of
the dual of $p_1^{\ast}\alpha$ at $(y,e)$.
Then it is the composite of
the inclusion of
$\Ker(T_{(y,e)}p_{2})
\subset T_{(y,e)}(Y\times G)$
and the natural projection
$T_{(y,e)}Y\times G\lrarr T_yY$.
Since we have 
$\Ker\bigl(T_{(y,e)}p_2\bigr)
\simeq \bigl\{
 (-Av,v)\,\big|\,v\in \gminig \bigr\}
\simeq
 \gminig$,
the map is same as $-A$.
Since $\alpha$ can be recovered
from $p_1^{\ast}\alpha$,
we are done.
\hfill\qed

\begin{rem}
Since $f^{\ast}L_{Y_G/S}$ is obtained as the decent of
$\widetilde{f}^{\ast}\Cone(A)[-1]$
for a morphism $f:Z\lrarr Y_G$,
Lemma {\rm\ref{lem;06.4.30.3}} can be used 
in the calculation.
\hfill\qed
\end{rem}

\begin{example}
\label{example;06.4.30.10}
Let $W_i$ $(i=-1,0)$ be $R_i$-dimensional vector spaces
over $k$.
Let $N(W_{-1},W_0)$ denote the vector space
of linear maps from $W_{-1}$ to $W_0$.
We have the right $GL(W_{-1})\times \GL(W_0)$-action
on $N(W_{-1},W_0)$ given by
$(g_{-1},g_0)\cdot f=g_0^{-1}\circ f\circ g_{-1}$.
Hence we obtain the quotient space
$Y(W_{\cdot}):=
 N(W_{-1},W_0)_{\GL(W_{-1})\times \GL(W_{0})}$.

Let $X$ and $U$ be stacks over $k$.
Let $V_i$ $(i=-1,0)$ be vector bundles
on $U\times X$ whose ranks are $R_i$.
Let $f:V_{-1}\lrarr V_0$ be a morphism of 
$\nbigo_{U\times X}$-modules.
Then we obtain the morphism
$\Phi_f:U\times X\lrarr Y(W_{\cdot})$.
The pull back of the cotangent complex
$\Phi_f^{\ast}L_{Y(W_{\cdot})/k}$ is quasi isomorphic to
the following complex:
\[
 \nhom\bigl(V_{0},V_{-1}\bigr)
\stackrel{\alpha}{\lrarr}
 \nhom\bigl(V_0,V_0\bigr)\oplus\nhom(V_{-1},V_{-1}).
\]
Here $\nhom\bigl(V_{0},V_{-1}\bigr)$ stands in
degree $0$,
and the map $\alpha$ is given by
$ \alpha(a)=\bigl(f\circ a,-a\circ f\bigr)$.
We remark that it is isomorphic to
$\nhom\bigl(V_{\cdot},V_{\cdot}\bigr)^{\lor}_{\leq 0}[-1]$.
(See the subsubsection {\rm\ref{subsubsection;06.4.29.30}}.)

Actually, we have only to care the signature.
We can see it formally.
Let $f$ be an element of $N(W_{-1},W_{0})$.
The differential of the action of $\GL(W_{-1})\times \GL(W_{0})$
gives the map:
\begin{equation}
\label{eq;06.6.21.1}
 \End(W_{-1})\oplus \End(W_{0})\lrarr
 T_fN(W_{-1},W_0)=N(W_{-1},W_0),
\quad
 (a_{-1},a_0)\longmapsto
 -a_0\circ f+f\circ a_{-1}
\end{equation}
If we regard $W_{-1}\stackrel{f}{\lrarr}W_0$
as a complex,
{\rm(\ref{eq;06.6.21.1})} can be regarded as
$\nhom(W_{\cdot},W_{\cdot})_{\geq\,0}$.
Then, Lemma {\rm\ref{lem;06.4.30.3}} says that
the cotangent complex corresponds
to $\bigl(
\nhom(W_{\cdot},W_{\cdot})_{\geq\,0}\bigr)^{\lor}[-1]$.
\hfill\qed
\end{example}

Let us consider the following diagram:
\[
 \begin{CD}
 Y @>{\psitilde}>> S \\
 @V{\pi}VV @V{\pi_1}VV \\
 Y_{G}@>{\psi}>> S_G
 \end{CD}
\]
We have the natural isomorphisms
$\pi^{\ast}L_{Y_G/S_G} \simeq L_{Y/S}$
and
$\pi^{\ast}\psi^{\ast}L_{S_G/S}[1]
\simeq
 \psitilde^{\ast}L_{S/S_G}
\simeq
 L_{Y/Y_G}$.
We identify them by the isomorphisms.
\begin{lem}
\label{lem;06.5.8.1}
Under the identification above,
the morphism
$\pi^{\ast}L_{Y_G/S_G}\lrarr \pi^{\ast}\psi^{\ast}L_{S_G/S}[1]$
is same as the natural morphism
$L_{Y/S}\lrarr L_{Y/Y_G}$ .
\end{lem}
\pf
We have the natural isomorphisms:
\[
 \pi^{\ast}L_{Y_G/S}\simeq
 \Cone\bigl(L_{Y/S}\lrarr L_{Y/Y_G}\bigr)[-1],
\quad
\pi^{\ast}\psi^{\ast}L_{S_G/S}\simeq
 \psitilde^{\ast}L_{S/S_G}[-1]
\]
The morphism
$\pi^{\ast}\psi^{\ast}L_{S_G/S}
\lrarr
 \pi^{\ast}L_{Y_G/S}$ is induced by
$\psitilde^{\ast}L_{S/S_G}\lrarr L_{Y/Y_G}$.
Hence, the distinguished triangle
$\pi^{\ast}\psi^{\ast}L_{S_G/S}
\lrarr \pi^{\ast}L_{Y_G/S}
\lrarr \pi^{\ast}L_{Y_G/S_G}
\lrarr \pi^{\ast}\psi^{\ast} L_{S_G/S}[1]$
is the following:
\[
 \psi^{\ast}L_{S/S_G}[-1]
\lrarr
 \Cone\bigl(
 L_{Y/S}\rarr L_{Y/Y_G}
 \bigr)[-1]
\lrarr
 L_{Y/S}
\lrarr 
 \psi^{\ast}L_{S/S_G}
\]
Then the claim of the lemma follows.
\hfill\qed

\subsubsection{Some more examples}
\label{subsubsection;06.5.25.2}

The result in this subsubsection
will be used in the subsections \ref{subsection;06.6.1.10}
and \ref{subsection;06.6.1.11}.
The author recommends the reader
to skip here. 
Let $X$ be a smooth projective surface,
and let $U_2$ be a quasi compact stack.
We consider a pair of $U_2$-coherent sheaf $\nbigf$
and a section $\varphi:\nbigo_{U_2\times X}\lrarr \nbigf$.
We assume that $p_{X\,\ast}\nbigf$ is locally free.
We have the induced section
$\nbigh(\varphi):\nbigo_{U_2\times X}\lrarr p_{X\,\ast}\nbigf$.
Let $U_1$ be a subscheme of $U_2$
contained in the $0$-set of $\nbigh(\varphi)$.

Assume we are given a datum 
$(V_{\cdot},P_{\cdot},\phi,\phitilde)$
as follows:
\begin{itemize}
\item
A locally free resolution of $\nbigf$:
\[
\begin{CD}
 V_{-1}@>{d_{-1}}>>
 V_0@>{d_{0}}>>
 V_1 @>{\epsilon}>> 
 \nbigf
\end{CD}
\]
Namely, the sequence
$0\lrarr V_{-1}\lrarr V_0\lrarr V_1$ is exact,
and $V_0/V_1\simeq \nbigf$.
\item
A morphism $\phi:\nbigo_{U_2\times X}\lrarr V_{1}$
such that $\epsilon\circ\phi=\varphi$.
Such a $\phi$ is called a lift of $\varphi$.
\item
A resolution
$P_{\cdot}=(P_{-1}\stackrel{\del}{\rarr} P_0)$ 
of $\nbigo_X$,
i.e., $P_{0}/P_{-1}\simeq\nbigo_X$.
\item
A morphism $\phitilde_{i}:P_{i}\lrarr V_i$ $(i=0,-1)$
such that the following diagram is commutative:
\[
 \begin{CD}
 V_{-1} @>{d_{-1}}>> V_0 @>{d_0}>> V_1  \\
 @A{\phitilde_{-1}}AA 
 @A{\phitilde_0}AA 
 @A{\phi}AA \\
 P_{-1} @>{\del}>> P_0 @>>> \nbigo
 \end{CD}
\]
\end{itemize}
We put 
$\gminik(V_{\cdot},P_{\cdot},\phi,\phitilde)
:=\nhom(p_{U_1}^{\ast}P_{\cdot},
 V_{\cdot|U_1\times X})^{\lor}$
on $U_1\times X$,
and we will construct a morphism
$\gminir(V_{\cdot},P_{\cdot},\phi,\phitilde):
\gminik(V_{\cdot},P_{\cdot},\phi,\phitilde)
\lrarr L_{U_1\times X/U_2\times X}$ on $U_1$,
which depends only on $(\nbigf,\varphi)$
in the derived category
$D(U_1\times X)$.

We have the vector bundles
$N(P_i,V_j)$ on $U_2\times X$.
We have the following smooth morphism:
\[
 h_1:N(P_0,V_0)\times_X N(P_{-1},V_{-1})
\lrarr N(P_{-1},V_0),
\quad
 h(a_0,a_{-1})=f\circ a_{-1}-a_0\circ \del
\]
We put $Z_1(V_{\cdot},P_{\cdot}):=h^{-1}(0)$.
We also have the following morphism:
\[
 h_2:N(P_{0},V_1)\lrarr N(P_{-1},V_1),
\quad
 h_2(a)=-a\circ \del.
\]
We put $Z_2(V_{\cdot},P_{\cdot}):=h_2^{-1}(0)$.
We remark $Z_2(V_{\cdot},P_{\cdot})$ is naturally isomorphic to
$N(\nbigo,V_1)$.

We have the naturally defined morphism
$\Gamma:Z_1(V_{\cdot},P_{\cdot})\lrarr Z_2(V_{\cdot},P_{\cdot})$.
The morphisms $\phi$ and $\phitilde$ give
the sections $\Phi_i:U_i\times X\lrarr Z_i$,
and we have the following commutative diagram:
\begin{equation}
\label{eq;06.5.24.1}
 \begin{CD}
 U_1\times X @>{\Phi_1}>> Z_1(V_{\cdot},P_{\cdot}) @>>>X\\
 @V{j_X}VV @VVV @VVV\\
 U_2\times X @>{\Phi_2}>> Z_2(V_{\cdot},P_{\cdot})@>>>X
 \end{CD}
\end{equation}
It is easy to see that
$\Phi_1^{\ast}
 L_{Z_1(V_{\cdot},P_{\cdot})/Z_2(V_{\cdot},P_{\cdot})}$ is expressed 
by $\gminik(V_{\cdot},P_{\cdot},\phitilde,\phi)_{\leq 0}$.
The composite
$\gminik(V_{\cdot},P_{\cdot},\phitilde,\phi)
\lrarr \Phi_1^{\ast}L_{Z_1(V_{\cdot},P_{\cdot})/Z_2(V_{\cdot},P_{\cdot})}
\lrarr L_{U_1\times X/U_2\times X}$
is denoted by
$\gminir(V_{\cdot},P_{\cdot},\phitilde,\phi)$.

\begin{lem}
$\gminir(V_{\cdot},P_{\cdot},\phitilde,\phi)$
and $\gminik(V_{\cdot},P_{\cdot},\phitilde,\phi)$
depend only on $(\nbigf,\varphi)$
in the derived category $D(U_1\times X)$.
\end{lem}
\pf
Let $(V^{(i)},P^{(i)},\phi^{(i)},\phitilde^{(i)})$
$(i=1,2)$ be data as above.
Let $V^{(3)}_1$ be the cokernel
of $(\phi^{(1)},-\phi^{(2)}):
 \nbigo_{U_2\times X}\lrarr V^{(1)}_1\oplus V^{(2)}_1$.
Then $\epsilon^{(i)}$ $(i=1,2)$ naturally induce
the morphism $\epsilon^{(3)}:V^{(3)}_1\lrarr \nbigf$.
We also have the morphism
$\phi^{(3)}:\nbigo_{U_2\times X}\lrarr V^{(3)}_1$,
naturally.
We take a surjection $A\lrarr \Ker(\epsilon^{(3)})$,
appropriately.
We have the naturally defined morphism
$V_0^{(1)}\oplus V_0^{(2)}\lrarr \Ker(\epsilon^{(3)})$.
We put $V_0^{(3)}:=A\oplus V_0^{(1)}\oplus V_0^{(2)}$,
and then we have the natural morphism
$d_0^{(3)}:V_0^{(3)}\lrarr V_1^{(3)}$.
The kernel of $d_0^{(3)}$ is locally free,
and we put $V_{-1}^{(3)}:=\Ker(d_0^{(3)})$.
We put $P_0^{(3)}:=P_0^{(1)}\oplus P_0^{(2)}$,
and $P_{-1}^{(3)}$ denotes the kernel of the natural morphism
$P_0^{(1)}\oplus P_0^{(2)}\lrarr \nbigo_{U_2\times X}$.
Then, we have the following compatible
diagrams on $U_2\times X$:
\[
 \begin{CD}
 V_{-1}^{(3)}@>>> V_{0}^{(3)} @>>> V_{1}^{(3)} @>>>\nbigf\\
 @A{a_3^{(i)}}AA @A{a_2^{(i)}}AA @A{a_1^{(i)}}AA @A{=}AA \\
 V_{-1}^{(i)} @>>> V_{0}^{(i)} @>>> V_1^{(i)} @>>>\nbigf
 \end{CD}
\quad\quad\quad\quad
 \begin{CD}
 P_{-1}^{(3)} @>>> P_0^{(3)} @>>> \nbigo\\
 @AAA @AAA @AAA \\
 P_{-1}^{(i)} @>>> P_0^{(i)} @>>>\nbigo
 \end{CD}
\]
The cokernels of $a_j^{(i)}$ are locally free.
The composite of $\phi^{(i)}$ and $a_1^{(i)}$ is same as
$\phi^{(3)}$.
On $U_1\times X$,
$a_j^{(i)}$ and $\phitilde_j$ are compatible.
Then, we have only to compare 
$\gminir(V^{(1)},P^{(1)},\phi^{(1)},\phitilde^{(1)})$
and $\gminir(V^{(3)},P^{(3)},\phi^{(3)},\phitilde^{(3)})$.

We give only an indication.
We regard
$V_j^{(1)}\subset V_j^{(3)}$ and
$P_j^{(1)}\subset P_j^{(3)}$ as the filtrations.
Let $N'(P_{j}^{(3)},V_{k}^{(3)})$ denote 
the vector bundle corresponding to the locally free sheaves
of filtration-preserving homomorphisms of 
$P_{j}^{(3)}$ to $V_{k}^{(3)}$.
We construct the vector bundles
$Z'_i(V_{\cdot},P_{\cdot})$ by using $N'(P^{(3)}_i,V^{(3)}_{j})$
instead of
$N(P_i^{(3)},V_j^{(3)})$.
Then, we have the morphisms:
$Z_i(V^{(3)}_{\cdot},P^{(3)}_{\cdot})
\lrarr
 Z'_i(V_{\cdot}^{(3)},P_{\cdot}^{(3)})$
and 
$Z_i(V_{\cdot}^{(3)},P_{\cdot}^{(3)})
\lrarr
Z_i(V_{\cdot}^{(1)},P_{\cdot}^{(1)})$.

Let $\nhom'(P_{\cdot},V_{\cdot})$
denote the subcomplex of 
$\nhom(P_{\cdot},V_{\cdot})$
which consists of filtration-preserving homomorphisms.
Let $\gminik'(V_{\cdot}^{(3)},P_{\cdot}^{(3)},
 \phi^{(3)},\phitilde^{(3)})$
denote the dual
$\nhom'(P^{(3)}_{\cdot},V_{\cdot})^{\lor}$.
We obtain the following diagram on $U$:
\[
 \begin{CD}
\gminik(V_{\cdot}^{(3)},P_{\cdot}^{(3)},\phi^{(3)},\phitilde^{(3)})
@>{b_1}>>
\gminik'(V_{\cdot}^{(3)},P_{\cdot}^{(3)},\phi^{(3)},\phitilde^{(3)})
@<{b_2}<<
\gminik(V_{\cdot}^{(1)},P_{\cdot}^{(1)},\phi^{(1)},\phitilde^{(1)})
 \\
@VVV @VVV @VVV \\
 L_{Z_1(V^{(3)}_{\cdot},P^{(3)}_{\cdot})/
 Z_2(V^{(3)}_{\cdot},P^{(3)}_{\cdot})}
@>>>
 L_{Z_1'(V^{(3)}_{\cdot},P^{(3)}_{\cdot})/
 Z_2(V^{(3)}_{\cdot},P^{(3)}_{\cdot})}
 @<<<
 L_{Z_1(V^{(1)}_{\cdot},P^{(1)}_{\cdot})/
 Z_2(V^{(1)}_{\cdot},P^{(1)}_{\cdot})} \\
 @VVV @VVV @VVV\\
 L_{U_1\times X/U_2\times X}
 @>{=}>>
 L_{U_1\times X/U_2\times X}
 @<{=}<<
 L_{U_1\times X/U_2\times X}
 \end{CD}
\]
Then, it is clear that the morphisms $b_i$ are quasi-isomorphic.
\hfill\qed

\vspace{.1in}

Now, we use the notation
$\gminik(\nbigf,\varphi)$ to denote
$\gminik(V_{\cdot},P_{\cdot},\phi,\phitilde)$,
in the following argument.
We also put as follows:
\[
 \Ob^H_{\rel}(\nbigf,\varphi):=
 Rp_{X\,\ast}\bigl(
\gminik(\nbigf,\varphi)\otimes\omega_X\bigr)
\]
Here, $\omega_X$ denotes the dualizing complex of $X$.
Let $\ob^H_{\rel}(\nbigf,\varphi)$
denote the composite of the morphisms
\[
\Ob^H_{\rel}(\nbigf,\varphi)
\lrarr
 Rp_{X\,\ast}\bigl(
 L_{U_1\times X/U_2\times X}\otimes\omega_X\bigr)
\lrarr
 L_{U_1/U_2}
\]

Let $(\nbigf_i,\varphi_i)$ $(i=1,2)$ be  as above.
Assume we are given a morphism
with a morphism $f:\nbigf_1\lrarr \nbigf_2$
such that $\varphi_2=f\circ \varphi_1$.

\begin{lem}
\label{lem;06.5.24.2}
We have the induced morphism
$\gminis(f):\gminik(\nbigf_2,\varphi_2)\lrarr
 \gminik(\nbigf_1,\varphi_1)$
such that
$\gminir(\nbigf_1,\varphi_1)\circ\gminis(f)
=\gminir(\nbigf_2,\varphi_2)$.
\end{lem}
\pf
Before going into the proof,
we give a rather canonical construction
of a datum for a given $(\nbigf,\varphi)$
as in the beginning of this subsubsection.
We take a sufficiently large integer $m_1$,
and we put 
$V_1:=p_X^{\ast}\bigl(
 p_{X\,\ast}\nbigf(m)\bigr)\otimes\nbigo(-m_1)
\oplus \nbigo_{U_2\times X}$.
The canonically defined morphism
$p_X^{\ast}\bigl(p_{X\,\ast}
 \nbigf(m_1)\bigr)\otimes\nbigo(-m_1)
\lrarr \nbigf$
and $\varphi:\nbigo_{U_2\times X}\lrarr \nbigf$
gives a surjection
$\epsilon:V_1\lrarr \nbigf$.
We also have the lift 
 $\phi:\nbigo_{U_2\times X}\lrarr V_1$.
Then, we take a sufficiently large integer $m_0$,
and we put
$V_0:=p_X^{\ast}\bigl(
 p_{X\,\ast}\Ker(\epsilon)(m_0)
 \bigr)\otimes\nbigo_X(-m_0)$.
Then, we have the surjection
$V_0\lrarr \Ker(\epsilon)$,
which induces $d_0:V_0\lrarr V_1$.
Since the kernel is locally free,
we put $V_{-1}:=\ker(d_0)$.
We take a resolution $P_{\cdot}$
such that $P_0$ is a direct sum of $\nbigo(-m_0)$.
Then we canonically obtain the morphism
$\phitilde_i:P_i\lrarr V_i$ $(i=0,-1)$ on $U_1\times X$.
Thus, we obtain a datum as above.

We take
$(V^{(i)}_{\cdot},P_{\cdot},\phi^{(i)},\phitilde^{(i)})$
for $(\nbigf_i,\varphi_i)$
by applying the construction explained above.
Then the morphism $f$ is canonically lifted
to $\widetilde{f}:V^{(1)}_{\cdot}\lrarr V^{(2)}_{\cdot}$
such that
$\phi^{(2)}=\widetilde{f}\circ \phi^{(1)}$
and $\phitilde^{(2)}_i=\widetilde{f}\circ \phitilde^{(1)}_i$.
Then, we obtain the following diagram:
\[
 \begin{CD}
 U_2\times X @>>> 
 Z_2(V^{(1)}_{\cdot},P_{\cdot})@>>>
 Z_2(V^{(2)}_{\cdot},P_{\cdot})\\
 @AAA @AAA @AAA \\
 U_1\times X @>>>
 Z_1(V^{(1)}_{\cdot},P_{\cdot})@>>>
 Z_1(V^{(2)}_{\cdot},P_{\cdot})
 \end{CD}
\]
Then, the claim of the lemma is clear.
\hfill\qed

\vspace{.1in}

Now, we assume $R^ip_{X\,\ast}\nbigf=0$ for $i>0$.
We put $\gbigv:=p_{X\,\ast}\nbigf$.
The morphism $\varphi:\nbigo_{U_2\times X}\lrarr \nbigf$
induces the section $\phibar$ of $\gbigv$.
We have the following diagram:
\begin{equation}
 \label{eq;06.5.24.3}
 \begin{CD}
 U_1 @>{j_2}>> U_2 \\
 @V{j_1}VV @V{i}VV \\
 U_2 @>{\phibar}>> \gbigv
 \end{CD}
\end{equation}
Here $i$ denotes the $0$-section.

\begin{prop}
\label{prop;06.5.24.5}
$\Ob_{\rel}^H(\nbigf,\varphi)$
is isomorphic to 
$j_2^{\ast}L_{U_2/\gbigv}
\simeq
 \gbigv^{\lor}[1]$,
and the morphism 
$\gminir(\nbigf,\varphi)$
is same as the morphism
$\kappa:j_2^{\ast}L_{U_2/\gbigv}\lrarr L_{U_1/U_2}$
induced from the diagram {\rm(\ref{eq;06.5.24.3})}.
\end{prop}
\pf
We have the naturally defined morphism
$a_1:p_X^{\ast}\gbigv\lrarr \nbigf$,
for which we have $\varphi=\phibar\circ a_1$.
Due to Lemma \ref{lem;06.5.24.2},
we have the following diagram:
\[
 \begin{CD}
 \gminik(\nbigf,\varphi)
@>{\gminis(a_1)}>>
 \gminik(p_X^{\ast}\gbigv,\phibar)
@>>>
 L_{U_1\times X/U_2\times X}
 \end{CD}
\]
It induces the following morphisms:
\[
 \begin{CD}
 \Ob^{H}_{\rel}(\nbigf,\varphi)
@>{b_0}>>
 \Ob^H_{\rel}(p_X^{\ast}\gbigv,p_X^{\ast}\phibar)
@>{b_1}>>
 L_{U_1/U_2}
 \end{CD}
\]

Let us see $b_1$ more closely.
In the construction for $(p_X^{\ast}\gbigv,p_X^{\ast}\phibar)$,
we can put $V_{1}=p_X^{\ast}\nbigf$, $V_i=U_2\times X$ $(i=0,-1)$,
$P_0=\nbigo_{U_1\times X}$ and $P_{-1}=U_2\times X$.
Then, $Z_1=U_2\times X$  and $Z_2=p_X^{\ast}\gbigv$.
The diagram (\ref{eq;06.5.24.1}) is given as follows:
\[
 \begin{CD}
 U_1\times X @>{j_1}>> Z_1 @=U_2\times X\\
 @V{j_2}VV @V{i}VV @VVV\\
 U_2\times X @>{\phibar}>> Z_2 @=\gbigv
 \end{CD}
\]
Here $i$ denotes the $0$-section.
We have 
$\gminik=j_1^{\ast}L_{Z_1/Z_2}\simeq p_X^{\ast}\gbigv_{|U_1\times X}[1]$,
and the morphism
$\gminir:\gminik\lrarr L_{U_1\times X/U_2\times X}$
is same as the pull back of $\kappa$.
In particular, we have the following factorization
of $b_1$:
\[
\begin{CD}
 \Ob^H(p_X^{\ast}\gbigv,p_X^{\ast}\phibar)
=\gbigv[1]\otimes Rp_{X\,\ast}\bigl(
p_{U_2}^{\ast}\omega_X\bigr)
@>{b_2}>>
 \gbigv[1]
@>{\kappa}>>
 L_{U_1/U_2}
\end{CD}
\]

It is easy to see that the composite
$b_2\circ a_0$ is isomorphic,
under the assumption $R^ip_{X\,\ast}\nbigf=0$ $(i>0)$.
Thus the proof of Proposition \ref{prop;06.5.24.5}
is finished.
\hfill\qed

\vspace{.1in}

We can obtain a similar result for a smooth projective curve.
The argument is similar and simpler,
and hence we omit to give a proof.

Let $D$ be a smooth projective curve.
Let $\nbigf$ be a $U_2$-coherent sheaf on $U_2\times D$
such that $p_{D\,\ast}\nbigf$ is locally free.
Let $\varphi$ be a morphism 
$\nbigo_{U_2\times D}\lrarr \nbigf$.
It induces the section $\nbigh(\varphi)$
of $p_{D\,\ast}\nbigf$.
Let $U_1$ be a subscheme of $U_2$
contained in $\nbigh(\varphi)^{-1}(0)$.

Assume that we have a locally free resolution
$(V_{0}\rarr V_1)$ of $\nbigf$
such that there exists a lift 
$\phi:\nbigo_{U_2\times D}\lrarr V_1$
of $\varphi$.
The morphism
$\phitilde:\nbigo_{U_1\times D}\lrarr V_{0\,|\,U_1\times D}$
is induced.
We put 
$\gminik(V_{\cdot},\phi):=
 \nhom(\nbigo_{U_1\times D},V_{\cdot|U_1\times D})^{\lor}$.

Let us construct a morphism
$\gminir(V_{\cdot},\phi):
 \gminik(V_{\cdot},\phi)\lrarr L_{U_1\times D/U_2\times D}$.
We put $Z_1:=N(\nbigo,V_0)$ and $Z_2:=N(\nbigo,V_1)$.
Then we have the naturally defined morphism
$Z_1\lrarr Z_2$.
The sections $\phi$ and $\phitilde$
induces the following commutative diagram:
\[
 \begin{CD}
 U_1\times D @>{j}>> Z_1 \\
 @VVV @VVV \\
 U_2\times D @>>> Z_2
 \end{CD}
\]
It induces the morphism
$\gminir(V_{\cdot},\phi):
 \gminik(V_{\cdot},\phi)
\simeq
 j^{\ast}L_{Z_1/Z_2}
\lrarr L_{U_1\times D/U_2\times D}$.
It can be shown that 
$\gminir(V_{\cdot},\phi)$
and $\gminik(V_{\cdot},\phi)$
depends only on $(\nbigf,\varphi)$.
Therefore, we use the notation
$\gminir(\nbigf,\varphi)$
and $\gminik(\nbigf,\varphi)$.
We put as follows:
\[
 \Ob^H_{\rel}(\nbigf,\varphi):=
 Rp_{D\,\ast}\bigl(
 \gminik(\nbigf,\varphi)\otimes\omega_D
 \bigr)
\]
Then, we have the induced morphism
$\ob_{\rel}^H(\nbigf,\varphi):
 \Ob^H_{\rel}(\nbigf,\varphi)
\lrarr L_{U_1/U_2}$.

We put $\gbigv:=p_{X\,\ast}\nbigf$.
We have the induced section $\phibar$.
Then, we obtain the diagram (\ref{eq;06.5.24.3}).
It induces the morphism
$\kappa:\gbigv^{\lor}[1]\lrarr L_{U_1/U_2}$.

\begin{prop}
\label{prop;06.5.25.200}
Assume $R^ip_{D\,\ast}\nbigf=0$ for $i>0$.
Then, we have the following commutative diagram:
\[
 \begin{CD}
\Ob^H_{\rel}(\nbigf,\varphi) 
@>{\ob_{\rel}(\nbigf,\varphi)}>> 
 L_{U_1/U_2}\\
 @V{\simeq}VV @V{=}VV\\
\gbigv^{\lor}[1]
 @>{\kappa}>> 
 L_{U_1/U_2}
 \end{CD}
\]
\end{prop}
\pf
It can be shown by an argument similar to 
the proof of Proposition \ref{prop;06.5.24.5}.
\hfill\qed

\subsection{Obstruction Theory}
\label{subsection;06.7.3.5}
\subsubsection{Definition and the foundational theorem of Behrend-Fantechi}

\label{subsubsection;06.6.8.3}

In the study of Gromov-Witten theory,
Li-Tian and Behrend-Fantechi introduced
the notion of virtual fundamental class
of moduli stacks with some good structure.
(See \cite{l-t} and \cite{bf}.)
In this paper,
we follow the work of Behrend-Fantechi.
Let us recall their notion of obstruction theory
of an algebraic stack
with a minor generalization.

\begin{df} 
\index{obstruction theory}
Let $\nbigx$ be an algebraic stack
over an algebraic stack $S$.
Let $E^{\cdot}$ be an object in $D(\nbigx)$
such that $\nbigh^i(E^{\cdot})$ are coherent $(i=-1,0,1)$.
A homomorphism
$\phi:E^{\cdot}\lrarr L_{\nbigx/S}$
is called an obstruction theory for $\nbigx/S$,
if $\nbigh^i(\phi)$ $(i\geq 0)$ are isomorphic
and $\nbigh^{-1}(\phi)$ is surjective.
In that case,
$E^{\cdot}$ is also called an obstruction theory
for $\nbigx/S$.
\hfill\qed
\end{df}

Since we have $\nbigh^i(L_{\nbigx})=0$ for $i>1$,
the condition implies $\nbigh^i(E^{\cdot})=0$ for $i>1$.
If $\nbigx$ is Deligne-Mumford,
we also have $\nbigh^1(E^{\cdot})=\nbigh^1(L_{\nbigx})=0$.

We will often use the following theorem.

\begin{prop}
[Behrend-Fantechi, Theorem 4.5, \cite{bf}]
 \label{prop;06.4.29.10}
Let $\nbigx$ be a Deligne-Mumford stack over $S$.
Let $\phi:E^{\cdot}\lrarr L_{\nbigx/S}$ be a morphism
in $D(\nbigx)$.
The following conditions are equivalent.
\begin{itemize}
\item
 $\phi$ is an obstruction theory.
\item
Let $T$ and $\overline{T}$ be $S$-schemes
such that $T$ is a closed subscheme of $\overline{T}$
whose ideal sheaf $J$ is square $0$.
Let $g:T\lrarr \nbigx$ be a morphism over $S$.
\begin{description}
\item[(A1)]
$g$ can be extended to a morphism
$\overline{g}:\overline{T}\lrarr \nbigx$
over $S$,
if and only if 
$\phi^{\ast}\bigl(o(g)\bigr)=0$
in $\Ext^1\bigl(g^{\ast}E^{\cdot},J\bigr)$,
where $o(g)$ is the obstruction class of $g$.
(See Proposition {\rm\ref{prop;06.4.29.5}}.)
\item[(A2)]
If $\phi^{\ast}\bigl(o(g)\bigr)=0$,
the set of the extension classes of $g$
is the torsor over the group
$\Ext^0\bigl(g^{\ast}E^{\cdot},J\bigr)$.
\hfill\qed
\end{description}
\end{itemize}
\end{prop}

We recall the perfectness of the obstruction theory
in the sense of Behrend-Fantechi
with a minor generalization.
\begin{df}
\label{df;06.6.13.12}
\index{perfect obstruction theory}
Let $\phi:E^{\cdot}\lrarr L_{\nbigx/S}$
be an obstruction theory 
of an algebraic stack $\nbigx$ over $S$.
It is called perfect,
if it is quasi isomorphic to a complex
of locally free sheaves $F^{-1}\rarr F^{0}\rarr F^1$
in the derived category $D(\nbigx)$.
\hfill\qed
\end{df}

In that case, the number $-\rank F^{1}+\rank F^0-\rank F^{-1}$
is well defined on each connected component of $\nbigx$.
The number is called the expected dimension of $\nbigx$
over $S$ with respect to $\phi$.

If $\nbigx$ is Deligne-Mumford,
we have $\nbigh^1(E^{\cdot})=\nbigh^1(L_{\nbigx})=0$
for the obstruction theory $E^{\cdot}$.
Hence, a perfect obstruction theory is
quasi isomorphic to $F^{-1}\rarr F^{0}$.
The important and foundational theorem 
of Behrend and Fantechi is the following.
(See also \cite{l-t}.)
Let $A_{\ast}(\nbigx)$ denote the Chow group of $\nbigx$
with rational coefficient.

\begin{prop}
 [Behrend-Fantechi, Section 5, \cite{bf}]
 \label{prop;06.6.13.20}
 Let $\nbigx$ be a Deligne-Mumford stack
 over a smooth scheme $S$.
 A perfect obstruction theory
 $\phi:E^{\cdot}\lrarr L_{\nbigx/S}$
 induces the virtual fundamental class
 $[\nbigx,\phi]\in A_d(\nbigx)$,
 where $d$ is the expected dimension
 with respect to $\phi$.

When $\nbigx$ is smooth,
$[\nbigx,\phi]$ is given by the Euler class
of the vector bundle $\nbigh^1\bigl(E^{\cdot\,\lor}\bigr)$.
\hfill\qed 
\end{prop}

Let $\nbigx_i$ $(i=1,2)$ be algebraic stacks over $S$
with obstruction theories
$\phi_{i}:E^{\cdot}_{i}\lrarr L_{\nbigx_i/S}$.
Assume we have the following commutative diagram:
\[
 \begin{CD}
 \nbigx_1 @>{f}>> \nbigx_2\\
 @V{g}VV @VVV \\
 \nbigy_1@>{h}>>\nbigy_2 @>>>S
 \end{CD}
\]
Recall the following definition in \cite{bf}.
\begin{df}
\label{df;06.5.23.110}
We say that $\phi_i$ are compatible over $h$,
if we have the following morphism
of distinguished triangles on $\nbigx_1$:
\[
 \begin{CD}
 f^{\ast}E^{\cdot}_2@>>>
 E^{\cdot}_1@>>>
 g^{\ast}L_{\nbigy_1/\nbigy_2}@>>>
 f^{\ast}E^{\cdot}_2[1]\\
 @VVV @VVV @VVV @VVV \\
 g^{\ast}L_{\nbigx_2/S}@>>>
 L_{\nbigx_1/S} @>>>
 L_{\nbigx_1/\nbigx_2} @>>>
 g^{\ast}L_{\nbigx_2/S}[1]
 \end{CD}
\]
\hfill\qed
\end{df}

We recall the following theorem for later use.
\begin{prop}
[Behrend-Fantechi, Proposition 7.5, \cite{bf}]
\label{prop;06.6.11.150}
Assume $\nbigx_i$ $(i=1,2)$ are Deligne-Mumford,
and the obstruction theories $\phi_i$ are perfect.
If $\phi_i$ are compatible over $h$,
then $h^{!}[\nbigx_2,\phi_2]=[\nbigx_1,\phi_1]$,
at least if $h$ is smooth
or $\nbigy_i$ are smooth over $S$
\hfill\qed
\end{prop}

See \cite{bf} for more detail about virtual fundamental classes.

\subsubsection{Easy example}
\label{subsubsection;06.6.22.1}

Let $X$ be a smooth variety over $k$.
We would like to construct an obstruction theory
of the moduli spaces $\nbigm$ of some objects on $X$.
Our naive strategy 
is summarized as follows
(See \cite{bf}, for example):
\begin{enumerate}
\item
Take the classifying stack $Y$ of such objects
over $X$.
It means that such objects over $U\times X$
bijectively correspond 
to morphisms $\Phi:U\times X\lrarr Y$ over $X$.
For example, recall that a vector bundle of rank $R$
over $U\times X$ corresponds to a map
$U\times X\lrarr X_{\GL(R)}$ over $X$.
\item
For any classifying map $\Phi:U\times X\lrarr Y$,
we obtain the morphism
$\Phi^{\ast}L_{Y/X}\lrarr L_{U\times X/X}$ on $U\times X$.
Let $\omega_X$ denote the dualizing complex on $X$,
i.e., it is the canonical sheaf shifted by the dimension of $X$.
Then, we obtain the morphisms on $U$:
\[
 \Ob_{U}:=Rp_{X\ast}\bigl(\Phi^{\ast}L_{Y/X}\otimes\omega_X\bigr)
 \lrarr Rp_{X\ast}\bigl(p^{\ast}_XL_{U/k}\otimes\omega_X\bigr)
 \lrarr L_{U/k}.
\]
In particular, we obtain the morphism
$\Ob_{\nbigm}\lrarr L_{\nbigm}$ on $\nbigm$.
\item
Then we hope that the morphism
$\Ob_{\nbigm}\lrarr L_{\nbigm}$ gives the obstruction theory,
in some cases.
Note that the property is local,
once the morphism is given globally.
Thus we have only to check the claim
for the sufficiently small etale open sets of $\nbigm$.
The tool for checking is Proposition \ref{prop;06.4.29.10}.
\end{enumerate}

\begin{rem}
In general,
we need some modification for the construction of $\Ob_{\nbigm}$
to obtain the good obstruction theory.
\hfill\qed
\end{rem}

Let us see the easiest example.
Let $F$ and $V$ be vector bundles defined on $X$.
Let $U$ be any scheme over $k$,
and let $f:p_U^{\ast}(F)\lrarr p_U^{\ast}(V)$ be a morphism
of $\nbigo_{U\times X}$-modules over  $U\times X$.
It is easy to see that such a morphism $f$ corresponds to
a morphism $\Phi_f:U\times X\lrarr N(F,V)$ over $X$.
Thus we obtain the complex
$\gminig(f):=\Phi_f^{\ast}L_{N(F,V)/X}$
and the morphism $\gminig(f)\lrarr L_{U\times X/X}$
in the derived category $D(U\times X)$.

\begin{lem}
The complex $\gminig(f)$ is quasi isomorphic to
$p_U^{\ast}\nhom(V,F)$.
\end{lem}
\pf
Let $\pi:N(F,V)\lrarr X$  denote the natural projection.
Since the morphism $N(F,V)\lrarr X$ is smooth,
the cotangent complex $L_{N(F,V)/X}$ is quasi isomorphic
to $\Omega_{N(F,V)/X}\simeq \pi^{\ast}\nhom(V,F)$.
Thus we obtain the quasi isomorphism
$\Phi_f^{\ast}L_{N(F,V)/X}\simeq p_U^{\ast}\nhom(V,F)$.
\hfill\qed

\vspace{.1in}

We put $\Ob(f):=Rp_{X\ast}\bigl(\gminig(f)\otimes\omega_X\bigr)$.
Then, we obtain the morphisms
$\Ob(f)\lrarr 
 Rp_{X\ast}\bigl(L_{U\times X/X}\otimes\omega_X\bigr)
 \lrarr L_U$ 
in the derived category  $D(U)$.
The composite is denoted by $\ob(f)$.

Now, let $M(F,V)$ denote the moduli scheme of
the morphisms $F\lrarr V$, i.e.,
maps $U\lrarr M(F,V)$
correspond to 
$f:p_U^{\ast}(F)\lrarr p_U^{\ast}(V)$ on $U\times X$.
It is easy to see that
$M(F,V)$ is isomorphic to the vector space
$H^0\bigl(X,\nhom(F,V)\bigr)$.
We have the universal morphism
$f^{u}:p_{M(F,V)}^{\ast}(F)\lrarr p_{M(F,V)}^{\ast}(V)$
over $M(F,V)\times X$.
It induces the morphism
$\ob(f^u):\Ob(f^u)\lrarr L_{M(F,V)}$.

\begin{lem}
\label{lem;06.5.3.1}
The morphism $ob(f^u)$ gives the obstruction theory
of $M(F,V)$.
\end{lem}
\pf
It is almost obvious
from the universal properties of
$N(F,V)$ and $M(F,V)$.
But, we give an argument for the explanation
of our later discussion.
We have only to check the conditions
(A1) and (A2) in Proposition \ref{prop;06.4.29.10}.

Since the claim is local,
we can check the claim for any sufficiently small
open subset $U$ of $M(F,V)$.
Let $T$ be an affine scheme over $k$.
A morphism $g:T\lrarr U$
induces the morphism
$g_X:T\times X\lrarr U\times X$
and $\gtilde_X=\Phi\circ g_X:T\times X\lrarr N(F,V)$
over $X$.
Let $\Tbar$ denote a scheme such that
$T$ is embedded in $\Tbar$
whose ideal  $J$ is square zero.
The deformation theory of the morphisms
$g$ and $\gtilde_X$ is controlled by the groups 
$\Ext^i\bigl(g^{\ast}L_{U/k},J\bigr)$
and $\Ext^i\bigl(\gtilde_X^{\ast}L_{N(F,V)/X},J_X\bigr)$
respectively.
We have the following commutative diagram:
\[
 \begin{CD}
 \Ext^i\bigl(g^{\ast}L_{U/k},J\bigr) 
 @>{h}>>
 \Ext^i\bigl(g^{\ast}\Ob(f^u),J\bigr)  \\
 @VVV @A{\simeq}AA \\
 \Ext^i\bigl(g_X^{\ast}L_{U\times X/X},J_X\bigr) 
 @>>>
 \Ext^i\bigl(g_X^{\ast}(\gminig),J\bigr)
 @=
 \Ext^i\bigl(\gtilde_X^{\ast}L_{N(F,V)/X},J\bigr)
 \end{CD}
\]
We have the obstruction classes
$o(g)\in \Ext^1\bigl(g^{\ast}L_{U/k},J\bigr)$
and $o(\gtilde_X)\in \Ext^1\bigl(g_X^{\ast}\gminig,J\bigr)$
of the morphisms $g$ and $\gtilde_X$
respectively.
By the functoriality of the cotangent complex,
the obstruction class $o(g)$
is mapped to the obstruction class $o(\gtilde_X)$
in the diagram above.

If the image $h\bigl(o(g)\bigr)$ is $0$,
the class $o(\gtilde_X)$ is $0$.
Thus $\gtilde_X$ can be extended to
a morphism $\Tbar\times X\lrarr N(F,V)$,
which induces a morphism
of $p^{\ast}_{\Tbar}(F)\lrarr p^{\ast}_{\Tbar}(V)$
on $\Tbar\times X$.
By the universal property of $M(F,V)$,
we obtain a morphism
$\Tbar\lrarr M(F,V)$,
which is the extension of $g$.
Therefore, the condition (A1) is satisfied.

Similarly, we know that the morphism
$\Ext^0\bigl( g^{\ast}L_{U/k},J\bigr)
 \lrarr 
 \Ext^0\bigl(\gtilde_X^{\ast}L_{N(F,V)/X},J\bigr)$
is isomorphic
from the universality of $M(F,V)$ and $N(F,V)$.
Hence the condition (A2) is also satisfied.
Thus we are done.
\hfill\qed

\subsubsection{Obstruction theory for locally free subsheaves}
\label{subsubsection;06.5.4.5}

Let $X$ be a smooth projective variety over $k$
with an ample line bundle $\nbigo_X(1)$.
Let $V$ be a locally free sheaf on $X$.
Let $W$ denote an $R$-dimensional $k$-vector space.
We denote $W\otimes\nbigo_X$ by $W_X$.
We have the natural right $\GL(W)$-action on $N(W_X,V)$.
The quotient stack is denoted by $\ywquo$.

We consider the deformation theory of
locally free subsheaves $F\subset V$ of rank $R$.
Let $U$ be any $k$-scheme.
Any locally free subsheaf $f:F\lrarr p_U^{\ast}V$ on $U\times X$
induces the morphism
$\Phi(F,f):U\times X\lrarr \ywquo$ over  $X$.
We put $\gminig(F,f):=\Phi(F,f)^{\ast}L_{\ywquo/X}$,
and $\Ob(F,f):=Rp_{X\ast}\bigl(\gminig(F,f)\otimes\omega_X\bigr)$.
Then, we have the morphism
$\gminig(F,f)\lrarr L_{U\times X/X}$ on $U\times X$,
which induces $\ob(F,f):\Ob(F,f)\lrarr L_U$ on $U$.
The following lemma can be shown by 
the argument explained in the subsection 
\ref{subsubsection;06.4.29.15}.

\begin{lem}
$\gminig(F,f)$ is represented by
$\Cone(\alpha)[-1]$ of the morphism
$\alpha:\nhom(p_U^{\ast}V,F)
\lrarr \nhom(F,F)$,
where $\alpha$ is given by $\alpha(a)=a\circ f$.
\hfill\qed
\end{lem}

\begin{rem}
We put $V_{-1}:=F$ and $V_0:=p_U^{\ast}V$,
and we regard $V_{\cdot}=(V_{-1}\rarr V_0)$ as a complex,
where $V_0$ stands in the degree $0$.
Then, $\Cone(\alpha)$ is naturally isomorphic
to $\nhom(V_{-1}[1],V_{\cdot})^{\lor}[-1]$.
\hfill\qed
\end{rem}

Let $H$ be a polynomial.
We have the moduli of quotients 
$(q:V\lrarr Q)$ of $V$
such that the Hilbert polynomials of $Q$
are $H$.
Let $M(V,H)$ denote the open subscheme
which consists of the points
$(q:V\lrarr Q)$ such that $\Ker(q)$ are locally free.
Then, we have the universal family
$f^u:F^u\lrarr p_{M(V,H)}^{\ast}(V)$ defined over
$M(V,H)\times X$.
We obtain the morphism 
$\ob(F^u,f^u):\Ob(F^u,f^u)\lrarr L_{M(V,H)}$.

\begin{prop} 
 \label{prop;06.5.2.1}
The morphism $\ob(F^u,f^u)$
gives the obstruction theory of $M(V,H)$.
\end{prop}
\pf
Let $N$ be a sufficiently large number satisfying the condition $O_N$
for the family $F^u$,
i.e., we have $H^i\bigl(X,F^u(N)_{|\{q\}\times X}\bigr)=0$
for any $q\in M(V,H)$ and $i>0$,
and $F^u_{|\{q\}\times X}(N)$ are globally generating for any $q\in M(V,H)$.
We put
$\Ftilde^u:=p_X^{\ast}p_{X\ast}\bigl(F^u(N)\bigr)
 \otimes\nbigo(-N)$.
We have the natural surjection
$g:\Ftilde^u\lrarr F^u$.

We put $\Fbar=\nbigo(-N)^{\oplus\,d}$,
where $d=\rank \Ftilde^u$.
We have the Grassmaniann bundle 
$\pi:Gr(\Fbar,R)\lrarr X$
associated to the vector bundle $\Fbar$,
i.e., the fiber of $\pi$ over a point $x\in X$ is
the Grassmann variety of the $R$-dimensional quotient spaces
of the vector space $\Fbar_{|x}$.
We denote the universal quotient bundle 
over  $Gr(\Fbar,R)$ by $Q$.
Then, we have the vector bundle
$\Ytilde_{\quo}:=N(Q,\pi^{\ast}V)$ over $Gr(\Fbar,R)$,
which is a variety smooth over $X$.
We have the natural morphism $\pi_1:\Ytilde_{\quo}\lrarr \ywquo$.

We would like to check the conditions (A1) and (A2)
in Proposition \ref{prop;06.4.29.10}.
Let $U$ be any sufficiently small open set of $M(V,H)$,
on which we can assume that there exists
an isomorphism $\Ftilde^u\simeq p_U^{\ast}\Fbar$.
Thus, the morphism
$\alpha:p_U^{\ast}\Fbar\lrarr F^u$
is given on $U\times X$.
From $\alpha$ and $f^u$,
we obtain the morphism
$\Phi(\alpha,F^u,f^u):U\times X\lrarr \Ytilde_{\quo}$
over $X$.
By the argument in the subsubsection
\ref{subsubsection;06.4.29.15},
we can show that the complex
$\Phi(\alpha,F^u,f^u)^{\ast}L_{\Ytilde_{\quo}/X}$
is quasi isomorphic to $\Cone(\beta)[-1]$
for the morphism
$\beta:\nhom(p_U^{\ast}V,F^u)\oplus 
 \nhom(F^u,p_U^{\ast}\Fbar)\lrarr\nhom(F^u,F^u)$,
where $\beta(b_1,b_2)=b_1\circ f^u-f^u\circ b_2$.
We can also show that the natural morphism
$\Cone(\alpha)[-1]\lrarr \Cone(\beta)[-1]\lrarr L_{U\times X/X}$
corresponds to the factorization
$\Phi(F^u,f^u)^{\ast}L_{Y_{\quo}/X}
\lrarr \Phi(\alpha,F^u,f^u)^{\ast}L_{\Ytilde_{\quo}/X}
\lrarr L_{U\times X/X}$ associated to
$U\times X\lrarr \Ytilde_{\quo}\lrarr Y_{\quo}(W_{\cdot})$.
We put as follows:
\[
 \gminig(\alpha,F^u,f^u)
 :=\Phi(\alpha,F^u,f^u)^{\ast}L_{\Ytilde_{\quo}/X},
\quad
\Ob(\alpha,F^u,f^u):=
 Rp_{X\ast}\bigl(\gminig(\alpha,F^u,f^u)\otimes\omega_X\bigr).
\]

Let $T$ be an affine scheme,
and let $g:T\lrarr U$ be a morphism.
We put $\widetilde{g}_X:=\Phi(F^u,f^u)\circ g$
and $\widehat{g}_X:=\Phi(\alpha,F^u,f^u)\circ g$.
For any coherent sheaf $J$ on $T$,
we have the following diagram:
\begin{equation}
 \label{eq;06.4.29.20}
 \begin{CD}
 \Ext^i(g^{\ast}L_{U/k},J) @>{h^i_1}>>
 \Ext^i\bigl(g^{\ast}\Ob(\alpha,F^u,f^u),J\bigr)
 @>{h^i_2}>>
 \Ext^i(g^{\ast}\Ob(F^u,f^u),J) \\
 @VVV @A{\simeq}AA @A{\simeq}AA \\
 \Ext^i(g_X^{\ast}L_{U\times X/X},J_X)
 @>>>
 \Ext^i\bigl(g^{\ast}_X\gminig(\alpha,F^u,f^u),J\bigr)
 @>>>
 \Ext^i\bigl(g^{\ast}_X\gminig(F^u,f^u),J\bigr)
 \end{CD}
\end{equation}
Let $\Tbar$ be an affine scheme into which
$T$ is embedded closedly
such that the corresponding ideal $J$ is square zero.
Due to the deformation theory of Illusie,
we have the obstruction classes
of the morphisms $g$ and $\widehat{g}_X$
in the groups $\Ext^1(g^{\ast}L_{U/k},J)$
and $\Ext^1\bigl(g^{\ast}_X\gminig(\alpha,F^u,f^u),J \bigr)$
respectively.
The classes are denoted by $o(g)$ and $o(\widehat{g}_X)$.
Due to the functoriality,
the class $o(g)$ is mapped to the class $o(\widehat{g}_X)$
in the diagram (\ref{eq;06.4.29.20}).
If $h^1_1(o(g))$ is $0$,
then the morphism $\widehat{g}_X$
can be extended.

Note that the cohomology sheaves
$R^ip_{X\ast}\bigl(\nhom(F^u,p_U^{\ast}\Fbar)\otimes\omega_X \bigr)$
vanish unless $i=0$,
because of our choice of $N$.
Thus, we have the isomorphism
$\Ext^i\bigl(g^{\ast}\Ob(\alpha,F^u,f^u),J\bigr)
\simeq 
 \Ext^i\bigl(g^{\ast}\Ob(F^u,f^u),J\bigr)$
for any $i>0$ and for any coherent sheaf $J$ on $T$. 
Hence $h^1_2\circ h^1_1(o(g))=0$
implies $h^1_1\bigl(o(g)\bigr)=0$.
Then the morphism $\widehat{g}_X$
can be extended over $\Tbar\times X$,
and hence $\widetilde{g}_X$ can also be extended
over $\Tbar\times X$.
Therefore, we obtain a locally free subsheaf
$\widehat{F}$ of $p_{\Tbar}^{\ast}(V)$ 
on $\Tbar\times X$,
which is the extension of $g_X^{\ast}F^u$.
Due to the universal property of $M(V,H)$,
the morphism $g$ can be extended over $\Tbar$.
Therefore, the condition (A1) is satisfied.

Let us check the condition (A2).
We put as follows:
\[
 H_0:=
 \Ext^0\Bigl(p_{X\ast}\bigl(
 g_X^{\ast}\nhom(F^u,p_U^{\ast}\Fbar)\otimes\omega_X
 \bigr),J
 \Bigr)
=H^0\Bigl(T,
 g^{\ast}\nEnd\bigl(
 p_{X\,\ast}\bigl(
 F^u(m)\bigr)
 \bigr)\otimes J
 \Bigr)
\]
\[
  H_1:= \Ext^0\bigl(g^{\ast}\Ob(g,F^u,f^u),J\bigr)
=\Ext^0\bigl(
 \widehat{g}_X^{\ast}L_{\Ytilde_{\quo}/X},J_X
 \bigr),
\quad
 H_2:=\Ext^0\bigl(g^{\ast}\Ob(F^u,f^u),J\bigr)
\]
We obtain the exact sequence
$0\lrarr H_0\lrarr H_1\lrarr H_2\lrarr 0$.
Due to the theory of Illusie, 
$H_1$ parameterizes the set of extensions
$\widehat{g}'_X:\Tbar\times X\lrarr \Ytilde_{\quo}$ 
of $\widehat{g}_X$.
The natural action of $H_0$ on $H_1$
determines the equivalence relation on $H_1$,
and it is easy to see that
$\widehat{g}'_X\sim \widehat{g}''_X$
if and only if 
$\pi_1\circ\widehat{g}_X'=\pi_1\circ\widehat{g}''_X$,
because $H_0$ parameterizes
the deformation of the morphisms
$\Fbar\lrarr F^u$.
Thus the set of the extensions
of the morphisms 
$T\times X\lrarr \ywquo$ over $\Tbar\times X$
is the torsor over the group $H_2$.

By the universal property of $M(V,H)$ and $\ywquo$,
the set of the extensions of 
$g$ over $\Tbar$ is also the torsor over $H_2$.
Namely the condition (A2) is satisfied.
Therefore we are done.
\hfill\qed

\vspace{.1in}
Usually, we consider the deformation theory
of quotients of $V$.
Let $H$ be a polynomial, and 
let $\Quot(V,H)$ denote the quot scheme
which parameterizes the quotient sheaves of $V$
whose Hilbert polynomials are $H$.
We have the universal quotient
$q^u:p_{\Quot(V,H)}^{\ast}(V)\lrarr Q^u$
on $\Quot(V,H)\times X$.
We denote the kernel of $q^u$
by $F^u$,
and the inclusion $F^u\lrarr p_{\Quot(V,H)}^{\ast}(V)$
is denoted by $f^u$.

Let us consider the case $\dim X=1$.
Let $H_V$ denote the Hilbert polynomial of $V$.
Then $\Quot(V,H)$ parameterizes
the locally free subsheaves of $V$
whose Hilbert polynomials is $H_V-H$.
Therefore,
we have obtained the obstruction theory
$\ob(f^u):\Ob(f^u)\lrarr L_{\Quot(V,H)}$.

\begin{prop}
In the case $\dim(X)=1$,
the obstruction theory
$\ob(f^u)$ is perfect.
The scheme $\Quot(V,H)$ is smooth,
if $H$ is a constant,
i.e., $H$ is a Hilbert polynomial of
sheaves of finite length.
\end{prop}
\pf
To show the perfectness of $\Ob(f^u)$,
we have only to show that
$Rp_{X\,\ast}\bigl(\gminig(F^u,f^u)^{\lor}\bigr)$
is perfect of amplitude contained in $[0,1]$.
Let $q$ be any point of $Q(V,H)$.
We put $F:=F^u_{|\{q\}\times X}$ and $Q:=V/F$.
The complex $\gminig(F^u,f^u)^{\lor}_{|\{q\}\times X}$
is $\Cone(\gamma)[-1]$
for the natural morphism $\gamma:\nHom(F,F)\lrarr \nHom(F,V)$,
which is quasi-isomorphic to $\nhom(F,Q)$.
Hence we have
$H^i\bigl(X,\gminig(F^u,f^u)^{\lor}_{|\{q\}\times X} \bigr)=0$
unless $i=0,1$ for any point $q\in \Quot(V,H)$.
Then the desired perfectness easily follows.

Let us show the second claim.
When $H$ is a constant,
i.e., $Q$ is a torsion sheaf,
we always have the vanishing
$H^1\bigl(X,
 \gminig(F^u,f^u)^{\lor}_{|\{q\}\times X}\bigr)=0$.
Let $T$ be any affine scheme over $k$,
and let $g:T\lrarr \Quot(V,H)$ be a morphism.
Then we obtain the vanishing
$\Ext^{1}\bigl( g^{\ast}\Ob(f^u),J\bigr)=0$
for any coherent $\nbigo_T$-module $J$,
and hence the vanishing
of any obstruction class.
Thus we obtain the smoothness.
\hfill\qed

\vspace{.1in}

\begin{rem}
\label{rem;06.5.3.30}
Let us consider the case $\dim X=2$.
Let $Q^{tf}(V,H)$ denote the open subset
of $Q(F,H)$ corresponding to the torsion-free quotient.
It gives an open subset of
the moduli stack of locally free subsheaves of $V$.
Then, $F^u_{|\{q\}\times X}$ are locally free
for any $q\in Q^{tf}(V,H)$.
Therefore,
we obtain the obstruction theory
$\ob(F^u,f^u):\Ob(F^u,f^u)\lrarr 
 L_{Q^{tf}(V,H)/k}$
from Proposition {\rm\ref{prop;06.5.2.1}}.
\hfill\qed
\end{rem}

\subsubsection{Obstruction theory for filtrations of a vector bundle
 on a curve}
\label{subsubsection;06.6.24.5}

Let $S$ be a scheme over $k$,
and let $\nbigd$ be a smooth projective 
curve over $S$ 
provided with an ample line bundle $\nbigo(1)$.
The projection $\nbigd\lrarr S$ is denoted by $p$.
For any point $s\in S$,
the fiber over $s$ is denoted by $\nbigd_s$.
Let $V$ and $F$ be a locally free sheaf
on $\nbigd$ provided with an
injective morphism $f:F\lrarr V$.
Assume that the quotient is $S$-flat.

Let $H_i$ be polynomials.
For an $S$-scheme $T$,
let $F(T)$ denote the set of 
the data $(g,V^{\ast})$ as follows:
\begin{itemize}
\item
 $g$ is a morphism $T\lrarr \nbigd$ over $S$.
\item
 $V^{\ast}$ denotes a filtration
 $g^{\ast}V=
 V^{(1)}\supset
 V^{(2)}\supset\cdots \supset
 V^{(l)}\supset V^{(l+1)}=g^{\ast}F$.
 We assume that the quotients
 $\Cok_{i}:=V^{(1)}/V^{(i+1)}$ are $T$-flat.
\item
 The Hilbert polynomials
 of $\Cok_{i|\,\nbigd_s}$ are $H_i$
 for any $i=1,\ldots,l$ and $s\in S$.
\end{itemize}
The functor is representable by a scheme,
which can be shown by the standard technique
using the quot schemes.
Let $M(H_{\ast})$ denote the moduli scheme.
Let $p_{M(H_{\ast})}$ denote the projection
$M(H_{\ast})\times_S \nbigd\lrarr\nbigd$.
We have the universal filtration
on $M(H_{\ast})\times_S \nbigd$:
\[
 p_{M(H_{\ast})}^{\ast}V
=\nbigv^{(1)}
\supset
 \nbigv^{(2)}\supset\cdots\supset
 \nbigv^{(l)}\supset\nbigv^{(l+1)}
=p_{M(H_{\ast})}^{\ast}F.
\]

To discuss an obstruction theory of $M(H_{\ast})$,
we introduce some stacks.
Take vector spaces $W_i$ $(i=2,\ldots,l)$
over $k$
such that $\rank W_i=\rank \nbigv^{(i)}=:r_i$.
We put $W^{(i)}:=W_i\otimes\nbigo_{\nbigd}$ $(i=2,\ldots,l)$.
We put $W^{(1)}=V$ and $W^{(l+1)}=F$.
We put 
$Y_0:=N(W^{(l+1)},W^{(1)})$
and 
$R_1:=\prod_{i=1}^l N(W^{(i+1)},W^{(i)})$.
We put $G(W_{\ast}):=\prod_{i=2}^{l}\GL(W_i)$.
We have the natural right
$G(W_{\ast})$-action on $R_1$.
Let $Y_1$ denote the quotient stack
of $R_1$ by the $G(W_{\ast})$-action.
By the composition of the maps,
we obtain the morphism
$\phi:R_1\lrarr Y_0$,
which induces $Y_1\lrarr Y_0$.
We also put $Y_2:=\nbigd$.
Then the morphism $F\lrarr V$ induces the morphism
$Y_2\lrarr Y_0$.
We put $Y:=Y_1\times_{Y_0}Y_2$.

Let $V^{\ast}$ denote
the above filtered vector bundle
on $T\times_S \nbigd$.
We have the naturally defined morphism:
\[
G(V^{\ast}):
 \Phi_0(V^{\ast})^{\ast} L_{Y_0/\nbigd}
\lrarr
 \bigoplus_{i=1,2}
 \Phi_i(V^{\ast})^{\ast} L_{Y_i/\nbigd}
\]
We use the notation in the subsubsection
\ref{subsubsection;06.5.23.10}.
We put
$\gminig(V^{\ast}):=
 C_2(V^{\ast},V^{\ast})^{\lor}[-1]$.
We obtain the morphism
$\Phi_i(V^{\ast}):T\times_S \nbigd\lrarr Y_i$.
By the argument in the subsubsection
\ref{subsubsection;06.4.29.15},
the cone of $G(V^{\ast})$ is expressed by the complex 
$\gminig(V^{\ast})$.
Thus, we have the naturally defined morphism
$\gminig(V^{\ast})\lrarr L_{T\times_S \nbigd/\nbigd}$.
We put $\Ob(V^{\ast}):=
 Rp_{\ast}\bigl(\gminig(V^{\ast})\otimes\omega_{\nbigd/S}\bigr)$,
and then we obtain the morphism
$\ob(V^{\ast}):\Ob(V^{\ast})\lrarr L_T$.

\begin{lem}
 \label{lem;06.5.3.15}
The morphism $\ob(\nbigv^{\ast})$ gives
an obstruction theory for 
$M(H_{\ast})$.
\end{lem}
\pf
Let us take locally free sheaves $J^{(i)}$
$(i=2,\ldots,l)$ on $\nbigd$ such that
$H^1\bigl(\nbigd_s, \nhom(J^{(i)},\nbigv^{(i)}_{|\nbigd_s})
\bigr)=0$ for any $s\in \nbigd$.
For any $S$-scheme $T$,
let $\Ftilde(T)$ denote the set of
the data $(g,V^{\ast},\varphi_{\ast})$
as follows:
\begin{itemize}
\item
 $g$ denotes a morphism $T\lrarr \nbigd$,
and $V^{\ast}$ denotes a filtration
 as above.
\item
 $\varphi_{\ast}$  denotes a tuple
 of surjections of $g^{\ast}J^{(i)}$
 onto $V^{(i)}$.
\end{itemize}
The functor $\Ftilde$ is representable
by the scheme which is denoted by
$\widetilde{M}(H_{\ast})$.
It is easily described.
We have the locally free sheaf
$N_i=\nhom\bigl(p_{M(H_{\ast})}^{\ast}J^{(i)},\nbigv^{(i)}\bigr)$
on $M(H_{\ast})\times_{S}\nbigd$.
Then $\widetilde{M}(H_{\ast})$
is isomorphic to an open subset of
$\bigoplus p_{\ast}N_i$.
On $\widetilde{M}(H_{\ast})\times_{S}\nbigd$,
we have the universal filtration
$\nbigv^{\ast}$
with the tuple of surjective morphisms
$\varphi^u_{\ast}$.

Let $Gr(J^{(i)},r_i)$ be the Grassmannian bundles
of $r_i$-dimensional quotient spaces
associated to the vector bundles $J^{(i)}$.
We have the universal quotient bundle $Q_i$.
We put $Z:=\prod_{i=2}^{l} Gr(J^{(i)},r_i)$,
where the fiber product is taken over $\nbigd$.
The pull back of $Q_i$
via the projection
$Z\lrarr Gr(J^{(i)},r_i)$ is denoted by
$\widetilde{W}^{(i)}$ $(i=2,\ldots,l)$.
The pull back of $V$ and $F$
via the projection
$Z\lrarr \nbigd$ are denoted by
$\widetilde{W}^{(1)}$
and $\widetilde{W}^{(l+1)}$
respectively.
Then we put
$\Ytilde_0:=
 N\bigl(\Wtilde^{(l+1)},\Wtilde^{(1)}\bigr)$,
$\Ytilde_1:=
 \prod_{i=1}^l N(\Wtilde^{(i+1)},\Wtilde^{(i)})$
and $\Ytilde_2:=Z$.
We have the natural morphisms
$\Ytilde_i\lrarr\Ytilde_0$ $(i=1,2)$ as above,
the fiber product
$\Ytilde_1\times_{\Ytilde_0}\Ytilde_2$
is denoted by $\Ytilde$.
The inclusions $\Ytilde\lrarr \Ytilde_i$
are denoted by $j_i$.
On $\Ytilde$, 
we have the natural morphism
$ j_0^{\ast}L_{\Ytilde_0/\nbigd}\lrarr
\bigoplus_{i=1,2}
  j_i^{\ast}L_{\Ytilde_i/\nbigd}$.
The cone of the morphism is denoted
by $\Ob(\Ytilde)$.
Then we have the naturally defined morphism
$\ob(\Ytilde):\Ob(\Ytilde)\lrarr L_{\Ytilde/\nbigd}$,
and it gives an obstruction theory
for $\Ytilde$ over $\nbigd$.
(Basic example in \cite{bf}).

Let $T$ be an $S$-scheme.
From $(V^{\ast},\varphi_{\ast})$,
we obtain the morphism
$\Phi_i(V^{\ast},\varphi_{\ast}):
 T\times_S \nbigd\lrarr \Ytilde_i$.
Therefore, we obtain
$\Phi(V^{\ast},\varphi_{\ast})^{\ast}
 \Ob(\Ytilde)\lrarr L_{T\times_S \nbigd/S}$.
We put 
$\Obtilde(V^{\ast},\varphi_{\ast}):=Rp_{\ast}\bigl(
 \Ob(\Ytilde)\otimes\omega_{\nbigd/S}\bigr)$,
and then we obtain the morphism
$\obtilde(V^{\ast},\varphi_{\ast}):
 \Obtilde(V^{\ast},\varphi_{\ast})
\lrarr L_{T/S}$.

Let us describe the complex
$\Obtilde(V^{\ast},\varphi_{\ast})$.
We have the morphism
$\nhom\bigl(V^{(i)},J^{(i)}\bigr)
\lrarr
 \nhom\bigl(V^{(i)},V^{(i)}\bigr)$
given by
$a_i\longmapsto \varphi_i\circ a_i$.
It induces the morphism of the complexes
$\alpha:\bigoplus_{i=2}^l\nhom\bigl(V^{(i)},J^{(i)}\bigr)[-1]
\lrarr \gminig(V^{\ast})$.
We put $\gminigtilde(V^{\ast}):=\cone(\alpha)$.
By using the argument in the subsubsection
\ref{subsubsection;06.4.29.15},
we can show that
$\gminigtilde(V^{\ast})$ expresses
$\Phi(V^{\ast},\varphi_{\ast})^{\ast}\Ob(\Ytilde)$.

Applying the above construction
to $(\nbigv^{\ast},\varphi_{\ast}^u)$,
we obtain the morphism
$\obtilde(\nbigv^{\ast},\varphi_{\ast}^u):
 \Obtilde(\nbigv^{\ast},\varphi_{\ast}^u)
\lrarr L_{\widetilde{M}/S}$.

\begin{lem}
 \label{lem;06.5.3.10}
The morphism
$\obtilde(\nbigv^{\ast},\varphi^u_{\ast})$
gives an obstruction theory of
$\Mtilde(H_{\ast})$ over $S$.
\end{lem}
\pf
Let $h:T\lrarr \Mtilde(H_{\ast})$ be a morphism,
and let $J$ be a coherent sheaf on $T$.
The pull back of $J$ via $T\times _S\nbigd\lrarr T$
is denoted by $J_{\nbigd}$.
We put
$\widehat{h}_{\nbigd}:=
 \Phi(\nbigv^{\ast},\varphi^u_{\ast})\circ
 h_{\nbigd}$.
We have the following commutative diagram:
\[
 \begin{CD}
 \Ext^1\bigl(h^{\ast}L_{\Mtilde(H_{\ast})/S},\,J \bigr)
@>{\psi}>>
 \Ext^1\bigl(
 h^{\ast}\Obtilde,\,J
 \bigr)\\
 @VVV @A{\simeq}AA \\
 \Ext^1\bigl(h_{\nbigd}^{\ast}
 L_{\Mtilde(H_{\ast})\times_S\nbigd/\nbigd},\,J_{\nbigd} \bigr)
 @>>>
 \Ext^1\bigl(h_{\nbigd}^{\ast}\gminigtilde(\nbigv^{\ast}),\,J_{\nbigd}\bigr)
 \end{CD}
\]
Let $\Tbar$ be an $S$-scheme
such that 
 $T$ is embedded as a closed subscheme
and that the corresponding ideal $J$ is square $0$.
We have the obstruction classes
$o(h)$ and $o(\widehat{h}_{\nbigd})$ in 
$\Ext^1\bigl(h^{\ast}L_{\Mtilde(H_{\ast})/S},\,J \bigr)$
and 
$\Ext^1\bigl(
 \widehat{h}_{\nbigd}^{\ast}L_{\Ytilde/S},J_X
 \bigr)$.
It is easy to see that
$\psi\bigl(\ob(h)\bigr)\in
 \Ext^1\bigl(h^{\ast}\Obtilde,J\bigr)$ 
is same as the image of $o(\widehat{h}_{\nbigd})$
via the composite of the following morphisms:
\[
\begin{CD}
 \Ext^1\bigl(
 \widehat{h}_{\nbigd}^{\ast}L_{\Ytilde/S},J_{\nbigd}
 \bigr)
 @>{b_1}>>
 \Ext^1\bigl(
 h_{\nbigd}^{\ast}\gminigtilde(\nbigv^{\ast}),J_{\nbigd}
 \bigr)
 @>{b_2}>>
 \Ext^1\bigl(h^{\ast}\Obtilde,J\bigr)
\end{CD}
\]
Hence the vanishing of $\psi\bigl(o(h)\bigr)$ implies
$b_1\bigl(o(\widehat{h}_{\nbigd})\bigr)=0$.
Since $\Obtilde$ gives the obstruction theory for $\Ytilde$,
it implies that $\widehat{h}$ can be extended
to a morphism $\Tbar\times_S\nbigd\lrarr \Ytilde$.
Then we obtain the extension of $h$
to the morphism $\Tbar\lrarr \Mtilde(H_{\ast})$
due to the universal property of $\Mtilde(H_{\ast})$.
Therefore, the condition (A1) of 
Proposition \ref{prop;06.4.29.10}
is checked.
The condition (A2) can also be checked
easily,
and the proof of Lemma \ref{lem;06.5.3.10}
is finished.
\hfill\qed

\vspace{.1in}

Let $\pi$ denote the projection
$\Mtilde(H_{\ast})\lrarr M(H_{\ast})$,
which is smooth.
We have the following commutative diagram:
\[
 \begin{CD}
 \widetilde{M}(H_{\ast})\times_S \nbigd
 @>{\Phi_i(\nbigv^{\ast},\varphi_{\ast}^u)}>> 
 \Ytilde_i\\
 @V{\pi}VV @VVV \\
 M(H_{\ast})\times_S\nbigd @>{\Phi_i(\nbigv^{\ast})}>> Y_i
 \end{CD}
\]
Then, we obtain the following morphism of
the distinguished triangles on $\Mtilde(H_{\ast})\times_S\nbigd$:
\[
\begin{CD}
 \pi^{\ast}G(\nbigv^{\ast})
@>>>
 \Phi(\nbigv^{\ast},\varphi_{\ast}^u)^{\ast}
 \Ob(\Ytilde)
@>>>
\bigoplus_{i=2}^l\nhom(\nbigv^{(i)},J^{(i)})
@>>>
 \pi^{\ast}G(\nbigv^{\ast})[1]\\
 @VVV @VVV @VVV @VVV \\
 \pi^{\ast}L_{M(H_{\ast})\times_S\nbigd/\nbigd}
@>>>
 L_{\Mtilde(H_{\ast})\times_S\nbigd/\nbigd}
@>>>
 L_{\Mtilde(H_{\ast})\times_S\nbigd/M(H_{\ast})\times_S\nbigd}
@>>>
 \pi^{\ast}L_{M(H_{\ast})\times_S\nbigd/\nbigd}[1]
\end{CD}
\]
Hence, we obtain the following morphism
of the distinguished triangles:
\[
 \begin{CD}
 \pi^{\ast}\Ob(\nbigv^{\ast})
 @>>>
 \Obtilde(\nbigv^{\ast},\varphi^{\ast})
 @>>>
 \bigoplus_i
\bigl(
 p_{\ast}\nhom(J^{(i)},\nbigv^{(i)})
\bigr)^{\lor}
 @>>>
 \pi^{\ast}\Ob(\nbigv^{\ast})[1]\\ 
 @VVV @VVV @V{\varphi}VV @VVV \\
 \pi^{\ast}L_{M(H_{\ast})}
 @>>>
 L_{\Mtilde(H_{\ast})/S}
 @>>>
 L_{\Mtilde(H_{\ast})/M(H_{\ast})}
 @>>>
 \pi^{\ast}L_{M(H_{\ast})/S}[1]
 \end{CD}
\]
It is easy to see that both of
$L_{\Mtilde(H_{\ast}/S)}$ and 
$\bigoplus_i\bigl(
p_{\ast}\nhom(J^{(i)},\nbigv^{(i)})\bigr)^{\lor}$
are isomorphic to the $0$-th cohomology sheaves,
and that the morphism $\varphi$ is isomorphic.
Then, the claim of Lemma \ref{lem;06.5.3.15}
immediately follows from Lemma  \ref{lem;06.5.3.10}.
\hfill\qed

\subsection{Equivariant Complexes on 
 Deligne-Mumford Stacks with GIT Construction}
\label{subsection;06.7.3.6}
The results in this subsection will be used
when we discuss the equivariant obstruction theory 
of the master space in the subsection 
\ref{subsection;06.6.21.2}.

\subsubsection{Locally free resolution}
\label{subsubsection;06.5.22.555}

Let $G_i$ $(i=1,2)$ be a linear reductive group over $k$.
Let $U$ be a quasi projective variety over $k$
provided with $G_1\times G_2$.
We assume that there exists a
$G_1\times G_2$-embedding into some
projective space $\proj^N$.
The closure of $U$ in $\proj^N$
is denoted by $\Ubar$.
The $G_1\times G_2$-equivariant polarization
is denoted by $\nbigo(1)$.
We assume that $U$ is contained in
the open subset of the stable points of $\Ubar$
with respect to the polarization $\nbigo(1)$
and the $G_2$-action.
We assume that $\nbigm=U/G_2$ is a separated
Deligne-Mumford stack.
The projection $U\lrarr\nbigm$ is denoted by
$\pi$.

\begin{lem}
Let $\nbigf$ be a $G_1$-equivariant
quasi coherent sheaf on $\nbigm$.
Then there exists a $G_1$-equivariant
locally free sheaf $\nbigv$ on $\nbigm$
and a $G_1$-equivariant surjection
$\phi:\nbigv\lrarr\nbigf$.
\end{lem}
\pf
There exists a coherent sheaf $\nbigg$
on $\Ubar$ such that
$\nbigg_{|U}=\pi^{\ast}\nbigf$.
There exists a large number $N$
such that
$\nbigg(N)$ is globally generating.
Then $\pi^{\ast}\nbigf(N)$ is also globally generating.
We may take a $G_1\times G_2$-equivariant
subspace $W$ of
$H^0\bigl(U,\pi^{\ast}\nbigf(N)\bigr)$
such that 
$W\otimes\nbigo(-N)\lrarr \pi^{\ast}\nbigf$
is surjective.
Then we have only to take the descent of 
$W\otimes\nbigo(-N)$ and the morphism.
\hfill\qed

\begin{cor}
Let $\nbigf_{\cdot}$ be a bounded
$G_1$-equivariant complex of coherent sheaves
on $\nbigm$.
Assume that there exist integers
$M_1$ and $M_2$ such that the following holds:
\begin{itemize}
\item For any point of $\nbigm$,
  there exists a neighbourhood $\nbigu$
 such that
 $\nbigf_{\cdot|\nbigu}$ is isomorphic to
 a $G_1$-equivariant coherent locally free complex
 $\nbigg^{\nbigu}_{\cdot}$ in $D(\nbigu)$
 where $\nbigg^{\nbigu}_i=0$ unless
 $M_1\leq i\leq M_2$.
\end{itemize}
Then there exists a global $G_1$-equivariant 
coherent locally free complex
$\nbigg_{\cdot}\simeq\nbigf$ in $D(\nbigm)$,
where $\nbigg_i=0$ unless
 $M_1\leq i\leq M_2$.
\hfill\qed
\end{cor}

\subsubsection{Equivariant representative}
\label{subsubsection;06.6.22.5}

We recall that the morphism 
of $\nbigm$ to the coarse scheme is finite
(Proposition \ref{prop;06.4.25.1}).

\begin{lem}
 \label{lem;06.5.4.555}
Let $C_{i\,\cdot}$ $(i=1,2)$ be $G_1$-equivariant
bounded complexes of coherent sheaves on $\nbigm$.
We assume that $C_{1\,\cdot}$ is perfect.
Let $\varphi$ be an element of
the $G_1$-invariant part of
$\Ext^0(C_{1\,\cdot},C_{2\,\cdot})$.
Then, we can take a
$G_1$-equivariant perfect complex $\Ctilde_{1\,\cdot}$
with a $G_1$-equivariant morphism
$\psi:\Ctilde_{1\,\cdot}\lrarr C_2$,
such that $\Ctilde_{1\,\cdot}$
is $G_1$-equivariantly quasi isomorphic
to $C_{1\,\cdot}$,
and that $\psi$ represents $\varphi$.
\end{lem}
\pf
We give only an indication.
We may assume that
$C_{2,i}=0$ unless $|i|<N$.
We take a sufficiently large number $N_1$,
and we replace $C_{1\,\cdot}$
with a $G_1$-equivariant quasi isomorphic
complex $\Ctilde_{1\,\cdot}$
with the property
$\Ext^k\bigl(C_{1,i},C_{2,j}\bigr)=0$
for any $k>0$ and $i>-N_1$,
and for any $j$.
Then, $\Ext^0(C_{1},C_2)\simeq \Ext^0(\Ctilde_1,C_2)$
is isomorphic to the first cohomology of the following:
\[
 \bigoplus_{-i+j=-1}
 \Ext^0(\Ctilde_{1,i},C_{2,j})
\lrarr
 \bigoplus_{-i+j=0}
 \Ext^0(\Ctilde_{1,i},C_{2,j})
\lrarr
 \bigoplus_{-i+j=1}
 \Ext^0(\Ctilde_{1,i},C_{2,j})
\]
Since $G_1$ is assumed to be reductive,
the claim is clear.
\hfill\qed

\vspace{.1in}
Let $B^{(i)}$ $(i=1,2)$ be $G_1$-equivariant
bounded complexes on $\nbigm$.
We assume that $B^{(1)}$ is perfect.
Let $\phi$ be an element of 
the $G_1$-invariant part of
$\Ext^0\bigl(B^{(1)}, B^{(2)}\bigr)$.
We take a $G_1$-equivariant perfect complex
$\Btilde^{(1)}$ with $G_1$-equivariant morphisms
$a_i:\Btilde^{(1)}\lrarr B^{(i)}$
such that $a_1$ is quasi isomorphic,
and that the diagram
$B^{(1)}\stackrel{a_1}{\llarr}
 \Btilde^{(1)}\stackrel{a_2}{\lrarr}
 B^{(2)}$ represents $\phi$.
We have the natural $G_1$-equivariant 
structure on the cone $\cone(a_2)$.
Assume we have other $G_1$-equivariant complex
$\Bhat^{(1)}$ with $G_1$-equivariant morphisms
$\widehat{a}_i:\Bhat^{(1)}\lrarr B^{(i)}$
such that the diagram
$B^{(1)}\stackrel{\widehat{a}_1}{\llarr}
 \Bhat^{(1)}\stackrel{\widehat{a}_2}{\lrarr} B^{(2)}$
represents $\phi$.
Then, 
there exists $G_1$-equivariant complex
$\Bbar^{(1)}$ with $G_1$-equivariant morphisms
with morphisms
$f:\Bbar^{(1)}\lrarr \Btilde^{(1)}$
and $g:\Bbar^{(1)}\lrarr \Bhat^{(1)}$
such that the following diagrams are commutative
up to homotopy for $i=1,2$:
\[
 \begin{CD}
 \Bbar^{(1)} @>{f}>> \Btilde^{(1)}\\
 @V{g}VV @V{a_i}VV \\
 \Bhat^{(1)}@>{\widehat{a}_i}>> B^{(i)}
 \end{CD}
\]
Due to an argument similar to the proof of 
Lemma \ref{lem;06.5.4.555},
we may assume that the homotopy is also $G_1$-equivariant.
Then, we have the $G_1$-equivariant quasi isomorphisms:
\[
 \cone(\widehat{a}_2)
\llarr
 \cone(\widehat{a}_2\circ g)
\simeq
 \cone(a_2\circ f)
\lrarr
 \cone(a_2)
\]
In this sense,
the $G_1$-equivariant complex
$\cone(a_2)$ is uniquely determined
up to $G_1$-equivariant quasi isomorphisms.
We denote it by $\cone(\varphi)$.

\subsection{Elementary Remarks on some Extreme Sets}
\label{subsection;06.7.3.7}

The results in this subsection will be used
when we discuss the geometric invariant
theory for the enhanced master space
in the subsections 
\ref{subsection;06.5.21.110}--\ref{subsection;06.6.6.2}.

\subsubsection{Preparation for the proof of Proposition \ref{prop;06.5.16.50}}
\label{subsubsection;06.5.15.150}

Let us consider a vector space
$\nbigu=\bigoplus_{i=1}^N\rnum\cdot e_i$.
We put
$f_j:=(j-N)\sum_{i\leq j}e_i+j\cdot \sum_{i>j}e_i$.
The following lemma is well known and easy to prove.

\begin{lem}
 \label{lem;06.4.22.4}
Take any element $\rho=\sum_{i=1}^N a_i\cdot e_i\in\nbigu$
satisfying $\sum_{j=1}^N a_j=0$ and 
$a_1\leq a_2\leq \cdots\leq a_N$.
Then there exist non-negative rational numbers
$b_j$ such that
$\rho=\sum b_j\cdot f_j$.
\hfill\qed
\end{lem}

Let $r_1,\ldots,r_s$ be positive integers
such that $\sum_{j=1}^s r_j=N$.
We put $R_j=\sum_{i\leq j}r_i$.
We put as follows:
\[
 v_j:=\sum_{R_{j-1}<i\leq R_j}e_i,
\quad
(j=1,\ldots,s).
\]
For an integer $j$ such that $1\leq j\leq s$,
we put as follows:
\[
 y(j):=
 -(N-R_j)\sum_{h\leq j}v_h
+R_j\cdot \sum_{h> j}v_h.
\]
For a pair of integers $(i_1,i_2)$
such that $1\leq i_1<i_2\leq s$,
we put as follows:
\[
 x(i_1,i_2):=
 -(N-R_{i_2})\sum_{h\leq i_1}v_h
+R_{i_1}\sum_{i_2<h}v_h.
\]
For an integer $i_0$ such that $1\leq i_0\leq s$,
we put as follows:
\[
 \nbigs(i_0):=\Bigl\{
 (i_1,i_2)\in\seisuu^2\,\Big|\,
 1\leq i_1<i_0<i_2\leq s
\Bigr\}.
\]

\begin{lem}
\label{lem;06.4.22.105}
Let $v=\sum_{j=1}^s a_j\cdot v_j$ be
an element of $\nbigu$ satisfying the following:
\begin{equation}
 \label{eq;06.4.22.1}
 a_1\leq a_2\leq\cdots\leq a_s,
\quad
 \sum_{j=1}^sr_j\cdot a_j=0.
\end{equation}
Take an integer $i_0$ such that $1\leq i_0\leq s$.
\begin{itemize}
\item
 Assume $a_{i_0}>0$.
 Then there exist the non-negative rational numbers
 $b(i_1,i_2)\in\rnum_{\geq 0}$
 for $(i_1,i_2)\in\nbigs(i_0)$ and
 the non-negative rational numbers
 $c_j$ $(1\leq j<i_0)$
 such that the following equality holds:
\begin{equation}
 \label{eq;06.4.22.2}
 v=\sum_{(i_1,i_2)\in\nbigs(i_0)}
 b(i_1,i_2)\cdot x(i_1,i_2)+
 \sum_{j=1}^{i_0-1} c_j\cdot y(j).
\end{equation}
\item
Assume $a_{i_0}=0$.
Then there exist the non-negative rational numbers
 $b(i_1,i_2)\in\rnum_{\geq 0}$
 for $(i_1,i_2)\in\nbigs(i_0)$
such that the following holds:
\[
v=\sum_{(i_1,i_2)\in\nbigs(i_0)}
 b(i_1,i_2)\cdot x(i_1,i_2).
\]
One of $b(i_1,i_2)$ is not $0$.
\item
Assume $a_{i_0}<0$.
 Then there exist the non-negative rational numbers
 $b(i_1,i_2)\in\rnum_{\geq 0}$
 for $(i_1,i_2)\in\nbigs(i_0)$ and
 the non-negative rational numbers
 $c_j$ $(i_0< j\leq N)$
 such that the following holds:
\[
 v=\sum_{(i_1,i_2)\in\nbigs(i_0)}
 b(i_1,i_2)\cdot x(i_1,i_2)+
 \sum_{j=i_0+1}^{N} c_j\cdot y(j).
\]
\end{itemize}
\end{lem}
\pf
We use an induction on the number
$d(v):=\#\bigl\{i\,\big|\,a_i\neq a_{i+1}\bigr\}$.
In the case $d(v)=0$, the claim is obvious.
Let $v$ be as in the lemma
such that $d(v)=m+1$.
Take the integers $h_1$ and $h_2$
satisfying the following:
\[
 a_1=a_2=\cdots =a_{h_1}<a_{h_1+1},
\quad
 a_s=a_{s-1}=\cdots =a_{h_2+1}>a_{h_2}.
\]
We remark the following:
\begin{itemize}
\item
In the case $a_{i_0}>0$,
we have $h_1<i_0$.
\item
In the case $a_{i_0}=0$,
we have $h_1<i_0<h_2$.
\item
In the case $a_{i_0}<0$,
we have $i_0<h_2$.
\end{itemize}

Let us discuss the case $a_{i_0}>0$.
If we have $i_0\leq h_2$,
we put as follows:
\[
v':=v-f\cdot x(h_1,h_2)=\sum a_i'\cdot v_i,
\quad
 f:=\min\left\{
 \frac{a_{h_1+1}-a_{h_1}}{N-R_{h_2}},
 \,\frac{a_{h_2+1}-a_{h_2}}{R_{h_1}}
 \right\}
\]
If we have $i_0\geq h_2+1$,
we put as follows:
\[
v':=v-g\cdot y(h_2),
\quad
 g:=\frac{a_{i_0}-a_{h_2}}{2R_{h_2}}
\]
It is easy to see that
the numbers $a_i'$ satisfy the condition
(\ref{eq;06.4.22.1}),
and that we have $d(v')\leq d(v)-1$.
Due to the hypothesis of the induction,
we have the expression for $v'$
as in (\ref{eq;06.4.22.2})
with the non-negative coefficients.
Then we obtain the desired expression for $v$.

The cases $a_{i_0}=0$ or $a_{i_0}<0$
can be discussed similarly.
\hfill\qed

\subsubsection{Preparation for the proof of Proposition \ref{prop;06.5.16.20}}
\label{subsubsection;06.4.23.15}

Let $N^{(\alpha)}$ $(\alpha=1,2)$ be positive integers.
Let us consider a vector space as follows:
\[
 \nbigu=\nbigu^{(1)}\oplus\nbigu^{(2)},
\quad
\nbigu^{(\alpha)}=
 \bigoplus_{i=1}^{N^{(\alpha)}} \rnum\cdot e^{(\alpha)}_i.
\]
Let $r_1^{(\alpha)},\ldots,r_{s(\alpha)}^{(\alpha)}$
$(\alpha=1,2)$ be positive integers
such that
$\sum_{j=1}^{s(\alpha)}r_j^{(\alpha)}=N^{(\alpha)}$.
We put
$R^{(\alpha)}_j=\sum_{h\leq j}r^{(\alpha)}_h$.
We put
$\Omega^{(\alpha)}=\sum_i e^{(\alpha)}_i$.
We also put as follows:
\[
 v^{(\alpha)}_j:=
 \sum_{R^{(\alpha)}_{j-1}<i\leq R^{(\alpha)}_j}
 e^{(\alpha)}_i,
\quad \bigl(j=1,\ldots, s(\alpha)\bigr).
\]
For each integer $j$ such that
$1\leq j\leq s(2)$,
we put as follows:
\[
 y^{(2)}(j)=
 -(N-R_j^{(2)})\cdot\sum_{h\leq j}v^{(2)}_h
+R_j^{(2)}\cdot \sum_{h>j}v^{(2)}_h.
\]
For each integer $j$ such that
$1\leq j\leq s(1)$,
we put as follows:
\[
 x_1(j):=
 -N^{(2)}\cdot\sum_{h\leq j}v^{(1)}_h
+R^{(1)}_j\cdot\Omega^{(2)},
\quad
 x_2(j):=
 N^{(2)}\cdot\sum_{h\geq j}v^{(1)}_h
+(R^{(1)}_{j-1}-N^{(1)})\cdot\Omega^{(2)}.
\]

\begin{lem}
 \label{lem;06.4.23.10}
Let 
$v=\sum_{\alpha=1,2}\sum_j
 a_{j}^{(\alpha)}\cdot v^{(\alpha)}_j$
be any element of $\nbigu$
satisfying the following conditions:
\begin{equation}
\label{eq;06.4.22.3}
 a_1^{(\alpha)}\leq a_2^{(\alpha)}\leq\cdots
 \leq a_{s(\alpha)}^{(\alpha)},
\quad
 \sum_{\alpha=1,2}\sum_{j}
 r^{(\alpha)}_j\cdot a^{(\alpha)}_{j}=0.
\end{equation}
Take an integer $i_0$
such that $1\leq i_0\leq s(1)$.
Then, there exist 
non-negative rational numbers
$c(j)\geq 0$ $(j=1,\ldots,s(2))$,
$d_1(i)\geq 0$ $(i=1,\ldots,i_0)$,
$d_2(i)\geq 0$ $(i=i_0+1,\ldots,s(1))$
and a rational number $A$
such that the following holds:
\begin{equation}
 \label{eq;06.4.22.4}
v=\sum_{j=1}^{s(2)}c(j)\cdot y^{(2)}(j)
 +\sum_{i<i_0}d_1(i)\cdot x_1(i)
 +\sum_{i>i_0}d_2(i)\cdot x_2(i)
 +A\cdot\bigl(N^{(2)}\Omega^{(1)}
 -N^{(1)}\cdot \Omega^{(2)}\bigr).
\end{equation}
\end{lem}
\pf
Due to Lemma \ref{lem;06.4.22.4},
we may assume 
$a^{(2)}_1=\cdots =a^{(2)}_{s(2)}$
from the beginning.
We use an induction on the number
$d(v)=\#\bigl\{
 i\,|\,a^{(1)}_i\neq a^{(1)}_{i+1}
 \bigr\}$.
In the case $d(v)=0$,
we have
$v=A\cdot \bigl(N^{(2)}\cdot \Omega^{(1)}
-N^{(1)}\cdot \Omega^{(1)}\bigr)$
for some $A$,
and hence the claim is clear.
Let $v$ be an element as in the lemma
such that $d(v)=m+1>0$.
Let us take the integer $h_1$
satisfying 
$a^{(1)}_1=a^{(1)}_{h_1}<a^{(1)}_{h_1+1}$.
In the case $i_0>h_1$,
we put as follows:
\[
 v':=v-
 \frac{a^{(1)}_{h_1+1}-a^{(1)}_{h_1}} {N^{(2)}}
 x_1(h_1)
\]
In the case $i_0\leq h_1$,
we put as follows:
\[
 v'=v-
 \frac{a^{(1)}_{h_1+1}-a^{(1)}_{h_1}}
 {N^{(2)}}
 \cdot x_2(h_1).
\]
Then $v'$ satisfies the condition (\ref{eq;06.4.22.3})
and $d(v')<d(v)$.
Due to the hypothesis of the induction,
we have the expression for $v'$
as in (\ref{eq;06.4.22.4}).
Hence we obtain the desired expression for $v$.
\hfill\qed

\subsection{Twist of Line Bundles}
\label{subsection;06.5.17.5}

This subsubsection is a preparation for 
the subsection \ref{subsection;06.6.21.3}.

\subsubsection{Construction}
\label{subsubsection;06.5.14.100}

Let $Y$ be an algebraic stack over 
a field $k$.
Let $G_m$ denote the one dimensional
algebraic torus $\Spec k[t,t^{-1}]$.
Let $I$ denote the trivial line bundle on $Y$.
A point of $I$ is denoted
by $(y,u)$ where $y\in Y$ and $u\in I_{|y}$.
For each integer $n$,
$\nbigt(n)$ denote the line bundle $I$
with the $G_m$-action
by $t\cdot (y,u):=(y,t^n\!\cdot\! u)$.

Let $L$ be any line bundle on $Y$.
Let $L^{\ast}$ denote the complement of
the image of the $0$-section,
i.e., $L^{\ast}:=L-Y$.
Let $\pi:L^{\ast}\lrarr Y$ denote 
the naturally defined projection.
A point of $L^{\ast}$ is also denoted by
$(y,v)$, where $y\in Y$ and $v\in \pi^{-1}(y)$.

Let us fix an integer $r$.
We consider the $G_m$-action on $L^{\ast}$
given by $t\cdot (y,v):=(y,t^r\!\cdot\! v)$.
We have the naturally defined $G_m$-action
on $\pi^{\ast}\nbigt(n)$.
It induces the line bundle $\nbigi_n$
on the algebraic stack $L^{\ast}/G_m$.
Let $\varphi:L^{\ast}/G_m\lrarr Y$ denote the
naturally induced morphism.

\begin{lem}
\label{lem;06.5.17.6}
We have the canonical isomorphism
$\nbigi_n\otimes\nbigi_m\simeq\nbigi_{n+m}$
and $\nbigi_{-n}\simeq\nbigi_n^{-1}$
and $\nbigi_0\simeq \nbigo_{Y/G_m}$.
We also have the canonical isomorphism
$\nbigi_{-r}\simeq \varphi^{\ast}L$.
\end{lem}
\pf
The first claim is obvious.
Let us show the second claim.
Let us denote a point of $\pi^{\ast}L$
by $(y,v,u')$,
where $y\in Y$, $v\in \pi^{-1}(y)$
and $u'\in L_{|y}$.
The trivial $G_m$-action on $L$
induces the $G_m$-action
on $\pi^{\ast}L$ over $L^{\ast}$,
which is given by
$t\cdot (y,u,v')=\bigl(y,t^r\cdot u,v'\bigr)$.

On the other hand,
let us denote a point of $\pi^{\ast}\nbigt(-r)$
by $(y,v,u)$
where $y\in Y$, $v\in \pi^{-1}(y)$
and $u\in \nbigt(-r)_{|y}$.
The action is denoted by
$t\cdot (y,v,u)=\bigl(y,t^rv,t^{-r}u\bigr)$.

We have the naturally defined isomorphism
$\pi^{\ast}\nbigt(-r)\lrarr \pi^{\ast}L$
given by
$(y,u,v)\longmapsto
 (y,u,u\cdot v)$,
which is $G_m$-equivariant.
Therefore,
we obtain the isomorphism
$\nbigi_{-r}\simeq \varphi^{\ast}L$.
\hfill\qed

\subsubsection{The weight of the induced action}

We use the notation in the previous subsubsection.
Let $G_m^{(i)}$ denote the torus
$\Spec k[t_i,t_i^{-1}]$.
Let us consider the action of
$G_m^{(1)}\times G_m^{(2)}$
on $L$ given by
$(t_1,t_2)\cdot (y,v):=
\bigl(y,t_1\!\cdot\! t_2\!\cdot\! v\bigr)$.

Let $\nbigt(n_1,n_2)$ denote 
the trivial line bundle $I$
with the $G_m^{(1)}\times G_m^{(2)}$
given by
$(t_1,t_2)\cdot (y,u)=
\bigl(y,t_1^{n_1}\!\cdot\! t_2^{n_2}\!\cdot\! u\bigr)$.
Then we have 
the $G_m^{(1)}\times G_m^{(2)}$-line bundle
$\pi^{\ast}\nbigt(n_1,n_2)$ on $L^{\ast}$.
We obtain the line bundle $\nbigi_{n_2}$
on $L^{\ast}/G_m^{(2)}$,
and we have the induced $G_m^{(1)}$-action
on $\nbigi_{n_2}$.

\begin{lem}
\label{lem;06.5.17.17}
The weight of the $G_m^{(1)}$-action
on $\nbigi_{n_2}$ is $n_1-n_2$.
\end{lem}
\pf
We put $\Gtilde_m^{(i)}:=\Spec k[s_i]$.
Let us consider the morphism
$\Gtilde_m^{(1)}\times\Gtilde_{m}^{(2)}
\lrarr
 G_m^{(1)}\times G_m^{(2)}$
given by
$(s_1,s_2)\longmapsto \bigl(s_1,s_1^{-1}\cdot s_2\bigr)$.
The induced $\Gtilde_m^{(1)}\times\Gtilde_m^{(2)}$-action
on $L^{\ast}$ and $\nbigt(n_1,n_2)$ is given by
$(s_1,s_2)\cdot (y,v)=\bigl(y,s_2^{r}\cdot v\bigr)$
and
$(s_1,s_2)\cdot (y,u)=
 \bigl(y,s_1^{n_1-n_2}\cdot s_2^{n_2}\cdot u\bigr)$.
Therefore, the weight of the $G_m^{(1)}$-action 
on $\nbigi_{n_2}$ is given by $n_1-n_2$.
\hfill\qed

\section{Parabolic $L$-Bradlow Pairs}
\label{section;06.5.17.20}

In this section,
we recall some definitions.
All of them are more or less familiar.
The purpose is to fix the meaning in this paper.
In the following of this section,
$X$ will denote a smooth irreducible projective variety
over an algebraically closed field $k$ of characteristic $0$.
Let $\Pic_X$ denote the Picard variety of $X$.
We fix a base point $x_0\in X$,
due to which we have the Poincar\'{e} bundle
$\poin_X$ on $\Pic_X\times X$.

\subsection{Sheaves with some Structure
   and their Moduli Stacks}
\label{subsection;06.5.17.21}
\subsubsection{Orientation}
\label{subsubsection;06.7.3.211}

Let $E$ be a $U$-coherent sheaf on $U\times X$.
Then we have the morphism $\det_E:U\lrarr \Pic_X$
induced by the determinant line bundle  $\det(E)$ of $E$,
which satisfies the condition
$\det(E)_{|\{u\}\times X}
\simeq
 \poin_{X\,|\,\{\det_E(u)\}\times X}$.
The morphism will be denoted by simply $\det$,
if there are no risk of confusion.
In general, the line bundles $\det^{\ast}_E\poin_X$
and $\det(E)$ are not isomorphic.

\begin{example}
 \label{example;06.4.9.1}
{\rm 
Let $c$ be an element of the second cohomology group
$H^2(X)$,
and let $\Pic_X(c)$ denote the Picard variety
of the line bundles whose first Chern classes are  $c$.
Assume $H^i(X,L)=0$ $(i>0)$ for any line bundle
$L\in \Pic_X(c)$.
Then, $p_{X\ast}(\poin_X)$ 
gives the vector bundle on $\Pic_X(c)$.
We obtain the associated projective space bundle
$\proj_c=\proj\bigl(p_{X\ast}(\poin_X)^{\lor}\bigr)$
on $\Pic_X(c)$.

Let $\pi$ denote the natural projection
$\proj_c\times X\lrarr \Pic_X\times X$.
We have the line bundle
$L(a)=\pi^{\ast}\poin_X\otimes p_X^{\ast}\nbigo_{\proj_c}(a)$
for each $a\in\seisuu$.
Here $\nbigo_{\proj_c}(1)$
denotes the tautological bundle
of the projective space bundle
$\proj_c\lrarr \Pic_X(c)$,
and $\nbigo_{\proj_c}(a)=\nbigo_{\proj_c}(1)^{\otimes\,a}$.
The determinant bundle
$\det\bigl(L(a)\bigr)$ is obviously $L(a)$ itself.
On the other hand,
$\det_{L(a)}$ is given by the projection $\pi$.
Thus, $L(a)$ and $\det^{\ast}_{L(a)}\poin_X$
are not isomorphic, if $a\neq 0$.
\hfill\qed}
\end{example}

\begin{df}[Orientation]
\index{orientation}
An orientation of a $U$-coherent sheaf $E$ on $U\times X$
is defined to be
an isomorphism $\rho:\det(E)\lrarr \det^{\ast}_E\poin_X$ on $U$.
A tuple $(E,\rho)$ is called an oriented $U$-coherent sheaf.

An isomorphism of two oriented sheaves $(E,\rho)$
and $(E',\rho')$ is defined to be
an isomorphism $\chi:E\lrarr E'$
satisfying $\rho'=\chi^{\ast}(\rho):=\rho\circ\chi$.
\hfill\qed
\end{df}

The restrictions $\bigl(\det_E^{\ast}\poin_{X}\bigr)_{|\{u\}\times X}$
and $\det(E)_{|\{u\}\times X}$ are isomorphic
for any point $u\in U$ by definition of $\det_E$,
so that the push-forward $p_{X\ast}\nhom\bigl(\det(E),\det_E^{\ast}\poin_X\bigr)$
is the line bundle on $U$.

\begin{df}[Orientation bundle]
\index{orientation bundle}
The line bundle $p_{X\ast}\nhom\bigl(\det(E),\det_E^{\ast}\poin_X\bigr)$
is called the orientation bundle of $E$.
We denote it by $\Or(E)$.
\hfill\qed
\end{df}

If $E$ is oriented, then the orientation bundle
$\Or(E)$ is naturally
isomorphic to the trivial line bundle $\nbigo_U$,
i.e.,
an orientation is equivalent to a trivialization of $\Or(E)$.

\begin{example}{\rm
Let $L(a)$ be given in 
Example \ref{example;06.4.9.1}.
The orientation bundle $\Or(L(a),\proj_c)$
induced by $L(a)$ 
is isomorphic to
$ p_{X\ast}\Bigl(\nhom\bigl(L(a),
 \det{}_{L(a)}^{\ast}\poin\bigr)\Bigr)
\simeq
 p_{X\ast}\circ p_X^{\ast}\bigl(
 \nbigo_{\proj_c}(-a)\bigr)
\simeq
 \nbigo_{\proj_c}(-a)$.
\hfill\qed
}
\end{example}

We have the following additive property
of the orientation bundles.
\begin{lem}
Let $E_i$ $(i=1,2)$ be $U$-coherent sheaves
on $U\times X$.
Then we have the natural isomorphism
$Or(E_1\oplus E_2)\simeq Or(E_1)\otimes Or(E_2)$.
\end{lem}
\pf
The natural isomorphism
$\det(E_1\oplus E_2)\simeq \det(E_1)\otimes \det(E_2)$
is given,
and hence
$ \det_{E_1}^{\ast}\poin_X\otimes\det_{E_2}^{\ast}\poin_X
\simeq
 \det^{\ast}_{E_1\oplus E_2}\poin_X$.
It induces the following isomorphism:
\[
 \nhom\Bigl(\det(E_1\oplus E_2),
   \det{}^{\ast}_{E_1\oplus E_2}\poin_X\Bigr)
\simeq
 \nhom\bigl(\det(E_1),\det{}^{\ast}_{E_1}\poin_X\bigr)
\otimes
 \nhom\bigl(\det(E_2),\det{}^{\ast}_{E_2}\poin_X\bigr).
\]
Therefore,
we obtain the natural morphism
$\Or(E_1)\otimes\Or(E_2)\lrarr \Or(E_1\oplus E_2)$,
which is isomorphic.
\hfill\qed

\subsubsection{Parabolic structure}
\label{subsubsection;06.6.9.25}

See \cite{my} and \cite{y} for detail on the notion
of parabolic sheaf.
Our terminology slightly differs from theirs.
We remark that it is different from
that in the author's other papers
(\cite{mochi2}, for example.)
Let $D$ be a Cartier divisor of $X$.
A $U$-parabolic sheaf,
or simply parabolic sheaf,
on $U\times (X,D)$ is defined to be a tuple
$\bigl(E,F_{\ast}(E),\alpha_{\ast}\bigr)$:
\begin{itemize}
\item $E$ is a $U$-coherent sheaf on $U\times X$.
\item
$F_{\ast}(E)$ denotes a filtration of $E$:
\[
 E=F_1(E)\supset F_2(E)\supset \cdots\supset F_l(E)
 \supset F_{l+1}(E)=E(-D).
\]
Here $E(-D)$ denotes $E\otimes p_U^{\ast}\nbigo(-D)$.
We assume that $\Cok_i(E):=E/F_{i+1}(E)$ are flat on $U$.
\item
$\alpha_{\ast}=(\alpha_1,\ldots,\alpha_l)$ is a tuple of numbers
$0< \alpha_1<\alpha_2<\cdots<\alpha_l\leq 1$.
 It is called a system of weights.
\end{itemize}
A tuple$(E,F_{\ast},\alpha_{\ast})$ will be often denoted simply by $E_{\ast}$.
The filtration $F_{\ast}$ is called a quasi-parabolic structure.
Isomorphisms of parabolic sheaves
are defined naturally.

The number $l$ is called the depth of the parabolic structure.
The tuple $\alpha_{\ast}$ will be called a weight
of the parabolic structure.
For any parabolic sheaf $E_{\ast}$,
we put $Gr_i(E):=F_i(E)/F_{i+1}(E)$,
which are the $U$-coherent sheaves
on $U\times D$.

\begin{rem}
We will often use the word ``parabolic''
even when a system of weights is not given.
\hfill\qed
\end{rem}

\noindent
{\bf (Subobject and quotient object)}\\
\index{subobject, quotient object of parabolic sheaves}
Let $E_{\ast}$ be a parabolic torsion-free sheaf defined over $U\times (X,D)$.
For any subsheaf $E'\subset E$
and any quotient sheaf $E\lrarr E''$,
we have the induced parabolic structures on $E'$
and $E''$.
Namely we put $F_i(E')=F_i(E)\cap E'$,
and $F_i(E'')=\Image(F_i(E)\lrarr E'')$.
The parabolic structures are called the induced parabolic structures.
We always consider the induced parabolic structures
on the subsheaves and the quotient sheaves.

\vspace{.1in}
\noindent{}
{\bf (The Condition $O_m$)}\\
\index{condition $O_m$}
Let $(E,F_{\ast})$ be a $U$-quasi-parabolic sheaf
on $U\times (X,D)$.
Let $m$ be a positive integer.
We say that $(E,F_{\ast})$ satisfies the condition $O_{m}$,
if the following holds:
\begin{itemize}
\item
 $F_i(E)(m)_{|\{u\}\times X}$
 and $\Cok_i(E)(m)_{|\{u\}\times X}$ are generated by its global sections,
 for all $i=1,\ldots,l+1$ and for all $u\in U$.
\item
 The higher cohomology groups
 of $F_i(E)(m)_{|\{u\}\times X}$
 and $\Cok_i(E)(m)_{|\{u\}\times X}$ vanish,
for all $i=1,\ldots,l+1$ and for all $u\in U$.
\end{itemize}

When we are given a $U$-quasi-parabolic sheaf
$(E,F_{\ast})$ on $U\times (X,D)$,
the open subset $U'$ is determined
by the condition $O_m$.

\vspace{.1in}
\noindent
{\bf (Twist)}\\
Let $m$ be an integer, and let $E$ be
a $U$-coherent sheaf defined over $U\times X$.
Recall that $E(m)$ denotes the coherent sheaf
$E\otimes p_U^{\ast}\nbigo_X(m)$.
If $E$ has a quasi parabolic structure $F_{\ast}(E)$,
we have the naturally induced parabolic structure
$F_{\ast}\bigl(E(m)\bigr)$ of $E(m)$.
The tuple $\bigl(E(m),F_{\ast}(E(m))\bigr)$
is denoted by $E_{\ast}(m)$.

\subsubsection{$L$-Bradlow pairs and reduced $L$-Bradlow pairs}
\label{subsubsection;06.7.3.212}

Let $L$ be a line bundle over $X$.

\begin{df}[$L$-section]
\index{$L$-section}
Let $E$ be a $U$-coherent sheaf on $U\times X$.
A morphism $\phi:p_{U}^{\ast}L\lrarr E$ is called an $L$-section.
An $\nbigo_X$-section is simply called a section.
\hfill\qed
\end{df}

\begin{df}
Let $(E,\phi)$ be a pair of a $U$-coherent sheaf on $X$
and an $L$-section.
We say that $\phi$ is non-trivial everywhere,
if $\phi_{|\{u\}\times X}\neq 0$ for every point $u\in U$.
\hfill\qed
\end{df}

\begin{df}
 \mbox{{}}
\begin{description}
\item[($L$-Bradlow pair)]
\index{$L$-Bradlow pair}
A parabolic $L$-Bradlow pair $(E_{\ast},\phi)$ on $U\times (X,D)$ is a pair
of a $U$-torsion-free parabolic sheaf $E_{\ast}$ and
an $L$-section $\phi:\pi_U^{\ast}L\lrarr E$.

An isomorphism between two such pairs $(E_{\ast},\phi)$ and $(E_{\ast}',\phi')$
is defined to be an isomorphism $\chi:E_{\ast}\lrarr E_{\ast}'$
satisfying $\phi'=\chi_{\ast}(\phi):=\chi\circ\phi$.

\item[(Oriented $L$-Bradlow pair)]
An oriented parabolic $L$-Bradlow pair
$(E_{\ast},\phi,\rho)$ on $U\times(X,D)$
is a pair of a parabolic $L$-Bradlow pair $(E_{\ast},\phi)$
and an orientation of $E$.
An isomorphism of two such pairs
is defined naturally.
\hfill\qed
\end{description}
\end{df}

\begin{rem}
We are mainly interested in 
parabolic $L$-Bradlow pairs $(E_{\ast},\phi)$
such that $\phi$ is non-trivial
everywhere.
Sometimes, we will assume it without mention.
\hfill\qed
\end{rem}

\begin{df}
\index{subobject of $L$-Bradlow pair}
Let $(E_{\ast},\phi)$ and $(E_{\ast}',\phi')$ be 
parabolic $L$-Bradlow pairs on $X$.
We say that $(E'_{\ast},\phi')$ is a subobject of
$(E_{\ast},\phi)$ if the following conditions hold:
\begin{itemize}
\item
 $E'$ is a subsheaf of $E$,
 and the parabolic structure is same as
 the induced one.
\item
If the image of $\phi$ is contained in $E$,
 we have $\phi'=\phi$.
 Otherwise, $\phi'=0$.
\hfill\qed
\end{itemize}
\end{df}

We also introduce the notion of reduced $L$-section.

\begin{df}[Reduced $L$-section]
\index{reduced $L$-section}
Let $L$ be a line bundle over $X$.
Let $E$ be a $U$-coherent sheaf on $U\times X$.
A reduced $L$-section of $E$
is defined to be a pair $(M,[\phi])$ of a line bundle $M$
on $U$ and a morphism
$[\phi]:p_{X}^{\ast}(M)\otimes p_U^{\ast}(L)
 \lrarr E$.

A reduced $L$-section is often 
denoted  simply by $[\phi]$
instead of $(M,[\phi])$,
if there are no risk of confusion.
\hfill\qed
\end{df}

\begin{df}
Let $(E,[\phi])$ be a pair of a $U$-coherent sheaf on
$U\times X$ and a reduced $L$-section.
We say that $[\phi]$ is non-trivial everywhere,
if $[\phi]_{|\{u\}\times X}\neq 0$  for each $u\in U$.
\hfill\qed
\end{df}

Let $(M,[\phi])$ be a reduced $L$-section of $E$
which is non-trivial everywhere.
Then, $(M,[\phi])$ induces the morphism
$\overline{[\phi]}:M\lrarr p_{X\ast}\nhom(L,E)$.
Under the vanishing
$H^i\bigl(X,\nhom(L,E)_{|\{u\}\times X}\bigr)=0$ $(i>0)$
for each $u\in U$,
we obtain the naturally induced section
$U\lrarr \proj\bigl(p_{X\,\ast}(\nhom(L,E))^{\lor}\bigr)$.
On the other hand,
such a section
$U\lrarr \proj\bigl(p_{X\,\ast}(\nhom(L,E))^{\lor}\bigr)$
induces a reduced $L$-section
which is non-trivial everywhere.

\begin{df}[Parabolic reduced $L$-Bradlow pair]
\index{parabolic reduced $L$-Bradlow pair}
\index{reduced $L$-Bradlow pair}
\mbox{{}}
\begin{itemize}
\item
A parabolic reduced $L$-Bradlow pair
$(E_{\ast},M,[\phi])$ on $U\times (X,D)$
is defined to be a pair
of a torsion-free parabolic sheaf $E_{\ast}$
on $U\times (X,D)$ and a reduced $L$-section $(M,[\phi])$
which is non-trivial everywhere.
It is often denoted by $(E_{\ast},[\phi])$
instead of $(E_{\ast},M,[\phi])$.

An isomorphism of two reduced $L$-Bradlow pairs
$(E_{i\,\ast},M_i,[\phi_i])$ $(i=1,2)$ is defined to be
a pair $(\chi,\eta)$ of
an isomorphism $\chi:E_{1\,\ast}\lrarr E_{2\,\ast}$
and $\eta:M_1\simeq M_2$
such that the following diagram is commutative:
\[
 \begin{CD}
 p_X^{\ast}M_1\otimes p_U^{\ast}L @>{\phi_1}>> E_1 \\
 @V{\eta\otimes\id_L}VV @V{\chi}VV\\
 p_X^{\ast}M_2\otimes p_U^{\ast}L @>{\phi_2}>> E_2
 \end{CD}
\]

\item
An oriented parabolic reduced $L$-Bradlow pair
is defined to be a tuple of
a parabolic reduced $L$-Bradlow pair $(E_{\ast},[\phi])$
with an orientation $\rho$ of $E$.

An isomorphism between two oriented reduced $L$-Bradlow pair
is naturally defined.
\hfill\qed
\end{itemize}
\end{df}

\begin{rem}
In the definition,
we assume that $[\phi]$ is non-trivial
everywhere,
contrast to the definition of non-reduced $L$-Bradlow pair.
\end{rem}

\begin{rem}
In the case $U=\Spec(k)$,
a parabolic reduced $L$-Bradlow pair
is just a parabolic $L$-Bradlow pair
whose $L$-section is non-trivial,
up to isomorphisms.
\hfill\qed
\end{rem}

\begin{rem}
We will often use the word
``$L$-Bradlow pair'' and
``reduced $L$-Bradlow pair'' 
instead of ``parabolic $L$-Bradlow pair''
and ``parabolic reduced $L$-Bradlow pair'',
if there are no risk of confusion.

We will also use the word
``quasi-parabolic $L$-Bradlow pair''
and ``quasi-parabolic reduced $L$-Bradlow pair'',
if a system of weights is not given.
\hfill\qed
\end{rem}

We will use the following notion subordinately.
\begin{df}

Let $\vecL=(L_1,L_2)$ be a tuple of line bundles on $X$.

\begin{itemize}
\item
\index{parabolic $\vecL$-Bradlow pair}
\index{$\vecL$-Bradlow pair}
A parabolic $\vecL$-Bradlow pair on $U\times (X,D)$
is defined to be a tuple
$(E_{\ast},\vecphi)$
of a $U$-parabolic torsion-free sheaf $E_{\ast}$ on $U\times (X,D)$
and a pair $\vecphi$
of $L_i$-sections $\phi_i$ $(i=1,2)$.
\item
\index{parabolic reduced $\vecL$-Bradlow pair}
\index{reduced $\vecL$-Bradlow pair}
An oriented parabolic reduced $\vecL$-Bradlow pair on $U\times (X,D)$
is defined to be a tuple
$(E_{\ast},[\vecphi],\rho)$
of a $U$-parabolic torsion-free sheaf $E_{\ast}$ on $U\times (X,D)$,
a pair $[\vecphi]$ of reduced $L_i$-sections $[\phi_i]$
which are non-trivial everywhere,
and an orientation $\rho$. 
Isomorphisms are defined naturally.
\hfill\qed
\end{itemize}
\end{df}

\subsubsection{Type and the moduli stacks}
\label{subsubsection;06.4.24.10}

Let $H^{\ast}$ denote an appropriate cohomology theory
with the appropriate Chern class for vector bundles.
We put $H^{\ev}(X):=\bigoplus_{i=0}^{\dim X} H^{2i}(X)$.
If $D$ is a smooth divisor,
we put as follows:
\[
 \Typetilde:=
\Bigl\{(y,y_1,y_2,\ldots \ldots)\in
H^{\ev}(X)\oplus
 \bigoplus_{i=1}^{\infty} H^{\ev}(Y)\,\Big|\,
 \sum_{i\geq 1}y_i= y_{|D}
\Bigr\}.
\]
In general,
we put as follows:
\[
 \Typetilde:=
\Bigl\{(y,y_1,y_2,\ldots \ldots)\in
H^{\ev}(X)\oplus
 \bigoplus_{i=1}^{\infty} H^{\ev}(X)\,\Big|\,
 \sum_{i\geq 1}y_i=
 y\cdot \Bigl(1-\ch\bigl(\nbigo(-D)\bigr)\Bigr)
\Bigr\}.
\]
In the following,
$y\cdot \Bigl(1-\ch\bigl(\nbigo(-D)\bigr)\Bigr)$
is denoted by $y_{|D}$
for simplicity of the notation,
even when $D$ is not necessarily smooth.

For any quasi-parabolic sheaf $(E,F_{\ast})$ on $X$
of depth $l$,
we obtain the element of $\Typetilde$:
\[
 \type\bigl(E,F_{\ast}\bigr):=
 \Bigl(\ch(E),\ch\bigl(\Gr_1(E)\bigr),
 \ldots,\ch\bigl(\Gr_{l-1}(E)\bigr),
 \ch\bigl(\Gr_{l}(E)\bigr),0,\ldots\Bigr)
 \in \Typetilde
\]
Let $\Type$ denote the subset of $\Typetilde$
which consists of 
$\type\bigl(E,F_{\ast}\bigr)$
for some quasi parabolic sheaves
$(E,F_{\ast})$.
\index{$\Type$}

Let $\vecy=(y,y_1,y_2,\ldots \ldots)$ be any element of $\Type$.
The number $\depth(\vecy):=\max\{i\,|\,y_i\neq 0\}$
is called the depth of $\vecy$.
The element $y$ is called the $H^{\ev}(X)$-component of $\vecy$,
and $(y_1,y_2,\ldots,)$ is called the parabolic part of $\vecy$.
The $H^0(X)$-component of $y$ is called the rank of $\vecy$,
and it is denoted by $\rank(\vecy)$ or $\rank y$.
\index{$\rank(\vecy)$, $\rank(y)$}
We put $\Type_r:=\bigl\{\vecy\in\Type\,\big|\,\rank(\vecy)=r\bigr\}$.
We denote by $\Type^{\circ}$
the set of types whose parabolic components are trivial,
i.e., $y_1=y_{|D}$ and $y_i=0$ for $i>1$.
We often identify $y$ and $\vecy$ in that case,
and we regard $\Type^{\circ}$ as the subset of
$H^{\ast}(X)$.
We put $\Type^{\circ}_r:=\Type^{\circ}\cap\Type_r$.
\index{$\Type_r$, $\Type^{\circ}$, $\Type_r^{\circ}$}

We have the sum
$\vecy^{(1)}+\vecy^{(2)}$
of two elements $\vecy^{(i)}\in\Type$ $(i=1,2)$
by taking the component summation.

For any quasi-parabolic sheaf $(E,F_{\ast})$ on $X$,
we obtain the element 
$\type\bigl(E,F_{\ast}\bigr)$
of $\Type$.
In general,
let $(E,F_{\ast})$ be a $U$-quasi-parabolic sheaf
on $U\times X$.
When we have the element $\vecy\in\Type$
such that $\type(E_{\ast|u})=\vecy$
for any closed point $u\in U$,
then $E_{\ast}$ is called of type $\vecy$.
When $U$ is connected,
such an element always exists.

\begin{df}
\index{type}
The type of an (oriented) parabolic $L$-Bradlow pairs
is defined to be the type
of the underlying quasi-parabolic sheaf.
\hfill\qed
\end{df}

We introduce the following notation.
\begin{notation}
In each line,
the left hand side denotes the moduli stack
of the object in the right hand side:

\begin{description}
\item[$\nbigm(\vecy)$:]
  Parabolic sheaves of type $\vecy$.
\item[$\nbigm(\vecyhat)$:]
 Oriented parabolic sheaves of type $\vecy$.
 \index{$\nbigm(\vecy)$, $\nbigm(\vecyhat)$}
\item[$\nbigm\bigl(\vecy,L\bigr)$:]
Parabolic $L$-Bradlow pairs of type $\vecy$
whose $L$-sections are non-trivial everywhere.
\item[$\nbigm\bigl(\vecyhat,L\bigr)$:]
 Oriented parabolic $L$-Bradlow pairs
 of type $\vecy$
 whose $L$-sections are non-trivial everywhere.
\index{$\nbigm(\vecy,L)$, $\nbigm(\vecyhat,L)$}
\item[\mbox{$\nbigm\bigl(\vecy,[L]\bigr)$}:]
Parabolic reduced $L$-Bradlow pairs
of type $\vecy$.
\item[\mbox{$\nbigm\bigl(\vecyhat,[L]\bigr)$}:]
Oriented parabolic reduced $L$-Bradlow pairs
of type $\vecy$.
\index{$\nbigm(\vecy,[L])$, $\nbigm(\vecyhat,[L])$}
\item[\mbox{$\nbigm\bigl(\vecyhat,[\vecL]\bigr)$}:]
Oriented parabolic reduced $\vecL$-Bradlow pairs
of type $\vecy$.
\index{$\nbigm(\vecyhat,[\vecL])$}
\hfill\qed
\end{description}
\end{notation}
The condition $O_m$ determines
the open substack of each moduli stack.
They are denoted by
$\nbigm(m,\vecy)$,
$\nbigm(m,\vecyhat)$,
\index{$\nbigm(m,\vecy)$, $\nbigm(m,\vecyhat)$}
$\nbigm(m,\vecy,L)$,
$\nbigm(m,\vecyhat,L)$
\index{$\nbigm(m,\vecy,L)$, $\nbigm(m,\vecyhat,L)$}
$\nbigm(m,\vecy,[L])$,
$\nbigm(m,\vecyhat,[L])$
\index{$\nbigm(m,\vecy,[L])$, $\nbigm(m,\vecyhat,[L])$}
and 
$\nbigm(m,\vecyhat,[\vecL])$
\index{$\nbigm(m,\vecyhat,[\vecL])$}
respectively.
When the parabolic part of $\vecy$ is trivial,
we often use the notation
$\nbigm(y)$, $\nbigm(\yhat)$,
$\nbigm(m,y)$, $\nbigm(m,\yhat)$, etc.

\subsubsection{
The tautological line bundle and
the relations among some moduli stacks}
\label{subsubsection;06.4.11.75}

Let $\vecy$ be an element of $\Type$,
and let $L$ be a line bundle on $X$.
Let $\nbigehat^u(L)$, $\nbige^u(L)$ and $\nbigehat^u[L]$
denote the universal sheaves over
$\nbigm(\vecyhat,L)\times X$,
$\nbigm(\vecy,L)\times X$
and $\nbigm(\vecyhat,[L])\times X$ 
respectively.
The universal $L$-sections of
$\nbigehat^u(L)$ and $\nbige^u(L)$
are denoted by $\phihat^u$ and $\phi^u$, respectively.
The universal reduced $L$-section
of $\nbigehat^u[L]$ is denoted by $[\phihat^u]$.

We have the $G_m$-action $\rho_1$ on 
$\nbigm(\vecyhat,L)$ given by
$\rho_1(t)\cdot (E,F_{\ast},\phi,\rho):=
 \bigl(E,F_{\ast},t\cdot\phi,\rho\bigr)$.
It is easy to observe that
the quotient stack is isomorphic to
$\nbigm(\vecyhat,[L])$.
Therefore, we can regard
$\nbigm(\vecyhat,L)$ as the $G_m$-torsor
over $\nbigm(\vecyhat,[L])$.
The associated line bundle is denoted by
$\nbigo_{\rel}(-1)$.
We put $\nbigo_{\rel}(1):=\nbigo_{\rel}(-1)^{\lor}$,
and we obtain $\nbigo_{\rel}(n)$ in the obvious manner.
\begin{df}
\index{relative tautological line bundle}
\index{tautological line bundle}
\index{$\nbigo_{\rel}(1)$}
The line bundle $\nbigo_{\rel}(1)$ is called
the relative tautological line bundle
of $\nbigm(\vecyhat,[L])$.
It is also called the tautological line bundle.
\hfill\qed
\end{df}

We can obtain the same line bundles 
by the way of the subsubsection \ref{subsubsection;06.5.14.100}.
Namely, let $\nbigt(n)$ denote the trivial line bundle
on $\nbigm(\vecyhat,[L])$ provided with the $G_m$-action
of the weight $n$.
Let $\pi_1:\nbigm(\vecyhat,L)\lrarr\nbigm(\vecyhat,[L])$
denote the natural projection.
Then we obtain the $G_m$-line bundle $\pi_1^{\ast}\nbigt(n)$.
By taking the descent, 
we obtain the line bundle $\nbigi_n$.
Then we have $\nbigi_n\simeq\nbigo_{\rel}(n)$.

We remark that $\nbigo_{\rel}(-1)$ appears
in the domain of the universal reduced $L$-section $[\phihat^u]$.
Namely, $[\phihat^u]$ is the morphism of
$p^{\ast}_{\nbigm(\vecyhat,[L])}L\otimes
 p_{X}^{\ast}\nbigo_{\rel}(-1)\lrarr \nbigehat^u[L]$.
To see that,
we have only to observe that
$\pi^{\ast}[\phihat^u]:p_{\nbigm(\vecyhat,L)}^{\ast}L
 \otimes\pi_{1,X}^{\ast}\nbigt(-1)
\lrarr \pi_{1,X}^{\ast}\nbigehat^u[L]$
is equivariant with respect to $\rho_1$,
and they give the universal objects 
over $\nbigm(\vecyhat,L)\times X$.

\begin{rem}
For example,
we have the projective space bundle
such as 
$\nbigm(m,\vecyhat,[\nbigo(-m)])
 \lrarr \nbigm(m,\vecyhat)$,
and the restriction of $\nbigo_{\rel}$
to $\nbigm(m,\vecyhat,[\nbigo(-m)])$
is the relative tautological line bundle of this bundle.
\hfill\qed
\end{rem}

\vspace{.1in}

On the other hand,
we have the $G_m$-action $\rho_2$
on $\nbigm(\vecyhat,L)$ given by
$\rho_2(t)\cdot(E,F_{\ast},\phi,\rho):=
 (E,F_{\ast},\phi,t\cdot \rho)$.
It is easy to observe that the quotient stack
is isomorphic to $\nbigm(\vecy,L)$.
Thus, we can regard $\nbigm(\vecyhat,L)$
as the $G_m$-torsor on $\nbigm(\vecy,L)$.
The associated line bundle
is clearly isomorphic to the orientation bundle
$\Or\bigl(\nbige^u(L)\bigr)$.

\vspace{.1in}

Let $r$ be the rank of $\vecy$.
The obvious multiplication of $E$
gives the isomorphism
$(E,t^{-1}\cdot\phi,t^r\rho)\simeq (E,\phi,\rho)$.
To see it,
note the following:
Let $f:E_1\lrarr E_2$ be an isomorphism.
Then the orientation $\rho$ of $E_2$ induces
the orientation of $E_1$ as follows:
\[
\begin{CD}
 \det(E_1)@>{\det(f)}>> \det(E_2)
 @>{\rho}>> \det^{\ast}\poin.
\end{CD}
\]
On the other hand, the $L$-morphism $\phi$ of $E_2$ induces
the $L$-morphism of $E_1$ as follows:
\[
 \begin{CD}
 L@>{\phi}>> E_2 @>{f^{-1}}>> E_1.
 \end{CD}
\]
Therefore, we have $\rho_1(t)=\rho_2(t^r)$.
Hence we have the 
naturally defined morphism
$\kappa:\nbigm(\vecyhat,[L])\lrarr\nbigm(\vecy,L)$,
which is etale and finite of degree $r^{-1}$.
However,
the morphism $\kappa$ does not preserve
the universal object.
We have the relation
$\kappa_X^{\ast}\nbige^u(L)
\simeq
 \nbigehat^u[L]\otimes\nbigo_{\rel}(1)$,
and hence
$\kappa^{\ast}\Or\bigl(\nbige^u(L)\bigr)
\simeq
 \nbigo_{\rel}(r)$.

\subsection{Hilbert Polynomials}
\label{subsection;06.5.17.22}
\subsubsection{Hilbert polynomials of coherent sheaves}

\index{Hilbert polynomial}
Let $\nbigo_X(1)$ be an ample line bundle on
a projective variety $X$.
Let $E$ be a coherent sheaf on $X$.
In this paper,
the non-reduced Hilbert polynomial of $E$
is denoted by $H_E$,
i.e.,
$H_E$ is the unique polynomial of $\rnum$-coefficients
such that $H_E(m)$ is same as the Euler number
$\sum (-1)^i \dim H^i\bigl(X,E(m)\bigr)$.
\index{$H_E$}
In the case $\rank(E)>0$,
the reduced Hilbert polynomial of $E$
is denoted by $P_E$,
i.e.,
$P_E:=H_E/\rank(E)$.
\index{$P_E$}

We also use the notation
$h^0(E)$ to denote $\dim H^0(X,E)$.
\index{$h^0(E)$}

\subsubsection{Hilbert polynomials of parabolic sheaves}

We recall the parabolic Hilbert polynomials
and the parabolic degree, which were introduced
by Maruyama and Yokogawa in \cite{my}.
See \cite{my} for more detail.

Let $E_{\ast}:=(E,F_{\ast},\alpha_{\ast})$ be a parabolic sheaves
of depth $l$.
We put $\epsilon_i:=\alpha_{i+1}-\alpha_{i}$
$(i=1,\ldots,l)$.
Recall that the non-reduced parabolic Hilbert polynomial
$H_{E_{\ast}}$ \index{$H_{E_{\ast}}$}
is defined to be as follows:
\[
 H_{E_{\ast}}(t):=
 H_{E(-D)}(t)+\sum_{i=1}^l\alpha_i \cdot H_{\Gr_i(E)}
=H_E(t)-\sum_{i=1}^l \epsilon_i \cdot H_{\Cok_i(E)}.
\]
The reduced parabolic Hilbert polynomial
$P_{E_{\ast}}(t)$ is defined to be 
$P_{E_{\ast}}(t):=
  H_{E_{\ast}}(t)\big/ \rank(E)$.
\index{$P_{E_{\ast}}$}

Since we have the equality
$H_{E}=H_{E(-D)}+\sum _{i=1}^l H_{\Gr_i(E)}$,
we obtain the following lemma.
\begin{lem}
We have the inequality
$H_{E_{\ast}}(t)\leq H_{E}(t)$
and $P_{E_{\ast}}(t)
\leq P_E(t)$ for any sufficiently large $t$.
\hfill\qed
\end{lem}

The parabolic degree is defined to be as follows:
\[
 \pardeg(E_{\ast})=
\deg(E)+
 (\dim X-1)!\times
 \Bigl(\mbox{the coefficient of $t^{\dim X-1}$ of the polynomial
 $\sum_{i=1}^l\alpha_i \cdot H_{\Gr_i}(t)$}\Bigr)
\]
\index{$\pardeg(E_{\ast})$}
Note that we have the inequality
$\deg(E)\leq \pardeg(E_{\ast})$.

The parabolic slope $\mu(E_{\ast})$ is defined to be
$\pardeg(E_{\ast})/\rank(E)$.
Then we have the inequality
$\mu(E)\leq \mu(E_{\ast})\leq \mu(E)+\deg(D)$
for the usual slope $\mu(E)$ of a torsion-free sheaf $E$.
\index{$\mu(E_{\ast})$}

We also put
$h^0(E_{\ast}):=
 \alpha_1h^0\bigl(E(-D)\bigr)
+\sum\epsilon_i h^0\bigl(F_{i+1}(E)\bigr)$.

\subsubsection{Hilbert polynomial for
  parabolic $L$-Bradlow pairs}

We recall the Hilbert polynomials for Bradlow pairs,
following \cite{th1}.

\begin{notation}
\index{$\nbigp^{\br}$}
Let $\nbigp^{\br}$ denote the set of polynomials $\delta$
of $\real$-coefficients
such that $\deg(\delta)\leq \dim X-1$
and that $\delta(t)>0$ for any sufficiently large $t$

For any element $\delta\in\nbigp^{\br}$,
the coefficient of $t^{d-1}$ in $\delta$ is denoted by $\delta_{\top}$,
which may be $0$.
\index{$\delta_{\top}$}

The total order $\leq$ on the set $\nbigp^{\br}$
is defined as follows:
Let $\delta$ and $\delta'$ be elements of $\nbigp^{\br}$.
Then
$\delta\leq\delta'$ if and only if
$\delta (t)\leq\delta'(t)$ for any sufficiently large $t$.
\hfill\qed
\end{notation}

Let $(E_{\ast},\phi)$ be a parabolic $L$-Bradlow pair
on $(X,D)$.
For any $\delta\in\nbigp^{\br}$,
the non-reduced $\delta$-Hilbert polynomial $H^{\delta}_{(E_{\ast},\phi)}$
of $(E_{\ast},\phi)$ is defined as follows:\index{$H^{\delta}_{(E_{\ast},\phi)}$}
\begin{equation}
 \label{eq;06.4.9.3}
 H^{\delta}_{(E_{\ast},\phi)}:=
 H_{E_{\ast}}+\epsilon(E_{\ast},\phi)\cdot \delta,
\quad
 \epsilon(E_{\ast},\phi):=\left\{
 \begin{array}{ll}
 1 & (\phi\neq 0)\\
 0 & (\phi=0)
 \end{array}
 \right.
\end{equation}
The reduced $\delta$-Hilbert polynomial
is defined to be
$P^{\delta}_{(E_{\ast},\phi)}:=
 H^{\delta}_{(E_{\ast},\phi)}/\rank E$.
\index{$P^{\delta}_{(E_{\ast},\phi)}$}
Similarly, the slope $\mu^{\delta}_{(E_{\ast},\phi)}$
is defined by
$\mu^{\delta} (E_{\ast},\phi):=
 \mu(E_{\ast})+\epsilon(\phi)\cdot \delta_{\top}/\rank(E)$.
\index{$\mu^{\delta}(E_{\ast},\phi)$}

Let $\vecL=(L_1,L_2)$ be a pair of line bundles on $X$.
We take $\vecdelta:=(\delta_1,\delta_2)\in
 \bigl(\nbigp^{\br}\bigr)^2$.
Let $(E_{\ast},\vecphi)$ on $(X,D)$
be a parabolic $\vecL$-Bradlow pair.
The $\vecdelta$-Hilbert polynomial is defined by 
$ H^{\vecdelta}_{(E_{\ast},\vecphi)}:=
 H_{E_{\ast}}+
\sum_{i=1,2} \epsilon(E_{\ast},\phi_i)\cdot\delta_i$,
where $\epsilon(E_{\ast},\phi_i)$ are given as in (\ref{eq;06.4.9.3}).
We also put
$P^{\vecdelta}_{(E_{\ast},\vecphi)}
=H^{\vecdelta}_{(E_{\ast},\vecphi)}/\rank E$.
\index{$P^{\vecdelta}_{(E_{\ast},\vecphi)}$}

\subsubsection{The Hilbert polynomials associated to a type}
\label{subsubsection;06.5.15.1}

Let $z$ be any element of $\bigoplus_{i\geq 1}H^{i}(X)$.
Since we have $z^k=0$  for a large integer $k$,
we have the polynomial
$\exp( t\!\cdot\! z)=\sum_{i=0}^{\infty}(k!)^{-1}(t\cdot z)^k$.
In the case $c=c_1\bigl(\nbigo_X(1)\bigr)$,
it will be denoted by $\ch\bigl(\nbigo_X(t)\bigr)$
in the following.
When we substitute $t=m$ for some integer $m$,
it is same as the ordinary meaning.

Let $\vecy=(y,y_1,y_2,\ldots,)$ be an element of $\Type$.
We put as follows:
\[
 H_{\vecy}(t):=\int_X \Td(X)\cdot\ch\bigl(\nbigo_X(t)\bigr)\cdot y,
\quad
 H_{\vecy,i}(t):=
 \int_X\Td(X)\cdot \ch\bigl(\nbigo_X(t)\bigr)\cdot
 \sum_{j\leq i}y_j.
\]
\index{$H_{\vecy}(t)$, $H_{\vecy,i}(t)$, $H_y(t)$}
When $D$ is smooth and we regard $y_i$ as an element of
$H^{\ast}(D)$,
we put $H_{\vecy,i}(t):=\int_{D}\Td(D)\cdot
   \ch\bigl(\nbigo_X(t)\bigr)\cdot\sum_{j\leq i}y_j$.
When the parabolic part of $\vecy$ is trivial,
we use the notation $H_{y}(t)$ instead of $H_{\vecy}(t)$.

If we are given a system of weights
$\alpha_{\ast}$,
we put $\epsilon_i:=\alpha_{i+1}-\alpha_i$.
And we put as follows: 
\[
 H^{\alpha_{\ast}}_{\vecy}:=
 H_{\vecy}-\sum \epsilon_i\cdot H_{\vecy,i},
\quad
 P^{\alpha_{\ast}}_{\vecy}:=
 \frac{H^{\alpha_{\ast}}_{\vecy}}{\rank \vecy}.
\]
\index{$H^{\alpha_{\ast}}_{\vecy}$,
 $P^{\alpha_{\ast}}_{\vecy}$}
We also put as follows:
\[
  \deg(\vecy,\alpha_{\ast}):=
 \int_{X}y\cdot c_1\bigl(\nbigo_X(1)\bigr)
-\sum\epsilon_i\cdot \int_X \sum_{j\leq i}y_j
\cdot c_1\bigl(\nbigo_X(1)\bigr),
\quad
 \mu(\vecy,\alpha_{\ast}):=
 \frac{\deg(\vecy,\alpha_{\ast})}{\rank\vecy}
\]
\index{$\deg(\vecy,\alpha_{\ast})$, $\mu(\vecy,\alpha_{\ast})$}

If we are given an element $\delta\in\nbigp^{\br}$,
we put as follows:
\[
 H^{\alpha_{\ast}\delta}_{\vecy}:=
 H^{\alpha_{\ast}}+\delta,
\quad
 P^{\alpha_{\ast}\delta}_{\vecy}:=
 \frac{H^{\alpha_{\ast}\delta}_{\vecy}}{\rank \vecy},
\quad
  \mu(\vecy,\alpha_{\ast},\delta):=
 \mu(\vecy,\alpha_{\ast})+
 \frac{\delta_{top}}{\rank(\vecy)}
\]
\index{$H^{\alpha_{\ast},\delta}_{\vecy}$,
 $P^{\alpha_{\ast},\delta}_{\vecy}$,
 $\mu(\vecy,\alpha_{\ast},\delta)$}

\subsection{Semistability}
\label{subsection;06.7.3.23}
\subsubsection{Semistability and the moduli stacks}
\label{subsubsection;06.6.8.2}

Let $(E_{\ast},\phi)$ be a parabolic $L$-Bradlow pair
on $(X,D)$,
and let $\delta$ be an element of $\nbigp^{\br}$.
Recall that $(E_{\ast},\phi)$ is called $\delta$-semistable,
if the following inequality holds
for each sub-objects $(E'_{\ast},\phi')$
of $(E,\phi)$:\index{$\delta$-semistable, $\delta$-stable}
\begin{equation}
 \label{eq;06.4.9.2}
 P^{\delta}_{(E'_{\ast},\phi')}(t)
\leq P^{\delta}_{(E_{\ast},\phi)}(t)
\quad
 \mbox{\rm
 ($t$ is sufficiently large)}
\end{equation}
If the strict inequality holds
in (\ref{eq;06.4.9.2}) for each subobject,
$(E_{\ast},\phi)$ is called $\delta$-stable.

Since parabolic reduced $L$-Bradlow pairs on $X$
are just parabolic $L$-Bradlow pairs whose $L$-section
is non-trivial,
the notion of $\delta$-(semi)stability
of parabolic reduced $L$-Bradlow pairs is also given.
The $\delta$-(semi)stability 
of oriented parabolic (reduced) $L$-Bradlow pairs on $X$
is defined by the $\delta$-(semi)stability
of the underlying parabolic (reduced) $L$-Bradlow pairs.

\begin{rem}
Let $E$ be a torsion-free sheaf on $X$.
We can regard it as the parabolic sheaf $E_{\ast}$
canonically.
The quasi-parabolic structure is given by
$F_{2}=E(-D)\subset F_1=E$.
The weight is given by $\alpha_1=1$.
Then we have $H_{E_{\ast}}=H_E$,
and hence the semistability of $(E_{\ast},\phi)$ is 
equivalent to the semistability of $(E,\phi)$.

If the weight is given by $\alpha_1=0$,
then we have $H_{E_{\ast}}=H_{E(-D)}$,
and hence the semistability 
of $(E_{\ast},\phi)$ and $(E,\phi)$
are not equivalent,
as remarked in Remark {\rm 1.1.1} in {\rm\cite{my}}.
\hfill\qed
\end{rem}

We also have the $\delta$-$\mu$-semistability
in the standard way.\index{$\delta$-$\mu$-semistability}
Namely,
a parabolic $L$-Bradlow pair $(E_{\ast},\phi)$
on $(X,D)$ is called $\delta$-$\mu$-semistable,
if the inequality
$\mu^{\delta}(E'_{\ast},\phi')
\leq
 \mu^{\delta}(E_{\ast},\phi)$
holds for any sub-objects
$(E'_{\ast},\phi')\subset (E_{\ast},\phi)$.
If the strict inequality holds
for any subjects,
$(E_{\ast},\phi)$ is called
$\delta$-$\mu$-stable.
It is easy to check the implication:
\[
 \mbox{\rm $\delta$-$\mu$-stable } 
\Longrightarrow
 \mbox{\rm $\delta$-stable }
\Longrightarrow
 \mbox{\rm $\delta$-semistable }
\Longrightarrow
 \mbox{\rm $\delta$-$\mu$-semistable}
\]
Similarly,
we have the notion of 
$\delta$-$\mu$-semistability
and $\delta$-$\mu$-stability
of parabolic reduced $L$-Bradlow pairs.

\vspace{.1in}

In the family case,
a $U$-parabolic $L$-Bradlow pair 
$(E_{\ast},\phi)$ on $U\times X$
is called $\delta$-(semi)stable,
if $(E_{\ast},\phi)_{|\{u\}\times X}$ is 
$\delta$-(semi)stable for each $u\in U$.
We obtain the $\delta$-(semi)stability
in the oriented case and in the reduced case
similarly.

\begin{rem}
As usual,
we have only to consider sub-objects
$(E'_{\ast},\phi')\subset (E_{\ast},\phi)$
such that $E'$ is saturated,
when we check the $\delta$-(semi)stability
of $(E_{\ast},\phi)$.
\hfill\qed
\end{rem}

Let $L_i$ $(i=1,2)$ be line bundles on $X$,
and let $\iota:L_1\lrarr L_2$ be a non-trivial morphism.
Let $(E_{\ast},\phi)$ be a parabolic $L_2$-Bradlow pair on $(X,D)$.
Then we obtain the parabolic $L_1$-Bradlow pair $(E,\phi\circ\iota)$.
When we consider the semistable condition for 
a parabolic $L$-Bradlow pair,
the choice of $L$ is not essential in the following sense.
\begin{lem}
A pair $(E_{\ast},\phi)$ is $\delta$-semistable,
if and only if $(E_{\ast},\phi\circ\iota)$ is $\delta$-semistable.
\end{lem}
\pf
Let $F$ be a saturated subsheaf of $E$.
The image $\phi\circ\iota(L_1)$ is contained in $F$
if and only if the image $\phi(L_2)$ is contained in $F$.
Thus the claim is clear.
\hfill\qed

\begin{lem}
\label{lem;06.4.24.6}
Let $(E_{i\,\ast},\phi_i)$ $(i=1,2)$
be $\delta$-semistable parabolic $L$-Bradlow pairs with
$P^{\delta}_{(E_{1\,\ast},\phi_1)}(t)
= P^{\delta}_{(E_{2\,\ast},\phi_2)}(t)$.
Let $f:(E_{1\,\ast},\phi_1)\lrarr (E_{2\,\ast},\phi_2)$
be a non-trivial morphism.
We have the induced $L$-Bradlow pairs
$(\Ker(f)_{\ast},\phi')$,
$(\Image(f)_{\ast},\phi'')$
and $(\Cok(f)_{\ast},\phi''')$.
Then they are also $\delta$-semistable.

A similar claim holds for $\delta$-$\mu$-semistability.
\end{lem}
\pf
We put 
$(E_{3\,\ast},\phi_3):=
\bigl(\Image (f)_{\ast},\phi''\bigr)$.
From the $\delta$-semistability of
$(E_{i\,\ast},\phi_i)$ $(i=1,2)$,
we obtain the inequality
for sufficiently large $t$:
\begin{equation}
 \label{eq;06.4.24.5}
 P^{\delta}_{(E_{1\,\ast},\phi_1)}(t)
\leq
 P^{\delta}_{(E_{3\,\ast},\phi_3)}(t)
\leq
 P^{\delta}_{(E_{2\,\ast},\phi_2)}(t).
\end{equation}
Due to the assumption
$P^{\delta}_{(E_{1\,\ast},\phi_1)}(t)
= P^{\delta}_{(E_{2\,\ast},\phi_2)}(t)$,
the equalities hold in (\ref{eq;06.4.24.5}).
Then, it is easy to derive the claims of the lemma
by definition of semistability.
\hfill\qed

\begin{cor}
\label{cor;06.5.21.35}
Any automorphisms of stable objects 
are constant multiplication.
\hfill\qed
\end{cor}

We introduce the following notation.
In each line,
the left hand side denotes the moduli stack 
of the object in the right hand side:

\begin{description}
\item[$\nbigm^s(\vecy,\alpha_{\ast})$:]
stable parabolic sheaves of type $\vecy$
with weight $\alpha_{\ast}$.
\item[$\nbigm^s(\vecyhat,\alpha_{\ast})$:]
stable oriented parabolic sheaves of type $\vecy$
with weight $\alpha_{\ast}$.
  \index{$\nbigm^s(\vecy,\alpha_{\ast})$,
 $\nbigm^s(\vecyhat,\alpha_{\ast})$}
\item[$\nbigm^s(\vecy,L,\alpha_{\ast},\delta)$:]
$\delta$-stable $L$-Bradlow pairs of type $\vecy$
with weight $\alpha_{\ast}$,
whose $L$-sections are non-trivial everywhere.
\item[$\nbigm^s(\vecyhat,L,\alpha_{\ast},\delta)$:]
$\delta$-stable oriented $L$-Bradlow pairs of type $\vecy$
with weight $\alpha_{\ast}$,
whose $L$-sections are non-trivial everywhere.
 \index{$\nbigm^s(\vecy,L,\alpha_{\ast}\delta)$,
 $\nbigm^s(\vecyhat,L,\alpha_{\ast},\delta)$}
\item[\mbox{$\nbigm^s(\vecy,[L],\alpha_{\ast},\delta)$}:]
$\delta$-stable reduced $L$-Bradlow pairs of type $\vecy$
with weight $\alpha_{\ast}$.
\item[\mbox{$\nbigm^s(\vecyhat,[L],\alpha_{\ast},\delta)$}:]
$\delta$-stable oriented reduced $L$-Bradlow pairs of type $\vecy$
with weight $\alpha_{\ast}$.
 \index{$\nbigm^s(\vecy,[L],\alpha_{\ast},\delta)$,
 $\nbigm^s(\vecyhat,[L],\alpha_{\ast},\delta)$}
\end{description}
For the moduli stack of semistable objects,
we use the notation
$\nbigm^{ss}(\vecy,\alpha_{\ast})$,
$\nbigm^{ss}(\vecyhat,\alpha_{\ast})$,
$\nbigm^{ss}(\vecy,L,\alpha_{\ast},\delta)$,
$\nbigm^{ss}(\vecyhat,[L],\alpha_{\ast},\delta)$, etc..
\index{
 $\nbigm^{ss}(\vecy,\alpha_{\ast})$,
$\nbigm^{ss}(\vecyhat,\alpha_{\ast})$}
\index{$\nbigm^{ss}(\vecy,L,\alpha_{\ast},\delta)$,
$\nbigm^{ss}(\vecyhat,[L],\alpha_{\ast},\delta)$}
When the parabolic structure of $\vecy$ is trivial,
we often use the notation
$\nbigm^{s}(\yhat)$,
$\nbigm^{s}(\yhat,[L],\delta)$,
$\nbigm^{s}(y,L,\delta)$, etc..
\index{
 $\nbigm^{s}(\yhat)$,
$\nbigm^{s}(\yhat,[L],\delta)$,
$\nbigm^{s}(y,L,\delta)$}

Similarly, we also have the notion of
$\vecdelta$-(semi)stability
for (oriented, reduced) $\vecL$-Bradlow pairs.
We also have
$\vecdelta$-$\mu$-(semi)stability.
We denote by $\nbigm^{s}(\vecyhat,[\vecL],\alpha_{\ast},\vecdelta)$
the moduli stack of $(\delta_1,\delta_2)$-stable
oriented reduced $\vecL$-Bradlow pairs
of type $\vecy$ with weight $\alpha_{\ast}$,
for example.

When any $E_{\ast}\in\nbigm^{s}(\vecyhat,\alpha_{\ast})$
satisfies the condition $O_m$,
we denote by $\nbigmtilde^{s}_m(\vecyhat,\alpha_{\ast})$
the full flag bundle associated
to the vector bundle $p_{X\,\ast}\Ehat^u(m)$,
where $\Ehat^u$ denotes the universal bundle
over $\nbigm^s_m(\vecyhat,\alpha_{\ast})\times X$.
We use the notation
$\nbigmtilde^s_m(\vecyhat,[L],\alpha_{\ast},\delta)$
and $\nbigmtilde^{s}_m(\vecyhat,[\vecL],\alpha_{\ast},\vecdelta)$,
etc. in similar ways.
If there are no risk of confusion, we will often omit to denote $m$.
\index{
 $\nbigmtilde^{s}_m(\vecyhat,\alpha_{\ast})$,
 $\nbigmtilde^{s}(\vecyhat,\alpha_{\ast})$}
\index{
 $\nbigmtilde^s_m(\vecyhat,[L],\alpha_{\ast},\delta)$,
 $\nbigmtilde^s(\vecyhat,[L],\alpha_{\ast},\delta)$}
\index{
$\nbigmtilde^{s}_m(\vecyhat,[\vecL],\alpha_{\ast},\vecdelta)$,
$\nbigmtilde^{s}(\vecyhat,[\vecL],\alpha_{\ast},\vecdelta)$}

\subsubsection{Harder-Narasimhan filtration
  and partial Jordan-H\"{o}lder filtration}
\index{Harder-Narasimhan filtration}
\index{partial Jordan-H\"{o}lder filtration}
\label{subsubsection;06.7.3.231}

Let $(E_{\ast},\phi)$ be a $\delta$-semistable
parabolic $L$-Bradlow pair.
In this paper, a filtration
\[
 (E_{1\,\ast},\phi_1)
\subset
 (E_{2\,\ast},\phi_2)
\subset\cdots\subset
 (E_{k\,\ast},\phi_k)
=(E_{\ast},\phi)
\]
is called a partial Jordan-H\"{o}lder filtration
with respect to $\delta$-semistability,
if each $(E_{i\,\ast},\phi_i)$
is $\delta$-semistable
such that
$P^{\delta}_{(E_{i\,\ast},\phi_i)}
=P^{\delta}_{(E_{\ast},\phi)}$.
Each $\Gr_i(E):=E_i/E_{i-1}$
has the induced parabolic structure
and the $L$-section $\Gr_i(\phi)$,
and the parabolic $L$-Bradlow pair
$\bigl(\Gr_{i\,\ast}(E),\Gr_i(\phi)\bigr)$
is $\delta$-semistable with
$P^{\delta}_{(\Gr_{i\,\ast}(E),\Gr_i(\phi))}=
 P^{\delta}_{(E_{\ast},\phi)}$.

If each $\bigl(\Gr_{i\,\ast}(E),\Gr_i(\phi)\bigr)$
is $\delta$-stable,
the filtration is called a Jordan-H\"older filtration
with respect to $\delta$-stability.
It can be shown that
the length of Jordan-H\"older filtration
and the collection of graded objects
$\bigl\{
\bigl(\Gr_{i\,\ast}(E),\Gr_i(\phi)\bigr)
\bigr\}$
are independent of the choice of 
Jordan-H\"older filtration.

Similarly,
we have the notion of
partial Jordan-H\"{o}lder filtration
with respect to $\delta$-$\mu$-semistability
and Jordan-H\"older filtration.

\begin{lem}
Let $(E_{\ast},\phi)$ be a parabolic $L$-Bradlow pair.
There exists the unique increasing filtration
$F=\bigl(F_i(E)\,|\,i=1,2,\ldots \bigr)$ of $E$
with the following properties:
\begin{itemize}
\item
The induced objects
$\bigl(
 \Gr^F_i(E),\phi_i'
 \bigr)$ are $\delta$-semistable.
\item
The inequalities
$P^{\delta}_{( \Gr^{F}_i(E) ,\phi_i' )}(t)>
 P^{\delta}_{(\Gr^{F}_{i+1}(E),\phi_{i+1}')}(t)$ hold,
for any sufficiently large $t$.
\end{itemize}
The filtration is called 
the Harder-Narasimhan filtration
with respect to the $\delta$-semistability.

We also have the unique Harder-Narasimhan filtration
with respect to the $\delta$-$\mu$-semistability.
\end{lem}
\pf
We give only an outline.
We use an induction on $\rank (E)$.
In the case $\rank(E)=1$,
the claim is trivial.
Take a sufficiently negative number $C$.
We know that
the family of saturated subsheaves $E'$ of $E$
such that
$\deg(E')\geq C$ is bounded
(Proposition \ref{prop;06.5.21.1}).
Therefore,
the family of saturated subsheaves $E'$
of $E$ such that
$P^{\delta}_{(E'_{\ast},\phi')}
\geq P^{\delta}_{(E_{\ast},\phi)}$ is bounded,
where $(E'_{\ast},\phi')$ denotes
the subobject of $(E_{\ast},\phi)$.
Hence we obtain the finiteness
of the set $S$ of the polynomials $P$
with the following properties:
\begin{itemize}
\item
$P(t)\geq P_{(E_{\ast},\phi)}(t)$
for any sufficiently large $t$.
\item
There exist a saturated subsheaf $E'$ of $E$
such that $P^{\delta}_{(E'_{\ast},\phi')}=P$.
\end{itemize}
We say $P\leq 'Q$ if $P(t)\leq Q(t)$ for any
sufficiently large $t$.
It gives the total order of $S$.
Let $P_0$ be the maximum.
Let $T(P_0)$ denote the family of saturated subsheaves
$E'$ of $E$
such that $P^{\delta}_{(E'_{\ast},\phi)}=P_0$.
Then it is easy to see that
$(E'_{\ast},\phi)$ is $\delta$-semistable
for any $E'\in T(P_0)$.
It is also easy to see
that $E'_1+E_2'\in T(P_0)$
for $E_i'\in T(P_0)$.
Therefore, we have the maximum $E_1$
in $T(P_0)$ with respect to the inclusion.

We put $\Ehat:=E/E_1$,
and let $\pi:E\lrarr \Ehat$ denote the natural projection.
Let $\Ehat'\subset \Ehat$ be any saturated subsheaf.
If 
$P^{\delta}_{(\Ehat',\phihat')}
\geq P^{\delta}_{(E_1,\phi_1)}$,
then we have
$P^{\delta}_{(\pi^{-1}(\Ehat'),\phi'')}
\geq P^{\delta}_{(E_1,\phi_1)}$,
which contradicts our choice of $E_1$.
Therefore, we have 
$P^{\delta}_{(\Ehat',\phihat')}
<P^{\delta}_{(E_1,\phi_1)}$.
By applying the hypothesis of the induction,
we have the Harder-Narasimhan filtration
of $(\Ehat_{\ast},\widehat{\phi})$
with respect to $\delta$-semistability.
Together with the above remark,
we obtain the Harder-Narasimhan filtration
of $(E_{\ast},\phi)$
with respect to $\delta$-semistability.
\hfill\qed

\subsubsection{$(\delta,\ell)$-Semistability}
\label{subsubsection;06.4.25.10}

Let $P$ be a polynomial, and let $r$ be a positive integer.
Let $m$ be a sufficiently large integer
satisfying the following condition:
\begin{itemize}
\item
 Let $(E_{\ast},\phi)$ be a $\delta$-parabolic $L$-Bradlow pair
 with $P^{\delta}_{(E_{\ast},\phi)}=P$
 and $\rank E\leq r$.
 Then $E_{\ast}$ satisfies the condition $O_m$.
\end{itemize}

\begin{df}
\index{$(\delta,\ell)$-semistability}
Let $(E_{\ast},[\phi])$ be
a parabolic reduced $L$-Bradlow pair on $X$
with $P^{\delta}_{(E_{\ast},\phi)}=P$,
and let $\nbigf$ be a full flag of $H^0\bigl(X,E(m)\bigr)$.
Let $\ell$ be any positive integer.
We say that $(E_{\ast},[\phi],\nbigf)$ is $(\delta,\ell)$-semistable,
if the following condition holds:
\begin{itemize}
\item
 $(E_{\ast},[\phi])$ is $\delta$-semistable.
\item
 Take any partial Jordan-H\"older filtration 
of $(E_{\ast},[\phi])$ with respect to $\delta$-semistability:
\[
 E^{(1)}_{\ast}\subset
 E^{(2)}_{\ast}\subset\cdots
 \subset E^{(i-1)}_{\ast}
 \subset (E^{(i)}_{\ast},\phi)
\subset\cdots\subset (E^{(k)}_{\ast},\phi)
\]
Then we have
$\nbigf_{\ell}\cap H^0\bigl(X,E^{(i-1)}(m)\bigr)=\{0\}$
and 
 $\nbigf_{\ell}\not\subset H^0\bigl(X,E^{(j)}(m)\bigr)$
 for $j<k$.
\end{itemize}

We denote by
$\nbigmtilde_m^{ss}\bigl(\vecy,[L],\alpha_{\ast},(\delta,\ell)\bigr)$,
the moduli stack of such tuples
$(E_{\ast},[\phi],\nbigf)$.
In the oriented case,
we use the notation
$\nbigmtilde_m^{ss}\bigl(\vecyhat,[L],\alpha_{\ast},(\delta,\ell)\bigr)$
as usual.
\index{
 $\nbigmtilde_m^{ss}(\vecy,[L],\alpha_{\ast},(\delta,\ell))$,
 $\nbigmtilde^{ss}(\vecy,[L],\alpha_{\ast},(\delta,\ell))$}
\index{
 $\nbigmtilde_m^{ss}(\vecyhat,[L],\alpha_{\ast},(\delta,\ell))$,
 $\nbigmtilde^{ss}(\vecyhat,[L],\alpha_{\ast},(\delta,\ell))$}

Similarly,
we have the notion of $(\delta,\ell)$-semistability 
for a $L$-Bradlow pair $(E_{\ast},\phi)$ such that
$\phi\neq 0$
and a full flag $\nbigf$ of $H^0(X,E(m))$.
The moduli stack is denoted by
$\nbigmtilde_m^{ss}\bigl(\vecy,L,\alpha_{\ast},(\delta,\ell)\bigr)$.
\index{$\nbigmtilde_m^{ss}\bigl(\vecy,L,\alpha_{\ast},(\delta,\ell)\bigr)$,
$\nbigmtilde^{ss}\bigl(\vecy,L,\alpha_{\ast},(\delta,\ell)\bigr)$}
When there are no risk of confusion,
we omit to denote $m$.
\hfill\qed
\end{df}

\begin{rem}
\label{rem;06.5.21.60}
When $\ell$ is sufficiently large,
the second condition is trivial.
The first condition is equivalent to
$i=1$ for any J-H filtration.
Therefore, $(q,E_{\ast},\phi,\nbigf)$ is
$(\delta,\ell)$-semistable
if and only if
$(E_{\ast},\phi)$ is $\delta_-$-semistable
for a parameter $\delta_-<\delta$
such that $|\delta-\delta_-|$ is sufficiently small.
\hfill\qed
\end{rem}

\begin{lem}
\label{lem;06.5.21.60}
Let $(E_{\ast},[\phi],\nbigf)$ be a reduced $L$-Bradlow pair
with a full flag $\nbigf$ of $ H^0\bigl(X,E(m)\bigr)$.
We assume that it is $(\delta,\ell)$-semistable.
Then the automorphism group of
$(E_{\ast},[\phi],\nbigf)$ is $G_m$.
\end{lem}
\pf
Let $f$ be an endomorphism 
of $(E_{\ast},[\phi],\nbigf)$.
We have the generalized eigen decomposition
$(E_{\ast},[\phi],\nbigf)=
 (E_{1\,\ast},[\phi_1],\nbigf^{(1)})
\oplus
 \bigoplus_{i=2}^k(E_{i\,\ast},\nbigf^{(i)})$.
But, the condition of $(\delta,\ell)$ is not satisfied
if $k\geq 2$.
Therefore, $f$ has the unique eigenvalue.
Let $N$ be the nilpotent part of $f$.
Let $h$ be the integer such that
$N^h\neq 0$ and $N^{h+1}=0$.
Assume $h\geq 0$, and we will derive the contradiction.
We put $E_1:=\Image N^h$
and $E_2:=\Ker N^h$.
We have the naturally induced parabolic 
structures and the $L$-sections 
$\phi_i$ of $E_{i}$.
Since we have $N(\phi)=0$, $\phi_1=0$.
Then, we obtain the partial Jordan-H\"older filtration
$E_{1\,\ast}
\subset
 (E_{2\,\ast},\phi_2)
\subset
 (E_{\ast},\phi)$
with respect to $\delta$-semistability,
due to Lemma \ref{lem;06.4.24.6}.

We have the induced filtrations
$\nbigf^{(1)}$ on $H^0\bigl(X,E_1(m)\bigr)$
and $\nbigf^{(3)}$ on $H^0\bigl(X,E/E_1(m)\bigr)$.
We have the induced isomorphism
$H^0\bigl(X, E/E_1(m) \bigr)
\lrarr
 H^0\bigl(X,E_1(m)\bigr)$,
which has to preserve the filtration
$\nbigf^{(i)}$ above.
However, we have 
$\nbigf_{\ell}^{(1)}=0$
and 
$\nbigf_{\ell}^{(3)}\neq 0$
due to the $(\delta,\ell)$-semistability.
Thus we arrive at the contradiction,
and hence we have $h=0$.
\hfill\qed

\subsection{Some Boundedness}
\label{subsection;06.7.3.24}
\subsubsection{Foundational theorems}
\label{subsubsection;06.7.3.241}

Let $\nbigo_X(1)$ be a very ample line bundle.
Let $D$ denote a Cartier divisor of $X$.
We recall several foundational theorems,
by following D. Huybrecht-M. Lehn (\cite{hl2}).

\begin{prop}
[A. Grothendieck, Lemma 2.5, \cite{gro}] 
\label{prop;06.5.21.1}
Let $\nbigd$ be a bounded family of coherent sheaf $F$
on $X$.
Let $C$ be any positive number.
Then we have the boundedness
of the family of torsion-free coherent sheaves $F''$
with the following property:
\begin{itemize}
\item
$\deg(F'')\leq C$.
\item
There exists a member $F\in\nbigd$
such that $F''$ is a quotient sheaf of $F$.
\hfill\qed
\end{itemize}
\end{prop}

For any torsion-free coherent sheaf $E$ on $X$,
we have the Harder-Narasimhan filtration
with respect to the standard semistability.
We denote the slope of the first (resp. last)
term of the Harder-Narasimhan filtration
by  $\mu_{\max}(E)$ (resp. $\mu_{\min}(E)$).
\index{$\mu_{\max}(E)$, $\mu_{\min}(E)$}

\begin{prop}
[M. Maruyama \cite{mar}] \label{prop;06.5.21.2}
Let $H$ be a polynomial, and let $C$ be a constant.
We have the boundedness
of the family of torsion-free coherent sheaves
$F$ on $X$ satisfying $\mu_{\max}(F)\leq C$
and $H_F=H$.
\hfill\qed
\end{prop}

We use the notation $[x]_+=\max\{x,0\}$
for any real number $x$.
\begin{prop}
[C. Simpson \cite{si}] \label{prop;06.5.21.3}
Let $r$ be a positive integer.
Then there is a positive constant $c$ such that
the following inequality holds 
for every $\mu$-semistable sheaf $F$
satisfying $\rank(F)<r$ and $\mu(F)<\mu$:
\[
 \frac{h^0(F)}{\rank(F)}\leq
 \frac{1}{g^{d-1}d!}([\mu+c]_+^d).
\]
Here $d$ denotes the dimension of $X$,
and $g$ denotes the number
$c_1\bigl(\nbigo_X(1)\bigr)^d\cap [X]$.
\hfill\qed
\end{prop}

\subsubsection{Boundedness of semistable $L$-Bradlow pairs}
\label{subsubsection;06.7.3.242}

Let $\vecy$ be an element of $\Type$,
and let $\alpha_{\ast}$ be a system of weights.
Let $L$ be a line bundle over $X$,
and let $\delta^{(0)}$ be an element of $\nbigp^{\br}$.
Let $\nbigss(\vecy,L,\alpha_{\ast},\delta^{(0)})$
denote the family of parabolic $L$-Bradlow pairs
$(E_{\ast},\phi)$ such that $\phi\neq 0$
with the following properties:
\begin{itemize}
\item
The type of $E_{\ast}$ is  $\vecy$,
and the weight of the parabolic structure is
given by  $\alpha_{\ast}$.
\item
$(E_{\ast},\phi)$ is 
 $\delta$-$\mu$-semistable
for some $\delta\leq \delta^{(0)}$ in $\nbigp^{\br}$.
\end{itemize}

\begin{lem}
 \label{lem;06.4.11.2}
The family $\nbigss(\vecy,L,\alpha_{\ast},\delta^{(0)})$ is bounded.
\end{lem}
\pf
Let $(E_{\ast},\phi)$ be a member of
$\nbigss(\vecy,L,\alpha_{\ast},\delta^{(0)})$.
Assume that it is $\delta$-$\mu$-semistable
for $\delta\in\nbigp^{\br}$.
Let $E'$ be the first member of
the Harder-Narasimhan filtration of $E$
with respect to the standard semistability.
Then we have the following inequalities:
\[
  \mu_{\max}(E)=\mu(E')\leq
 \mu(E'_{\ast})+
\frac{ \epsilon(E',\phi)\cdot\delta_{top}}{\rank E'}
\leq
 \mu(E_{\ast})+\frac{\epsilon(E,\phi)\cdot\delta_{top}}{\rank E}
\leq \mu(E_{\ast})+\frac{\delta_{top}^{(0)}}{\rank E}.
\]
The last term depends only on $(\vecy,\delta^{(0)})$.
Thus we obtain the boundedness from
Proposition \ref{prop;06.5.21.2}.
\hfill\qed

\vspace{.1in}

Recall the following important observation
due to Thaddeus \cite{th1}.
\begin{prop}
 \label{prop;06.4.11.1}
Take an element $\vecy\in\Type$
and a line bundle $L$ on $X$.
Assume $r=\rank(\vecy)>1$.
Let $\delta$ be an element of $\nbigp^{\br}$
satisfying the following condition:
\begin{equation}
 \label{eq;06.5.21.5}
 \delta_{top}>\frac{r}{r-1}
 \bigl(\mu(\vecy,\alpha_{\ast})-\deg(L)\bigr).
\end{equation}
Then, there does not exist $\delta$-semistable
parabolic $L$-Bradlow pair
$(E_{\ast},\phi)$  of type $\vecy$
such that $\phi\neq 0$.
\end{prop}
\pf
Let $(E_{\ast},\phi)$ be a $\delta$-semistable
$L$-Bradlow pair such that $\phi\neq 0$.
The $L$-section $\phi$ generates
the subsheaf $E'$ of $E$ with $\rank E'=1$.
Due to the stability,
we have the following inequality:
\[
 \mu(E'_{\ast})+\delta_{top}
\leq
 \mu(E_{\ast})+\frac{\delta_{top}}{r}.
\]
We also have
$\mu(E'_{\ast})\geq \mu(E')=\deg(E')\geq \deg(L)$.
Therefore,
we obtain the following:
\[
 \left(
 1-\frac{1}{r}
 \right)\cdot \delta_{top}
\leq
 \mu(E_{\ast})-\deg(L).
\]
Thus we are done.
\hfill\qed

\begin{cor}
The family
$\nbigss\bigl(\vecy,L,\alpha_{\ast}\bigr):=
\bigcup_{\delta\in \nbigp^{\br}}
 \nbigss\bigl(\vecy,L,\alpha_{\ast},\delta\bigr)$
is bounded.
\end{cor}
\pf
It follows from Proposition \ref{prop;06.4.11.1}
and Lemma \ref{lem;06.4.11.2}.
\hfill\qed

\vspace{.1in}

Let $\vecL=(L_1,L_2)$ be a pair of line bundles on $X$.
Let $\vecdelta^{(0)}=(\delta_1^{(0)},\delta^{(0)}_2)$
be an element of $\nbigp^{\br}$.
Let $\nbigss(\vecy,\vecL,\alpha_{\ast},\vecdelta^{(0)})$
denote the family of parabolic $\vecL$-Bradlow pairs
$(E_{\ast},\vecphi)$ of type $\vecy$ with weight $\alpha_{\ast}$
such that $\phi_i\neq 0$,
which is $(\delta_1,\delta_2)$-$\mu$-semistable
for some $\delta_i\leq \delta_i^{(0)}$.
By an argument similar to the proof of
Lemma \ref{lem;06.4.11.2},
we can show the following:
\begin{lem}
The family $\nbigss(\vecL,\alpha_{\ast},\vecdelta)$
is bounded.
\hfill\qed
\end{lem}

By an argument similar to the proof of
Proposition \ref{prop;06.4.11.1},
we can show the following lemma.
\begin{lem}
\label{lem;06.5.21.10}
There exists $\delta^{(0)}_2$
such that the following holds
for any $\delta_2\geq \delta_2^{(0)}$
and for any $\delta_1$:
\begin{itemize}
\item
There does not exist any
$(\delta_1,\delta_2)$-semistable
$\vecL$-Bradlow pair $(E_{\ast},\phi_1,\phi_2)$
of type $\vecy$ such that $\phi_2\neq 0$.
\end{itemize}

In particular, we obtain the boundedness of
the family
$\nbigss(\vecL,\alpha_{\ast},\delta_1):=
 \bigcup_{\delta_2}
 \nbigss\bigl(\vecL,\alpha_{\ast},(\delta_1,\delta_2)\bigr)$.
\hfill\qed
\end{lem}

\subsubsection{Boundedness of Yokogawa family}
\label{subsubsection;06.7.3.243}

Let $\vecy$ be an element of $\Type$,
and let $L$ be a line bundle on $X$.
Let us fix a system of weights $\alpha_{\ast}$,
and we put $\epsilon_i:=\alpha_{i+1}-\alpha_i$.
For each positive integer $m$,
let us take an $H_{\vecy}(m)$-dimensional
vector space $V_m$.
Let $\proj_m$ denote the projectivization of $V_m$,
i.e.,
$\proj_m:=\proj(V_m^{\lor})$.
We also take an inclusion
$\iota:\nbigo(-m)\lrarr L$.

\begin{df}
\label{df;06.6.21.11}
A Yokogawa datum of type $(\vecy,m)$
is defined to be a tuple
$\bigl(q,E_{\ast},\phi,W_{\ast},[\phitilde]\bigr)$
as follows:
\begin{itemize}
\item
 $(E_{\ast},\phi)$ is a parabolic $L$-Bradlow pair
 over $X$,  such that  $\type(E_{\ast})=\vecy$.
\item
 $q$ is a generically surjective morphism 
 $V_{m,X}\lrarr E(m)$,
 where we put $V_{m,X}:=V_m\otimes\nbigo_X$.
\item
 $W_{\ast}=(W_1,\ldots,W_{l})$ is a tuple of
 subspaces of $V_m$ such that
 $\dim W_i=H_{\vecy}(m)-H_{\vecy,i}(m)$.
 We assume that
 $H^0(q)(W_i)$ is contained in
$H^0\bigl(X,F_{i+1}(E)(m)\bigr)$.
\item
$[\widetilde{\phi}]$ is a point of $\proj_m$,
and there exists a non-zero scalar $c$
such that
 $H^0(q)(\widetilde{\phi})=c\cdot\iota(\phi)$,
where $\iota(\phi)$ denotes the section
of $E(m)$ induced by $\iota$ and $\phi$.
\hfill\qed
\end{itemize}
\end{df}

\begin{df}
\label{df;06.5.15.20}
\index{$\YOKtilde(m,K,\vecy,L,\delta)$}
Let $\delta$ be an element of $\nbigp^{\br}$.
Let $K$ be any non-negative number.
Let $\YOKtilde(m,K,\vecy,L,\delta)$
denotes the set of Yokogawa data
$\bigl(q,E_{\ast},\phi,W_{\ast},[\phitilde]\bigr)$
such that the following inequality holds
for any subspace $W\subset V_m$:
\begin{equation}
 \label{eq;06.4.19.5}
 P^{\delta,\alpha_{\ast}}_{\vecy}(m)\cdot
 \rank \nbige_W
-\epsilon(W,[\phitilde])\cdot\delta(m)
-\sum_{i=1}^l\epsilon_i\cdot \dim(W_i\cap W)
-\alpha_1\cdot \dim(W)+K\geq 0.
\end{equation}
Here, $\nbige_W$ denotes the subsheaf of $E(m)$
generated by $W$ and $q$,
and 
$\epsilon(W,[\phitilde])$ is defined to be
$1$ $([\phitilde]\in\proj_W)$ or
$0$ $([\phitilde]\not\in\proj_W)$.
We remark
$\epsilon(W,[\phitilde])=
 \epsilon(\nbige_W(-m),\phi)$,
where $\epsilon(\nbige_W(-m),\phi)$ is given as
in {\rm(\ref{eq;06.4.9.3})}.
\hfill\qed
\end{df}

For each positive integer $N$,
we put $\nYOKtilde(N,K,\vecy,L,\delta)
:=\bigcup_{m\geq N}\YOKtilde(m,K,\vecy,L,\delta)$.
\index{$\nYOKtilde(N,K,\vecy,L,\delta)$}
Following Yokogawa \cite{y},
we consider the family $\nYOK(N,K,\vecy,L,\delta)$
\index{$\nYOK(N,K,\vecy,L,\delta)$}
of quasi-parabolic $L$-Bradlow pairs
$(E,F_{\ast},\phi)$ of type $\vecy$
such that there exists a lift
$\bigl(q,E,F_{\ast},\phi,W_{\ast},[\widetilde{\phi}]\bigr)
\in \nYOKtilde(N,K,\vecy,L,\delta)$.
The family will be called the Yokogawa family.

\begin{prop}
\label{prop;06.5.15.10}
There exists a small positive number 
$K_0=K_0(\vecy,L,\delta)$
and a large  integer $N_0=N_0(\vecy,L,\delta)$
such that the following holds:
\begin{itemize}
\item
The family $\nYOK\bigl(N_0,K_0,\vecy,L,\delta\bigr)$
is bounded,
and it satisfies the condition $O_m$
for any $m\geq N_0$.
\item
For any $(q,E_{\ast},\phi,W_{\ast},[\phitilde])
\in \nYOK\bigl(N_0,K_0,\vecy,L,\delta\bigr)$,
the morphism $q$ is surjective.
In particular,
$H^0(q):V_m\lrarr H^0(X,E(m))$ is isomorphic.
\item
All members of $\nYOK(N_0,K_0,\vecy,L,\delta)$ are
$\delta$-semistable.
\end{itemize}
\end{prop}
\pf
We follow the argument of Yokogawa \cite{y}.
We begin with the following lemma.

\begin{lem}
Let $\bigl(E_{\ast},\phi\bigr)$ be 
a parabolic $L$-Bradlow pair
contained in the family $\nYOK(N,K,\vecy,L,\delta)$
with a lift
$\bigl(q,E_{\ast},\phi,W_{\ast},[\widetilde{\phi}]\bigr)$.
Let $E'$ denote a quotient sheaf of $E$.
Then the following inequality holds:
\begin{equation}
 \label{eq;06.4.19.15}
 P_{\vecy}^{\delta,\alpha_{\ast}}(m)\leq
 \frac{h^0\bigl(E'_{\ast}(m)\bigr)
     +\epsilon(E',\phi)\cdot\delta(m) +K }
  {\rank(E')}.
\end{equation}
\end{lem}
\pf
Let $W$ denote the kernel
of the composite of the morphisms:
\[
\begin{CD}
 V_m@>{H^0(q)}>>
H^0\bigl(X,E(m)\bigr) @>>> H^0\bigl(X,E'(m)\bigr)
\end{CD}
\]
Let $\nbige_W$ denote the subsheaf of $E(m)$
generated by the image of $W$ via $q$.
Due to (\ref{eq;06.4.19.5}),
we obtain the following inequality:
\[
 P_{\vecy}^{\delta,\alpha_{\ast}}(m)\cdot\rank(\nbige_W)
+\epsilon(W,\widetilde{\phi})\cdot\delta(m)
-\sum_i\epsilon_i\dim(W_i\cap W)-\alpha_1\dim(W)+K
\geq 0.
\]
We have the inequalities
$\dim(W_i\cap W)\geq \dim(W_i)-h^0\bigl(F_{i+1}(E')(m)\bigr)$
for $i=1,\ldots,l$.
We also have the equalities
$\epsilon(W,[\widetilde{\phi}])+\epsilon(E',[\phi'])=1$,
where $\phi'$ is the induced $L$-section of $E'$ by $\phi$.
Since $q$ is generically surjective,
we have $\rank(E')+\rank(\nbige_W)=\rank(E)$.
Thus, we obtain the following inequality:
\begin{multline}
 0\leq
 P_{\vecy}^{\delta,\alpha_{\ast}}(m)\cdot\rank E
-\delta(m)
-\sum_i\epsilon_i\cdot \dim (W_i)
 +\sum_i\epsilon_i \cdot h^0\bigl(F_{i+1}(E')(m)\bigr)
 -\alpha_1\dim(W)\\
 -P_{\vecy}^{\delta,\alpha_{\ast}}(m)\cdot \rank E'
+\epsilon(E',\phi')\cdot\delta(m)+K.
 \end{multline}
We have the following equality:
\begin{multline}
P_{\vecy}^{\delta,\alpha_{\ast}}(m)\cdot\rank E
-\delta(m)
-\sum_i\epsilon_i \cdot \dim (W_i)
 -\alpha_1\dim(W)\\
=H(m)-\sum \epsilon_i \cdot H_i(m)
 -\sum_i\epsilon_i \cdot \dim (W_i)
 -\alpha_1\dim(V)
 +\alpha_1\dim(V/W).
\end{multline}
Since we have the equalities
$\dim(W_i)=H(m)-H_i(m)$ and $\dim(V)=H(m)$,
the right hand side equals to $\alpha_1\dim (V/W)$.
We have the inequality
$\dim(V/W)\leq \dim H^0\bigl(X,E'(m)\bigr)$.
Hence we obtain the following inequality:
\[
 0\leq -\rank(E')\cdot P^{\delta,\alpha_{\ast}}(m)
 +h^0\bigl(E_{\ast}'(m)\bigr)
+\epsilon(E',\phi)\cdot \delta(m)+K.
\]
Then (\ref{eq;06.4.19.15}) immediately follows.
\hfill\qed

\begin{lem}
 \label{lem;06.4.19.10}
There exists an integer $N_1$ such that
the Yokogawa family $\nYOK(N_1,K,\vecy,L,\delta)$
is bounded.
\end{lem}
\pf
We put $d:=\dim X$
and $g:=c_1\bigl(\nbigo_X(1)\bigr)^d\cap [X]$.
Take a sufficiently negative number $C$
satisfying the following inequality
for any sufficiently large $t$:
\[
\frac{1}{g^{d-1}d!}
 \bigl(C+tg+c\bigr)^{d}
+\delta(t)+K
<P^{\delta,\alpha_{\ast}}(t).
\]
Take a large integer $N_1$
such that the following inequalities hold
for any $m>N_1$:
\[
 C+mg+c>0,\quad
 \delta(m)>0,\quad
 \frac{\bigl(C+mg+c\bigr)^d}{g^{d-1}\cdot d!}
+\delta(m)+K
 <P^{\alpha_{\ast},\delta}(m).
\]
Then we will show that $\nYOK(N_1,K,\vecy,L,\delta)$ is bounded.

Let $\bigl(E_{\ast},\phi\bigr)$ be 
a parabolic $L$-Bradlow pair
contained in the family $\nYOK(N_1,K,\vecy,L,\delta)$.
Let $E'$ denote the last member of
the Harder-Narasimhan filtration of $E$
with respect to the standard semistability.
Assume $\mu(E')<C$, and we will derive the contradiction.
Due to Proposition \ref{prop;06.5.21.3},
we have the following inequality:
\[
 \frac{
 h^0\bigl(E'_{\ast}(m)\bigr)}{\rank(E')}
 \leq
 \frac{h^0(E'(m))}{\rank(E')}
\leq
\frac{\bigl[\mu(E')+mg+c\bigr]_+^d}{g^{d-1}d!}.
\]
By the assumption $\mu(E')<C$,
we obtain the following inequality:
\begin{equation}
 \label{eq;06.4.11.40}
 \frac{h^0\bigl(E'_{\ast}(m)\bigr)
 +\epsilon(E',\phi)\cdot\delta(m)+K}
 {\rank(E')}\\
\leq
\frac{(C+mg+c)^d}{g^{d-1}d!}+\delta(m)+K
 < P^{\delta,\alpha_{\ast}}(m).
\end{equation}
However, it contradicts with
(\ref{eq;06.4.19.15}).
Thus we obtain $\mu(E')>C$.
It implies $\mu_{\max}(E)<C'$ for some constant $C'$,
and thus we obtain that
the family $\nYOK(N_1,K,\vecy,\delta)$ is bounded
due to Proposition \ref{prop;06.5.21.2}.
\hfill\qed

\vspace{.1in}
Then, there exists an integer $N_2$ such that
 the family $\nYOK(N_2,K,\vecy,L,\delta)$
 satisfies the condition $O_m$ for any $m\geq N_2$.

\begin{lem}
Assume $K_1$ is strictly smaller than $\alpha_1$.
Then,  the map $H^0(q):V_m\lrarr H^0\bigl(X,E(m)\bigr)$
is isomorphic for any
$\bigl(q,E,F_{\ast},\phi,W_{\ast},[\widetilde{\phi}]\bigr)
 \in \nYOKtilde(N_2,K_1,\vecy,L,\delta)$.
In particular, $q$ is surjective.
\end{lem}
\pf
We have only to check that $H^0(q)$ is injective.
Let $W$ denote the kernel of $H^0(q)$.
Then we obtain the following inequality
from (\ref{eq;06.4.19.5}):
\[
 -\sum \epsilon_i\cdot \dim (W_i\cap W)
-\alpha_1\cdot \dim(W)+K_1\geq 0.
\]
Since $K_1$ is strictly smaller than $\alpha_1$,
we obtain $\dim(W)=0$,
i.e., $H^0(q)$ is injective.
\hfill\qed

\vspace{.1in}
Let us finish the proof of Proposition \ref{prop;06.5.15.10}.
Let us consider the family $\nbigs_1$ of
parabolic $L$-Bradlow pairs $(E'_{\ast},\phi')$
with the following property:
\begin{itemize}
\item
There exists 
some $(E_{\ast},\phi)\in \nYOK(N_2,K_1,\vecy,L,\delta)$
such that $(E'_{\ast},\phi')$ is isomorphic to
the last member of Harder Narasimhan filtration
of $(E_{\ast},\phi)$
with respect to $\delta$-semistability.
\end{itemize}
Since $\nbigs_1$ is bounded,
the number of the Hilbert polynomials of members 
in $\nbigs_1$ is finite.
Therefore,
there exists a small positive number $K_0<K_1$
and a large integer $N_0\geq N_2$
such that the following holds  for any $m\geq N_0$.
\begin{itemize}
\item
  Let $(E_{\ast},\phi)$ be a member of
 $\nYOK(N_0,K_0,\vecy,L,\delta)$,
 which is not $\delta$-semistable.
 Let $(E'_{\ast},\phi')$ denote last member 
 of the Harder-Narasimhan filtration 
 of $(E_{\ast},\phi)$
 with respect to $\delta$-semistability.
 Then,  the inequality
 $P^{\delta_0}_{(E'_{\ast},\phi')}(t)+K_0
 <P^{\delta_0,\alpha_{\ast}}_{\vecy}(t)$
 holds for any $t\geq m$.
\item
 The family $\nbigs_1$ satisfies the condition $O_m$.
\end{itemize}
Let $(E_{\ast},\phi)$ be a member of
$\nYOK(N_0,K_0,\vecy,L,\delta)$,
and let $(q,E_{\ast},\phi,W_{\ast},[\widetilde{\phi}])$
be its lift in $\nYOKtilde(N_0,K_0,\vecy,L,\delta)$.
Assume $(E_{\ast},\phi)$ is not $\delta$-semistable,
and let $(E'_{\ast},\phi')$  be the last member 
of the Harder-Narasimhan filtration of $(E_{\ast},\phi)$
with respect to the $\delta$-semistability.
Then, we obtain the following inequality from
(\ref{eq;06.4.19.15}) and the second condition:
\[
 P_{\vecy}^{\delta,\alpha_{\ast}}(m)\leq
 \frac{h^0\bigl(E'_{\ast}(m)\bigr)
     +\epsilon(E',\phi)\cdot\delta(m) +K_0}
  {\rank(E')}
\leq P^{\delta}_{(E'_{\ast},\phi')}(m)+K_0
\]
It contradicts with the first condition.
Thus we are done.
\hfill\qed

\subsection{$1$-Stability Condition
 and $2$-Stability Condition}
\label{subsection;06.7.3.25}
\subsubsection{Parabolic sheaf}

Let $\vecy=(y,y_1,\ldots,y_l)$ be an element of $\Type$,
and let $\alpha_{\ast}=(\alpha_1,\ldots,\alpha_l)$
be a system of weights.

\begin{df}
\mbox{{}}
\index{$1$-stability condition}
\index{$2$-stability condition}
\begin{itemize}
\item
We say that the $1$-stability condition
holds for $(\vecy,\alpha_{\ast})$,
if $\nbigm^s(\vecy,\alpha_{\ast})=\nbigm^{ss}(\vecy,\alpha_{\ast})$
holds.
\item
We say that the $2$-stability condition
holds for $(\vecy,\alpha_{\ast})$,
if the automorphism group of
$E_{\ast}\in\nbigm^{ss}(\vecy,\alpha_{\ast})$
is $G_m$ or $G_m^2$.
\hfill\qed
\end{itemize}
\end{df}

\begin{lem}
 \label{lem;06.4.10.20}
The $2$-stability condition 
for $(\vecy,\alpha_{\ast})$ is equivalent
to the following condition:
\begin{itemize}
\item
Let $E_{\ast}$ be a parabolic sheaf of type $\vecy$
with weight $\alpha_{\ast}$,
which is polystable but not stable.
Then we have the unique decomposition
$E_{\ast}=E_{1\,\ast}\oplus E_{2\,\ast}$,
where $E_{i\,\ast}$ are stable parabolic sheaves
with weight $\alpha_{\ast}$
such that $E_{1\,\ast}\not\simeq E_{2\,\ast}$.
\end{itemize}
\end{lem}
\pf
Assume that the $2$-stability condition holds.
Let $E_{\ast}$ be a polystable parabolic sheaf
of type $\vecy$ with weight $\alpha_{\ast}$.
We have a decomposition
$E_{\ast}=E_{1\,\ast}\oplus E_{2\,\ast}$.
If one of $E_{i\,\ast}$ is not stable,
then by taking the graduation of a Jordan-H\"{o}lder filtration,
we obtain a polystable parabolic sheaf
$E'_{\ast}=E'_{1\,\ast}\oplus E'_{2\,\ast}\oplus E'_{3\,\ast}$
of type $\vecy$ with weight $\alpha_{\ast}$.
However, the automorphism of $E'_{\ast}$
is $G_m^3$, which contradicts with
the $2$-stability condition.

Assume that the condition above holds.
Let $E_{\ast}$ be a semistable parabolic sheaf
of type $\vecy$ with weight $\alpha_{\ast}$.
Let $f:E_{\ast}\lrarr E_{\ast}$ be an endomorphism.
The eigenvalues of $f$ are constant.
Let $E_{\ast}=\bigoplus_{i=1}^N E_{i,\ast}$
be the generalized eigen decomposition of $f$.
We have the decomposition $f=\bigoplus_{i=1}^Nf_{i}$.
If $N\geq 3$,
the length of a Jordan-H\"{o}lder filtration is longer than $3$,
and hence we have a polystable object
which has more than three stable components.
Hence $N\leq 2$.

In the case $N=2$,
it can be shown that $E_{i\,\ast}$ are stable
by the same argument.
Hence the automorphism group is $G_m^2$.

Let us consider the case $N=1$.
If $E_{\ast}$ is stable,
the automorphism group is $G_m$.
In the case $E_{\ast}$ is not stable,
the length of the Harder-Narasimhan filtration
is $2$.
Moreover, the graded components
are not mutually isomorphic.
Hence the automorphism group is $G_m$.
\hfill\qed

\vspace{.1in}

In the above argument,
the following corollary is proved.
\begin{cor}
Assume that the $2$-stability condition holds
for $(\vecy,\alpha_{\ast})$.
Let $E_{\ast}$ be a semistable parabolic
sheaf of type $\vecy$ with weight $\alpha_{\ast}$.
Then one of the following holds:
\begin{itemize}
\item
 $E_{\ast}$ is stable.
\item
 $E_{\ast}$ is uniquely decomposed into
 $E_{1\,\ast}\oplus E_{2\,\ast}$,
 where $E_{i\,\ast}$ are stable such that
 $E_{1\,\ast}\not\simeq E_{2\,\ast}$.
\item
 We have the non-split exact sequence
 $0\lrarr E_{1\,\ast}\lrarr E_{\ast}\lrarr E_{2\,\ast}\lrarr 0$,
 where $E_{i\,\ast}$ are stable such that
 $E_{1\,\ast}\not\simeq E_{2\,\ast}$.
\hfill\qed
\end{itemize}
\end{cor}

\subsubsection{Parabolic $L$-Bradlow pair}

Let $L$ be a line bundle on $X$,
and let $\delta$ be an element of $\nbigp^{\br}$.
\begin{df}
\mbox{{}}
\index{$1$-stability condition}
\index{$2$-stability condition}
\begin{itemize}
\item
We say that the $1$-stability condition
holds for $(\vecy,\alpha_{\ast},L,\delta)$,
if $\nbigm^{s}(\vecy,[L],\alpha_{\ast},\delta)=
 \nbigm^{ss}(\vecy,[L],\alpha_{\ast},\delta)$ holds.
\item
We say that the $2$-stability condition
holds for $(\vecy,\alpha_{\ast},[L],\delta)$,
if the automorphism group of
$(E_{\ast},\phi)\in\nbigm^{ss}(\vecy,[L],\alpha_{\ast},\delta)$
is $G_m$ or $G_m^2$.
\hfill\qed
\end{itemize}
\end{df}

\begin{df}
\label{df;06.4.24.1}
\index{critical}
Fix a type $\vecy\in\Type$,
a system of weights $\alpha_{\ast}$
and a line bundle $L$.
A parameter $\delta\in\nbigp^{\br}$ 
is called critical
for $(\vecy,\alpha_{\ast},L)$,
if the $1$-stability condition does not hold
for $(\vecy,\alpha_{\ast},L,\delta)$.
The set of such critical parameters
is denoted by
$\Cr(\vecy,\alpha_{\ast},L)$.
\index{$\Cr(\vecy,\alpha_{\ast},L)$}
\hfill\qed
\end{df}

\begin{lem}
\label{lem;06.5.21.11}
The set $\Cr(\vecy,L,\alpha_{\ast})$ 
of critical values is finite.
\end{lem}
\pf
We may assume $\rank(\vecy)>1$
and $\mu(\vecy,\alpha_{\ast})\geq 0$.
Recall Proposition \ref{prop;06.4.11.1}.
We can take a sufficiently negative number $C$
such that there does not exist
$\delta$-semistable $L$-Bradlow pairs
for any $\delta_{top}\geq -C$.

Let $\nbigs_1$ denote the family of
$L$-Bradlow pairs $(E'_{\ast},\phi')$
with the following property:
\begin{itemize}
\item
There exists a member $(E_{\ast},\phi)$
of $\nbigss(\vecy,L,\alpha_{\ast})$
such that $(E'_{\ast},\phi')$ is a saturated
subobject of $(E_{\ast},\phi)$.
\item
$\deg(E'_{\ast})\geq C$.
\end{itemize}
Since $\nbigs_1$ is bounded
(Proposition \ref{prop;06.5.21.1}),
we obtain the finiteness of the set $\nbigt_1$
of the polynomials $H$ such that
 There exists a member
 $(E'_{\ast},\phi')\in \nbigs_1$
satisfying $H_{E'_{\ast}}=H$.

Let $\delta$ be an element of
$\Cr(\vecy,L,\alpha_{\ast})$.
Then there exists 
a $\delta$-semistable $(E_{\ast},\phi)$ 
such that $\phi\neq 0$,
which has a non-trivial partial Jordan-H\"older filtration
$(E'_{\ast},\phi')\subset (E_{\ast},\phi)$.
We have the following equality:
\[
  \frac{\deg(E'_{\ast})+\epsilon(E',\phi')\cdot \delta_{top}}
  {\rank E'}
=\frac{\deg(E_{\ast})+\delta_{top}}{\rank E}
\]
Since we have $\delta_{top}\leq -C$,
we obtain the following inequality:
\[
 \deg(E'_{\ast})\geq \rank(E')\cdot \mu(E_{\ast})
 -\delta_{top}
\geq C
\]
Hence, $(E'_{\ast},\phi')$
is a member of $\nbigs_1$.
We have the equality
$P_{(E'_{\ast},\phi')}^{\delta}
=P_{(E_{\ast},\phi)}^{\delta}$.
Hence there exist $H\in\nbigt_1$,
an integer $r_1$
and an integer $\epsilon$
satisfying the following:
\[
 \left(
 \frac{1}{r}-\frac{\epsilon}{r_1}
 \right)\cdot\delta
=\frac{H_{E_{\ast}}}{r}
-\frac{H}{r_1},
\quad
 0<r_1<r,
\quad
 \epsilon=0,\mbox{ or }1
\]
Thus we obtain the finiteness
of $\Cr(\vecy,L,\alpha_{\ast})$
from the finiteness of $\nbigt_1$
\hfill\qed

\begin{cor}
If $\delta'\neq\delta$ is sufficiently close
to some element $\delta\in\nbigp^{\br}$,
the $1$-stability condition holds
for $(\vecy,L,\alpha_{\ast},\delta')$.
If $\delta$ is sufficiently close to $0$,
the $1$-stability condition
for $(\vecy,L,\alpha_{\ast},\delta')$.
\hfill\qed
\end{cor}

\begin{lem}
\label{lem;06.5.21.6}
Let $\delta_0$ be an element of $\nbigp^{\br}$.
If $\delta_1\in\nbigp^{\br}$ is sufficiently close to
$\delta_0$,
any $\delta_1$-semistable $L$-Bradlow pair
is also $\delta_0$-semistable.
\end{lem}
\pf
Let $\nbigs$, $\nbigs_1$ and $\nbigt_1$
be as in the proof of
Lemma \ref{lem;06.5.21.11}.
 Let $r_1$ and $\epsilon$ be integers
 such that $0<r_1<r$ and $\epsilon=0,1$.
 We put
 $P(H,r_1,\epsilon,\delta):=
 r_1^{-1}\bigl(H+\epsilon\cdot\delta\bigr)$
 for any $\delta\in\nbigp^{\br}$
and $H\in\nbigt_1$.
If $\delta_1$ is sufficiently close to $\delta_0$,
the following holds:
\begin{description}
\item[(A)]
 $P(H,r_1,\epsilon,\delta_1)(t)>0$
implies $P(H,r_1,\epsilon,\delta_0)(t)\geq 0$
for any sufficiently large $t$.
\end{description}

Let $\delta_1$ be as above.
We may assume that
the $1$-stability condition holds
for $(\vecy,L,\alpha_{\ast},\delta_1)$.
Let $(E_{\ast},\phi)$ be 
$\delta_1$-stable $L$-Bradlow pair
of type $\vecy$,
and we assume that it is not $\delta$-semistable.
Then there exists a saturated subobject
$(E'_{\ast},\phi')$ of $(E_{\ast},\phi)$
with the following property:
\begin{itemize}
\item
$P^{\delta_0}_{(E'_{\ast},\phi')}(t)\geq 
 P^{\delta_0}_{(E_{\ast},\phi)}(t)$
 for any sufficiently large $t$.
\item
 $P^{\delta_1}_{(E'_{\ast},\phi')}(t)\leq
 P^{\delta_1}_{(E_{\ast},\phi)}(t)$
 for any sufficiently large $t$.
\end{itemize}
Then $(E'_{\ast},\phi')$ is a member of $\nbigs_1$
due to the first inequality.
Therefore, the two inequalities
contradict with the condition (A) above.
Thus we are done.
\hfill\qed

\vspace{.1in}

By a similar argument,
we can show the following.
\begin{lem}
\label{lem;06.5.21.20}
Let $\vecy$ be an element of $\Type$,
and let $\alpha_{\ast}$ be a system of weights.
Assume $\delta$ is sufficiently small.
Then, 
$E_{\ast}$ is semistable
if $(E_{\ast},\phi)$ be a $\delta$-stable
$L$-Bradlow pair.
\end{lem}
\pf
Let $\nbigs$ be a family of 
$\mu$-semistable parabolic torsion-free sheaves
of type $\vecy$.
Let $\nbigsbar$ denote the family of
parabolic torsion-free sheaves $E'_{\ast}$
with the following property:
\begin{itemize}
\item
There exists $E_{\ast}\in\nbigs$
such that $E'_{\ast}$ is a saturated subobject
of $E_{\ast}$.
Moreover, we have $\mu(E'_{\ast})=\mu(\vecy,\alpha_{\ast})$.
\end{itemize}
Let $\nbigt$ denote the set of
the polynomials $P$
such that $P\neq P^{\alpha_{\ast}}_{\vecy}$
and $P=P_{E'_{\ast}}$ for some
$E'_{\ast}\in \nbigsbar$.
Since the families $\nbigs$ and $\nbigsbar$ are bounded,
the set $\nbigt$ is finite.
We take a positive number $\delta_1$
satisfying the following:
\[
 0<\delta_1<
 \frac{1}{10\cdot \rank(\vecy)}
 \min\bigl\{
 |P-P^{\alpha_{\ast}}_{\vecy}|\,\big|\,
 P\in\nbigt
 \bigr\}
\]
We regard $\delta_1$ as the polynomial
of degree $0$.

Let $(E_{\ast},\phi)$ be a $\delta_1$-semistable
$L$-Bradlow pair of type $\vecy$
with weight $\alpha_{\ast}$.
It is easy to observe that
$E_{\ast}$ is $\mu$-semistable.
Let $E'_{\ast}$ be a subobject of $E_{\ast}$.
In the case $\mu(E'_{\ast})<\mu(E_{\ast})$,
we obviously have
$P_{E'_{\ast}}<P_{E_{\ast}}$.
Assume $\mu(E'_{\ast})=\mu(E_{\ast})$
and $P_{E'_{\ast}}\neq P_{E_{\ast}}$.
Then $E'_{\ast}$ is a member of
$\nbigsbar$,
and $P_{E'_{\ast}}$ is a member of
$\nbigt$.
Due to $\delta_1$-semistability of $(E_{\ast},\phi)$,
we have the following:
\[
 P_{E'_{\ast}}+\frac{\epsilon(E'_{\ast},\phi)\cdot \delta_1}{\rank E'}
\leq
 P_{E_{\ast}}+\frac{\delta_1}{\rank E}
\]
It implies $P_{E'_{\ast}}< P_{E_{\ast}}$
due to our choice of $\delta_1$.
\hfill\qed

\vspace{.1in}

By an argument similar to the
proof of Lemma \ref{lem;06.4.10.20},
we obtain the following lemma.
\begin{lem}
The $2$-stability condition 
for $(\vecy,\alpha_{\ast},L,\delta)$ is equivalent to
the following:
\begin{itemize}
\item
Let $(E_{\ast},\phi)$ be a parabolic $L$-Bradlow pair
 of type $\vecy$
with weight $\alpha_{\ast}$ such that $\phi\neq 0$,
which is $\delta$-polystable but not $\delta$-stable.
Then we have the unique decomposition
$(E_{\ast},\phi)=(E_{1\,\ast},\phi_1)\oplus E_{2\,\ast}$,
where $(E_{1\,\ast},\phi_1)$ is $\delta$-stable
and $E_{2\,\ast}$ is stable.
\end{itemize}

Moreover,
when the $2$-stability condition holds
for $(\vecy,\alpha_{\ast},L,\delta)$,
one of the following holds
for any $\delta$-semistable parabolic $L$-Bradlow pair
$(E_{\ast},\phi)$ of type $\vecy$ with weight $\alpha_{\ast}$
such that  $\phi\neq 0$.
\begin{itemize}
\item
 $(E_{\ast},\phi)$ is $\delta$-stable.
\item
 $(E_{\ast},\phi)$ is uniquely decomposed into
 $(E_{\ast},\phi)=(E_{1\,\ast},\phi_1)\oplus E_{2\,\ast}$,
 where $(E_{1\,\ast},\phi)$ is $\delta$-stable
 and $E_{2\,\ast}$ is stable.
\item
 We have the non-split exact sequence
 $0\lrarr (E_{1\,\ast},\phi_1)\lrarr (E_{\ast},\phi)
 \lrarr E_{2\,\ast}\lrarr 0$
 or 
 $0\lrarr  E_{2\,\ast}\lrarr
 (E_{\ast},\phi)\lrarr  (E_{1\,\ast},\phi_1) \lrarr 0$,
 where $(E_{1\,\ast},\phi)$ is $\delta$-stable
 and $E_{2\,\ast}$ is stable.
\hfill\qed
\end{itemize}
\end{lem}

\subsubsection{Parabolic $\vecL$-Bradlow pair}
\label{subsubsection;06.6.15.5}

Let $\vecL=(L_1,L_2)$ be a pair of line bundles on $X$,
and let $\vecdelta=(\delta_1,\delta_2)$ be elements of 
$\bigl(\nbigp^{\br}\bigr)^2$.
Similarly, we have the notion of
$1$-stability and $2$-stability conditions.
\begin{df}
\mbox{{}}
\index{$1$-stability condition}
\index{$2$-stability condition}
\label{df;06.6.8.15}
\begin{itemize}
\item
We say that the $1$-stability condition
holds for $(\vecy,\alpha_{\ast},\vecL,\vecdelta)$,
if $\nbigm^{s}(\vecy,[\vecL],\alpha_{\ast},\vecdelta)
=\nbigm^{ss}(\vecy,[\vecL],\alpha_{\ast},\vecdelta)$.
\item
We say that the $2$-stability condition
holds for $(\vecy,\alpha_{\ast},\vecL,\vecdelta)$,
if the automorphism group of any
$(E_{\ast},\vecphi)\in
 \nbigm^{ss}(\vecy,[\vecL],\alpha_{\ast},\vecdelta)$
is $G_m$ or $G_m^2$.
\hfill\qed
\end{itemize}
\end{df}

\begin{df}
\mbox{{}}
\index{critical}
\begin{itemize}
\item
Fix a type $\vecy\in\Type$,
a system of weights $\alpha_{\ast}$
and a pair of line bundles $\vecL$.
A parameter $\vecdelta\in\nbigp^{\br\,2}$ 
is called critical
for $(\vecy,\alpha_{\ast},\vecL)$,
if the $1$-stability condition does not hold
for $(\vecy,\alpha_{\ast},\vecL,\vecdelta)$.
The set of such critical parameters
is denoted by
$\Cr(\vecy,\alpha_{\ast},\vecL)$.
\item
We also put as follows:
\[
 \Cr(\vecy,\alpha_{\ast},\vecL,\delta_1):=
 \bigl\{\delta_2\in\nbigp^{\br}\,\big|\,
 (\delta_1,\delta_2)\in\Cr(\vecy,\alpha_{\ast},\vecL)
 \bigr\}
\]
Any element $\delta_2\in\Cr(\vecy,\alpha_{\ast},\vecL,\delta_1)$
is called critical for
$(\vecy,\alpha_{\ast},\vecL,\delta_1)$.
\hfill\qed
\end{itemize}
\end{df}

We can show the following lemma 
by using Lemma \ref{lem;06.5.21.10}
and an argument similar to the proof 
of Lemma \ref{lem;06.5.21.11}.
\begin{lem}
The set $\Cr(\vecy,\vecL,\alpha_{\ast},\delta_1)$
is finite.
\hfill\qed
\end{lem}

By an argument similar to the proof of
Lemma \ref{lem;06.5.21.6},
we can show the following lemma.
\begin{lem}
\label{lem;06.5.21.15}
If $\delta_2'$ are sufficiently close to $\delta_2$,
any $(\delta_1,\delta_2')$-semistable $\vecL$-Bradlow pair
is also $(\delta_1,\delta_2)$-semistable.

If $\delta_2'$ is sufficiently close to $0$,
the $L_1$-Bradlow pair
$(E_{\ast},\phi_1)$ is $\delta_1$-semistable
for any $(\delta_1,\delta_2)$-semistable
$\vecL$-Bradlow pair
$(E_{\ast},\phi_1,\phi_2)$.
\hfill\qed
\end{lem}

\begin{lem}
Assume that $\delta_1+\delta_2$ is sufficiently small
as in Lemma {\rm\ref{lem;06.5.21.20}}.
Then, 
if $(E_{\ast},\phi_1,\phi_2)$is 
a $(\delta_1,\delta_2)$-semistable
$\vecL$-Bradlow pair,
$E_{\ast}$ is semistable.
\hfill\qed
\end{lem}

\begin{lem}
\label{lem;06.6.13.11}
If both of $\delta_i$ are sufficiently small,
the following claims hold:
\begin{itemize}
\item
 If the $1$-stability condition holds
 for $(\vecy,\alpha_{\ast})$,
 then the $1$-stability condition holds also for
 $(\vecy,\alpha_{\ast},\vecL,\vecdelta)$.
\item
 Even if  the $1$-stability condition does not hold
 for $(\vecy,\alpha_{\ast})$,
 the $2$-stability condition holds for
 $(\vecy,\alpha_{\ast},\vecL,\vecdelta)$.
\item
 If the $1$-stability condition does not hold
 for $\vecdelta=(\delta_1,\delta_2)$,
 the equality
 $\delta_1/r_1=\delta_2/r_2$ holds 
 for some decomposition $r_1+r_2=\rank(\vecy)$.
\end{itemize}
\end{lem}
\pf
We take $\delta_i$ as in the proof of 
Lemma \ref{lem;06.5.21.20}.
The first claim is clear.
Let us see the second claim.
Let $(E_{\ast},\phi_1,\phi_2)$ be 
a $(\delta_1,\delta_2)$-polystable
$(L_1,L_2)$-Bradlow pair of type $\vecy$
with weight $\alpha_{\ast}$
such that $\phi_i\neq 0$.
We remark that $E_{\ast}$ is semistable.
Since $\delta_i$ are sufficiently small,
the number of stable components are at most $2$.
If it is decomposed,
it is of the form $(E_{1\,\ast},\phi_1,0)\oplus (E_{2\,\ast},0,\phi_2)$,
where $(E_{i\,\ast},\phi_i)$ are $\delta_i$-stable
with $\phi_i\neq 0$.
The components are not isomorphic.
Hence the $2$-stability condition holds.
The third claim also immediately follows.
\hfill\qed

\subsection{Quot Schemes}
\label{subsection;06.5.15.15}
\subsubsection{Preliminary}

Let $y$ be an element of $H^{\ev}(X)$.
We will denote $y\cdot \ch\bigl(\nbigo(m)\bigr)$
by $y(m)$ for any integer $m$.
We also have the element $\det(y)=y_1$,
where $y_1$ denote the $H^2(X)$-component of $y$.
\index{$y(m)$}
\index{$\det(y)$}

Let $H_y$ be the associated Hilbert polynomial to $y$
(the subsubsection \ref{subsubsection;06.5.15.1}).
Take a large integer $m$ such that $H_y(m)>0$.
We also assume that
any line bundles $M$ with $c_1(M)=\det\bigl(y(m)\bigr)$
satisfy $H^i(X,M)=0$ $(i>0)$.

We take an $H_y(m)$-dimensional vector space $V_{m}$ over $k$.
We denote $V_{m}\otimes\nbigo_X$ by $V_{m,X}$.

\subsubsection{Quotient sheaves}
\label{subsubsection;06.5.9.50}

A pair of $U$-coherent sheaf $\nbige$ on $U\times X$
and a surjection $q:p_U^{\ast}V_{m,X}\lrarr \nbige$
is called a $U$-quotient sheaf of $V_{m,X}$,
which is denoted by $(q,\nbige)$ or simply by $q$.
A $U$-quotient sheaf  $(q,\nbige)$
with an orientation $\rho$
is called an oriented $U$-quotient sheaf of $V_{m,X}$.
The type of $(q,\nbige)$ or $(q,\nbige,\rho)$
is defined to be the type of $\nbige(-m)$.
(See the subsubsection \ref{subsubsection;06.4.24.10}
for the type.)

Recall that 
the moduli functor of quotient sheaves of
$V_{m,X}$ with type $y$
is representable,
and the moduli scheme is projective. (\cite{gro}).
We denote it by $Q(m,y)$. \index{$Q(m,y)$}
Let $(q^u,\nbige^u)$ denote
the universal quotient sheaf of
$V_{m,X}$ on $Q(m,y)\times X$.
A point of $Q(m,y)$ is denoted
by the corresponding quotient
$\bigl(q,\nbige\bigr)$.
We have the right action of $GL(V_m)$
on $Q(m,y)$
given by $g\cdot (q,\nbige):=(q\circ g,\nbige)$.

Let $(q,\nbige)$ be a $U$-quotient sheaf
of $V_{m,X}$ with type $y$
defined over $U\times X$.
We say that $\bigl(q,\nbige\bigr)$ satisfies
the (TFV)-condition,
if the following holds for any $u\in U$:
\begin{description}
\index{(TFV)-condition}
\item[(TFV):]
The sheaf $\nbige_{|\{u\}\times X}$ is torsion-free,
the induced map
$V_m\lrarr H^0\bigl(X,\nbige_{|\{u\}\times X}\bigr)$
is isomorphic,
and 
$H^i\bigl(X,\nbige_{|\{u\}\times X}\bigr)$ vanish for any $i>0$.
\end{description}
In general,
the condition determines the maximal open subset of $U$
on which the (TFV)-condition holds.
In particular,
it determines the open subset of $Q(m,y)$,
which is denoted by $Q^{\circ}(m,y)$.
\index{$Q^{\circ}(m,y)$}

Let $\Or(\nbige^u)$ denote the orientation bundle,
which is the line bundle on $Q(m,y)$.
The moduli functor of oriented quotient sheaves
of $V_{m,X}$ with type $y$ is representable
by $Q(m,\widehat{y}):=\Or(\nbige^u)^{\ast}$.
We have the induced right action of $GL(m)$
on $Q(m,\widehat{y})$.
We put $Q^{\circ}(m,\yhat):=
 Q^{\circ}(m,y)\times_{Q(m,y)}Q(m,\yhat)$.
\index{$Q(m,\yhat)$, $Q^{\circ}(m,\yhat)$}

\subsubsection{Quotient quasi-parabolic sheaves
 and the Maruyama-Yokogawa construction}
\label{subsubsection;06.5.14.1}

Let $D$ be a Cartier divisor of $X$.
Let $y$ and $V_{m,X}$ be given as above.
A $U$-quasi-parabolic quotient sheaf $(q,\nbige,F_{\ast})$
of $V_{m,X}$ on $U\times (X,D)$ is 
defined to be a $U$-quotient sheaf $(q,\nbige)$
of $V_{m,X}$
with a $U$-quasi-parabolic structure $F_{\ast}$
of $\nbige$ at $D$.
The type of $U$-quasi parabolic quotient sheaf
$(q,\nbige,F_{\ast})$ of $V_{m,X}$
is defined to be the type
of the underlying $U$-quasi-parabolic sheaf
$\bigl(\nbige(-m),F_{\ast}\bigr)$.
(See the subsubsection \ref{subsubsection;06.4.24.10}
for the type.)

Let $\vecy$ be an element of $\Type$,
whose $H^{\ast}(X)$-component is $y$.
Let $\bigl(q, \nbige,F_{\ast}\bigr)$ be
a $U$-quotient quasi-parabolic sheaf of $V_{m,X}$
with type $\vecy$ on $U\times (X,D)$.
We say that $\bigl(q,\nbige,F_{\ast}\bigr)$
satisfies the (TFV)-condition
for quasi-parabolic shaves,
if the following holds for any $u\in U$ and for any $i$:
\begin{description}
\index{(TFV)-condition}
\item[(TFV):]
$(q,\nbige)$ satisfies the (TFV)-condition.
Moreover,
the sheaves
$F_i(\nbige)_{|\{u\}\times X}$
and $\Cok_i(\nbige)_{|\{u\}\times X}$ 
are generated by global sections,
and the higher cohomology groups
of $F_i(\nbige)_{|\{u\}\times X}$
and $\Cok_i(\nbige)_{|\{u\}\times X}$
vanish.
\end{description}
In general, the condition determines
the maximal open subset of $U$
on which the (TFV)-condition holds.

\vspace{.1in}

We put $H_i:=H_{\vecy,i}$.
Let $Q_{m,i}$ denote the scheme
representing the moduli functor of
the quotient sheaves of $V_{m,X}$
whose Hilbert polynomials are $H_i$.
We have the natural right $GL(V_m)$-action
on $Q_{m,i}$.
We have the open subset $U_{m,i}$ of $Q_{m,i}$
given by the conditions
that $H^j\bigl(X,\nbige_i\bigr)=0$ for any $j>0$
and that $V_m\lrarr H^0\bigl(X,\nbige_i\bigr)$ is surjective.
Let $G_{m,i}$ denote the Grassmannian variety,
which is the moduli of $H_i(m)$-dimensional
quotient space of the vector space $V_m$.
We have the $GL(V_m)$-equivariant morphism of
$U_{m,i}$ to $G_{m,i}$
by the correspondence:
\[
 \bigl(q_i,\nbige_i\bigr)\longmapsto
 \Bigl(H^0(q_i):V_{m}\rarr H^0\bigl(X,E_i(m)\bigr)\Bigr).
\]

Let $Q^{tf}(m,y)$ denote the open subset
of $Q(m,y)$ consisting of the points
corresponding to the torsion-free quotients.
For the construction of the moduli of semistable parabolic sheaves,
Maruyama and Yokogawa constructed the scheme
$\Gamma$ which is obtained as the subscheme of
$Q^{tf}(m,y)\times\prod_i^l U_{m,i}$.
(See section 3 of \cite{my}.)
The scheme $\Gamma$ is the moduli of 
quotient quasi-parabolic sheaves 
$(q,\nbige,F_{\ast})$ of $V_{m,X}$ with type $\vecy$
satisfying the following conditions:
(i) $\nbige$ is torsion-free.
(ii) The higher cohomology groups of $\Cok_i(\nbige)$
 vanish, and $V_{m}\lrarr H^0\bigl(X,\Cok_i(\nbige)\bigr)$ are surjective
 for any $i$.

Moreover, the (TFV)-condition determines the open subset
of $\Gamma$,
which is denoted by $Q^{\circ}(m,\vecy)$.\index{$Q^{\circ}(m,\vecy)$}
The scheme $Q^{\circ}(m,\vecy)$ represents
the moduli functor of quotient quasi-parabolic sheaves
of $V_{m,X}$ with type $\vecy$
satisfying (TFV)-conditions.
We have the universal objects on 
$Q^{\circ}(m,\vecy)\times X$,
which is denoted by
$\bigl(q^u,\nbige^u,F^u_{\ast}\bigr)$.
We have the right $GL(V_m)$-action
on $Q^{\circ}(m,\vecy)$
given by $g\cdot(q,\nbige,F_{\ast}):=(q\circ g,\nbige,F_{\ast})$.

\vspace{.1in}

Let $\Pic_X\bigl(c\bigr)$ denote 
the component of the Picard variety of $X$
such that any line bundle $M\in \Pic_X\bigl(c\bigr)$
satisfy $c_1(M)=c$.
Let $\poin_X\bigl(c\bigr)$ denote the Poincar\'{e} bundle
on $\Pic_X\bigl(c\bigr)\times X$.
Then we obtain the locally free sheaf:
\begin{equation}
\label{eq;06.5.21.30}
\Zhat_m:=
 p_{X\,\ast}
 \nhom\Bigl(\bigwedge^{\rank y} V_{m,X},
 \poin_X\bigl(\det y(m)\bigr)\Bigr) .
\end{equation}
\index{$\Zhat_m$, $Z_m$}
The projectivization $Z_m$
is called the Gieseker space.
We have the natural right action of $GL(V_m)$
on $Z_m$.

It is known that
the correspondence
$(q,\nbige)
\longmapsto H^0\bigl(\bigwedge^{r}q\bigr)$
gives the $GL(V_m)$-equivariant
morphism of $Q^{\circ}(m,y)$ to $Z_m$,
which is known to be an immersion.
Therefore, we obtain the morphism
$Q^{\circ}(m,y)\times \prod_{i}U_{m,i}
 \lrarr Z_m\times \prod_iG_{m,i}$.
By the composition of the inclusion
$Q^{\circ}(m,\vecy)
 \lrarr Q^{\circ}(m,y)\times \prod_iU_{m,i}$,
we obtain the $GL(V_m)$-equivariant morphism
$Q^{\circ}(m,\vecy)\lrarr Z_m\times\prod_i G_{m,i}$.
The following lemma was shown in \cite{my}.

\begin{lem}
 [\cite{my}, Proposition 3.2] \label{lem;06.4.10.1}
 The morphism
 $Q^{\circ}(m,\vecy)\lrarr Z_m\times\prod_i G_{m,i}$
 is an immersion.
\end{lem}
\pf
We give only some remarks.
Since the morphism
$Q^{\circ}(m,x)\lrarr Z_m$ is an immersion,
we have only to care the morphism
of $Q^{\circ}(m,\vecy)$ 
to $Q^{\circ}(m,y)\times \prod G_{m,i}$.
Recall the precise result of Maruyama and Yokogawa:
Let $\alpha_{\ast}$ be a system of weights.
Let $\Gamma^{ss}$ denote the open subscheme
of $\Gamma$ such that
the corresponding parabolic sheaves
$(\nbige,F_{\ast},\alpha_{\ast})$ are semistable.
We may assume that 
the (TFV)-condition holds
for each member $\Gamma^{ss}$.

They showed that
the morphism of $\Gamma^{ss}$ to
$Q^{\circ}(m,y)\times \prod_{i}G_{m,i}$ is immersion.
Their argument can be summarized as follows:
\begin{description}
\item
[(i)]
Construct a subscheme $\nbigz$ of
$Q(m,y)\times \prod_i Q_{m,i}\times \prod_iG_{m,i}$
such that
the projection of $\nbigz$
to $\qmy\times\prod_i G_{m,i}$
is immersive.
($\nbigz$ is denoted as
 $\Delta_1^0\times_{Q}
  \times_Q\Delta_2^0\times_{Q}
   \cdots\times_{Q}\Delta_l^0$ in \cite{my}.)
\item
[(ii)]
The morphism $\Gamma^{ss}\lrarr \prod_iG_{m,i}$ 
gives the graph $\nbigg$, which is a subset of
$\qmy\times \prod_i Q_{m,i}\times \prod_{i}G_{m,i}$.
It can be shown that $\nbigg$ is contained in $\nbigz$.
Then the projection of $\nbigg$
to $Q(m,y)\times\prod_i G_{m,i}$
is immersive.
Hence the morphism
$\Gamma^{ss}\lrarr Q^{tf}(m,y)\times\prod_i G_{m,i}$ is immersion.
\end{description}
For the argument,
we need the only fact that each member of
$\Gamma^{ss}$ satisfies the (TFV)-condition.
(It is given as the conditions (3.0.1) and (3.0.2)
in \cite{my}.)
Thus the morphism
$Q^{\circ}(m,\vecy)
 \lrarr Z_m\times \prod_iG_{m,i}$ is immersive.
See \cite{my} for more detail.
\hfill\qed

\subsubsection{Quotient parabolic $L$-Bradlow pair}  

Let $L$ be a line bundle on $X$.
Fix a non-trivial morphism
$\iota:\nbigo(-m)\lrarr L$.
Let $q:p_U^{\ast}V_{m,X}\lrarr \nbige$ be
a $U$-quotient sheaf with type $y$
defined over $U\times X$ satisfying (TFV)-condition.
We have the morphism
$\nbige(-m)\otimes L^{-1}\lrarr \nbige$ induced by $\iota$.
Let $\phi$ be an $L$-section of $\nbige(-m)$.
Then $\iota$ and $\phi$
induce the $\nbigo$-section $\iota(\phi)$
of $\nbige$.

\begin{df}[Quotient $L$-section]
Let $(q,\nbige)$ be a quotient sheaf of $\vmx$
defined over $U\times X$.
An $\nbigo_X$-section $\overline{\phi}$
of $p_U^{\ast}\vmx$ is called
a quotient $L$-section of $\bigl(q,\nbige\bigr)$
with respect to $\iota$,
if there exists an $L$-section of $\nbige(-m)$
such that $\iota(\phi)=q\circ\overline{\phi}$.

The condition means that
the element $q\circ\overline{\phi}\in
 \Gamma\bigl(U\times X,\nbige\bigr)$
is contained in $\Gamma\bigl(U\times X,\nbige(-m)\otimes L^{-1}\bigr)$,
where the inclusion $\nbige(-m)\otimes L^{-1}\lrarr \nbige$
is given by $\iota$.
\hfill\qed
\end{df}

\begin{df}
A quotient $U$-quasi-parabolic $L$-Bradlow pair 
of $V_{m,X}$ with type $\vecy$
on $U\times (X,D)$ is defined to be 
a pair $\bigl(q,\nbige,F_{\ast},\overline{\phi}\bigr)$
of quotient $U$-quasi parabolic sheaf 
 $(q,\nbige,F_{\ast})$ of $V_{m,X}$
with type $\vecy$ over $U\times (X,D)$
and a quotient $L$-section $\overline{\phi}$
with respect to $\iota$.
\hfill\qed
\end{df}

Let us construct the moduli scheme
$Q^{\circ}(m,\vecy,L)$ \index{$Q^{\circ}(m,\vecy,L)$}
of quotient quasi-parabolic $L$-Bradlow pairs
of $V_{m,X}$ with type $\vecy$
satisfying (TFV)-condition,
whose $L$-section is non-trivial everywhere.
We put
$Q^{\circ}\bigl(m,\vecy,\nbigo(-m)\bigr):=
 Q^{\circ}\bigl(m,\vecy\bigr)\times V_m^{\ast}$.
Let $\pi$ denote the projection of
$Q^{\circ}\bigl(m,\vecy,\nbigo(-m)\bigr)\times X$
onto $Q^{\circ}\bigl(m,\vecy\bigr)\times X$.
On $Q^{\circ}\bigl(m,\vecy,\nbigo(-m)\bigr)\times X$,
we have the quotient quasi-parabolic sheaf
$\pi^{\ast}(q^u,\nbige^u,F^u_{\ast})$
of $V_{m,X}$ with type $\vecy$.
We also have the canonical 
$\nbigo_X$-section $\overline{\phi}^u$
of $p^{\ast}_{Q^{\circ}(m,\vecy,\nbigo(-m))}V_{m,X}$,
which is induced by the identity of $V_{m}$.

We have the composite $\Lambda$ of the following morphisms
on $Q^{\circ}\bigl(m,\vecy,\nbigo(-m)\bigr)\times X$,
where the last quotient sheaf is induced by $\iota$:
\[
 \nbigo_{Q^{\circ}(m,\vecy,\nbigo(-m))}
  \stackrel{\overline{\phi}^u}{\lrarr}
 p_{Q^{\circ}(m,\vecy,\nbigo(-m))}^{\ast}V_{m,X}
\stackrel{\pi^{\ast}q^u}{\lrarr}
\pi^{\ast}\nbige^u
\lrarr
 \frac{\pi^{\ast}\nbige^u}{\pi^{\ast}\nbige^u(-m)\otimes L^{-1}}.
\]

Recall the following result (Lemma 4.3 of \cite{y3})
due to Yokogawa.
\begin{lem}
 \label{lem;4.26.1}
Let $f:\nbigx\lrarr S$ be a proper morphism of Noetherian
schemes and $\phi:I\lrarr F$ be an $\nbigo_{\nbigx}$-morphism
of coherent sheaves.
Assume that $F$ is flat over $S$.
Then there exists a unique closed subscheme $Z$ of $S$
such that for all morphism $g:T\lrarr S$, $g^{\ast}\phi=0$
if and only if $g$ factors through $Z$.
\hfill\qed
\end{lem}

Remark the following easy lemma.
\begin{lem}
Let $E$ be a $U$-coherent sheaf on $X\times U$
such that $E_{|X}$ is torsion-free for any $u\in U$.
Let $D$ be a Cartier divisor of $X$.
Then $E\otimes\nbigo_D$ is flat over $U$.
\end{lem}
\pf
Let $f$ be any local section of $\nbigo_X$,
and let $\nbigf$ be any $\nbigo_U$-coherent sheaf.
We have only to show the injectivity of
the endomorphism of $E\otimes_{\nbigo_U}\nbigf$
induced by $f$.
By considering the support of $\pi_X\bigl(\Ker(f)\bigr)$,
we can reduce the case $U=\Spec(K)$
and $\nbigf=\nbigo_U$ for some field $K$.
In the case, the claim is trivial.
\hfill\qed

\vspace{.1in}
Due to the above two lemmas,
the vanishing condition of $\Lambda$
gives the closed subscheme of $Q^{\circ}\bigl(m,\vecy,\nbigo(-m)\bigr)$,
which is $Q^{\circ}(m,\vecy,L)$.
By the construction,
it is easy to see that
$Q^{\circ}(m,\vecy,L)$ has the desired property.
We denote the universal pair
on $Q^{\circ}(m,\vecy,L)\times X$
by $\bigl(q^u,\nbige^u,F^u_{\ast},\overline{\phi}^u\bigr)$.
We remark that we have the unique $L$-section
$\phi^u$ of $\nbige^u(-m)$
such that $\iota(\phi^u)=q_{\ast}(\overline{\phi}^u)$.

\vspace{.1in}

The right $GL(V_m)$-action
on $Q^{\circ}\bigl(m,\vecy,\nbigo(-m)\bigr)$ is given by
$g\cdot(q,\nbige,F_{\ast},\phi)=
\bigl(q\circ g,\nbige,F_{\ast},g^{-1}\circ \phi\bigr)$.
The action can be naturally lifted to the action
on the universal object.
Note that $Q^{\circ}\bigl(m,\vecy,L\bigr)$ 
is a closed $GL(V_m)$-stable subscheme
of $Q^{\circ}\bigl(m,\vecy,\nbigo(-m)\bigr)$.

\vspace{.1in}

\subsubsection{Quotient reduced $L$-Bradlow pair}

Let $\iota:\nbigo\lrarr L(m)$ be the fixed non-trivial morphism.
Let $q:p_U^{\ast}\vmx\lrarr \nbige$ be
a $U$-quotient sheaf of $\vmx$
defined over $U\times X$.
A reduced $\nbigo_X$-section $[\overline{\phi}]$
of $p_U^{\ast}\vmx$ and $q$ induce
the reduced $\nbigo_X$-section
$q_{\ast}([\phi])$ of $\nbige$.
On the other hand,
a reduced $L$-section $[\phi]$ of $\nbige(-m)$
and $\iota$ induce the reduced
$\nbigo_X$-section $\iota([\phi])$
of $\nbige$.

\begin{df}
\mbox{{}} \label{df;06.4.10.10}
\begin{itemize}
 \item
A reduced $\nbigo_X$-section $[\bar{\phi}]$ of $p_U^{\ast}\vmx$
is called a quotient reduced $L$-section
if there exists a reduced $L$-section of $\nbige(-m)$
such that $\iota([\phi])=q_{\ast}([\bar{\phi}])$.

\item
A quotient $U$-quasi-parabolic reduced $L$-Bradlow pair
of $V_{m,X}$
with type $\vecy$ on $U\times X$
is defined to be 
a tuple $\bigl(q,\nbige,F_{\ast},[\overline{\phi}]\bigr)$
of quotient $U$-quasi-parabolic sheaf $(q,\nbige,F_{\ast})$
of $V_{m,X}$ with type $\vecy$
and a reduced $L$-section $[\overline{\phi}]$ of $\nbige(-m)$.
We also assume that $[\overline{\phi}]$
is non-trivial everywhere.
\hfill\qed
\end{itemize}
\end{df}

Let us construct the scheme
$Q^{\circ}(m,\vecy,[L])$ \index{$Q^{\circ}(m,\vecy,[L])$}
representing the moduli functor
of quotient quasi parabolic reduced $L$-Bradlow pair
of $V_{m,X}$ with type $\vecy$
satisfying (TFV)-condition.
We have the free $G_m$-action on the scheme 
$Q^{\circ}(m,\vecy,L)$ defined by 
$t\cdot(q,\nbige,F_{\ast},\overline{\phi})=
 \bigl(q,\nbige,F_{\ast},
 t\cdot\overline{\phi}\bigr)$.
Then we put
$Q^{\circ}(m,\vecy,[L]):=
 Q^{\circ}(m,\vecy,L)/G_m$.
It is the closed subscheme
of $Q^{\circ}(m,\vecy)\times\proj_m$.
Due to the naturally defined morphism
$\pi:Q^{\circ}(m,\vecy,[L])\times X
 \lrarr Q^{\circ}(m,\vecy)\times X$,
we have the quotient quasi parabolic sheaf
$\bigl(\widehat{q}^u,\widehat{\nbige}^u,
 \widehat{F}^{u}_{\ast}\bigr)
:=(\pi^{\ast}q^u,\pi^{\ast}\nbige^u,\pi^{\ast}F^u_{\ast})$.
The morphism $\phibar^u$ naturally induces
the reduced $\nbigo_X$-section
$[\phibar^u]: p_{Q^{\circ}(m,\vecy,[L])}^{\ast}\nbigo_X
\otimes
 \nbigo_{\proj_m}(-1)
\lrarr 
p_{Q^{\circ}(m,\vecy,[L])}^{\ast}V_{m,X}$.
We also have the reduced $L$-section
$[\phi^u]: p_{Q^{\circ}(m,\vecy,[L])}^{\ast}L
\otimes
 \nbigo_{\proj_m}(-1)
\lrarr 
p_{Q^{\circ}(m,\vecy,[L])}^{\ast}\widehat{\nbige}(-m)$.

\begin{lem}
$Q^{\circ}(m,\vecy,[L])$
has the desired universal property,
and $\bigl(\widehat{q}^u,\widehat{\nbige}^u,
 \widehat{F}^u_{\ast},[\phi^u]\bigr)$
is the universal object.
\end{lem}
\pf
We give only an outline.
Due to Lemma \ref{lem;4.26.1},
we can reduce the case $L=\nbigo(-m)$
and $\iota$ is the identity.
Let $(q,\nbige,F_{\ast},[\phibar])$ denote
a $U$-quotient quasi-parabolic reduced $L$-Bradlow pair,
satisfying (TFV)-condition.
We have the map
$F:U\lrarr Q^{\circ}(m,\vecy)$ corresponding to
$(q,\nbige,F_{\ast})$.

We have the locally free sheaf 
$p_{X\,\ast}\nbige$
and $p_{X\,\ast}\widehat{\nbige}^u$ on $U$
and $Q^{\circ}(m,\vecy)$ respectively.
Let $\proj_1$ and $\proj_2$ denote the projectivization
of them.
We have $\Psi^{\ast}\proj_2\simeq\proj_1$,
and we have the natural morphism
$\Psitilde:\proj_1\lrarr \proj_2$.
We remark that $\proj_2$ is naturally isomorphic to
$Q^{\circ}\bigl(m,\vecy,[\nbigo(-m)]\bigr)
=Q^{\circ}\bigl(m,\vecy\bigr)\times\proj_m$.
The pull back of the naturally defined reduced $L$-section
$[\phibar^u]$ on $\proj_2\times X$
is same as the naturally defined reduced $L$-section
on $\proj_1\times X$.

On the other hand,
the reduced $L$-morphism $[\phibar]$ induces
the section $f:U\lrarr \proj_1$.
It is easy to see that 
$[\phibar]$ is the pull back of
the naturally defined reduced $L$ section
on $\proj_1\times X$.
Thus we are done.
\hfill\qed

\vspace{.1in}
Since the above $G_m$-action
and the $GL(V_m)$-action
on $Q^{\circ}(m,\vecy,L)$ are commutative,
we have the induced
right $GL(V_m)$-action on $Q^{\circ}(m,\vecy,[L])$
and the universal object.
We also have the $GL(V_m)$-equivariant immersion
of $Q^{\circ}(m,\vecy,[L])$ to
$Z_m\times\prod_i G_{m,i}\times \proj(V_m^{\lor})$.

\subsubsection{Quotient reduced $\vecL$-Bradlow pair}

Let $\vecL=(L_1,L_2)$ be a pair of line bundles on $X$.
Quotient quasi-parabolic reduced $\vecL$-Bradlow pairs 
of $V_{m,X}$ with type $\vecy$ can be 
given as in Definition \ref{df;06.4.10.10}.
It is easy to construct 
the scheme $Q^{\circ}(m,\vecy,[\vecL])$
\index{$Q^{\circ}(m,\vecy,[\vecL])$}
representing the moduli functor
of quotient quasi-parabolic reduced $\vecL$-Bradlow pairs 
of $V_{m,X}$ with type $\vecy$
satisfying (TFV)-conditions.
In fact,
we have only to take the fiber product
of $Q^{\circ}(m,\vecy,[L_i])$ $(i=1,2)$
over $Q^{\circ}(m,\vecy)$.
We have the natural right $GL(V_m)$-action
on $Q^{\circ}(m,\vecy,[\vecL])$
and the $GL(V_m)$-equivariant immersion
$Q^{\circ}(m,\vecy,[\vecL])
\lrarr Z_m\times \prod_iG_{m,i}\times\proj_m\times\proj_m$.

\subsubsection{Oriented objects}
\label{subsubsection;06.5.4.20}

We put
$Q^{\circ}(m,\vecyhat,[L]):=
 Q^{\circ}(m,\vecy,[L])\times_{Q^{\circ}(m,y)}
 Q^{\circ}(m,\yhat)$, \index{$Q^{\circ}(m,\vecyhat,[L])$}
which represents the moduli functor
of quotient oriented quasi-parabolic
reduced $L$-Bradlow pairs 
of $V_{m,X}$ with type $\vecy$
satisfying (TFV)-condition.
We naturally have the universal object.
Similarly,
the schemes $Q^{\circ}(m,\vecyhat)$,
$Q^{\circ}(m,\vecyhat,L)$,
$Q^{\circ}(m,\vecyhat,[L_1],[L_2])$, etc.
are given,
which represent appropriate functors
respectively.

Let $\Zhat_m$ be as in (\ref{eq;06.5.21.30}).
We have the following naturally induced
Cartesian diagram:
\[
 \begin{CD}
 Q^{\circ}(m,\yhat) @>>> \Zhat_m\\
 @VVV @VVV \\
 Q^{\circ}(m,y)@>>> Z_m
 \end{CD}
\]
The morphisms are $\GL(V_m)$-equivariant.
Hence, we obtain the $\GL(V_m)$-equivariant morphism:
\[
 Q^{\circ}(m,\vecyhat)\lrarr \Zhat_m\times\prod_i G_{m,i}
\]

\subsubsection{Quotient stacks}

We have the universal quotient quasi-parabolic sheaf
$\bigl(q^u,\nbige^u,F_{\ast}\bigr)$ of
$p_{Q^{\circ}(m,\vecy)}V_{m,X}$
over $Q^{\circ}(m,\vecy)\times X$.
The $\GL(V)$-action on $Q^{\circ}(m,\vecy)$
is naturally lifted to the action
on $\bigl(q^u,\nbige^u,F_{\ast}\bigr)$.
Then, we obtain the quasi-parabolic sheaf
$(E^u,F_{\ast})$ on 
$\bigl(Q^{\circ}(m,\vecy)/\GL(V_m)\bigr)
\times (X,D)$
by taking the descent of
$\bigl(\nbige^u(-m),F_{\ast}\bigr)$.
The following lemma is well known.
\begin{lem}
Let $\vecy$ be an element of $\Type$.
The quotient stack $Q^{\circ}(m,\vecy)/\GL(V_m)$
is isomorphic to $\nbigm(m,\vecy)$,
and the quasi parabolic sheaf $(E^u,F_{\ast})$
is the universal one.
\end{lem}
\pf
Let $g:T\lrarr \nbigm(m,\vecy)$ be a morphism.
Then we have the corresponding $T$-quasi-parabolic sheaf
$(E,F_{\ast})$ on $T\times X$ of type $\vecy$,
satisfying the condition $O_m$.
Then we obtain the vector bundle
$\nbigv:=p_{X\,\ast}E(m)$ on $T$.
Let $P$ denote the associated principal $\GL(V_m)$-bundle,
and let $\pi:P\lrarr T$ denote the projection.
On $P$,
we have the equivariant trivialization 
$\pi_X^{\ast}\nbigv\simeq V_{m}\otimes\nbigo_T$.
Thus we obtain the equivariant morphism
$q:p_{P}^{\ast}V_{m,X}\lrarr \pi_X^{\ast}E(m)$.
Therefore, we obtain the quotient quasi-parabolic sheaf
$\bigl(q,\pi_X^{\ast}E(m),F_{\ast}\bigr)$ on $T\times X$
which is naturally $\GL(V_m)$-equivariant.
It also satisfies the (TFV)-condition.
Therefore, we obtain the $\GL(V_m)$-equivariant morphism
$P\lrarr Q^{\circ}(m,\vecy)$,
in other words,
the morphism $T\lrarr Q^{\circ}(m,\vecy)/\GL(V_m)$.
In particular,
we obtain $\nbigm(m,\vecy)\lrarr Q^{\circ}(m,\vecy)/\GL(V_m)$.

On the other hand, 
let $g:T\lrarr Q^{\circ}(m,\vecy)/\GL(V_m)$.
We have the corresponding $\GL(V_m)$-torsor $P(g)$
over $T$.
On $P(g)\times X$,
we have the quotient
$g_X^{\ast}(q^u):
 p_{T}^{\ast}V_{m,X}\lrarr g_X^{\ast}\nbige^u$
with a quasi parabolic structure $t_X^{\ast}F$,
which is $\GL(V_m)$-equivariant.
By taking the descent with respect to the action,
we obtain the $T$-quasi parabolic sheaf
$\bigl(\nbige^u(-m),F_{\ast}\bigr)$ on $T\times X$.
It satisfies the condition $O_m$.
Therefore, we obtain
the morphism of $T$ to $\nbigm(m,\vecy)$.
In particular, we obtain 
$Q^{\circ}(m,\vecy)/\GL(V_m)\lrarr\nbigm(m,\vecy)$.

It is easy to check that they are mutually inverse.
\hfill\qed

\vspace{.1in}
By the same argument,
we obtain the description of the moduli stacks
$\nbigm(m,\vecyhat)$,
$\nbigm(m,\vecy,L)$,
$\nbigm(m,\vecy,[L])$,
$\nbigm(m,\vecyhat,[L])$,
$\nbigm(m,\vecyhat,[\vecL])$, etc.
as the quotient stacks of
$Q^{\circ}(m,\vecyhat)$,
$Q^{\circ}(m,\vecy,L)$,
$Q^{\circ}(m,\vecy,[L])$,
$Q^{\circ}(m,\vecyhat,[L])$,
$Q^{\circ}(m,\vecyhat,[\vecL])$, etc..
The universal objects are constructed
by the same procedure.

\section{Geometric Invariant Theory and Enhanced Master Space}
\label{section;06.7.3.10}

\subsection{$\delta$-Semistability and Numerical Criterion}
\label{subsection;06.7.3.11}
\subsubsection{Statement}
\label{subsubsection;06.5.14.30}

Let $X$ be a smooth $d$-dimensional projective variety over 
an algebraically closed field $k$ of characteristic $0$,
and let $\nbigo_X(1)$ be a very ample line bundle.
We denote the number $c_1(\nbigo_X(1))^d\cap[X]$ by $g$.
Let $D$ denote a Cartier divisor of $X$.
We do not have to assume the smoothness of $D$
in this section.
We use the notation in the subsection 
\ref{subsection;06.5.15.15}.

Let $\vecy=(y,y_1,y_2,\ldots,y_l)$ be an element of $\Type$.
(See the subsubsection \ref{subsubsection;06.4.24.10}.)
We have the associated Hilbert polynomials
$H_{\vecy,\ast}=
 (H_{\vecy},H_{\vecy,1},H_{\vecy,2},\ldots,H_{\vecy,l})$.
(See the subsubsection \ref{subsubsection;06.5.15.1}.)
Let $V_m$ be an $H_{\vecy}(m)$-dimensional vector space over $k$,
and let $V_{m,X}$ denote $V_m\otimes\nbigo_X$.
Let $Z_m$ denote the Gieseker space
over $\Pic\bigl(\det(y(m))\bigr)$,
and let $G_{m,i}$ denote the Grassmann variety
of $H_{\vecy}(m)-H_{\vecy,i}(m)$-dimensional
quotients of $V_m$.
We put as follows:
\[
 \nbiga_m(\vecy):=
 Z_m\times\prod_i G_{m,i},
\quad
\nbiga_m(\vecy,[L]):=
\nbiga_m(\vecy)\times \proj_m,
\quad
 \nbiga_m(\vecy,[\vecL]):=
 \nbiga_m(\vecy)\times\proj_m^{(1)}\times\proj_m^{(2)}.
\]
Here $L$ denotes a line bundle on $X$,
and $\vecL=(L_1,L_2)$ denotes a pair of line bundles on $X$.
We also put $\proj_m^{(i)}=\proj_m$.
Recall that we have obtained the $SL(V_m)$-equivariant
immersions $\Psi_m$ of
$Q^{\circ}(m,\vecy)$ to $\nbiga_m(\vecy)$.
(See the subsubsection \ref{subsubsection;06.5.14.1}).
We take an inclusion $\iota:\nbigo(-m)\lrarr L$.
Then we have the closed immersion
$Q^{\circ}(m,\vecy,[L])\lrarr
Q^{\circ}\bigl(m,\vecy,[\nbigo(-m)]\bigr)
=Q^{\circ}(m,\vecy)\times\proj_m$.
Therefore,
we obtain the equivariant immersion
of $Q^{\circ}(m,\vecy,[L])$
to $\nbiga_m(\vecy,[L])$.
Similarly, we can obtain the equivariant immersion
of $Q^{\circ}(m,\vecy,[\vecL])$
to $\nbiga_m(\vecy,[\vecL])$
by taking an inclusion $\iota_i:\nbigo(-m)\lrarr L_i$.
The immersions are denoted by $\Psi_m$.

Since $Z_m$ is the projective space bundle on 
$\Pic_X\bigl(\det(y(m))\bigr)$,
we have the relative tautological bundle
$\nbigo_{Z_m}(1)$.
We also have the canonical polarizations
$\nbigo_{G_{m,i}}(1)$ and $\nbigo_{\proj_m}(1)$
of $G_{m,i}$ and $\proj_m$ respectively.
They are $SL(V_m)$-equivariant line bundles.

For a positive number $A$
and a tuple of positive numbers
$B_{\ast}=(B_i\,|\,i=1,\ldots,l)$,
we formally consider the line bundles
on $\nbiga_m(\vecy)$ given as follows,
although it is not precisely a line bundle
when the numbers are not integers:
\[
 \nbigl_{\vecy}(A,B_{\ast}):=
 \nbigo_{Z_m}(A)\otimes
 \bigotimes_{i=1}^l\nbigo_{G_{m,i}}(B_i).
\]
Similarly, for positive numbers $C$ and $C_j$ $(j=1,2)$,
we formally consider
$ \nbigl_{\vecy,L}(A,B_{\ast},C):=
\nbigl_{\vecy}(A,B_{\ast})
\otimes
 \nbigo_{\proj_m}(C) $
on $\nbiga_m(\vecy,L)$
and
$ \nbigl_{\vecy,\vecL}(A,B_{\ast},C_1,C_2):=
\nbigl_{\vecy}(A,B_{\ast})
\otimes
 \nbigo_{\proj_m^{(1)}}(C_1)
\otimes
 \nbigo_{\proj_m^{(2)}}(C_2)$
on $\nbiga_m(\vecy,\vecL)$.

When the numbers are positive rational numbers,
$\nbigl_{\vecy}(A,B_{\ast})$
gives $SL(V_m)$-equivariant 
$\rnum$-polarizations of $\nbiga_m(\vecy)$.
Even if the numbers are not rational,
the Hilbert-Mumford criterion formally
provides us the semistability condition
with respect to $\nbigl_{\vecy}(A,B_{\ast})$.
So, let $\nbiga^{ss}_m(\vecy,A,B_{\ast})$
denote the set of semistable points 
of the $SL(V_m)$-action on $\nbiga_m(\vecy)$
with respect to $\nbigl_{\vecy}(A,B_{\ast})$.
Similarly,
$\nbiga^{ss}_m(\vecy,L,A,B_{\ast},C)$
and $\nbiga^{ss}_m(\vecy,\vecL,A,B_{\ast},C_1,C_2)$
are given.

\vspace{.1in}

The semistability (resp. stability) condition
with respect to the system of weights $\alpha_{\ast}$
determines the open subset
$Q^{ss}(m,\vecy,\alpha_{\ast})$
(resp. $Q^{s}(m,\vecy,\alpha_{\ast})$)
\index{$Q^{ss}(m,\vecy,\alpha_{\ast})$,
 $Q^s(m,\vecy,\alpha_{\ast})$}
of $Q^{\circ}(m,\vecy)$.
Namely, it is the maximal open subset of
$Q^{\circ}(m,\vecy)$
which consists of the points $(q,\nbige,F_{\ast})$
such that the parabolic sheaf
$(\nbige(-m),F_{\ast},\alpha_{\ast})$
is semistable.
Similarly,
we have the 
open subset $\Qmss$
(resp. $Q^{s}(m,\vecy,[L],\alpha_{\ast},\delta)$)
\index{$Q^{ss}(m,\vecy,[L],\alpha_{\ast},\delta)$,
 $Q^{s}(m,\vecy,[L],\alpha_{\ast},\delta)$}
of $Q^{\circ}(m,\vecy,[L])$
determined by
the semistability (stability) condition with respect to
the system of weights $\alpha_{\ast}$
and a parameter $\delta$.
We also have the open subsets
$Q^{ss}(m,\vecy,[\vecL],\alpha_{\ast},\vecdelta)$
and $Q^{s}(m,\vecy,[\vecL],\alpha_{\ast},\vecdelta)$
of $Q^{\circ}(m,\vecy,[\vecL])$
in a similar way.
\index{$Q^{ss}(m,\vecy,[\vecL],\alpha_{\ast},\vecdelta)$,
 $Q^{s}(m,\vecy,[\vecL],\alpha_{\ast},\vecdelta)$}

The first claim of the following proposition
was proved by Maruyama-Yokogawa \cite{my}
and the second claim was proved by Yokogawa in \cite{y}.
\begin{prop}
\label{prop;06.5.21.31}
There exists an integer
$N(\vecy,\alpha_{\ast})$ such that
the following holds for any $m\geq N(\vecy,\alpha_{\ast})$:
\begin{enumerate}
\item
The image of $Q^{ss}(m,\vecy,\alpha_{\ast})$
via the morphism $\Psi_m$
is contained in
$\nbiga^{ss}\bigl(m,\vecy,P_{\vecy}^{\alpha_{\ast}}(m),
 \epsilon_{\ast}\bigr)$.
Thus we obtain the morphism
$\Psihat_m:
Q^{ss}(m,\vecy,\alpha_{\ast})\lrarr
\nbiga^{ss}\bigl(m,\vecy,P_{\vecy}^{\alpha_{\ast}}(m),\epsilon_{\ast}\bigr)$.
\item
The morphism $\Psihat_m$ above is proper.
 In particular, it is a closed immersion.
\hfill\qed
\end{enumerate}
\end{prop}

By the same argument,
we can show the following proposition.

\begin{prop}
 \label{prop;06.4.11.50}
There exists an integer $N_1(\vecy,L,\alpha_{\ast},\delta)$
such that the following claims hold
for any $m\geq N_1(\vecy,L,\alpha_{\ast},\delta)$:
\begin{enumerate}
\item
The image of $\Qmss$
via the morphism $\Psi_m$
is contained in 
$\Amss$.
Thus we obtain the morphism
$ \Psihat_m:
\Qmss
\lrarr
\Amss$.
\item
The morphism $\Psihat_m$ above is proper.
 In particular, it is a closed immersion.
\end{enumerate}
Similarly, there exists a large integer
$N_1(\vecy,\vecL,\alpha_{\ast},\vecdelta)$
such that we have the $SL(V_m)$-equivariant
closed immersion for any $m\geq N_1(\vecy,\vecL,\alpha_{\ast},\vecdelta)$:
\[
 \Psihat_m:
 Q^{ss}(m,\vecy,[\vecL],\alpha_{\ast},\vecdelta)
\lrarr
 \nbiga_m^{ss}\bigl(\vecy,[\vecL],
 P_{\vecy}^{\alpha_{\ast},\vecdelta}(m),\epsilon_{\ast},
 \vecdelta(m) \bigr).
\]
Here, we put $\vecdelta(m):=(\delta_1(m),\delta_2(m))$.
\end{prop}

Although we need only a minor modification,
we will later give a rather detailed proof of the claims
for $(\vecy,L,\alpha_{\ast},\delta)$
in Proposition \ref{prop;06.4.11.50},
for the convenience of the reader.
We closely follow the arguments of \cite{my} and \cite{y}.
We also use the argument in \cite{hl2}.
Since the claim for $(\vecy,\vecL,\alpha_{\ast},\vecdelta)$
can be shown similarly, we omit to give the proof of it.

We also obtain the following.
\begin{prop}
 \label{prop;06.4.11.85}
The image of 
$Q^{s}(m,\vecy,[L],\alpha_{\ast},\delta)$
via $\Psi_m$ is contained
in $\nbiga^s_m(\vecy,[L],\alpha_{\ast},\delta)$.
Similar claims hold for
$Q^s(m,\vecy,\alpha_{\ast})$
and $Q^s(m,\vecy,[\vecL],\alpha_{\ast},\vecdelta)$.
\end{prop}
The claim for $Q^{s}(m,\vecy)$ was proved 
by Maruyama-Yokogawa in \cite{my}.
Since the argument is similar to the proof of
the first claim of Proposition \ref{prop;06.4.11.50},
we give just some remarks in the proof of 
Proposition \ref{prop;06.4.11.50}.

\begin{rem}
When $\alpha_i$  and the coefficients of $\delta$ and 
$\delta_i$ $(i=1,2)$ are rational,
it is standard to obtain the projective coarse moduli scheme
of $\delta$-semistable parabolic $L$-Bradlow pairs
or $\vecdelta$-semistable 
parabolic $\vecL$-Bradlow pairs
from Proposition {\rm\ref{prop;06.4.11.50}}.

Even if the numbers are not rational,
we can obtain the coarse moduli schemes
of $\delta$-stable parabolic reduced $L$-Bradlow pairs,
if the $1$-stability condition holds
for $(\vecy,L,\alpha_{\ast},\delta)$.
Take a sufficiently large $N>0$,
and take rational numbers
$A$, $B_{\ast}=(B_i\,|\,i=1,\ldots,l)$
and $C$
which are close to 
$N\cdot P_{\vecy}^{\delta,\alpha_{\ast}}(m)$,
$N\cdot \epsilon_{\ast}$
and $N\cdot \delta(m)$ respectively.
If a point of $\nbiga(m,\vecy,[L])$
is stable with respect to
$\nbigl_{\vecy,L}\bigl(P_{\vecy}^{\delta,\alpha_{\ast}}(m),
 \epsilon_{\ast},\delta(m) \bigr)$,
then it is stable with respect to
$\nbigl_{\vecy,L}(A,B_{\ast},C)$.
Thus we can obtain the coarse scheme
by the geometric invariant theory
(Proposition {\rm\ref{prop;06.6.21.10}}).
Therefore,
when the $1$-stability condition holds,
we obtain the projective coarse moduli scheme
of semistable ones.
\hfill\qed
\end{rem}

Before going into the proof of 
Proposition \ref{prop;06.4.11.50},
we give some consequence about the property
of the moduli stacks.

\begin{cor}
When the $1$-stability condition holds
for $(\vecy,\alpha_{\ast})$,
the moduli stack $\nbigm^s(\vecyhat,\alpha_{\ast})$
is Deligne-Mumford and proper.
\end{cor}
\pf
Under the assumption,
we obtain
$Q^{ss}(\vecy,\alpha_{\ast})=
Q^{s}(\vecy,\alpha_{\ast})$.
It is easy to see the finiteness of the stabilizer of
any point $z\in Q^{ss}(\vecy,\alpha_{\ast})$
with respect to the $\SL(V_m)$-action
from Corollary \ref{cor;06.5.21.35}.
Then, the quotient stack
$Q^{s}(\vecy,\alpha_{\ast})/\SL(V_m)$
is Deligne-Mumford and proper,
due to Proposition \ref{prop;06.5.21.31},
Proposition \ref{prop;06.4.11.85}
and Proposition \ref{prop;06.4.25.1}.

Since we have the etale finite morphisms
$Q^s(\vecy,\alpha_{\ast})/\SL(V_m)
\lrarr Q^{s}(\vecy,\alpha_{\ast})/\PGL(V_m)$
and 
$\nbigm(\vecyhat,\alpha_{\ast})
\lrarr Q^{s}(\vecy,\alpha_{\ast})/\PGL(V_m)$,
it is easy to derive the claimed property
of $\nbigm(\vecyhat,\alpha_{\ast})$.
\hfill\qed

\vspace{.1in}
Similarly, we obtain the following.
\begin{prop}
\mbox{{}}
\begin{itemize}
\item
If the $1$-stability condition holds for
$(\vecy,L,\alpha_{\ast},\delta)$,
the moduli stacks
$\nbigm^s(\vecyhat,[L],\alpha_{\ast},\delta)$
and $\nbigm^s(\vecy,L,\alpha_{\ast},\delta)$
are Deligne-Mumford and proper.
\item
If the $1$-stability condition holds for
$(\vecy,\vecL,\alpha_{\ast},\vecdelta)$,
the moduli stack
$\nbigm^{s}(\vecyhat,[\vecL],\alpha_{\ast},\vecdelta)$
is Deligne-Mumford and proper.
\hfill\qed
\end{itemize}
\end{prop}

\subsubsection{The numerical criterion}

Let us start the proof of Proposition \ref{prop;06.4.11.50}.
Let $(q,E_{\ast},\phi,W_{\ast},[\phitilde])$
be a Yokogawa datum of type $(\vecy,m)$
(Definition \ref{df;06.6.21.11}).
From $W_{\ast}$,
we obtain the tuple of the quotients
$\bigl(V_m/W_i\,|\,i=1,\ldots,l\bigr)$,
which gives the point of $\prod_{i}G_{m,i}$.
We denote the point by $V/W_{\ast}$.
Then we obtain the following element
of $\nbiga_m(\vecy,[L])$:
\begin{equation}
\label{eq;06.5.21.50}
 \Psi\bigl(q,E_{\ast},\phi,W_{\ast},[\phitilde]\bigr):=
 \left(
 H^0\Bigl(\bigwedge^r q\Bigr),\,
 V_m/W_{\ast},
 [\phitilde]\right)
\in\nbiga_m(\vecy,[L]).
\end{equation}

Let $W$ be a subspace of $V_m$.
The number $\epsilon(W,[\phitilde])$
is given as in Definition \ref{df;06.5.15.20}.
Let $\nbige_{W}$ denote the subsheaf of $E(m)$
generated by the image of $W$ via $q$.

\begin{lem} 
 \label{lem;06.4.11.35}
The point $\Psi\bigl(q,E_{\ast},\phi,W_{\ast},[\phitilde]\bigr)$
is contained in $\Amss$,
if and only if the following inequalities hold
for any non-trivial subspace $W\subset V_m$:
\begin{equation} 
 \label{eq;06.4.11.3}
 P_{\vecy}^{\alpha_{\ast},\delta}(m)\cdot\rank(\nbige_W)
-\epsilon(W,[\phi])\cdot\delta(m)
-\sum \epsilon_i \cdot\dim(W_i\cap W)
-\alpha_1\cdot\dim(W)\geq 0.
\end{equation}
In other words,
$\Psi(q,E_{\ast},\phi,W_{\ast},[\phitilde])$
is semistable point,
if and only if 
$(q,E_{\ast},\phi,W_{\ast},[\phitilde])
\in \YOKtilde(m,0,\vecy,L,\delta)$.

The point $\Psi\bigl(q,E_{\ast},\phi,W_{\ast},[\phitilde]\bigr)$
is contained in $\Ams$,
if and only if the strict inequalities hold
in {\rm(\ref{eq;06.4.11.3})}
for any subspace $W\subset V_m$.
\end{lem}
\pf
We give only an indication.
We put
$\nbigl_1:=
 \nbigo_{Z}\bigl(P_{\vecy}^{\alpha_{\ast}}(m)\bigr)\otimes
 \bigotimes_{i=1}^l
 \nbigo_{G_{m\,i}}(\epsilon_{i})$
and 
$ \nbigl_{2}:=
 \nbigo_{Z}\bigl(r^{-1}\cdot\delta(m)\bigr)\otimes
 \nbigo_{\proj_m}(\delta(m))$.
Then we have
$\nbigl\bigl(P^{\alpha_{\ast},\delta}_{\vecy}(m), 
 \epsilon_{\ast},\delta(m)\bigr)=\nbigl_1\otimes\nbigl_2$.

We put $N=\dim V_m$.
Let $u_1,\ldots,u_N$ be the vectors
satisfying the following:
\begin{equation}
 \label{eq;06.4.22.6}
\lambda(t)\cdot u_i=t^{w_i}\cdot u_i,
\quad
w_1\leq w_2\leq\cdots\leq w_N,
\quad
\sum w_i=0.
\end{equation}
We remark that
we can regard {\rm(\ref{eq;06.4.22.6})}
as the condition for $\lambda$,
when we fix $u_1,\ldots,u_N$.

The number $\mu_{\lambda}(P,\nbigl_1)$
is calculated in \cite{my} (see also \cite{bhosle}):
\begin{multline}
 \label{eq;06.4.22.10}
 \mu_{\lambda}(P,\nbigl_1)=
-P_{\vecy}^{\alpha_{\ast}}(m)
\sum_{i=1}^N
\bigl(
 \rank\nbige^{(i)}-\rank\nbige^{(i-1)}
\bigr) w_i 
-\sum_{j=1}^l
 \epsilon_j
 \sum_{i=1}^N
 \Bigl(
 \rank W_{j}\cap \nbigv^{(i-1)}
-\rank W_{j}\cap \nbigv^{(i)}+1\Bigr)\cdot
 w_i.
\end{multline}
Here, $\nbige^{(i)}$ denote the subsheaves
of $E(m)$
generated by $u_1,\ldots,u_i$ via $q$,
and $\nbigv^{(i)}$ denote the subspaces
of $V_m$ generated by $u_1,\ldots,u_i$.
The number $\mu_{\lambda}(P,\nbigl_2)$
is as follows:
\begin{equation}
 \label{eq;06.4.22.11}
 \delta(m)\cdot\left(
-\frac{1}{r}\sum_{i=1}^N
\bigl(
 \rank\nbige^{(i)}-\rank\nbige^{(i-1)}
\bigr)\cdot w_i
+\sum_{i=1}^N
 \Bigl(
 \dim \langle \phitilde\rangle\cap \nbigv^{(i)}
-\dim\langle\phitilde\rangle\cap\nbigv^{(i-1)}
 \Bigr)\cdot w_i
\right)
\end{equation}
Here, $\langle\phitilde\rangle$ denotes
the subspace of $V_m$ generated by $\phitilde$.
The number $\mu_{\lambda}(P,\nbigl)$
can be obtained as the sum of
(\ref{eq;06.4.22.10}) and (\ref{eq;06.4.22.11}).
We write it as the reference for the later argument:
\begin{multline}
 \label{eq;06.4.22.5}
\mu_{\lambda}(P,\nbigl)=
-P_{\vecy}^{\alpha_{\ast}}(m)\cdot
\sum_{i=1}^N
\bigl(
 \rank\nbige^{(i)}-\rank\nbige^{(i-1)}
\bigr)\cdot w_i
-\sum_{j=1}^l
 \epsilon_j
 \sum_{i=1}^N
 \Bigl(
 \rank W_{j}\cap \nbigv^{(i-1)}
-\rank W_{j}\cap \nbigv^{(i)}+1\Bigr)
 \cdot w_i
 \\
+ \delta(m)
\left(
-\frac{1}{r}\sum_{i=1}^N
 \bigl(\rank\nbige^{(i)}-\rank\nbige^{(i-1)}\bigr)\cdot w_i
+\sum_{i=1}^N
 \bigl(
 \dim \langle \phitilde\rangle\cap \nbigv^{(i)}
-\dim\langle\phitilde\rangle\cap\nbigv^{(i-1)}
 \bigr)\cdot w_i
\right).
\end{multline}
The right hand side of  (\ref{eq;06.4.22.5})
is linear with respect to $(w_1,\ldots,w_N)$.
Therefore,
$\mu_{\lambda}(P,\nbigl)\geq 0$
holds for any $\lambda$
satisfying (\ref{eq;06.4.22.6}),
if and only if 
$\mu_{f_k}(P,\nbigl)\geq 0$
for any $k=1,\ldots,N-1$,
where $f_k=(\overbrace{k-N,\ldots,k-N}^k,k\ldots,k)$.
Due to the calculation of Maruyama and Yokogawa,
we have the following:
\begin{equation}
\label{eq;06.5.21.150}
 \mu_{f_k}(P,\nbigl_1)
=
 H(m)\cdot
\left(
 P_{\vecy}^{\alpha_{\ast}}(m)\cdot \rank \nbige^{(k)}
-\sum_{i=1}^l\epsilon_i\cdot \dim W_i\cap \nbigv^{(k)}
-\alpha_1\cdot\dim \nbigv^{(k)}.
\right)
\end{equation}
By a direct calculation,
we have the following:
\begin{equation}
 \label{eq;06.5.21.151}
 \mu_{f_k}(P,\nbigl_2)
=H(m)\cdot\delta(m)\cdot
 \left(
 \frac{\rank \nbige^{(k)}}{r}
-\epsilon\bigl(\nbigv^{(k)},[\phitilde]\bigr)
 \right).
\end{equation}
Hence $\mu_{f_k}(P,\nbigl)$
depends only on $\nbigv^{(k)}$.
Therefore, we obtain 
the function 
$F_{(P,\nbigl)}$ on 
the sets of the non-trivial subspaces
$\{0\neq W\subsetneq V\}$,
and 
$P$ is semistable with respect to
$\nbigl$
if and only if $F_{(P,\nbigl)}(W)\geq 0$
for any subspace $W$.
Since
$F_{(P,\nbigl)}(W)$
is the left hand side of 
(\ref{eq;06.4.11.3})
multiplied with the positive number $H(m)$,
we are done.
\hfill\qed

\subsubsection{A lemma}

Following Huybrecht-Lehn \cite{hl2},
we show the following lemma.
\begin{lem}
\label{lem;06.4.11.90}
Let $\vecy$ be an element of $\Type$.
There exists an integer $N_4$ with the following property:
\begin{itemize}
\item
Let $(E_{\ast},\phi)$ be a $\delta$-semistable
parabolic $L$-Bradlow pair of type $\vecy$
such that $\phi\neq 0$.
Then the following inequality holds
for any integer $m\geq N_4$
and for any subsheaf $F\subset E$:
\begin{equation}
 \label{eq;06.4.11.80}
 \frac{h^0\bigl(F_{\ast}(m)\bigr)
   +\epsilon(\phi,F)\cdot\delta(m)
  }
 {\rank (F)}
\leq
 P^{\delta}_{(E_{\ast},\phi)}(m).
\end{equation}
\item
 If the equality holds in {\rm (\ref{eq;06.4.11.80})},
 $(F_{\ast},\phi')$ is $\delta$-semistable
 with $P^{\delta}_{(F_{\ast},\phi')}=P^{\delta}_{(E_{\ast},\phi)}$.
\item
If $(E_{\ast},\phi)$ is $\delta$-stable,
the strict inequality holds in {\rm(\ref{eq;06.4.11.80})}.
\end{itemize}
\end{lem}
\pf
We have the inequality
$\mu(F)\leq \mu^{\delta}(F_{\ast})
\leq\mu^{\delta}(E_{\ast})$
for any subsheaf $F\subset E$,
and hence
$\mu_{\max}(F)\leq 
 \mu(E_{\ast})+\delta_{top}=C_0$.
On the other hand, we have the inequality
$\mu_{\min}(F)\leq \mu(F)$ by definition.
Hence we obtain the following inequality
by using Proposition \ref{prop;06.5.21.3}
(\cite{hl2}):
\[
 \frac{h^0(F(m))}{\rank(F)}
 \leq
 \frac{1}{g^{d-1}d!}
 \left(
  \Bigl(1-\frac{1}{\rank(F)}\Bigr)
 \bigl[C_0+mg+c\bigr]^d_+
  +\frac{1}{\rank(F)}\bigl[\mu(F)+mg+c\bigr]_+^d
 \right).
\]

We can take a sufficiently negative number $C$
with the following property:
\begin{itemize}
\item
For any positive integer $r'$ such that $r'<\rank(E)$
and for any sufficiently large $t>0$,
the following inequalities hold:
\begin{equation}
 \label{eq;06.4.11.30}
 \frac{1}{g^{d-1}d!}
\left( 
\Bigl(1-\frac{1}{r^{\prime}}\Bigr)(C_0+tg+c)^d
+\frac{1}{r^{\prime}}(C+tg+c)^d 
 \right)
+\frac{\delta(t)}{r'}
\leq
 P^{\alpha_{\ast}}_{\vecy}(t).
\end{equation}
\end{itemize}
Note that the coefficient of $t^d$ of the both sides
are same by the construction,
and thus we can take such $C$.
Let $N_5$ be a large number
such that $-C+mg+c>0$ for any $m>N_5$.

Let $\nbigs$ be the family of the sheaves $F$
with the following property:
\begin{itemize}
\item
There exists
a $\delta$-semistable parabolic reduced $L$-Bradlow pairs
$(E_{\ast},\phi)$ of type $\vecy$ with
weight $\alpha_{\ast}$
such that $F$ is a saturated subsheaf of $E$.
\end{itemize}

We divide $\nbigs$ into two families
by the following conditions:
(i) $\mu(F)<C$ or (ii) $\mu(F)\geq C$.
If $F$ is contained in the family (i),
then the desired inequality holds 
for any $m\geq N_5$ because of our choices of $C$ and $N_5$.
On the other hand,
the family (ii) is bounded (Proposition \ref{prop;06.5.21.1}).
Thus we can take a large number $N_6$
such that the family (ii) satisfies the condition $(O_m)$
for any $m\geq N_6$.
Then the desired inequalities hold
for any $m\geq N_6$,
because of the $\delta$-semistability of
$(E_{\ast},\phi)$.
Thus we have only to put
$N_4:=\max\{N_5,N_6\}$.
\hfill\qed

\subsubsection{Proof of the claim 1 in Proposition \ref{prop;06.4.11.50}}

Let $N_1(\vecy,L,\alpha_{\ast},\delta)$
be an integer larger than $N_4$
in Lemma \ref{lem;06.4.11.90}
and $N_0(\vecy,L,\delta)$ in Proposition \ref{prop;06.5.15.10}.
Let $\bigl(q,E(m),F_{\ast},[\phitilde]\bigr)$ be a point of $\Qmss$.
The reduced $L$-section $[\phi]$ of $E$ is induced by
$[\phitilde]$ and $\iota$.
By definition,
$H^0(q)$ gives the isomorphism
$V_m\simeq H^0\bigl(X,E(m)\bigr)$.
Let $W$ be a subspace of $V_m$,
which generates the subsheaf $\nbige_W$ of $E(m)$
via $q$.
Due to Lemma  \ref{lem;06.4.11.90},
we have the following inequality:
\begin{equation}
\label{eq;06.4.11.55}
 \frac{h^0(\nbige_{W,\ast})
 +\epsilon\bigl(\nbige_W(-m),\phi\bigr)\cdot \delta(m) }
 {\rank(\nbige_W)}
 \leq
 P_{\vecy}^{\delta,\alpha_{\ast}}(m).
\end{equation}
We have the equalities
$\epsilon([\phi],W)=\epsilon([\phi],\nbige_W)$
and the inequalities
$\dim(W)\leq h^0(\nbige_{W})$ and
$\dim(W\cap W_i)\leq h^0(F_i(\nbige_W))$.
Thus we obtain the following inequality:
\begin{equation} 
\label{eq;06.4.11.60}
 \alpha_1\dim(W)+\sum \epsilon_i\dim(W_i\cap W)
\leq
 \alpha_1\cdot h^0(\nbige_W)+
 \sum\epsilon_i \cdot h^0(F_i(\nbige_W))
=h^0\bigl(\nbige_{W\ast}\bigr).
\end{equation}
By substituting (\ref{eq;06.4.11.60})
to (\ref{eq;06.4.11.55}),
we obtain the inequality (\ref{eq;06.4.11.3}).
Hence $\Psi_m(q,E(m),F_{\ast},[\phitilde])$ gives 
the point contained in $\Amss$.
\hfill\qed

\begin{rem}
We obtain Proposition {\rm\ref{prop;06.4.11.85}}
by using the same argument and
the second claims in Lemma {\rm\ref{lem;06.4.11.35}}
and Lemma {\rm\ref{lem;06.4.11.90}}.
\hfill\qed
\end{rem}

\subsubsection{Proof of the claim 2
of Proposition \ref{prop;06.4.11.50}}
\label{subsubsection;06.5.15.80}

Take a discrete valuation ring $R$ over $k$.
We denote the quotient field by $K$.
Assume that we have the following diagram:
\begin{equation}\label{eq;6.19.32}
 \begin{CD}
 \Spec(K) @>{f}>>
 \Qmss
 \\
 @VVV @V{\widehat{\Phi}}VV \\
 \Spec(R) @>{g}>> 
 \Amss
 \end{CD}
\end{equation}
We have only to show the existence of
a lift $\Spec(R)\lrarr Q^{ss}(m,\vecy,[L],\alpha_{\ast},\delta)$.
Let $X_K$ and $X_R$ denote
$X\times\Spec K$ and $X\times \Spec R$
respectively.

The morphism $f$ corresponds to the tuple 
$(q_K,E_{K\ast }(m),[\phitilde_K])$
of the quotient parabolic sheaf
$\bigl(q_K,E_{K\ast }(m)\bigr)$
and the quotient reduced $L$-section $[\phitilde_K]$
defined over $X_K$.
The tuple $(q_K,E_{K\ast}(m))$ satisfies the (TFV)-condition.
The $L$-section $[\phi_K]$ of $E_K$ is induced 
by $[\phitilde_K]$ and $\iota$,
and the parabolic $L$-Bradlow pair
$(E_{K\,\ast},\phi_K)$ with weight $\alpha_{\ast}$
is $\delta$-semistable.

As in \cite{y} pp. 502--503,
$E_{K\ast}(m)$ can be extended
to the parabolic torsion-free sheaf $E_{R\ast}(m)$
over $X_R$.
The morphism $q_K$ can be extended to the morphism
$q_R:V_{m}\otimes\nbigo_{X_R}\lrarr E_{R}(m)$
such that the restriction of $q_R$
to the closed fiber is generically surjective.
Since $\proj_m$ is proper,
we can extend $[\phitilde_K]$ to $[\phitilde_R]$ over $R$.
We put
$W_{K,i}:=H^0\bigl(X\otimes K,F_{i+1}(E(m))\bigr)$ $(i=1,\ldots,l)$
which give the subspaces of $V_{m}\otimes K$.
Since $\prod_i G_{m,i}$ are proper,
we obtain the subbundle
$W_{R,i}$ of $V_m\otimes R$ over $\Spec R$.
We put $W_{R,\ast}=(W_{R,i}\,|\,i=1,\ldots,l)$.
The family
$\bigl(q_R,E_{R\,\ast},W_{R\,\ast},[\phitilde_R]\bigr)$
induces the morphism $\Spec(R)\lrarr \nbiga_m(\vecy,[L])$
as in (\ref{eq;06.5.21.50}).
By separatedness of $\nbiga_m(\vecy,[L])$,
it coincides with $g$ in the diagram (\ref{eq;6.19.32}).

Let $\bigl(q_0,E_{0\ast},W_{0\ast},[\phitilde_0]\bigr)$
denote the specialization of
 $\bigl(q_R,E_{R\ast},W_{R\ast}, [\phitilde_R] \bigr)$ 
to the closed point of $\Spec R$.
We also have the induced $L$-section $\phi_0$ of $E_0$.
The tuple
$\bigl(q_0,E_{0\ast},\phi_0,W_{0\ast},[\phitilde_0]\bigr)$
is the Yokogawa datum such that 
$\Psi\bigl(q_0,E_{0\,\ast},\phi_0,W_{0\,\ast},[\phitilde_0]\bigr)$
is contained in $\Amss$.
Due to Lemma \ref{lem;06.4.11.35},
$\bigl(q_0,E_{0\ast},\phi_0,W_{0\ast},[\phitilde_0]\bigr)$
is contained in 
$\YOKtilde(m,0,\vecy,L,\delta)$.
Recall $m\geq N_1(\vecy,L,\alpha_{\ast},\delta)
 \geq N_0(\vecy,L,\delta)$,
where $N_0(\vecy,L,\delta)$ is as in 
Proposition \ref{prop;06.5.15.10}.
Hence, we know that $q_0$ is surjective,
that $E_{0,\ast}(m)$ satisfies the $TFV$-condition,
and that the parabolic $L$-Bradlow pair
$(E_{0,\ast},\phi_0)$ is $\delta$-semistable.
Hence 
$\bigl(q_R,E_{R\,\ast},W_{R\ast},[\phi_R]\bigr)$
gives a map $f_R:\Spec(R)\lrarr
 Q^{ss}\bigl(m,\vecy,[L],\alpha_{\ast},\delta \bigr)$,
whose restriction to $\Spec(K)$ is $f$.
It is clear that $f_R$ gives the lift of $g$.
Hence the claim $2$ of Proposition \ref{prop;06.4.11.50}
is proved.
\hfill\qed

\subsubsection{Complement}

We give a consequence of the proof.
We put $V:=V_m$,
$Q:=Q^{ss}(m,\vecy,[L],\alpha_{\ast},\delta)$.
We also put as follows:
\[
 \nbigl:=\nbigl_{\vecy,L}(P^{\delta,\alpha_{\ast}}_{\vecy}(m),
 \epsilon_{\ast},\delta(m))
\]

\begin{lem}
 \label{lem;06.4.24.65}
Let $z=(q,E,F_{\ast},[\phi])$ be a point of $Q$.
Let $V=V'\oplus V''$ be a decomposition,
and let $\lambda$ be the one-parameter subgroup
of $\SL(V)$ given by
$t^{-\rank V''}\cdot \id_{V'}
\oplus t^{\rank V'}\cdot \id_{V''}$.
Let $E'(m)$ denote the subsheaf generated by $V'$.
We have the induced parabolic structure
and the $L$-Bradlow pair $\phi'$ of $E'$.
Then $(E'_{\ast},\phi')$ is $\delta$-semistable
with $P^{\delta}_{(E'_{\ast},\phi')}=
   P^{\alpha_{\ast},\delta}_{\vecy}$,
if and only if $\mu_{\lambda}(z,\nbigl)=0$ hold.
\end{lem}
\pf
Assume $\mu_{\lambda}(z,\nbigl)=0$.
We put $W:=H^0\bigl(X,E'(m)\bigr)$.
Let $W_i$ denote the kernel 
of $H^0\bigl(X,E(m)\bigr)\lrarr H^0\bigl(X,\Cok_{i+1}(m)\bigr)$.
From the calculation in the proof of Lemma 
\ref{lem;06.4.11.35},
we have the following:
\begin{equation}
 \label{eq;06.4.24.50}
0=\mu_{\lambda}(z,\nbigl)
=H(m)\cdot
\left(
 P_{\vecy}^{\alpha_{\ast},\delta}(m)\cdot\rank E'
-\epsilon(E',\phi')\cdot\delta(m)
-\sum \epsilon_i \cdot\dim(W_i\cap W)
-\alpha_1\cdot\dim(W)
\right).
\end{equation}
Let $\nbige_{W}$ denote the subsheaf of $E(m)$
generated by $W$ via $q$.
Therefore, we obtain the following inequality
from (\ref{eq;06.4.24.50}):
\begin{multline}
 \label{eq;06.4.24.60}
 P^{\alpha_{\ast},\delta}_{\vecy}(m)
=\frac{\alpha_1\cdot\dim(W)+\sum\epsilon_i\cdot \dim (W_i\cap
W)+\epsilon(E',\phi')\cdot
 \delta(m)}{\rank E'}
\leq
 \frac{h^0(\nbige_{W\,\ast})+\epsilon(\nbige_{W},\phi)\cdot\delta(m)}
 {\rank \nbige_W} \\
\leq
 \frac{h^0(E'_{\ast}(m))+\epsilon(E',\phi)\cdot \delta(m)}
 {\rank E'}
\leq P^{\alpha_{\ast},\delta}_{\vecy}(m)
\end{multline}
Here, the first inequality is obtained
in (\ref{eq;06.4.11.60}),
the second inequality follows from
$\nbige_W\subset E'(m)$,
and the third inequality follows from 
Lemma \ref{lem;06.4.11.90}.
We can conclude that the equality holds
in (\ref{eq;06.4.24.60}).
Then $(E'_{\ast},\phi')$ is $\delta$-semistable
due to the second claim of Lemma \ref{lem;06.4.11.90}.

Assume $(E'_{\ast},\phi')$ is $\delta$-semistable.
Note that the condition $O_m$ holds for
$(E'_{\ast},\phi')$,
because it holds for 
$(E'_{\ast},\phi')\oplus \bigl((E/E')_{\ast},\phi''\bigr)$,
where $\phi''$ denotes the induced $L$-section
on the quotient $E/E'$.
Hence we have
$h^0(E'_{\ast})=
 \alpha_1\cdot \dim V'+\sum \epsilon_i\cdot \dim W_i\cap V'$.
Then $\mu_{\lambda}(z,\nbigl)=0$
follows 
from the calculation of $\mu_{\lambda}(z,\nbigl)$
in the proof of Lemma \ref{lem;06.4.11.35}.
See (\ref{eq;06.5.21.150})
and (\ref{eq;06.5.21.151}).
We remark $\lambda=f_k$ and $V'=V^{(k)}$ in this case.
\hfill\qed

\begin{cor}
\label{cor;06.5.21.155}
Let $z=(q,E,F_{\ast},[\phi])$ be a point of $Q$.
Let $\lambda$ be a one parameter subgroup of $\SL(V)$.
Let $V=\bigoplus\nbigv_i$ be the weight decomposition
of $\lambda$, i.e.,
$\lambda$ preserves the decomposition,
and the weight on $\nbigv_i$ is $i$.
Let $\nbige^{(i)}$ be the subsheaf of
$E(m)$ generated by $\nbigv_j$ $(j\leq i)$ via $q$.
We have the induced $L$-section
$\phi_i$ and the parabolic structure
of $\nbige^{(i)}(-m)$.
Then all $\bigl(\nbige^{(i)}(-m)_{\ast},\phi_i\bigr)$
are $\delta$-semistable
with $P^{\delta}_{(\nbige^{(i)}(-m)_{\ast},\phi_i)}
=P^{\delta}_{(E_{\ast},\phi)}$,
if and only if
$\mu_{\lambda}(z,\nbigl)=0$ holds.
\end{cor}
\pf
We put $U_i=\bigoplus_{j\leq i}\nbigv_j$
and $U_i':=\bigoplus_{j>i}\nbigv_i$.
Let $\lambda_i$ be the one-parameter subgroup
of $\SL(V)$ given by
$t^{-\rank U_i''}\cdot id_{U_i'}\oplus
 t^{\rank U_i'}\cdot \id_{U_i''}$.
It is easy to see that $\lambda$ can be 
expressed as $\prod \lambda_i^{a_i}$
with $a_i\in\rnum_{> 0}$.
The condition $\mu_{\lambda}(z,\nbigl)=0$
implies $\mu_{\lambda_i}(z,\nbigl)=0$.
Therefore, the claim immediately follows from
Lemma \ref{lem;06.4.24.65}.
\hfill\qed

\subsection{Perturbation of $\delta$-Semistability}
\label{subsection;06.5.15.100}
\subsubsection{Preliminary}
\label{subsubsection;06.6.21.35}

We continue to use the notation in the subsubsection 
\ref{subsubsection;06.5.14.30}.
We put $V:=V_m$,
$Q:=Q^{ss}(m,\vecy,[L],\alpha_{\ast},\delta)$ and
$\nbigl:=\nbigl_{\vecy,L}(P^{\delta,\alpha_{\ast}}_{\vecy}(m),
 \epsilon_{\ast},\delta(m))$.
For simplicity,
we assume that
$\alpha_i$ and the coefficients of $\delta$
are rational.

Take a sufficiently large number $k$
such that
$\nbigl^{\otimes\,k}$ is a line bundle
on $\nbiga$.
For a rational number $\gamma$,
we put
$\nbigl_{\gamma}:=
 \nbigl^{\otimes\,k}\otimes\nbigo_{\proj_m}(\gamma)$.
Let $\nbiga^{ss}(\nbigl_{\gamma})$ 
(resp. $\nbiga^s(\nbigl_{\gamma})$)
denote the set of the semistable (resp. stable)
points of $\nbiga$ with respect to $\nbigl_{\gamma}$.

Let $\Flag(V,\Nbar)$ denote the full flag variety:
\index{$\Flag(V,\Nbar)$}
\begin{equation}
\label{eq;06.6.15.1}
 \Bigl\{
 \nbigf_{\ast}=
 \bigl(0\subset\nbigf_1\subset\nbigf_2\subset\cdots\subset
 \nbigf_N=V\bigr)\,\Big|\,
 \dim \nbigf_i/\nbigf_{i-1}=1
 \Bigr\}.
\end{equation}

Let  $G_l(V)$ denote the Grassmann variety
of $l$-dimensional subspaces of $V$.
We have the natural morphism
$\rho_l:\Flag(V,\Nbar)\lrarr G_l(V)$.
Let $\nbigo_{G_l(V)}(1)$ denote 
the canonical polarization of $G_l(V)$.
For a tuple of positive rational numbers
$n_{\ast}=(n_1,n_2,\ldots,n_N)$,
we put as follows:
\[
 \nbigo_{\Flag}(n_{\ast}):=
 \bigotimes_{i=1}^N\rho_i^{\ast}\nbigo_{G_i(V)}(n_i).
\]

We put $\Qtilde:=Q\times\Flag(V,\Nbar)$
and $\nbigatilde:=\nbiga\times\Flag(V,\Nbar)$.
We have the induced map
$\Psitilde_m:\Qtilde\lrarr \nbigatilde$.
For a tuple $n_{\ast}$
and a rational number $\gamma$,
let us consider the following $\rnum$-line bundle:
\[
 \nbigltilde(\gamma,n_{\ast}):=
 \nbigl_{\gamma}\otimes
 \nbigo_{\Flag}(n_{\ast})
=\nbigl^{\otimes\,k}
 \otimes\nbigo_{\proj_m}(\gamma)
\otimes\nbigo_{\Flag}(n_{\ast}).
\]
Let $\nbigatilde^{ss}(\gamma,n_{\ast})$
denote the set of the semistable points
with respect to $\nbigltilde(\gamma,n_{\ast})$.
We will be interested in the open subset
$\Psitilde_m^{-1}\bigl(\nbigatilde^{ss}(\gamma,n_{\ast})\bigr)$.

\subsubsection{$\delta_+$-semistability and $\delta_-$-semistability}

Let $\delta_+$ and $\delta_-$
denote elements of $\nbigp^{\br}$
with $\delta_-<\delta<\delta_+$
such that
$\delta_+$ and $\delta_-$ are sufficiently close to
$\delta$.
The following lemma is clear from 
Lemma \ref{lem;06.5.21.6}.
\begin{lem}
When $\delta_+$ (resp. $\delta_-$)
is sufficiently close to $\delta$,
a parabolic $L$-Bradlow pair $(E_{\ast},\phi)$
is $\delta_+$-semistable ($\delta_-$-semistable)
if and only if
the following condition holds:
\begin{itemize}
\item
Take any partial Jordan-H\"{o}lder filtration
of $(E_{\ast},\phi)$
 with respect to $\delta$-semistability:
\[
  \bigl(
 E^{(1)}_{\ast},\phi^{(1)}
 \bigr)
\subset
 \bigl(E^{(2)}_{\ast},\phi^{(2)}\bigr)
\subset\cdots\subset
 \bigl(E^{(k)}_{\ast},\phi^{(k)}\bigr)
=\bigl(E_{\ast},\phi\bigr).
\]
Then we have $\phi^{(i)}(0)=0$
for $i<k$
(resp. $\phi^{(1)}\neq 0$).
\end{itemize}
Moreover,
any $\delta_+$-semistable 
(resp. $\delta_-$-semistable)
$L$-Bradlow pair
is also $\delta_+$-stable
(resp. $\delta_-$-stable).
\hfill\qed
\end{lem}

We put
$Q_+:=Q^{ss}\bigl(m,\vecy,[L],\alpha_{\ast},\delta_+\bigr)$
and 
$Q_-:=Q^{ss}\bigl(m,\vecy,[L],\alpha_{\ast},\delta_-\bigr)$.
They are independent of choices of
$\delta_+$ and $\delta_-$
when $\delta_+$ and $\delta_-$ are sufficiently close to
$\delta$ due to the previous lemma.
We denote the signature of $\gamma$
by $\sign(\gamma)$.
The absolute value of $\gamma$ is denoted by $|\gamma|$.

\begin{prop}
\label{prop;06.5.15.60}
Assume that $|\gamma|$ is sufficiently small,
and that $n_i$ are sufficiently smaller than $|\gamma|$.
Then, we have
$\Psitilde_m^{-1}\bigl(
 \nbigatilde^{ss}(\gamma,n_{\ast}) \bigr)
=Q_{\sign(\gamma)}\times\Flag(V,\Nbar)$.
In particular,
we have the closed immersion
$Q_{\sign(\gamma)}\times\Flag(V,\Nbar)\lrarr
 \nbigatilde^{ss}(\gamma,n_{\ast})$.
Moreover,
the image is contained
in $\nbigatilde^{s}(\gamma,n_{\ast})$.
\end{prop}
\pf
Let us begin with the following lemma.

\begin{lem}
\label{lem;06.5.15.55}
Assume that the absolute value of $\gamma\neq 0$
is sufficiently small.
\begin{itemize}
\item
Then, we have
$\Psi_m^{-1}\bigl(\nbiga^{ss}(\nbigl_{\gamma})\bigr)
=Q^{ss}_{\sign(\gamma)}$.
\item
The induced morphism
$\Psi_m:Q_{\sign(\gamma)}^{ss}\lrarr
 \nbiga^{ss}(\nbigl_{\gamma})$
is a  closed immersion.
Moreover, the image
is contained in
$\nbiga^{s}(\nbigl_{\gamma})$.
\end{itemize}
\end{lem}
\pf
Let us show the first claim.
Let $z$ denote a point $(q,E_{\ast},\phi)\in Q$.
As a preparation,
we consider the following two cases:
\begin{description}
\item[(A)]
 There exists a partial Jordan-H\"{o}lder filtration
$E_{\ast}'\subset (E_{\ast},\phi)$
 with respect to $\delta$-semistability.
\item[(B)]
 There exists a partial Jordan-H\"{o}lder filtration
 $(E''_{\ast},\phi)\subset (E_{\ast},\phi)$
 with respect to $\delta$-semistability.
\end{description}
In the case (A),
we put $V':=H^0\bigl(X,E'(m)\bigr)$,
and take a complement $V''$ of $V'$
in $V$.
We consider the one parameter subgroup
$\lambda$
given by 
$ t^{-\rank V''}\!\cdot\!\id_{V'}
\,\oplus\,
 t^{\rank V'}\!\cdot\!\id_{V''}$.
Since we have $\mu_{\lambda}(z,\nbigl)=0$,
the equality 
$ \mu_{\lambda}\bigl(z,\nbigl_{\gamma}\bigr)
=\gamma\cdot \mu_{\lambda}\bigl(z,\nbigo_{\proj_m}(1)\bigr)
=\gamma\cdot\rank V'$
holds.

In the case (B),
we put $V'':=H^0\bigl(X,E''(m)\bigr)$,
and take a complement $V'$ of $V''$ in $V$.
Let us consider the one parameter subgroup
$\lambda$ given by
$ t^{-\rank V'}\!\cdot\!\id_{V''}
\,\oplus\, t^{\rank V''}\!\cdot\! \id_{V'}$.
As before,
we have
$ \mu(z,\nbigl_{\gamma})
=-\gamma\cdot \rank V'$.

From the above considerations,
we easily obtain
$\Psi_m^{-1}(\nbiga^{ss}(\nbigl_{\gamma})) 
\subset Q^{ss}_{\sign(\gamma)}$.

\vspace{.1in}

Let us show the reverse implication.
We use the argument in the proof of 
Lemma \ref{lem;06.4.11.35}.
Let $u_1,\ldots,u_N$ be any base of $V$.
Let $(w_1,\ldots,w_N)$ be an element of $\seisuu^N$
such that $w_i\leq w_{i+1}$
and $\sum_{i=1}^N w_i=0$.
Let $\lambda$ be the one parameter
subgroup of $\SL(V)$
given by $\lambda(t)\cdot u_i=t^{w_i}\cdot u_i$.
We have 
$\mu_{\lambda}(z,\nbigl_{\gamma})
=k\cdot\mu_{\lambda}(z,\nbigl)
+\gamma\cdot \mu_{\lambda}\bigl(z,\nbigo_{\proj_m}(1)\bigr)$.
As seen in the proof of Lemma \ref{lem;06.4.11.35},
$\mu_{\lambda}(z,\nbigl)$ is linear
with respect to $(w_1,\ldots,w_N)$.
It is also well known
that $\mu_{\lambda}\bigl(z,\nbigo_{\proj_m}(1)\bigr)$
is linear with respect to $(w_1,\ldots,w_N)$.
(See Lemma \ref{lem;06.5.15.50}, for example.)
Therefore,
we have only to show
$\mu_{f_h}\bigl(z,\nbigl_{\gamma}\bigr)>0$
for any $h=1,\ldots,N-1$.

In the case $\mu_{f_h}\bigl(z,\nbigl\bigr)>0$,
we have 
$k\cdot \mu_{f_h}\bigl(z,\nbigl\bigr)\geq 1$.
On the other hand,
the absolute value of
$\mu_{f_h}\bigl(z,\nbigo_{\proj_m}(1)\bigr)$
is dominated by $\dim V$.
Therefore,
if $\gamma$ is sufficiently small,
we have $\mu_{f_h}\bigl(z,\nbigl_{\gamma}\bigr)>0$.

Let us consider the case $\mu_{f_h}(z,\nbigl)=0$.
Let $W$ denote the subspace of $V$
generated by $u_1,\ldots,u_h$.
Let $\nbige_W$ denote the subsheaf of $E(m)$
generated by $W$ and $q$.
We put $\phi':=\phi$ if the image of $\phi$
is contained in $\nbige_W$,
and $\phi':=0$ otherwise.
Then $(\nbige'(-m),\phi')$ is $\delta$-semistable
(Lemma \ref{lem;06.4.24.65}).
Due to the considerations (A) and (B),
we obtain $\mu_{f_h}(z,\nbigl_{\gamma})>0$
in the case $(E_{\ast},\phi)$ is $\delta_{\sign(\gamma)}$-semistable.
Thus the first claim is proved.

\vspace{.1in}

In the above argument,
we showed
$\mu_{\lambda}(\Psi_m(z),\nbigl_{\gamma})>0$
for any point $z$ of $Q_{\sign(\gamma)}$.
Thus we have already obtained
$\Psi_m(Q_{\sign(\gamma)})\subset
 \nbiga^{s}(\nbigl_{\gamma})$.
Since the morphism $\Psi_m:Q\lrarr \nbiga^{ss}(\nbigl_0)$
is proper  (Proposition \ref{prop;06.4.11.50}),
the properness
of $Q_{\sign(\gamma)}\lrarr \nbiga^{ss}(\nbigl_{\gamma})$
follows from the first claim.
Thus the proof of Lemma \ref{lem;06.5.15.55}
is finished.
\hfill\qed

\vspace{.1in}
Let us return to the proof of Proposition \ref{prop;06.5.15.60}.
Let $z$ be any point of $Q\times\Flag(V,\Nbar)$.
Let $\lambda$ be the one parameter subgroup
as in the proof of Lemma \ref{lem;06.5.15.55}.
Due to Lemma \ref{lem;06.4.22.15},
$\mu_{\lambda}\bigl(z,\nbigo_{\Flag}(n_{\ast})\bigr)$
is linear with respect to $w_1,\ldots w_N$.
Hence
$\mu_{\lambda}\bigl(z,\nbigl(\gamma,n_{\ast})\bigr)$
is also linear.
Therefore,
we have only to show
$\mu_{f_h}\bigl(z,\nbigl(\gamma,n_{\ast})\bigr)>0$
for any $h=1,\ldots,N-1$.
We have the following:
\[
 \mu_{f_h}\bigl(z,\nbigl_{\gamma}\bigr)
+\sum_j n_h\cdot 
 \mu_{f_h}\bigl(z,\rho_j^{\ast}\nbigo_{G_j(V)}(1)\bigr).
\]
If we take a large integer $k'$
such that $\nbigl_{\gamma}^{\otimes k'}$
is a line bundle,
then we have $\mu_{f_h}\bigl(z,
 \nbigl_{\gamma}^{\otimes\,k'}\bigr)\geq 1$,
because 
$\mu_{f_h}\bigl(z,
 \nbigl_{\gamma}^{\otimes\,k'}\bigr)$
is a positive integer.
On the other hand,
$\mu_{f_h}\bigl(z,\rho_i^{\ast}\nbigo_{G_i(V)}(1)\bigr)$
is dominated by $2\dim V^2$.
Therefore, if $n_{i}$ are sufficiently small,
the contribution of $\nbigo_{\Flag}(n_{\ast})$
to $\mu_{f_h}\bigl(z,\nbigl(\gamma,n_{\ast})\bigr)$
is sufficiently small.
Thus we are done.
\hfill\qed

\subsubsection{$(\delta,\ell)$-semistability}

Let $\ell$ be a positive integer.
Let $\Qtilde^{ss}(\delta,\ell)$ denote
the maximal subset of $\Qtilde$,
which consists of the points 
$\bigl(q,E_{\ast},[\phi],\nbigf\bigr)$
such that $(E_{\ast},[\phi],\nbigf)$
is $(\delta,\ell)$-semistable.
(See the subsubsection \ref{subsubsection;06.4.25.10}
 for $(\delta,\ell)$-semistability.)

\begin{prop}
\label{prop;06.5.15.70}
There exist negative rational number $\gamma$
and a tuple of positive rational numbers $n_{\ast}$,
for which the following holds:
\begin{itemize}
\item
We have $\Qtilde^{ss}(\delta,\ell)=
\Psitilde_m^{-1}\bigl(\nbigatilde^{ss}(\gamma,n_{\ast})\bigr)$.
\item
The induced morphism
$\Qtilde^{ss}(\delta,\ell)\lrarr
 \nbigatilde^{ss}(\gamma,n_{\ast})$
is a closed immersion.
The image is contained in 
$\nbigatilde^{s}(\gamma,n_{\ast})$.
\end{itemize}
\end{prop}

Before going into the proof of Proposition 
\ref{prop;06.5.15.70},
we give a consequence.
\begin{cor}
The moduli stack
$\nbigmtilde^{ss}\bigl(\vecy,[L],\alpha_{\ast},(\delta,\ell)\bigr)$
of the $(\delta,\ell)$-semistable objects
is Deligne-Mumford and proper.
\end{cor}
\pf
From Lemma \ref{lem;06.5.21.60},
Proposition \ref{prop;06.5.15.70}
and Proposition \ref{prop;06.4.25.1},
the quotient stack 
$\Qtilde^{ss}(\delta,\ell)/\SL(V_m)$ is 
Deligne-Mumford and proper.
Since we have the etale and finite morphisms
$\Qtilde^{ss}(\delta,\ell)/\SL(V_m)
\lrarr \Qtilde^{ss}(\delta,\ell)/\PGL(V_m)$
and
$\nbigmtilde\bigl(\vecyhat,[L],
 \alpha_{\ast},(\delta,\ell)\bigr)
\lrarr
 \Qtilde^{ss}(\delta,\ell)/\PGL(V_m)$,
we easily obtain the desired properties
for $\nbigmtilde\bigl(\vecyhat,[L],\alpha_{\ast},(\delta,\ell)\bigr)$
by using the valuative criterion.
\hfill\qed

\vspace{.1in}

Let us start the proof of Proposition
\ref{prop;06.5.15.70}.
We put $N:=\dim V_m=H_{\vecy}(m)$.
When $\ell$ is larger than $N$,
$(\delta,\ell)$-semistability is same as
$\delta_-$-semistability
(Remark \ref{rem;06.5.21.60}).
Hence the claim follows from Proposition \ref{prop;06.5.15.60}.
Therefore, we will assume $\ell<N$ in the following argument.

We take $\gamma$ and $n_{\ast}$
satisfying the following condition:

\begin{condition}
 \label{condition;06.4.20.1}
Let $K_0(\vecy,L,\delta)$ be the number as in Proposition
{\rm \ref{prop;06.5.15.10}}.
Take a small rational number $\epsilon$
satisfying the following:
\[
 0<\epsilon<\frac{K_0(\vecy,L,\delta)}{100\cdot N^{100}}
\]
Take an irrational number $a>0$
satisfying the following:
\[
 \left(\ell-\frac{1}{\ell}\right)\cdot a<\epsilon
 < \ell\cdot a.
\]
Take mutually distinct prime numbers
$p_1,\ldots,p_{\ell}$ such that
$p_1>100\cdot N^{100}$
and $p_i>100\cdot N^{100}\cdot p_{i-1}$.
We also assume the following:
\[
 \sum\frac{i}{p_i}
<\min\left\{
 \bigl|\ell\cdot a-\epsilon\bigr|,\,\,
 \bigl|
 \epsilon-(\ell-\ell^{-1})\cdot a
 \bigr|
 \right\}
\]
Take rational numbers $q_1,\ldots,q_{\ell}$
such that the following holds:
\[
 \left|
 \frac{q_i}{p_i}-\frac{a}{i}
 \right|
<\frac{1}{p_i},
\quad
 (i=1,\ldots,\ell).
\]
We put $\gamma:=-\epsilon$
and $n_i:=q_i/p_i$ $(i=1,\ldots,\ell)$.
We remark $n_1>n_2>\cdots >n_{\ell}$
and the following inequalities:
\[
 \gamma+\sum_{i=1}^{\ell}i\cdot n_i>0,
\quad
\gamma+ \sum_{\substack{1\leq i\leq \ell\\ i\neq i_0}}
 i\cdot n_i
+n_{i_0}\cdot (i_0-1)<0
\]

We also take prime numbers
$p_i$ $(i=\ell+1,\ldots,N)$
satisfying 
$p_i>100\cdot N^{100}\cdot p_{i-1}$
 $(i=\ell+1,\ldots,N)$
and the following:
\[
 \sum_{i=\ell+1}^{N}
 \frac{i}{p_i}
\cdot 100\cdot N^{100}
<\min\left\{
\Bigl|
 \sum_{i=1}^{\ell}i\cdot n_i-\epsilon
\Bigr|,
\,\,\,
\Bigl|
\sum_{\substack{1\leq i\leq \ell\\ i\neq i_0}}
 i\cdot n_i
+n_{i_0}\cdot (i_0-1)
-\epsilon
\Bigr|
 \right\}
\]
We put $n_i=p_i^{-1}$ $(i=\ell+1,\ldots,N)$.
\hfill\qed
\end{condition}

We remark the following elementary fact.
\begin{lem}
 \label{lem;06.4.22.90}
Let $p_1,\ldots,p_N$ be mutually distinct prime numbers.
Let $q_i\neq 0$ be an integer which is coprime to $p_i$,
for each $i$.
Then the sum
$\sum_{i=1}^N q_i/p_i$ cannot be an integer.
\end{lem}
\pf
Assume that $\sum_{i=1}^Nq_i/p_i=a$
is an integer.
Then we obtain the relation:
\[
 q_1\cdot \prod_{j>1}p_j
+p_1\cdot\left(
 a\cdot\prod_{i=2}^Np_i
+ \sum_{j=2}^N
 q_j\cdot\prod_{\substack{i\geq j,\\ i\neq j}}p_i
 \right)=0.
\]
It contradicts with the assumption
that $q_1$ and $p_1$ are coprime.
\hfill\qed

\vspace{.1in}

Let us start the proof of Proposition \ref{prop;06.5.15.70}.
Let $z=(q,E_{\ast},\phi,\nbigf)$ be a point of $\Qtilde$.
We give preliminary considerations.

\vspace{.1in}
(A) If there exists a partial Jordan-H\"{o}lder
filtration $E'_{\ast}\subset (E_{\ast},\phi)$
with respect to $\delta$-semistability,
we put $V':=H^0\bigl(X,E'(m)\bigr)$,
and we take a complement $V''$ of $V'$ in $V$.
Let $\lambda$ denote the one parameter subgroup of
$\SL(V)$
given by
$ t^{-\rank V''}\!\cdot\! \id_{V'}
\,\oplus\,
 t^{\rank V'}\!\cdot\! \id_{V''}$.
Then we have the following:
\begin{multline}
\mu_{\lambda}\bigl(
 z,\nbigltilde(\gamma,n_{\ast})
 \bigr)
=k\cdot\mu_{\lambda}(z,\nbigl)
+\gamma\cdot \mu_{\lambda}\bigl(z,\nbigo_{\proj}(1)\bigr)
+\mu_{\lambda}\bigl(z,\nbigo_{\Flag}(n_{\ast})\bigr)\\
=\gamma\cdot \rank V'
+\sum n_i\cdot \left(
 -\dim\bigl(\nbigf_i\cap V'\bigr)\cdot\rank V''
+\dim\Bigl(\nbigf_i\big/\nbigf_i\cap V'
 \Bigr)\cdot\rank V'
 \right)\\
=
\left( \gamma+
\sum_{i=1}^{\ell}
 n_i\cdot \dim\bigl(\nbigf_i/\nbigf_i\cap V'\bigr)
\right)\cdot\rank V'
-\sum_{i=1}^{\ell} n_i\cdot
 \dim \bigl(\nbigf_i\cap V'\bigr)\cdot\rank V'' \\
+\sum_{i=\ell+1}^N
 n_i\cdot \left(
 -\dim (\nbigf_i\cap V')\cdot\rank V''
 +\dim\bigl(\nbigf_i/\nbigf_i\cap V'
 \bigr)\cdot\rank V'
 \right).
\end{multline}
Due to our choice of $\epsilon$ and $n_{\ast}$,
the right hand side is larger than $0$,
if and only if $\nbigf_{\ell}\cap V_{1}'=\{0\}$.
Moreover,
if it is larger than $0$,
it is strictly larger than $0$
due to Lemma \ref{lem;06.4.22.90}.

\vspace{.1in}
(B)
Let us consider the case
where there exists a partial Jordan-H\"{o}lder filtration
$(E_{\ast}'',\phi)\subset (E_{\ast},\phi)$
with respect to $\delta$-semistability.
We put $V'':=H^0\bigl(X,E''(m)\bigr)$
and take a complement $V'$ of $V''$ in $V$.
Consider the one parameter subgroup $\lambda$
given by
$\bigl\{
 t^{-\rank V''}\cdot \id_{V'}
\oplus
 t^{\rank V'}\cdot\id_{V''}
 \bigr\}$.
Then we have the following:
\begin{multline}
\mu_{\lambda}\bigl(z,\nbigltilde(\gamma,n_{\ast})\bigr)
=-\gamma\cdot\rank V'
+\sum_{i=1}^Nn_i\cdot
 \left(
 -\dim\bigl(\nbigf_i\cap V''\bigr)\cdot\rank V'
+\dim\bigl(
 \nbigf_i/\nbigf_i\cap V''
 \bigr)\cdot \rank V''
 \right) \\
=\left(
 -\gamma-\sum_{i=1}^{\ell}
 n_i\cdot \dim \bigl(\nbigf_i\cap V''\bigr)
 \right)\cdot\rank V'
+\sum_{i=1}^{\ell}
 n_i\cdot
 \dim\bigl(\nbigf_i/\nbigf_i\cap V''\bigr)
 \cdot\rank V'' \\
+\sum_{i=\ell+1}^{N}
n_i\cdot\Bigl(
 -\dim(\nbigf_i\cap V'')\cdot\rank V'
+\dim\bigl(\nbigf_i/\nbigf_i\cap V''\bigr)\cdot \rank V''
 \Bigr).
\end{multline}
Due to our choice of $\gamma$ and $n_{\ast}$,
it is strictly smaller than $0$,
if and only if $\nbigf_{\ell}\subset V''$.
Namely,
it is larger than $0$
if and only if $\nbigf_{\ell}\not\subset V''$.
Moreover,
if it is larger than $0$,
it is strictly larger than $0$
due to Lemma \ref{lem;06.4.22.90}.

From the above preliminary consideration,
we obtain
$\Psitilde^{-1}\bigl(\nbigatilde^{ss}(\delta,n_{\ast})\bigr)
\subset \Qtilde^{ss}(\delta,\ell)$.
Let us show the reverse implication.
We use the standard argument as in the proof
of Lemma \ref{lem;06.4.11.35}.
Let $z=(q,E_{\ast},\phi,\nbigf)$ be a point of
$\Qtilde^{ss}(\delta,\ell)$.
Let $u_1,\ldots,u_N$ be a base of $V$,
and let $(w_1,\ldots,w_N)$ be an element of
$\seisuu^N$
such that $w_i\leq w_{i+1}$
and $\sum w_i=0$.
Let $\lambda$ be the one-parameter subgroup
of $\SL(V)$
given by $\lambda(t)\cdot u_i=t^{w_i}\cdot u_i$.
Then $\mu_{\lambda}\bigl(z,\nbigl(\gamma,n_{\ast})\bigr)$
is linear with respect to $(w_1,\ldots,w_N)$.
Hence we have only to show 
$\mu_{f_h}\bigl(z,\nbigl(\gamma,n_{\ast})\bigr)>0$
for any $h$.

First, let us consider the case $\mu_{f_h}\bigl(z,\nbigl\bigr)>0$.
We have $\mu_{f_h}\bigl(z,\nbigl^{\otimes\,k}\bigr)\geq 1$.
Since $|\gamma|$ and $n_i$
are smaller than $100^{-1}\cdot (\dim V)^{-100}$,
we have
$\mu_{f_h}\bigl(z,\nbigo_{\proj}(\gamma)\otimes
 \nbigo_{\Flag}(n_{\ast})\bigr)<10^{-1}$.
Hence we obtain
$\mu_{f_h}\bigl(z,\nbigl(\gamma,n_{\ast})\bigr)>0$.
Next, let us consider the case
$\mu_{f_h}\bigl(z,\nbigl\bigr)=0$.
Let $\nbige'$ denote the subsheaf
of $E(m)$ generated by $u_1,\ldots,u_h$.
We put $\phi':=\phi$ if the image of
$\phi$ is contained in $\nbige'(-m)$,
and $\phi':=0$ otherwise.
Then $\bigl(\nbige'(-m)_{\ast},\phi'\bigr)\subset (E_{\ast},\phi)$
is a partial Jordan-H\"{o}lder filtration
as in (A) or (B).
Therefore,
we have $\mu_{f_h}\bigl(z,\nbigl(\gamma,n_{\ast})\bigr)>0$
in this case, too.
Hence the first claim of Proposition \ref{prop;06.5.15.70}
is obtained.

Since we have shown
that $\mu_{f_h}\bigl(z,\nbigl(\gamma,n_{\ast})\bigr)>0$
for any $h$,
we obtain that the image
$\widetilde{\Psi}(\Qtilde^{ss}(\delta,\ell))$
is contained in $\nbigatilde^{s}(\gamma,n_{\ast})$.
Let us show the properness
of $\widetilde{\Psi}:\Qtilde^{ss}(\delta,\ell)
\lrarr\nbigatilde^{ss}(\gamma,n_{\ast})$.
We use the argument in the subsubsection
\ref{subsubsection;06.5.15.80}.
Assume that we have the following diagram:
\[
 \begin{CD}
 \Spec K @>{f}>>\Qtilde^{ss}(\delta,\ell)\\
 @VVV @VVV \\
 \Spec R @>{g}>> 
 \nbigatilde^{ss}(\gamma,n_{\ast})
 \end{CD}
\]
Let $X_K$ and $X_R$ denote
$X\times\Spec K$ and $X\times \Spec R$
respectively.
We denote the projection
$\nbigatilde\lrarr\nbiga$ by $\pi$.

Let $\bigl(q_K,E_{K\,\ast},[\phitilde_K],\nbigf_{K,\ast}\bigr)$
be the objects on $X_K$ corresponding to $f$.
As in the subsubsection \ref{subsubsection;06.5.15.80},
we obtain the objects
$(q_R,E_{R\,\ast},\phi_R,W_{R\,\ast},[\phitilde_R])$
on $X_R$.
It induces the morphism
$f_1:\Spec R\lrarr \nbiga$,
which is same as $\pi\circ g$.
We also obtain the Yokogawa datum
$(q_0,E_{0\,\ast},\phi_0,W_{0,\ast},[\phitilde_0])$.

Since $|\gamma|$ and $n_{i}$ are sufficiently small,
the tuple
$\bigl(q_0,E_{0\ast},\phi_0,W_{0\,\ast},[\phitilde_0]\bigr)$
is contained in the Yokogawa family
$\nYOK(N_0,K_0,\vecy,L,\delta)$
for $N_0=N_0(\vecy,L,\delta)$ and 
$K_0=K_0(\vecy,L,\delta)$.
Hence 
$\bigl(q_0,E_{0\ast},[\widehat{\phi}]\bigr)$
gives a point of $Q=Q^{ss}(m,\vecy,[L],\alpha_{\ast},\delta)$,
due to Proposition \ref{prop;06.5.15.10}.
It implies the image of $f_1=\pi\circ g$ is contained 
in $\nbiga^{ss}(\nbigl_0)$.
Then the desired properness immediately follows
from the first claim.
\hfill\qed

\subsection{Enhanced Master Space}
\label{subsection;06.5.21.110}
\subsubsection{The construction}
\label{subsubsection;06.5.22.369}

We use the notation in the subsection
\ref{subsection;06.5.15.100}.
Take a rational number $\gamma_2<0$,
whose absolute value is sufficiently small.
We will consider the following two situations:
\begin{description}
\item[(I)]
 $\gamma_1$ is a sufficiently small positive rational number,
 and $n_{i}$  are sufficiently smaller than $\gamma_1$.
\item[(II)]
 $\gamma_1$ and $n_{\ast}$ are as in 
 Condition \ref{condition;06.4.20.1}.
\end{description}
In the both cases,
we assume $|\gamma_1|$ is sufficiently 
smaller than $|\gamma_2|$.
We also assume the following:
\begin{equation}
\label{eq;06.5.15.111}
\sum_{i}i\cdot n_i+|\gamma_2|\leq K_0(\vecy,L,\delta).
\end{equation}
Here $K_0(\vecy,L,\delta)$ denotes the constant
in Proposition \ref{prop;06.5.15.10}.

Let $k'$ be a number such that
$k'\cdot(\gamma_1-\gamma_2)=1$.
We consider the following $\rnum$-line bundles
on $\Qtilde$:
\[
 \nbigltilde_1:=\nbigltilde(\gamma_1,n_{\ast})^{\otimes\,k'},
\quad
 \nbigltilde_2:=\nbigltilde(\gamma_2,n_{\ast})^{\otimes\,k'}.
\]
Then we have
$\nbigltilde_2=\nbigltilde_1\otimes\nbigo_{\proj_m}(-1)$.

Let $\pi_1:\Qtilde\lrarr \proj_m$
denote the projection.
We put
$\nbigbtilde:=
\proj\bigl(
 \pi_1^{\ast}\nbigo_{\proj_m}(0)\oplus\pi_1^{\ast}\nbigo_{\proj_m}(1)
 \bigr)$ over $\Qtilde$.
We put $\nbigo_{\nbigbtilde}(1):=
 \nbigo_{\proj_m}(1)\otimes\nbigltilde_1$,
where $\nbigo_{\proj_m}(1)$ denotes
the tautological bundle of
$\proj\bigl(\pi_1^{\ast}\nbigo_{\proj_m}(0)
  \oplus\pi_1^{\ast}\nbigo_{\proj_m}(1)\bigr)$
over $\Qtilde$.
Let $\nbigbtilde^{ss}$ (resp. $\nbigbtilde^s$)
denote the set of the semistable points
(resp. the stable points)
with respect to $\nbigo_{\nbigbtilde}(1)$.

We put $\nbigbtilde_1:=\proj\bigl(\pi_1^{\ast}\nbigo_{\proj_m}(0)\bigr)$
and $\nbigbtilde_2:=\proj\bigl(\pi_1^{\ast}\nbigo_{\proj_m}(1)\bigr)$.
We naturally regard $\nbigbtilde_i$
as the closed subscheme of
$\nbigbtilde$.
The following lemma is clear from the construction.
\begin{lem}
The restriction of $\nbigo_{\nbigbtilde}(1)$ to
$\nbigbtilde_i$
is same as $\nbigltilde_i$.
The $\rnum$-line bundle $\nbigo_{\nbigbtilde}(1)$
gives a $\GL(V_m)$-equivariant polarization.
\hfill\qed
\end{lem}

We put 
$\TH:=\nbigbtilde\times_{\nbigatilde}\Qtilde$,
$\TH_i:=\nbigbtilde_i\times_{\nbigatilde}\Qtilde$
and $\TH^{\ast}:=\TH-\bigl(\TH_1\cup\TH_2\bigr)$.
We remark that
$\TH^{\ast}$ is isomorphic to
$Q^{ss}(m,\vecy,L,\alpha_{\ast},\delta)
\times\Flag(V,\Nbar)$.
We also put
$\TH^{ss}=\nbigbtilde^{ss}\times_{\nbigatilde}\Qtilde$.
The following lemma is obvious
from Proposition \ref{prop;06.5.15.60},
Proposition \ref{prop;06.5.15.70}
and our choice of the constants.

\begin{lem}
\mbox{{}}
\begin{itemize}
\item
In the case (I),
we have 
$\TH^{ss}\times_{\TH}\TH_1
=Q^{ss}_+\times\Flag(V,\Nbar)$.
\item
In the case (II), we have 
$\TH^{ss}\times_{\TH}\TH_1
=\Qtilde(\delta,\ell)$.
\item
In the both cases,
we have
$\TH^{ss}\times_{\TH}\TH_2
=Q^{ss}_{-}\times\Flag(V,\Nbar)$.
\hfill\qed
\end{itemize}
\end{lem}

The quotient stack $\TH^{ss}/\SL(V)$
is called the enhanced master space.
In the rest of this subsection,
we will show the following proposition.
\begin{prop}
\label{prop;06.5.16.50}
The stack $\TH^{ss}/\SL(V)$ is Deligne-Mumford
and proper.
\end{prop}

Due to Corollary \ref{cor;06.4.25.100},
the proposition is obtained from the following
three lemmas.

\begin{lem}
 \label{lem;06.4.22.101}
Let us consider the $\SL(V)$-action on $\TH^{ss}$.
The stabilizer of any point $z\in \TH^{ss}$ is finite
and reduced.
As a result,
$\TH^{ss}/\SL(V)$ is Deligne-Mumford.
\end{lem}

\begin{lem}
 \label{lem;06.4.25.102}
If $m$ is sufficiently large,
then the morphism
$\TH^{ss}\lrarr \nbigbtilde^{ss}$
is proper.
\end{lem}

\begin{lem}
\label{lem;06.4.22.103}
The image of $\TH^{ss}\lrarr\nbigbtilde^{ss}$
is contained in $\nbigbtilde^{s}$.
\end{lem}

\subsubsection{Proof of Lemma \ref{lem;06.4.22.101}}

The claim is obvious for any point
$z\in \TH^{ss}\cap\bigl(\TH_1\cup \TH_2\bigr)$.
So we discuss the stabilizer of a point 
$z=\bigl(q,E_{\ast},[\phi],\nbigf_{\ast},u\bigr)
\in \TH^{ss}\cap \TH^{\ast}$.
Let $g\in\SL(V)$ be any element such that $g\cdot z=z$.
Let $V=\bigoplus V^{(i)}$
denote the generalized eigen decomposition of $g$.
Correspondingly,
we have the decomposition:
\[
 \bigl(q,E_{\ast},[\phi],\nbigf\bigr)
=\bigl(q^{(1)},E^{(1)}_{\ast},[\phi],\nbigf^{(1)}\bigr)
\oplus
 \bigoplus_{i=2}^{l}
 \bigl(q^{(i)},E^{(i)}_{\ast},\nbigf^{(i)}\bigr).
\]

\begin{lem}
\label{lem;06.4.22.50}
$l\leq 2$.
\end{lem}
\pf
Assume $l>2$, and we will derive a contradiction.
Let us consider the one parameter subgroup
$\lambda$ given by
$ t^{-\rank V^{(3)}}\id_{V^{(2)}}\,
\oplus\,
 t^{\rank V^{(2)}}\id_{V^{(3)}}$.
It is easy to see that
$\lambda$ fixes the point $z$.
Since we have
$\mu_{\lambda}(z,\nbigl)
=\mu_{\lambda}\bigl(z,\nbigo_{\proj_m}(1)\bigr)=0$,
we have the following equality:
\begin{equation}
 \label{eq;06.4.22.30}
 \mu_{\lambda}\bigl(z,\nbigo_{\nbigbtilde}(1)\bigr)
=\sum_{i=1}^N
 n_i\cdot
\Bigl(
 -\rank V^{(3)}
 \cdot\rank\nbigf^{(2)}_i
+\rank V^{(2)}\cdot\rank\nbigf^{(3)}_i
\Bigr).
\end{equation}
Due to our choice of $n_{\ast}$,
the right hand side of (\ref{eq;06.4.22.30})
can be $0$,
if and only if
the following equality holds for any $i$:
\[
 -\rank V^{(3)}\cdot
 \rank\nbigf^{(2)}_i
+\rank V^{(2)}\cdot \rank\nbigf^{(3)}_i=0.
\]
However,
there exists a number $i_0$
such that
$\rank\nbigf^{(2)}_{i_0+1}
=\rank\nbigf^{(2)}_{i_0}+1$
and 
$\rank\nbigf^{(3)}_{i_0+1}
=\rank\nbigf^{(3)}_{i_0}$.
Thus
$\mu_{\lambda}\bigl(z,\nbigo_{\nbigbtilde}(1)\bigr)
 \neq 0$.
Let $\lambda^{-1}$ denote
the one-parameter subgroup given by
$\lambda^{-1}(t)=\lambda(t)^{-1}$.
Then one of 
$\mu_{\lambda}\bigl(z,\nbigo_{\nbigbtilde}(1)\bigr)$
or $\mu_{\lambda^{-1}}\bigl(z,\nbigo_{\nbigbtilde}(1)\bigr)$
is negative.
Hence $z$ cannot be semistable.
\hfill\qed

\begin{lem}
\label{lem;06.4.22.51}
Let $N$ be a nilpotent endomorphism
of $(E_{\ast},[\phi],\nbigf)$.
Then $N=0$.
\end{lem}
\pf
Assume $N\neq 0$,
and we will derive a contradiction.
There exists the integer such that
$N^j\neq 0$ and $N^{j+1}=0$.
It is easy to obtain $N(\phi)=0$.
We obtain the subsheaves
$\Image N^j$ and $\Ker(N^j)$ of $E$.
Then we obtain the naturally induced
parabolic $L$-Bradlow pairs
$(\Image N^j_{\ast},\phi')
\subset
 (\Ker N^j_{\ast},\phi'')
\subset
 (E_{\ast},\phi)$,
which gives the partial Jordan-H\"{o}lder
filtration with respect to $\delta$-semistability,
due to Lemma \ref{lem;06.4.24.6}.
We take subspaces $\nbigk_i$ $(i=1,2,3)$
of $V$
satisfying the following condition:
\[
 \nbigk_1=H^0\bigl(X,\Image N^j\bigr),
\quad
 \nbigk_1\oplus\nbigk_2
=H^0\bigl(X,\Ker N^j\bigr),
\quad
 \nbigk_1\oplus\nbigk_2\oplus\nbigk_3=V.
\]
We remark that
$N^j$ induces the isomorphism
$\nbigk_3\lrarr \nbigk_1$.
From the inclusion
$\nbigk_1\subset V$,
we have the induced filtration $\nbigf_{\nbigk_1}$.
From the isomorphism
$\nbigk_3\simeq V\big/(\nbigk_1\oplus\nbigk_2)$,
the filtration $\nbigf_{\nbigk_3}$ is induced.
We remark
$N\bigl(\nbigf_{\nbigk_3\,h}\bigr)
\subset 
 \nbigf_{\nbigk_1\,h}$.
Let us consider the one-parameter subspace
$\lambda$ given by
$t^{-1}\id_{\nbigk_1}\,
\oplus\,
 t^{}\id_{\nbigk_3}$,
for which we have
$\mu_{\lambda}(z,\nbigl)=
 \mu_{\lambda}\bigl(z,\nbigo_{\proj_m}(1)\bigr)=0$.
From
$N\bigl(\nbigf_{\nbigk_3\,h}\bigr)
\subset 
 \nbigf_{\nbigk_1\,h}$
and $n_i>n_{i+1}$,
we obtain
$\mu_{\lambda}\bigl(z,\nbigo_{\Flag}(n_{\ast})\bigr)<0$.
Thus the point $z$ cannot be semistable.
\hfill\qed

\vspace{.1in}

From the lemmas \ref{lem;06.4.22.50}
and \ref{lem;06.4.22.51},
we obtain the following lemma.
\begin{lem}
\label{lem;06.5.21.100}
Let $(q,E_{\ast},[\phi],\nbigf)$ be a point of $\Qtilde$
such that $z=(q,E_{\ast},[\phi],\nbigf,u)\in\TH^{ss}\cap \TH^{\ast}$.
Then, either one of the following holds:
\begin{enumerate}
\item
 The automorphism group of
 $(E_{\ast},[\phi],\nbigf)$ is
 $G_m$.
\item
 There exists the unique decomposition
 $(q,E_{\ast},[\phi],\nbigf)
 =(q^{(1)},E^{(1)}_{\ast},[\phi^{(1)}],\nbigf^{(1)})
 \oplus
 (q^{(2)},E^{(2)}_{\ast},\nbigf^{(2)})$,
 and the automorphism group of
 $(q,E_{\ast},[\phi],\nbigf)$
 is $G_m^2$.
\hfill\qed
\end{enumerate}
\end{lem}

In the first case,
the stabilizer of $(q,E_{\ast},[\phi],\nbigf,u)$
with respect to the $SL(V)$-action
is trivial.
In the second case,
we put $V^{(i)}:=H^0\bigl(X,E^{(i)}(m)\bigr)\subset V$.
Then the intersection $G_m^2\cap \SL(V)$ 
consists of the elements 
$\rho(t)=t^{a}\cdot \id_{V^{(1)}}
\oplus
 t^b\cdot\id_{V^{(2)}}$
satisfying
$a\cdot \rank V^{(1)}+b\cdot \rank V^{(2)}=0$.
By considering the action along the direction of
the fiber $\TH/\Qtilde$,
which is given by $\rho(t)u=t^{a}\cdot u$,
we obtain that the stabilizer
is finite.
\hfill\qed

\subsubsection{Proof of Lemma \ref{lem;06.4.25.102}}

Let $(q,E_{\ast},\phi,W_{\ast},[\phitilde])$
be a Yokogawa datum.
Recall that 
we obtain the element
$\Psi(q,E_{\ast},\phi,W_{\ast},[\phitilde])
 \in \nbiga$.

\begin{lem}
\label{lem;06.5.15.121}
Assume that $m$ is larger than 
the constant $N(\vecy,L,\delta)$
in Proposition {\rm\ref{prop;06.5.15.10}}.
Let $z$ be a point of $\nbigbtilde^{ss}$
such that
$\pi(z)=
 \Psi\bigl(q,E_{\ast},[\phi],W_{\ast},[\phitilde]\bigr)$,
where $\pi$ denotes the naturally defined projection
$\nbigbtilde\lrarr\nbiga$.
Then, 
$(E_{\ast},[\phi])$ is $\delta$-semistable,
$q$ is onto,
and the condition $O_m$ holds for $(E_{\ast},[\phi])$.
\end{lem}
\pf
Let $W$ be any subspace of $V$.
Take a complement $W'$ of $W$ in $V$.
Consider the one-parameter subgroup
of $\SL(V)$ given by
$t^{-\rank W'}\id_W\,
\oplus\, t^{\rank W}\id_{W'}$.
We have the following:
\begin{equation}
 \mu_{\lambda}\bigl(z,\nbigo_{\nbigbtilde}(1)\bigr)
=k\cdot k'\mu_{\lambda}\bigl(z,\nbigl\bigr)
+k'\sum n_j\cdot \mu_{\lambda}\bigl(z,
 \nbigo_{G_l(V)}(1)\bigr)
+k'\max\Bigl\{
 \gamma_1\cdot \mu_{\lambda}\bigl(z,\nbigo_{\proj_m}(1)\bigr),
\,\,
 \gamma_2\cdot\mu_{\lambda}\bigl(z,\nbigo_{\proj_m}(1)\bigr)
 \Bigr\}.
\end{equation}
The first term in the right hand side is as follows:
\[
 k\cdot k'\cdot H(m)
 \cdot\Bigl(
 P^{\alpha_{\ast},\delta}(m)
\cdot\rank\nbige_W
-\sum \epsilon_i\dim (W_i\cap W)
-\alpha_1\dim (W)
-\epsilon(W,[\phi])\cdot\delta(m)
 \Bigr).
\]
The absolute value of the second term
can be dominated by 
$ k'\sum_{i=1}^n n_i\cdot i\cdot\dim V$.
The absolute value of the third term
can be dominated by
$k'\cdot |\gamma_2|\cdot\dim V$.
Recall $\dim V=H(m)$.
Since we have assumed
that $|\gamma_i|$ and $n_j$ are sufficiently small
as in (\ref{eq;06.5.15.111}),
we obtain the following inequality:
\begin{equation}
 \label{eq;06.4.22.60}
 P_{\vecy}^{\alpha_{\ast},\delta}(m)
\cdot\rank\nbige_{W}
-\sum\epsilon_i\cdot\dim(W_i\cap W)
-\alpha_1\dim(W)-\epsilon(W,[\phi])\cdot\delta(m)
+K\geq 0.
\end{equation}
Here $K=K_0(\vecy,L,\delta)$ denotes the constant
in Proposition \ref{prop;06.5.15.10}.
Namely,
$(q,E_{\ast},\phi,W_{\ast},[\phitilde])$
is contained in
$\YOKtilde(m,K,\vecy,L,\delta)$.
Therefore,
the claim of the lemma follows from Proposition \ref{prop;06.5.15.10}.
\hfill\qed

\vspace{.1in}

Now, we use the same argument as the last part
of the proof of Proposition \ref{prop;06.5.15.70}.
Assume that we have a diagram:
\[
 \begin{CD}
 \Spec K @>{f}>> \TH^{ss}\\
 @VVV @VVV\\
 \Spec R@>{g}>>
 \nbigbtilde^{ss}
 @>>> \nbiga
 \end{CD}
\]
Then we can show the composite $\pi\circ g$
is contained in $\nbiga^{ss}(\nbigl_0)$,
by using Lemma \ref{lem;06.5.15.121}.
Then the desired properness
follows from the definition of $\TH^{ss}$.
\hfill\qed

\subsubsection{Proof of Lemma \ref{lem;06.4.22.103},
 Step 1}

Let us show that the image of
$\TH^{ss}$ is contained in $\nbigbtilde^{s}$.
Let $z=(q,E_{\ast},[\phi],\nbigf,u)$ 
be a point of $\TH^{ss}$.
We have only to consider the case $u\neq 0$.

Let $u_1,\ldots,u_N$ be a base of $V$,
and let $(w_1,\ldots,w_N)$ be an element of
$\seisuu^N$
such that $w_i\leq w_{i+1}$ and $\sum w_i=0$.
Let $\lambda$ be the one-parameter subgroup
of $\SL(V)$ given by
$\lambda(t)\cdot u_i=t^{w_i}\cdot u_i$.
We will not distinguish the elements
$w=(w_1,\ldots,w_N)$ and $\lambda$.
We have the following:
\[
 \mu_{\lambda}\bigl(z,\nbigo_{\nbigbtilde}(1)\bigr)
=k\cdot k'\cdot \mu_{\lambda}\bigl(z,\nbigl\bigr)
+k'\cdot \mu_{\lambda}\bigl(z,\nbigo_{\Flag}(n_{\ast})\bigr)
+k'\max_{i=1,2}\Bigl\{
 \gamma_i\cdot\mu_{\lambda}
 \bigl(z,\nbigo_{\proj_m}(1)\bigr)
 \Bigr\}.
\]

Recall that we have the expression
$\lambda=\sum a_j\cdot f_j$
for $a_j\geq 0$,
where $f_j=\bigl(\overbrace{j-N,\ldots,j-N}^j,j\ldots,j\bigr)$.
\begin{lem}
If $\mu_{\lambda}(z,\nbigl)=0$,
then $\mu_{f_h}(z,\nbigl)=0$
for any $h$ such that $a_h\neq 0$.
\end{lem}
\pf
Since $(E_{\ast},[\phi])$ is $\delta$-semistable,
the claim immediately follows.
\hfill\qed

\vspace{.1in}

We put
$S_1:=\bigl\{
 j\,\big|\, \mu_{f_j}(z,\nbigl)=0
 \bigr\}$
and
$S_2:=\bigl\{
 j\,\big|\,\mu_{f_j}(z,\nbigl)>0
 \bigr\}$.

\begin{lem}
\label{lem;06.5.16.2}
For any element
$0\neq \rho=\sum_{j\in S_2}a_j\cdot f_j$,
we have the following:
\[
 k\cdot\mu_{\rho}\bigl(z,\nbigl\bigr)
+\mu_{\rho}\bigl(z,\nbigo_{\Flag}(n_{\ast})\bigr)
+\min_{i=1,2}\Bigl\{
 \gamma_i\cdot
 \mu_{\rho}\bigl(z,\nbigo_{\proj_m}(1)\bigr)
 \Bigr\}>0.
\]
\end{lem}
\pf
We put
$F_i(\rho):=k\cdot\mu_{\rho}(z,\nbigl)
+\mu_{\rho}(z,\nbigo_{\Flag}(n_{\ast}))
+\gamma_i\cdot\mu_{\rho}(z,\nbigo_{\proj_m}(1))$
for $i=1,2$.
We have only to show $F_i(\rho)>0$.
Since $F_i$ are linear with respect to $\rho$,
we have only to show
$F_i(f_j)>0$ for any $j\in S_2$.
We remark
$k\cdot \mu_{f_j}(z,\nbigl)\geq 1$.
The number $\mu_{f_j}(z,\nbigo_{\Flag}(n_{\ast}))$
is dominated by 
$\dim(V)$ and $n_{i}$ $(i=1,\ldots,N)$.
The number $\gamma_i\cdot\mu_{f_j}(z,\nbigo_{\proj_m}(1))$
is dominated by
$\dim(V)$ and $\gamma_i$.
Since $n_j$ and $\gamma_i$ are sufficiently small,
the claim is clear.
\hfill\qed

\begin{lem}
 \label{lem;06.4.22.150}
To show $\mu_{\lambda}\bigl(z,\nbigo_{\nbigbtilde}(1)\bigr)>0$
for any $\lambda$,
we have only to show
the following inequality
for $0\neq \rho=\sum_{j\in S_1}a_j\cdot f_j$:
\begin{equation}
 \label{eq;06.4.22.70}
 k\cdot \mu_{\rho}(z,\nbigl)
+\mu_{\rho}\bigl(z,\nbigo_{\Flag}(n_{\ast})\bigr)
+\max_{i=1,2}\Bigl\{
 \gamma_i\cdot\mu_{\rho}(z,\nbigo_{\proj_m}(1))
 \Bigr\}>0.
\end{equation}
\end{lem}
\pf
We have the decomposition
$\lambda=\lambda^{(1)}+\lambda^{(2)}$,
where $\lambda^{(i)}=\sum_{j\in S_i}a_j\cdot f_j$.
Assume (\ref{eq;06.4.22.70}) holds for
$\lambda^{(1)}$.
Then we obtain the following:
\begin{multline}
\frac{1}{k'}
\mu_{\rho}\bigl(z,\nbigo_{\nbigbtilde}(1)\bigr)
=\sum_{i=1,2}\Bigl(
 k\cdot\mu_{\rho^{(i)}}(z,\nbigl)
+\mu_{\rho^{(i)}}(z,\nbigo_{\Flag}(n_{\ast}))
 \Bigr)
+\max_{j=1,2}\bigl\{
 \gamma_j\cdot\mu_{\rho^{(1)}}(z,\nbigo_{\proj_m}(1))
+\gamma_j\cdot\mu_{\rho^{(2)}}(z,\nbigo_{\proj_m}(1))
 \bigr\} \\
\geq
 k\cdot \mu_{\rho^{(1)}} (z,\nbigl)
+\mu_{\rho^{(1)}}(z,\nbigo_{\Flag}(n_{\ast}))
+\min_{j=1,2}\bigl\{
 \gamma_j\cdot \mu_{\rho^{(1)}}(z,\nbigo_{\proj_m}(1))
 \bigr\} \\
+k\cdot\mu_{\rho^{(2)}}\bigl(z,\nbigl\bigr)
+\mu_{\rho^{(2)}}\bigl(z,\nbigo_{\Flag}(n_{\ast})\bigr)
+\max_{j=1,2}\bigl\{
 \gamma_j\cdot\mu_{\rho^{(2)}}\bigl(
 z,\nbigo_{\proj_m}(1) \bigr)
\bigr\}>0.
\end{multline}
Thus we are done.
\hfill\qed

\subsubsection{Proof of Lemma \ref{lem;06.4.22.103},
 Step 2}

To show (\ref{eq;06.4.22.70}),
we give some preliminary consideration.

\vspace{.1in}

(A) If  there exists a partial Jordan-H\"{o}lder filtration
 $E'_{\ast}\subset (E_{\ast},\phi)$
with respect to $\delta$-semistability,
 take a decomposition $V=V'\oplus V''$
 such that
 $V'=H^0\big(X,E'(m)\bigr)$
 and $V=V'\oplus V''$.
 Consider the one parameter subgroup
 $\lambda$ given by
 $t^{-\rank V''}\id_{V'}\oplus\, t^{\rank V'}\cdot\id_{V''}$.
 Then we have
$\mu_{\lambda}(z,\nbigl)=0$,
$\mu_{\lambda}\bigl(z,\nbigo_{\proj_m}(1)\bigr)
=\rank V'>0$,
and the following equality:
\[
 \mu_T(z,\nbigo_{\Flag}(n_{\ast}))
=\sum n_j\cdot \Bigl(
 -\rank V''\cdot\dim\nbigf_j\cap V'
+\rank V'\cdot\dim\frac{\nbigf_j}{\nbigf_j\cap V'}
 \Bigr).
\]
From the semistability of $z$,
we have the following:
\begin{equation}
 \label{eq;06.4.22.80}
 \gamma_1\cdot\rank V'+
\sum n_j\cdot\Bigl(
-\rank V''\cdot\dim(\nbigf_j\cap V')
+\rank V'\cdot\dim\frac{\nbigf_j}{\nbigf_j\cap V'}
 \Bigr)\geq 0.
\end{equation}
We remark that the strict inequality
holds in (\ref{eq;06.4.22.80}).
In the case (I), it is obvious.
In the case (II),
it follows from our choice of $n_{\ast}$
and Lemma \ref{lem;06.4.22.90}.
Hence 
$\mu_{\lambda}\bigl(z,\nbigo_{\nbigbtilde}(1)\bigr)>0$
in this case.

\vspace{.1in}

(B)
If there exists a partial Jordan-H\"{o}lder filtration
$(E''_{\ast},\phi'')\subset (E_{\ast},\phi)$
with respect to $\delta$-semistability
such that $\phi''\neq 0$,
let us take a decomposition
$V=V'\oplus V''$ such that
$V''=H^0\bigl(X,E''(m)\bigr)$.
Consider the one-parameter subgroup
$\lambda$ given by
$t^{-\rank V'}\id_{V''}\,
\oplus\,
 t^{\rank V''}\id_{V'}$.
We have
$\mu_{\lambda}(z,\nbigl)=0$,
$\mu_{\lambda}(z,\nbigo_{\proj_m}(1))=-\rank V'$
and the following:
\[
 \mu_{\lambda}(z,\nbigo_{\Flag}(n_{\ast}))
=\sum n_j\cdot
\left(
 -\rank V'\cdot \rank\bigl(\nbigf_j\cap V''\bigr)
+\rank V''\cdot\rank\frac{\nbigf_j}{\nbigf_j\cap V''}
\right).
\]
From the semistability of $z$,
we obtain the following:
\begin{equation}
 \label{eq;06.4.22.95}
 -\gamma_2\cdot\dim V'
+\sum n_j\cdot
\left(
 -\rank V'\cdot \rank\bigl(\nbigf_j\cap V''\bigr)
+\rank V''\cdot\rank\frac{\nbigf_j}{\nbigf_j\cap V''}
\right)\geq 0
\end{equation}
In the both cases of (I) and (II),
the strict inequality holds
in (\ref{eq;06.4.22.95}).
Hence we have
$\mu_{\lambda}\bigl(
 z,\nbigo_{\nbigbtilde}(1)\bigr)>0$
in this case.

\vspace{.1in}
(C)
Let us consider the case
that there exists a partial Jordan-H\"{o}lder filtration
$E'_{\ast}\subset (E''_{\ast},\phi)\subset (E_{\ast},\phi)$
with respect to $\delta$-semistability.
We take a decomposition
$V=V'\oplus V''\oplus V'''$
such that
$V'=H^0(X,E'(m))$ and
$V'\oplus V''=H^0(X,E''(m))$.
Consider the one-parameter subgroup
$\lambda$
given by
$t^{-\rank V'''}\cdot\id_{V'}\,
\oplus\,
 t^{\rank V'}\cdot \id_{V'''}$.
Then we have
$\mu_{\lambda}(z,\nbigl)
=\mu_{\lambda}(z,\nbigo_{\proj_m}(1))=0$.
Hence,
we obtain the following inequality
from the semistability of the point $z$:
\begin{equation}
 \label{eq;06.4.22.100}
 \sum n_j\cdot\left(
 -\rank V''\cdot\rank(\nbigf_j\cap V')
+\rank V'\cdot\rank\left(\frac{\nbigf_j}{\nbigf_j\cap V''}\right)
\right)
\geq 0.
\end{equation}
Due to our choice of $n_{\ast}$
and Lemma \ref{lem;06.4.22.90},
the strict inequality holds
in (\ref{eq;06.4.22.100}).
Therefore,
we have $\mu_{\lambda}(z,\nbigo_{\nbigbtilde}(1))>0$
in this case.

\subsubsection{Proof of Lemma \ref{lem;06.4.22.103},
 Step 3}

For $\rho=\sum_{f_j\in S_1}a_j\cdot f_j$,
let $V=\bigoplus V_i$ be the weight decomposition.
We have the number $i_0$
such that $\phi\in\nbigg_{i_0}-\nbigg_{i_0-1}$.
We put $r_i=\dim V_i$.

We use the notation in the subsubsection
\ref{subsubsection;06.5.15.150}.
We regard $\rho$ as an element of
$\nbigu=\bigoplus_{i=1}^n\rnum\cdot e_i$.
Then we have the expression
$\rho=\sum a'_j\cdot v_j$
such that
$a_1'\leq a_2'\leq \cdots \leq a_s'$
and $\sum r_i\cdot a_i'=0$.

Assume $a'_{i_0}>0$.
Then we have the expression,
due to Lemma \ref{lem;06.4.22.105}:
\[
 \rho=\sum_{(i_1,i_2)\in\nbigs(i_0)}
 b(i_1,i_2)\cdot x(i_1,i_2)+
 \sum_{j=1}^{i_0-1} c_j\cdot y(j).
\]
Here the coefficients $b(i_1,i_2)$ and $c_j$
are non-negative,
and one of $b(i_1,i_2)$ or $c_j$ is positive.
We remark that we have the following linearity:
\[
 \max_{i=1,2}\bigl\{
 \gamma_i\cdot \mu_{\rho}(z,\nbigo_{\proj_m}(1))
 \bigr\}
=\gamma_1\cdot
 \sum_{j=1}^{i_0-1}c_j\cdot
 \mu_{y(j)}(z,\nbigo_{\proj_m}(1)).
\]
Hence the following holds:
\[
 \mu_{\rho}\bigl(z,\nbigo_{\nbigbtilde}(1)\bigr)
=\sum b(i_1,i_2)\cdot
 \mu_{x(i_1,i_2)}\bigl(z,\nbigo_{\nbigbtilde}(1)\bigr)
+\sum_{j=1}^{i_0-1}c_j\cdot
 \mu_{y(j)}\bigl(z,\nbigo_{\nbigbtilde}(1)\bigr).
\]
We have the positivity
$\mu_{x(i_1,i_2)}\bigl(z,\nbigo_{\nbigbtilde}(1)\bigr)>0$
and $\mu_{y(j)}\bigl(z,\nbigo_{\nbigbtilde}(1)\bigr)>0$
from (C) and (A), respectively.
Therefore,
we obtain the positivity
$\mu_{\rho}\bigl(z,\nbigo_{\nbigbtilde}(1)\bigr)>0$.

By similar arguments,
we can show the desired positivity
in the cases $a'_{i_0}=0$ and $a'_{i_0}<0$.
Thus the image of $\TH^{ss}$
is contained in $\nbigbtilde^{s}$.
Therefore,
we obtain Lemma \ref{lem;06.4.22.103}
and hence Proposition \ref{prop;06.5.16.50}.
\hfill\qed

\subsection{Fixed Point Set of the Torus Action
 on Enhanced Master Space}
\label{subsection;06.6.6.2}
\subsubsection{Preliminary}

We continue to use the notation
in the subsection \ref{subsection;06.5.21.110}.
We have the $G_m$-action $\rho$ on
$\nbigbtilde=
\proj\bigl(\nbigo_{\proj_m}(0)\oplus
 \nbigo_{\proj_m}(1) \bigr)$
given by
$\rho(t)\cdot [u_1:u_2]=[t\cdot u_1:u_2]$.
It induces the $G_m$-action on $\TH$,
which is also denoted by $\rho$.
Since it commutes with the $\SL(V)$-action,
we obtain the induced action $\rhobar$
on $\TH^{ss}/\SL(V)$.

We would like to discuss the fixed point set
of the enhanced master space.
The stack theoretic fixed point set 
(see \cite{gp})
is given in our case, as follows:
We have the $\SL(V)\times G_m$-equivariant
closed immersion
$\TH^{ss}\lrarr \nbigbtilde^s$.
Therefore, we have an open neighbourhood
$\nbigu$ of $\TH^{ss}/\SL(V)$
in $\nbigbtilde^s/\SL(V)$,
which is $G_m$-invariant,
Deligne-Mumford and smooth.
The embedding
$\TH^{ss}/\SL(V)\lrarr \nbigu$
is $G_m$-equivariant.
The fixed point set 
$\nbigu^{G_m}$
of $\nbigu$
is defined to be the $0$-set of 
the vector field induced by the $G_m$-action.
Then, the stack theoretic fixed point set $M^{G_m}$
of $M$ is defined to be the intersection
$M\cap\nbigu^{G_m}$.

However, we restrict ourselves
to the set theoretic fixed point set in this subsection.
In other words, we will consider
only the closed points of the fixed point set,
although it is not difficult to prove the result
for the stack theoretic fixed point set.
We will later discuss the stack theoretic fixed point
set of the enhanced master space
in the oriented case.

\vspace{.1in}

Let $\pi$ be the projection
$\TH^{ss}\lrarr \TH^{ss}/\SL(V)$.
We will use the notation
$\bigl(q,E_{\ast},[\phi],\nbigf,u\bigr)$
to denote a point of $\TH$,
where $(q,E_{\ast},[\phi],\nbigf)$ denotes
a point of $\Qtilde$,
and $u$ denotes a point of the fiber
of $\TH\lrarr \Qtilde$ over $(q,E_{\ast},[\phi],\nbigf)$.

\begin{lem}
\label{lem;06.5.16.1}
Let $z=\bigl(q,E_{\ast},[\phi],\nbigf,u\bigr)$
be a point of $\TH^{ss}$.
The point $\pi(z)$ is contained in the fixed point set
if and only if one of the following holds:
\begin{enumerate}
\item
$z\in \TH_1\cup \TH_2$.
\item
We have the unique decomposition
$(q,E_{\ast},[\phi],\nbigf)
=(q^{(1)},E^{(1)}_{\ast},[\phi^{(1)}],\nbigf^{(1)})
\oplus
 (q^{(2)},E^{(2)}_{\ast},\nbigf^{(2)})$.
\end{enumerate}
\end{lem}
\pf
We have only to consider the case
$z\in \TH^{ss}\cap \TH^{\ast}$.
Assume that the condition $2$ holds.
We put $V^{(i)}:=H^0\bigl(X,E^{(i)}(m)\bigr)$,
and we consider the one parameter subgroup
$\lambda$ of $\SL(V)$ given by
$ t^{-\rank V^{(2)}} \id_{V^{(1)}}
\oplus 
 t^{\rank V^{(1)}}\cdot \id_{V^{(2)}}$.
It fixes
$\bigl(q,E_{\ast},[\phi],\nbigf\bigr)$,
and 
it acts non-trivially along the direction
of the fiber $\TH\lrarr \Qtilde$,
as $\lambda(t)u=t^{-\rank V^{(2)}}u$.
Therefore,
the action $\rhobar$ fixes $[z]$.

On the other hand,
if $\pi(z)$ is a fixed point with respect to $\rhobar$,
then we obtain the one-parameter subgroup
of $\SL(V)$ which fixes $(q,E_{\ast},[\phi],\nbigf)$
due to Lemma \ref{lem;06.5.21.100}.
Hence, it has the decomposition.
The uniqueness follows from 
Lemma \ref{lem;06.4.22.101}.
\hfill\qed

\vspace{.1in}

Let $z=(q,E_{\ast},[\phi],\nbigf)$ be
a point of $\TH^{ss}\cap\TH^{\ast}$
such that $\pi(z)\in M^{G_m}$.
We have the decomposition
as in Lemma \ref{lem;06.5.16.1}.
Then, we obtain the types
$\vecy_1=\type(E^{(1)}_{\ast})$ and
$\vecy_2=\type(E^{(2)}_{\ast})$.
We also obtain the decomposition
$I_1\sqcup I_2=\Nbar=\{1,\ldots,N\}$:
\[
 I_j:=\bigl\{
 i\in \Nbar\,\big|\,
 \nbigf^{(j)}_i/\nbigf^{(j)}_{i-1}\neq 0
 \bigr\}
\]
The datum $(\vecy_1,\vecy_2,I_1,I_2)$
is called the decomposition type of $z$.
Thus, we prepare the following definition.
\begin{df}
\label{df;06.5.22.1}
\index{decomposition type}
A decomposition type is defined 
to be a datum $\gbigi:=(\vecy_1,\vecy_2,I_1,I_2)$
as follows:
\begin{itemize}
\item
 $\vecy=\vecy_1+\vecy_2$ in $\Type$
 such that
 $P_{\vecy_1}^{\alpha_{\ast},\delta}
 =P_{\vecy}^{\alpha_{\ast},\delta}$.
\item
 $\Nbar=I_1\sqcup I_2$
 such that  $|I_i|=H_{y_i}(m)$.
\end{itemize}
The set of the decomposition types 
is denoted by $\Dec(m,\vecy,\alpha_{\ast},\delta)$.
\index{$\Dec(m,\vecy,\alpha_{\ast},\delta)$}
\hfill\qed
\end{df}

We remark that the condition $O_m$-holds
for $\nbigm^{ss}(\vecy_1,L,\alpha_{\ast},\delta)$
and $\nbigm^{ss}(\vecy_2,\alpha_{\ast})$
for a decomposition type $(\vecy_1,\vecy_2,I_1,I_2)$,
if $\nbigm^{ss}(\vecy_1,L,\alpha_{\ast},\delta)
\times\nbigm^{ss}(\vecy_2,\alpha_{\ast})\neq \emptyset$.

\subsubsection{Statement}

\label{subsubsection;06.5.16.110}

Let $\gbigi:=(\vecy_1,\vecy_2,I_1,I_2)$ be a decomposition type.
Take a decomposition
$V=V^{(1)}\oplus V^{(2)}$
such that
$\dim V^{(i)}=H_{\vecy_i}(m)$.
We put $\proj^{(i)}:=\proj(V^{(i)\lor})$.
Then, we put
$Q^{(1)}:=
Q^{ss}\bigl(m,\vecy_1,[L],\alpha_{\ast},\delta\bigr)$
and
$Q^{(2)}:=Q^{ss}(m,\vecy_2,\alpha_{\ast})$.
We put as follows: 
\[
 \Flag^{(i)}:=\Flag(V^{(i)},I_i):=
\Bigl\{
 \nbigf^{(i)}_{\ast}\,\Big|\,
 \mbox{\rm filtration indexed by $\Nbar$},\,\,
 \dim\Gr^{\nbigf^{(i)}}_j=1\,\, (j\in I_i),
\mbox{\rm or } =0\,\,(j\not\in I_i)
\Bigr\}.
\]
Clearly, $\Flag^{(i)}$ are isomorphic to
the full flag varieties of $V^{(i)}$.
We put $\Qtilde^{(1)}:=Q^{(1)}\times\Flag^{(1)}$,
$\Qtilde^{(2)}:=Q^{(2)}\times\Flag^{(2)}$,
and $\Qtilde_{\spl}(\gbigi)
:=\Qtilde^{(1)}\times\Qtilde^{(2)}$.
Then, we have the naturally defined morphism
$\Qtilde_{\spl}(\gbigi)\lrarr \Qtilde$.
We put 
$\TH_{\spl}(\gbigi):=
 \TH\times_{\Qtilde}\Qtilde_{\spl}(\gbigi)$,
$\TH^{\ast}_{\spl}(\gbigi):=
 \TH^{\ast}\times_{\Qtilde}\Qtilde_{\spl}(\gbigi)$
and $\TH_{i,\spl}(\gbigi):=
 \TH_i\times_{\Qtilde}\Qtilde_{\spl}(\gbigi)$.
We have the closed immersion
$\iota:\TH_{\spl}(\gbigi)\lrarr \TH$.

Let $z=\bigl(
 (q_1,E^{(1)}_{\ast},[\phi],\nbigf^{(1)}),
 (q_2,E^{(2)}_{\ast},\nbigf^{(2)}),u \bigr)$
be a point of $\TH_{\spl}(\gbigi)$.
Let $\min(I_2)$ denote the minimum of $I_2$.
We will prove the following lemma, later
(the subsubsection \ref{subsubsection;06.6.21.20}).

\begin{lem}
\label{lem;06.5.16.11}
In the case (II) (the subsubsection {\rm\ref{subsubsection;06.5.22.369}}),
if $\iota(z)$ is contained in $\TH^{ss}$,
then we have $\min(I_2)> \ell$.
\end{lem}

We put $\nbigf^{(2)}_{\min}:=\nbigf^{(2)}_{\min(I_2)}$.
We remark $\nbigf^{(2)}_{\min-1}=0$,
and $\dim\nbigf^{(2)}_{\min}=1$.
We also remark that the pair
$\bigl(E^{(2)}_{\ast},\nbigf_{\min}^{(2)}\bigr)$
can be regarded as a reduced $\nbigo_X(-m)$-Bradlow pair
on $X$.
We will also prove the following proposition.

\begin{prop}
\label{prop;06.5.16.20}
$\iota(z)$ is contained in $\TH^{ss}$,
if and only if the following conditions hold:
\begin{itemize}
\item $z\in \TH^{\ast}_{\spl}$
\item
 $\bigl(E^{(2)}_{\ast},\nbigf^{(2)}_{\min}\bigr)$
 is an $\epsilon$-semistable
 reduced $[\nbigo(-m)]$-Bradlow pair
 for any sufficiently small $\epsilon>0$.
\item
 $(E^{(1)}_{\ast},[\phi],\nbigf^{(1)})$
 is $\bigl(\delta,\min(I_2)-1\bigr)$-semistable.
\end{itemize}
\end{prop}

\subsubsection{Step 1}

Let $G_1$ denote the subgroup of
$\GL(V_1)\times\GL(V_2)$
determined by $\det(g)=1$,
i.e.,
$G=\{(g_1,g_2)\in\GL(V_1)\times\GL(V_2)\,|\,
 \det(g_1)\cdot \det(g_2)=1\}$.

\begin{lem}
\label{lem;06.5.21.121}
The following two conditions are equivalent:
\begin{itemize}
\item
$\mu_{\lambda}\bigl(\iota(z),\nbigo_{\nbigbtilde}(1)\bigr)
\geq 0$ for any one parameter subgroup
$\lambda$ of $\SL(V)$.
\item
$\mu_{\lambda}\bigl(\iota(z),\nbigo_{\nbigbtilde}(1)\bigr)
\geq 0$ for any one parameter subgroup
$\lambda$ of $G_1$.
\end{itemize}
\end{lem}
\pf
The first condition clearly implies the second condition.
Let us show the reverse implication.
Assume the second condition holds.
Let $\lambda:G_m\lrarr\SL(V)$ be 
a one-parameter subgroup.
We have the decomposition
$\lambda=\lambda^{(1)}+\lambda^{(2)}$
such that
$\mu_{\lambda^{(1)}}\bigl(z,\nbigl\bigr)=0$
and 
$\mu_{\lambda^{(2)}}\bigl(z, \nbigl\bigr)>0$.
By the same argument as the proof of 
Lemma \ref{lem;06.4.22.150},
we have only to show the following inequalities:
\begin{equation}
 k\cdot \mu_{\lambda^{(1)}}(z,\nbigl)
+\max_{j=1,2}\bigl\{\gamma_j\cdot 
 \mu_{\lambda^{(1)}}(z,\nbigo_{\proj}(1))\bigr\}
+\mu_{\lambda^{(1)}}(z,\nbigo_{\Flag}(n_{\ast}))
>0.
\end{equation}
\begin{equation}
 \label{eq;06.4.22.160}
 k\cdot \mu_{\lambda^{(2)}}(z,\nbigl)
+\min_{j=1,2}\bigl\{\gamma_j\cdot 
 \mu_{\lambda^{(2)}}(z,\nbigo_{\proj}(1))\bigr\}
+\mu_{\lambda^{(2)}}(z,\nbigo_{\Flag}(n_{\ast}))
> 0.
\end{equation}
Since $n_i$ and $\gamma_j$ are sufficiently small,
we can show that
the inequality (\ref{eq;06.4.22.160}) always holds
by the same argument as the proof of
Lemma \ref{lem;06.5.16.2}.
So we may and will assume 
$\mu_{\lambda}(z,\nbigl)=0$.

Let $V=\bigoplus_i \nbigv_i$ denote 
the weight decomposition of $\lambda$.
We put $\nbigg_j=\bigoplus_{i\leq j}\nbigv_i$.
We have the number $i_0$
determined by
$\phi\in \nbigg_{i_0}-\nbigg_{i_0-1}$.
We have only to show the following:
\[
\max_{j=1,2}\bigl\{\gamma_j\cdot i_0\bigr\}
+\sum_j n_j\sum_{i}
 \dim\frac{\nbigf_j\cap\nbigg_i}{\nbigf_j\cap\nbigg_{i-1}}
\geq 0.
\]

We put $\nbigh_0=V^{(1)}$ and $\nbigh_1=V$,
which gives the filtration of $V$.
We have the natural identification
$\Gr^{\nbigh}(V)\simeq V$.
Since $\nbigf$ is compatible with
the decomposition $V=V^{(1)}\oplus V^{(2)}$,
the induced filtration by $\nbigf$
is same as $\nbigf$.
Let $\nbigg'$ denote the induced filtration
by $\nbigg$ on $\Gr^{\nbigh}(V)\simeq V$:
\[
 \nbigg_i'=\nbigg_i'\cap V^{(1)}\oplus \nbigg_i'\cap V^{(2)},
\quad
 \nbigg_i'\cap V^{(1)}=\nbigg_i\cap V^{(1)},
\quad
 \nbigg_i\cap V^{(2)}=\frac{\nbigg_i}{\nbigg_i\cap V^{(1)}}
\]
Because of $\phi\in V^{(1)}$,
we have $\phi\in\nbigg_{i_0}'-\nbigg_{i_0-1}'$.

Let us take any decomposition
$V^{(a)}=\bigoplus\nbigk_i^{(a)}$ $(a=1,2)$
such that
$\nbigg_j'\cap V^{(a)}=\bigoplus_{i\leq j}\nbigk_i^{(a)}$.
We put $\nbigk_i=\nbigk_i^{(1)}\oplus\nbigk_i^{(2)}$.
Then we obtain the one parameter subgroup
$\lambda'$ of $G_1$
whose weight decomposition is
$\bigoplus\nbigk_i$.

\begin{lem}
\label{lem;06.5.21.122}
We have $\mu_{\lambda'}\bigl(z,\nbigl\bigr)=0$.
\end{lem}
\pf
Let $\nbige_i(m)$ denote the subsheaf of $E(m)$
generated by $\nbigf_i$.
Each $\nbige_i$ has the induced parabolic structure
and the reduced $L$-section $[\phi_i]$.
Then the filtration 
$\cdots\subset
 (\nbige_{i\,\ast},[\phi_i])
\subset
 (\nbige_{i+1\,\ast},[\phi_{i+1}])
\subset\cdots$
is a partial Jordan-H\"older filtration
due to Corollary \ref{cor;06.5.21.155}.

Let us consider the filtration
$\nbigh$ of $E$
given by
$\nbigh_0=E^{(1)}$
and $\nbigh_1=E$.
On $\Gr^{\nbigh}(\nbige_i)$,
we have the induced parabolic structure
and the reduced $L$-section $[\Gr^{\nbigh}(\phi_i)]$.
Then, the tuple
$\bigl(\Gr^{\nbigh}(\nbige_i)_{\ast},\Gr^{\nbigh}(\phi_i)\bigr)$
is $\delta$-semistable.

We have the canonical isomorphism
$\Gr^{\nbigh}(E)\simeq E$.
We regard $\Gr^{\nbigh}(\nbige_i)(m)$
as the subsheaf of $E(m)$.
Then,
$\nbigg_i'$ generates $\Gr^{\nbigh}(\nbige_i)(m)$.
Then, we obtain $\mu_{\lambda'}(z,\nbigl)=0$
from Corollary \ref{cor;06.5.21.155}
\hfill\qed

\vspace{.1in}

Let us return to the proof of 
Lemma \ref{lem;06.5.21.121}.
Since we have assumed the second condition
in the lemma,
we have $\mu_{\lambda'}\bigl(z,\nbigo_{\nbigbtilde}(1)\bigr)\leq 0$.
Then, we obtain the following inequality
from Lemma \ref{lem;06.5.21.122}:
\[
\max_{j=1,2}\bigl\{\gamma_j\cdot i_0\bigr\}
+\sum_j n_j\sum_{i}
 \dim\frac{\nbigf_j\cap\nbigg'_i}{\nbigf_j\cap\nbigg'_{i-1}}
\geq 0.
\]
Hence we have only to show the following inequality:
\[
\sum_i i\cdot \dim
 \frac{\nbigf_j\cap \nbigg_i'}{\nbigf_j\cap\nbigg_{i-1}'}
\leq 
\sum_i i\cdot \dim
 \frac{\nbigf_j\cap \nbigg_i}{\nbigf_j\cap\nbigg_{i-1}}.
\]
It can be shown by the geometric argument
on the $\SL(V)$-action on the Grassmann variety
of $N-H_j(m)$-dimensional quotients of $V$.
But, we give a more elementary proof.
We put $M:=\max\{i\,|\,\nbigv^{(i)}\neq 0\}$,
and we put $H:=\nbigf_j$.
Then we have the following:
\begin{multline}
 \sum i\cdot\dim\frac{H\cap \nbigg_i}{H\cap\nbigg_{i-1}}
=
\sum_{i\leq M} 
 i\cdot\dim\frac{H\cap \nbigg_i}{H\cap\nbigg_{i-1}}
=\sum_{i\leq M}
 i\cdot \dim(H\cap \nbigg_i)
-\sum_{i\leq M-1}
(i+1)\dim H\cap \nbigg_i \\
=\dim H\cdot M-\sum_{i\leq M-1}\dim H\cap \nbigg_i.
\end{multline}
Hence we have only to show 
$\dim H\cap \nbigg_i\leq \dim H\cap\nbigg_i'$
for each $i$.
But we have the equality
$H\cap\nbigg_i\cap V^{(1)}
=H\cap \nbigg'_i\cap V^{(1)}$
and the following inclusion:
\[
 \frac{H\cap\nbigg_i}{H\cap\nbigg_i\cap V^{(1)}}
\subset
 H\cap \nbigg_i'\cap V^{(2)}.
\]
Hence, we obtain the desired inequality,
and we are done.
\hfill\qed

\subsubsection{Step 2}

We give some preliminary consideration.

(O1)
Assume that 
there exists a partial Jordan-H\"older filtration
$E^{(2)}_{1\,\ast}\subset E^{(2)}_{2\,\ast}$.
We take a decomposition
$V^{(2)}=V^{(2)}_1\oplus V^{(2)}_2$
such that
$V^{(2)}_1=H^0\bigl(X,E^{(2)}_1(m)\bigr)$,
and we consider the one-parameter subgroup
$\lambda$ given by
$t^{-\rank V_2^{(2)}}\id_{V_1^{(2)}}\,
\oplus\,
 t^{\rank V_1^{(2)}}\id_{V^{(2)}}$.
In the case, we have the following:
\begin{equation}
 \label{eq;06.5.16.10}
 \mu_{\lambda}\bigl(z,\nbigo_{\nbigbtilde}(1)\bigr)
=
 \sum_{j=\min(I_2)}^Nn_j\cdot
 \Bigl(
-\rank V^{(2)}_2\cdot
 \dim \nbigf^{(2)}_j\cap V^{(2)}_1
+\rank V^{(2)}_1\cdot
 \dim\frac{\nbigf^{(2)}_j}{\nbigf^{(2)}_j\cap V^{(2)}_1}
 \Bigr).
\end{equation}
Since $n_{\min(I_2)}$ is sufficiently larger than
$n_j$ $(j>\min(I_2))$,
(\ref{eq;06.5.16.10}) is larger than $0$,
if and only if
$\nbigf^{(2)}_{\min(I_2)}\cap V^{(2)}_1=\{0\}$.
We also remark that 
(\ref{eq;06.5.16.10}) cannot be $0$.

\vspace{.1in}

In particular,
we obtain the following:
\begin{lem}
\label{lem;06.5.16.25}
If $\iota(z)$ is contained in $\nbigbtilde^{ss}$,
then the reduced $\nbigo(-m)$-Bradlow pair
$(E_{\ast}^{(2)},\nbigf_{\min})$
is $\epsilon$-stable for any sufficiently small
number $\epsilon>0$.
\hfill\qed
\end{lem}

\subsubsection{Step 3}
\label{subsubsection;06.6.21.20}
We give some preliminary consideration.

(O2)
Let us consider the one-parameter subgroup $\lambda$
of $G_1$ given by
$t^{\rank V^{(2)}}\id_{V^{(1)}}\,
\oplus\, t^{-\rank V^{(1)}}\id_{V^{(2)}}$.
Let $\pi_i$ denote the projection of
$\TH_{\spl}$
onto 
$\TH_{\spl}\cap \TH_i$.
Let us consider points
$z_i:=\pi_i(z)$
which are fixed with respect to $\lambda$.
We have
$\mu_{\lambda}(z_i,\nbigl)=0$ for $i=1,2$.
We have $\mu_{\lambda}\bigl(z_2,
 \nbigo_{\proj_m}(\gamma)\bigr)
=\gamma_2\cdot\rank V^{(2)}$.
Since $n_{j}$ are sufficiently smaller
than $|\gamma_2|$,
we always have the inequality
$\mu_{\lambda}\bigl(z_2,
 \nbigo_{\nbigbtilde_{\spl}}(1)\bigr)<0$.

Let us consider the condition
$\mu_{\lambda}\bigl(z_1,\nbigo_{\nbigbtilde}(1)\bigr)\geq 0$.
It is equivalent to the following inequality:
\begin{equation}
 \label{eq;06.4.23.1}
 \mu_{\lambda}\bigl(z_1,\nbigo_{\proj_m}(\gamma_1)\bigr)
+\mu_{\lambda}\bigl(z_1,\nbigo_{\Flag}(n_{\ast})\bigr)\geq 0.
\end{equation}

In the case (I),
we have $\gamma_1>0$,
and $n_{i}$ are sufficiently smaller
than $\gamma_1$.
Hence, the inequality (\ref{eq;06.4.23.1})
always holds.
More strongly, the strict inequality holds.

In the case (II),
the inequality (\ref{eq;06.4.23.1}) can be rewritten
as follows:
\begin{multline}
 \label{eq;06.4.23.2}
 \gamma_1\cdot\rank V^{(2)}
+\sum_{i}n_i\cdot \left(
 \dim(\nbigf_i\cap V^{(1)})\cdot\rank V^{(2)}
-\dim\left(\frac{\nbigf_i}{\nbigf_i\cap V^{(1)}}\right)
 \cdot \rank V^{(1)}
 \right) \\
=\left(
 \gamma_1+\sum_{i=1}^{\ell}
 (\dim \nbigf_i\cap V^{(1)})\cdot n_i
 \right)\cdot\rank V^{(2)}
-\sum_{i=1}^{\ell}
 \dim\left(\frac{\nbigf_i}{\nbigf_i\cap V^{(1)}}\right)
 \cdot\rank V^{(1)}\cdot n_i \\
+\sum_{i\geq \ell+1}
 n_i\cdot \left(
 \dim(\nbigf_i\cap V^{(1)})\cdot\rank V^{(2)}
-\dim\left(\frac{\nbigf_i}{\nbigf_i\cap V^{(1)}}\right)
 \cdot \rank V^{(1)}
 \right)\geq 0
\end{multline}
The inequality (\ref{eq;06.4.23.2}) 
is equivalent to the condition
$\nbigf_{\ell}\subset V^{(1)}$,
due to our choice of $\gamma_1$ and $n_{\ast}$.
Moreover,
if the inequality holds,
the strict inequality holds due to
our choice of $n_{\ast}$.

\vspace{.1in}
Now, we give a proof of Lemma \ref{lem;06.5.16.11}.
If $\iota(z)$ is contained in $\nbigbtilde^{ss}$,
we have
$\mu_{\lambda}\bigl(z,\nbigo_{\nbigbtilde}(1)\bigr)
\geq 0$ for $\lambda$ as above.
Therefore, we obtain
the inequality
$\mu_{\lambda}\bigl(z_1,\nbigo_{\nbigbtilde}(1)\bigr)\geq 0$,
and hence 
$\nbigf_{\ell}\subset V^{(1)}$.
It means $\min(I_2)>\ell$.
Thus Lemma \ref{lem;06.5.16.11} is proved.

\subsubsection{Step 4}

We put $k:=\min(I_2)-1$
in the following argument.
We give some more preliminary considerations.

(O3)
Assume there exists a partial Jordan-H\"older filtration
$E^{(1)}_{1\,\ast}\subset (E^{(1)}_{\ast},\phi)$
with respect to $\delta$-semistability.
We take a decomposition
$V^{(1)}=V^{(1)}_1\oplus V^{(1)}_2$
such that
$V^{(1)}_1=H^0\bigl(X,E^{(1)}_1(m)\bigr)$,
and we consider the one-parameter subgroup
$\lambda$ given by
$t^{-\rank V^{(2)}}\id_{V_1^{(1)}}\,
\oplus\,
 t^{\rank V_1^{(1)}}\id_{V^{(2)}}$.
We have
$\mu_{\lambda}(z,\nbigl)
=\mu_{\lambda}\bigl(z,\nbigo_{\proj_m}(1)\bigr)
=0$ for such one parameter subgroups.
Therefore,
the condition
$\mu_{\lambda}\bigl(z,\nbigo_{\nbigbtilde}(1)\bigr)\geq 0$
is equivalent to the following inequality:
\begin{multline}
 \label{eq;06.4.23.3}
0\leq
\mu_{\lambda}\bigl(z,\nbigo_{\Flag}(n_{\ast})\bigr) 
=\sum_{i}n_i\cdot\left(
-\rank V^{(2)}\cdot
 \dim\bigl(\nbigf^{(1)}_i\cap V_1^{(1)}\bigr)
+\rank V_1^{(1)}\cdot
 \dim\nbigf^{(2)}_i
 \right)\\
=-\sum_{i=1}^{k-1}n_i\cdot
 \rank V^{(2)}\cdot
 \dim\bigl(\nbigf^{(1)}_i\cap V_1^{(1)}\bigr)
+n_k\cdot
 \left(
-\rank V^{(2)}\cdot
 \dim\bigl(\nbigf^{(1)}_k\cap V_1^{(1)}\bigr)
+\rank V_1^{(1)}\cdot
 \dim\nbigf^{(2)}_k
 \right)\\
+\sum_{i>k}n_i\cdot\left(
-\rank V^{(2)}\cdot
 \dim\bigl(\nbigf^{(1)}_i\cap V_1^{(1)}\bigr)
+\rank V_1^{(1)}\cdot
 \dim\nbigf^{(2)}_i
 \right).
\end{multline}
Recall that $n_{i}$ are sufficiently smaller than $n_{i-1}$.
Hence the inequality (\ref{eq;06.4.23.3}) holds,
if and only if
$\nbigf^{(1)}_{k-1}\cap V_1^{(1)}=\{0\}$.

\vspace{.1in}
(O4)
Assume that there exists a partial Jordan-H\"older filtration
$(E_{2\,\ast}^{(1)},\phi)\subset (E^{(1)}_{\ast},\phi)$
with respect to $\delta$-semistability.
We take a decomposition
$V^{(1)}=V^{(1)}_2\oplus V^{(1)}_1$
such that
$V^{(1)}_2=H^0\bigl(X,E^{(1)}_2(m)\bigr)$,
and we take the one-parameter subgroup $\lambda$
given by
$t^{\rank V^{(2)}}\id_{V^{(1)}_1}\,
\oplus\,
 t^{-\rank V_1^{(1)}}\cdot\id_{V^{(2)}}$.
Then, we have
$\mu_{\lambda}(z,\nbigl)
=\mu_{\lambda}\bigl(z,\nbigo_{\proj_m}(1)\bigr)
=0$ for such one-parameter subgroups.
Therefore,
the condition 
$\mu_{\lambda}\bigl(z,\nbigo_{\nbigbtilde}(1)\bigr)\geq 0$
is equivalent to the following conditions:
\begin{multline}
 \label{eq;06.4.23.5}
0\leq
 \sum_i n_i\left(
 \rank V^{(2)}\cdot \rank\left(
 \frac{\nbigf^{(1)}_i}{\nbigf^{(1)}_i\cap V_2^{(1)}}
 \right)
-\rank V_1^{(1)}\cdot \rank\nbigf^{(2)}_i
 \right) \\
=\sum_{i=1}^{k-1}
 n_i\cdot \rank V^{(2)}
\rank\left(
 \frac{\nbigf^{(1)}_i}{\nbigf^{(1)}_i\cap V_2^{(1)}}
 \right)
+n_k\left(
 \rank V^{(2)}\cdot \rank\left(
 \frac{\nbigf^{(1)}_k}{\nbigf^{(1)}_k\cap V_2^{(1)}}
 \right)
-\rank V_1^{(1)}\cdot \rank\nbigf^{(2)}_k
 \right) \\
+\sum_{i>k} n_i\left(
 \rank V^{(2)}\cdot \rank\left(
 \frac{\nbigf^{(1)}_i}{\nbigf^{(1)}_i\cap V_2^{(1)}}
 \right)
-\rank V_1^{(1)}\cdot \rank\nbigf^{(2)}_i
 \right) \\
\end{multline}
Since $n_i$ is sufficiently smaller than $n_{i-1}$,
the inequality (\ref{eq;06.4.23.5}) holds
if and only if 
$\nbigf^{(1)}_{k-1}\not\subset V^{(1)}_2$.

\vspace{.1in}
We obtain the following claim
from the above preliminary considerations.
\begin{lem}
\label{lem;06.5.16.26}
If $\iota(z)$ is contained in $\nbigbtilde^{ss}$,
the tuple $(E_{\ast},\phi,\nbigf^{(1)})$
is $\bigl(\delta,\min(I_2)-1\bigr)$-semistable.
\hfill\qed
\end{lem}

\subsubsection{End of the proof of 
 Proposition \ref{prop;06.5.16.20}}

When $\iota(z)$ is contained in $\nbigbtilde^{ss}$,
it is easy to see $z\in\TH^{\ast}_{\spl}$.
We have already seen the other two conditions
are satisfied
(Lemma \ref{lem;06.5.16.25} and
Lemma \ref{lem;06.5.16.26}).

Let $z\in\TH_{\spl}$ be a point which satisfies
the conditions in Proposition \ref{prop;06.5.16.20}.
Let $u_1,\ldots,u_{N^{(1)}}$ be a base of $V^{(1)}$,
and let $u_{N^{(1)}+1},\ldots,u_{N}$ be 
a base of $V$.
Let $(w_1,\ldots,w_N)$ be an element of
$\seisuu^N$
such that $w_i\leq w_{i+1}$
and $\sum w_i=0$.
Let $\lambda$ be the one-parameter subgroup
of $G$, given by
$\lambda(t)\cdot u_i=t^{w_i}\cdot u_i$.
We do not distinguish $\lambda$ and $(w_1,\ldots,w_N)$.
We have the following:
\[
 \mu_{\lambda}\bigl(z,\nbigo_{\nbigbtilde}(1)\bigr)
=k\cdot k'\cdot \mu_{\lambda}\bigl(z,\nbigl\bigr)
+k'\mu_{\lambda}\bigl(z,\nbigo_{\Flag}(n_{\ast})\bigr)
+k'\max_{i=1,2}\Bigl\{
 \gamma_i\cdot \mu_{\lambda}(z,\nbigo_{\proj_m}(1))
 \Bigr\}.
\]
We put
$S_1:=\bigl\{
 h\,\big|\,\mu_{f_h}(z,\nbigl)=0 \bigr\}$
and $S_2:=\bigl\{
 h\,\big|\,\mu_{f_h}(z,\nbigl)>0
 \bigr\}$,
where
$f_h=\bigl(\overbrace{h-N,\ldots,h-N}^h,h,\ldots,h\bigr)$.
Since $n_{j}$ and $\gamma_i$
are sufficiently small,
we can show that 
the following inequality holds
for  any $h\in S_2$
by the same argument as the proof of 
Lemma \ref{lem;06.5.16.2}:
\begin{equation}
 \label{eq;06.4.26.1}
 k\cdot \mu_{f_h}(z,\nbigl)
+\mu_{f_h}(z,\nbigo_{\Flag}(n_{\ast}))
+\min_{i=1,2}\bigl\{
 \gamma_i\cdot\mu_{f_h}(z,\nbigo_{\proj_m}(1))
 \bigr\}>0.
\end{equation}
By the same argument as the proof of
Lemma \ref{lem;06.4.22.150},
it can be shown that
we have only to show the following 
inequalities,
for any $0\neq \rho=\sum_{j\in S_1}a_j\cdot f_j$
 with $a_j\geq 0$:
\begin{equation}
\label{eq;06.4.26.2}
F(\rho):=
k\cdot \mu_{\rho}(z,\nbigl)
+\mu_{\rho}(z,\nbigo_{\Flag}(n_{\ast}))
+\max_{i=1,2}\bigl\{
 \gamma_i\cdot\mu_{\rho}(z,\nbigo_{\proj_m}(1))
 \bigr\}>0.
\end{equation}
Let us show (\ref{eq;06.4.26.2}).
We have the weight decomposition
$V^{(\alpha)}=
 \bigoplus_j^{s(\alpha)} \nbigv^{(\alpha)}_j$ of $\rho$.
We put $N^{(\alpha)}:=\dim V^{(\alpha)}$.
We put $r^{(\alpha)}_j:=\dim \nbigv^{(\alpha)}_j$.
We use the notation in the subsubsection
\ref{subsubsection;06.4.23.15}.
Then $\rho$ can be expressed
as $\sum a_j^{(\alpha)}\cdot v^{(\alpha)}_j$
satisfying (\ref{eq;06.4.22.3}).
Let $i_0$ be the number determined by
$\phi\in \bigoplus_{i\leq i_0} \nbigv_{i}
 -\bigoplus_{i\leq i_0-1}\nbigv_i$.
Due to Lemma \ref{lem;06.4.23.10},
we have the expression:
\begin{equation}
\rho=\sum_{j=1}^{s(2)}c(j)\cdot y^{(2)}(j)
 +\sum_{i<i_0}d_1(i)\cdot x_1(i)
 +\sum_{i>i_0}d_2(i)\cdot x_2(i)
 +A\cdot\bigl(N^{(2)}\Omega^{(1)}
 -N^{(1)}\cdot \Omega^{(2)}\bigr).
\end{equation}
Here $c(j)$, $d_1(i)$ and $d_2(i)$
are non-negative numbers,
and $A$ is a rational number.
One of $c(j)$, $d_1(i)$, $d_2(i)$ or $A$
is not zero.
Due to $\mu_{\kappa}(z,\nbigo_{\proj_m}(1))=0$
for $\kappa=y^{(2)}(j), x_1(i), x_2(i)$,
we have the following linearity:
\[
 F(\rho)
=\sum_{j=1}^{s(2)}c(j)\cdot F\bigl(y^{(2)}(j) \bigr)
+\sum_{i<i_0}d_1(i)\cdot F\bigl(x_1(i)\bigr)
+\sum_{i>i_0}d_2(i)\cdot F\bigl(x_2(i)\bigr)
+F\Bigl(
A\cdot\bigl(N^{(2)}\Omega^{(1)}
-N^{(1)}\cdot\Omega^{(2)}\bigr)
 \Bigr).
\]
We obtain $F\bigl(y^{(2)}(j)\bigr)>0$,
$F\bigl(x_1(i)\bigr)>0$
and 
$F\bigl(x_2(i)\bigr)>0$
from the preliminary considerations
(O1), (O3) and (O4) respectively.
We have
$F\Bigl(
A\cdot\bigl(N^{(2)}\Omega^{(1)}
-N^{(1)}\cdot\Omega^{(2)}\bigr)
 \Bigr)>0$ in the case $A\neq 0$,
from the preliminary consideration (O2).
Therefore,
we obtain the desired positivity,
and the proof of Proposition \ref{prop;06.5.16.20}
is finished.
\hfill\qed

\subsection{Enhanced Master Space in the Oriented Case}
\label{subsection;06.5.11.100}
\subsubsection{The construction}

\label{subsubsection;06.5.17.10}

We put
$\widehat{Q}:=
 Q^{ss}(\vecyhat,[L],\alpha_{\ast},\delta)$
and 
$\widehat{\widetilde{Q}}:=
 \widehat{Q}\times\Flag(V,\Nbar)$.
Then, we put
$\THhat:=\TH\times_{\Qtilde}\Qtildehat$
and
$\THhat^{\ast}:=\TH^{\ast}\times_{\TH}\THhat$.
We remark that
$\THhat^{\ast}$ is isomorphic to
$Q^{ss}(m,\vecyhat,L,\alpha_{\ast},\delta)
\times\Flag(V,\Nbar)$.
We also put
$\THhat^{ss}=\TH^{ss}\times_{\TH}\THhat$.
We have the natural $\GL(V)$-action
on $\THhat^{ss}$.
The quotient stack is called
the enhanced master space
\index{enhanced master space}
\index{master space}
in the oriented case,
and it is denoted by $\Mhat$.
\index{$\Mhat$}
It is also called the master space
for abbreviation.

\begin{prop}
\label{prop;06.6.9.10}
$\Mhat$ is Deligne-Mumford and proper.
\end{prop}
\pf
From Proposition \ref{prop;06.5.16.50},
the stack $\TH^{ss}/\PGL(V)$ is Deligne-Mumford
and proper.
We have the naturally defined morphism
$\Mhat\lrarr \TH^{ss}/\PGL(V)$,
which is etale and finite.
Then, we obtain that $\Mhat$ is
also Deligne-Mumford and proper.
\hfill\qed

\vspace{.1in}
We have the universal quotient objects 
$\bigl(q^u,\nbige^u,F_{\ast},[\phi^u],\rho^u\bigr)$
on $\Qhat\times X$.
It induces the oriented reduced $L$-Bradlow pair
$\pi_X^{\ast}\bigl(\nbige^u(-m),
 F_{\ast},[\phi^u],\rho^u \bigr)$
on $\THhat\times X$,
where $\pi_X:\THhat\times X\lrarr \Qhat\times X$
denotes the natural projection.
By taking the descent with respect to the $\GL(V)$-action,
we obtain the oriented reduced $L$-Bradlow pair
$(\Ehat^{\Mhat},F_{\ast}^{\Mhat},[\phi^{\Mhat}],\rho^{\Mhat})$.
We also have the induced full flag
$\nbigf^{\Mhat}$ of $p_{X\,\ast}\Ehat^{\Mhat}(m)$.

The $G_m$ action $\rho$ on $\TH$
induces the $G_m$-action on $\THhat$,
which is also denoted by $\rho$.
It induces the $G_m$-action on $\Mhat$,
which is denoted by $\rhobar$.
Let us see the stack theoretic fixed point set
of $\Mhat$ in the following subsubsections.

\subsubsection{The obvious fixed point sets}
\label{subsubsection;06.5.17.25}

We put $\THhat^{ss}_i=\TH_i\times_{\TH}\THhat^{ss}$.
We have the substacks $\Mhat_i:=\THhat_i^{ss}/\GL(V)$
$(i=1,2)$. \index{$\Mhat_i$}
The stack $\Mhat_2$ can be easily described.
By construction,
$\Mhat_2$ gives the moduli stack 
$\nbigmtilde^{s}(\vecyhat,[L],\alpha_{\ast},\delta_-)$
of the tuples
$(E_{\ast},[\phi],\rho,\nbigf)$,
where 
 $(E_{\ast},[\phi],\rho)$ is an oriented
 $\delta_-$-stable reduced $L$-Bradlow pairs 
 of type $\vecy$,
and $\nbigf$ is a full flag of
 $H^0(X,E(m))$.
It is easily related with 
the moduli stack $\nbigm^s(\vecyhat,[L],\alpha_{\ast},\delta_-)$
of $\delta_-$-semistable oriented
reduced $L$-Bradlow pairs.
We have the universal sheaf
$\Ehat^{u}$ over
$\nbigm^s(\vecyhat,[L],\alpha_{\ast},\delta_-)\times X$.
We obtain the locally free sheaf
$p_{X\,\ast}\Ehat^u(m)$ on
$\nbigm^s(\vecyhat,[L],\alpha_{\ast},\delta_-)$.
The associated full flag bundle
is isomorphic to $\Mhat_2$.

Let us see $\Mhat_1$.
In the case (I),
$\Mhat_1$ is isomorphic to 
the moduli stack of the tuples
$(E_{\ast},[\phi],\rho,\nbigf)$,
where
 $(E_{\ast},[\phi],\rho)$
 is a $\delta_+$-semistable oriented 
 reduced $L$-Bradlow pair of type $\vecy$,
and  $\nbigf$ is a full flag of $H^0(X,E(m))$.
It is related with 
the moduli stack
$\nbigm^s\bigl(\vecyhat,[L],\alpha_{\ast},\delta_+\bigr)$,
as in the case of $\Mhat_2$.
In the case (II),
$\Mhat_1$ is the moduli stack of 
$(\delta,\ell)$-semistable tuples
$(E_{\ast},[\phi],\rho,\nbigf)$.

We have the restriction of the oriented reduced $L$-Bradlow pair
$(\Ehat^{\Mhat},F_{\ast}^{\Mhat},[\phi^{\Mhat}],\rho^{\Mhat})$
to $\Mhat_i\times X$.
We also have the universal objects on $\Mhat_i\times X$
by the moduli theoretic meaning of $\Mhat_i$.
It is clear that they are isomorphic,
by the construction of $\Ehat^{\Mhat}$.
It is also easy to observe that
the weight of $\rhobar$ on
$\Ehat^{\Mhat}_{|\Mhat_i\times X}$ is $0$.

\subsubsection{Fixed point sets associated to
 decomposition types}
\label{subsubsection;06.5.17.1}

Let us describe the components of
the fixed point set contained in 
$\Mhat^{\ast}:=\Mhat-(\Mhat_1\cup\Mhat_2)$.
\index{$\Mhat^{\ast}$}
We use the notation in the subsubsection 
\ref{subsubsection;06.5.16.110}.
Let $\gbigi=(\vecy_1,\vecy_2,I_1,I_2)$ be
a decomposition type.
In the case (II) (the subsubsection \ref{subsubsection;06.5.22.369}),
we assume $\ellbar\subset I_1$.
We put $k:=\min(I_2)-1$.
Let $\Qtilde^{(1)}(\delta,k)$ denote
the maximal open subset of $\Qtilde^{(1)}$
determined by the $(\delta,k)$-semistability.
(It is open due to Proposition \ref{prop;06.5.15.70}.)
Let $\Qtilde^{(2)}_+$ denote the maximal open subset 
determined by the condition that
the reduced $\nbigo(-m)$-Bradlow pair
$\bigl(E^{(2)}_{\ast},\nbigf_{k+1}\bigr)$
is $\epsilon$-stable for any sufficiently small $\epsilon>0$.
We have the naturally defined morphism
$\Qtilde^{(1)}(\delta,k)\times
 \Qtilde^{(2)}_+\lrarr
 \Qtilde$.
Let $\Qtildehat_{\spl}(\gbigi)$
denote the fiber product of
$\Qtilde^{(1)}(\delta,k)\times\Qtilde^{(2)}_+$
and $\Qtildehat$ over $\Qtilde$.
We also put
$\THhat^{ss}_{\spl}(\gbigi):=
 \THhat^{\ast}\times_{\Qtildehat}\Qtildehat_{\spl}(\gbigi)$.
We have the natural 
$\GL(V^{(1)})\times\GL(V^{(2)})$-action
on $\THhat^{ss}_{\spl}(\gbigi)$.
The quotient stack is denoted by
$\Mhat^{G_m}(\gbigi)$. \index{$\Mhat^{G_m}(\gbigi)$}
Due to Proposition \ref{prop;06.5.16.20},
we have the naturally defined morphism
$\varphi_{\gbigi}:\Mhat^{G_m}(\gbigi)\lrarr \Mhat$.

\begin{lem}
\label{lem;06.5.17.22}
\mbox{{}}
\begin{enumerate}
\item
$\Mhat^{G_m}(\gbigi)$ is isomorphic to the moduli stack
of the objects
$\bigl(
 (E^{(1)}_{\ast},\phi,\nbigf^{(1)}),
 (E^{(2)}_{\ast},\nbigf^{(2)}),\rho
 \bigr)$ as follows:
\begin{itemize}
\item
 $(E^{(1)}_{\ast},\phi)$ is $\delta$-semistable
 parabolic $L$-Bradlow pair of type $\vecy_1$,
and $\nbigf^{(1)}$ is a full flag of 
 $p_{X\,\ast}E^{(1)}(m)$ such that
 $(E^{(1)}_{\ast},\phi,\nbigf^{(1)})$ is 
 $(\delta,k)$-semistable.
\item
 $E^{(2)}_{\ast}$ is semistable parabolic sheaf
 of type $\vecy_2$,
and
 $\nbigf^{(2)}$ is a full flag of
 $p_{X\,\ast}E^{(2)}(m)$ such that
 $(E^{(2)},\nbigf^{(2)}_{\min})$ is
 a reduced $\nbigo(-m)$-Bradlow pair. 
\item
 $\rho$ is an orientation of 
 $E^{(1)}\oplus E^{(2)}$.
\end{itemize}

\item
We have the decomposition
$\varphi_{\gbigi}^{\ast}\Ehat^{\Mhat}_{\ast}
=E^{\Mhat}_1\oplus E^{\Mhat}_2$
and $\varphi_{\gbigi}^{\ast}\nbigf^{\Mhat}
=\nbigf^{\Mhat}_1\oplus\nbigf^{\Mhat}_2$.
The pull back of
the reduced $L$-section $[\phi^{\Mhat}]$ of $\Ehat^{\Mhat}$
naturally gives the $L$-section $\phi_1^{\Mhat}$ of $E_1^{\Mhat}$.

Then, the object
$\bigl(
 (E^{\Mhat}_{1\,\ast},\phi^{\Mhat}_1,\nbigf^{\Mhat}_1),
 (E^{\Mhat}_{2\,\ast},\nbigf^{\Mhat}_2),
 \varphi_{\gbigi}^{\ast}\rho^{\Mhat}
\bigr)$ gives the universal object
over $\Mhat^{G_m}\times X$
in the moduli theoretic meaning above.
\end{enumerate}
\end{lem}
\pf
We recall the fibration $\TH^{\ast}\lrarr \Qtilde$ 
is $\nbigo_{\proj_m}(-1)^{\ast}$.
Then, the claim is clear by construction.
\hfill\qed

\vspace{.1in}

We will later give a more convenient description
of $\Mhat^{G_m}(\gbigi)$ and the restriction
of $\bigl(\Ehat^{\Mhat},F_{\ast},[\phi],\rho\bigr)$
to $\Mhat^{G_m}(\gbigi)$.

\subsubsection{Statement}

We put as follows:
\[
 \nbigs(m,\vecy):=\left\{
 \begin{array}{ll}
 \Dec(m,\vecy,\alpha_{\ast},\delta) & \mbox{\rm (the case (I))} \\
 \mbox{{}}\\
 \bigl\{\gbigi=(\vecy_1,\vecy_2,I_1,I_2)\in\Dec(m,\vecy,\alpha_{\ast},\delta)
 \,\big|\,\ellbar\subset I_1
 \bigr\} & (\mbox{\rm the case (II)})
 \end{array}
 \right.
\]

\begin{prop}
\label{prop;06.5.16.150}
Let $\Mhat^{G_m}$ denote the stack theoretic
fixed point set of $\Mhat$
with respect to the action $\rhobar$.
Then we have the following isomorphism:
\[
 \Mhat^{G_m}\simeq
 \Mhat_1\sqcup\Mhat_2\sqcup
 \coprod_{\gbigi\in\nbigs(m,\vecy)}
 \Mhat^{G_m}(\gbigi).
\]
\end{prop}

We will prove the proposition 
in the following subsubsections.

\subsubsection{Ambient stack}

We have the Poincar\'{e} bundle
$\Poin$ on $\Pic\bigl(\det(y(m))\bigr)\times X$.
Then,
we obtain the vector bundle
$\Zhat_m:=p_{X\,\ast}\bigl(
 \nhom\bigl(\bigwedge^r V_{m,X},\Poin\bigr)
 \bigr)$ 
on $\Pic\bigl(\det(y(m))\bigr)$.
Recall that the Gieseker space
is the projectivization of $\Zhat_m$.
We have the natural projection
$\nbigbtilde\lrarr Z_m$.
Thus, we put
$\nbigbtildehat^{ss}:=
 \Zhat_m\times_{Z_m}\nbigbtilde^{ss}$.
We have the naturally defined $\GL(V_m)$-action
on $\nbigbtildehat^{ss}$.
The quotient stack is denoted by
$\gbigb'$.
The $G_m$-action $\rho$ on $\nbigbtilde^{ss}$
induces the action on $\nbigbtildehat^{ss}$.
Since it commutes with the action of $\GL(V)$,
we obtain the $G_m$-action $\rhobar$
on $\gbigb'$.

We have the naturally defined closed immersion
$\THhat^{ss}\lrarr \nbigbtildehat^{ss}$,
which is $\GL(V_m)\times G_m$-equivariant.
Therefore,
we obtain the closed $G_m$-equivariant
immersion $\Mhat\lrarr \gbigb'$.
Since $\Mhat$ is Deligne-Mumford and separated,
we can take an equivariant open neighbourhood
$\gbigb$ of $\Mhat$ in $\gbigb'$,
which is Deligne-Mumford and separated.
We remark that $\gbigb$ is also smooth.
Let $\gbigb^{G_m}$ denote the fixed point set of $\gbigb$
with respect to $G_m$,
i.e., the $0$-set of the vector field
induced by the $G_m$-action.
Recall that $\Mhat^{G_m}$
is defined to be
$\gbigb^{G_m}\times_{\gbigb}\Mhat$.

\subsubsection{Fixed point set of the ambient space}

We also put
$\nbigbtildehat_i^{ss}:=
 \nbigbtilde_i\times_{\nbigbtilde}\nbigbtildehat^{ss}$.
The quotient stack
$\nbigbtildehat_i^{ss}$ is denoted by
$\gbigb_i$.
It gives an open subset of $\gbigb^{G_m}$.
It is easy to see
$\Mhat_i=\gbigb_i\times_{\gbigb}\Mhat$.

We put $N=H_y(m)$.
Let $\vecN=(N_1,N_2)$, $\vecr=(r_1,r_2)$,
$\veck_{\ast}:=\bigl(k_{1,j},k_{2,j}\,\big|\,j=1,\ldots,l\bigr)$,
$\vecI:=(I_1,I_2)$ be a datum as follows:
\begin{itemize}
\item
 $N_i$ are positive integers such that
 $N_1+N_2=N$.
\item
 $r_i$ are positive integers such that $r_1+r_2=r$.
\item
 $k_{i,j}$ are positive integers such that
 $k_{1,j}+k_{2,j}=H_{y}(m)-H_{\vecy,j}(m)$.
\item
 $I_1\sqcup I_2=\Nbar$.
\end{itemize}

Such a tuple
$\gbigq:=(\vecN,\vecr,\veck_{\ast},\vecI)$
is called a decomposition type for $\nbigatilde$.
We remark that a decomposition type
$\gbigi=(\vecy_1,\vecy_2,I_1,I_2)$
induces the decomposition type $\gbigq(\gbigi)$
for $\nbigatilde$, as follows:
\[
 N_i=|I_i|,\quad
 r_i=\rank(\vecy_i),\quad
 k_{i,j}=H_{y_i}(m)-H_{\vecy_i,j}(m).
\]

We put $\nbigatildehat:=\nbigatilde\times_{Z_m}\Zhat_m$.
For a decomposition type
$\gbigq=(\vecN,\vecr,\veck_{\ast},\vecI)$ for $\nbigatilde$,
the locally closed regular subvariety 
$\gbigc_1(\gbigq)$ of $\nbigatildehat$ is 
the set of the points
$(\gminif,K_{\ast},[\phi],\nbigf)$ satisfying
the following conditions:
\begin{itemize}
\item
 There exists the unique
 decomposition $V=V^{(1)}\oplus V^{(2)}$.
\item
 $\gminif\in \Zhat_m$ is contained in
 $H^0\Bigl(
 X, \nhom\Bigl(
 \bigwedge^{r_1}V_X^{(1)}\otimes\bigwedge^{r_2}V_X^{(2)},\,
 L \Bigr)
 \Bigr)$
 for some line bundle $L\in \Pic\bigl(\det(y(m))\bigr)$,
 where $V^{(i)}_X:=V^{(i)}\otimes\nbigo_X$.
\item
 $K_{\ast}=\bigl(K_i\,\big|\,i=1,\ldots,l\bigr)\in
 \prod_{i} Q_{m,i}$ is compatible
 with the decomposition $V=V^{(1)}\oplus V^{(2)}$,
 i.e., we have the decomposition
 $K_j=K_{j}^{(1)}\oplus K_j^{(2)}$
 such that $K_j^{(i)}$ are quotients of $V^{(i)}$.
 We also assume
 $\dim K_j^{(i)}=k_{i,j}$.
\item
 $[\phi]$ is contained in
 the projectivization of $V^{(1)}$.
\item
 We have the decomposition
 $\nbigf=\nbigf^{(1)}\oplus \nbigf^{(2)}$
 compatible with the decomposition
 $V=V^{(1)}\oplus V^{(2)}$.
 Moreover,
 $\nbigf^{(i)}_j/\nbigf^{(i)}_{j-1}\neq 0$
 if and only if
 $j\in I^{(i)}$.
\end{itemize}

We put $\nbigbtilde^{\ast}:=
 \nbigbtilde-\bigcup\nbigbtilde_i$.
We put
$\gbigc_2(\gbigq):=\nbigbtilde^{\ast}\times_{\nbigatilde}
 \gbigc_1(\gbigq)$.
We have the natural $\GL(V)$-action
on $\gbigc_2(\gbigq)$,
and the quotient stack is denoted by
$\gbigc'_3(\gbigq)$,
which is the closed substack of $\gbigb'$.
We put $\gbigc_3(\gbigq)=\gbigc_3'(\gbigq)\cap\gbigb$.
The following lemma can be checked easily.
\begin{lem}
 $\gbigc_3(\gbigq)$ are open subsets 
of $\gbigb^{G_m}$.
\hfill\qed
\end{lem}

\subsubsection{Proof of Proposition \ref{prop;06.5.16.150}}

We put as follows:
\[
 \gbigb^{G_m}_0:=
 \coprod_{\gbigi\in\nbigs(m,\vecy)}
 \gbigc_3\bigl(\gbigq(\gbigi)\bigr)
\]
We have the natural morphism
$\Mhat^{G_m}(\gbigi)\lrarr
 \gbigc_3\bigl(\gbigq(\gbigi)\bigr)$.
Therefore,
we have the morphism
$\Mhat^{G_m}(\gbigi)
\lrarr \gbigb_0^{G_m}$.
We obtain the following morphism:
\[
\psi_1:
 \coprod \Mhat^{G_m}(\gbigi)
\lrarr
 \Mhat\times_{\gbigb}
 \gbigb_0^{G_m} 
\]

\begin{lem}
$\psi_1$ is isomorphic.
\end{lem}
\pf
Let $g:T\lrarr \Mhat\times_{\gbigb}\gbigb_0^{G_m}$
be a morphism.
Then we have the $\GL(V)$-equivariant torsor $P(g)$.
We also obtain the following data
from $g:T\lrarr \gbigb_0^{G_m}$.
\begin{itemize}
\item
We have the $\GL(V)$-equivariant decomposition
$\nbigv_1\oplus\nbigv_2$
of $p_{P(g)}^{\ast}V_{X}$ over $X\times P(f)$.
\item
 $\vecf:\bigwedge^r\nbigv\lrarr\det^{\ast}\Poin$,
 which is contained in 
$\Hom\bigl(
\bigwedge^{r_1}\nbigv_1\otimes
 \bigwedge^{r_2}\nbigv_2,\,
\det^{\ast}\Poin\bigr)$.
\end{itemize}
The composite
 $T\lrarr \Mhat\times_{\gbigb}\gbigb_0^{G_m}
 \lrarr \Mhat$
is contained in $\Mhat^{\ast}$,
and hence
in $Q^{ss}(m,\vecyhat,L,\alpha_{\ast},\delta)/\GL(V_m)$.
Therefore,
we obtain the following data:
\begin{itemize}
\item 
 The oriented quotient parabolic 
 $L$-Bradlow pair
 $\bigl(q,E_{\ast}(m),\rho,\phi\bigr)$
\item
 $q:p_{P(g)}^{\ast}V_{m,X}\lrarr E(m)$
 satisfies (TFV)-condition.
\end{itemize}
We have 
$\rho\circ\bigwedge^rq=\gminif$.
Then, we obtain the decomposition
$E=E^{(1)}\oplus E^{(2)}
=q\bigl(\nbigv^{(1)}\bigr)
\oplus q\bigl(\nbigv^{(2)}\bigr)$.
The claim is clear on the open subset
where $E$ is locally free.
Since $E$ is torsion-free and 
$q$ is surjective,
the decomposition is obtained on the whole space.

By taking the decent with respect to
the $\GL(V)$-action,
we obtain 
$\bigl((q_1,E^{(1)}_{\ast},\phi_1,\nbigf^{(1)}),
 (q_2,E^{(2)}_{\ast},\nbigf^{(2)}),\rho\bigr)$
on $T$.
We remark that the decomposition data
is determined on each connected component of $T$.
The conditions in Proposition \ref{prop;06.5.16.20}
is satisfied for the specialization of 
$(q_1,E^{(1)}_{\ast},\phi_1,\nbigf^{(1)})$
and $ (q_2,E^{(2)}_{\ast},\nbigf^{(2)})$
to any closed fibers $\{u\}\times X$.
Hence we obtain the morphism
$T\lrarr \coprod \Mhat^{G_m}(\gbigi)$.
In particular,
we obtain the morphism
$\psi_2:\Mhat\times_{\gbigb}\gbigb^{G_m}_0
\lrarr \coprod\Mhat^{G_m}(\gbigi)$.
It is easy to see that
$\psi_1$ and $\psi_2$ are mutually inverse.
\hfill\qed

\vspace{.1in}

Then, we obtain the following:
\[
 \Mhat_1\sqcup\Mhat_2\sqcup
\coprod\Mhat^{G_m}(\gbigi)
=
 \Mhat^{G_m}\times_{\gbigb}
 \bigl(
 \gbigb_1\cup\gbigb_2\cup\gbigb^{G_m}_0
 \bigr)
=\Mhat^{G_m}\times_{\gbigb}
 \gbigb^{G_m}
\]

In particular,
$\Mhat_1\sqcup\Mhat_2\sqcup
 \coprod \Mhat^{G_m}(\gbigi)$
is the closed substack of $\Mhat$.
Since $\gbigb(\gbigq)$ and $\gbigb_i$ 
are open in $\gbigb^{G_m}$,
it is easy to see that
$\Mhat_i$ and $\Mhat^{G_m}(\gbigi)$
are unions of connected components of
the fixed point set.
Thus, the proof of Proposition \ref{prop;06.5.16.150}
is finished.
\hfill\qed

\subsection{Decomposition into Product of Two Moduli Stacks}
\label{subsection;06.6.21.3}
\subsubsection{Statement}
\label{subsubsection;06.5.17.20}

Let $\gbigi=(\vecy_1,\vecy_2,I_1,I_2)$ be
a decomposition type.
We would like to decompose
$\Mhat^{G_m}(\gbigi)$ into the product of
two moduli stacks
up to etale finite morphisms.
We put $k:=\min(I_2)-1$.
We use the notation in the subsubsection 
\ref{subsubsection;06.5.17.1}.
We introduce some more moduli stacks.

Let $\nbigmtilde^{ss}\bigl(
 \vecyhat_1,[L],\alpha_{\ast},(\delta,k)\bigr)$
denote the moduli stack of the objects
$(E^{(1)}_{\ast},[\phi^{(1)}],\rho^{(1)},\nbigf^{(1)})$
as follows:
\begin{itemize}
\item
 $(E^{(1)}_{\ast},[\phi^{(1)}],\rho^{(1)})$ 
is a $\delta$-semistable oriented reduced
 $L$-Bradlow pair of type $\vecy_1$,
 and 
 $\nbigf^{(1)}$ is a full flag of
 $p_{X\,\ast}E^{(1)}(m)$
 such that
$(E^{(1)}_{\ast},[\phi^{(1)}],\nbigf)$ 
is $(\delta,k)$-semistable.
\end{itemize}
The universal object over
$\nbigmtilde^{ss}\bigl(\vecyhat_1,[L],\alpha_{\ast},(\delta,k)\bigr)\times X$
is denoted by
$\bigl(\Ehat_{1\ast}^u,\phi_1^u,\rho^u_1,\nbigf^u_1\bigr)$.

\vspace{.1in}
Let $\nbigmtilde^{ss}\bigl(
 \vecy_1,L,\alpha_{\ast},(\delta,k)\bigr)$
denote the moduli stack of the objects
$(E^{(1)}_{\ast},\phi^{(1)},\nbigf^{(1)})$
as follows:
\begin{itemize}
\item
 $(E^{(1)}_{\ast},\phi^{(1)})$ 
is a $\delta$-semistable $L$-Bradlow pair of type $\vecy_1$
such that $\phi^{(1)}$ is non-trivial everywhere,
and  $\nbigf^{(1)}$ is a full flag of $p_{X\,\ast}E^{(1)}(m)$
such that
$(E^{(1)}_{\ast},\phi^{(1)},\nbigf)$ 
is $(\delta,k)$-semistable.
\end{itemize}
The universal object over
$\nbigmtilde^{ss}\bigl(\vecy_1,L,\alpha_{\ast},(\delta,k)\bigr)\times X$
is denoted by $(E_{1\,\ast}^u,\phi^u,\nbigf^u_1)$.

\vspace{.1in}

Let $\nbigmtilde^{ss}\bigl(\vecyhat_2,
 \alpha_{\ast},+\bigr)$ denote the moduli stack
of the objects $(E^{(2)}_{\ast},\rho^{(2)},\nbigf^{(2)})$
as follows:
\begin{itemize}
\item
 $(E^{(2)}_{\ast},\rho^{(2)})$ is semistable
 oriented parabolic sheaf of type $\vecy_2$,
and 
 $\nbigf^{(2)}$ is a full flag of
 $p_{X\,\ast}E^{(2)}(m)$.
\item
 $(E^{(2)}_{\ast},\nbigf^{(2)}_{\min})$
 is an $\epsilon$-semistable
 reduced $\nbigo(-m)$-Bradlow pair
 for any sufficiently small $\epsilon>0$.
\end{itemize}
The universal object over
$\nbigmtilde^{ss}(\vecyhat_2,\alpha_{\ast},+)\times X$
is denoted by
$(\Ehat^u_2,\rho_2^u,\nbigf^{u}_2)$.

\begin{prop}
\label{prop;06.5.17.30}
We put $r_i:=\rank\vecy_i$.
There exists the algebraic stack $\nbigs$
with the following properties:
\mbox{{}}
\begin{itemize}
\item
There exist the following diagram:
\begin{equation}
 \label{eq;06.5.23.20}
\begin{CD}
 \Mhat^{G_m}(\gbigi)@<{F}<<
 \nbigs
 @>{G}>>
 \nbigmtilde^{ss}\bigl(\vecyhat_1,[L],\alpha_{\ast},(\delta,k)\bigr)
 \times
 \nbigmtilde^{ss}\bigl(\vecyhat_2,\alpha_{\ast},+\bigr)
\end{CD}
\end{equation}
The morphisms $F$ and $G$ are etale and finite of degree 
$(r_1\cdot r_2)^{-1}$ and $r_2^{-1}$, respectively.
\item
Let $\nbigo_{1,\rel}(1)$ denote the tautological line bundle
of $\nbigmtilde^{ss}(\vecyhat_1,[L],\alpha_{\ast},\delta)$.
There exists the line bundle
$\nbigo_{1,\rel}(1/r_2)$ on $\nbigs$ such that
     $\nbigo_{1,\rel}(1/r_2)^{r_2}=G^{\ast}\nbigo_{1,\rel}(1)$,
and we have the following relations:
\begin{equation}
 \label{eq;06.5.17.15}
 F^{\ast}E_1^{\Mhat}
\simeq G^{\ast}\Ehat_1^{u}\otimes\nbigo_{1,\rel}(1),
\quad
 F^{\ast}E_2^{\Mhat}
\simeq G^{\ast}\Ehat_2^u\otimes\nbigo_{1,\rel}(-r_1/r_2)
\end{equation}
\end{itemize}
The weight of the $G_m$-action $\rhobar$ on $E_1^{\Mhat}$
and $E_2^{\Mhat}$ are 
$-1$ and $r_1/r_2$, respectively.
\end{prop}

\begin{cor}
\label{cor;06.5.22.50}
We have the following diagram:
\begin{equation}
 \label{eq;06.5.17.21}
\begin{CD}
 \Mhat^{G_m}(\gbigi)@<{F}<<
 \nbigs
 @>{G'}>>
 \nbigmtilde^{ss}\bigl(\vecy_1,L,\alpha_{\ast},(\delta,k)\bigr)
 \times
 \nbigmtilde^{ss}\bigl(\vecyhat_2,\alpha_{\ast},+\bigr)
\end{CD}
\end{equation}
Here $G'$ is etale and finite of degree $(r_1\cdot r_2)^{-1}$.
We have the following relations
\[
 F^{\ast}E_1^{\Mhat}\simeq
 G^{\prime\ast}E_1^u,\quad
 F^{\ast}E_2^{\Mhat}\simeq
 G^{\prime\ast}\Ehat_2^u\otimes\Or(E_1^u)^{-1/r_2}
\]
Here, we put $\Or(E_1^u)^{-1/r_2}:=\nbigo_{1,\rel}(-r_1/r_2)$.
\hfill\qed
\end{cor}

Before going into the proof,
we give some remark on the inductive process
to investigate the transition of the moduli stacks,
heuristically.
Let $\delta$ be an element of 
$\Cr(\vecy,L,\alpha_{\ast})$.
Let $\delta_-$ and $\delta_+$ be sufficiently close to
$\delta$ such that $\delta_-<\delta<\delta_+$.
We will be interested in the difference of
$\nbigm^{ss}(\vecyhat,[L],\alpha_{\ast},\delta_-)$
and $\nbigm^{ss}(\vecyhat,[L],\alpha_{\ast},\delta_+)$.
We make the enhanced master space.
Then $\Mhat_i$ are isomorphic to the full flag variety bundles
over $\nbigmtilde^{ss}(\vecyhat,[L],\alpha_{\ast},\delta_-)$
and $\nbigm^{ss}(\vecyhat,[L],\alpha_{\ast},\delta_+)$.
So we can derive some information
from the fixed point sets $\Mhat^{G_m}(\gbigi)$
by the localization technique,
and $\Mhat^{G_m}(\gbigi)$ is isomorphic to
$\nbigmtilde^{ss}(\vecyhat_1,[L],\alpha_{\ast},\delta,\ell)
\times
 \nbigmtilde^{ss}(\vecyhat_2,\alpha_{\ast},+)$
up to finite and etale morphisms,
where $\ell=\min(I_2)-1$.

The structure of
$\nbigmtilde^{s}(\vecyhat_2,\alpha_{\ast},+)$
can be easily related with the moduli stack
$\nbigm=
 \nbigmtilde^{s}(\vecyhat_2,[\nbigo(-m)],\alpha_{\ast},\epsilon)$,
where $\epsilon$ is a sufficiently small number.
We have the universal oriented sheaf
$\Ehat^u$ over $\nbigm\times X$
with the universal reduced section $[\phi_1^u]$.
Then we obtain the vector bundle
$p_{X\,\ast}\bigl(\Ehat^u(m)\bigr)$
with the line subbundle
$\nbigq\subset p_{X\,\ast}\bigl(
\Ehat^u(m)\bigr)$ induced by $[\phi^u]$.
We put $\nbigc:=p_{X\,\ast}\bigl(\Ehat^u(m)\bigr)/\nbigq$.
Then the associated full flag bundle to $\nbigc$
is isomorphic to
$\nbigmtilde^{ss}
 \bigl(\vecyhat_2,[\nbigo(-m)],\alpha_{\ast},\epsilon\bigr)$.

On the other hand,
the structure of
$\nbigmtilde^{ss}\bigl(\vecyhat_1,[L],\alpha_{\ast},\delta,\ell\bigr)$
is not so easy to describe.
However, we can make the enhanced master space $\Mhat^{(1)}$
again,
so that $\Mhat^{(1)}_1$ and $\Mhat^{(1)}_2$
are isomorphic to
$\nbigmtilde^{ss}\bigl(\vecyhat_1,[L],\alpha{\ast},\delta,\ell\bigr)$
and the full flag variety bundle over
$\nbigmtilde^{ss}\bigl(
\vecyhat_1,[L],\alpha_{\ast},\delta_-\bigr)$.
Thus, we can proceed inductively.

We remark $\rank(\vecy_1)<\rank(\vecy)$.
Therefore, the process will stop,
and we will arrive the description of the difference of
$\nbigm^{ss}(\vecyhat,[L],\alpha_{\ast},\delta_+)$
and $\nbigm^{ss}(\vecyhat,[L],\alpha_{\ast},\delta_-)$
in terms of the products of the moduli stacks
of semistable objects with lower ranks.
We use such an argument in the subsection 
\ref{subsection;06.6.21.30}.

\subsubsection{Preliminary}

We use the notation
in the subsection \ref{subsection;06.6.6.2}.
Let $\Qtildehat^{(1)}$ denote the maximal open subset
of $Q^{\circ}(m,\vecyhat_1,L)\times\Flag^{(1)}$
determined by the $(\delta,k)$-semistability.
Let $\Qtildehat^{(2)}$ denote the maximal open subset
of $Q^{\circ}(m,\vecyhat_2)\times \Flag^{(2)}$,
which consists of the points
$(q_2,E_{2\,\ast},\rho_2,\nbigf^{(2)})$
such that $(E_{2\,\ast},\nbigf^{(2)}_{\min})$
is $\epsilon$-semistable for any small $\epsilon>0$.
We also put $T_1:=\Spec k[t_1,t_1^{-1}]$.
We have the $T_1$-action
on $\Qtildehat^{(1)}\times \Qtildehat^{(2)}$, 
given as follows:
\[
 t_1\cdot\Bigl(
 \bigl(E_{1\,\ast},\phi,\rho,\nbigf^{(1)}\bigr),
 \bigl(E_{2\,\ast},\phi,\rho,\nbigf^{(2)}\bigr)
\Bigr)
=\Bigl(
 \bigl(E_{1\,\ast},\phi,\,t_1\!\cdot\! \rho,\,\nbigf^{(1)}\bigr),
 \bigl(E_{2\,\ast},\phi,\,t_1^{-1}\!\cdot\!\rho,\,\nbigf^{(2)}\bigr)
\Bigr)
\]
By construction,
we have the isomorphism
$ \THhat^{ss}_{\spl}(\gbigi)
\simeq
\bigl(\Qtildehat^{(1)}\times \Qtildehat^{(2)}\bigr) /T_1$.
We have the naturally defined actions
$\GL(V_i)$ on $\Qtildehat^{(i)}$.
We put $\nbigm_i:=\Qtildehat^{(i)}/\GL(V_i)$.
Then we obtain the following description:
\[
 \Mhat^{G_m}(\gbigi)
\simeq
\frac{\Qtildehat^{(1)}\times \Qtildehat^{(2)}}
{\GL(V_1)\times \GL(V_2)\times T_1}
\simeq
\frac{\nbigm_1\times\nbigm_2}{T_1}
\]

Let us see the right hand side more closely.
The stack $\nbigm_1$ is isomorphic to
the moduli stack of the objects
$(E^{(1)}_{\ast},\phi^{(1)},\rho^{(1)},\nbigf^{(1)})$
as follows:
\begin{itemize}
\item
 $(E^{(1)}_{\ast},\phi^{(1)},\rho^{(1)})$ 
is a $\delta$-semistable oriented
 $L$-Bradlow pair of type $\vecy_1$,
and 
 $\nbigf^{(1)}$ is a full flag of
 $p_{X\,\ast}E^{(1)}(m)$
such that
 $(E^{(1)}_{\ast},\phi^{(1)},
 \nbigf^{(1)})$ is $(\delta,k_0)$-semistable.
\end{itemize}
The quotient stack $\nbigm_2$ is
isomorphic to the moduli stack 
$\nbigmtilde^{ss}\bigl(\vecyhat_2,
 \alpha_{\ast},+\bigr)$.

The $T_1$-action on $\nbigm_1$ is given by
$t_1\cdot
 \bigl(E^{(1)}_{\ast},\phi^{(1)},\rho^{(1)},\nbigf^{(1)}\bigr)
=\bigl(E^{(1)}_{\ast},\phi^{(1)},
 t_1\!\cdot\! \rho^{(1)},\nbigf^{(1)}\bigr)$.
The $T_1$-action on $\nbigm_2$
is given by
$t_1\cdot \bigl(E^{(2)}_{\ast},\rho^{(2)},\nbigf^{(2)}\bigr)
=\bigl(E^{(2)}_{\ast},t_1^{-1}\!\cdot\!\rho^{(2)},
 \nbigf^{(2)}\bigr)$.

\subsubsection{Construction of $\nbigs$}

We put  $\Ttilde_1:=\Spec k[s_1,s_1^{-1}]$.
Let $\Ttilde_1\lrarr T_1$ be the morphism
given by $t_1=s_1^{r_1r_2}$,
where $r_i=\rank \vecy_i$.
We have the naturally induced $\Ttilde_1$-action
on $\nbigm_1\times\nbigm_2$.
Let $\nbigs$ denote the quotient stack
$(\nbigm_1\times\nbigm_2)/\Ttilde_1$.
Then, we have the following morphism:
\[
\begin{CD}
 \nbigs
 @>{F}>>
 (\nbigm_1\times\nbigm_2)/T_1
\simeq
 \Mhat^{G_m}(\gbigi)
\end{CD}
\]
Here $F$ is etale and finite of degree $(r_1\cdot r_2)^{-1}$.

Let us see the $\Ttilde_1$-action on $\nbigm_i$ $(i=1,2)$.
The induced $\Ttilde_1$-action on
$\nbigm_2$ is trivial,
i.e.,
$s_1\cdot \bigl(E^{(2)}_{\ast},\rho^{(2)},\nbigf^{(2)}\bigr)
=\bigl(E^{(2)}_{\ast},s_1^{-r_1\cdot r_2}\!\cdot\!\rho^{(2)},\nbigf^{(2)}\bigr)
\simeq
 \bigl(E^{(2)}_{\ast},\rho^{(2)},\nbigf^{(2)}\bigr)$.
The isomorphism is given by the following diagram:
\begin{equation}
 \label{eq;06.4.26.100}
\begin{CD}
 E^{(2)} \\
 @V{s_1^{-r_1}}VV\\
 E^{(2)}
\end{CD}
\mbox{{}}
\quad\quad\quad\quad
 \begin{CD}
 \det (E^{(2)}) @>{s_1^{-r_1r_2}\!\cdot\!\rho^{(2)}}>> \det^{\ast}\poin\\
 @V{s_1^{-r_1r_2}}VV @V{\id}VV\\
 \det(E^{(2)})@>{\rho^{(2)}}>> \det^{\ast}\poin
 \end{CD}
\end{equation}

The induced $\Ttilde_1$-action on $\nbigm_1$
is given as follows:
\[
 s_1\cdot\bigl(E^{(1)}_{\ast},
 \phi^{(1)},\rho^{(1)},\nbigf^{(1)}\bigr)
=\bigl(E^{(1)}_{\ast},\phi^{(1)},
 s_1^{r_1r_2}\!\cdot\!\rho^{(1)}, \nbigf^{(1)}\bigr)
\simeq
 \bigl(E^{(1)}_{\ast},s_1^{r_2}\!\cdot\!\phi^{(1)},\rho^{(1)},
 \nbigf^{(1)}\bigr).
\]
The isomorphism is given by the following diagram:
\begin{equation}
 \label{eq;06.4.26.101}
 \begin{CD}
 L @>>> E^{(1)} @. \mbox{{\hspace{1cm}}}
@.\det(E^{(1)})@>{s_1^{r_1r_2}\rho^{(1)}}>>
 \det^{\ast}\poin\\
 @V{\id}VV @V{s_1^{r_2}}VV@. @V{s_1^{r_1r_2}}VV
 @V{\id}VV\\
 L @>{s_1^{r_2}\phi}>> E^{(1)} @.\mbox{{}}@.
 \det(E^{(1)})@>{\rho^{(1)}}>>\det^{\ast}\poin
 \end{CD}
\end{equation}

In particular,
since $\Ttilde_1$-action on $\nbigm_2$ is trivial,
we obtain the following morphism $G$:
\[
\begin{CD}
 \nbigs=
 (\nbigm_1/\Ttilde_1)\times\nbigm_2
@>{G}>>
 (\nbigm_1/T_1)\times\nbigm_2
=
 \nbigmtilde^{ss}\bigl(\vecyhat_1,[L],\alpha_{\ast},
 (\delta,k_0)\bigr)
\times
 \nbigmtilde^{ss}(\vecy_2,\alpha_{\ast},I_2) 
\end{CD}
\]
The morphism $G$ is etale and finite
of degree $1/r_2$.

\subsubsection{The universal sheaf}

We put $\nbigmbar_1:=
 \nbigmtilde^{ss}\bigl(\vecyhat_1,[L],
 \alpha_{\ast},(\delta,k_0)\bigr)$.
We remark that
$\pi:\nbigm_1\times\nbigm_2
\lrarr \nbigmbar_1\times\nbigm_2$
is isomorphic to 
$\nbigo_{1,\rel}(-1)^{\ast}
\lrarr \nbigmbar_1\times\nbigm_2$,
and the $\Ttilde_1$-action
on $\nbigo_{1,\rel}(-1)^{\ast}$
is given by $s\cdot v=s^{r_2}v$
on each fiber.
We use the argument
in the subsection \ref{subsection;06.5.17.5}.
Let $\nbigt(n)$ denote the trivial line bundle
on $\nbigmbar_1\times\nbigm_2$
with the $\Ttilde_1$-action of weight $n$.
It induces the $\Ttilde_1$-line bundle
$\pi^{\ast}\nbigt(n)$ over $\nbigm_1\times\nbigm_2$.
By the descent,
we obtain the line bundle
$\nbigi_n$ over $\nbigm_1/\Ttilde_1\times\nbigm_2$.
We put $\nbigo_{1,\rel}(1/r_2):=\nbigi_1$.
It satisfies
 $\nbigo_{1,\rel}(1/r_2)^{r_2}\simeq 
 G^{\ast}\nbigo_{1,\rel}(1)$,
due to Lemma \ref{lem;06.5.17.6}.

Let $\Ehat^{\prime\,u}_1$
and $\Ehat^{\prime\,u}_2$
denote the pull back of $\Ehat_1^u$
via the morphism
$\nbigm_1\times\nbigm_2\times X\lrarr
 \nbigmbar_1\times\nbigm_2\times X$.
By the construction,
we have the $\Ttilde_1$-equivariant
sheaves $\Ehat^{\,\prime\TH}_1\oplus \Ehat_2^{\prime\,\TH}$ on
$\nbigm_1\times \nbigm_2\times X$,
which is induced by the sheaf on $\THhat^{ss}$
as in the subsubsection \ref{subsubsection;06.5.17.10}.
When we take it into account of the $G_m$-action $\rhobar$,
we have
$\Ehat_1^{\prime\TH}\simeq\Ehat^{\prime\,u}_1\otimes\pi^{\ast}\nbigt(r_2)$
and 
$\Ehat_2^{\prime\TH}\simeq\Ehat^{\prime\,u}_2\otimes\pi^{\ast}\nbigt(-r_1)$
due to the diagrams (\ref{eq;06.4.26.100}) and (\ref{eq;06.4.26.101}).
Therefore,
we obtain (\ref{eq;06.5.17.15}).

We put $\Ttilde_2:=\Spec k[s_2,s_2^{-1}]$,
and let $\Ttilde_2\lrarr G_m$ is a morphism
given by $t=s_2^{r_2}$.
We have the induced action $\rhotilde$
of $\Ttilde_2$ on $\THhat^{ss}$.
On $\Qtildehat^{(1)}\times \Qtildehat^{(2)}$,
the induced action is given as follows:
\[
 s_2\cdot \bigl(E_{1\,\ast},\phi,\rho_1,\nbigf_1,
 E_{2\,\ast},\rho_2,\nbigf_2\bigr)
=\bigl(
 E_{1\,\ast},s_2^{r_2}\phi,\rho_1,\nbigf_1,
E_{2\,\ast},\rho_2,\nbigf_2 \bigr).
\]
Therefore,
the action of $\Ttilde_1$ and $\Ttilde_2$
on $\nbigm_1\times\nbigm_2$ are same.
On the other hand,
the induced bundles over $\nbigm_1\times\nbigm_2\times X$
are $\Ehat_1^{\prime\,u}$ and $\Ehat_2^{\prime\,u}$
with respect to the $\rhotilde$-action.
Therefore,
the weight of $\rhotilde$
on $\varphi^{\ast}\Ehat_1^u\otimes\nbigi_{r_2}$
is $-r_2$,
and the weight of $\rhotilde$
on $\varphi^{\ast}\Ehat_2^u\otimes\nbigi_{-r_1}$
is $r_1$, due to Lemma \ref{lem;06.5.17.17}.
Thus, the proof of Proposition \ref{prop;06.5.17.30}
is finished.
\hfill\qed

\subsection{Simple Cases}
\label{subsection;06.7.3.15}
\subsubsection{The case where the $2$-stability condition is satisfied}
\label{subsubsection;06.5.17.100}

We give some indication
about what happens when the $2$-stability condition
is satisfied for $(\vecy,L,\alpha_{\ast},\delta)$,
without proof.
In this case,
we do not have to consider the enhanced master space
and $(\delta,\ell)$-semistability.
Hence the problem is simpler.

We use the notation in the subsubsections
\ref{subsubsection;06.5.14.30}
and \ref{subsubsection;06.6.21.35}.
We take a positive rational number $\gamma_1$
and a negative rational number $\gamma_2$
such that $|\gamma_i|$ are sufficiently small.
We take a large rational number $k'$
such that $k'\cdot(\gamma_1-\gamma_2)=1$.
We put $\nbigl_i:=\nbigl_{\gamma_i}^{\otimes k'}$.
We have $\nbigl_2=\nbigl_1\otimes\nbigo_{\proj_m}(-1)$.

Let us consider $\nbigb:=\proj\bigl(\nbigo_{\proj_m}(0)
 \oplus\nbigo_{\proj_m}(1)\bigr)$ over $\nbiga$. 
We put $\nbigb_1=\proj\bigl(\nbigo_{\proj_m}(0)\bigr)$
and $\nbigb_2=\proj\bigl(\nbigo_{\proj_m}(1)\bigr)$,
which are naturally regarded as the closed subscheme of $\nbigb$.
We have the tautological line bundle
$\nbigo_{\proj}(1)$,
and we put $\nbigo_{\nbigb}(1):=\nbigo_{\proj}(1)\otimes\nbigl_1$.
We have the natural $\SL(V)$-action on $\nbigb$,
and $\nbigo_{\nbigb}(1)$ gives the equivariant polarization.
Let $\nbigb^{ss}$ denote the set of the semistable points
of $\nbigb$ with respect to $\nbigo_{\nbigb}(1)$.
We put $\TH^{ss}:=Q\times_{\nbiga}\nbigb^{ss}$.
We have the natural $\SL(V)$-action on $\TH^{ss}$.
The following proposition can be shown
by an argument similar to the proof of
Proposition \ref{prop;06.5.16.50}.
In fact, it is much simpler.
\begin{prop}
\label{prop;06.5.17.50}
The quotient stack $\TH^{ss}/\SL(V)$
is Deligne-Mumford.
\hfill\qed
\end{prop}

We put $\THhat^{ss}:=\TH^{ss}\times_{Q}\Qhat$,
where we put $\Qhat:=Q^{ss}(\vecyhat,[L],\alpha_{\ast},\delta)$.
We put $\Mhat:=\THhat^{ss}/\GL(V)$,
which is called the master space in the oriented case.
\index{master space}
\index{$\Mhat$}
We have the $G_m$-action $\rhobar$
as in the enhanced case.
From Proposition \ref{prop;06.5.17.50},
we obtain the following.
\begin{prop}
$\Mhat$ is Deligne-Mumford and proper.
\hfill\qed
\end{prop}

We put $\THhat^{ss}_i:=\nbigb^{ss}_i\times_{\nbiga}\THhat^{ss}$,
and we put $\Mhat_i:=\THhat^{ss}_i/\GL(V)$.
\index{$\Mhat_i$}
We put $\Mhat^{\ast}:=\Mhat-\Mhat_1\cup\Mhat_2$.
\index{$\Mhat^{\ast}$}
It is easy to observe that
$\Mhat^{\ast}$ is an open substack of
$\nbigm(m,\vecyhat,L)$.

Due to Lemma \ref{lem;06.5.15.55},
$\Mhat_1$ and $\Mhat_2$ are isomorphic to
$\nbigm^{ss}(\vecyhat,[L],\alpha_{\ast},\delta_+)$
and $\nbigm^{ss}(\vecyhat,[L],\alpha_{\ast},\delta_-)$
respectively.
They give the obvious fixed point sets of
$\Mhat$ with respect to $\rhobar$.

Let us see the other components of the fixed point set.
A decomposition type 
is defined to be $\gbigi:=(\vecy_1,\vecy_2)\in\Type^2$
satisfying $\vecy_1+\vecy_2=\vecy$
and $P^{\alpha_{\ast},\delta}_{\vecy}=
 P^{\alpha_{\ast},\delta}_{\vecy_1}=P^{\alpha_{\ast}}_{\vecy_2}$.
For such a decomposition type $\gbigi:=(\vecy_1,\vecy_2)$,
we consider the moduli stack $\Mhat^{G_m}(\gbigi)$
of objects
$\bigl(E_{\ast}^{(1)},\phi,E_{\ast}^{(2)},\rho\bigr)$
as follows:
\begin{itemize}
\item
 $(E_{\ast}^{(1)},\phi)$ is $\delta$-stable
 $L$-Bradlow pair of type $\vecy_1$.
\item
 $E_{\ast}^{(2)}$ is stable parabolic sheaf
 of type $\vecy_2$.
\item
 $\rho$ is an orientation of $E^{(1)}\oplus E^{(2)}$.
\end{itemize}

Note that the $2$-stability condition
for $(\vecy,L,\alpha_{\ast},\delta)$
implies the $1$-stability conditions
for $(\vecy_1,L,\alpha_{\ast},\delta)$
and $(\vecy_2,\alpha_{\ast})$
if $\nbigm^{ss}(\vecy_1,L,\alpha_{\ast},\delta)
\times\nbigm^{ss}(\vecy_2,\alpha_{\ast})\neq \emptyset$.
We have the naturally defined morphism
$\Mhat^{G_m}(\gbigi)\lrarr \Mhat$
as in the subsubsection \ref{subsubsection;06.5.17.1}.

Let $S(\vecy,\alpha_{\ast},\delta)$ denote the set of decomposition types.
Then we can show the following 
by the same argument as the proof of 
Proposition \ref{prop;06.5.16.150}.
\begin{prop}
$\Mhat_1\sqcup \Mhat_2\sqcup
 \coprod_{\gbigi \in S(\vecy,\alpha_{\ast},\delta)}
 M^{G_m}(\gbigi)$
is the stack theoretic fixed point set
of $\Mhat$ with respect to $\rhobar$.
\hfill\qed
\end{prop}

We naturally have the oriented reduced $L$-Bradlow pair
$(\Ehat^{\Mhat}_{\ast},[\phi^{\Mhat}],\rho^{\Mhat})$
on $\Mhat\times X$, as in the subsubsection
\ref{subsubsection;06.5.17.10}.
The restriction of 
$(\Ehat^{\Mhat}_{\ast},[\phi^{\Mhat}],\rho^{\Mhat})$ to
$\Mhat_i\times X$ has the universal property
with respect to the moduli theoretic meaning above.
Let $\varphi:\Mhat\lrarr \nbigm(m,\vecyhat,[L])$ denote 
the naturally defined morphism.
Then the restriction of 
$\varphi^{\ast}\nbigo_{\rel}(1)$ to $\Mhat^{\ast}$
is canonically trivialized.
Therefore, the restriction of $[\phi^{\Mhat}]$
to $\Mhat^{\ast}$ gives the $L$-section,
which we denote by $\phi^{\Mhat}$.
Then, the restriction of $(\Ehat^{\Mhat},\phi^{\Mhat},\rho^{\Mhat})$
to $\Mhat^{G_m}(\gbigi)$ has the universal property
with respect to the moduli theoretic meaning above.
Correspondingly, we have the decomposition
$\Ehat^{\Mhat}_{|\Mhat^{G_m}(\gbigi)}=
 E_1^{\Mhat}\oplus E_2^{\Mhat}$.

\vspace{.1in}

It is convenient to decompose
$\Mhat^{G_m}(\gbigi)$ into the product
of two moduli stacks up to etale finite morphisms.
By the same argument as the proof of
Proposition \ref{prop;06.5.17.30},
we obtain the following description
of $\Mhat^{G_m}(\gbigi)$
up to etale finite morphisms.
\begin{prop}
\label{prop;06.6.10.5}
We put $r_i:=\rank\vecy_i$.
There exists the algebraic stack $\nbigs$
with the following properties:
\mbox{{}}
\begin{itemize}
\item
There exists the following diagram:
\[
\begin{CD}
 \Mhat^{G_m}(\gbigi)@<{F}<<
 \nbigs
 @>{G}>>
 \nbigm^{ss}\bigl(\vecyhat_1,[L],\alpha_{\ast},\delta\bigr)
 \times
 \nbigm^{ss}\bigl(\vecyhat_2,\alpha_{\ast}\bigr)
\end{CD}
\]
The morphisms $F$ and $G$ are etale and finite of degree 
$(r_1\cdot r_2)^{-1}$ and $r_2^{-1}$, respectively.
We also have the following diagram:
\begin{equation}
\label{eq;06.5.22.300}
\begin{CD}
  \Mhat^{G_m}(\gbigi)@<{F}<<
 \nbigs
 @>{G'}>>
 \nbigm^{ss}\bigl(\vecy_1,L,\alpha_{\ast},\delta\bigr)
 \times
 \nbigm^{ss}\bigl(\vecyhat_2,\alpha_{\ast}\bigr)
\end{CD}
\end{equation}
Here, $G'$ is etale and finite of degree $(r_1r_2)^{-1}$.
\item
Let $\nbigo_{1,\rel}(1)$ denote the tautological line bundle
of $\nbigm^{ss}(\vecyhat_1,[L],\alpha_{\ast},\delta)$.
We use the same notation to denote the pull back via 
an appropriate morphism.
Then, there exists the line bundle
$\nbigo_{1,\rel}(1/r_2)$ on $\nbigs$
such that $\nbigo_{1,\rel}(1/r_2)^{r_2}=G^{\ast}\nbigo_{1,\rel}(1)$,
and we have the following relation:
\begin{equation}
 F^{\ast}E_1^{\Mhat}
\simeq G^{\ast}\Ehat_1^{u}\otimes\nbigo_{1,\rel}(1),
\quad
 F^{\ast}E_2^{\Mhat}
\simeq G^{\ast}\Ehat_2^u\otimes\nbigo_{1,\rel}(-r_1/r_2)
\end{equation}
Here $\Ehat_1^u$ and $\Ehat_2^u$
are induced by the universal sheaves
over $\nbigm^{ss}(\vecyhat_1,[L],\alpha_{\ast},\delta)\times X$
and $\nbigm^{ss}(\vecyhat_2,\alpha_{\ast})\times X$,
respectively.
We also have the following relation:
\[
 F^{\ast}E_1^{\Mhat}
\simeq
 G^{\prime\,\ast}E_1^u,
\quad
 F^{\ast}E_2^{\Mhat}
\simeq
 G^{\prime\ast}\Ehat_2^u\otimes\Or(E_1^u)^{-1/r_2}
\]
Here, $E_1^u$ denotes the pull back 
of the universal sheaf over
$\nbigm^{ss}(\vecy_1,L,\alpha_{\ast},\delta)\times X$,
and 
$\Or(E_1^u)^{-1/r_2}$ denotes
$\nbigo_{1,\rel}(-r_1/r_2)$.
\item
The weights of the $G_m$-action
$\rhobar$ on $E_1^{\Mhat}$
and $E_2^{\Mhat}$ are 
$-1$ and $r_1/r_2$, respectively.
\hfill\qed
\end{itemize}
\end{prop}

\subsubsection{The case of oriented reduced $\vecL$-Bradlow pairs}
\label{subsubsection;06.5.21.555}

Let $\vecL=(L_1,L_2)$ be a pair of line bundles on $X$.
Let $\vecdelta=(\delta_1,\delta_2)$ be an element of
$\bigl(\nbigp^{\br}\bigr)^2$.
We can discuss the transition of the moduli stacks
$\nbigm^{ss}(\vecy,\vecL,\alpha_{\ast},\vecdelta)$,
when $\delta_1$ is moved.
For simplicity,
we restrict ourselves to the case
where both of $\delta_i$ are sufficiently small.
Recall the results in the subsubsection
\ref{subsubsection;06.6.15.5}.
Then, the $2$-stability condition is always satisfied,
and the problem is simple as in the subsubsection 
\ref{subsubsection;06.5.17.100}.
We give only an indication without proof.

Let $\vecdelta=(\delta_1,\delta_2)$ be an element
of $\nbigp^{\br\,2}$
such that $\delta_i$ are sufficiently small.
Assume that the $1$-stability condition
does not hold for $(\vecy,\vecL,\alpha_{\ast},\vecdelta)$.
In that case,
we have the positive integers $r_i$ $(i=1,2)$
satisfying the following:
\[
 r_1+r_2=r,
\quad
 \frac{\delta_1}{r_1}=\frac{\delta_2}{r_2}
\]
We take $\delta_{1,-}$ and $\delta_{1,+}$
such that $\delta_{1,-}<\delta_1<\delta_{1,+}$,
which are sufficiently close to $\delta_1$.
We put $\vecdelta_{\kappa}:=(\delta_{1,\kappa},\delta_2)$
for $\kappa=\pm$.
We would like to compare the moduli stacks
$\nbigm^{ss}(\vecy,\vecL,\alpha_{\ast},\vecdelta_-)$
and $\nbigm^{ss}(\vecy,\vecL,\alpha_{\ast},\vecdelta_+)$.

We use the notation in the subsubsection
\ref{subsubsection;06.5.14.30}.
We put $\nbiga:=\nbiga_m(\vecy,[\vecL])$ and
$\nbigl:=\nbigl_{\vecy,\vecL}\bigl(P^{\alpha_{\ast},\vecdelta}_{\vecy}(m),
 \epsilon_{\ast},\vecdelta(m)\bigr)$,
which gives a $\GL(V)$-polarization on $\nbiga$.
We put
$\nbigl_{\gamma}:=\nbigl\otimes\nbigo_{\proj_m^{(1)}}(\gamma)$
for rational number $\gamma$.
Let $\nbiga^{ss}(\nbigl_{\gamma})$
denote the set of the semistable points of $\nbiga$
with respect to $\nbigl_{\gamma}$.

We put $Q:=Q^{ss}(m,\vecy,[\vecL],\alpha_{\ast},\vecdelta)$.
We have the $\GL(V)$-actions on $\nbiga$ and $Q$.
We also have the equivariant morphism
$\Psi_m:Q\lrarr\nbiga$.
The $\vecdelta_{\kappa}$-semistability condition
determines the open subset $Q^{ss}_{\kappa}$.
The following lemma can be shown by the same argument
as the proof of Lemma \ref{lem;06.5.15.55}.
\begin{lem}
\label{lem;06.5.17.155}
Assume that the absolute value of $\gamma\neq 0$
is sufficiently small.
\begin{itemize}
\item
Then, we have
$\Psi_m^{-1}\bigl(\nbiga^{ss}(\nbigl_{\gamma})\bigr)
=Q^{ss}_{\sign(\gamma)}$.
\item
The induced morphism
$\Psi_m:Q_{\sign(\gamma)}^{ss}\lrarr
 \nbiga^{ss}(\nbigl_{\gamma})$
is a  closed immersion.
Moreover, the image
is contained in
$\nbiga^{s}(\nbigl_{\gamma})$.
\hfill\qed
\end{itemize}
\end{lem}

We take a positive rational number $\gamma_1$
and a negative rational number $\gamma_2$
such that $|\gamma_i|$ are sufficiently small.
We take large number $k'$
such that $k'(\gamma_1-\gamma_2)=1$.
We put $\nbigl_i:=\nbigl_{\gamma_i}^{\otimes k'}$.
We have $\nbigl_2=\nbigl_1\otimes\nbigo_{\proj^{(1)}_m}(-1)$.

Let us consider $\nbigb:=\proj\bigl(\nbigo_{\proj^{(1)}_m}(0)
 \oplus\nbigo_{\proj^{(1)}_m}(1)\bigr)$ over $\nbiga$.
We put $\nbigb_1=\proj\bigl(\nbigo_{\proj^{(1)}_m}(0)\bigr)$
and $\nbigb_2=\proj\bigl(\nbigo_{\proj^{(1)}_m}(1)\bigr)$,
which are naturally regarded as the closed subscheme of $\nbigb$.
We have the tautological line bundle
$\nbigo_{\proj}(1)$,
and we put $\nbigo_{\nbigb}(1):=\nbigo_{\proj}(1)\otimes\nbigl_1$.
We have the natural $\SL(V)$-action on $\nbigb$,
and $\nbigo_{\nbigb}(1)$ gives the equivariant polarization.
Let $\nbigb^{ss}$ denote the set of the semistable points
of $\nbigb$ with respect to $\nbigo_{\nbigb}(1)$.
We put $\TH^{ss}:=Q\times_{\nbiga}\nbigb^{ss}$.
We have the natural $\SL(V)$-action on $\TH^{ss}$.
As in the case of Proposition \ref{prop;06.5.17.50},
the following proposition can be shown easily.
\begin{prop}
\label{prop;06.5.17.150}
The quotient stack $\TH^{ss}/\SL(V)$
is Deligne-Mumford.
\hfill\qed
\end{prop}

We put $\THhat^{ss}:=\TH^{ss}\times_{Q}\Qhat$,
where we put $\Qhat:=Q(\vecyhat,[\vecL],\alpha_{\ast},\vecdelta)$.
We put $\Mhat:=\THhat^{ss}/\GL(V)$,
which is called the master space.
\index{master space}
\index{$\Mhat$}
We have the $G_m$-action $\rhobar$
as usual.
From Proposition \ref{prop;06.5.17.150},
we obtain the following.
\begin{prop}
$\Mhat$ is Deligne-Mumford and proper.
\hfill\qed
\end{prop}

We put $\THhat^{ss}_i:=\nbigb^{ss}_i\times_{\nbiga}\THhat^{ss}$,
and we put $\Mhat_i:=\THhat^{ss}_i/\GL(V)$.
\index{$\Mhat_i$}
We put $\Mhat^{\ast}:=\Mhat-\Mhat_1\cup\Mhat_2$.
\index{$\Mhat^{\ast}$}
It is easy to observe that
$\Mhat^{\ast}$ is an open substack
of the moduli stack
$\nbigm(m,\vecyhat,L_1,[L_2])$
of the tuple
$(E_{\ast},\rho,\phi_1,[\phi_2])$ as follow:
\begin{itemize}
\item $E_{\ast}$ is a parabolic sheaf
 with an orientation $\rho$
 satisfying the condition $O_m$.
\item $\phi_1$ is an $L_1$-section of $E_{\ast}$
 such that $\phi_1\neq 0$.
\item
 $[\phi_2]$ is a reduced $L_2$-section
 of $E_{\ast}$ such that $[\phi_2]\neq 0$.
\end{itemize}

From Lemma \ref{lem;06.5.17.155},
$\Mhat_1$ and $\Mhat_2$ are isomorphic to
$\nbigm^{ss}(\vecyhat,[\vecL],\alpha_{\ast},\vecdelta_+)$
and $\nbigm^{ss}(\vecyhat,[\vecL],\alpha_{\ast},\vecdelta_-)$
respectively.
They give the obvious fixed point sets of
$\Mhat$ with respect to $\rhobar$.

Let us see the other components of the fixed point set.
A decomposition type is defined to be 
a datum $\gbigi:=(\vecy_1,\vecy_2)\in\Type^2$
satisfying the following:
\[
\vecy_1+\vecy_2=\vecy,
\quad
P^{\alpha_{\ast}}_{\vecy}=
 P^{\alpha_{\ast}}_{\vecy_1}=P^{\alpha_{\ast}}_{\vecy_2},
\quad
r_i=\rank(\vecy_i)
\]
For a decomposition type $\gbigi:=(\vecy_1,\vecy_2)$,
let $\Mhat^{G_m}(\gbigi)$ be the moduli stack 
of the objects
$\bigl(E_{\ast}^{(1)},\phi_1,E_{\ast}^{(2)},[\phi_2],\rho\bigr)$
as follows:
\begin{itemize}
\item
 $(E_{\ast}^{(1)},\phi_1)$ is $\delta_1$-stable
 $L_1$-Bradlow pair of type $\vecy_1$.
\item
 $(E_{\ast}^{(2)},[\phi_2])$ is $\delta_2$-stable
 reduced $L_2$-Bradlow pair  of type $\vecy_2$.
\item
 $\rho$ is an orientation of $E^{(1)}\oplus E^{(2)}$.
\end{itemize}

Note that $\delta_i$ are sufficiently small,
and hence the $1$-stability conditions
for $(\vecy_1,L_1,\alpha_{\ast},\delta_1)$
and $(\vecy_2,L_2,\alpha_{\ast},\delta_2)$
are satisfied.
We have the naturally defined morphism
$\Mhat^{G_m}(\gbigi)\lrarr \Mhat$,
as in the subsubsection \ref{subsubsection;06.5.17.1}.

Let $S(\vecy,\alpha_{\ast},\vecdelta)$
denote the set of decomposition type.
Then, we can show the following
by the same argument as the proof of 
Proposition \ref{prop;06.5.16.150}.
\begin{prop}
$\Mhat_1\sqcup \Mhat_2\sqcup
 \coprod_{\gbigi \in S(\vecy,\alpha_{\ast},\vecdelta)}
 M^{G_m}(\gbigi)$
is the stack theoretic fixed point set
of $\Mhat$ with respect to $\rhobar$.
\hfill\qed
\end{prop}

We naturally have the oriented reduced $\vecL$-Bradlow pair
$(\Ehat^{\Mhat}_{\ast},[\phi^{\Mhat}_1],[\phi^{\Mhat}_2],\rho)$
on $\Mhat\times X$, as in the subsubsection
\ref{subsubsection;06.5.17.10}.
The restriction of 
$(\Ehat^{\Mhat}_{\ast},[\phi^{\Mhat}_1],[\phi^{\Mhat}_2],\rho)$ to
$\Mhat_i\times X$ has the universal property
with respect to the moduli theoretic meaning of $\Mhat_i$.

Let $\nbigo_{\rel}^{(i)}(1)$ denote the line bundle
on $\nbigm(m,\vecyhat,[\vecL])$
which is the pull back
of the relatively tautological line bundle
on $\nbigm(m,\vecyhat,[L_i])$
via the natural morphism
$\nbigm(m,\vecyhat,[\vecL])\lrarr \nbigm(m,\vecyhat,[L_i])$.
Let $\varphi:\Mhat\lrarr \nbigm(m,\vecyhat,[\vecL])$
denote the naturally defined morphism.
The restriction of
$\varphi^{\ast}\nbigo^{(1)}_{\rel}(1)$ to $\Mhat^{\ast}$
is canonically trivialized.
Thus, the restriction
$[\phi_1^{\Mhat}]_{|\Mhat^{\ast}}$
induces the $L_1$-section of $\Ehat^{\Mhat}$,
which we denote by $\phi_1^{\Mhat}$.
We put $\nbigi^{(2)}:=\varphi^{\ast}\nbigo^{(2)}_{\rel}(-1)$.
Then, $[\phi^{\Mhat}_2]$
gives the morphism
$\nbigi^{(2)}\otimes L_2\lrarr \Ehat^{\Mhat}$.
The restriction of
$(\Ehat^{\Mhat},\phi^{\Mhat}_1,[\phi_2^{\Mhat}],\rho)$
to $\Mhat^{G_m}(\gbigi)$ has the universal property
with respect to the moduli theoretic meaning above.
Correspondingly, we have the decomposition
$\Ehat^{\Mhat}_{|\Mhat^{G_m}(\gbigi)}
=E_1^{\Mhat}\oplus E^{\Mhat}_2$.

\vspace{.1in}

It is convenient to decompose
$\Mhat^{G_m}(\gbigi)$ into the product
of two moduli stacks up to etale finite morphisms.
By the same argument as the proof of
Proposition \ref{prop;06.5.17.30},
we obtain the following description
of $\Mhat^{G_m}(\gbigi)$
up to etale finite morphisms.
\begin{prop}
We put $r_i:=\rank\vecy_i$.
There exists the algebraic stack $\nbigs$
with the following properties:
\mbox{{}}
\begin{itemize}
\item
There exist the following diagram:
\begin{equation}
\begin{CD}
 \Mhat^{G_m}(\gbigi)@<{F}<<
 \nbigs
 @>{G}>>
 \nbigm^{ss}\bigl(\vecyhat_1,[L_1],\alpha_{\ast},\delta_1\bigr)
 \times
 \nbigm^{ss}\bigl(\vecyhat_2,[L_2],\alpha_{\ast},\delta_2\bigr)
\end{CD}
\end{equation}
The morphisms $F$ and $G$ are etale and finite of degree 
$(r_1\cdot r_2)^{-1}$ and $r_2^{-1}$, respectively.
We also have the following diagram:
\begin{equation}
\label{eq;06.6.22.10}
\begin{CD}
 \Mhat^{G_m}(\gbigi)@<{F}<<
 \nbigs
 @>{G'}>>
 \nbigm^{ss}\bigl(\vecy_1,L_1,\alpha_{\ast},\delta_1\bigr)
 \times
 \nbigm^{ss}\bigl(\vecyhat_2,[L_2],\alpha_{\ast},\delta_2\bigr)
\end{CD}
\end{equation}
Here $G'$ is etale and finite of degree $(r_1r_2)^{-1}$.
\item
Let $\nbigo_{i,\rel}(1)$ denote the pull back of
the relative tautological line bundle
of $\nbigm^{ss}(\vecy_i,[L_i],\alpha_{\ast},\delta)$.
There exists the line bundle
$\nbigo_{1,\rel}(1/r_2)$ on $\nbigs$
such that  $\nbigo_{1,\rel}(1/r_2)^{r_2}\simeq G^{\ast}\nbigo_{1,\rel}(1)$,
and we have the following relation:
\begin{equation}
 F^{\ast}E_1^{\Mhat}
\simeq G^{\ast}\Ehat_1^{u}\otimes\nbigo_{1,\rel}(1),
\quad
 F^{\ast}E_2^{\Mhat}
\simeq G^{\ast}\Ehat_2^u\otimes\nbigo_{1,\rel}(-r_1/r_2),
\quad
 F^{\ast}\nbigi^{(2)}\simeq
 G^{\ast}\nbigo_{2,\rel}(-1)\otimes
 \nbigo_{1,\rel}(-r_1/r_2).
\end{equation}
Here $\Ehat_i^u$
denotes the pull back of the universal sheaf
over $\nbigm^{ss}(\vecyhat_i,[L_i],\alpha_{\ast},\delta_i)\times X$.
We also have the following:
\[
 F^{\ast}E_1^{\Mhat}
\simeq
 G^{\prime\,\ast}E_1^u,
\quad
 F^{\ast}E_2^{\Mhat}
\simeq
 G^{\prime\,\ast}\Ehat_2^u\otimes\Or(E_1^u)^{-1/r_2}
\]
Here, $E_1^u$ denotes the pull back of the universal sheaf
over $\nbigm^{ss}(\vecy_1,L_1,\alpha_{\ast},\delta_1)\times X$,
and $\Or(E_1^u)^{-1/r_2}$ denotes $\nbigo_{1,\rel}(-r_1/r_2)$.
\item
The weights of $\rhobar$ on $E_1^{\Mhat}$,
$E_2^{\Mhat}$
and $\nbigi^{(2)}_{|\Mhat^{G_m}(\gbigi)}$
are $-1$, $r_1/r_2$ and $r_1/r_2$, respectively.
\hfill\qed
\end{itemize}
\end{prop}

\section{Obstruction Theory of the Moduli Stacks and the Master spaces}
\label{section;06.6.4.250}
In this section,
we discuss the obstruction theory
of the moduli stack of parabolic reduced Bradlow pairs
on a smooth projective {\em surface} $X$.
We also assume the smoothness of the divisor $D$,
which is the support of the parabolic structure.
The naive strategy is explained in 
the subsubsection \ref{subsubsection;06.6.22.1}.
We will also discuss the obstruction theory
for the master spaces.

\vspace{.1in}

\noindent
{\bf Notation}
Let $S$ be a scheme, $Z$ be a stack over $S$,
and $G$ be a smooth group scheme over $S$.
When $G$ acts on $Z$,
the quotient stack is denoted by $Z_G$.

\subsection{Deformation of Torsion-free Sheaves}
\label{subsection;06.7.3.31}
\subsubsection{Construction of the basic complex}
\label{subsubsection;06.5.2.3}

Let $U$ be any algebraic stack over $k$,
and let $E$ be a torsion-free $U$-coherent sheaf
defined over $U\times X$.
Assume that we have a locally free resolution
$V_{\cdot}=(V_{-1}\rarr V_0)$ of $E$ on $U$,
i.e.,
$V_i$ are locally free sheaves of finite ranks,
and we have the surjection $V_{0}\lrarr E$
whose kernel is $V_{-1}$.
The inclusion $V_{-1}\subset V_0$
is denoted by $f$.
We put as follows (see the subsubsection \ref{subsubsection;06.4.29.30}):
\[
 \gminig(V_{\cdot}):=\nhom(V_{\cdot},V_{\cdot})^{\lor}[-1]
\]
\index{$\gminig(V_{\cdot})$}

Let $W_{-1}$ and $W_0$ be vector spaces over $k$
such that $\rank W_{i}=\rank V_i$.
We denote $W_{i}\otimes\nbigo_X$ by $W_{i\,X}$.
We put $\GL(W_{\cdot}):=\GL(W_{-1})\times\GL(W_0)$.
We have the natural right $\GL(W_{\cdot})$-action
on the vector bundle $N(W_{-1\,X},W_{0\,X})$
given by
$(g_{-1},g_{0})\cdot f=g_{0}^{-1}\circ f\circ g_{-1}$.
Here $g_i$ denotes an element of $GL(W_i)$,
and $f$ denotes an element of $N(W_{-1\,X},W_{0\,X})$.
The quotient stack is denoted by $\yw$.
\index{$Y(W_{\cdot})$}

Then, we have the classifying map
$\Phi(V_{\cdot}):U\times X\lrarr \yw$
over $X$,
and thus 
$\Phi(V_{\cdot})^{\ast}L_{\yw/X}\lrarr L_{U\times X/X}$.
As explained in Example \ref{example;06.4.30.10},
$\Phi(V_{\cdot})^{\ast}L_{\yw/X}$ is represented by
the complex $\gminig(V_{\cdot})_{\leq 1}$.
We have the naturally defined morphism
$\gminig(V_{\cdot})\lrarr \gminig(V_{\cdot})_{\leq 1}$,
and hence
$\gminig(V_{\cdot})\lrarr L_{U\times X/X}$.

Let $\omega_X$ denote the dualizing complex of $X$,
and we put as follows:
\[
 \Ob(V_{\cdot}):=Rp_{X\,\ast}\bigl(
 \gminig(V_{\cdot})\otimes\omega_X \bigr)
\]\index{$\Ob(V_{\cdot})$, $\ob(V_{\cdot})$}
Then, we have the naturally defined morphism
$\ob(V_{\cdot}):\Ob(V_{\cdot})\lrarr L_U$.

\begin{lem}
\label{lem;6.14.6}
The object $\Ob(V_{\cdot})$
and the morphism $\ob(V_{\cdot})$ 
depends only on $E$
in the derived category $D(U)$,
in the sense that it is independent of
the choice of a resolution $V_{\cdot}$.
\end{lem}
\pf
Let $V^{(1)}_{\cdot}$ be another resolution.
We would like to compare the two morphisms
$\ob(V^{(1)}_{\cdot})$ and $\ob(V_{\cdot})$
in the derived category $D(U)$.
We put $V^{(2)}_0:=V_0\oplus V^{(1)}_0$
and $V^{(2)}_{-1}:=\ker(V^{(2)}_0\lrarr E)$.
We have the morphism $V_{\cdot}\lrarr V^{(2)}_{\cdot}$
and $V_{\cdot}\lrarr V^{(1)}_{\cdot}$.
Therefore,
we have only to compare
the morphisms $\ob(V^{(2)}_{\cdot})$
and $\ob(V_{\cdot})$.
Note that $V_{i}$ is a subbundle of $V^{(2)}_i$,
i.e.,
we have the filtration
$V_{\cdot}\subset V^{(2)}_{\cdot}$.
Let $\nhom'\bigl(V^{(2)}_i,V^{(2)}_j\bigr)$
 $(i,j=0,-1)$
denote the sheaf of $\nbigo_X$-morphisms
$V^{(2)}_i\lrarr V^{(2)}_j$
preserving the filtrations.
They naturally form the complex of sheaves
$\nhom'(V^{(2)}_{\cdot},V^{(2)}_{\cdot})$.
We put 
$\gminig(V_{\cdot},V^{(2)}_{\cdot}):=
 \nhom'(V^{(2)}_{\cdot},V^{(2)}_{\cdot})^{\lor}[-1]$.
We have the naturally defined quasi isomorphisms
$\nhom'(V^{(2)}_{\cdot},V^{(2)}_{\cdot})
 \lrarr \nhom(V_{\cdot},V_{\cdot})$
and $\nhom'(V^{(2)}_{\cdot},V^{(2)}_{\cdot})
 \lrarr \nhom(V^{(2)}_{\cdot},V^{(2)}_{\cdot})$.
They induce the quasi isomorphisms
$\gamma_1:\gminig(V_{\cdot})\lrarr
 \gminig(V_{\cdot},V^{(2)}_{\cdot})$
and $\gamma_2:\gminig(V^{(2)})\lrarr \gminig(V_{\cdot},V^{(2)}_{\cdot})$.

Let $W_i^{(2)}$ be vector spaces
such that $\rank W_i^{(2)}=\rank V_i^{(2)}$ $(i=1,2)$.
We fix inclusions $W_i\subset W_i^{(2)}$.
We denote $W_i^{(2)}\otimes\nbigo_X$ by $W^{(2)}_{i\,X}$.
We have the filtration
$W_{i\,X}\subset W^{(2)}_{i\,X}$.
Let $\nhom'(W_{i\,X}^{(2)},W_{j\,X}^{(2)})$ be
the sheaf of $\nbigo_X$-morphisms
$W^{(2)}_{i\,X}\lrarr W^{(2)}_{j\,X}$
preserving the filtration.
The corresponding vector bundle is denoted
by $N'(\wx{i}^{(2)},\wx{j}^{(2)})$.
We have the natural morphisms
$N'\bigl(\wx{i}^{(2)},\wx{j}^{(2)}\bigr)\lrarr N(\wx{i},\wx{j})$
and
$N'\bigl(\wx{i}^{(2)},\wx{j}^{(2)}\bigr)
 \lrarr N\bigl(\wx{i}^{(2)},\wx{j}^{(2)}\bigr)$.
Let $\GL'(W^{(2)}_i)$ be the subgroup
of $\GL(W^{(2)}_i)$,
which consists of the elements of $\GL(W^{(2)}_i)$
preserving the filtration.
Then we have the natural right
$\GL'(W^{(2)}_1)\times \GL'(W^{(2)}_2)$-action
on $N'(W_{j}^{(2)},W_{l}^{(2)})$.
The quotient stack is denoted by $Y'(W^{(2)}_{\cdot})$.
We have the homomorphisms
$\GL'(W^{(2)}_i)\lrarr \GL(W^{(2)}_i)$
and $\GL'(W^{(2)}_i)\lrarr \GL(W_i)$.
Thus we have the morphisms
$Y'(W^{(2)})\lrarr Y(W_{\cdot})$
and $Y'(W^{(2)})\lrarr Y(W^{(2)}_{\cdot})$.

From the tuple $\bigl(V_{\cdot},V^{(2)}_{\cdot}\bigr)$,
we obtain the morphism
$\Phi\bigl(V_{\cdot},V^{(2)}_{\cdot}\bigr):
 U\times X\lrarr Y'(W^{(2)}_{\cdot})$.
Then, it can be shown 
by the argument in the subsubsection  
\ref{subsubsection;06.4.29.15},
that the complex $\Phi\bigl(V_{\cdot},V^{(2)}_{\cdot}\bigr)$
is represented by $\gminig(V_{\cdot},V^{(2)}_{\cdot})$.
We also have the following commutative diagram:
\begin{equation} 
\begin{CD}
 \Phi(V_{\cdot})^{\ast}
 L_{\yw/X} @>>>
 \Phi(V_{\cdot},V^{(2)}_{\cdot})^{\ast}
 L_{Y'(W^{(2)}_{\cdot})/X }
 @<<<
 \Phi(V^{(2)}_{\cdot})^{\ast}
 L_{Y(W^{(2)}_{\cdot})/X }\\
 @AAA @AAA @AAA \\
 \gminig(V_{\cdot})_{\leq 1}
 @>{\gamma_1}>>
 \gminig(V_{\cdot},V^{(2)}_{\cdot})_{\leq 1}
 @<{\gamma_2}<<
 \gminig(V^{(2)}_{\cdot})_{\leq 1}
\end{CD}
\end{equation}
Then, we obtain the following diagram:
\[
\begin{CD}
 L_{U\times X/X}
 @>{=}>>
 L_{U\times X/X}
 @<{=}<< L_{U\times X/X}\\
 @AAA @AAA @AAA \\
 \Phi(V_{\cdot})^{\ast}
 L_{\yw/X} @>>>
 \Phi(V_{\cdot},V^{(2)}_{\cdot})^{\ast}
 L_{Y'(W^{(2)}_{\cdot})/X }
 @<<<
 \Phi(V^{(2)}_{\cdot})^{\ast}
 L_{Y(W^{(2)}_{\cdot})/X }\\
 @AAA @AAA @AAA \\
 \gminig(V_{\cdot}) @>{\simeq}>>
 \gminig(V_{\cdot},V_{\cdot}^{(2)})
 @<{\simeq}<<
 \gminig(V_{\cdot}^{(2)}).
\end{CD}
\]
We put $\Ob\bigl(V^{(2)}_{\cdot},V_{\cdot}\bigr):=
 Rp_{X\,\ast}\bigl(
 \gminig(V^{(2)}_{\cdot},V_{\cdot})\otimes\omega_X\bigr)$.
Then we obtain the following diagram in $D(U)$:
\[
 \begin{CD}
 L_{U} @>{=}>> L_U @<{=}<< L_U \\
 @AAA @AAA @AAA \\
 \Ob(V_{\cdot})@>{\simeq}>>
 \Ob\bigl(V^{(2)}_{\cdot},V_{\cdot}\bigr)
 @<{\simeq}<<
 \Ob(V^{(2)}_{\cdot})
 \end{CD}
\]
Thus we are done.
\hfill\qed

\vspace{.1in}
If $E$ is a vector bundle of rank $R$,
we may take $V_0=E$ and $V_{-1}=0$.
In this case, 
the construction can be reworded
as follows:
We have the classifying map
$\Phi(E):U\times X\lrarr X_{\GL(R)}$.
It induces the morphism
$\Phi(E)^{\ast}L_{X_{\GL(R)/X}}\lrarr 
 L_{U\times X/X}$.
We have the expression
$\Phi(E)^{\ast}L_{X_{\GL(R)/X}}
\simeq \nhom(E,E)$.
We put
$\Ob(E):=
 Rp_{X\,\ast}\bigl(\nhom(E,E)\otimes\omega_X\bigr)$,
and we obtain the morphism
$\ob(E):\Ob(E)\lrarr L_{U/k}$.

\subsubsection{The trace-free part and the diagonal part}
\label{subsubsection;06.5.9.20}

We have the homomorphism
$\GL(W_{\cdot})\lrarr G_m$
given by $(f_{-1},f_0)\longmapsto
 \det(f_{-1})^{-1}\cdot \det (f_0)$.
It induces the morphism
$\gminiw:Y(W_{\cdot})\lrarr X_{G_m}$.
The composite
$\gminiw\circ\Phi(V_{\cdot}):U\times X\lrarr X_{G_m}$
is same as the classifying map
of the determinant bundle
$\det(E)\simeq\det(V_0)\otimes\det(V_{-1})^{-1}$,
which is denoted by $\Phi(\det(E))$.
Then, we obtain the following commutative diagram:
\[
 \begin{CD}
 \Phi\bigl(\det(E)\bigr)^{\ast}L_{X_{G_m}/X}
@>>>
 \Phi(V_{\cdot})^{\ast}L_{Y(W_{\cdot})/X}
@>>>
 L_{U\times X/X}\\
 @A{\simeq}AA @AAA \\
 \nbigo[-1]
@>{i}>>
 \gminig(V_{\cdot})
 \end{CD}
\]
Here, the map $i:\nbigo[-1]\lrarr\gminig(V_{\cdot})$
is given as follows:
\[
 \nbigo\lrarr\nhom(V_0,V_0)\oplus\nhom(V_{-1},V_{-1}),
\quad
 f\longmapsto (f\cdot \id_{V_0},-f\cdot \id_{V_{-1}})
\]
On the other hand,
we have the trace map 
$\tr:\gminig(V_{\cdot})\lrarr\nbigo[-1]$:
\[
 \nhom(V_0,V_0)\oplus\nhom(V_{-1},V_{-1})
\lrarr
 \nbigo,
\quad
 (f_0,f_{-1})\longmapsto
 \tr(f_0)+\tr(f_{-1})
\]
We put $\Ker(\tr):=\gminig^{\circ}(V_{\cdot})$,
and $\gminig^d(V_{\cdot}):=\Image(i)$.
\index{$\gminig^{\circ}(V_{\cdot})$, $\gminig^d(V_{\cdot})$}
We have the decomposition
$\gminig(V_{\cdot})=\gminig^{\circ}(V_{\cdot})
\oplus \gminig^{d}(V_{\cdot})$,
which induces the decomposition
$\Ob(V_{\cdot})=\Ob^{\circ}(V_{\cdot})\oplus 
 \Ob^{d}(V_{\cdot})$.
\index{$\Ob^{\circ}(V_{\cdot})$, $\Ob^d(V_{\cdot})$}
The complexes
$\gminig^{\circ}(V_{\cdot})$ and
$\Ob^{\circ}(V_{\cdot})$
(resp. $\gminig^d(V_{\cdot})$
and $\Ob^d(V_{\cdot})$)
are called the trace-free part
(resp. the diagonal part).

The determinant bundle induces the morphism
$\dettilde_E:U\lrarr\nbigm(1)$.
We also have the following commutative diagram:
\[
 \begin{CD}
 U\times X @>{\Phi(V_{\cdot})}>>
 Y(W_{\cdot}) \\
 @V{\dettilde_{E,X}}VV @V{\gminiw}VV\\
 \nbigm(1)\times X@>>> X_{G_m}
 \end{CD}
\]
The composite $\gminiw\circ\Phi(V_{\cdot})$
is same as the classifying map $\Phi(\det(E))$
of $\det(E)$.
Thus we obtain the following commutative diagram:
\[
 \begin{CD}
 L_{U\times X/X} @<<<
 \Phi(V_{\cdot})^{\ast}L_{Y(W_{\cdot})/X} 
 @<<< \gminig(V_{\cdot})\\
 @AAA @AAA @AAA \\
 \dettilde_{E,X}^{\ast}L_{\nbigm(1)\times X/X}
 @<<< \Phi(\det(E))^{\ast}L_{X_{G_m}/X}
 @=
 \gminig^d(V_{\cdot})
 \end{CD}
\]
Therefore, we obtain the following commutative diagram:
\begin{equation}
 \label{eq;06.5.13.1}
 \begin{CD}
 L_{U} @<<< \Ob(V_{\cdot})\\
 @AAA @AAA \\
 \dettilde_{E}^{\ast}L_{\nbigm(1)} @<<<
 \Ob^d(V_{\cdot})
 \end{CD}
\end{equation}

\subsubsection{Preparation for  Master space}
\label{subsubsection;06.5.4.100}

We also put $A(W_{\cdot}):=X_{\GL(W_0)}$.
We have the natural morphism
$\Gamma:\yw\lrarr A(W_{\cdot})$.
We put $\Psi(V_{\cdot}):=\Gamma\circ\Phi(V_{\cdot})$.
Then $\Psi(V_{\cdot})^{\ast}L_{A(W_{\cdot})/X}$
is represented by 
$\nhom(V_0,V_0)^{\lor}[-1]$.
For these representatives,
the natural morphism
$\Psi(V_{\cdot})^{\ast}L_{A(W_{\cdot})/X}
\lrarr
 \Phi(V_{\cdot})^{\ast}L_{Y(W_{\cdot})/X}$
is expressed by the obvious inclusion
of $\nhom(V_0,V_0)$ 
to $\nhom(V_0,V_0)\oplus \nhom(V_{-1},V_{-1})$.
We put
$\gminih(V_{\cdot}):=
 \nhom(V_{0},V_{\cdot})^{\lor}[-1]$
and 
$\Ob^G(V_{\cdot}):=
 Rp_{X\,\ast}\bigl(
 \gminih(V_{\cdot})\otimes\omega_X
 \bigr)$.
\index{$\gminih(V_{\cdot})$}
\index{$\Ob^G(V_{\cdot})$}
Then we obtain the following diagram:
\[
 \begin{CD}
 \gminih(V_{\cdot}) @>>>
 \gminig(V_{\cdot}) @.\\
@VVV @VVV \\
 \Psi(V_{\cdot})^{\ast}L_{A(W_{\cdot})/X}
 @>>>
 \Phi(V_{\cdot})^{\ast}L_{Y(W_{\cdot})/X}
 @>>>
 L_{U\times X/X}
 \end{CD}
\]
Therefore, we obtain the morphisms
$ \Ob^G(V_{\cdot})\lrarr
 \Ob(V_{\cdot})\lrarr
 L_{U/k}$.

\vspace{.1in}
Now, we assume the following condition (C)
for $E$ and $V_{\cdot}$:
\begin{description}
\item[(C1)]
 For any point $u\in U$,
 the higher cohomology groups
 $H^i\bigl(X,E_{|\{u\}\times X}\bigr)$
 vanish.
\item[(C2)]
 We put $V':=p_{X\,\ast}(E)$.
 Then we have $V_0=p_X^{\ast}V'$.
\end{description}

We put $B(W_{\cdot}):=\Spec(k)_{\GL(W_0)}$.
\index{$B(W_{\cdot})$}
We remark $A(W_{\cdot})=X\times B(W_{\cdot})$.
Because of $V_0=p_X^{\ast}V'$,
we have $\rank W_0=\rank V'$,
and hence we have the classifying map
$\Psi(V'):U\lrarr B(W_{\cdot})$ of $V'$.
We have the following commutative diagram:
\begin{equation}
\begin{CD}
U\times X 
 @>{\Phi(V_{\cdot})}>> Y(W_{\cdot}) @>>> X\\
@V{\Psi(V')_X}VV @V{\Gamma}VV @VVV\\
B(W_{\cdot})\times X @>{=}>> A(W_{\cdot}) @>>> X
\end{CD}
\end{equation}
Thus we obtain the following diagram on $U\times X$:
\begin{equation}
\begin{CD}
L_{U\times X/X} @<<<
 \Phi(V_{\cdot})^{\ast}L_{Y(W_{\cdot})/X} @<<< 
 \gminig(V_{\cdot})\\
 @AAA @AAA @AAA \\
\Psi(V')_X^{\ast}L_{B(W_{\cdot})\times X/X} 
 @<{\simeq}<<
 \Psi(V_{\cdot})^{\ast}L_{A(W_{\cdot})/X}
 @<<<
 \gminih(V_{\cdot})
\end{CD}
\end{equation}
Hence, we obtain the following diagram on $U$:
\begin{equation}
\label{eq;06.5.9.1}
\begin{CD}
L_U @<<< \Ob(V_{\cdot}) \\
@AAA @AAA \\
\Psi(V')^{\ast}L_{B(W_{\cdot})}
 @<{\tau_1}<< 
 \Ob^G(V_{\cdot}).
\end{CD}
\end{equation}

It is easy to see that
both of $\Ob^G(V_{\cdot})$
and $\Psi(V')^{\ast}L_{B(W_{\cdot})}$
are isomorphic to
$\nhom(V',V')[-1]$
under the condition (C).

\begin{lem}
\label{lem;06.5.4.201}
The morphism
$\tau_1$ in {\rm (\ref{eq;06.5.9.1})} is isomorphic.
\end{lem}
\pf
The composite of the following naturally defined morphisms
is isomorphic:
\[
  Rp_{X\,\ast}\Bigl(
 \bigl(
 \nhom(V_0,V_0)\rarr \nhom(V_{-1},V_0)
\bigr)\otimes\omega_X
 \Bigr)
\lrarr
 Rp_{X\,\ast}\bigl(
 \nhom(V_0,V_0)\otimes\omega_X
 \bigr)
\lrarr \nhom(V',V')
\]
Then, the claim of the lemma immediately follows.
\hfill\qed

\subsubsection{Basic complex on the moduli stack $\nbigm(m,y)$}
\label{subsubsection;06.7.3.311}

Let $y\in H^{\ast}(X)$ be a Chern character
of a coherent sheaf on $X$.
Let $H$ denote the polynomial associated to $y$.
We take an $H(m)$-dimensional vector space $V_m$.
We have the scheme $Q^{\circ}(m,y)$.
(See the subsubsection \ref{subsubsection;06.5.9.50}.)
We consider the universal quotient
$q^u:p_{Q^{\circ}(m,y)}^{\ast}V_{m,X}\lrarr E^u(m)$ 
defined over $Q^{\circ}(m,y)\times X$.
We put $V^u_0:=p_{Q^{\circ}(m,y)}^{\ast}V_{m,X}$ and
$V^u_{-1}:=\ker\bigl(V^u_0\lrarr E^u(m)\bigr)$.
The inclusion
$V^u_{-1}\lrarr V_{0}^u$ is denoted by $f^u$.
We put
$V':=V_m\otimes\nbigo_{Q^{\circ}(m,y)}=p_{X\,\ast}V^u_0$.
We have the morphism
$\pi:Q^{\circ}(m,y)\lrarr\nbigm(m,y):=Q^{\circ}(m,y)/GL(V_m)$.
The latter is an open subset
of the moduli stack of torsion-free sheaves
of type $y$,
determined by the condition $O_m$.
The descent of $E^u$,
$V_{\cdot}$ and $V'$
with respect to the $\GL(V_m)$-action
are denoted by $\nbige^u$,
$\nbigv_{\cdot}$ and $\nbigv'$.
The sheaf $\nbige^u$ is the universal sheaf.

We put $W_0:=V_m$,
and we take a vector space $W_{-1}$
such that $\dim W_{-1}=H(m)-\rank(y)$.
Applying the result in the subsubsection
\ref{subsubsection;06.5.2.3},
we obtain the complex $\Ob(\nbigv_{\cdot})$
and the morphism
$\ob(\nbigv_{\cdot}):
 \Ob(\nbigv_{\cdot})\lrarr L_{\nbigm(m,y)/k}$.
We obviously have
$\pi^{\ast}\Ob(\nbigv_{\cdot})
=\Ob(V^u_{\cdot})$.

\begin{lem}
 \label{lem;06.5.2.2}
We have the following morphism of the distinguished triangles
on $Q^{\circ}(m,y)$:
\[
 \begin{CD}
 \pi^{\ast}\Ob(\nbigv_{\cdot})
 @>>> 
 \Ob(V^u_{-1},f^u) @>>>
 \nhom(V',V') @>>>
 \pi^{\ast}\Ob(\nbigv_{\cdot})[1]\\
 @VVV @VVV @V{\simeq}VV @VVV \\
 \pi^{\ast}L_{\nbigm(m,y)}@>>>
 L_{Q^{\circ}(m,y)}@>>>
 L_{Q^{\circ}(m,y)/\nbigm(m,y)}@>>>
 \pi^{\ast}L_{\nbigm(m,y)}[1].
 \end{CD}
\]
\end{lem}
\pf
We use the notation in the subsubsection
\ref{subsubsection;06.5.4.5}.
We have the following commutative diagram:
\[
 \begin{CD}
 Q^{\circ}(m,y) @>{\Phi(V^u_{-1},f^u)}>> \yw_{quo} @>>> X \\
 @VVV @V{\pi_0}VV @VVV \\
 \nbigm(m,y) @>>> \yw @>>> X
 \end{CD}
\]
The composite $\pi_0\circ \Phi(V^u_{-1},f^u)$
is same as $\Phi(V^u_{\cdot})$.
Thus we obtain the following morphism of distinguished triangles
on $Q^{\circ}(m,y)\times X$:
\[
 \begin{CD}
 \Phi(V^u_{\cdot})^{\ast}L_{ \yw/X}
 @>{a}>>
 \Phi(V^u_{-1},f^u)^{\ast}L_{\ywquo/X}
 @>>>
 \Phi(V^u_{-1},f^u)^{\ast}L_{\ywquo/\yw}\\
 @VVV @VVV @V{\simeq}VV \\
 \pi^{\ast}L_{\nbigm(m,y)\times X/X} @>>>
 L_{Q^{\circ}(m,y)\times X/X} @>>> L_{Q^{\circ}(m,y)\times X/\nbigm(m,y)\times X}.
 \end{CD}
\]
Recall that
$\Phi(V^u_{\cdot})^{\ast}L_{ \yw/X}$
and $\Phi(V^u_{-1},f^u)^{\ast}L_{\ywquo/X}$
are expressed by
$\gminig(V^u_{\cdot})_{\leq 1}$
and $\gminig(V^u_{-1},f^u)_{\leq 1}$,
respectively.
It is easy to see that the morphism $a$
is expressed by the naturally defined morphism 
$\gminig(V_{\cdot}^u)_{\leq 1}\lrarr
\gminig(V^u_{-1},f^u)_{\leq 1}$
We put $\gminik(V^u_{\cdot}):=
 \cone\bigl(\gminig(V_{\cdot}^u)\lrarr
 \gminig(V^u_{-1},f^u)\bigr)$.
Then, obtain the following morphism
of the distinguished triangles on $Q^{\circ}(m,y)\times X$:
\[
 \begin{CD}
 \gminig(V^u_{\cdot})
 @>>>
 \gminig(V^u_{-1},f^u)
 @>>>
 \gminik(V^u_{\cdot})
 @>>>
 \gminig(V^u_{\cdot})[1]\\
 @VVV @VVV @VVV @VVV \\
 \pi^{\ast}L_{\nbigm\times X/X} @>>>
 L_{Q^{\circ}(m,y)\times X/X} @>>> 
 L_{Q^{\circ}(m,y)\times X/\nbigm(m,y)\times X} @>>>
 \pi^{\ast}L_{\nbigm(m,y)\times X/X}[1].
 \end{CD}
\]
Hence we obtain the following morphism of distinguished triangles
on $Q^{\circ}(m,y)$:
\[
 \begin{CD}
 \Ob(V_{\cdot})@>>>
 \Ob(V_{-1},f) @>>>
 Rp_{X\ast}\bigl(\gminik(V_{\cdot})\otimes\omega_X\bigr)
 @>>>
 \Ob(V_{\cdot})[1]\\
 @VVV @VVV @V{\varphi}VV @VVV \\
 \pi^{\ast}L_{\nbigm}@>>>
 L_{Q^{\circ}(m,y)}@>>>
 L_{Q^{\circ}(m,y)/\nbigm(m,y)}@>>>
 \pi^{\ast}L_{\nbigm(m,y)}[1].
 \end{CD}
\]
Recall
$\gminig(V^u_{\cdot})=
 \nhom\bigl(V^u_{\cdot},\,V^u_{\cdot}\bigr)^{\lor}[-1]$ and 
$\gminig(V^u_{-1},f^u)=
 \nhom\bigl(V^u_{-1}[1],\,V^u_{\cdot}\bigr)^{\lor}[-1]$.
It is easy to observe that
$\gminik(V^u_{\cdot})$ is expressed
by the complex
$\nhom(V_0^u,V_0^u)\rarr \nhom(V_{-1}^u,V_{0}^u)$,
where the first term stands at the degree $0$.
Under the identification,
the morphism $\varphi$ is given
by the identity of $\nhom(V_0^u,V_0^u)$.
Then, it is easy to check that
$Rp_{X\ast}(\gminik(V_{\cdot}))$
and $L_{Q^{\circ}(m,y)/\nbigm(m,y)}$ 
are quasi isomorphic to
their $0$-th cohomology sheaves $\nhom(V',V')$,
and that the morphism $\varphi$
in the diagram is isomorphic,
as in Lemma \ref{lem;06.5.4.201}.
\hfill\qed

\begin{cor}
 \label{cor;06.5.10.12}
The morphism
$\ob(\nbigv_{\cdot})$ gives an obstruction
theory for $\nbigm(m,y)$.
\hfill\qed
\end{cor}

We also use the notation
$\Ob(m,y)$ and $\ob(m,y)$
to denote $\Ob(\nbigv_{\cdot})$
and $\ob(\nbigv_{\cdot})$.
\index{$\Ob(m,y)$, $\ob(m,y)$}

\subsubsection{The case of the moduli of line bundles}
\label{subsubsection;06.5.11.1}

Let $\Poin$ denote the Poincar\'e bundle 
on $\Pic\times X$.
Then we have the classifying map
$\Phi(\Poin):\Pic\times X\lrarr X_{G_m}$.
We put $\gminig(\Poin):=\Phi(\Poin)^{\ast}L_{X_{G_m}/X}$
and
$\Ob(\Poin):=
 Rp_{X\,\ast}\bigl(\gminig(\Poin)\otimes\omega_X\bigr)$.
Then we have the morphism
$\ob(\Poin):\Ob(\Poin)\lrarr L_{\Pic}$
on $\Pic$.

Since $\gminig(\Poin)\simeq \nbigo[-1]$,
we have an isomorphism
in the derived category $D(\Pic)$:
\[
  \Ob(\Poin)\simeq 
\Bigl(
 H^0(X,\nbigo)^{\lor}\otimes\nbigo_{\Pic}[-1]
\Bigr)
\oplus
\Bigl(
 H^1(X,\nbigo)^{\lor}\otimes\nbigo_{\Pic}[0]
\Bigr)
\oplus
\Bigl(
 H^2(X,\nbigo)^{\lor}\otimes\nbigo_{\Pic}[1]
\Bigr)
\]
We have to be careful
that the decomposition is not canonical,
but the morphism
$\Ob(\Poin)\lrarr L_{\Pic}\lrarr
 H^1(X,\nbigo)^{\lor}\otimes\nbigo_{\Pic}[0]
\simeq\nbigh^1\bigl(\Ob(\Poin)\bigr)$
induces the canonical decomposition:
\[
 \tau_{\geq 0}\Ob(\Poin)
\simeq
\Bigl(
 H^1(X,\nbigo)^{\lor}\otimes\nbigo_{\Pic}[0]
\Bigr)
\oplus
\Bigl(
 H^0(X,\nbigo)^{\lor}\otimes\nbigo_{\Pic}[-1]
\Bigr)
\]
We also remark that the composite
$\tau_{\leq -1}\Ob(\Poin) \lrarr L_{\Pic}$ is trivial.

\subsection{Relative Obstruction Theory for Orientations}
\label{subsection;06.7.3.32}
\subsubsection{Construction of the complex}
\label{subsubsection;06.5.5.10}

We use the notation 
in the subsubsections
\ref{subsubsection;06.5.2.3}--\ref{subsubsection;06.5.9.20}.
Let $U_1$ be an Artin stack.
Let $F_1:U_1\lrarr U$ be a morphism.
Assume that we have an orientation
$\rho$ of the sheaf $F_{1X}^{\ast}(E)$
over $U_1\times X$.
We have the morphism
$\det_E:U_1\lrarr \Pic$ induced by $\det(E)$.
We denote $\det_E\times\id_X:U_1\times X\lrarr \Pic\times X$
by $\det_{E,X}$.
We have 
$\det_{E,X}^{\ast}\Poin\simeq \det(E)$.
Then we have the following commutative diagram:
\[
 \begin{CD}
 U_1\times X  @>{\Phi(V_{\cdot})}>> Y(W_{\cdot})\\
 @V{\det_{E,X}}VV @V{\gminiw}VV \\
\Pic\times X
@>{\Phi(\Poin)}>> X_{G_m}
 \end{CD}
\]
It induce the following commutative diagram:
\[
 \begin{CD}
 L_{U_1\times X/X}
 @<<<
 \Phi(V_{\cdot})^{\ast} L_{Y(W_{\cdot})/X}
 @<<<
 \gminig(V_{\cdot}) \\
@AAA @AAA @AAA\\
 \det_{E,X}^{\ast}L_{\Pic\times X/X}
 @<<<
 \Phi(\det(E))^{\ast}L_{X_{G_m}/X}
 @<<< 
 \Phi(\det(E))^{\ast}L_{X_{G_m}/X}
 \end{CD}
\]
Since we have $\Phi(\det(E))^{\ast}L_{X_{G_m}/X}=\gminig^d(V_{\cdot})$,
we obtain the following:
\[
 \begin{CD}
 L_{U_1} @<<< \Ob(V_{\cdot})\\
 @AAA @AAA \\
 \det_E^{\ast}L_{\Pic}
 @<<< \Ob^d(V_{\cdot})
 \end{CD}
\]
We put
$\Ob_{\rel}(V_{\cdot},\rho):=
\Cone\bigl(\Ob^d(V_{\cdot})\lrarr
 \det_E^{\ast}L_{\Pic}\bigr)$.
\index{$\Ob_{\rel}(V_{\cdot},\rho)$, 
 $\ob_{\rel}(V_{\cdot},\rho)$}
We have the morphisms
$\Ob_{\rel}(V_{\cdot},\rho)[-1]
\lrarr \Ob^d(V_{\cdot})
\lrarr \Ob(V_{\cdot})$.

Let $\nbigm(1)$ denote the moduli of line bundles,
i.e.,
$\nbigm(1)=\Pic/G_m$.
Let $\pi$ denote the projection $\Pic\lrarr\nbigm(1)$.
We have the following commutative diagrams
induced by $\det(E)$:
\[
 \begin{CD}
 U_1 @>{\det_E}>> \Pic \\
 @V{F_1}VV @V{\pi}VV \\
 U @>{\dettilde_E}>> \nbigm(1)
 \end{CD}
\quad\quad\quad\quad
 \begin{CD}
 U_1\times X @>{\det_{E,X}}>> \Pic\times X\\
 @V{F_{1X}}VV @V{\pi_X}VV \\
 U\times X @>>> \nbigm(1)\times X
 @>>> X_{G_m}
 \end{CD}
\]
Hence, we have the following diagram on $U_1\times X$:
\begin{equation}
 \begin{CD}
 L_{U_1\times X/X}@<<< \det_{E,X}^{\ast}L_{\Pic\times X/X}\\
 @AAA @AAA \\
 F_{1X}^{\ast}L_{U\times X/X} @<<<
 \det_{E,X}^{\ast}\pi_X^{\ast}L_{\nbigm(1)\times X/X} @<<<
 \Phi(\det(E))^{\ast}L_{X_{G_m}/X}
 \end{CD}
\end{equation}
Therefore,
we obtain the following diagram on $U_1$:
\begin{equation}
\label{eq;06.5.10.11}
 \begin{CD}
 L_{U_1}@<<< \det_E^{\ast}L_{\Pic}\\
 @AAA @AAA \\
 F_1^{\ast}L_{U}@<<<
 \det_E^{\ast}\pi^{\ast}L_{\nbigm(1)} @<<<
 \Ob^d(V_{\cdot})
 \end{CD}
\end{equation}
Thus, we obtain the following commutative
diagram:
\begin{equation}
\label{eq;06.5.10.1}
 \begin{CD}
 \Ob_{\rel}(V_{\cdot},\rho)[-1] @>>>
 \Ob(V_{\cdot})\\
 @VVV @VVV \\
 L_{U_1/U}[-1] @>>> F_1^{\ast}L_U
 \end{CD}
\end{equation}

The following lemma can be shown
by an argument similar to Lemma \ref{lem;6.14.6}.
\begin{lem}
\label{lem;6.16.20}
The diagram {\rm(\ref{eq;06.5.10.1})}
depends only on $(E,\rho)$
in the derived category $D(U_1)$.
\hfill\qed
\end{lem}

The following lemma is easy to see by construction
and the argument in the subsubsection
\ref{subsubsection;06.5.11.1}.
\begin{lem}
\label{lem;06.5.11.20}
$\Ob_{\rel}(V_{\cdot},\rho)$ is isomorphic to
$\Bigl(H^0(X,\nbigo)^{\lor}\otimes\nbigo_{U_1}[0]\Bigr)
\oplus
 \Bigl(H^2(X,\nbigo)^{\lor}\otimes\nbigo_{U_1}[2]  \Bigr)$.
The composite of the morphisms
$\bigl(
 \tau_{\leq -2}
 \Ob_{\rel}(V_{\cdot},\rho)
\bigr)[-1]\lrarr \Ob(V_{\cdot})\lrarr L_{U_1}$
is trivial.
\hfill\qed
\end{lem}

\subsubsection{Relative obstruction property}

For any $U$-scheme $g:T\lrarr U$,
let $F_1(T)$ denote the set of 
orientations of $g^{\ast}E$.
Then, we obtain the functor $F_1$ of
the category of $U$-schemes
to the category of sets.
The functor $F_1$ is representable
by the scheme $\Or(E)^{\ast}$.
Let $\pi$ denote the projection
$\Or(E)^{\ast}\lrarr U$.
On $\Or(E)^{\ast}\times X$,
we have the universal orientation
$\rho^u$ of $\pi^{\ast}E$.
From the resolution $V_{\cdot}$
and the orientation $\rho^u$,
we obtain the morphism:
\[
 \ob_{\rel}(V_{\cdot},\rho^u):
 \Ob_{\rel}(V_{\cdot},\rho^u)
\lrarr
 L_{\Or(E)^{\ast}/U}
\]

\begin{lem}
\label{lem;06.5.3.42}
The morphism
$\ob_{\rel}(V_{\cdot},\rho^u)$
gives the relative obstruction theory
for $\Or(E)^{\ast}$ over $U$.
\end{lem}
\pf
We have only to show that
$\nbigh^0(\ob_{\rel}(V_{\cdot},\rho^u))$ is isomorphic.
From the diagram (\ref{eq;06.5.10.11}),
we obtain the following morphisms:
\[
\begin{CD}
 \Ob_{\rel}(V_{\cdot},\rho^u)
@>{\varphi_1}>>
 \det_E^{\ast}\Cone\bigl(
 \pi^{\ast}L_{\nbigm(1)}\lrarr L_{\Pic}
 \bigr)
@>{\varphi_2}>>
L_{\Or(E)^{\ast}/U}
\end{CD}
\]
Since we have the isomorphism
$\Or(E)^{\ast}\simeq
 U\times_{\nbigm(1)}\Pic$,
the morphism $\varphi_2$ is isomorphic.
We have the following commutative diagram:
\[
 \begin{CD}
 \nbigh^0(\Ob^d) @>>> 
\det_E^{\ast}
 \Omega_{\Pic} @>>>
 \nbigh^0\bigl(\Ob_{\rel}(V_{\cdot},\rho^u)\bigr)
 @>>>
 \nbigh^1(\Ob^d)\\
 @V{a_1}VV @V{\simeq}VV @V{a_2}VV 
 @V{a_3}VV \\
 \nbigh^{0}(
 \det_E^{\ast}\pi^{\ast}L_{\nbigm(1)})
 @>>>
 \det_E^{\ast}\Omega_{\Pic}
 @>>>
 \det_E^{\ast}L_{\Pic/\nbigm(1)}
 @>>>
 \det_E^{\ast}\pi^{\ast}\nbigh^1(L_{\nbigm(1)})
 \end{CD}
\]
We also have $\nbigh^1(\det_E^{\ast}L_{\Pic})=0$.
The morphisms
$a_i$ $(i=1,3)$ are isomorphic,
by applying Corollary \ref{cor;06.5.10.12}
to $\nbigm(1)$.
Thus $a_2$ is isomorphic.
Then, we obtain the claim of the lemma.
\hfill\qed

\subsection{Relative Obstruction Theory for $L$-Section}
\label{subsection;06.7.3.33}
\subsubsection{Construction of the complex}
\label{subsubsection;06.5.4.1}

We use the notation in the subsubsection
\ref{subsubsection;06.5.2.3}.
Let $L$ be a line bundle on $X$.
Let $P_{\cdot}:=
 \bigl(P_{-1}\stackrel{\del_L}{\lrarr}P_0\bigr)\simeq L$
be a locally free resolution,
where $P_0$ stands in the degree $0$.
We have the natural right
$\GL(W_{\cdot})$-action
on $N\bigl(P_i,W_j\bigr)$
given by
$(g_{-1},g_0)\cdot f=
 g_j^{-1}\circ f$.
It induces the $\GL(W_{\cdot})$-actions
on the varieties
$N\bigl(P_{-1},W_{0\,X}\bigr)$,
$X$ and
$N(W_{-1\,X},W_{0\,X})\times_X
 N(P_0,W_{0\,X})\times_XN(P_{-1},W_{-1\,X})$.
The quotient stacks are denoted by
$Y_{0}(W_{\cdot},P_{\cdot})$,
$Y_1(W_{\cdot},P_{\cdot})$,
and $Y_2(W_{\cdot},P_{\cdot})$
respectively.

We have the equivariant map 
$h:N(W_{-1\,X},W_{0\,X})\times_X
 N(P_{-1},W_{-1\,X})\times_X
 N(P_{0},W_{0\,X})\lrarr
 N(P_{-1},W_{0\,X})$
given by
$h(f,a_{-1},a_0)
=f\circ a_{-1}-a_0\circ\del_L$.
Since $\del_{L|x}$ is injective
for any point $x\in X$,
the map $h$ is smooth.
It induces the smooth morphism
$Y_{2}(W_{\cdot},P_{\cdot})
\lrarr Y_0(W_{\cdot},P_{\cdot})$.
We also have the morphism
$Y_1(W_{\cdot},P_{\cdot})\lrarr
 Y_0(W_{\cdot},P_{\cdot})$
induced by the $0$-section
$X\lrarr N(P_{-1},W_{0,X})$.
We denote the fiber product
$Y_1\times_{Y_0}Y_2$
by $Y(W_{\cdot},P_{\cdot})$.
\index{$Y(W_{\cdot},P_{\cdot})$}

\vspace{.1in}

Let $U_2$ be an algebraic stack
with a morphism $F_2:U_2\lrarr U$
and an $L$-section
$\phi:p_{U_2}^{\ast}(L)\lrarr F_{2\,X}^{\ast}(E)$.
We assume that
we have a morphism of complexes
$\phitilde_{\cdot}=(\phitilde_{-1},\phitilde_{0}):
 p_{U_2}^{\ast}(P_{\cdot})
\lrarr F_{2\,X}^{\ast}(V_{\cdot})$
which gives $\phi$ in the cohomology level.
Such a $\phitilde$ is called a lift of $\phi$.
We put as follows:
\[
\gminig_{\rel}(V_{\cdot},\phitilde)
=
\nhom\bigl(
p_{U_2}^{\ast}P_{\cdot},\,
 F_{2,X}^{\ast}V_{\cdot}
\bigr)^{\lor}
\] \index{$\gminig_{\rel}(V_{\cdot},\phitilde)$}
We have the naturally induced morphism
$\gamma(\phitilde_{\cdot}):
 \gminig_{\rel}(V_{\cdot},\phitilde)[-1]\lrarr 
 \gminig(V_{\cdot})$
and 
$\gamma(\phitilde_{\cdot})_{\leq 1}:
 \gminig_{\rel}(V_{\cdot},\phitilde)[-1]_{\leq 1}
 \lrarr \gminig(V_{\cdot})_{\leq 1}$.

We have the classifying map
$\Phi(V_{\cdot},\phitilde):
 U_2\times X\lrarr Y(W_{\cdot},P_{\cdot})$
which induces the maps
$\Phi_i(V_{\cdot},\phitilde):
 U_2\times X\lrarr Y_i(W_{\cdot},P_{\cdot})$.

\begin{lem}
\label{lem;06.5.1.10}
$\Phi(V_{\cdot},\phitilde)^{\ast}
 L_{Y(W_{\cdot},P_{\cdot})/X}$
is expressed by the complex
$\Cone\bigl(\gamma(\phitilde_{\cdot})_{\leq 1}\bigr)$.
\end{lem}
\pf
We have the induced morphisms
$\kappa_i:\Phi_0(V_{\cdot},\phitilde)^{\ast}
 L_{Y_0(W_{\cdot},P_{\cdot})/X}
\lrarr
\Phi_i(V_{\cdot},\phitilde)^{\ast}
 L_{Y_i(W_{\cdot},P_{\cdot})/X}$
$(i=1,2)$.
Since $Y_2(W_{\cdot},P_{\cdot})\lrarr
 Y_0(W_{\cdot},P_{\cdot})$
is smooth,
$\Phi(V_{\cdot},\phitilde)^{\ast}
 L_{Y(W_{\cdot},P_{\cdot})/X}$
is isomorphic to
the cone of the induced morphism
$\Phi_0(V_{\cdot},\phitilde)^{\ast}
 L_{Y_0(W_{\cdot},P_{\cdot})/X}
\lrarr
 \bigoplus_{i=1,2}
 \Phi_i(V_{\cdot},\phitilde)^{\ast}
 L_{Y_i(W_{\cdot},P_{\cdot})/X}$.

Let us see 
$\Phi_i(V_{\cdot},\phitilde)^{\ast}
 L_{Y_i(W_{\cdot},P_{\cdot})/X}$.
We use the argument explained
in the subsubsection
\ref{subsubsection;06.4.29.15}.
We will omit to denote $p_{U_2}^{\ast}$
and $F_{2\,X}^{\ast}$.
In the case $i=2$,
it is expressed by the following complex:
\begin{equation}
 \label{eq;06.5.1.1}
\begin{array}{c}
 \nhom(V_0,V_{-1})
\oplus
 \nhom(V_{-1},P_{-1})
\oplus
 \nhom(V_{0},P_0)
\lrarr
 \nhom(V_0,V_0)
\oplus
 \nhom(V_{-1},V_{-1})\\
\mbox{{}}\\
 (b,c_{-1},c_0)
\longmapsto
 \bigl(
 f\circ b+\phitilde_0\circ c_0,\,\,\,
 -b\circ f +\phitilde_{-1}\circ c_{-1} 
 \bigr)
\end{array}
\end{equation}
Here, the first term stands in the degree $0$.
In the case $i=0$,
it is expressed by the following complex:
\begin{equation}
 \label{eq;06.5.1.2}
 \nhom(V_0,P_{-1})
\longmapsto 
 \nhom\bigl(V_{0},V_0\bigr)
\oplus
 \nhom\bigl(V_{-1},V_{-1}\bigr),
\quad
 c_2\longmapsto
 \bigl(0,0
 \bigr)
\end{equation}
Again, the first term stands in the degree $0$.
In the case $i=1$,
it is expressed by the following:
\begin{equation}
 \label{eq;06.5.1.3}
0\lrarr
 \nhom(V_0,V_0)\oplus \nhom(V_{-1},V_{-1})
\end{equation}
Here the term $0$ stands in the degree $0$.
For the description
(\ref{eq;06.5.1.1})
and (\ref{eq;06.5.1.2}),
the degree $0$-part of $\kappa_2$ is 
given as follows:
\[
 \nhom(V_0,P_{-1})\lrarr
 \nhom(V_0,V_{-1})\oplus
 \nhom(V_{-1},P_{-1})\oplus
 \nhom(V_0,P_0),
\quad
 a\longmapsto
\bigl(
 \phi_{-1}\circ a,\,\,
 a\circ f,\,\,
  -\del_L\circ a
\bigr).
\]
The degree $1$-part of $\kappa_2$
is given by the identity.
On the other hand,
the morphism $\kappa_1$ is
the obvious one for the descriptions
(\ref{eq;06.5.1.2}) and (\ref{eq;06.5.1.3}).
Then the claim of the lemma can be checked directly.
\hfill\qed

\vspace{.1in}

We put $\gminig(V_{\cdot},\phitilde_{\cdot}):=
 \Cone\bigl(\gamma(\phitilde_{\cdot})\bigr)$.
Then, we obtain the following commutative diagram:
\[
 \begin{CD}
 \gminig(V_{\cdot}) @>>>
 \Phi\bigl(F_{2\,X}^{\ast}V_{\cdot}\bigr)^{\ast}L_{Y(W_{\cdot})/X}
 @>>> F_{2,X}^{\ast}L_{U\times X/X}\\
 @VVV @VVV @VVV\\
 \gminig(V_{\cdot},\phitilde_{\cdot})
 @>>>
 \Phi(V_{\cdot},\phitilde_{\cdot})^{\ast}
  L_{Y(W_{\cdot},P_{\cdot})/X}
 @>>>
 L_{U_2\times X/X}
 \end{CD}
\]
We obtain the following morphism
of the distinguished triangles on $U_2\times X$:
\[
  \begin{CD}
 \gminig(V_{\cdot})@>>>
 \gminig(V_{\cdot},\phitilde_{\cdot}) @>>>
 \gminig_{\rel}(V_{\cdot},\phitilde_{\cdot}) @>>>
 \gminig(V_{\cdot})[1]\\
 @VVV @VVV @VVV @VVV \\
 F_{2\,X}^{\ast}L_{U\times X/X} @>>>
 L_{U_2\times X/X} @>>>
 L_{U_2\times X/U\times X} @>>>
 F_{2\,X}^{\ast}L_{U\times X/X}[1].
 \end{CD}
\]

We put
$\Ob_{\rel}\bigl(V_{\cdot},\phitilde_{\cdot}\bigr):=
 Rp_{X\ast}\bigl(
\gminig_{\rel}(V_{\cdot},\phitilde_{\cdot})
 \otimes\omega_X\bigr)$.
\index{$\Ob_{\rel}\bigl(V_{\cdot},\phitilde_{\cdot}\bigr)$}
Then we obtain the following diagram on $U_2$:
\begin{equation} 
 \label{eq;6.14.7}
\begin{CD}
\Ob_{\rel}\bigl(V_{\cdot},\phitilde_{\cdot}\bigr)[-1] 
 @>>> \Ob(V_{\cdot})\\
@VVV @VVV \\
L_{U_2/U}[-1] @>>> F_{2\,X}^{\ast}L_{U/k}
\end{CD}
\end{equation}

\begin{lem}
 \label{lem;6.14.14}
The diagram {\rm(\ref{eq;6.14.7})} 
depends only on $(E,\phi)$
in the derived category
$D(U_2)$ .
\end{lem}
\pf
Let 
$\bigl(
 V^{(1)}_{\cdot},P^{(1)}_{\cdot},\phitilde^{(1)}
\bigr)$ be another choice.
We take the resolution $V^{(2)}_{\cdot}$ of $E$
as in the proof of Lemma \ref{lem;6.14.6}.
We put $P^{(2)}_0=P^{(1)}_0\oplus P_0$
and $P^{(2)}_{-1}=\ker(P^{(2)}_0\lrarr L)$.
Then the lift
$\phitilde^{(2)}:P^{(2)}_{\cdot}\lrarr 
 F_{2,X}^{\ast}V^{(2)}_{\cdot}$
is naturally obtained from
the lifts $\phitilde^{(1)}$ and $\phitilde$.
We have the compatible inclusions
$P_{\cdot}\lrarr P^{(2)}_{\cdot}$
and $V_{\cdot}\lrarr V^{(2)}_{\cdot}$.
Then we can show the claim of the lemma
by using the filtered objects
as in Lemma \ref{lem;6.14.6}.
\hfill\qed

\subsubsection{Relative obstruction property}
\label{subsubsection;06.5.4.50}

For any $U$-scheme $g:T\lrarr U$,
let $F(T)$ denote the set of
$L$-sections of $g_X^{\ast}E$.
Thus we obtain the functor of the category
of $U$-schemes to the category of sets.
The functor $F$ is representable
by a scheme $M(L)$.
Let $M(L)\lrarr U$ denote the projection.
On $M(L)\times X$,
we have the universal $L$-section $\phi$
of $\pi_X^{\ast}E$.
Assume that we have a locally free
resolution $P_{\cdot}$ of $L$
for which we have a lift
$\phitilde:P_{\cdot}\lrarr V_{\cdot}$
of $\phi$.
Then we obtain the morphism:
\[
 \ob_{\rel}(V_{\cdot},\phitilde):
 \Ob_{\rel}(V_{\cdot},\phitilde)\lrarr L_{M(L)/U}
\]

\begin{lem}
 \label{lem;06.5.3.2}
$\ob_{\rel}(V_{\cdot},\phitilde)$
gives the relative obstruction theory
for $M(L)$ over $U$.
\end{lem}
\pf
We have only to show the claim
on any sufficiently small open subsets $\nbigu$ of $M(L)$.
Recall that $\ob_{\rel}(V_{\cdot},\phitilde)$
is independent of the choice of $P_{\cdot}$
and $\phitilde$.
Therefore,
we may assume
$H^i\bigl(X,P_{0}^{-1}\otimes
 V_{l\,|\,\{u\}\times X}\bigr)=0$
for $i=0,1$, $l=0,-1$
and for any $u\in \nbigu$.

For simplicity of the notation,
we put $C_0:=\gminig_{\rel}(V_{\cdot},\phitilde)$
and $C_1:=\gminig_{\rel}(V_{\cdot},\phitilde)_{\leq 0}$.
We have the obvious map $C_0\lrarr C_1$.
By the argument in the subsubsection
\ref{subsubsection;06.4.29.15},
the complex $C_1$ expresses
$\Psi(V_{\cdot},\phitilde)^{\ast}
 L_{Y(W_{\cdot},P_{\cdot})/Y(W_{\cdot})}$.
Thus we have the natural morphism
$C_1\lrarr L_{\nbigu\times X/U\times X}$.
We put
$\Ob_1:=Rp_{X\,\ast}\bigl(C_1\otimes\omega_X\bigr)$.
Then we have the induced morphisms
$\Ob_{\rel}(V_{\cdot},\phitilde)
\lrarr \Ob_{1}$
and $\ob_1:\Ob_1\lrarr L_{\nbigu/U}$.
It is easy to see that the composite
of the morphisms
is same as $\ob_{\rel}(V_{\cdot},\phitilde)$.
We use the following lemma.

\begin{lem}
The morphism 
$\nbigh^0\bigl(\Ob_{\rel}(V_{\cdot},\phitilde)\bigr)
\lrarr \nbigh^0\bigl(\Ob_1\bigr)$
is isomorphic,
and 
$\nbigh^1\bigl(\Ob_{\rel}(V_{\cdot},\phitilde)\bigr)
\lrarr \nbigh^1\bigl(\Ob_1\bigr)$
is surjective.
\end{lem}
\pf
We have the exact sequence of the complexes
$ 0\lrarr \nhom\bigl(V_{-1},P_0\bigr)[-1]
 \lrarr C_0\lrarr C_1\lrarr 0$.
Due to our choice of $P_{\cdot}$,
we have
$H^i\bigl(X, P_0^{-1}\otimes V_{-1\,|\,\{u\}\times X}\bigr)=0$
for any $u\in \nbigu$.
Then the claim can be easily shown.
\hfill\qed

\vspace{.1in}

Therefore,
we have only to show that
$\ob_1$ gives the obstruction theory 
for $\nbigu$ over $U$.
For that purpose,
we have only to check the conditions
(A1) and (A2) in
Proposition \ref{prop;06.4.29.10}.
We use the same argument as that
in the proof of Lemma \ref{lem;06.5.3.1}.
Let $T$ be a scheme with a morphism
$h:T\lrarr \nbigu$.
Then, we have the following commutative
diagram,
for any coherent sheaf $J$ on $T$:
\[
 \begin{CD}
 \Ext^i\bigl(h^{\ast}L_{\nbigu/U},J\bigr)
 @>{\phi}>>
 \Ext^i\bigl(h^{\ast}\Ob_1,J\bigr)\\
 @VVV @A{\simeq}AA \\
 \Ext^i\bigl(h_X^{\ast}
 L_{\nbigu\times X/U\times X},J_X\bigr)
 @>>>
 \Ext^i\bigl(h_X^{\ast}C_1,\,J_X \bigr)
 \end{CD}
\]
Let $\Tbar$ be a scheme
such that
$T$ is a closed subscheme of $\Tbar$
whose corresponding ideal sheaf $J$
is square $0$.
Let $h'$ be a morphism
$\Tbar\lrarr U$ such that
the restriction $h'_{|T}$ is same
as $\pi\circ h$, where $\pi$ denotes
the projection $\nbigu\lrarr U$.
We have the obstruction class
$o(h,h')\in \Ext^1\bigl(h^{\ast}L_{\nbigu/U},\,J\bigr)$.
We put
$\widehat{h}_X:=
 \Phi(V_{\cdot},\phitilde_{\cdot})\circ h_X:
 T\times X\lrarr Y(W_{\cdot},\phitilde)$ 
and 
$\widehat{h}_X':=
 \Phi(V_{\cdot})\circ h'_X:\Tbar\times X\lrarr
 Y(W_{\cdot})$.
Since the complex $C_1$ expresses
$\Psi(V_{\cdot},\phitilde)^{\ast}
 L_{Y(W_{\cdot},P_{\cdot})/Y(W_{\cdot})}$,
we have the obstruction class
$o(\widehat{h}_X,\widehat{h}'_X)
 \in \Ext^1\bigl(h^{\ast}C_1,J_X\bigr)$.
By the functoriality,
$o(h,h')$ is mapped to 
$o(\widehat{h}_X,\widehat{h}'_X)$ in the diagram.
Hence, $\phi\bigl(o(h,h')\bigr)=0$ implies
that $\widehat{h}_X$ can be extended
to a morphism
$\widehat{h}_{1,X}:\Tbar\lrarr
 Y(W_{\cdot},P_{\cdot})$,
which is a lift of $\widehat{h}'_X$.
Hence, we obtain an $L$-section of
$h_X^{\prime\,\ast}E$.
Then we obtain an extension
of $h$ to a morphism
$h_1:\Tbar\lrarr \nbigu$
which is a lift of $h'$,
due to the universal property of 
$\nbigm(L)$.
Thus the condition (A1) is checked.
The condition (A2) can also be checked easily,
and the proof of Lemma \ref{lem;06.5.3.2}
is finished.
\hfill\qed

\subsubsection{Preparation for the obstruction theory
 of Master space}
\label{subsubsection;06.5.21.250}

We will use the notation in the subsubsections
\ref{subsubsection;06.5.4.100}
and \ref{subsubsection;06.5.4.1}.
We have the naturally defined right
$\GL(W_0)$-action on $N\bigl(\nbigo_X,W_{0,X}\bigr)$.
The quotient stack is denoted by 
$A(W_{\cdot},P_{\cdot})$.
Assume that we have a morphism
$\iotatilde:\nbigo_X\lrarr P_0$.
Then we have the induced morphism
$Y_1(W_{\cdot},P_{\cdot})\lrarr A(W_{\cdot},P_{\cdot})$,
and thus
$\Gamma_L:Y(W_{\cdot},P_{\cdot})
\lrarr A(W_{\cdot},P_{\cdot})$.
From $(F_{2X}^{\ast}V_{\cdot},\phitilde)$ on 
$U_2\times X$ and $\iotatilde$,
we obtain the morphism
$\phibar:p_{U_2}^{\ast}\nbigo_X
\lrarr F_{2X}^{\ast}V_0$.
We put
$\gminih_{\rel}(V_{\cdot},\phitilde)
:=\nhom\bigl(p_{U_2}^{\ast}\nbigo_X,\,
 F_{2X}^{\ast}V_{\cdot}\bigr)^{\lor}$.
\index{$\gminih_{\rel}(V_{\cdot},\phitilde)$}
Then, we have the induced maps
$\gamma(\phibar):
 \gminih_{\rel}(V_{\cdot},\phitilde)[-1]\lrarr
 \gminih(V_{\cdot})$.

We put 
$\Psi(V_{\cdot},\phitilde_{\cdot}):=
 \Gamma_L\circ\Phi(V_{\cdot},\phitilde)$,
which is the classifying map
of $\bigl(F_{2\,X}^{\ast}V_0,\phibar\bigr)$.
Then, 
$\Psi(V_{\cdot},\phitilde_{\cdot})^{\ast}
 L_{A(W_{\cdot},P_{\cdot})/X}$
is expressed by 
$\cone\bigl(\gamma(\phibar)_{\leq 1}\bigr)$.
The induced morphism
$\Psi(V_{\cdot},\phitilde_{\cdot})^{\ast}
 L_{A(W_{\cdot},P_{\cdot})/X}
\lrarr 
 \Phi(V_{\cdot},\phitilde_{\cdot})^{\ast}
 L_{Y(W_{\cdot},P_{\cdot})/X}$
is expressed by the naturally given morphism
$\Cone\bigl(\gamma(\phibar)_{\leq 1}\bigr)
\lrarr
 \Cone\bigl(\gamma(\phitilde_{\cdot})_{\leq 1}\bigr)$.

We put 
$\gminih(V_{\cdot},\phitilde):=
 \cone\bigl(\gamma(\phibar)\bigr)$.
\index{$\gminih(V_{\cdot},\phitilde)$}
Then, we obtain the following commutative diagram:
\[
 \begin{CD}
 \gminih(V_{\cdot},\phitilde) @>>> 
 \gminig(V_{\cdot},\phitilde) @>>>
 L_{U_2\times X/X}\\
 @AAA @AAA @AAA \\
 \gminih(V_{\cdot}) @>>>
 \gminig(V_{\cdot}) @>>>
 F_{2\,X}^{\ast}L_{U\times X/X}
 \end{CD}
\]
We put
$ \Ob^G_{\rel}(V_{\cdot},\phitilde):=
 Rp_{X\,\ast}\bigl(
 \gminih_{\rel}(V_{\cdot},\phitilde)\otimes\omega_X \bigr)$.
Then, we have the following diagram on $U_2$:
\begin{equation}
 \begin{CD}
 L_{U_2/U}[-1] @<<<
 \Ob_{\rel}(V_{\cdot},\phitilde)[-1] @<<< 
 \Ob^G_{\rel}(V_{\cdot},\phitilde)[-1]\\
 @VVV @VVV @VVV \\
 F_{2}^{\ast}L_{U/k}@<<<
 \Ob(V_{\cdot}) @<<< 
 \Ob^G(V_{\cdot})
 \end{CD}
\end{equation}

\vspace{.1in}

Now, we assume the condition (C)
in the subsubsection \ref{subsubsection;06.5.4.100}.
We have the natural right $\GL(W_0)$-action
on $N(k,W_0)$.
The quotient stack is denoted by
$B(W_{\cdot},P_{\cdot})$.
\index{$B(W_{\cdot},P_{\cdot})$}
We have the natural isomorphism
$J_L:B(W_{\cdot},P_{\cdot})\times X\lrarr A(W_{\cdot},P_{\cdot})$.
From $\phibar:p_{U_2}^{\ast}\nbigo_X
\lrarr V_0$,
we obtain the morphism
$\Xi(V_{\cdot},\phibar):
 U_2\lrarr B(W_{\cdot},P_{\cdot})$.
Note that the composite
of $\Xi(V_{\cdot},\phibar)_X$
and $J_L$ is same as
$\Psi(V_{\cdot},\phibar)$.
Therefore,
we have the following commutative diagram:
\[
 \begin{CD}
 U_2\times X @>>>
 B(W_{\cdot},P_{\cdot})\times X
 @>>>
 A(W_{\cdot},P_{\cdot}) \\
 @VVV @VVV @VVV \\
 U\times X @>>> 
 B(W_{\cdot})\times X@>>> 
 A(W_{\cdot})
 \end{CD}
\]
Thus, we obtain the following commutative
diagram:
\[
 \begin{CD}
 F_2^{\ast}L_{U/k} @<<<
 \Phi(V_{\cdot})^{\ast}L_{B(W_{\cdot})/k} @<{\tau_1}<{\simeq}<
 \Ob^G(V_{\cdot}) \\
 @AAA @AAA @AAA \\
 L_{U_2/U}[-1] @<<<
 \Xi(V_{\cdot},\phibar)^{\ast}
 L_{B(W_{\cdot},P_{\cdot})/B(W_{\cdot})}[-1]
 @<{\tau_2}<<
 \Ob^G_{\rel}(V_{\cdot},\phitilde)
 \end{CD}
\]
\begin{lem}
\label{lem;06.5.2.10}
The morphism $\tau_2$ is isomorphic.
\end{lem}
\pf
The complex 
$\gminih_{\rel}(V_{\cdot},\phitilde)$
is quasi isomorphic to
$\bigl(\nhom(V_0,\nbigo)\rarr 
 \nhom(V_{-1},\nbigo)\bigr)$,
where the first term stands at the degree $0$.
And the degree $0$-part of 
the morphism
$\gminih_{\rel}(V_{\cdot},\phitilde)_{\leq 1}
\lrarr
 \Xi(V,\phibar)_X^{\ast}L_{B(W_{\cdot},P_{\cdot})/B(W_{\cdot})}
\simeq
 \nhom(V_0,\nbigo)$
is given by the identity of $\nhom(V_0,\nbigo)$.
Hence, the claim can be shown
as in Lemma \ref{lem;06.5.4.201}.
\hfill\qed

\subsubsection{Preparation for Proposition
       \ref{prop;06.5.10.50}}

We have the morphism of the complexes
$\gminif:\gminig_{\rel}(V_{\cdot},\phi)
\lrarr \nbigo_{U_2\times X}$
given as follows:
\[
 \nhom(F_{2X}^{\ast}V_0,p_{U_2}^{\ast}P_0)
\oplus
 \nhom(F_{2X}^{\ast}V_{-1},p_{U_2}^{\ast}P_{-1})
\lrarr \nbigo_{U_2\times X},
\quad
 (a_0,a_{-1})\longmapsto
 \tr\bigl(\phitilde_0\circ a_0\bigr)
+\tr\bigl(\phitilde_{-1}\circ a_{-1}\bigr)
\]
It is easy to check that the following diagram
is commutative:
\[
 \begin{CD}
 \gminig_{\rel}(V_{\cdot},\phitilde)[-1]
 @>>> \gminig(V_{\cdot})\\
 @V{\gminif}VV @V{\tr}VV \\
 \nbigo_{U_2\times X}[-1]@>{\id}>>\nbigo_{U_2\times X}[-1]
 \end{CD}
\]
It induces the following commutative diagram:
\begin{equation}
\label{eq;06.5.10.20}
\begin{CD}
 \Ob_{\rel}(V_{\cdot},\phitilde)[-1]
@>>>
 \Ob(V_{\cdot})\\
 @VVV @VVV \\
\nbigo_{U_2\times X}[-1]@>>>\nbigo_{U_2\times X}[-1]
\end{CD}
\end{equation}

\subsection{Relative Obstruction Theory for Reduced $L$-Section}
\label{subsection;06.7.3.34}
\subsubsection{Construction of the complex}
\label{subsubsection;06.5.4.110}

We use the notation 
in the subsubsection \ref{subsubsection;06.5.4.1}.
The weight $(-1)$-action of $G_m$
on $P_i$ induces the $G_m$-action
on $Y_i(W_{\cdot},P_{\cdot})$
$(i=0,1,2)$
and $Y(W_{\cdot},P_{\cdot})$.
The quotient stacks are denoted by
$Y_i(W_{\cdot},[P_{\cdot}])$
$(i=0,1,2)$
and $Y(W_{\cdot},[P_{\cdot}])$.
\index{$Y(W_{\cdot},[P_{\cdot}])$}
We have the morphism
$\pi_3:Y(W_{\cdot},[P_{\cdot}])
\lrarr X_{G_m}$
induced by $Y(W_{\cdot},P_{\cdot})\lrarr X$.

Let $U_3$ be a stack
with a morphism $F_3:U_3\lrarr U$
and a reduced $L$-section
$[\phi]:p_{U_3}^{\ast}(L)\otimes p_X^{\ast}(M)
\lrarr F_{3\,X}^{\ast}(E)$,
where $M$ denotes a line bundle on $U_3$.
Let $\Phi_M:U_3\lrarr \Spec(k)_{G_m}$ 
denote the classifying map for $M$.

\begin{notation}
When we have a map
$g:U_3\lrarr T$,
then $\Phi_M$ induces the morphism
$U_3\lrarr T\times \Spec(k)_{G_m}$,
which is denoted by $g_M$
in the following argument.
\hfill\qed
\end{notation}

We assume that
we have a morphism of complexes
$[\phitilde_{\cdot}]=
\bigl([\phitilde_{-1}],[\phitilde_{0}]\bigr):
 p_{U_3}^{\ast}(P_{\cdot})\otimes p_X^{\ast}M
\lrarr F_{3\,X}^{\ast}(V_{\cdot})$
which gives $[\phi]$ in the cohomology level.
Such a map $[\phitilde]$ is called a lift of $[\phi]$.
We put as follows:
\[
 \gminig'_{\rel}(V_{\cdot},[\phitilde])
=
\nhom\bigl(
p_U^{\ast}P_{\cdot}\otimes p_X^{\ast}M,\,\,
F_{3,X}^{\ast}V_{\cdot}
\bigr)^{\lor}
\] \index{$\gminig'_{\rel}(V_{\cdot},[\phitilde])$}
We have the naturally induced morphism
$\gamma[\phitilde_{\cdot}]:
 \gminig'_{\rel}(V_{\cdot},[\phitilde])[-1]
\lrarr \gminig'(V_{\cdot})$
and 
$\gamma[\phitilde_{\cdot}]_{\leq 1}:
 \gminig'(V_{\cdot},[\phitilde])[-1]_{\leq 1}
 \lrarr \gminig'(V_{\cdot})_{\leq 1}$.

We have the classifying map
$\Phi\bigl(V_{\cdot},[\phitilde]\bigr):
 U_3\times X\lrarr Y(W_{\cdot},[P_{\cdot}])$.
We consider the trivial $G_m$-actions
on $U\times X$,  $Y(W_{\cdot})$ and $X$,
and the quotient stacks are denoted by
$(U\times X)_{G_m}$,
$Y(W_{\cdot})_{G_m}$
and $X_{G_m}$ respectively.
Then, we have the following commutative diagram:
\[
 \begin{CD}
 U_3\times X @>{\Phi(V_{\cdot},[\phitilde_{\cdot}])}>>
 Y(W_{\cdot},[P_{\cdot}])
 @>>> X_{G_m}\\
 @V{F_{3XM}}VV @V{\pi_3}VV @VVV \\
 (U\times X)_{G_m}@>>> Y(W_{\cdot})_{G_m}
 @>>> X_{G_m}
 \end{CD}
\]
The composite
 $\pi_3\circ\Phi(V_{\cdot},[\phitilde_{\cdot}])$
is same as $\Phi(F_{3,X}^{\ast}V_{\cdot})_M$.
The following lemma can be shown
by the same argument
as the proof of Lemma \ref{lem;06.5.1.10}.

\begin{lem}
$\Phi(V_{\cdot},[\phitilde_{\cdot}])^{\ast}
 L_{Y(W_{\cdot},[P_{\cdot}])/X_{G_m}}$
is expressed by the complex
$\Cone\bigl(\gamma[\phitilde_{\cdot}]_{\leq 1}\bigr)$.
\hfill\qed
\end{lem}

We put $\gminig'(V_{\cdot},[\phitilde_{\cdot}]):=
 \Cone\bigl(\gamma[\phitilde_{\cdot}]\bigr)$.
\index{$\gminig'(V_{\cdot},[\phitilde_{\cdot}])$}
Then we obtain the following commutative diagram:
\[
 \begin{CD}
 \gminig(V_{\cdot}) @>>>
 \Phi\bigl(F_3^{\ast}V_{\cdot}\bigr)_{M}^{\ast}
 L_{Y(W_{\cdot})_{G_m}/X_{G_m}}
 @>>> F_{3,X,M}^{\ast}L_{(U\times X)_{G_m}/X_{G_m}}\\
 @VVV @VVV @VVV\\
 \gminig'(V_{\cdot},[\phitilde_{\cdot}])
 @>>>
 \Phi(V_{\cdot},\phitilde_{\cdot})^{\ast}
  L_{Y(W_{\cdot},[P_{\cdot}])/X_{G_m}}
 @>>>
 L_{U_3\times X/X_{G_m}}
 \end{CD}
\]
We remark
$F_{3,X,M}^{\ast}L_{(U\times X)_{G_m}/X_{G_m}}$
and 
$\Phi\bigl(F_{3\,X}^{\ast}
 V_{\cdot}\bigr)_{M}^{\ast}L_{Y(W_{\cdot})_{G_m}/X_{G_m}}$
are naturally isomorphic to
$F_{3\,X,M}^{\ast}L_{U\times X/X}$
and $\Phi\bigl(F_{3\,X}^{\ast}V_{\cdot}\bigr)^{\ast}
 L_{Y(W_{\cdot})/X}$,
respectively.
Then, we obtain the following morphism
of the distinguished triangles on $U_3\times X$:
\[
  \begin{CD}
 \gminig(V_{\cdot})@>>>
 \gminig'(V_{\cdot},[\phitilde_{\cdot}]) @>>>
 \gminig'_{\rel}(V_{\cdot},[\phitilde_{\cdot}]) @>>>
 \gminig(V_{\cdot})[1]\\
 @VVV @VVV @VVV @VVV \\
 F_{3\,X}^{\ast}L_{(U\times X)_{G_m}/X_{G_m}} @>>>
 L_{U_3\times X/X_{G_m}} @>>>
 L_{U_3\times X/(U\times X)_{G_m}} @>>>
 F_{3\,X}^{\ast}L_{(U\times X)_{G_m}/X_{G_m}}[1].
 \end{CD}
\]

We put
$\Ob_{\rel}'\bigl(V_{\cdot},[\phitilde]\bigr):=
 Rp_{X\,\ast}\bigl(
 \gminig'_{\rel}(V_{\cdot},[\phitilde])
 \bigr)$.
\index{$\Ob_{\rel}'\bigl(V_{\cdot},[\phitilde]\bigr)$}
Then we obtain the following diagram
on $U_3$:
\begin{equation} \label{eq;6.14.12}
\begin{CD}
\Ob_{\rel}'\bigl(V_{\cdot},[\phitilde]\bigr)[-1] @>>>
L_{U_3/U_{G_m}}[-1] @.\\
@VVV @VVV @.\\
\Ob(V_{\cdot})@>>>
F_{3,M}^{\ast}L_{U_{G_m}/k}@>>>
F_{3,M}^{\ast}L_{U_{G_m}/U}
\end{CD}
\end{equation}
We obtain the morphism
$\Ob_{\rel}'\bigl(V_{\cdot},[\phitilde_{\cdot}]\bigr)[-1]
\lrarr F_{3,M}^{\ast}L_{U_{G_m}/U}$,
by taking the composite of the morphisms
in the diagram above.
We put as follows:
\[
 \Ob_{\rel}(V_{\cdot},[\phitilde]):=
 \Cone\bigl(
 \Ob_{\rel}'(V_{\cdot},[\phitilde])
\lrarr F_{3,M}^{\ast}L_{U_{G_m}/U}[1]
 \bigr)[-1].
\] \index{$\Ob_{\rel}(V_{\cdot},[\phitilde])$}
We obtain the following morphism
of the distinguished triangles:
\[
 \begin{CD}
 F_{3,M}^{\ast}L_{U_{G_m}/U}
@>>>
 \Ob_{\rel}(V_{\cdot},[\phitilde])
@>>>
 \Ob'_{\rel}(V_{\cdot},[\phitilde])
@>>>
 F_{3,M}^{\ast}L_{U_{G_m}/U}[1] \\
 @VVV @VVV @VVV @VVV \\
F_{3,M}^{\ast}L_{U_{G_m}/U}
@>>>
 L_{U_3/U}
@>>>
 L_{U_3/U_{G_m}}
@>>>
F_{3,M}^{\ast}L_{U_{G_m}/U}[1]
 \end{CD}
\]
Therefore, we obtain
the following commutative diagram
on $U_3$:
\begin{equation} \label{eq;6.14.10}
\begin{CD}
\Ob_{\rel}\bigl(V_{\cdot},[\phitilde]\bigr)[-1]
@>>> \Ob(V_{\cdot})\\
@V{\ob_{\rel}(V_{\cdot},[\phitilde])}VV @VVV \\
L_{U_3/U}[-1] @>>> F_3^{\ast}L_{U}.
\end{CD}
\end{equation}

\begin{lem}
\label{lem;6.14.27}
The diagram {\rm(\ref{eq;6.14.10})}
depends only on $(E,[\phi])$.
\end{lem}
\pf
By an argument similar to the proof of Lemma \ref{lem;6.14.14},
we can show that the diagram (\ref{eq;6.14.12})
is independent of a choice of 
$(V_{\cdot},[P_{\cdot}],[\phitilde])$.
Then the diagram (\ref{eq;6.14.10}) is also independent.
\hfill\qed

\subsubsection{Relative obstruction property}

For any $U$-scheme $g:T\lrarr U$,
let $F(T)$ denote the set of
the reduced $L$-sections of $g_X^{\ast}E$.
Thus, we obtain the functor of
the category of $U$-schemes to
the category of sets.
The functor is representable by a scheme $M[L]$.
Let $\pi:M[L]\lrarr U$ denote the projection.
We have the line bundle $\nbigi$ on $M[L]$,
and the universal reduced $L$-section
$[\phi^u]:p_X^{\ast}\nbigi\otimes p_{M[L]}^{\ast}L
\lrarr  \pi_X^{\ast}E$
over $M[L]\times X$.
Assume that we have a locally free
resolution $P_{\cdot}$ of $L$
for which we have a lift 
$[\phitilde]:p_{X}^{\ast}\nbigi\otimes
 p_{M[L]}^{\ast}P_{\cdot} \lrarr \pi_X^{\ast}V_{\cdot}$
of $[\phi]$.
Then,
we obtain the morphism
$\ob_{\rel}\bigl(V_{\cdot},[\phitilde_{\cdot}]\bigr):
 \Ob_{\rel}\bigl(V_{\cdot},[\phitilde_{\cdot}]\bigr)
\lrarr L_{M[L]/U}$.
We also have the complex
$\Ob'_{\rel}\bigl(V_{\cdot},[\phitilde_{\cdot}]\bigr)$
and the morphism 
$\Ob'_{\rel}\bigl(V_{\cdot},[\phitilde]\bigr)
\lrarr L_{M[L]/U_{G_m}}$.

Let $M(L)$ be as in the subsubsection
\ref{subsubsection;06.5.4.50}.
It is easy to observe that
$M(L)$ is isomorphic to $\nbigi^{\ast}$.
We have the smooth projection
$\pi_2:M(L)\lrarr M[L]$.
The pull back of $[\phitilde]$ via
$\pi_{2\,X}$ is denoted by $\phitilde$.

\begin{lem}
We have the following commutative diagram:
\begin{equation}
 \label{eq;06.5.9.12}
 \begin{CD}
 \pi_2^{\ast}\Ob_{\rel}'(V,[\phitilde])
 @>{\simeq}>>
 \Ob_{\rel}(V,\phitilde)\\
 @VVV @VVV \\
 \pi_2^{\ast}L_{M[L]/U_{G_m}}
 @>{\simeq}>> L_{M(L)/U}
 \end{CD}
\end{equation}
\end{lem}
\pf
We have the following commutative diagram:
\[
 \begin{CD}
 M(L)\times X @>{\Phi(V_{\cdot},\phitilde)}>>
 Y(W_{\cdot},P_{\cdot})
 @>{\psitilde}>> Y(W_{\cdot}) @>>> X\\
 @VVV @VVV @VVV @VVV \\
 M[L]\times X @>{\Phi(V_{\cdot},[\phitilde])}>>
 Y(W_{\cdot},[P_{\cdot}])
 @>{\psi}>> Y(W_{\cdot})_{G_m}
 @>>> X_{G_m}
 \end{CD}
\]
Therefore, we obtain the following:
\[
\begin{CD}
 L_{M(L)\times X/U\times X}@<<<
 \Phi(\phi)^{\ast}L_{Y(W_{\cdot},P_{\cdot})/Y(W_{\cdot})}
 @<<< \gminig_{\rel}(V_{\cdot},\phitilde)\\
 @A{\simeq}AA @A{\simeq}AA @A{\simeq}AA\\
 \pi_{2,X}^{\ast}L_{M[L]\times X/U\times X_{G_m}}
 @<<<
 \Phi([\phi])^{\ast}L_{Y(W_{\cdot},[P_{\cdot}])/Y(W_{\cdot})_{G_m}}
 @<<<
 \pi_{2,X}^{\ast}\gminig'_{\rel}(V_{\cdot},[\phitilde])
\end{CD}
\]
Then the claim is clear.
\hfill\qed

\vspace{.1in}
We have the morphisms
$\Ob_{\rel}(V_{\cdot},\phitilde)\lrarr
 L_{M(L)/U}\lrarr L_{M(L)/M[L]}$,
where the latter is the projection.

\begin{lem}
 \label{lem;06.5.9.100}
We have the commutative diagram:
\begin{equation}
 \label{eq;06.5.9.10}
 \begin{CD}
 \pi_2^{\ast}\Ob_{\rel}\bigl(V_{\cdot},[\phitilde]\bigr)
 @>>> 
 \Ob_{\rel}(V_{\cdot},\phitilde) 
 @>>> 
 L_{M(L)/M[L]}
 @>>> 
 \pi_2^{\ast}\Ob_{\rel}(V_{\cdot},[\phitilde])[1]\\
 @VVV @VVV @V{\simeq}VV @VVV \\
 \pi_2^{\ast}L_{M[L]/U}
 @>>> 
 L_{M(L)/U}
 @>>> 
 L_{M(L)/M[L]}
 @>>> 
 \pi_2^{\ast}L_{M[L]/U}[1].
 \end{CD} 
\end{equation}
In particular,
$\pi_2^{\ast}\Ob_{\rel}(V_{\cdot},[\phitilde])
\simeq\Cone\bigl(
 \Ob_{\rel}(V_{\cdot},[\phitilde])
\lrarr L_{M(L)/M[L]} \bigr)$.
\end{lem}
\pf
Let $F_{3,\nbigi}$ denote the morphism
$M[L]\lrarr U_{G_m}$, induced by $\nbigi$.
Then we have the following isomorphism
of the distinguished triangles on $M(L)$:
\begin{equation}
 \label{eq;06.5.9.11}
 \begin{CD}
 \pi_2^{\ast}L_{M[L]/U}
 @>>>
 \pi_2^{\ast}L_{M[L]/U_{G_m}}
 @>>>
 \pi_2^{\ast}F_{3,\nbigi}^{\ast}L_{U_{G_m}/U}[1]
 @>>>
 \pi_2^{\ast}L_{M[L]/U}[1]\\
 @V{=}VV @V{\simeq}VV @V{\simeq}VV @V{=}VV\\
 \pi_2^{\ast}L_{M[L]/U}
 @>>>
 L_{M(L)/U}
 @>>>
 L_{M(L)/M[L]}
 @>>>
 \pi_2^{\ast}L_{M[L]/U}[1]
 \end{CD}
\end{equation}
We obtain (\ref{eq;06.5.9.10})
from (\ref{eq;06.5.9.12}) and (\ref{eq;06.5.9.11}).
\hfill\qed

\begin{lem}
\label{lem;06.5.3.41}
$\ob_{\rel}\bigl(V_{\cdot},[\phitilde_{\cdot}]\bigr)$
gives a relative obstruction theory
for $M[L]$ over $U$.
\end{lem}
\pf
The complex $L_{M(L)/M[L]}$
is quasi isomorphic to
the $0$-th cohomology sheaf,
and we know that 
$\ob_{\rel}(V_{\cdot},\phitilde)$
gives the obstruction theory
for $M(L)$ over $U$
(Lemma \ref{lem;06.5.3.2}).
Then, the claim of the lemma follows
from the diagram (\ref{eq;06.5.9.10}).
\hfill\qed

\subsubsection{Preparation for the obstruction theory
 of Master space}
\label{subsubsection;06.5.4.150}

We will use the notation in the subsubsections
\ref{subsubsection;06.5.4.110}
and \ref{subsubsection;06.5.4.100}.
We have the weight $(-1)$-action
of $G_m$ on $\nbigo_X$.
It induces the $G_m$-action
on $A(W_{\cdot},P_{\cdot})$.
The quotient stack is denoted by
$A(W_{\cdot},[P_{\cdot}])$.
We have the induced morphism
$\Gamma_{[L]}:Y(W_{\cdot},[P_{\cdot}])
\lrarr A(W_{\cdot},[P_{\cdot}])$.
From $(F_{3,X}^{\ast}V_{\cdot},[\phitilde])$ on
$U_3\times X$ and $[\iotatilde]$,
we obtain the morphism
$[\phibar]:p_{U_3}^{\ast}\nbigo_X\otimes p_{X}^{\ast}M
\lrarr F_{3,X}^{\ast}V_0$.
We put
$\gminih'_{\rel}(V_{\cdot},[\phitilde]):=
 \nhom\bigl(
 p_{U_3}^{\ast}\nbigo_X\otimes p_X^{\ast}M,
 F_{3X}^{\ast}V_{\cdot}
 \bigr)^{\lor}$.
\index{$\gminih'_{\rel}(V_{\cdot},[\phitilde])$}
We have the induced map
$\gamma([\phibar]):
\gminih'_{\rel}(V_{\cdot},[\phitilde])[-1]
\lrarr
\gminih(V_{\cdot})$,
and 
$\gamma([\phibar])_{\leq 1}:
\gminih'_{\rel}(V_{\cdot},[\phitilde])[-1]_{\leq 1}
\lrarr
\gminih(V_{\cdot})_{\leq 1}$.

We put 
$\Psi(V_{\cdot},[\phitilde_{\cdot}]):=
 \Gamma_{[L]}\circ\Phi(V_{\cdot},[\phitilde])$,
which give the classifying map of the tuple
$(F_{X,D}^{\ast}V_0,[\phibar])$.
Then, 
$\Psi(V_{\cdot},[\phitilde_{\cdot}])^{\ast}
 L_{A(W_{\cdot},[P_{\cdot}])/X_{G_m}}$
is expressed by the complex
$\cone\bigl(\gamma([\phibar])_{\leq 1}\bigr)$,
and the naturally defined morphism
$\cone\bigl(\gamma([\phibar])_{\leq 1}\bigr)
\lrarr
 \cone\bigl(\gamma([\phitilde_{\cdot}])_{\leq 1}\bigr)$
expresses
$\Psi(V_{\cdot},[\phitilde_{\cdot}])^{\ast}
 L_{A(W_{\cdot},[P_{\cdot}])/X_{G_m}}
\lrarr 
 \Phi(V_{\cdot},[\phitilde_{\cdot}])^{\ast}
 L_{Y(W_{\cdot},[P_{\cdot}])/X_{G_m}}$.

We put 
$\gminih'(V_{\cdot},\phitilde):=
\cone\bigl(\gamma([\phibar])\bigr)$.
Then, we obtain the following commutative diagram:
\[
 \begin{CD}
 \gminih'(V_{\cdot},[\phitilde]) @>>> 
 \gminig'(V_{\cdot},[\phitilde]) @>>>
 L_{U_3\times X/X_{G_m}}\\
 @AAA @AAA @AAA \\
 \gminih(V_{\cdot}) @>>>
 \gminig(V_{\cdot}) @>>>
 F_{3\,X}^{\ast}L_{U\times X/X}
 \end{CD}
\]
We put 
$ \Ob^{\prime\,G}_{\rel}(V_{\cdot},[\phitilde]):=
 Rp_{X\,\ast}\bigl(
 \gminih'_{\rel}(V_{\cdot},[\phitilde])\otimes\omega_X 
 \bigr)$. 
\index{$\Ob^{\prime\,G}_{\rel}(V_{\cdot},[\phitilde])$}
Then, we have the following diagram on $U_3$:
\begin{equation}
 \begin{CD}
 L_{U_3/U_{G_m}}[-1] @<<<
 \Ob'_{\rel}(V_{\cdot},[\phitilde])[-1] @<<<
 \Ob^{\prime\,G}_{\rel}(V_{\cdot},[\phitilde])[-1]\\
 @VVV @VVV @VVV \\
 F_{3}^{\ast}L_{U/k}@<<<
 \Ob(V_{\cdot}) @<<< 
 \Ob^G(V_{\cdot})
 \end{CD}
\end{equation}
We obtain the composite
$\Ob^{\prime\,G}_{\rel}(V,[\phitilde])
\lrarr L_{U_3/U_{G_m}}
\lrarr F_{3,M}^{\ast}L_{U_{G_m}/U}[1]$.
Let $\Ob^{G}_{\rel}(V_{\cdot},[\phitilde])$
denote the cone of the morphism.
Then, we obtain the following diagram:
\begin{equation}
 \label{eq;06.5.4.300}
 \begin{CD}
 L_{U_3/U}[-1] @<<<
 \Ob_{\rel}(V_{\cdot},[\phitilde])[-1] @<<< 
 \Ob^G_{\rel}(V_{\cdot},[\phitilde])[-1]\\
 @VVV @VVV @VVV \\
 F_{3}^{\ast}L_{U/k}@<<<
 \Ob(V_{\cdot}) @<<< 
 \Ob^G(V_{\cdot})
 \end{CD}
\end{equation}

Now, we assume the condition (C) in the subsubsection
\ref{subsubsection;06.5.4.100}.
We have the weight $(-1)$-action of $G_m$
on the one dimensional vector space $k$.
It induces the $G_m$-action
on $B(W_{\cdot},P_{\cdot})$.
The quotient stack is denoted by
$B(W_{\cdot},[P_{\cdot}])$.
\index{$B(W_{\cdot},[P_{\cdot}])$}
Similarly,
we obtain $A(W_{\cdot},[P_{\cdot}])$
from $A(W_{\cdot},P_{\cdot})$.
We have the natural isomorphism
$J_{[L]}:B(W_{\cdot},[P_{\cdot}])\times X
 \lrarr A(W_{\cdot},[P_{\cdot}])$.
From 
$[\phibar]:p_X^{\ast}M\otimes 
 p_{U_3}^{\ast}\nbigo_X\lrarr V_0$,
we obtain the morphism
$\Xi(V_{\cdot},[\phibar]):
 U_3\lrarr B(W_{\cdot},[P_{\cdot}])$.
Note that the composite
of $\Xi(V_{\cdot},[\phibar])_X$
and $J_{[L]}$ is same as
$\Psi(V_{\cdot},[\phibar])$.
Therefore,
we have the following commutative diagram:
\[
 \begin{CD}
 U_3\times X @>>>
 B(W_{\cdot},[P_{\cdot}])\times X
 @>>>
 A(W_{\cdot},[P_{\cdot}]) \\
 @VVV @VVV @VVV \\
 (U\times X)_{G_m} @>>> 
\bigl(B(W_{\cdot})\times X\bigr)_{G_m} @>>> 
 A(W_{\cdot})_{G_m}
 \end{CD}
\]
Therefore, we obtain the following commutative
diagram:
\[
 \begin{CD}
 F_{3,M}^{\ast}L_{U_{G_m}/k_{G_m}} @<<<
 \Phi(V_{\cdot})_M^{\ast}L_{B(W_{\cdot})_{G_m}/k_{G_m}} @<<<
 \Ob^G(V_{\cdot}) \\
 @AAA @AAA @AAA \\
 L_{U_3/U_{G_m}}[-1] @<<<
 \Xi(V_{\cdot},[\phibar])^{\ast}
 L_{B(W_{\cdot},[P_{\cdot}])/B(W_{\cdot})_{G_m}}[-1]
 @<{\tau_3'}<<
 \Ob^{\prime\,G}_{\rel}(V_{\cdot},[\phitilde])
 \end{CD}
\]
The following lemma is similar to
Lemma \ref{lem;06.5.2.10}.
\begin{lem}
The morphism $\tau_3'$ is isomorphic.
\hfill\qed
\end{lem}

By the standard modification,
we obtain the following commutative diagram:
\begin{equation}
 \label{eq;06.5.4.250}
 \begin{CD}
 F_{3}^{\ast}L_{U/k} @<<<
 \Phi(V_{\cdot})^{\ast}L_{B(W_{\cdot})/k} @<<<
 \Ob^G(V_{\cdot}) \\
 @AAA @AAA @AAA \\
 L_{U_3/U}[-1] @<<<
 \Xi(V_{\cdot},[\phibar])^{\ast}
 L_{B(W_{\cdot},[P_{\cdot}])/B(W_{\cdot})}[-1]
 @<{\tau_3}<<
 \Ob^{G}_{\rel}(V_{\cdot},[\phitilde])[-1]
 \end{CD}
\end{equation}
\index{$\Ob^{G}_{\rel}(V_{\cdot},[\phitilde])$}
The following lemma immediately follows from
the previous lemma.
\begin{lem}
\label{lem;06.5.4.200}
The morphism $\tau_3$ is isomorphic.
\hfill\qed
\end{lem}

\subsubsection{Preparation for Proposition \ref{prop;06.5.10.50}}

We use the notation in the subsubsection 
\ref{subsubsection;06.5.4.110}.
As a preparation for the proof of Proposition \ref{prop;06.5.10.50},
let us see the morphism
$\varphi_1:\Ob'_{\rel}(V_{\cdot},[\phitilde])
\lrarr F_{3}^{\ast}L_{U_{G_m}/U}$
on $U_3$, more closely.
We recall
$ \Phi(V,[\phitilde])^{\ast}
 L_{Y(W_{\cdot},[P_{\cdot}])/Y(W_{\cdot})_{G_m}}
\simeq
 \gminig'_{\rel}(V_{\cdot},[\phitilde])_{\leq 0}$.

\begin{lem}
\label{lem;06.5.10.30}
We have the following commutative diagram:
\[
 \begin{CD}
 \Phi(V,[\phitilde])^{\ast}
 L_{Y(W_{\cdot},[P_{\cdot}])/Y(W_{\cdot})_{G_m}}
@>>>
 \Phi(F_{3,X}^{\ast}V_{\cdot})_{M}^{\ast}
 L_{Y(W_{\cdot})_{G_m}/Y(W_{\cdot})}[1]\\
 @AAA @AA{\simeq}A\\
 \gminig_{\rel}'(V_{\cdot},[\phitilde])
 @>{\gminif'}>>
 \nbigo
 \end{CD}
\]
Here, the morphism $\gminif'$ is given as follows:
\[
 \nhom\bigl(F_{3,X}^{\ast}V_{0},\,\,
p_U^{\ast}P_{0}\otimes p_X^{\ast}M\bigr)
\oplus
 \nhom\bigl(F_{3,X}^{\ast}V_{0},\,\,
p_U^{\ast}P_{0}\otimes p_X^{\ast}M\bigr)
\lrarr
 \nbigo,
\quad
 (a_0,a_{-1})
\longmapsto
 \tr\bigl([\phitilde_0]\circ a_0\bigr)
+\tr\bigl([\phitilde_{-1}]\circ a_{-1}\bigr).
\]
\end{lem}
\pf
It follows from Corollary \ref{cor;06.5.9.2}.
\hfill\qed

\vspace{.1in}
The morphism $\gminif'$ induces
the following morphisms:
\[
 \Ob'_{\rel}(V_{\cdot},[\phitilde])=
 Rp_{X\,\ast}\bigl(
 \gminig_{\rel}(V_{\cdot},[\phitilde])\otimes\omega_X
\bigr)
 \lrarr
 Rp_{X\,\ast}\bigl(p_{U_3}^{\ast}\omega_X\bigr)
\lrarr \nbigo_{U_3}.
\]
It is easy to check that
the composite is same as
$\varphi_1$.

\subsection{Relative Obstruction Theory for Parabolic Structure}
\label{subsection;06.7.3.35}
\subsubsection{Construction of the complex}
\label{subsubsection;06.5.5.5}

We use the notation in the subsubsection
\ref{subsubsection;06.5.2.3}.
Let $F_4:U_4\lrarr U$ be a morphism,
and let $F_{\ast}$ be a quasi-parabolic structure
of $F_{4\,X}^{\ast}E$ at $D$.
We denote the kernel of
$F_{4,X}^{\ast}V_{0\,|\,D}\lrarr 
\Cok_{h-1}(E)$ by $V_D^{(h)}$.
Since the smoothness of $D$ is assumed,
$V_D^{(h)}$ are locally free.
The filtered vector bundle
$V^{(1)}_D\supset V^{(2)}_D\supset
\cdots\supset V_D^{(l+1)}$
is denoted by $V_D^{\ast}$.

We put $W^{(1)}:=W_{0}$
and $W^{(l+1)}:=W_{-1}$.
Let $W^{(h)}$ $(h=2,\ldots,l)$
be vector spaces over $k$
such that $\dim W^{(h)}=\rank V_D^{(h)}$.
Let $D$ be a smooth divisor of $X$.
We denote $W^{(h)}\otimes\nbigo_D$
and $W_i\otimes\nbigo_D$
by $W_{i\,D}$ and $W^{(h)}_D$ respectively.
We have the natural right
$\GL(W_{\cdot})$-action
on $N(W_{-1,D},W_{0,D})$ given by
$(g_0,g_{-1})\cdot f=g_0^{-1}\circ f\circ g_{-1}$.
The quotient stack is denoted by
$Y_D(W_{\cdot})$.
\index{$Y_D(W_{\cdot})$}
Similarly, we have the natural right
$\GL(W^{(h)})
\times\GL(W^{(h+1)})$-action
on $N(W_{D}^{(h+1)},W_{D}^{(h)})$
given by
$(g^{(h)},g^{(h+1)})\cdot f
=g^{(h)\,-1}\circ f\circ g^{(h+1)}$.
Thus we obtain the right
$\prod_{h=1}^{l+1} \GL(W^{(h)})$-action
on $\prod_{h=1}^{l}
 N\bigl(W_D^{(h+1)},W_D^{(h)}\bigr)$,
where the latter fiber product is taken
over $D$.
The quotient stack is denoted by
$Y_D(W_{\cdot},W^{\ast})$.
\index{$Y_D(W_{\cdot},W^{\ast})$}
The composition of the morphisms induce the map
$\prod_{h=1}^{l}
  N\bigl(W_D^{(h+1)}, W_D^{(h)}\bigr)
\lrarr N(W_{-1\,D},W_{0\,D})$.
It induce the morphism
$Y_D(W_{\cdot},W^{\ast})
\lrarr Y_D(W_{\cdot})$.

We have the classifying map
$\Phi_D(V_{\cdot},F_{\ast}):
 U_4\times D\lrarr Y_D(W_{\cdot},W^{\ast})$
over $D$
obtained from the tuple $V_D^{\ast}$.
We also have the classifying map
$\Phi(V_{\cdot|D}):
 U\times D\lrarr Y_D(W_{\cdot})$.
Thus, we obtain the following diagram on $U_4\times D$:
\[
 \begin{CD}
 L_{U_4\times D/D} @<<<
 \Phi_D(V_{\cdot},F_{\ast})^{\ast}
 L_{Y_D(W_{\cdot},W^{\ast})/D} \\
 @AAA @AAA \\
 F_{4,D}^{\ast}L_{U\times D/D} @<<<
 \Phi\bigl(F_{4\,D}^{\ast}(V_{\cdot|D})\bigr)^{\ast}
 L_{Y_D(W_{\cdot})/D}.
 \end{CD}
\]

We use the notation in the subsubsection
\ref{subsubsection;06.5.23.10}.
We put $\gminig_D(V_{\cdot},F_{\ast}):=
 C_1(V_D^{\ast},V_D^{\ast})^{\lor}[-1]$
and
$\gminig_{\rel}(V_{\cdot},F_{\ast}):=
 C_2(V_{D}^{\ast},V_D^{\ast})^{\lor}[-1]$.
\index{$\gminig_D(V_{\cdot},F_{\ast})$,
 $\gminig_{\rel}(V_{\cdot},F_{\ast})$}
We also put
$\gminig(V_{\cdot|D}):=
 \nhom(V_{\cdot|D},V_{\cdot|D})^{\lor}[-1]$.
We have the morphism
$\gamma_D:\gminig(V_{\cdot|D})\lrarr
 \gminig_D(V_{\cdot},F_{\ast})$
induced by $\varphi$ given in \ref{eq;06.6.14.1}.
It is easy to see that
$\gminig_{\rel}(V_{\cdot},F_{\ast})$
is quasi isomorphic to $\cone(\gamma_D)$.

By an argument explained in the subsubsection
\ref{subsubsection;06.4.29.15},
we can show that
$\Phi_D(V_{\cdot},F_{\ast})^{\ast}
 L_{Y_D(W_{\cdot},W^{\ast})/D}$
and $\Phi(V_{\cdot|D})^{\ast}L_{Y_D(W_{\cdot})/D}$
are expressed by
$\gminig_D(V_{\cdot},F_{\ast})_{\leq 1}$
and $\gminig(V_{\cdot|D})_{\leq 1}$.
Under the identification,
the natural morphism
$\Phi_D(V_{\cdot},F_{\ast})^{\ast}
 L_{Y_D(W_{\cdot},W^{\ast})/D}
\lrarr \Phi(V_{\cdot\,|\,D})^{\ast}L_{Y_D(W_{\cdot})/D}$
is given by the restriction $\gamma_{D\leq 1}$.
Then we obtain the following commutative diagram:
\[
 \begin{CD}
 \gminig(V_{\cdot|D}) @>>>
 \gminig(V_{\cdot|D})_{\leq 1}
 @>>>
 F_{4}^{\ast}L_{U/k}\\
 @VVV @VVV @VVV \\
 \gminig_D(V_{\cdot},F_{\ast}) @>>>
 \gminig_D(V_{\cdot},F_{\ast})_{\leq 1}
 @>>>
 L_{U_4/k}
 \end{CD}
\]
Then we obtain the following morphism
of the distinguished triangles
on $U_4$:
\[
 \begin{CD}
 \gminig(V_{\cdot|D})
@>>>
  \gminig_D(V_{\cdot},F_{\ast})
 @>>>
 \gminig_{\rel}(V_{\cdot|D},F_{\ast})
 @>>>
 \gminig(V_{\cdot|D})[1]\\
 @VVV @VVV @VVV @VVV \\
F_{4\,X}^{\ast}L_{U\times X/X} @>>>
 L_{U_4\times X/X} @>>>
 L_{U_4\times X/ U\times X} @>>>
F_{4\,X}^{\ast} L_{U\times X/X}
 \end{CD}
\]

Let $\omega_D$ denote the dualizing complex of $D$.
We put $\Ob_{\rel}(V_{\cdot},F_{\ast}):=
Rp_{D\,\ast}\bigl(
 \gminig_{\rel}\bigl(V_{\cdot},F_{\ast}\bigr)
 \otimes\omega_D\bigr)$ 
\index{$\Ob_{\rel}(V_{\cdot},F_{\ast})$}
and $\Ob(V_{\cdot\,|\,D}):=
 Rp_{D\,\ast}\bigl(\gminig(V_{\cdot\,|\,D})\otimes\omega_D\bigr)$.
Then, we obtain the following commutative diagram:
\begin{equation}
 \label{eq;6.14.30}
\begin{CD}
\Ob_{\rel}\bigl(V_{\cdot},F_{\ast}\bigr)[-1] 
@>>> L_{U_4/U}[-1]\\
@VVV @VVV \\
\Ob(V_{\cdot|D})
@>{\ob(V_{\cdot|D})}>> F_{4}^{\ast}L_{U/k}.
\end{CD}
\end{equation}

We have the following exact sequence
of complexes:
\[
 0\lrarr \gminig(V_{\cdot})\otimes\omega
 \lrarr \gminig(V_{\cdot})\otimes\omega(D)
 \lrarr \gminig(V_{\cdot|D})\otimes\omega_D[1]
 \lrarr 0.
\]
Thus we obtain the morphism
$Rp_{D\,\ast}\bigl(\gminig(V_{\cdot|D})
 \otimes\omega_D\bigr)
 \lrarr Rp_{X\,\ast
 }\bigl(\gminig(V_{\cdot})\otimes\omega_X\bigr)$,
that is,
$\eta:\Ob(V_{\cdot|D})
 \lrarr \Ob(V_{\cdot})$.
We have the morphism
$\ob(V_{\cdot}):\Ob(V_{\cdot}) \lrarr L_U$,
and hence the composite
$\ob(V_{\cdot})\circ \eta:
 \Ob(V_{\cdot|D})\lrarr L_U$.
On the other hand,
we have the morphism $\ob(V_{\cdot|D})$
in the diagram (\ref{eq;6.14.30}).

\begin{lem}
\label{lem;06.5.13.10}
Two morphisms $\ob(V_{\cdot})\circ \eta$ and 
$\ob(V_{\cdot\,|\,D})$ are same in the derived category.
\end{lem}
\pf
It is easy to observe that
$\gminig(V_{\cdot|D})\lrarr 
\Phi(V_{\cdot|D})^{\ast}L_{U\times D/D}$
 is the restriction of
$\gminig(V_{\cdot})\lrarr 
 \Phi(V_{\cdot})^{\ast}L_{U\times X/X}$
to $D$.
Then the coincidence
$\ob(V_{\cdot})\circ\eta=
\ob(V_{D\,\cdot})$ follows from
the compatibility of of the traces
for $\omega_X$ and $\omega_D$.
(See \cite{altman-kleiman}, for example).
\hfill\qed

\vspace{.1in}
Thus, we obtain the following commutative diagram:
\begin{equation}
 \label{eq;6.14.31}
\begin{CD}
\Ob_{\rel}\bigl(V_{\cdot},F_{\ast}\bigr)[-1]
@>>>
 \Ob(V_{\cdot|D}) 
@>>>
\Ob(V_{\cdot})\\
@V{\ob_{\rel}(V_{\cdot},F_{\ast})}VV @VVV @VVV\\ 
L_{U_4/U}[-1]
@>>>
 F_4^{\ast}L_{U/k}
@>>> 
F_4^{\ast}L_{U/k}.
\end{CD}
\end{equation}

The following lemma can be shown by an argument 
to use the filtered objects as in 
the proof of Lemma \ref{lem;6.14.6}.
\begin{lem} 
 \label{lem;6.16.30}
The diagram {\rm(\ref{eq;6.14.31})}
depends only on $\bigl(E,F_{\ast}(E)\bigr)$.
\hfill\qed
\end{lem}

\subsubsection{Relative obstruction property}

For any $U$-scheme $g:T\lrarr U$,
let $F(T)$ denote the set of 
the parabolic structure of $g_X^{\ast}E$
of a fixed type.
Thus, we obtain the functor of
the category of $U$-schemes
to the category of sets.
It is easy to see that
the functor $F$ is representable
by a scheme $M$
via the method of the quot schemes.

Let $\pi:M\lrarr U$ denote the projection.
On $M\times X$,
we have the universal parabolic structure $F^u_{\ast}$
of $\pi_X^{\ast}E$.
From the resolution $V_{\cdot}$
and $F^u_{\ast}$,
we obtain the complex
$\Ob_{\rel}(V_{\cdot},F^u_{\ast})$
and the morphism
$\ob_{\rel}(V_{\cdot},F^u_{\ast}):
 \Ob_{\rel}(V_{\cdot},F^u_{\ast})
\lrarr
 L_{M/U}$.

\begin{lem}
\label{lem;06.5.3.40}
$\ob_{\rel}(V_{\cdot},F_{\ast}^u)$
gives an obstruction theory of $M$ over $U$.
\end{lem}
\pf
It follows from Lemma \ref{lem;06.5.3.15}.
Note that $\gminig_{\rel}(V_{\cdot},F^u_{\ast})$
is naturally isomorphic as the complex considered there.
\hfill\qed

\subsubsection{The decomposition into the trace-free part
 and the diagonal part}
\label{subsubsection;06.5.11.200}

\begin{lem}
The morphism
$\Ob_{\rel}(V,F_{\ast})[-1]\lrarr \Ob(V_{\cdot})$
factors through 
the trace free part $\Ob^{\circ}(V_{\cdot})$.
\end{lem}
\pf
We have the trace map 
$\tr:\gminig_D(V_{\cdot},F_{\ast})\lrarr\nbigo_D[-1]$
given as follows:
\[
 \bigoplus_{i=1}^l\nhom(V_D^{(i)},V_{D}^{(i)})
\lrarr
 \nbigo,
\quad
 (f_i)\longmapsto
 \sum \tr(f_i)
\]
We also have the map $i:\nbigo_D[-1]\lrarr\gminig_D(V,F_{\ast})$
given as follows:
\[
 \nbigo\lrarr \bigoplus_{i=1}^{l}\nhom(V_D^{(i)},V_D^{(i)}),
\quad
 t\longmapsto (t\cdot \id_{V^{(1)}},0,\ldots,0,-t\cdot \id_{V^{(l+1)}}).
\]
We put $\gminig_D^{\circ}(V_{\cdot},F_{\ast}):=\Ker(\tr)$
and $\gminig_D^{d}(V_{\cdot},F_{\ast}):=\Image(i)$.
\index{$\gminig_D^{\circ}(V_{\cdot},F_{\ast})$,
 $\gminig_D^d(V_{\cdot},F_{\ast})$}
We also have the decomposition
$\gminig(V_{\cdot|D})=
 \gminig^{\circ}(V_{\cdot|D})\oplus
 \gminig^{d}(V_{\cdot|D})$
as in the subsubsection \ref{subsubsection;06.5.9.20}.
We have the following commutative diagrams:
\[
 \begin{CD}
 \nbigo[-1] @>>> \nbigo[-1]\\
 @V{i}VV @V{i}VV \\
 \gminig(V_{\cdot|D})
 @>>>
 \gminig_D(V_{\cdot},F_{\ast})
 \end{CD}
\quad\quad\quad\quad\quad\quad
 \begin{CD}
 \gminig(V_{\cdot|D})
 @>>>
 \gminig_D(V_{\cdot},F_{\ast}) \\
 @V{\tr}VV @V{\tr}VV \\
 \nbigo[-1] @>>> \nbigo[-1]
 \end{CD}
\]
Therefore, the decomposition is compatible
with $\gminig(V_{\cdot|D})\lrarr\gminig_D(V_{\cdot},F_{\ast})$.
It follows from 
$\gminig_{\rel}(V_{\cdot},F_{\ast})[-1]
\lrarr
 \gminig(V_{\cdot|D})$ factors through
$\gminig^{\circ}(V_{\cdot|D})$.
Therefore,
$\Ob_{\rel}(V,F_{\ast})[-1]
\lrarr \Ob(V_{\cdot|D})$ factors through
$\Ob^{\circ}(V_{\cdot|D})$.
It is easy to see that the morphism
$\Ob(V_{\cdot|D})\lrarr \Ob(V_{\cdot})$
is compatible with the decomposition
into the trace-free part and the diagonal part.
Thus we are done.
\hfill\qed

\vspace{.1in}
As an immediate corollary,
we obtain the decomposition of the cone:
\begin{equation}
 \label{eq;06.5.9.25}
 \Cone\Bigl(
 \Ob_{\rel}(V_{\cdot},F_{\ast})
\lrarr
 \Ob(V_{\cdot})
\Bigr)
\simeq
 \Cone\Bigl(
 \Ob_{\rel}(V_{\cdot},F_{\ast})
\lrarr \Ob^{\circ}(V_{\cdot})
 \Bigr)
\oplus
 \Ob_{\rel}^d(V_{\cdot})
\end{equation}

\subsection{Obstruction Theory for Moduli Stacks
 of Stable Objects}
\label{subsection;06.6.13.5}
\subsubsection{Relative complexes}
\label{subsubsection;06.5.22.5}

Let $y$ be an element of $H^{\ast}(X)$.
Let $\vecy$ be an element of $\Type$
whose $H^{\ast}(X)$-component is $y$.
Let $\nbigm(m,\vecy)$ be the open subset
of $\nbigm(\vecy)$ determined by
the condition $O_m$.
We have the natural morphism
$\gminip_1:\nbigm(m,\vecy)\lrarr\nbigm(m,y)$.
On $\nbigm(m,\vecy)\times X$,
we have the universal sheaf
$\gminip_{1\,X}^{\ast}\nbige^u$ 
over $\nbigm(m,\vecy)\times X$
with the parabolic structure
$F_{\ast}^u$ at $D$.
From the resolution
$\nbigv_{\cdot}$ of $\nbige^u(m)$
and the parabolic structure $F$,
we obtain the complex
$\Ob_{\rel}(\nbigv_{\cdot},F_{\ast})$
and the morphism:
\[
 \ob_{\rel}(m,\vecy):
 \Ob_{\rel}(m,\vecy)
\lrarr
 L_{\nbigm(m,\vecy)/\nbigm(m,y)}
\]
\index{$\Ob_{\rel}(m,\vecy)$, $\ob_{\rel}(m,\vecy)$}

\vspace{.1in}

Let $\nbigm(m,y,L)$ denote the open subset
of $\nbigm(y,L)$ determined by the condition $O_m$.
The natural morphism
$\nbigm(m,y,L)\lrarr\nbigm(m,y)$ is denoted by
$\gminip_2$.
We have the universal $L$-section
$\phi^u$ of $\gminip_{2\,X}^{\ast}\nbige^u$.
It induces the $L(m)$-section
of $\gminip_{2\,X}^{\ast}\nbige^u(m)$,
which is also denoted by $\phi^u$.
We fix an inclusion $\iota:\nbigo(-m)\lrarr L$.
If $m$ is sufficiently large,
we may assume that $L(m)$ has a locally free resolution
$P_{\cdot}=(P_{-1}\rarr P_0)$
such that $P_0$ is a direct sum of some $\nbigo_X$.
Since we have the isomorphism
$p_{X\,\ast}\bigl(\nbige^u(m)\bigr)
 \simeq p_{X\,\ast}(\nbigv_0)$,
the $L(m)$-section $\phi^u$ is canonically lifted to
the morphism 
$\phitilde^u_{\cdot}:
 p_{\nbigm(m,y,L)}^{\ast}P_{\cdot}\lrarr 
 \gminip_{2\,X}^{\ast}\nbigv_{\cdot}$.
Then we obtain the morphism:
\[
 \ob_{\rel}(m,y,L):
 \Ob_{\rel}(m,y,L)\lrarr
 L_{\nbigm(m,y,L)/\nbigm(m,y)}
\]
\index{$\Ob_{\rel}(m,y,L)$, $\ob_{\rel}(m,y,L)$}

\vspace{.1in}

Let $\nbigm(m,y,[L])$ denote the open subset
of $\nbigm(y,[L])$ determined by the condition $O_m$.
The natural morphism
$\nbigm(m,y,[L])\lrarr\nbigm(m,y)$
is denoted by $\gminip_3$.
Then we have the universal reduced $[L]$-section
$[\phi^u]$ of $\gminip_{3\,X}^{\ast}\nbige^u$
over $\nbigm(m,y,[L])\times X$.
As before,
we obtain the morphism:
\[
 \ob_{\rel}\bigl(m,y,[L]\bigr):
 \Ob_{\rel}\bigl(m,y,[L]\bigr)
\lrarr
 L_{\nbigm(m,y,[L])/\nbigm(m,y)}.
\]
\index{$\Ob_{\rel}(m,y,[L])$, $\ob_{\rel}(m,y,[L])$}

Let $\nbigm(m,\yhat)$ denote the open subset
of $\nbigm(\yhat)$ determined by
the condition $O_m$.
The natural morphism
$\nbigm(m,\yhat)\lrarr\nbigm(m,y)$
is denoted by $\gminip_4$.
Then we have the universal orientation
$\rho^u$ of $\gminip_{4\,X}^{\ast}\nbige^u$
over $\nbigm(m,\yhat)\times X$.
From the resolution
$\nbigv_{\cdot}$  and the orientation $\rho^u$,
we obtain the morphism:
\[
 \ob_{\rel}(m,\yhat):
 \Ob_{\rel}(m,\yhat)
\lrarr
 L_{\nbigm(m,\yhat)/\nbigm(m,y)}.
\]
\index{$\Ob_{\rel}(m,\yhat)$, $\ob_{\rel}(m,\yhat)$}

\subsubsection{Construction of the complexes and the morphisms}
\label{subsubsection;06.5.4.450}

We have the naturally defined morphisms of
$\nbigm(m,\vecyhat,[L])$
to $\nbigm(m,y)$, $\nbigm(m,\vecy)$,
$\nbigm(m,y,[L])$
and $\nbigm(m,\yhat)$.
The pull back of the complexes
$\Ob(m,y)$,
$\Ob_{\rel}\bigl(m,\vecy\bigr)$,
$\Ob_{\rel}\bigl(m,y,[L]\bigr)$
and 
$\Ob_{\rel}\bigl(m,\yhat\bigr)$
are denoted by the same notation.
Then, we put as follows
on $\nbigm(m,\vecyhat,[L])$:
\[
  \Ob(m,\vecyhat,[L]):=
 \cone\Bigl(
  \Ob_{\rel}\bigl(m,\yhat\bigr)[-1]
\oplus
 \Ob_{\rel}\bigl(m,y,[L]\bigr)[-1]
\oplus
 \Ob_{\rel}\bigl(m,\vecy\bigr)[-1]
\lrarr
 \Ob\bigl(m,y\bigr)
 \Bigr)
\]
\index{$\Ob(m,\vecyhat,[L])$, $\ob(m,\vecyhat,[L])$}

\begin{prop}
We have the naturally defined morphism
$\ob(m,\vecyhat,[L]):
 \Ob(m,\vecyhat,[L])
\lrarr L_{\nbigm(m,\vecyhat,[L])}$.
\end{prop}
\pf
We put
$C:=\Ob_{\rel}\bigl(m,\yhat\bigr)
\oplus
 \Ob_{\rel}\bigl(m,y,[L]\bigr)
\oplus
 \Ob_{\rel}\bigl(m,\vecy\bigr)$.
Due to the diagrams (\ref{eq;06.5.10.1}), (\ref{eq;6.14.10})
and (\ref{eq;6.14.31}),
we have the following commutative diagram:
\[
 \begin{CD}
 C[-1] @>>> \Ob(\nbigv_{\cdot})\\
 @VVV @VVV \\
 L_{\nbigm(m,\vecyhat,[L])/\nbigm(m,y)}[-1]
 @>>>
 L_{\nbigm(m,y)/k}
 \end{CD}
\]
It induces the desired morphism.
\hfill\qed

\vspace{.1in}

Similarly,
on the moduli stack $\nbigm(m,\vecy,L)$,
we put as follows:
\[
 \Ob\bigl(m,\vecy,L\bigr):=
\Cone\Bigl(
 \Ob_{\rel}(m,y,L)[-1]
\oplus
 \Ob_{\rel}\bigl(m,\vecy\bigr)[-1]
\lrarr
 \Ob\bigl(m,y\bigr)
 \Bigr).
\]
By using the sequence of the morphisms,
$\nbigm(m,\vecy,L)
\lrarr \nbigm(m,\vecy)
\lrarr\nbigm(m,y)$,
we obtain the morphism:
$\ob\bigl(m,\vecy,L\bigr):
\Ob\bigl(m,\vecy,L\bigr)
\lrarr
 L_{\nbigm(m,\vecy,L)/k}$.

\vspace{.1in}

We put as follows on the moduli
$\nbigm(m,\vecyhat)$:
\[
 \Ob\bigl(m,\vecyhat\bigr):=
\Cone\Bigl(
 \Ob_{\rel}(m,\yhat)[-1]
\oplus
 \Ob_{\rel}(m,\vecy)[-1]
\lrarr
 \Ob(m,y)
 \Bigr).
\]
Then we obtain the morphism
$ \ob\bigl(m,\vecyhat\bigr):
\Ob\bigl(m,\vecyhat\bigr)
\lrarr
 L_{\nbigm(m,\vecyhat)/k}$.
\index{
 $\Ob\bigl(m,\vecyhat\bigr)$,
 $\ob\bigl(m,\vecyhat\bigr)$}

\vspace{.1in}

Let $\vecL=(L_1,L_2)$ be a tuple of
line bundles on $X$.
We put as follows, on $\nbigm(m,\vecyhat,[\vecL])$:
\[
 \Ob\bigl(m, \vecyhat,[\vecL]\bigr):=
\Cone\Bigl(
 \Ob_{\rel}(m,\yhat)[-1]
\oplus
 \Ob_{\rel}(m,\vecy)[-1]
\oplus\bigoplus_{i=1,2}
 \Ob_{\rel}(m,y,[L_i])[-1].
\lrarr
 \Ob(m,y)
 \Bigr).
\]
Then we obtain the morphism
$\ob\bigl(m,\vecyhat,[\vecL]\bigr):
\Ob\bigl(m,\vecyhat,[\vecL]\bigr)
\lrarr
 L_{\nbigm(m,\vecyhat,[\vecL])/k}$.
\index{$\Ob\bigl(m, \vecyhat,[\vecL]\bigr)$,
 $\ob\bigl(m,\vecyhat,[\vecL]\bigr)$}

\subsubsection{Obstruction theory of the quot scheme
 and the moduli stacks}

Let $Q^{\circ}(m,\vecyhat,[L])$ be as in
the subsubsection \ref{subsubsection;06.5.4.20}.
Let $\pi:Q^{\circ}(m,\vecyhat,[L])\lrarr
 \nbigm(m,\vecyhat,[L])$ denote 
the projection.
On $Q^{\circ}(m,\vecyhat,[L])$,
we put as follows:
\[
 \Ob^Q(m,\vecyhat,[L]):=
\Cone\Bigl(
 \pi^{\ast}\Ob_{\rel}(m,\yhat)
\oplus
 \pi^{\ast}\Ob_{\rel}(m,y,[L])
\oplus
 \pi^{\ast}\Ob_{\rel}(m,\vecy)
\lrarr
 \Ob(V_{-1},f)
 \Bigr)
\]
Then we obtain the morphism
$\ob^Q(m,\vecyhat,[L]):
 \Ob^Q(m,\vecyhat,[L])
\lrarr L_{Q^{\circ}(m,\vecyhat,[L])/k}$
by an argument as in the subsubsection
\ref{subsubsection;06.5.4.450}.

\begin{prop}
 \label{prop;06.5.4.500}
The morphism $\ob^Q(m,\vecyhat,[L])$ 
gives an obstruction theory
of $Q^{\circ}(m,\vecyhat,[L])$.
\end{prop}
\pf
It follows from 
Proposition \ref{prop;06.5.2.1},
Remark \ref{rem;06.5.3.30},
Lemma \ref{lem;06.5.3.40},
\ref{lem;06.5.3.41},
and Lemma \ref{lem;06.5.3.42}.
\hfill\qed

\begin{prop}
\label{prop;06.5.11.2}
The morphism
$\ob(m,\vecyhat,[L])$ gives an obstruction
theory of $\nbigm(m,\vecyhat,[L])$ over $k$.
\end{prop}
\pf
We have
$\pi^{\ast}\nbigh^i\bigl(
 \Ob(m,\vecyhat,[L])
 \bigr)
\simeq
  \nbigh^i\bigl(
 \Ob^Q(m,\vecyhat,[L]) \bigr)$
for $i<0$.
We also have the following diagram
from Lemma \ref{lem;06.5.2.2}
and the construction of the complexes:
\[
 \begin{CD}
 \pi^{\ast}
 \nbigh^0\bigl(
 \Ob(m,\vecyhat,[L])
 \bigr)
 @>>>
 \nbigh^0\bigl(
 \Ob^Q(m,\vecyhat,[L])
 \bigr)
 @>>>
 \nhom(V',V')
 @>>>
 \pi^{\ast}
 \nbigh^1\bigl(
 \Ob(m,\vecyhat,[L])
 \bigr)\\
 @VVV @V{\simeq}VV @V{\simeq}VV @VVV \\
 \pi^{\ast}\nbigh^0(L_{\nbigm(m,\vecyhat,[L])/k})
 @>>>
 \nbigh^0(L_{Q^{\circ}(m,\vecyhat,[L])/k})
 @>>>
 \nhom(V',V')
 @>>>
 \pi^{\ast}\nbigh^1(L_{\nbigm(m,\vecyhat,[L])/k})
 \end{CD}
\]
We also remark that
$\nbigh^1\bigl(\Ob^Q(m,\vecy,[L])\bigr)
=\nbigh^1\bigl(L_{Q^{\circ}(m,\vecy,[L])/k}\bigr)=0$
and 
$\nbigh^{-1}(L_{Q(m)/\nbigm(m)})=0$.
Therefore,
the claim follows from Proposition \ref{prop;06.5.4.500}.
\hfill\qed

\vspace{.1in}
By a similar argument,
we obtain the following:
\begin{prop}
$\ob(m,\vecy,L)$,
$\ob(m,\vecyhat)$
and $\ob(m,\vecyhat,[\vecL])$
give obstruction theories of
$\nbigm(m,\vecy,L)$,
$\nbigm(m,\vecyhat)$
and $\nbigm(m,\vecyhat,[\vecL])$
respectively.
\hfill\qed
\end{prop}

\subsubsection{Obstruction theory of 
 the moduli stacks of stable objects}

Let $\alpha_{\ast}$ be a system of weight,
and let $\delta$ be an element of $\nbigp^{\br}$.
Take a sufficiently large integer $m$.
Then $\nbigm^s(\vecyhat,[L],\alpha_{\ast},\delta)$
is the open substack of $\nbigm(m,\vecyhat,[L])$.

\begin{prop}
\label{prop;06.5.3.50}
The restriction of
$\ob(m,\vecyhat,[L])$
gives an obstruction theory of 
$\nbigm^s(\vecyhat,[L],\delta,\alpha_{\ast})$.
It is independent of a choice of $m$.
\end{prop}
\pf
The first claim follows from 
Proposition \ref{prop;06.5.11.2}.
The second claim follows from 
Lemma \ref{lem;6.16.20},
Lemma \ref{lem;6.14.27}
and Lemma \ref{lem;6.16.30}.
\hfill\qed

\vspace{.1in}

By the same argument,
we can show the following proposition:
\begin{prop}
Let $\vecy$ be an element of $\Type$.
Let $\alpha_{\ast}$ be a system of weights.
We take a large integer $m$ appropriately
in the following claims.
\begin{itemize}
\item
The morphism
$\ob(m,\vecyhat):\Ob(m,\vecyhat)\lrarr
 L_{\nbigm^s(\vecyhat,\alpha_{\ast})/k}$
gives an obstruction theory.
\item
Let $\delta$ be an element of $\nbigp^{\br}$.
The morphism
$ob\bigl(m,\vecy,L\bigr):
 \Ob(m,\vecy,L)\lrarr
 L_{\nbigm^s(\vecy,L,\alpha_{\ast},\delta)}$
gives an obstruction theory.
\item
Let $\vecdelta=(\delta_1,\delta_2)$ be an 
element of $\nbigp^{\br}$.
The morphism
$\ob\bigl(m,\vecy,[\vecL]\bigr):
 \Ob(m,\vecy,[\vecL])\lrarr
 L_{\nbigm^s(\vecyhat,[\vecL],\alpha_{\ast},\vecdelta)/k}$
gives an obstruction theory.
\hfill\qed
\end{itemize}
\end{prop}

\subsection{Obstruction Theory for the Enhanced Master Space
 and the Related Stacks}
\label{subsection;06.5.23.1}
\subsubsection{The enhanced master space}
\label{subsubsection;06.5.8.5}

Take a sufficiently large number $m$.
We put $N:=H_y(m)$.
We take an $N$-dimensional vector space $V_m$.
We put $\proj_m:=\proj(V_m^{\lor})$.
We put 
$Z_1:=\proj\bigl(
 \nbigo_{\proj_m}(0)\oplus \nbigo_{\proj_m}(1) \bigr)$
over $\proj_m$.
We have the natural right $\GL(V_m)$-action
on $Z_2:=Z_1\times\Flag(V_m,\Nbar)$,
where $\Flag(V_m,\Nbar)$ denotes the full flag variety
of $V_m$
as in (\ref{eq;06.6.15.1}).
The quotient stack is denoted by
$\nbigqtilde$.
The quotient stack 
$(\proj_m)_{\GL(V_m)}$ is same as
$B(W_{\cdot},[P_{\cdot}])$
in the subsubsection \ref{subsubsection;06.5.4.150},
when $W_0=V_m$.
So we use the notation.

Let us fix an inclusion
$\iota:\nbigo(-m)\lrarr L$.
Since reduced $L$-sections induce
$\nbigo(-m)$-sections,
we obtain the morphism
$\Xi(\nbigv_{\cdot},[\phibar]):
 \nbigm(m,\vecyhat,[L])\lrarr B(W_{\cdot},[P_{\cdot}])$.
For simplicity of the notation,
we use $\Psi_1$ instead of 
$\Xi(\nbigv_{\cdot},[\phibar])$.
We use the following lemma
in the construction of the deformation theory
of the enhanced master space.

\begin{lem}
 \label{lem;06.5.4.350}
We have the morphism
$\varphi:\Psi_1^{\ast}L_{B(W_{\cdot},[P_{\cdot}])/k}\lrarr 
 \Ob\bigl(m,\vecyhat,[L]\bigr)$
such that 
the composite of $\varphi$ and
$\ob\bigl(m,\vecyhat,[L]\bigr)$
 is same as the naturally 
defined morphism
$\Psi_1^{\ast}L_{B(W_{\cdot},[P_{\cdot}])/k}
\lrarr L_{\nbigm(m,\vecyhat,[L])/k}$.
\end{lem}
\pf
Recall the complexes $\Ob^G(\nbigv_{\cdot})$
and $\Ob^{G}_{\rel}(\nbigv_{\cdot},[\phitilde])$
constructed in the subsubsection 
\ref{subsubsection;06.5.4.100}
and the subsubsection \ref{subsubsection;06.5.4.150}
respectively.
We put
$\Ob^G(\nbigv_{\cdot},[\phitilde])=
 \Cone\bigl(\Ob^G_{\rel}(\nbigv_{\cdot},[\phitilde])[-1]
\lrarr\Ob^G(\nbigv)\bigr)$.
Then we obtain the following 
commutative diagram
from (\ref{eq;06.5.4.300})
and (\ref{eq;06.5.4.250}):
\[
 \begin{CD}
 \Ob^G(\nbigv_{\cdot},[\phitilde])
 @>>> 
 \Ob\bigl(m,\vecyhat,[L]\bigr) \\
 @V{\varphi_1}VV @VVV \\
 \Psi_1^{\ast}L_{B(W_{\cdot},[P_{\cdot}])/k}
 @>>> L_{\nbigm(m,\vecyhat,[L])/k}.
 \end{CD}
\]
It is easy to show that $\varphi_1$
is isomorphic due to Lemma \ref{lem;06.5.4.201}
and Lemma \ref{lem;06.5.4.200}.
\hfill\qed

\vspace{.1in}

The fiber product
$\nbigm(m,\vecyhat,[L])\times_{B(W_{\cdot},[P_{\cdot}])}
 \nbigqtilde$ is denoted by $\nbign$.
By construction, the enhanced master space
$\Mhat$ is an open subset of $\nbign$.
The induced morphism $\nbign\lrarr \nbigqtilde$
is denoted by $\Psi_2$.
Let $\gminip$ denote the naturally defined morphism
of $\Mhat$ to $\nbigm(m,\vecyhat,[L])$.
Then we obtain the following morphism
of the distinguished triangles on $\Mhat$:
\begin{equation}
 \label{eq;06.5.4.360}
 \begin{CD}
  \Psi_2^{\ast}L_{\nbigqtilde/B(W_{\cdot},[P_{\cdot}])}[-1]
 @>>>
 \gminip^{\ast}\Psi_1^{\ast} L_{B(W_{\cdot},[P_{\cdot}])/k}
 @>>>
 \Psi_2^{\ast}L_{\nbigqtilde/k}
 @>>>
 \Psi_2^{\ast}L_{\nbigqtilde/B(W_{\cdot},[P_{\cdot}])}
 \\
 @V{\simeq}VV @VVV @VVV @VVV \\
 L_{\Mhat/\nbigm(m,\vecyhat,[L])}[-1]
 @>>>
 \gminip^{\ast}L_{\nbigm(m,\vecyhat,[L])/k}
 @>>>
 L_{\Mhat/k}
 @>>>
 L_{\Mhat/\nbigm(m,\vecyhat,[L])}
 \end{CD}
\end{equation}

From Lemma \ref{lem;06.5.4.350}
and the diagram (\ref{eq;06.5.4.360}),
we obtain the morphism
$\Psi_2^{\ast}L_{\nbigqtilde/B(W_{\cdot},[P_{\cdot}])}[-1]\lrarr 
 \gminip^{\ast}\Ob(m,\vecy,[L])$.
The cone is denoted by $\Ob(\Mhat)$.
We have the induced morphism
$\ob(\Mhat):\Ob(\Mhat)\lrarr L_{\Mhat/k}$.

\begin{prop}
\label{prop;06.5.6.1}
The morphism $\ob(\Mhat)$ gives
an obstruction theory of the master space $\Mhat$.
\end{prop}
\pf
We put $\nbigm:=\nbigm(m,\vecyhat,[L])$.
By construction,
we have the following morphism of distinguished triangles:
\begin{equation}
 \label{eq;06.6.13.1}
 \begin{CD}
 L_{\Mhat/\nbigm}[-1]
 @>>>
 \gminip^{\ast}\Ob(m,\vecyhat,[L])
 @>>>
 \Ob(\Mhat)
 @>>>
 L_{\Mhat/\nbigm}\\
 @V{\simeq}VV @VVV @VVV @V{\simeq}VV \\
 L_{\Mhat/\nbigm}[-1]
 @>>>
 \gminip^{\ast}L_{\nbigm/k}
 @>>>
 L_{\Mhat/k}
 @>>>
 L_{\Mhat/\nbigm}
 \end{CD}
\end{equation}
Then the claim follows from
Proposition \ref{prop;06.5.11.2}.
\hfill\qed

\subsubsection{The substack $\Mhat^{\ast}$}
\label{subsubsection;06.5.22.1000}

We have the natural $\GL(V_m)$-action
on $V_m^{\ast}\times\Flag(V_m,\Nbar)$.
The quotient stack is denoted by $\nbigqtilde^{\ast}$.
We put $B(W_{\cdot},P_{\cdot}):=(V_m^{\ast})_{\GL(V_m)}$.
We have the natural morphism
$\nbigqtilde^{\ast}\lrarr B(W_{\cdot},P_{\cdot})$
which is a full flag bundle.
Since $V_m^{\ast}$ naturally gives
the open subset 
$Z_1-\proj\bigl(\nbigo_{\proj}(0)\bigr)\cup
 \proj\bigl(\nbigo_{\proj}(1)\bigr)$ of $Z_1$,
the stack $\nbigqtilde^{\ast}$ naturally gives the open subset
of $\nbigqtilde$.

Recall that
we put $\Mhat^{\ast}:=\Mhat-\bigl(\Mhat_1\cup\Mhat_2\bigr)$
(the subsection \ref{subsection;06.5.21.110}).
The stack $\Mhat^{\ast}$
is an open subset of
$\nbigm(m,\vecyhat,L)\times_{B(W_{\cdot},P_{\cdot})}
 \nbigqtilde^{\ast}$ by construction.
We have the commutative diagrams:
\[
 \begin{CD}
 \Mhat^{\ast} @>{\Psi_4}>> \nbigqtilde^{\ast}\\
 @V{\gminip_2}VV @VVV \\
 \nbigm(m,\vecyhat,L) @>{\Psi_3}>> B(W_{\cdot},P_{\cdot})
 \end{CD}
\quad\quad\quad\quad
 \begin{CD}
 \Psi^{\ast}_4L_{\nbigqtilde^{\ast}/B(W_{\cdot},P_{\cdot})}[-1]
@>>>
 \gminip_2^{\ast}\Psi_3^{\ast}L_{B(W_{\cdot},P_{\cdot})}\\
 @V{\simeq}VV @V{\varphi}VV \\
 L_{\Mhat^{\ast}/\nbigm(m,\vecyhat,L)}[-1]
 @>>>
 \gminip_2^{\ast}L_{\nbigm(m,\vecyhat,L)}
 \end{CD} 
\]

\begin{lem}
$\varphi$ factors through $\gminip_2^{\ast}\Ob(m,\vecyhat,L)$.
\end{lem}
\pf
Recall the complexes
$\Ob^{G}(\nbigv_{\cdot})$
and $\Ob^{G}_{\rel}(\nbigv_{\cdot},\phitilde)$
constructed in the subsubsection
\ref{subsubsection;06.5.4.100}
and the subsubsection \ref{subsubsection;06.5.21.250},
respectively.
We put $\Ob^G(\nbigv_{\cdot},\phitilde):=\cone\bigl(
 \Ob_{\rel}^G(\nbigv_{\cdot},\phitilde)[-1]
\lrarr \Ob^G(\nbigv_{\cdot}) \bigr)$.
We obtain the following commutative diagram
on $\nbigm(m,\vecy,L)$:
\[
 \begin{CD}
 \Ob^{G}(\nbigv_{\cdot},\phitilde)
 @>>> \Ob(m,\vecy,L)\\
 @V{\varphi_2}VV @VVV \\
 \Psi_3^{\ast}L_{B(W,P)/k} @>>>
 L_{\nbigm(m,\vecy,L)}
 \end{CD}
\]
It is easy to show $\varphi_2$ is isomorphic
by using Lemma \ref{lem;06.5.4.201}
and Lemma \ref{lem;06.5.2.10}.
\hfill\qed

\vspace{.1in}

Since $\varphi$ factors through
$\Ob(m,\vecyhat,L)$,
we obtain the morphism
$\Psi^{\ast}_4L_{\nbigqtilde/B(W,P)}[-1]\lrarr\Ob(m,\vecyhat,L)$.
The cone is denoted by
$\Ob(\Mhat^{\ast})$.
Then we obtain the morphism
$\ob(\Mhat^{\ast}):
 \Ob(\Mhat^{\ast})\lrarr L_{\Mhat^{\ast}}$.

\begin{lem}
\label{lem;06.5.8.7}
There exists
the quasi isomorphism 
$u:\Ob(\Mhat)_{|\Mhat^{\ast}}\lrarr\Ob(\Mhat^{\ast})$
such that
 $\ob(\Mhat^{\ast})\circ u=\ob(\Mhat)$.
\end{lem}
\pf
Let $\pi_1$ denote the natural morphism
$\nbigm(m,\vecyhat,L)\lrarr\nbigm(m,\vecyhat,[L])$.
By construction of $\Ob^G(\nbigv_{\cdot},[\phitilde])$
and $\Ob^G(\nbigv_{\cdot},\phitilde)$,
we have the following commutative diagram
on $\nbigm(m,\vecyhat,L)$:
\begin{equation}
 \label{eq;06.5.21.500}
 \begin{CD}
 \Phi_2^{\ast}L_{B(W,P)/k} @<{\simeq}<<
 \Ob^G(\nbigv,\phitilde) @>>>
 \Ob(m,\vecyhat,L)\\
 @A{\psi_1}AA @A{\psi_2}AA @A{\psi_3}AA \\
 \pi_1^{\ast}\Phi_1^{\ast}L_{B(W,[P])/k}
 @<{\simeq}<< \pi_1^{\ast}\Ob^G(\nbigv,[\phitilde])
 @>>> \pi_1^{\ast}\Ob(m,\vecyhat,[L])
 \end{CD}
\end{equation}
Moreover, the induced morphisms
$\cone(\psi_2)\lrarr \cone(\psi_i)$ $(i=1,3)$
are isomorphisms.
Hence,
we obtain the following commutative diagram
on $\Mhat^{\ast}$:
\[
 \begin{CD}
 \Psi^{\ast}_4L_{\nbigqtilde^{\ast}/B(W,P)}[-1] @>>>
 \gminip_2^{\ast}\Psi_3^{\ast}L_{B(W,P)} @>>>
 \gminip_2^{\ast}\Ob(m,\vecyhat,L) @>>>
 \gminip_2^{\ast}L_{\nbigm(m,\vecyhat,L)}\\
 @A{\mu_1}AA @A{\mu_2}AA @A{\mu_3}AA @AAA \\
 \Psi_4^{\ast}L_{\nbigqtilde^{\ast}/B(W,[P])}[-1] @>>>
 \gminip^{\ast}\Psi_1^{\ast}L_{B(W,[P])} @>>>
 \gminip^{\ast}\Ob(m,\vecyhat,[L]) @>>>
 \gminip^{\ast}L_{\nbigm(m,\vecyhat,[L])}
 \end{CD}
\]
Moreover, the induced morphisms
$\cone(\mu_1)\lrarr \cone(\mu_2)\lrarr\cone(\mu_3)$
are quasi isomorphic.
Then the claim of the lemma is clear.
\hfill\qed

\vspace{.1in}

We have another description of $\Ob(\Mhat^{\ast})$.
When $V_m=W_0$,
we have $B(W_{\cdot}):=\Spec(k)_{\GL(V_m)}$.
We put $\Fbar:=\Flag(V_m,\Nbar)_{\GL(V_m)}$.
We have the following commutative diagram:
\begin{equation}
 \label{eq;06.6.15.3}
 \begin{CD}
 \Mhat^{\ast}@>{\Psi_4}>>\nbigqtilde^{\ast} @>{\Gamma_1}>> \Fbar\\
 @V{\gminip_2}VV  @VVV @VVV \\
 \nbigm(m,\vecyhat,L) @>{\Psi_2}>> B(W_{\cdot},P_{\cdot})
 @>{\Gamma_2}>> B(W_{\cdot})
 \end{CD}
\end{equation}
Then, we have the isomorphism
$\Psi_4^{\ast}L_{\nbigqtilde^{\ast}/B(W_{\cdot},P_{\cdot})}
\simeq
 \Psi_4^{\ast}\Gamma_1^{\ast}L_{\Fbar/B(W_{\cdot})}$.
We also obtain the following morphisms:
\[
\begin{CD}
 \Psi_2^{\ast} \Gamma_2^{\ast}L_{B(W_{\cdot})}
 @>>>
 \Psi_2^{\ast}L_{B(W_{\cdot},P_{\cdot})}
 @>>>
 \Ob(m,\vecyhat,L)
\end{CD}
\]
Therefore, we obtain the morphism
$\alpha:\Psi_4^{\ast}\Gamma_1^{\ast}
 L_{\Fbar/B(W_{\cdot})}[-1]
\lrarr 
 \gminip_2^{\ast}\Ob(m,\vecyhat,L)$.
We naturally obtain the following quasi-isomorphism:
\begin{equation}
\label{eq;06.5.22.250}
 \cone(\alpha)\simeq \Ob(\Mhat^{\ast})
\end{equation}

\subsubsection{The moduli stack $\nbigmtilde(m,\vecyhat,[L])$}
\label{subsubsection;06.6.13.40}

Let $\nbigmtilde(m,\vecyhat,[L])$ be the moduli stack of
the tuple $(E_{\ast},\rho,[\phi],\nbigf)$ as follows:
\begin{itemize}
\item
 $(E_{\ast},\rho,[\phi])$ is an oriented reduced $L$-Bradlow pair
 of type $\vecy$,
 satisfying the condition $O_m$.
\item
 $\nbigf$ is a full flag of $H^0(X,E(m))$.
\end{itemize}

By an argument in the subsubsection 
\ref{subsubsection;06.5.8.5},
we can obtain the obstruction theory of
$\nbigmtilde(m,\vecyhat,[L])$.
We also use the notation
$\nbigmtilde$ and $\nbigm$
to denote
$\nbigmtilde(m,\vecyhat,[L])$
and $\nbigm(m,\vecyhat,[L])$,
respectively.
When $V_m=W_0$,
we have $B(W_{\cdot}):=\Spec(k)_{\GL(V_m)}$.
We also put $\Fbar:=\Flag(V_m)_{\GL(V_m)}$.
The following diagram is Cartesian:
\[
  \begin{CD}
 \nbigmtilde(m,\vecyhat,[L]) @>{\Psi_{11}}>> \Fbar\\
 @V{\gminip_1}VV @VVV \\
 \nbigm(m,\vecyhat,[L]) @>{\Psi_{12}}>> B(W_{\cdot})
 \end{CD}
\]
Then, we obtain the following morphism of distinguished triangles
on $\nbigmtilde(m,\vecyhat,[L])$:
\[
 \begin{CD}
 \Psi_{11}^{\ast}L_{\Fbar/B(W_{\cdot})}[-1]
@>>>
 \gminip_1^{\ast}\Psi_{12}^{\ast}
 L_{B(W_{\cdot})/k}
@>>>
 \Psi_{11}^{\ast}L_{\Fbar/k}
@>>>
 \Psi_{11}^{\ast}L_{\Fbar/B(W_{\cdot})}\\
 @V{\simeq}VV @V{\varphi}VV @VVV @V{\simeq}VV \\
 L_{\nbigmtilde/\nbigm}[-1]
 @>>>
 \gminip_1^{\ast}L_{\nbigm/k}
 @>>>
 L_{\nbigmtilde/k}
 @>>>
 L_{\nbigmtilde/\nbigm} \\
 \end{CD}
\]

\begin{lem}
$\varphi$ factors through
$\gminip_1^{\ast}\Ob(m,\vecyhat,[L])$
\end{lem}
\pf
It can be shown by using the complex $\Ob^G(V_{\cdot})$
and the argument in the proof of 
Lemma \ref{lem;06.5.4.350}.
\hfill\qed

\vspace{.1in}

Then, we obtain the morphism
$\Psi_{11}^{\ast}L_{\Fbar/B(W_{\cdot})}[-1]
\lrarr \gminip_1^{\ast}\Ob(m,\vecyhat,[L])$.
The cone is denoted by 
$\Obtilde(m,\vecyhat,[L])$.
We obtain the natural morphism:
\[
 \obtilde\bigl(m,\vecyhat,[L]\bigr):
\Obtilde(m,\vecyhat,[L])
\lrarr L_{\nbigmtilde(m,\vecyhat,[L])/k}
\]
\index{$\Obtilde(m,\vecyhat,[L])$,
 $\obtilde(m,\vecyhat,[L])$}

By the same argument as the proof of Proposition
\ref{prop;06.5.6.1},
we can show the following.
\begin{prop}
$\obtilde\bigl(m,\vecyhat,[L]\bigr)$
gives an obstruction theory 
for
$\nbigmtilde\bigl(m,\vecy,[L]\bigr)$.
\hfill\qed
\end{prop}

We have the equivalent obstruction theory.
Let $\nbigq_1$ denote the quotient stack
of $\proj_m\times\Flag(V,\Nbar)$
via the natural $\GL(V_m)$-action.
Then, the following diagram is Cartesian:
\[
 \begin{CD}
 \nbigmtilde(m,\vecyhat,[L]) @>{\Psi_{13}}>> \nbigq_1\\
 @V{\gminip_1}VV @VVV \\
 \nbigm(m,\vecyhat,[L]) @>{\Psi_{14}}>>  B(W,[P])
 \end{CD}
\]
Then, we have the following diagram
on $\nbigmtilde(m,\vecyhat,[L])$:
\[
 \begin{CD}
 \Psi_{13}^{\ast}L_{\nbigq_1/B(W,[P])}[-1]
@>>>
 \gminip_1^{\ast}\Psi_{14}^{\ast}L_{B(W,[P])}
@>>>
 \Psi_{13}^{\ast}L_{\nbigq_1/k}
@>>>
 \Psi_{13}^{\ast}L_{\nbigq_1/B(W,[P])}\\
 @VVV @V{\varphi_1}VV @VVV @VVV\\
L_{\nbigmtilde/\nbigm}[-1]
@>>>
 \gminip_1^{\ast}L_{\nbigm/k}
@>>>
 L_{\nbigmtilde/k}
@>>>
L_{\nbigmtilde/\nbigm}
 \end{CD}
\]
By an argument similar to the proof of
Lemma \ref{lem;06.5.4.350},
we can show that $\varphi_1$
factors through $\gminip_1^{\ast}\Ob(m,\vecyhat,[L])$.
Let $\Obtilde_2(m,\vecyhat,[L])$
denote the cone of
$\Psi_{13}^{\ast}L_{\nbigq_1/B(W,[P])}[-1]
\lrarr
 \gminip_1^{\ast}\Ob(m,\vecyhat,[L])$,
and then we have the naturally defined morphism:
\[
 \obtilde_2(m,\vecyhat,[L]):
 \Obtilde_2(m,\vecyhat,[L])
\lrarr L_{\nbigmtilde(m,\vecyhat,[L])}
\]

\begin{lem}
We have the natural quasi isomorphism
$\psi:\Ob(m,\vecyhat,[L])\lrarr
 \Ob_2(m,\vecyhat,[L])$
such that the composite
$\ob_2(m,\vecyhat,[L])\circ
\psi$ is same as $\ob(m,\vecyhat,[L])$.
\end{lem}
\pf
Let $\nbigv_{\cdot}$ be the canonical resolution
of $\nbige^u(m)$ over $\nbigm(m,\vecyhat,[L])$.
We obtain the complexes
$\Ob^G(\nbigv_{\cdot})$
and $\Ob^{G}_{\rel}(\nbigv_{\cdot},[\phitilde])$
by the constructions in the subsubsection
\ref{subsubsection;06.5.4.100}
and the subsubsection \ref{subsubsection;06.5.4.150}.
Let $\Ob^G(\nbigv_{\cdot},[\phitilde])$
denote the cone of the morphism
$\Ob^G_{\rel}(\nbigv_{\cdot},[\phitilde])[-1]
 \lrarr \Ob^G(\nbigv_{\cdot})$.
We have the following commutative diagram
on $\nbigm(m,\vecyhat,[L])$:
\[
 \begin{CD}
 \Psi_{12}^{\ast}L_{B(W)} @>>>
 \Ob^G(\nbigv_{\cdot})@>>>
 \Ob(m,\vecyhat)\\
 @VVV @VVV @VVV \\
 \Psi_{14}^{\ast}L_{B(W,[P])} @>>>
 \Ob^{G}(\nbigv_{\cdot},[\phitilde]) @>>>
 \Ob(m,\vecyhat,[L])
 \end{CD}
\]
We have the commutative diagram:
\[
 \begin{CD}
 \nbigmtilde(m,\vecyhat,[L])@>>>\nbigq_1 @>>>\Fbar\\
 @VVV @VVV @VVV \\
 \nbigm(m,\vecyhat,[L]) @>>>  B(W,[P]) @>>> B(W)
 \end{CD}
\]
We obtain the following diagram
on $\nbigmtilde(m,\vecyhat,[L])$:
\[
 \begin{CD}
 \Psi_{13}^{\ast}L_{\nbigq_1/B(W,[P])}[-1]
 @>>>
 \gminip_1^{\ast}\Psi_{14}^{\ast}L_{B(W,[P])} @>>>
 \gminip_1^{\ast}\Ob(m,\vecyhat,[L])
 @>>>
 \gminip_1^{\ast}L_{\nbigm(m,\vecyhat,[L])}\\
 @A{\simeq}AA @AAA @AAA @AAA\\
 \Psi_{11}^{\ast}L_{\Fbar/B(W)}[-1] @>>>
 \gminip_1^{\ast}\Psi_{12}^{\ast}L_{B(W)}@>>>
 \gminip_1^{\ast}\Ob(m,\vecyhat,[L])
 @>>>
 \gminip_1^{\ast}L_{\nbigm(m,\vecyhat,[L])}
 \end{CD}
\]
Then, the claim is clear.
\hfill\qed

\vspace{.2in}

Recall that the moduli stack 
$\nbigmtilde^s\bigl(\vecyhat,[L],\alpha_{\ast},(\delta,\ell)\bigr)$
of $(\delta,\ell)$-stable objects
is the open substack of $\nbigmtilde(m,\vecyhat,[L])$.
By restricting 
$\obtilde(m,\vecyhat,[L])$,
we obtain the obstruction theory
of $\nbigmtilde^{s}\bigl(\vecyhat,[L],\alpha_{\ast},(\delta,\ell)\bigr)$.

\subsubsection{The moduli stack $\nbigmtilde(m,\vecy,L)$}
\label{subsubsection;06.6.13.50}

Let $\nbigmtilde(m,\vecy,L)$ be the moduli stack of
the tuple $(E_{\ast},\phi,\nbigf)$ as follows:
\begin{itemize}
\item
 $(E_{\ast},\phi)$ is an $L$-Bradlow pair
 of type $\vecy$,
 satisfying the condition $O_m$.
\item
 $\nbigf$ is a full flag of $H^0(X,E(m))$.
\end{itemize}

We use the notation in the subsubsection
\ref{subsubsection;06.6.13.40}.
We have the following Cartesian diagram:
\[
  \begin{CD}
 \nbigmtilde(m,\vecy,L) @>{\Psi_{11}}>> \Fbar\\
 @V{\gminip_1}VV @VVV \\
 \nbigm(m,\vecy,L) @>{\Psi_{12}}>> B(W_{\cdot})
 \end{CD}
\]
By using the construction in 
the subsubsection \ref{subsubsection;06.6.13.40},
we obtain the obstruction theory:
\[
 \obtilde(m,\vecy,L):
 \Obtilde(m,\vecy,L)\lrarr
 L_{\nbigmtilde(m,\vecy,L)}
\]
\index{$\Obtilde(m,\vecy,L)$, $\obtilde(m,\vecy,L)$}

The moduli stack
$\nbigmtilde^{ss}\bigl(\vecy,L,\alpha_{\ast},(\delta,\ell)\bigr)$
(the subsubsection \ref{subsubsection;06.4.25.10})
is the open substack of $\nbigmtilde(m,\vecy,L)$.
By restricting 
$\obtilde(m,\vecy,L)$,
we obtain the obstruction theory
of $\nbigmtilde^{ss}\bigl(\vecy,L,\alpha_{\ast},(\delta,\ell)\bigr)$.

\subsubsection{The moduli stack $\nbigmtilde^s(\vecyhat,\alpha_{\ast},+)$}
\label{subsubsection;06.6.13.51}

Let $\nbigmtilde(m,\vecyhat)$
denote the moduli stack of the objects $(E_{\ast},\nbigf)$
as follows:
\begin{itemize}
\item
 $E_{\ast}$ is a parabolic torsion-free sheaf
 satisfying the condition $O_m$.
\item
 $\nbigf$ is a full flag of $H^0(X,E(m))$.
\end{itemize}
We use the same notation in the subsubsection
\ref{subsubsection;06.6.13.40}.
In this case, we have the following Cartesian diagram:
\[
 \begin{CD}
 \nbigmtilde(m,\vecyhat) @>{\Psi_{11}}>> \Fbar\\
 @VVV @VVV \\
 \nbigm(m,\vecyhat)@>{\Psi_{12}}>> B(W_{\cdot}).
 \end{CD}
\]
By using the construction 
in the subsubsection \ref{subsubsection;06.6.13.40},
we obtain the obstruction theory:
\[
 \obtilde(m,\vecyhat):
 \Obtilde(m,\vecyhat)\lrarr
 L_{\nbigmtilde(m,\vecyhat)}
\]
\index{$\Obtilde(m,\vecyhat)$, $\obtilde(m,\vecyhat)$}

Recall $\nbigmtilde^s(\vecyhat,\alpha_{\ast},+)$
denotes the moduli stack of the objects
$(E_{\ast},\nbigf)$ as follows
(the subsubsection \ref{subsubsection;06.5.17.20}):
\begin{itemize}
\item
 $E_{\ast}$ is a parabolic torsion-free sheaf
 of type $\vecy$ with weight $\alpha_{\ast}$.
\item
 $\nbigf$ is a full flag of $H^0(X,E(m))$.
\item
 $(E_{\ast},\nbigf_{\min})$ is
 $\epsilon$-semistable reduced $\nbigo(-m)$-Bradlow pair,
 where $\epsilon$ denotes any sufficiently small positive number.
\end{itemize}

Since $\nbigmtilde^s(\vecyhat,\alpha_{\ast},+)$
is the open substack of $\nbigmtilde(m,\vecyhat,\alpha_{\ast})$,
we obtain the obstruction theory
of $\nbigmtilde^s(\vecyhat,\alpha_{\ast},+)$
by restricting $\obtilde(m,\vecyhat)$.

\subsubsection{The case where the $2$-stability condition is satisfied}
\label{subsubsection;06.5.21.600}

Let us construct the obstruction theory of the master space
in the case where the $2$-stability condition
is satisfied for $(\vecy,L,\alpha_{\ast},\delta)$.
It can be done by the way given in the subsubsection 
\ref{subsubsection;06.5.8.5},
so we give only an indication.
We use the notation in the subsubsections
\ref{subsubsection;06.5.17.100}
and \ref{subsubsection;06.5.8.5}.

We have the natural $\GL(V)$-action
on the $\proj^1$-bundle
$\proj\bigl(\nbigo_{\proj_m}(0)\oplus\nbigo_{\proj_m}(1)\bigr)$
over $\proj_m$.
The quotient stack is denoted by $\nbigq$.
We have the map
$\Psi_2:\nbigm(m,\vecyhat,[L])\lrarr B(W_{\cdot},[P])$,
and $\Mhat$ is an open substack of
$\nbigm(m,\vecyhat,[L])\times_{B(W,[P])} \nbigq$.
\[
 \begin{CD}
 \Mhat @>{\Psi_1}>> \nbigq \\
 @V{\gminip}VV @VVV\\
 \nbigm(m,\vecyhat,[L])@>{\Psi_2}>> B(W,[P])
 \end{CD}
\]
We obtain the following morphism
of distinguished triangles on $\Mhat$:
\[
 \begin{CD}
 \Psi_1^{\ast}L_{\nbigq/B(W,[P])}[-1]  @>>>
 \gminip^{\ast}\Psi_2^{\ast}L_{B(W,[P])} @>>>
 \Psi_1^{\ast}L_{\nbigq} @>>>
  \Psi_1^{\ast}L_{\nbigq/B(W,[P])} \\
 @VVV @V{\varphi}VV @VVV @VVV \\
 L_{\Mhat/\nbigm(m,\vecyhat,[L])}[-1] @>>>
 \gminip^{\ast}L_{\nbigm(m,\vecyhat,[L])} @>>>
 L_{\Mhat} @>>>
 L_{\Mhat/\nbigm(m,\vecyhat,[L])}
 \end{CD}
\]
Since $\varphi$ factors through
$\gminip^{\ast}\Ob(m,\vecyhat,[L])$
(Lemma \ref{lem;06.5.4.350}),
the morphism
$\Psi_1^{\ast}L_{\nbigq/B(W,[P])}[-1]
\lrarr \Ob(m,\vecyhat,[L])$
is obtained.
The cone is denoted by
$\Ob(\Mhat)$.
We have the naturally defined morphism
$\ob(\Mhat):\Ob(\Mhat)\lrarr L_{\Mhat}$.
By an argument similar to the proof of
Proposition \ref{prop;06.5.6.1},
we can show that $\ob(\Mhat)$
gives an obstruction theory of
$\Mhat$.

\vspace{.1in}

Recall we put $\Mhat^{\ast}:=\Mhat-\Mhat_1\cup\Mhat_2$.
It is an open substack of $\nbigm(m,\vecyhat,L)$.
We put $\Ob(\Mhat^{\ast}):=\ob(m,\vecyhat,L)$,
and then we have the obstruction theory
$\ob(\Mhat^{\ast}):\Ob(\Mhat^{\ast})\lrarr L_{\Mhat^{\ast}}$.

\begin{lem}
\label{lem;06.5.21.654}
There exists the quasi isomorphism
$u:\Ob(\Mhat)_{|\Mhat^{\ast}}\lrarr \Ob(\Mhat^{\ast})$
such that
$\ob(\Mhat^{\ast})\circ u=\ob(\Mhat)$.
\end{lem}
\pf
We have the commutative diagram:
\[
 \begin{CD}
 \Mhat^{\ast} @>{\Psi_1}>> B(W,P)\\
 @V{\pi}VV @VVV \\
 \nbigm(m,\vecyhat,[L])
 @>{\Psi_2}>> B(W,[P])
 \end{CD}
\]
From the diagram (\ref{eq;06.5.21.500}),
we obtain the following diagram:
\[
 \begin{CD}
 \Psi_1^{\ast}L_{B(W,P)} @>>>
 \Ob(m,\vecyhat,L) @>>> L_{\Mhat^{\ast}}\\
 @A{\mu_1}AA @A{\mu_2}AA @A{\mu_3}AA \\
 \pi^{\ast}\Psi_2^{\ast}L_{B(W,[P])} @>>>
 \pi^{\ast}\Ob(m,\vecyhat,[L]) @>>>
 \pi^{\ast}L_{\nbigm(m,\vecyhat,[L])}
 \end{CD}
\]
Moreover, the induced morphisms
$\cone(\mu_1)\lrarr\cone(\mu_2)\lrarr \cone(\mu_3)$
are quasi isomorphic.
Then the claim of the lemma is clear.
\hfill\qed

\subsubsection{The case of oriented reduced $\vecL$-Bradlow pairs}
\label{subsubsection;06.6.13.16}

Let $\vecL=(L_1,L_2)$ be a pair of line bundles over $X$.
Let us construct the obstruction theory
of the master space for the moduli stacks of
the oriented reduced $\vecL$-Bradlow pairs,
under the setting
in the subsubsection \ref{subsubsection;06.5.21.555}.
We give only an indication.
We use the notation in the subsubsection
\ref{subsubsection;06.5.21.600}.

We construct the master space $\Mhat$
as in the subsubsection \ref{subsubsection;06.5.21.555}.
Let us take an inclusion $\iota_1:\nbigo(-m)\lrarr L_1$.
Then the universal reduced $L_1$-section
$[\phi_1^u]$ induces the reduced $\nbigo_X$-section $[\phibar_1]$.
Therefore,
we obtain the morphism
$\Psi_2:\nbigm(m,\vecyhat,[\vecL])\lrarr B(W,[P])$.
By construction,
$\Mhat$ is an open subset of
$\nbigm(m,\vecyhat,[\vecL])\times_{B(W,[P])}\nbigq$.
\[
\begin{CD}
 \Mhat @>{\Psi_1}>> \nbigq \\
 @V{\gminip}VV @VVV\\
 \nbigm(m,\vecyhat,[\vecL])@>{\Psi_2}>> B(W,[P])
 \end{CD}
\]
We obtain the following morphism
of distinguished triangles on $\Mhat$:
\[
 \begin{CD}
 \Psi_1^{\ast}L_{\nbigq/B(W,[P])}[-1]  @>>>
 \gminip^{\ast}\Psi_2^{\ast}L_{B(W,[P])} @>>>
 \Psi_1^{\ast}\nbigq @>>>
  \Psi_1^{\ast}L_{\nbigq/B(W,[P])} \\
 @VVV @V{\varphi}VV @VVV @VVV \\
 L_{\Mhat/\nbigm(m,\vecyhat,[\vecL])}[-1] @>>>
 \gminip^{\ast}L_{\nbigm(m,\vecyhat,[\vecL])} @>>>
 L_{\Mhat} @>>>
 L_{\Mhat/\nbigm(m,\vecyhat,[\vecL])}
 \end{CD}
\]
By an argument similar to the proof of Lemma 
\ref{lem;06.5.4.350},
it can be shown that
$\varphi$ factors through
$\gminip^{\ast}\Ob(m,\vecyhat,[\vecL])$.
Hence we obtain the morphism
$\Psi_1^{\ast}L_{\nbigq/B(W,[P])}[-1]
\lrarr \Ob(m,\vecyhat,[\vecL])$.
The cone is denoted by
$\Ob(\Mhat)$.
We have the naturally defined morphism
$\ob(\Mhat):\Ob(\Mhat)\lrarr L_{\Mhat}$.
By an argument similar to the proof of Proposition
\ref{prop;06.5.6.1},
we can show that $\ob(\Mhat)$
gives an obstruction theory of
$\Mhat$.

\vspace{.1in}

Recall we put $\Mhat^{\ast}:=\Mhat-\Mhat_1\cup\Mhat_2$.
It is an open substack of the moduli stack
$\nbigm(m,\vecyhat,L_1,[L_2])$.
(See the subsubsection \ref{subsubsection;06.5.21.555}
for $\nbigm(m,\vecyhat,L_1,[L_2])$.)
On $\nbigm(m,\vecyhat,L_1,[L_2])$,
we have the following morphism:
\[
 \Ob_{\rel}(m,\vecy)[-1]
\oplus
 \Ob_{\rel}(m,y,L_1)[-1]
\oplus
 \Ob_{\rel}(m,y,[L_2])[-1]
\oplus
 \Ob_{\rel}(m,\yhat)[-1]
\lrarr
 \Ob(m,y)
\]
The cone is denoted by $\Ob(m,\vecyhat,L_1,[L_2])$.
As in the subsubsection \ref{subsubsection;06.5.4.450},
we can naturally construct the morphism
$\ob(m,\vecyhat,L_1,[L_2]):
\Ob(m,\vecyhat,L_1,[L_2])
\lrarr L_{\nbigm(m,\vecyhat,L_1,[L_2])}$.
By an argument similar to the proof of Proposition 
\ref{prop;06.5.3.50},
it can be shown that
$\ob(m,\vecyhat,L_1,[L_2])$ gives 
an obstruction theory.

\begin{lem}
There exists the quasi isomorphism
$u:\Ob(\Mhat)_{|\Mhat^{\ast}}\lrarr \Ob(\Mhat^{\ast})$
such that
$\ob(\Mhat^{\ast})\circ u=\ob(\Mhat)$.
\end{lem}
\pf
It can be shown by
the same argument as the proof of Lemma 
\ref{lem;06.5.21.654}.
\hfill\qed

\subsection{Moduli Theoretic Obstruction Theory 
of the Fixed Point Set}
\label{subsection;06.6.21.2}
\subsubsection{Statement}

Let $\gbigi=(\vecy_1,\vecy_2,I_1,I_2)$
be a decomposition type as in Definition \ref{df;06.5.22.1}.
We use the notation in the subsubsection
\ref{subsubsection;06.5.17.20}.
We put
$\nbigmtilde_{\spl}:=
 \nbigmtilde^{ss}\bigl(\vecy_1,L,\alpha_{\ast},(\delta,k_0)\bigr)
\times
 \nbigmtilde^{ss}(\vecyhat_2,\alpha_{\ast},+)$.
Recall that we constructed the obstruction theory
$\Obtilde(m,\vecy_1,L)$ of
$\nbigmtilde^{ss}\bigl(\vecy_1,L,\alpha_{\ast},(\delta,k_0)\bigr)$
(the subsubsection \ref{subsubsection;06.6.13.50})
and the obstruction theory
$\Obtilde(m,\vecy_2)$
of $\nbigmtilde^{ss}(\vecyhat_2,\alpha_{\ast},+)$
(the subsubsection \ref{subsubsection;06.6.13.51}).
The direct sum $\Ob(\nbigmtilde_{\spl})$ gives
the obstruction theory of $\nbigmtilde_{\spl}$.
The affine line $\Spec k[t]$ is denoted by $A^1$.

\begin{prop}
\label{prop;06.5.22.100}
We have the obstruction theory
$\ob(\Mhat^{G_m}(\gbigi)):\Ob(\Mhat^{G_m}(\gbigi))
\lrarr L_{\Mhat^{G_m}(\gbigi)}$
and the deformation
$\obtilde(\Mhat^{G_m}(\gbigi)):
 \Obtilde(\Mhat^{G_m}(\gbigi))\lrarr L_{\Mhat^{G_m}(\gbigi)\times A^1/A^1}$
with the following property:
\begin{itemize}
\item
We have the following commutative diagram:
\begin{equation}
 \label{eq;06.5.22.151}
 \begin{CD}
 \varphi_{\gbigi}^{\ast}\Ob(\Mhat)
 @>>>
 \varphi_{\gbigi}^{\ast}L_{\Mhat}\\
 @VVV @VVV \\
 \Ob(\Mhat^{G_m}(\gbigi))@>>>
 L_{\Mhat^{G_m}(\gbigi)}
 \end{CD}
\end{equation}
\item
 Let
 $\obtilde_a(\Mhat^{G_m}(\gbigi)):
 \Obtilde_a(\Mhat^{G_m}(\gbigi))
 \lrarr L_{\Mhat^{G_m}(\gbigi)}$
 denote the specialization
 of $\obtilde(\Mhat^{G_m}(\gbigi))$
 at $t=a$.
 At $t=1$,
 we have $\obtilde_1(\Mhat^{G_m}(\gbigi))= \ob(\Mhat^{G_m}(\gbigi))$.
 At $t=0$, we have the following commutative diagram
 in the diagram {\rm(\ref{eq;06.5.17.21})}:
\begin{equation}
\label{eq;06.5.22.160}
 \begin{CD}
 F^{\ast}\Obtilde_0(\Mhat^{G_m}(\gbigi)) @>>>
 F^{\ast}L_{\Mhat^{G_m}(\gbigi)}\\
 @V{\simeq}VV @V{\simeq}VV \\
 G^{\prime\ast}\Ob(\nbigmtilde_{\spl})
 @>>>
 G^{\prime\ast}L_{\nbigmtilde_{\spl}}
\end{CD}
\end{equation}
On each $a$, we have the following distinguished triangle:
\begin{equation}
\label{eq;06.5.22.165}
\begin{CD}
 G^{\prime\,\ast}\Obtilde(m,\vecy_1,L)
 @>>>
 F^{\ast}\Obtilde_a(\Mhat^{G_m}(\gbigi))
 @>>>
 G^{\prime\,\ast}\Obtilde(m,\vecyhat_2)
 @>>>
 G^{\prime\,\ast}\Obtilde(m,\vecy_1,L)[1]
\end{CD}
\end{equation}
\end{itemize}
\end{prop}

We will prove Proposition \ref{prop;06.5.22.100}
in the subsubsection \ref{subsubsection;06.5.22.123},
after some preparation.

\subsubsection{The moduli stack of split objects with an orientation}
\label{subsubsection;06.5.22.350}

We put $\nbigm_1:=\nbigm(m,\vecy_1,L)$
and $\nbigm_2:=\nbigm(m,\vecy_2)$.
We put $\nbigm_3:=\nbigm_1\times\nbigm_2$.

Let us consider the moduli stack $\nbigmhat_3$
of the objects
$(E_1,F_{1\,\ast},\phi,E_2,F_{2\,\ast},\rho)$
as follows:
\begin{itemize}
\item
 $(E_1,F_{1\,\ast},\phi)\in\nbigm_1$
and
 $(E_2,F_{2\,\ast})\in\nbigm_2$.
\item
$\rho$ denotes an orientation of $E_1\oplus E_2$.
\end{itemize}

We have the obstruction theory
$\Ob(\nbigm_3):=\Ob(m,\vecy_1,L)\oplus
 \Ob(m,\vecy_2)$ of $\nbigm_3$.
The relative obstruction theory 
of $\nbigmhat_3$
over $\nbigm_3$ is constructed in
the standard manner,
which we explain in the following.
Let $\pi:\nbigmhat_3\lrarr\nbigm_3$ denote
the projection.
We have the universal objects
$(\nbige_1^u,F_{1\ast}^u,\phi^u)$
over $\nbigm_1\times X$
and $(\nbige_2^u,F_{2\ast}^u)$
over $\nbigm_2\times X$.
We also have the canonical resolutions
$\nbigv^{(i)}_{\cdot}$ of $\nbige^u_i(m)$.
We denote the induced objects
over $\nbigmhat_3\times X$ by the same notation.
Then we have the orientation
of $\nbige_1^u\oplus\nbige_2^u$
over $\nbigmhat_3$.
Therefore, we obtain the following diagram:
\[
 \begin{CD}
 \gminig^d(\nbigv^{(1)}_{\cdot}\oplus\nbigv^{(2)}_{\cdot})
 @>>>
 \gminig(\nbigv^{(1)}_{\cdot})\oplus
 \gminig(\nbigv^{(2)}_{\cdot})\\
 @VVV @VVV \\
 \det_{\nbige^u_1\oplus\nbige_2^u,X}^{\ast}
 L_{\Pic\times X/X}
@>>> 
 L_{\nbigmhat_3\times X/X}
 \end{CD}
\]
Therefore, we obtain the following:
\[
 \begin{CD}
 \Ob^d(\nbigv^{(1)}_{\cdot}\oplus 
 \nbigv^{(2)}_{\cdot}) @>>>
 \Ob(m,\vecy_1,L)\oplus\Ob(m,\vecy_2) \\
 @VVV @VVV \\
 \det_{\nbige_1^u\oplus\nbige_2^u}^{\ast}
 L_{\Pic}@>>> L_{\nbigmhat_3}
 \end{CD}
\]
We put
$\Ob_{\rel}(\nbigmhat_3/\nbigm_3):=\Cone\bigl(
 \Ob^d(\nbigv_{\cdot}^{(1)}\oplus 
 \nbigv_{\cdot}^{(2)})\lrarr 
 \det_{\nbige_1^u\oplus\nbige_2^u}^{\ast}
 L_{\Pic} \bigr)$.
We have the natural morphism:
\[
\gamma:\Ob_{\rel}(\nbigmhat_3/\nbigm_3)[-1]
 \lrarr 
\Ob(m,\vecy_1,L)\oplus\Ob(m,\vecy_2)
=\pi^{\ast}\Ob(\nbigm_3)
\]
The cone of $\gamma$
is denoted by $\Ob(\nbigmhat_3)$.
Then we have the natural morphism
$\ob(\nbigmhat_3):\Ob(\nbigmhat_3)\lrarr L_{\nbigmhat_3}$.

We have the following commutative diagram
as in the subsubsection \ref{subsubsection;06.5.5.10}:
\[
\begin{CD}
 L_{\nbigmhat_3}
 @<<<
 \det^{\ast}_{\nbige^u_1\oplus\nbige_2^u}L_{\Pic} \\
 @AAA @AAA \\
 \pi^{\ast}L_{\nbigm_3}
 @<<<
\Phi\bigl(\det(\nbige^u_1\oplus\nbige_2^u)
 \bigr)^{\ast}
 \pi_1^{\ast}
 L_{\nbigm(1)}
 @<<<
 \Ob^d(\nbigv_{\cdot}^{(1)}\oplus\nbigv_{\cdot}^{(2)})
\end{CD}
\]
Here, $\pi_1$ denotes the projection
$\Pic\lrarr\nbigm(1)$.
Therefore, we have the following commutative diagram:
\[
 \begin{CD}
 \Ob_{\rel}(\nbigmhat_3/\nbigm_3)[-1]
 @>{\gamma}>>
 \pi^{\ast}\Ob(\nbigm_3) \\
 @V{\ob_{\rel}(\nbigmhat_3/\nbigm_3)}VV @VVV\\
 L_{\nbigmhat_3/\nbigm_3}[-1]
 @>>>
 \pi^{\ast}L_{\nbigm_3}
 \end{CD}
\]
By the same argument as the proof of
Lemma \ref{lem;06.5.3.42},
it can be shown that
$\ob_{\rel}(\nbigmhat_3/\nbigm_3)$ gives a relative
obstruction theory of $\nbigmhat_3$ over $\nbigm_3$.
Therefore,
$\Ob(\nbigmhat_3)$ gives an obstruction theory of
$\nbigmhat_3$.

\subsubsection{The embedding into the moduli stack of non-split objects}
\label{subsubsection;06.5.23.5}

We put $\nbigm_0:=\nbigm(m,\vecy,L)$
and $\nbigmhat_0:=\nbigm(m,\vecyhat,L)$.
The projection $\nbigmhat_0\lrarr\nbigm_0$ is denoted by $\pi_0$
Let $\nbigv_{\cdot}$ denote the canonical resolution
of $p_{X\,\ast}\nbige^u(m)$ on $\nbigm_0\times X$.
We put $\nbigv_0':=p_{X\,\ast}\nbigv_{0}=p_{X\,\ast}\nbige^u(m)$.
We have the naturally defined morphism
$\gminif:\nbigm_3\lrarr \nbigm_0$.
We have the decomposition
$\gminif_X^{\ast}\nbige^u=
 \nbige^u_1\oplus\nbige_2^u$,
$\gminif_X^{\ast}\nbigv_{\cdot}=
 \nbigv_{\cdot}^{(1)}\oplus\nbigv_{\cdot}^{(2)}$
and $\gminif^{\ast}\nbigv_0^{\prime}=
 \nbigv_0^{\prime\,(1)}\oplus \nbigv_0^{\prime\,(2)}$.
We have the naturally defined projections:
\[
 \gminif_X^{\ast}\gminig(\nbigv_{\cdot})
\lrarr
 \gminig(\nbigv_{\cdot}^{(1)})\oplus
 \gminig(\nbigv_{\cdot}^{(2)}),
\quad
\gminif_X^{\ast}\gminig_{\rel}(\nbigv_{\cdot},\phitilde)
\lrarr
 \gminig_{\rel}(\nbigv^{(1)}_{\cdot},\phitilde),
\]
\[
\gminif_X^{\ast}\gminig(\nbigv_{\cdot|D})
\lrarr
 \gminig(\nbigv^{(1)}_{\cdot|D})\oplus
 \gminig(\nbigv^{(2)}_{\cdot|D}),
\quad
\gminif_X^{\ast}\gminig_{\rel}(\nbigv_{\cdot},F_{\ast}^u)
\lrarr
 \gminig_{\rel}(\nbigv^{(1)}_{\cdot},F_{\ast}^{u\,(1)})
\oplus\gminig_{\rel}(\nbigv^{(2)}_{\cdot},F_{\ast}^{u\,(2)}). 
\]
They induce the following morphisms:
\[
 \gminif^{\ast}\Ob(\nbigv_{\cdot})
\lrarr
 \Ob(\nbigv^{(1)}_{\cdot})\oplus \Ob(\nbigv^{(2)}_{\cdot}),
\quad
 \gminif^{\ast}\Ob_{\rel}(\nbigv_{\cdot},\phitilde)
\lrarr
 \Ob_{\rel}(\nbigv^{(1)},\phitilde)
\]
\[
 \gminif^{\ast}\Ob(\nbigv_{\cdot|D})
\lrarr
 \Ob(\nbigv^{(1)}_{\cdot|D})\oplus
 \Ob(\nbigv^{(2)}_{\cdot|D}),
\quad
 \gminif^{\ast}\Ob(\nbigv_{\cdot},F^u_{\ast})
\lrarr
 \Ob_{\rel}(\nbigv^{(1)}_{\cdot},F^{u\,(1)}_{\ast})
\oplus
 \Ob_{\rel}(\nbigv^{(2)}_{\cdot},F^{u\,(2)}_{\ast})
\]
Therefore, we obtain the morphism:
\[
 \mu_1:
 \gminif^{\ast}\Ob(m,\vecy,L)
\lrarr 
 \Ob(m,\vecy_1,L)\oplus \Ob(m,\vecy_2)
\]

\begin{lem}
\label{lem;06.5.22.16}
The following diagram is commutative:
\begin{equation}
\label{eq;06.5.22.15}
 \begin{CD}
 \gminif^{\ast}\Ob(m,\vecy,L)
 @>>>\gminif^{\ast}L_{\nbigm_0}\\
 @VVV @VVV \\
 \Ob(m,\vecy_1,L)\oplus
 \Ob(m,\vecy_2)
 @>>>
 L_{\nbigm_3}
 \end{CD}
\end{equation}
\end{lem}
\pf
We take an $H_{\vecy}(m)$-dimensional vector space
$W_0=V_m$ with a decomposition
$W_0=W_0^{(1)}\oplus W_0^{(2)}$,
where $\dim W_0^{(i)}=H_{\vecy_i}(m)$.
We also take a $\bigl(H_{\vecy}(m)-\rank(\vecy)\bigr)$-dimensional
vector space $W_{-1}$ with a decomposition
$W_{-1}=W_{-1}^{(1)}\oplus W_{-1}^{(2)}$,
where $\dim W_{-1}^{(i)}=H_{\vecy_i}(m)-\rank(\vecy_i)$.
We put 
$Y(W^{(1)}_{\cdot},W_{\cdot}^{(2)}):=
 Y(W^{(1)}_{\cdot})\times Y(W^{(2)}_{\cdot})$.
We have the naturally defined morphism
$Y(W^{(1)}_{\cdot},W^{(2)}_{\cdot}) \lrarr Y(W_{\cdot})$.
By considering the classifying map of 
$\nbigv_{\cdot}$ and $\nbigv_{\cdot}^{(1)}\oplus\nbigv_{\cdot}^{(2)}$,
we obtain the following commutative diagram:
\[
 \begin{CD}
 \nbigm_0\times X
 @>{\Phi(\nbigv_{\cdot})}>> 
 Y(W_{\cdot})\\
 @A{\gminif_X}AA @AAA\\
 \nbigm_3\times X 
 @>{\Phi(\nbigv_{\cdot}^{(1)},\nbigv_{\cdot}^{(2)})}>> 
 Y(W_{\cdot}^{(1)},W_{\cdot}^{(2)})
 \end{CD}
\]
By the argument in the subsubsection
\ref{subsubsection;06.4.29.15},
we can show that 
$\gminif_X^{\ast}\gminig(\nbigv_{\cdot})_{\leq 1}
\lrarr
 \gminig(\nbigv_{\cdot}^{(1)})_{\leq 1}\oplus 
 \gminig(\nbigv_{\cdot}^{(2)})_{\leq 1}$
expresses
the morphism
$\gminif_X^{\ast}\Phi(\nbigv_{\cdot})^{\ast}L_{Y(W_{\cdot})/X}
\lrarr
 \Phi(\nbigv_{\cdot}^{(1)},\nbigv_{\cdot}^{(2)})^{\ast}
L_{Y(W_{\cdot}^{(1)},W_{\cdot}^{(2)})/X}$.
Therefore, we obtain the following commutative diagram:
\[
 \begin{CD}
 \gminif_X^{\ast}L_{\nbigm_0\times X/X}
 @<<<
 \gminif_X^{\ast}\Phi(\nbigv_{\cdot})^{\ast}
 L_{Y(W_{\cdot})/X}
 @<<<
 \gminig(\nbigv_{\cdot})\\
 @VVV @VVV @VVV \\
 L_{\nbigm_3\times X/X}
 @<<<
 \Phi(\nbigv_{\cdot}^{(1)},\nbigv_{\cdot}^{(2)})^{\ast}
 L_{Y(W_{\cdot}^{(1)},W_{\cdot}^{(2)})/X}
 @<<<
 \gminig(\nbigv_{\cdot}^{(1)}) \oplus
 \gminig(\nbigv_{\cdot}^{(2)})
 \end{CD}
\]
Therefore, we obtain the following commutative diagram:
\begin{equation}
\label{eq;06.5.22.11}
 \begin{CD}
 \gminif^{\ast}L_{\nbigm_0} @<<<
 \gminif^{\ast}\Ob(m,y)\\
 @VVV @VVV \\
 L_{\nbigm_3} @<<<
 \Ob(m,y_1)\oplus \Ob(m,y_2).
 \end{CD}
\end{equation}

Let $P_{\cdot}$ be a locally free resolution of $L(m)$,
as in the subsubsection \ref{subsubsection;06.5.22.5}.
We put $Y(W_{\cdot}^{(1)},W_{\cdot}^{(2)},P_{\cdot}):=
 Y(W^{(1)}_{\cdot},P_{\cdot})\times Y(W_{\cdot}^{(2)})$.
Then, we have the naturally defined morphism
$Y(W_{\cdot}^{(1)},W_{\cdot}^{(2)},P_{\cdot})
\lrarr Y(W_{\cdot},P_{\cdot})$.
Then we have the following diagram:
\[
 \begin{CD}
 \nbigm_0\times X@>{\Phi(\nbigv_{\cdot},\phitilde)}>> 
 Y(W_{\cdot},P_{\cdot})\\
 @A{\gminif_X}AA @AAA \\
 \nbigm_3\times X 
 @>{\Phi(\nbigv_{\cdot}^{(1)},\nbigv_{\cdot}^{(2)},\phitilde)}>>
 Y(W_{\cdot}^{(1)},W_{\cdot}^{(2)},P_{\cdot})
 \end{CD}
\]
Then, we obtain the following commutative diagram:
\[
 \begin{CD}
 \gminif_X^{\ast}L_{\nbigm_0\times X/X}
 @<<<
 \gminif_X^{\ast}\Phi(\nbigv_{\cdot},\phitilde)^{\ast}
 L_{Y(W_{\cdot},P_{\cdot})/X}
 @<<<
 \gminif_X^{\ast}\gminig(\nbigv_{\cdot},\phitilde)\\
 @VVV @VVV @VVV \\
 L_{\nbigm_3\times X/X}
 @<<<
 \Phi(\nbigv_{\cdot}^{(1)},\nbigv_{\cdot}^{(2)},\phitilde)^{\ast}
 L_{Y(W_{\cdot}^{(1)},W_{\cdot}^{(2)},P_{\cdot})}
 @<<<
 \gminig(\nbigv_{\cdot}^{(1)},\phitilde)
 \oplus
 \gminig(\nbigv_{\cdot}^{(2)})
 \end{CD}
\]
It is easy to observe that
$Rp_{X\,\ast}\gminig(\nbigv_{\cdot},\phitilde)\otimes\omega_X$
is naturally isomorphic to $\Ob(m,y,L)$.
Then, we obtain the following diagram:
\begin{equation}
 \label{eq;06.5.22.10}
 \begin{CD}
 \gminif^{\ast}L_{\nbigm_0} @<<< \gminif^{\ast}\Ob(m,y,L)\\
 @VVV @VVV\\
 L_{\nbigm_3} @<<<
 \Ob(m,y_1,L)\oplus \Ob(m,y_2)
 \end{CD}
\end{equation}
We have the natural morphism $\Ob(m,y)\lrarr\Ob(m,y,L)$
and $\Ob(m,y_1)\oplus \Ob(m,y_2)\lrarr \Ob(m,y_1,L)\oplus \Ob(m,y_2)$.
The diagrams (\ref{eq;06.5.22.11}) and (\ref{eq;06.5.22.10})
are compatible for the natural morphisms
in the sense of the subsubsection \ref{subsubsection;06.5.22.15}.

\vspace{.1in}
We put $\nbigv^{(i)(j)}_D:=
 \Ker(\nbigv_{0|D}^{(i)})\lrarr \Cok^{(i)}_{j-1}$
which are locally free sheaves on $\nbigm_3\times D$.
Let $\nbigv_D^{(i)\ast}$ denote the vector bundle
$\nbigv_D^{(i)}$ with the filtration
$\nbigv_{D}^{(i)(1)}\supset\nbigv_D^{(i)(2)}
 \cdots \supset \nbigv_D^{(i)(l+1)}$.
Similarly,
we have the filtered vector bundle $\nbigv_D^{\ast}$
on $\nbigm_0\times D$.

We put $W^{(i)(1)}:=W_0^{(i)}$ and
$W^{(i)(l+1)}=W_{-1}^{(i)}$.
We take vector spaces $W^{(i)}$ $(i=2,\ldots,l)$
with decompositions
$W^{(i)}=W^{(i)(1)}\oplus W^{(i)(2)}$,
where $\rank W^{(i)}=\rank \nbigv_D^{(i)}$
and $\rank W^{(i)(j)}=\rank\nbigv_D^{(i)(j)}$.
We use the notation in the subsubsection
\ref{subsubsection;06.5.5.5}.
We put $Y_D(W_{\cdot},W^{(1)\,\ast},W^{(2)}):=
 Y_D(W^{(1)}_{\cdot},W^{(1)\ast})
\times
 Y_D(W^{(2)}_{\cdot},W^{(2)\ast})$,
and $Y_D(W_{\cdot}^{(1)},W_{\cdot}^{(2)}):=
 Y_D(W_{\cdot}^{(1)})\times Y_D(W_{\cdot}^{(2)})$.
We have the naturally defined commutative diagram:
\[
 \begin{CD}
 Y_D(W_{\cdot},W^{\ast}) @>>> Y_D(W_{\cdot}) \\
 @AAA @AAA \\
 Y_D(W_{\cdot},W^{(1)\ast},W^{(2)\ast})
 @>>>
 Y_D(W_{\cdot}^{(1)},W_{\cdot}^{(2)})
 \end{CD}
\]
By considering the classifying maps
of $\nbigv_D^{\ast}$ and 
$(\nbigv_D^{(1)\ast}, \nbigv_D^{(2)\ast})$,
we obtain the following commutative diagram:
\[
 \begin{CD}
 \nbigm_0\times D
 @>{\Phi_D(\nbigv_{\cdot}^{\ast},F_{\ast})}>> 
 Y_D(W_{\cdot},W^{\ast}) @>>>
 Y_D(W_{\cdot})\\
 @AAA @AAA @AAA \\
 \nbigm_3\times D 
 @>>>\Phi_1
 Y_D(W_{\cdot},W^{(1)\ast},W^{(2)\ast})
 @>>>
 Y_D(W_{\cdot}^{(1)},W_{\cdot}^{(2)})
 \end{CD}
\]
Therefore, we obtain the following:
\[
 \begin{CD}
 \gminif_D^{\ast}L_{\nbigm_0\times D/D} @<<<
 \gminif_D^{\ast}\Phi_D(\nbigv_{\cdot},F_{\ast}^u)^{\ast}
 L_{Y_D(W_{\cdot},W^{\ast})/D}
 @<<<
 \gminif_D^{\ast}\Phi(\nbigv_{\cdot|D})^{\ast}
 L_{Y_D(W_{\cdot})/D}\\
 @VVV @VVV @VVV \\
 L_{\nbigm_3\times D/D} @<<< 
 \Phi_1^{\ast}
 L_{Y_D(W_{\cdot},W^{(1)\ast},W^{(2)\ast})/D}
 @<<<
 \Phi(\nbigv_{\cdot|D}^{(1)}\oplus\nbigv_{\cdot|D}^{(2)})^{\ast}
 L_{Y_D(W_{\cdot}^{(1)},W_{\cdot}^{(2)})/D}
 \end{CD}
\]
Then, we obtain the following:
\[
 \begin{CD}
 \gminif_D^{\ast}L_{\nbigm_0\times D/D}
 @<<< 
 \gminif_D^{\ast}\gminig(\nbigv_{\cdot},F^u_{\ast}) 
 @<<<
 \gminif_D^{\ast}\gminig(\nbigv_{\cdot|D})\\
 @VVV @VVV @VVV\\
 L_{\nbigm_3\times D/D}
 @<<<
 \gminig_D(\nbigv^{(1)}_{\cdot},F^{u(1)}_{\ast})
 \oplus
 \gminig_D(\nbigv^{(2)}_{\cdot},F^{u(2)}_{\ast})
 @<<<
 \gminig(\nbigv^{(1)}_{\cdot|D})
\oplus
 \gminig(\nbigv^{(2)}_{\cdot|D})
 \end{CD}
\]
We put 
$\Ob_D(\nbigv,F^u_{\ast}):=
 Rp_{D\,\ast}\bigl(
 \gminig(\nbigv_{\cdot},F^u_{\ast})\otimes\omega_D
 \bigr)$.
Similarly,
we obtain
$\Ob_D(\nbigv^{(i)},F^{u(i)}_{\ast})$.
Then, we obtain the following commutative diagram:
\[
 \begin{CD}
 \gminif^{\ast}L_{\nbigm_0} @<<<
 \gminif^{\ast}\Ob_D(\nbigv_{\cdot},F^{u}_{\ast})
 @<<< \gminif^{\ast}\Ob(\nbigv_{\cdot|D})\\
 @VVV @VVV @VVV\\
 L_{\nbigm_3} @<<<
 \Ob_D(\nbigv_{\cdot}^{(1)},F^{u(1)}_{\ast})
\oplus 
 \Ob_D(\nbigv_{\cdot}^{(2)},F^{u(2)}_{\ast})
 @<<<
 \Ob(\nbigv^{(1)}_{\cdot|D})
\oplus
 \Ob(\nbigv^{(2)}_{\cdot|D})
 \end{CD}
\]
We remark that the cone of
$\Ob(\nbigv_{\cdot|D})\lrarr
 \Ob(\nbigv_{\cdot})\oplus
 \Ob_D(\nbigv,F_{\ast})$ is naturally isomorphic 
to $\Ob(m,\vecy)$.
Thus, we obtain the following commutative diagram,
which is compatible with (\ref{eq;06.5.22.11}):
\begin{equation}
\label{eq;06.5.22.12}
\begin{CD}
 \gminif^{\ast}L_{\nbigm_0}
 @<<<
 \gminif^{\ast}\Ob(m,\vecy)\\
 @VVV @VVV \\
 L_{\nbigm_3} @<<<
 \Ob(m,\vecy_1)\oplus
 \Ob(m,\vecy_2)
\end{CD}
\end{equation}
From (\ref{eq;06.5.22.11}),
(\ref{eq;06.5.22.10}) and (\ref{eq;06.5.22.12}),
we obtain the desired diagram
(\ref{eq;06.5.22.15}).
Thus the proof of Lemma \ref{lem;06.5.22.16}
is finished.
\hfill\qed

\vspace{.1in}
We have the naturally defined morphism
$\gminifhat:\nbigmhat_3\lrarr \nbigmhat_0$.
By construction of the obstruction theories,
we have the following commutative diagram:
\[
 \begin{CD}
 \gminifhat^{\ast}\Ob_{\rel}(m,\yhat)[-1] @>>> 
 \gminifhat^{\ast}\pi_0^{\ast}\Ob(m,\vecy,L) \\
 @V{\simeq}VV @V{\mu_1}VV \\
 \Ob_{\rel}(\nbigmhat_3/\nbigm_3)[-1] @>>>
 \pi^{\ast}\Ob(m,\vecy_1,L)\oplus\Ob(m,\vecy_2)
 \end{CD}
\]
Therefore, we obtain the following commutative diagram:
\[
 \begin{CD}
 \gminifhat^{\ast}\Ob(m,\vecyhat,L)
 @>{\mu_2}>>
 \Ob(\nbigmhat_3) \\
 @VVV @VVV \\
\gminifhat^{\ast}L_{\nbigmhat_0}
@>>>
 L_{\nbigmhat_3}
 \end{CD}
\]

\subsubsection{Some compatibility}
\label{subsubsection;06.5.23.2}

We take $H_{\vecy_i}(m)$-dimensional vector spaces $W_0^{(i)}$.
We put $B(W_0^{(1)},W_0^{(2)}):=
 k_{\GL(W_0^{(1)})\times \GL(W_0^{(2)})}$.
We put $\nbigv_0^{\prime\,(i)}:=p_{X\,\ast}\nbigv_0^{(i)}$.
Then we have the classifying map
$\Phi(\nbigv_0^{\prime(1)}\oplus\nbigv_0^{\prime(2)}):
 \nbigm_3\lrarr B(W_0^{(1)},W_0^{(2)})$.
Therefore, we have the morphism
$\varphi:\Phi(\nbigv_0^{\prime(1)}\oplus\nbigv_0^{\prime(2)})^{\ast}
 L_{B(W_0^{(1)},W_0^{(2)})}
\lrarr L_{\nbigm_3}$.

\begin{lem}
\label{lem;06.5.22.150}
The morphism $\varphi$ factors through
$\Ob(m,y_1,L)\oplus \Ob(m,y_2)$.
In particular, it factors through
$\Ob(m,\vecy_1,L)\oplus\Ob(m,\vecy_2)$.
\end{lem}
\pf
We have the following commutative diagram:
\[
\begin{CD}
 \gminig(\nbigv_{\cdot}^{(1)})\oplus\gminig(\nbigv_{\cdot}^{(2)})
 @<<<
 \gminih(\nbigv_{\cdot}^{(1)})\oplus\gminih(\nbigv_{\cdot}^{(2)}) \\
 @VVV @VVV \\
 L_{\nbigm_3\times X}
 @<<<
 \Phi(\nbigv_0^{\prime(1)},\nbigv_0^{\prime(2)})^{\ast}
 L_{B(W_0^{(1)},W_0^{(2)})\times X/X}
\end{CD}
\]
It induces the following diagram:
\[
 \begin{CD}
 \Ob(m,y_1)\oplus\Ob(m,y_2) @<<<
 \Ob^{G}(\nbigv^{(1)}_{\cdot})\oplus
 \Ob^{G}(\nbigv^{(2)}_{\cdot})\\
 @VVV @V{\varphi_1}VV \\
 L_{\nbigm_3}@<<< 
 \Phi(\nbigv_0^{\prime\,(1)},
 \nbigv_0^{\prime\,(2)})^{\ast}
 L_{B(W_0^{(1)},W_0^{(2)})}
 \end{CD}
\]
It is easy to check $\varphi_1$ is isomorphic.
Then the claim of the lemma is clear.
\hfill\qed

\vspace{.1in}

We put $W_0:=W_0^{(1)}\oplus W_0^{(2)}$
and $B(W_0):=k_{\GL(W_0)}$.
By considering the classifying maps
of the vector bundles
$\nbigv_0$ and $\nbigv^{(i)}_0$,
we obtain the following commutative diagram:
\begin{equation}
\label{eq;06.5.22.2}
 \begin{CD}
\nbigmhat_3@>{\pi}>>
 \nbigm_3 @>{\Phi(\nbigv_0^{(1)},\nbigv_0^{(2)})}>>
  B(W_0^{(1)},W_0^{(2)})
 \\
@V{\gminifhat}VV @V{\gminif}VV @VVV \\
\nbigmhat_0@>{\pi_2}>>
 \nbigm_0 @>{\Phi(\nbigv_0)}>>
 B(W_0)
 \end{CD}
\end{equation}

\begin{lem}
We have the following commutative diagram:
\begin{equation}
\label{eq;06.5.11.10}
 \begin{CD}
 L_{\nbigmhat_3} @<<<
 \Ob(\nbigmhat_3) @<<<
\pi^{\ast} 
\Phi(\nbigv_0^{(1)},\nbigv_0^{(2)})^{\ast}
 L_{B(W_0^{(1)},W_0^{(2)})}
\\
 @AAA @A{\mu_2}AA @AAA \\
 \gminifhat^{\ast}L_{\nbigmhat_0} 
 @<<<
 \gminifhat^{\ast}\Ob(m,\vecyhat,L) @<<<
 \gminifhat^{\ast}\pi_2^{\ast}
 \Phi(\nbigv_0)^{\ast}
 L_{B(W_0)}
 \end{CD}
\end{equation}
Here, the composite of the horizontal arrows
are the naturally defined morphisms
from the diagram {\rm(\ref{eq;06.5.22.2})}.
\end{lem}
\pf
We have the following commutative diagram:
\[
 \begin{CD}
 \gminif_X^{\ast}\gminih(\nbigv_{\cdot})
 @>>>
 \gminif_X^{\ast}\gminig(\nbigv_{\cdot})\\
 @VVV @VVV \\
 \gminih(\nbigv^{(1)}_{\cdot})
\oplus
 \gminih(\nbigv^{(2)}_{\cdot})
 @>>>
 \gminig(\nbigv^{(1)}_{\cdot})
\oplus
 \gminig(\nbigv^{(2)}_{\cdot})
 \end{CD}
\]
Therefore, we obtain the following:
\[
 \begin{CD}
 \gminif^{\ast}\Ob^{G}(\nbigv_{\cdot})
 @>>>
 \gminif^{\ast}\Ob(m,y)\\
 @VVV @VVV \\
 \Ob^G(\nbigv_{\cdot}^{(1)})
\oplus
 \Ob^{G}(\nbigv_{\cdot}^{(2)})
 @>>>
 \Ob(m,y_1)\oplus\Ob(m,y_2)
 \end{CD}
\]
On the other hand,
we have the following commutative diagram:
\[
\begin{CD}
 \gminif_X^{\ast}\gminih(\nbigv_{\cdot})
 @>>>
 \gminif_X^{\ast}\Phi(\nbigv_0)^{\ast}_X
 L_{B(W_0)\times X/X}\\
 @VVV @VVV \\
 \gminih(\nbigv_{\cdot}^{(1)})
\oplus
 \gminih(\nbigv_{\cdot}^{(2)})
 @>>>
 \Phi(\nbigv_0^{\prime\,(1)},
 \nbigv_0^{\prime\,(2)})^{\ast}_X
 L_{B(V^{(1)}_m,V^{(2)}_m)\times X/X}
\end{CD}
\]
Therefore, we obtain the following:
\[
 \begin{CD}
 \gminif^{\ast}\Ob^G(\nbigv_{\cdot})
 @>{\tau_1}>>
 \gminif^{\ast}\Phi(\nbigv_{0})^{\ast}L_{B(W_0)}\\
 @VVV @VVV \\
 \Ob^G(\nbigv_{\cdot}^{(1)})
\oplus
 \Ob^G(\nbigv_{\cdot}^{(2)})
 @>{\tau_2}>>
 \Phi(\nbigv_0^{\prime\,(1)},
 \nbigv_0^{\prime\,(2)})^{\ast}
 L_{B(W_0^{(1)},W_0^{(2)})}
 \end{CD}
\]
The morphisms $\tau_i$ are isomorphic.
Thus we obtain the claim of the lemma.
\hfill\qed

\subsubsection{Deformation}
\label{subsubsection;06.5.22.360}

As explained in the subsubsection
\ref{subsubsection;06.5.11.200},
we have the decomposition
$\Ob(m,\vecy_2)=\Ob^{\circ}(m,\vecy_2)\oplus
 \Ob^d(m,\vecy_2)$.
We put
$\Obcheck(m,\vecy_2):=
 \Ob^{\circ}(m,\vecy_2)\oplus
 \tau_{\leq 0}\Ob^d(m,\vecy_2)$.
We have the following commutative diagram:
\[
 \begin{CD}
 \Ob_{\rel}(\nbigv_{\cdot}^{(1)}\oplus\nbigv_{\cdot}^{(2)})[-1]
 @>>>
 \Ob(m,\vecy_1,L)\oplus \Ob(m,\vecy_2)\\
 @AAA @AAA \\
 \tau_{\leq -1}
 \Ob^d(\nbigv_{\cdot}^{(1)}\oplus\nbigv_{\cdot}^{(2)})
 @>{\lambda_1}>>
 \Ob(m,\vecy_1,L)
\oplus
 \Obcheck(m,\vecy_2)
 \end{CD}
\]
We put $\Ob_1(\nbigmhat_3):=\Cone\bigl(\lambda_1\bigr)$.
Then we have the morphisms
$\Ob_1(\nbigmhat_3)\lrarr\Ob(\nbigmhat_3)\lrarr L_{\nbigmhat_3}$.
Since the first morphism is quasi isomorphic,
the composite $\ob_1(\nbigmhat_3)$
of the morphisms gives an obstruction theory
of $\nbigmhat_3$.
We remark we have the following commutative diagram
(We put
$C_1:=\tau_{\leq -1}\Ob^d(\nbigv_{\cdot}^{(1)}
 \oplus\nbigv^{(2)}_{\cdot})$
and $C_2:=\Ob(m,\vecy_1,L)\oplus\Obcheck(m,\vecy_2)$
in the diagram, to save the space.):
\begin{equation}
 \label{eq;06.6.22.3}
 \begin{CD}
 \nbigh^{-1}\bigl(C_1 \bigr)
 @>>>
 \nbigh^{-1}(C_2)
 @>>>
 \nbigh^{-1}\bigl(\Ob(\nbigmhat_3)\bigr)
 @>>> 0\\
 @VVV @VVV @VVV @VVV\\
 0 @>>>
 \nbigh^{-1}(\pi^{\ast}L_{\nbigm_3})
 @>>>
 \nbigh^{-1}\bigl( L_{\nbigmhat_3}\bigr)
 @>>>0
 \end{CD}
\end{equation}
\begin{equation}
 \label{eq;06.6.22.4}
 \begin{CD}
\nbigh^0(C_2)
 @>{\simeq}>>
 \nbigh^0\bigl(\Ob_1(\nbigmhat_3)\bigr)
 @>>> 0 @>>>
 \nbigh^1(C_2) @>>>
 \nbigh^1(\Ob(\nbigmhat_3))\\
 @VVV @V{\simeq}VV @VVV 
 @VVV @V{\simeq}VV\\
 \nbigh^0(\pi^{\ast}L_{\nbigm_3})
 @>>>
 \nbigh^0(L_{\nbigmhat_3})
 @>>>
 \nbigh^0(L_{\nbigmhat_3/\nbigm_3})
 @>{\simeq}>>
 \nbigh^1(\pi^{\ast}L_{\nbigm_3})
 @>>>
 \nbigh^1(L_{\nbigmhat_3})
 \end{CD}
\end{equation}

We would like to deform $\ob_1(\nbigmhat_3)$.
Let $i_1$, $i_2$ and $\eta$ denote the following naturally defined
morphisms:
\[
 i_1:\tau_{\leq -1}\Ob^d(\nbigv_1)\lrarr \Ob(m,\vecy_1,L),
\quad
 i_2:\tau_{\leq -1}\Ob^d(\nbigv_2)\lrarr\Obcheck(m,\vecy_2),
\]
\[
  \eta:\tau_{\leq -1}\Ob^d(\nbigv_1\oplus\nbigv_2)
\lrarr \tau_{\leq -1}\Ob^{d}(\nbigv_1)\oplus
 \tau_{\leq -1}\Ob^d(\nbigv_2)
\]

The following is a special case of Lemma \ref{lem;06.5.11.20}.
\begin{lem}
 \label{lem;06.5.11.50}
The composite $i_1\circ\ob(m,\vecy_1,L)$
and $i_2\circ\ob(m,\vecy_2)$ are trivial.
\hfill\qed
\end{lem}

For any $a\in k$,
let $\varphi_a:\tau_{\leq -1}\Ob^d(\nbigv_{\cdot}^{(1)}
 \oplus\nbigv_{\cdot}^{(2)})
\lrarr \Ob(m,\vecy_1,L)\oplus\Obcheck(m,\vecy_2)$
be the morphism given by
$\varphi_a:=\bigl(a\cdot i_1,i_2\bigr)\circ \eta$.
Then the following diagram is commutative
for any $a$, due to Lemma \ref{lem;06.5.11.50}:
\begin{equation}
 \label{eq;06.6.22.2}
 \begin{CD}
 \tau_{\leq -1}\Ob^d(\nbigv_{\cdot}^{(1)}
 \oplus \nbigv_{\cdot}^{(2)})
 @>{\varphi_a}>>
 \Ob(m,\vecy_1,L)\oplus\Obcheck(m,\vecy_2)\\
 @VVV @VVV \\
 L_{\nbigmhat_3/\nbigm_3}[-1] @>>>
 \pi^{\ast} L_{\nbigm_3}
 \end{CD}
\end{equation}
We put $\Obtilde_a(\nbigmhat_3):=\Cone(\varphi_a)$.
From the commutativity of the diagram
(\ref{eq;06.6.22.2}),
we obtain the morphism
$\obtilde_a(\nbigmhat_3):\Obtilde_a(\nbigmhat_3)
\lrarr L_{\nbigmhat_3}$ for any $a$.
The choice of $a$ does not have any effect
on the diagram (\ref{eq;06.6.22.4}).
Hence, it is easy to observe that
$\bigl\{\obtilde_a(\nbigmhat_3)\,\big|\,a\in k\bigr\}$ gives
an obstruction theory of 
$\obtilde(\nbigmhat_3):\Obtilde(\nbigmtilde_3)\lrarr
 L_{\nbigmhat_3\times A^1/A^1}$
of $\nbigmhat_3\times A^1$ over $A^1$.
At $a=1$,
$\Obtilde_1(\nbigmhat_3)$ is same as
$\Ob(\nbigmhat_3)$ in the derived category.

\begin{lem}
We put $\nbigmhat_2:=\nbigm(m,\vecyhat_2)$.
There exists an algebraic stack $\nbigs'$
with the following diagram:
\[
 \begin{CD}
 \nbigmhat_3 @<{F_1}<< \nbigs' @>{G_1'}>> \nbigm_1\times\nbigmhat_2\\
 @AAA @AAA @AAA \\
 \Mhat^{G_m}(\gbigi)@<{F}<< \nbigs@>{G'}>>
 \nbigmtilde_{\spl}
 \end{CD}
\]
Here the bottom horizontal arrows are given
in {\rm (\ref{eq;06.5.17.21})}.
The morphisms $F_1$ and $G_1'$ are etale and finite
of degree $(r_1\cdot r_2)^{-1}$.
\end{lem}
\pf
Similar to Proposition \ref{prop;06.5.17.30}
and Corollary \ref{cor;06.5.22.50}.
\hfill\qed

\begin{lem}
\label{lem;06.5.22.200}
For each $a$,
we have the following distinguished triangle:
\begin{equation}
 \label{eq;06.5.22.111}
 \begin{CD}
 G_1^{\prime\,\ast} \Ob(m,\vecy_1,L)
@>>>
 F_1^{\ast}\Obtilde_a(\nbigmhat_3)
@>>>
 G_1^{\prime\,\ast}\Ob(m,\vecyhat_2)
@>>>
 G_1^{\prime\,\ast}\Ob(m,\vecy_1,L)[1]
 \end{CD}
\end{equation}
We also have
$F_1^{\ast}\Obtilde_0(\nbigmhat_3)=
 G_1^{\prime\,\ast}\Ob(m,\vecy_1,L)
\oplus
 G_1^{\prime\,\ast}\Ob(m,\vecyhat_2)$.
\end{lem}
\pf
We have the following naturally defined 
distinguished triangle on $\nbigmhat_3$:
\begin{equation}
\label{eq;06.5.22.110}
 \begin{CD}
 \Ob(m,\vecy_1,L)
@>>>
 \Obtilde_a(\nbigmhat_3)
@>>>
 \cone(i_2)
@>>>
 \Ob(m,\vecy_1,L)[1]
 \end{CD}
\end{equation}
In the case $a=0$, it splits.
It is easy to see
$F_1^{\ast}\Ob(m,\vecy_1,L)
\simeq G_1^{\prime\,\ast}\Ob(m,\vecy_1,L)$
and 
$F_1^{\ast}\cone(i_2)
\simeq \Ob(m,\vecyhat_2)$.
Hence we obtain (\ref{eq;06.5.22.111})
from (\ref{eq;06.5.22.110}).
\hfill\qed

\begin{lem}
\label{lem;06.5.22.201}
The composite of the following morphisms
is independent of the choice of $a$,
and it is same as the naturally defined one:
\begin{equation}
\label{eq;06.5.22.20}
 \pi^{\ast}
 \Phi(\nbigv_0^{(1)}\oplus\nbigv_0^{(2)})^{\ast}
 L_{B(V_m^{(1)},V_m^{(2)})}
\lrarr
 \pi^{\ast}\bigl(
 \Ob(m,\vecy_1,L)\oplus\Ob(m,\vecy_2)
 \bigr)
\lrarr
 \Ob_a(\nbigmhat_3)
\stackrel{\varphi_a}
{\lrarr}
 L_{\nbigmhat_3}
\end{equation}
\end{lem}
\pf
It is clear from the construction.
\hfill\qed

\subsubsection{Proof of Proposition \ref{prop;06.5.22.100}}
\label{subsubsection;06.5.22.123}

Let us construct an obstruction theory
of $\Mhat^{G_m}(\gbigi)$.
We take $H_{\vecy_i}(m)$-dimensional
vector spaces $V_{m}^{(i)}$.
Let $F_i$ denote the full flag variety of $V_m^{(i)}$.
We put $\Fbar_i:=F_{i\,\GL(V_m^{(i)})}$.
Then, $\Mhat^{G_m}(\gbigi)$
is an open subset of the fiber product of
$\Fbar_1\times\Fbar_2$
and $\nbigmhat_3$ over
$B(W_0^{(1)},W_0^{(2)})$.
Hence we have the following commutative diagram:
\begin{equation}
 \label{eq;06.6.15.2}
\begin{CD}
 \Mhat^{G_m}(\gbigi)
 @>{g}>>
 \Fbar_1\times\Fbar_2\\
 @V{\pi_1}VV @V{q}VV \\
 \nbigmhat_3 @>>>
 B(W_0^{(1)},W_0^{(2)})
\end{CD}
\end{equation}
Then we have the isomorphism
$g^{\ast}L_{\Fbar_1\times\Fbar_2/
  B(W_0^{(1)},W_0^{(2)}) }
\simeq
 L_{\Mhat^{G_m}(\gbigi)/\nbigmhat_3}$.
Due to Lemma \ref{lem;06.5.22.150},
we have the following morphisms:
\[
\begin{CD}
 g^{\ast}L_{\Fbar_1\times\Fbar_2/
 B(W_0^{(1)},W_0^{(2)}) }[-1]
@>{\varphi_1}>>
 g^{\ast}q^{\ast}
 L_{B(W_0^{(1)},W_0^{(2)})}
@>{\varphi_2}>>
 \pi_1^{\ast}\Ob(\nbigmhat_3)
\end{CD}
\]
We put 
$\Ob\bigl(\Mhat^{G_m}(\gbigi)\bigr)
:=\Cone(\varphi_2\circ\varphi_1)$.
Then we obtain the morphism
$\ob\bigl(\Mhat^{G_m}(\gbigi)\bigr):
 \Ob\bigl(\Mhat^{G_m}(\gbigi)\bigr)
\lrarr L_{\Mhat^{G_m}(\gbigi)}$.
By the same argument as that in the subsubsection
\ref{subsubsection;06.5.8.5},
it can be shown that
$\ob\bigl(\Mhat^{G_m}(\gbigi)\bigr)$ gives an obstruction theory
for $\Mhat^{G_m}(\gbigi)$.

We put $N:=H_{y}(m)$.
Let $\Flag(V_m^{(i)},I_i)$ denote the moduli
of filtrations $\nbigf^{(i)}_{\ast}$ of $V_m^{(i)}$
as follows:
\[
 \nbigf^{(i)}_1\subset
 \nbigf^{(i)}_2\subset\cdots\subset
 \nbigf^{(i)}_{N}=V^{(i)}_m,
\quad
\dim \nbigf^{(i)}_{j}/\nbigf^{(i)}_{j-1}
=\left\{
 \begin{array}{ll}
 1 & (j\in I_i)\\
 0 & (j\not\in I_i)
 \end{array}
 \right.
\]
We put $V_m:=V_m^{(1)}\oplus V_m^{(2)}$.
Recall that $\Flag(V_m,\Nbar)$ denotes the full flag variety
of $V_m$.
We have the naturally defined inclusion
$\Flag(V_m^{(1)},I_1)\times \Flag(V_m^{(2)},I_2)
\lrarr \Flag(V_m,\Nbar)$.
Clearly, $\Flag(V^{(i)}_m,I_i)$ is naturally isomorphic $F_i$.
We use the notation in the subsubsection
\ref{subsubsection;06.5.22.1000}.
Then, the diagram (\ref{eq;06.6.15.2})
is compatible with the following diagram,
which is given by (\ref{eq;06.6.15.3}):
\[
 \begin{CD}
 \Mhat^{\ast} @>>> \Fbar \\
 @VVV @VVV \\ 
 \nbigmhat_0 @>>> B(W_0)
 \end{CD}
\]
By using the description (\ref{eq;06.5.22.250})
of $\Ob(\Mhat^{\ast})$,
we also obtain (\ref{eq;06.5.22.151})
from (\ref{eq;06.5.11.10}).

We would like to construct 
$\Obtilde\bigl(\Mhat^{G_m}(\gbigi)\bigr)$.
Corresponding to the equivalence
$\Obtilde_1(\nbigmhat_3)\simeq\Ob(\nbigmhat_3)$,
we have the equivalent obstruction theory
$\Obtilde_1\bigl(\Mhat^{G_m}(\gbigi)\bigr)\simeq
 \Ob\bigl(\Mhat^{G_m}(\gbigi)\bigr)$.
Due to Lemma \ref{lem;06.5.22.201},
we obtain the deformation
$\obtilde\bigl(\Mhat^{G_m}(\gbigi)\bigr):
\Obtilde\bigl(\Mhat^{G_m}(\gbigi)\bigr)
\lrarr L_{\Mhat^{G_m}(\gbigi)\times A^1/A^1}$
from 
$\obtilde(\nbigmhat_3):
 \Obtilde(\nbigmhat_3)\lrarr
 L_{\nbigmhat_3\times A^1/A^1}$.
We also obtain
the distinguished triangle (\ref{eq;06.5.22.165})
and the splitting (\ref{eq;06.5.22.160}) at $t=0$ 
from Lemma \ref{lem;06.5.22.200}.
Thus the proof of Proposition \ref{prop;06.5.22.100}
is finished.
\hfill\qed

\subsubsection{The case where the $2$-stability condition
 is satisfied}
\label{subsubsection;06.7.3.36.1}

Let us describe the obstruction theory
of the fixed point set of the master space,
when the $2$-stability condition is satisfied.
We use the notation in the subsubsection
\ref{subsubsection;06.5.17.100}.
We give only the statement.
We put
$\nbigm_{\spl}:=
 \nbigm^{s}(\vecy_1,L,\alpha_{\ast},\delta)
 \times \nbigm^s(\vecyhat_2,\alpha_{\ast})$.

\begin{prop}
\label{prop;06.6.10.11}
We have the obstruction theory
$\ob(\Mhat^{G_m}(\gbigi)):\Ob(\Mhat^{G_m}(\gbigi))
\lrarr L_{\Mhat^{G_m}(\gbigi)}$
and the deformation
$\obtilde(\Mhat^{G_m}(\gbigi)):
 \Obtilde(\Mhat^{G_m}(\gbigi))\lrarr L_{\Mhat^{G_m}(\gbigi)\times A^1/A^1}$
with the following property:
\begin{itemize}
\item
We have the following commutative diagram:
\begin{equation}
 \begin{CD}
 \varphi_{\gbigi}^{\ast}\Ob(\Mhat)
 @>>>
 \varphi_{\gbigi}^{\ast}L_{\Mhat}\\
 @VVV @VVV \\
 \Ob(\Mhat^{G_m}(\gbigi))@>>>
 L_{\Mhat^{G_m}(\gbigi)}
 \end{CD}
\end{equation}
\item
 Let
 $\obtilde_a(\Mhat^{G_m}(\gbigi)):
 \Obtilde_a(\Mhat^{G_m}(\gbigi))
 \lrarr L_{\Mhat^{G_m}(\gbigi)}$
 denote the specialization
 of $\obtilde(\Mhat^{G_m}(\gbigi))$
 at $t=a$.
 At $t=1$,
 we have $\obtilde_1(\Mhat^{G_m}(\gbigi))= \ob(\Mhat^{G_m}(\gbigi))$.
 At $t=0$, we have the following commutative diagram
 in the diagram {\rm(\ref{eq;06.5.22.300})}:
\begin{equation}
 \begin{CD}
 F^{\ast}\Obtilde_0(\Mhat^{G_m}(\gbigi)) @>>>
 F^{\ast}L_{\Mhat^{G_m}(\gbigi)}\\
 @V{\simeq}VV @V{\simeq}VV \\
 G^{\prime\ast}\Ob(\nbigm_{\spl})
 @>>>
 G^{\prime\ast}L_{\nbigm_{\spl}}
\end{CD}
\end{equation}
On each $a$, we have the following distinguished triangle:
\begin{equation}
\begin{CD}
 G^{\prime\,\ast}\Ob(m,\vecy_1,L)
 @>>>
 F^{\ast}\Obtilde_a(\Mhat^{G_m}(\gbigi))
 @>>>
 G^{\prime\,\ast}\Ob(m,\vecyhat_2)
 @>>>
 G^{\prime\,\ast}\Ob(m,\vecy_1,L)[1]
\end{CD}
\end{equation}
\end{itemize}
\end{prop}
\pf
It can be shown by an argument
similar to the proof of Proposition \ref{prop;06.5.22.100}.
In this case,
$\Mhat^{G_m}(\gbigi)$ is an open substack of
$\nbigmhat_3$.
The obstruction theory
$\ob(\nbigmhat_3):\Ob(\nbigmhat_3)\lrarr \Mhat^{G_m}(\gbigi)$ 
and the deformation
$\obtilde(\nbigmhat_3):
\Obtilde(\nbigmhat_3)\lrarr L_{\Mhat^{G_m}(\gbigi)\times A^1/A^1}$
gives the desired objects.
\hfill\qed

\subsubsection{The case of oriented reduced $\vecL$-Bradlow pair}
\label{subsubsection;06.6.12.1}

Let $\vecL=(L_1,L_2)$ be a pair of line bundles over $X$.
Let us describe the obstruction theory
of the fixed point set of the master space
for the moduli stack of the oriented 
$\vecL$-Bradlow pairs,
under the setting
in the subsubsection \ref{subsubsection;06.5.21.555}.
We give only the statement.
We put
$\nbigm_{\spl}:=
 \nbigm^{s}(\vecy_1,L_1,\alpha_{\ast},\delta_1)
 \times \nbigm^{s}(\vecyhat_2,[L_2],\alpha_{\ast},\delta_2)$.

\begin{prop}
\label{prop;06.6.10.15}
We have the obstruction theory
$\ob(\Mhat^{G_m}(\gbigi)):\Ob(\Mhat^{G_m}(\gbigi))
\lrarr L_{\Mhat^{G_m}(\gbigi)}$
and the deformation
$\obtilde(\Mhat^{G_m}(\gbigi)):
 \Obtilde(\Mhat^{G_m}(\gbigi))\lrarr 
 L_{\Mhat^{G_m}(\gbigi)\times A^1/A^1}$
with the following property:
\begin{itemize}
\item
We have the following commutative diagram:
\begin{equation}
 \begin{CD}
 \varphi_{\gbigi}^{\ast}\Ob(\Mhat)
 @>>>
 \varphi_{\gbigi}^{\ast}L_{\Mhat}\\
 @VVV @VVV \\
 \Ob(\Mhat^{G_m}(\gbigi))@>>>
 L_{\Mhat^{G_m}(\gbigi)}
 \end{CD}
\end{equation}
\item
 Let
 $\obtilde_a(\Mhat^{G_m}(\gbigi)):
 \Obtilde_a(\Mhat^{G_m}(\gbigi))
 \lrarr L_{\Mhat^{G_m}(\gbigi)}$
 denote the specialization
 of $\obtilde(\Mhat^{G_m}(\gbigi))$
 at $t=a$.
 At $t=1$,
 we have $\obtilde_1(\Mhat^{G_m}(\gbigi))= \ob(\Mhat^{G_m}(\gbigi))$.
 At $t=0$, we have the following commutative diagram
 in the diagram {\rm(\ref{eq;06.5.22.300})}:
\begin{equation}
 \begin{CD}
 F^{\ast}\Obtilde_0(\Mhat^{G_m}(\gbigi)) @>>>
 F^{\ast}L_{\Mhat^{G_m}(\gbigi)}\\
 @V{\simeq}VV @V{\simeq}VV \\
 G^{\prime\ast}\Ob(\nbigm_{\spl})
 @>>>
 G^{\prime\ast}L_{\nbigm_{\spl}}
\end{CD}
\end{equation}
On each $a$, we have the following distinguished triangle:
\begin{equation}
\begin{CD}
 G^{\prime\,\ast}\Ob(m,\vecy_1,L_1)
 @>>>
 F^{\ast}\Obtilde_a(\Mhat^{G_m}(\gbigi))
 @>>>
 G^{\prime\,\ast}\Ob(m,\vecyhat_2,[L_2])
 @>>>
 G^{\prime\,\ast}\Ob(m,\vecy_1,L_1)[1]
\end{CD}
\end{equation}
\end{itemize}
\end{prop}
\pf
We put $\nbigm_1:=\nbigm(m,\vecy_1,L_1)$,
$\nbigm_2:=\nbigm(m,\vecyhat_2,[L_2])$
and $\nbigm_3=\nbigm_1\times\nbigm_2$.
We consider the moduli stack $\nbigmhat_3$
of the objects
$(E_1,F_{1\,\ast},\phi_1,E_2,F_{2\,\ast},[\phi_2],\rho)$
as follows:
\begin{itemize}
\item
 $(E_1,F_{1\,\ast},\phi_1)\in\nbigm_1$
and
 $(E_2,F_{2\,\ast},[\phi_2])\in\nbigm_2$.
\item
$\rho$ denotes an orientation of $E_1\oplus E_2$.
\end{itemize}
Then $\Mhat^{G_m}(\gbigi)$ is the open substack
of $\nbigmhat_3$.

We can construct the obstruction theory
$\ob(\nbigmhat_3):\Ob(\nbigmhat_3)
\lrarr L_{\nbigmhat_3}$ by the argument
in the subsubsection \ref{subsubsection;06.5.22.350}.
We can also construct the deformation
$\obtilde(\nbigmhat_3):\Obtilde(\nbigmhat_3)
\lrarr L_{\Mhat^{G_m}(\gbigi)\times A^1/A^1}$
by the argument in the subsubsection 
\ref{subsubsection;06.5.22.360}.
It can be checked that they give the desired objects
by an argument similar to the proof of
Proposition \ref{prop;06.5.22.100}.
\hfill\qed

\subsection{Equivariant Obstruction Theory of the Master Space}
\label{subsection;06.6.17.2}
\subsubsection{Statements}
\label{subsubsection;06.5.31.11}

In this subsection,
we would like to explain the following claim, first.

\begin{prop}
\label{prop;06.6.10.2}
We have the $G_m$-equivariant lift of
the obstruction theory $\ob(\Mhat)$
of the enhanced master space $\Mhat$.
\end{prop}
The construction is explained in the subsubsection
\ref{subsubsection;06.5.23.30}--\ref{subsubsection;06.5.23.36}.
We explain that $\Ob(\Mhat^{\ast})$
can also be lifted equivariantly,
in the subsubsection \ref{subsubsection;06.5.23.31}.

We will later apply the localization formula
of Graber and Pandharipande (\cite{gp}).
For that purpose,
we have to see the induced obstruction theory
and the virtual normal bundle at the fixed point set.
Let $\iota_i:\Mhat_i\lrarr \Mhat$ denote the inclusion
$(i=0,1)$.
We have the decomposition of
$\iota_i^{\ast}\Ob(\Mhat)$
into the invariant part 
$\iota_i^{\ast}\Ob(\Mhat)^{\inv}$
and the moving part
$\iota_i^{\ast}\Ob(\Mhat)^{\mov}$.
We put
$\gbign(\Mhat_i):=\bigl(
\iota_i^{\ast}\Ob(\Mhat)^{\mov}
\bigr)^{\lor}[1]$,
which is called the virtual normal bundle.
We have the induced obstruction theory
$\ob_1(\Mhat_i):
 \iota_i^{\ast}\Ob(\Mhat)^{\inv}\lrarr L_{\Mhat_i}$
(\cite{gp}).
We will prove the following proposition
in the subsubsection \ref{subsubsection;06.5.23.32}
\begin{prop}
\label{prop;06.5.23.40}
The induced obstruction theory
$\ob_1(\Mhat_i)$ is equivalent to
the moduli theoretic obstruction theory
$\ob(\Mhat_i):\Ob(\Mhat_i)\lrarr L_{\Mhat_i}$
given in the subsubsection 
{\rm\ref{subsubsection;06.5.13.40}}.
The virtual normal bundle
$\gbign(\Mhat_i)$ is isomorphic to
$\nbigo_{\rel}\bigl((-1)^{i-1}\bigr)$,
and the weight of the induced $G_m$-action
is $(-1)^i$.
\end{prop}

Let $\gbigi=(\vecy_1,\vecy_2,I_1,I_2)$ be a decomposition type.
Let $\varphi_{\gbigi}:\Mhat^{G_m}(\gbigi)\lrarr \Mhat$
denote the inclusion.
Similarly,
we have the decomposition of
$\varphi_{\gbigi}^{\ast}\Ob(\Mhat)$ 
into the invariant part 
$\varphi_{\gbigi}^{\ast}\Ob(\Mhat)^{\inv}$
and the moving part
$\varphi_{\gbigi}^{\ast}\Ob(\Mhat)^{\mov}$.
We obtain the virtual normal bundle
$\gbign(\Mhat^{G_m}(\gbigi)):=
\bigl(
 \varphi_{\gbigi}^{\ast}\Ob(\Mhat)^{\mov}
\bigr)^{\lor}[1]$,
and the obstruction theory
$\ob_1(\Mhat^{G_m}(\gbigi)):
 \varphi_{\gbigi}^{\ast}\Ob(\Mhat)^{\inv}
\lrarr L_{\Mhat^{G_m}(\gbigi)}$.

To describe 
$\gbign(\Mhat^{G_m}(\gbigi))$,
we prepare some notation.
In the $K$-theory of coherent sheaves on $\Mhat^{G_m}(\gbigi)$,
we put as follows:
\begin{equation}
 \label{eq;06.5.23.15}
 \gbign(E^{\Mhat}_i,E^{\Mhat}_j)
:=-\sum_{i=0,1,2}
(-1)^{i}R^ip_{X\ast}\nrhom(E^{\Mhat}_i,E^{\Mhat}_j)
\end{equation}
\begin{equation}
 \label{eq;06.5.23.16}
 \gbign(L,E^{\Mhat}_2)
:=\sum_{i=0,1,2}
(-1)^iR^ip_{X\,\ast}\nhom(L,E^{\Mhat}_2)
\end{equation}
\begin{equation}
 \label{eq;06.5.23.17}
 \gbign_D(E^{\Mhat}_{i\,\ast},E^{\Mhat}_{j\,\ast}):=
-\sum_{i=0,1}
(-1)^iR^ip_{D\,\ast}\nrhom'_2\bigl(
 E^{\Mhat}_{i|D\,\ast},E^{\Mhat}_{j|D\,\ast}
 \bigr)
\end{equation}
Here $E^{\Mhat}_{i|D\,\ast}$
denotes the restriction
$E^{\Mhat}_i\otimes\nbigo_{\Mhat^{G_m}(\gbigi)\times D}$
with the induced filtration.
(See the subsubsection \ref{subsubsection;06.5.23.10}
for $\nrhom'_2$.)
Let $\nbigmhat_0$ and $\nbigmhat_3$ be as 
in the subsubsections
\ref{subsubsection;06.5.22.350}--\ref{subsubsection;06.5.23.5}.
It is easy to observe that
$\Mhat^{G_m}(\gbigi)\lrarr 
 \Mhat\times_{\nbigmhat_0}\nbigmhat_3$
is a regular embedding.
The normal bundle is denoted by $N_0$.

We will prove the following proposition in the subsubsection
\ref{subsubsection;06.5.23.35}.
\begin{prop}
\label{prop;06.5.23.41}
\mbox{{}}
\begin{itemize}
\item
The induced obstruction theory
$\varphi_{\gbigi}^{\ast}\Ob(\Mhat^{\ast})^{\inv}$
is naturally isomorphic to
$\Ob(\Mhat^{G_m}(\gbigi))$,
and the induced obstruction theory
$\varphi_{\gbigi}^{\ast}\Ob(\Mhat^{\ast})^{\inv}
\lrarr
 L_{\Mhat^{G_m}(\gbigi)}$
is equivalent to
$\ob(\Mhat^{G_m}(\gbigi)):
\Ob(\Mhat^{G_m}(\gbigi))\lrarr L_{\Mhat^{G_m}(\gbigi)}$.
\item
The virtual normal bundle $\gbign(\Mhat^{G_m}(\gbigi))$
is $K$-theoretically given as follows:
\begin{multline}
 \gbign(E^{\Mhat}_1,E^{\Mhat}_2)\otimes I_{1+r_1/r_2}
+\gbign(E^{\Mhat}_2,E^{\Mhat}_1)\otimes I_{-1-r_1/r_2}
+\gbign(L,E^{\Mhat}_2)\otimes I_{1+r_1/r_2} \\
+\gbign_D(E^{\Mhat}_{1\ast},
 E^{\Mhat}_{2\ast})\otimes I_{1+r_1/r_2}
+\gbign_D(E^{\Mhat}_{2\ast},
 E^{\Mhat}_{1\ast})\otimes I_{-1-r_1/r_2}
+N_0
\end{multline}
Here, $I_{n}$ denote the trivial line bundle on
$\Mhat^{G_m}(\gbigi)$
with the $G_m$-action of weight $n$.
\end{itemize}
\end{prop}

Before going into the proof,
we give a remark.
Recall the diagram (\ref{eq;06.5.23.20}).
We put $A=1+r_1/r_2$.
Let $\nbigo_{1,\rel}(1)$ denote the tautological line bundle
obtained from
$\nbigmtilde\bigl(\vecyhat_1,[L],\alpha_{\ast},(\delta,k)\bigr)$.
On $\nbigs$, we have the following:
\begin{equation}
\label{eq;06.5.23.25}
 F^{\ast}\gbign(E^{\Mhat}_1,E^{\Mhat}_2)
=G^{\ast}
\gbign(\Ehat^u_1,\Ehat^u_2)
 \otimes \nbigo_{1,\rel}(-A),
\quad
 F^{\ast}\gbign(E^{\Mhat}_2,E^{\Mhat}_1)
=G^{\ast}
 \gbign(\Ehat^u_2,\Ehat^u_1)
 \otimes\nbigo_{1,\rel}(A)
\end{equation}
\begin{equation}
\label{eq;06.5.23.26}
 F^{\ast}\gbign(L,E^{\Mhat}_2)
=G^{\ast}\gbign(L,\Ehat^u_2)\otimes\nbigo_{1,\rel}(-r_1/r_2)
\end{equation}
\begin{equation}
\label{eq;06.5.23.27}
 F^{\ast}\gbign_D(E^{\Mhat}_{1\ast},E^{\Mhat}_{2\ast})
=G^{\ast}\bigl(
 \gbign_D(\Ehat^u_{1\ast},\Ehat^u_{2\ast})\otimes\nbigo_{1,\rel}(-A)
 \bigr),
\,\,\,
 F^{\ast}\gbign_D(E^{\Mhat}_{2\ast},E^{\Mhat}_{1\ast})
=G^{\ast}\bigl(
 \gbign_D(\Ehat^u_{2\ast},\Ehat^u_{1\ast})
 \otimes\nbigo_{1,\rel}(A)\bigr).
\end{equation}
Here, $\gbign(\Ehat^u_{a},\Ehat^u_{b})$,  $\gbign(L,\Ehat^u_2)$ and
$\gbign_D(\Ehat^u_{a\ast},\Ehat^u_{b\ast})$
are elements of the $K\bigl(
\nbigmtilde\bigl(\vecyhat_1,[L],\alpha_{\ast},(\delta,\ell)\bigr)
 \times\nbigmtilde(\vecyhat_2,\alpha_{\ast})\bigr)$,
given as follows:
\begin{equation}
\label{eq;06.5.31.15}
 \gbign(\Ehat^{u}_i,\Ehat^{u}_j)
:=-\sum_{i=0,1,2}
(-1)^{i}R^ip_{X\ast}\nrhom(\Ehat^{u}_i,\Ehat^{u}_j)
\end{equation}
\begin{equation}
 \gbign(L,\Ehat^{u}_2)
:=\sum_{i=0,1,2}
(-1)^iR^ip_{X\,\ast}\nhom(L,\Ehat^{u}_2)
\end{equation}
\begin{equation}
 \gbign_D(\Ehat^{u}_{i\,\ast},\Ehat^{u}_{j\,\ast}):=
-\sum_{i=0,1}
(-1)^iR^ip_{D\,\ast}\nrhom'_2\bigl(
 \Ehat^{u}_{i|D\,\ast},\Ehat^{u}_{j|D\,\ast}
 \bigr)
\end{equation}
Here $\Ehat^{u}_{i|D\,\ast}$
denotes the restriction
$\Ehat^{u}_i\otimes\nbigo_{\Mhat^{G_m}(\gbigi)\times D}$
with the induced filtration.

\subsubsection{$G_m$-equivariant lift of $\Ob(\Mhat)$ and $\ob(\Mhat)$}
\label{subsubsection;06.5.23.30}

We would like to obtain the $G_m$-equivairant structure
of $\Ob(\Mhat)$.
We use the notation in the subsubsection
\ref{subsubsection;06.5.8.5}.
We have the following commutative diagram:
\[
 \begin{CD}
 \Mhat @>{\Psi_2}>> \nbigqtilde\\
 @V{\gminip}VV @V{\gminip'}VV \\
 \nbigm(m,\vecyhat,[L]) @>{\Psi_1}>>
 B(W_{\cdot},[P_{\cdot}])
 \end{CD}
\]
We have the $G_m$-action on $Z_1$
given by $t[u_0:u_1]=[t\cdot u_0:u_1]$.
It induces the $G_m$-action on $\nbigqtilde$.
We remark that $\Psi_2$ is $G_m$-equivariant map.
We have the natural $G_m$-equivariant structure
of $\gminip^{\ast}\Ob(m,\vecyhat,[L])$
and $\Psi_2^{\ast}L_{\nbigqtilde/B(W_{\cdot},[P_{\cdot}])}$.

Let $C_i$ $(i=1,2)$ be a bounded $G_m$-complex
on $\Mhat$.
We have the induced $G_m$-action
on $\Ext^0(C_1,C_2)$,
where $\Ext^0(C_1,C_2)$ denotes 
the vector space of the homomorphisms
of $C_1$ to $C_2$ in the derived category $D(\Mhat)$.
In the following,
we say that a morphism $\varphi:C_1\lrarr C_2$
is contained in the $G_m$-invariant part
of the $\Ext^0$-group,
if $\varphi$ is contained in the $G_m$-invariant
part of $\Ext^0(C_1,C_2)$.

\begin{lem}
\label{lem;06.5.22.556}
The morphism
$\Psi_2^{\ast}L_{\nbigqtilde/B(W_{\cdot},[P_{\cdot}])}[-1]
\lrarr \gminip^{\ast}\Ob(m,\vecyhat,[L])$
is contained in the $G_m$-invariant part of
the $\Ext^0$-group.
\end{lem}
\pf
By using Remark \ref{rem;06.5.22.500},
we obtain the natural $G_m$-equivariant representatives
$C'(\nbigqtilde)$ and $C'(B(W_{\cdot},[P_{\cdot}]))$, $C_3$
of $\Psi_2^{\ast}L_{\nbigqtilde}$,
$\Psi_2^{\ast}\gminip^{\prime\,\ast}L_{B(W_{\cdot},[P_{\cdot}])}$
and 
$\Psi_2^{\ast}L_{\nbigqtilde/B(W_{\cdot},[P_{\cdot}])}$
respectively.
We regard them as $G_m$-equivariant complex.
We have the natural $G_m$-equivariant morphism
$\alpha:C(B(W_{\cdot},[P_{\cdot}]))\lrarr C(\nbigqtilde)$.
We put $C'(\nbigqtilde/B(W_{\cdot},[P_{\cdot}])):=\cone(\alpha)$.
Then, we obtain
 the $G_m$-equivariant quasi isomorphism
$C_3\simeq C'(\nbigqtilde/B(W_{\cdot},[P_{\cdot}]))$.
Then the morphism
$\Psi_2^{\ast}L_{\nbigqtilde/B(W_{\cdot},[P_{\cdot}])}[-1]
\lrarr
 \Psi_2^{\ast}\gminip^{\prime\,\ast}
 L_{B(W_{\cdot},[P_{\cdot}])}$
can be expressed by the $G_m$-equivariant morphism:
\[
 C'(\nbigqtilde/B(W_{\cdot},[P_{\cdot}]))[-1]\lrarr 
 C'(B(W_{\cdot},[P_{\cdot}]))
\]

On the other hand,
the morphism
$\gminip^{\ast}u:\gminip^{\ast}\Psi_2^{\ast}L_{B(W_{\cdot},[P_{\cdot}])}
\lrarr \gminip^{\ast}\Ob(m,\vecyhat,[L])$
is contained in the $G_m$-invariant part of 
the $\Ext^0$-group.
Then the claim of the lemma is clear.
\hfill\qed

\vspace{.1in}

By applying the general non-sense
in the subsubsection \ref{subsubsection;06.6.22.5},
we obtain $G_m$-equivariant representative $C(\Mhat)$
of $\Ob(\Mhat)$.

\vspace{.1in}

Recall that we need only
the $(-1)$-truncated cotangent complexes
for the construction of the virtual classes,
and the complex $C(\Mhat)$ above are used for the localization theory
by Graber and Pandharipande (\cite{gp}).
Since $L_{\Mhat/\nbigm(m,\vecyhat,[L])}$
is isomorphic to the $0$-th cohomology sheaf,
we have the distinguished triangle:
\[
\tau_{\geq\,-1}
 \gminip^{\ast}L_{\nbigm(m,\vecyhat,[L])}
\lrarr
\tau_{\geq\,-1}
 L_{\Mhat}
\lrarr 
 L_{\Mhat/\nbigm(m,\vecyhat,[L])}
\lrarr
\tau_{\geq\,-1}
 \gminip^{\ast}L_{\nbigm(m,\vecyhat,[L])}[1]
\]
Since the morphism
$\gminip^{\ast}\Ob(m,\vecyhat,[L])\lrarr 
\tau_{\geq -1}\gminip^{\ast}L_{\nbigm(m,\vecyhat,[L])}$
is contained in the $G_m$-invariant part
of the $\Ext^0$-group,
the morphism
$L_{\Mhat/\nbigm(m,\vecyhat,[L])}[-1]
\lrarr
 \gminip^{\ast}L_{\nbigm(m,\vecyhat,[L])}$
is also contained in the $G_m$-invariant part of
the $\Ext^0$-group,
due to Lemma \ref{lem;06.5.22.556}.
Therefore,
we also obtain the $G_m$-equivariant structure
of $\tau_{\geq\,-1}L_{\Mhat}$.

By construction,
$\ob(\Mhat):\Ob(\Mhat)\lrarr \tau_{\geq\,-1}L_{\Mhat}$
is contained in the $G_m$-invariant part of
the $\Ext^0$-group.
Therefore, we can take the $G_m$-equivariant
representatives of 
$\Ob(\Mhat)$, $\tau_{\geq\,-1}L_{\Mhat}$
and $\ob(\Mhat)$.

\subsubsection{Equivalent $G_m$-equivariant structure of $\tau_{\geq -1}L_{\Mhat}$ }
\label{subsubsection;06.5.23.36}

On the other hand,
we have the $G_m$-equivariant structure
of $\tau_{\geq\,-1}L_{\Mhat}$,
obtained from the $G_m$-equivariant embedding
into a smooth Deligne-Mumford stack.
We use the notation in the subsubsection
\ref{subsubsection;06.5.22.369}.

Let $\Zhat_m$ be the vector bundle
as in (\ref{eq;06.5.21.30}).
Then we put
$\nbigbtildehat:=\nbigbtilde\times_{Z_m}\Zhat_m$.
We have the natural morphism
$\THhat^{ss}\lrarr \nbigbtildehat$
which is $\GL(V_m)\times G_m$-equivariant.
Therefore,
we obtain the $G_m$-equivariant immersion
$\iota:\Mhat\lrarr \nbigbtildehat_{\GL(V_m)}$.
Since $\Mhat$ is Deligne-Mumford,
we can take a smooth Deligne-Mumford 
open substack $P$ of $\nbigbtildehat_{\GL(V_m)}$,
which contains $\Mhat$.

Let $I$ denote the ideal sheaf of $P$ corresponding to $\Mhat$.
We put $C(\Mhat):=\Cone\bigl(I/I^2\lrarr \iota^{\ast}C(P)\bigr)$
on $\Mhat$,
where we put $C(P):=\Omega_{P/k}$.
It is naturally $G_m$-complex,
and it is the representatives
of the $(-1)$-truncated cotangent complex
$\tau_{\geq\,-1}L_{\Mhat/k}$.

\begin{lem}
The above two $G_m$-equivariant structures 
of $\tau_{\geq\,-1}L_{\Mhat}$ are equivalent.
\end{lem}
\pf
We have the complex
$C'(\nbigbtildehat):=
 \Cone(\Omega_{\nbigbtildehat}\lrarr \Omega_{\nbigbtildehat/P})[-1]$
on $\nbigbtildehat$.
It is naturally $\GL(V_m)\times G_m$-equivariant,
and thus 
it induces $G_m$-equivariant complex $C'(P)$ on $P$.
We have the natural $G_m$-equivariant quasi isomorphism
$C(P)\lrarr C'(P)$ on $P$.
We put $C'(\Mhat):=\Cone\bigl(I/I^2\lrarr \iota^{\ast}C'(P)\bigr)$
on $\Mhat$,
then we have the natural $G_m$-equivariant
quasi isomorphism $C(\Mhat)\lrarr C'(\Mhat)$.

We put
$\nbigahat:=\nbiga \times_{Z_m}\Zhat_m$.
We have the natural $\GL(V_m)$-action
on $\nbigahat$.
The quotient stack is denoted by $\gbigq$.
We have the $\GL(V_m)$-equivariant map
$Q^{\circ}(m,\vecyhat,[L])\lrarr  \nbigahat$,
and hence
$\nbigm(m,\vecyhat,[L])\lrarr \gbigq$.

Let $\gminip_3:\nbigbtildehat\lrarr\nbigahat$
denote the projection.
We put
$C'(\nbigahat):=
\Cone\bigl(\gminip_{3}^{\ast}\Omega_{\nbigahat}
 \lrarr \gminip_3^{\ast}\Omega_{\nbigahat/\gbigq}\bigr)[-1]$.
The complex is provided with the natural $\GL(V_m)$-action.
Hence, it induces the complex $C'(\gbigq)$ on $P$.
We have the natural morphism $C'(\gbigq)\lrarr C'(P)$.

It is clear that the morphism $I/I^2\lrarr \iota^{\ast}C(P)$
factors through $\iota^{\ast}C'(\gbigq)$.
We put $C'(\nbigm)=
\Cone\bigl(I/I^2\lrarr\iota^{\ast}C'(\gbigq)\bigr)$.
We have the exact sequences of the $G_m$-equivariant
complexes:
\[
\begin{array}{ll}
 0\lrarr C'(\gbigq)\lrarr C'(P)
 \lrarr \Omega_{P/\gbigq}\lrarr 0
 & \mbox{ on $P$ }
 \\
\mbox{{}}\\
 0\lrarr C'(\nbigm)\lrarr
 C'(\Mhat)\lrarr \Omega_{\Mhat/\nbigm}\lrarr 0
 & \mbox{ on $\Mhat$}
\end{array}
\]

We put $C'(P/\gbigq):=\Cone\bigl(C'(\gbigq)\lrarr C'(P)\bigr)$.
We have the $G_m$-equivariant morphism
$C'(P/\gbigq)[-1]\lrarr C'(\gbigq)$ on $P$.
We have the natural $G_m$-equivariant morphism
$\iota^{\ast}C'(P/\gbigq)[-1]
 \lrarr 
 \iota^{\ast}C'(\gbigq)
 \lrarr C'(\nbigm)$ of the $G_m$-complexes on $\Mhat$.
We put 
$C_1(\Mhat):=
\Cone\bigl(\iota^{\ast}C(P/\gbigq)[-1]\lrarr
 C'(\nbigm)\bigr)$.

Let us show that $C_1(\Mhat)$ is $G_m$-equivariant
quasi-isomorphic to $C'(\Mhat)$.
We put $C_0:=\Cone(I/I^2\lrarr I/I^2)$.
We have the composite of the morphisms
$C'(\Mhat)\lrarr I/I^2[1]\lrarr C_0[1]$.
We have the naturally defined $G_m$-equivariant
quasi isomorphism
$\Cone\bigl(C'(\Mhat)\rarr C_0[1]\bigr)[-1]
 \lrarr C'(\Mhat)$.
We have the morphism $I/I^2\lrarr \iota^{\ast}C'(\gbigq)$.
It induces the $G_m$-equivariant quasi isomorphism
$\Cone\bigl(C'(\Mhat)\rarr C_0[1]\bigr)[-1]
 \lrarr C_1(\Mhat)$.

We have the following commutative diagram:
\[
 \begin{CD}
 \Mhat @>>> P @>>> \nbigqtilde\\
 @VVV @VVV @VVV \\
 \nbigm @>>> \gbigq @>>> B(W_{\cdot},[P_{\cdot}])
 \end{CD}
\]
It induces the following commutative diagram:
\[
\begin{CD}
C'(\nbigqtilde)@>>> C'(P) \\
 @AAA @AAA \\
C'(B(W_{\cdot},[P_{\cdot}]))
@>>> C'(\gbigq) 
\end{CD}
\]
Hence,
we obtain the isomorphism
$C'(\nbigqtilde/B(W_{\cdot},[P_{\cdot}]))
\simeq  C'(P/\gbigq)$.
We also have
$C'(\nbigm)\simeq\gminip^{\ast}\tau_{\leq -1}L_{\nbigm}$.
Therefore, 
two $G_m$-equivariant structures are equivalent.
\hfill\qed

\subsubsection{$G_m$-equivariant structure of
 $\Ob(\Mhat^{\ast})$}
\label{subsubsection;06.5.23.31}

Since we have the quasi isomorphism
$\Ob(\Mhat)_{|\Mhat^{\ast}}\simeq
 \Ob(\Mhat^{\ast})$,
we have already obtained the $G_m$-equivariant
structure of $\Ob(\Mhat^{\ast})$.
We give another description of 
the $G_m$-equivariant structure.
We use the notation in the subsubsection 
\ref{subsubsection;06.5.22.1000}.

We have the natural $G_m$-equivariant
structure on the sheaves $\Ehat^{\Mhat}$
over $\Mhat\times X$.
It induces the $G_m$-equivariant structure
on $\gminip_2^{\ast}\Ob(m,\vecyhat,L)$.
On the other hand,
the morphisms $\Psi_3$ and $\gminip_2$
are $G_m$-equivariant.
Therefore, we have the $G_m$-equivariant
representative of $L_{B(W_{\cdot},P_{\cdot})}$.

\begin{lem}
The morphism
$\gminip_2^{\ast}\Psi_3^{\ast}
 L_{B(W_{\cdot},P_{\cdot})}
\lrarr
 \gminip_2^{\ast}\Ob(m,\vecyhat,L)$
is $G_m$-equivariant.
\end{lem}
\pf
Let $\nbigv_{\cdot}$ denote the canonical resolution
of $\Ehat^{\Mhat}(m)$.
By the remark \ref{rem;06.5.22.500},
$\gminih(\nbigv_{\cdot},\phitilde)_{\leq 1}$
gives the $G_m$-equivariant representative
of $\Phi_{3X}^{\ast}L_{Y(W_{\cdot},P_{\cdot})/X}$,
where the $G_m$-equivariant structure of
$\gminih(\nbigv_{\cdot},\phitilde)$
is induced by the $G_m$-equivariant structure
of $\Ehat^{\Mhat}$.
Therefore, 
the $G_m$-equivariant representative
of $\Phi_X^{\ast}L_{Y(W_{\cdot},P_{\cdot})/X}$
is given by
$\Ob^G(\nbigv_{\cdot},\phitilde)$.
On the other hand,
the morphism
$\Ob^{G}(\nbigv_{\cdot},\phitilde)\lrarr 
 \Ob(m,\vecyhat,L)$ is $G_m$-equivariant,
because their $G_m$-equivariant structures
are induced by that of $\Ehat^{\Mhat}$.
Thus we are done.
\hfill\qed

\vspace{.1in}

Then, the morphism
$\Psi_4^{\ast}L_{\nbigqtilde^{\ast}/B(W_{\cdot},P_{\cdot})}[-1]
\lrarr
 \gminip_2^{\ast}\Ob(m,\vecyhat,L)$
is contained in the $G_m$-invariant part
of the $\Ext^0$-group.
Therefore, we can take the $G_m$-equivariant
representative of $\Ob(\Mhat^{\ast})$.

\begin{lem}
Two $G_m$-equivariant structure
of $\Ob(\Mhat^{\ast})$
are equivalent.
\end{lem}
\pf
We have the following diagram:
\[
 \begin{CD}
 \Psi_4^{\ast}L_{\nbigqtilde^{\ast}/Y(W_{\cdot},P_{\cdot})}[-1]
 @>>>
 \Ob(m,\vecyhat,L)\\
 @A{a}AA @A{b}AA \\
 \Psi_4^{\ast}L_{\nbigqtilde^{\ast}/Y(W_{\cdot},[P_{\cdot}])}[-1]
 @>>>
 \Ob(m,\vecyhat,[L])
 \end{CD}
\]
The morphisms are contained
in the $G_m$-invariant part of
the $\Ext^0$-groups.
Therefore, the induced $G_m$-equivariant structures
on $\cone(a)$ and $\cone(b)$ are equivalent.
\hfill\qed

\subsubsection{Proof of Proposition \ref{prop;06.5.23.40}}
\label{subsubsection;06.5.23.32}

We use the notation in the subsection
\ref{subsection;06.5.11.100}
and the subsection \ref{subsection;06.5.23.1}.
We put $\Fbar:=\Flag(V_m,\Nbar)_{\GL(V_m)}$
and $\Zbar_1:=Z_{1\,\GL(V_m)}$.
Since we have
$\nbigqtilde=\Zbar_{1}\times_{B(W_{\cdot})}\Fbar$,
we have
$L_{\nbigqtilde/B(W,[P])}
=L_{\Zbar_1/B(W,[P])}\oplus L_{\Fbar/B(W)}$.

Recall $\Ob(\Mhat)$ is given by
$\cone\bigl(\Psi_2^{\ast}L_{\nbigqtilde/B(W,[P])}\lrarr
 \gminip^{\ast}\Ob(m,\vecyhat,[L])\bigr)$.
The moving part $\iota_i^{\ast}\Ob(\Mhat)^{\mov}$
is given by the following:
\[
 \iota_i^{\ast}\Psi_2^{\ast}L_{\Zbar_1/B(W,[P])}
\simeq
 L_{\Mhat_i/\Mhat}
\simeq
 \nbigo_{\rel}\bigl((-1)^i\bigr)[1]
\]
Then, it is easy to see that
the virtual normal bundle
$\gbign_{\Mhat}(\Mhat_i)$ is given by 
the line bundle $\nbigo_{\rel}\bigl((-1)^{i-1}\bigr)$,
and that the weight of the induced $G_m$-action is $(-1)^i$.
By construction,
we have the following:
\[
\iota_i^{\ast}\Ob(\Mhat)^{\inv}:=
\Cone\bigl( \Psi_{13}^{\ast}L_{\Fbar/B(W)}[-1]
\lrarr \Ob(m,\vecyhat,[L]) \bigr) 
\simeq 
\Ob(\Mhat_i)
\]
It is easy to observe that
we have the  following commutative diagram:
\[
 \begin{CD}
 \iota_i^{\ast}\Ob(\Mhat)@>{\varphi_1}>> \Ob(\Mhat_i)\\
 @V{\iota_i^{\ast}\ob(\Mhat)}VV @V{\ob(\Mhat_i)}VV \\
 \iota_i^{\ast}L_{\Mhat} @>{\varphi_2}>> L_{\Mhat_i}
 \end{CD}
\]
Here, $\varphi_2$ is the naturally defined one,
and $\varphi_1$ is the projection onto the invariant part.
Thus the induced obstruction theory
is given by $\ob(\Mhat_i)$.

\subsubsection{Proof of Proposition \ref{prop;06.5.23.41}}
\label{subsubsection;06.5.23.35}

We use the notation in the subsubsections
\ref{subsubsection;06.5.22.1000} and
\ref{subsubsection;06.5.23.2}.
Recall the expression of $\Ob(\Mhat^{\ast})$
as in (\ref{eq;06.5.22.250}).
Let us describe the decomposition
of $\Psi_4^{\ast}\Gamma^{\ast}L_{\Fbar/B(W_{0})}$.
We take a decomposition
$W_0=W_0^{(1)}\oplus W_0^{(2)}$
as in the subsubsection \ref{subsubsection;06.5.23.2}.
We put $\Fbar_i:=\Flag(W_0^{(i)},I_i)_{\GL(W_0^{(i)})}$.
We also put
$\Fbar':=\Flag(W_0,\Nbar)_{\GL(W_0^{(1)})\times\GL(W_0^{(2)})}$.
Then we have the regular immersion
$\Fbar_{1}\times\Fbar_2\lrarr \Fbar'$.
We have the following commutative diagram:
\[
 \begin{CD}
 \Mhat^{G_m}(\gbigi)@>{\Psi_{21}}>> 
 \Fbar_{1}\times\Fbar_2@>>>
 \Fbar \\
 @VVV @VVV @VVV \\
 \nbigmhat_3@>>>
 B(W_0^{(1)},W_0^{(2)}) @>>>
  B(W_0)
 \end{CD}
\]
Therefore, we have the isomorphism:
\[
 \Psi_4^{\ast}\Gamma^{\ast}
 L_{\Fbar/B(W_{\cdot})}\simeq
 \Psi_{21}^{\ast}L_{\Fbar'/B(W_0^{(1)},W_0^{(2)})}
\]
The invariant part of
$\Psi_4^{\ast}\Gamma^{\ast}
 L_{\Fbar/B(W_{\cdot})}$
is isomorphic to the pull back of
the relative cotangent bundle of $\Fbar_{1}\times\Fbar_2$
over $B(W_0^{(1)},W_0^{(2)})$.
The moving part is same as
the pull back of the conormal bundle of
$\Fbar_{1}\times\Fbar_2$ in $\Fbar'$,
which is naturally isomorphic to $N_0^{\lor}[1]$.

Let us see the decomposition of
$\Ob(m,\vecy,L)$.
Corresponding to the decomposition
$\varphi_{\gbigi}^{\ast}\Ehat^{\Mhat}=
 E^{\Mhat}_1\oplus E^{\Mhat}_2$,
we obtain the decomposition of the resolution
$\nbigv_{\cdot}=\nbigv_{1\,\cdot}\oplus\nbigv_{2\,\cdot}$.
It induces the following decompositions:
\[
 \gminig(\nbigv_{\cdot})^{\inv}=
 \gminig(\nbigv^{(1)}_{\cdot})
 \oplus\gminig(\nbigv^{(2)}_{\cdot}),
\quad
 \gminig(\nbigv_{\cdot})^{\mov}=
 \nhom\bigl(\nbigv^{(1)}_{\cdot},
 \nbigv^{(2)}_{\cdot}\bigr)^{\lor}[-1]
\oplus
 \nhom\bigl(\nbigv^{(2)}_{\cdot},
 \nbigv^{(1)}_{\cdot}\bigr)^{\lor}[-1]
\]

\[
 \gminig_{\rel}(\nbigv_{\cdot},\phi)^{\inv}=
 \gminig_{\rel}(\nbigv^{(1)}_{\cdot},\phi),
\quad
 \gminig_{\rel}(\nbigv_{\cdot},\phi)^{\mov}
=\nhom\bigl(P_{\cdot},\nbigv^{(2)}_{\cdot} \bigr)^{\lor}
\]

\[
 \gminig_D(\nbigv,F_{\ast})^{\inv}=
 \gminig_D(\nbigv^{(1)},F^{(1)}_{\ast})\oplus
 \gminig_D(\nbigv^{(2)},F^{(2)}_{\ast}),
\quad
  \gminig_D(\nbigv,F_{\ast})^{\mov}
=C_1\bigl(\nbigv_D^{(1)\ast},\nbigv_D^{(2)\ast}\bigr)^{\lor}[-1]
\oplus
  C_1\bigl(\nbigv_D^{(2)\ast},\nbigv_D^{(1)\ast}\bigr)^{\lor}[-1]
\]
(See the subsubsection \ref{subsubsection;06.5.23.10}
for $C_1$.)
\[
 \gminig(\nbigv_{\cdot\,|\,D})^{\inv}
=\gminig(\nbigv^{(1)}_{\cdot|D})\oplus
 \gminig(\nbigv^{(2)}_{\cdot|D}),
\quad
 \gminig(\nbigv_{\cdot|D})^{\mov}=
 \nhom_{\nbigo_D}(\nbigv^{(1)}_{\cdot|D},
 \nbigv^{(2)}_{\cdot|D})^{\lor}[-1]
\oplus
 \nhom_{\nbigo_D}(\nbigv^{(2)}_{\cdot|D},
 \nbigv^{(1)}_{\cdot|D})^{\lor}[-1]
\]
We also have 
$\gminig^d(\nbigv_{\cdot})^{\inv}=\gminig^d(\nbigv_{\cdot})$.
The contribution to the virtual normal bundle can be calculated
formally.
We can also easily observe that
$\Ob(m,\vecy,L)^{\inv}$ is naturally isomorphic
to $\Ob(\nbigmhat_3)$ given in the subsubsection 
\ref{subsubsection;06.5.22.350}.
Then $\iota^{\ast}\Ob(\Mhat^{\ast})^{\inv}$
is obtained as the cone of the composite
of the following morphisms:
\[
\begin{CD}
 L_{\Fbar_1\times\Fbar_2/B(W^{(1)},W^{(2)})}[-1]
 @>>>
 L_{B(W^{(1)},W^{(2)})}
@>>>
 \Ob(\nbigmhat_3)
\end{CD}
\]
Namely, it is same as $\Ob(\Mhat^{G_m}(\gbigi))$.
(See the subsubsection \ref{subsubsection;06.5.22.123}.)
We also have  the diagram (\ref{eq;06.5.22.151}).
Then it is easy to observe that the induced
obstruction theory is given by
$\ob(\Mhat^{G_m}(\gbigi))$.

\subsubsection{The case where the $2$-stability condition is satisfied}
\label{subsubsection;06.7.3.37.1}

We give only the statement about the $G_m$-equivariant
obstruction theory of the master space
in the case where the $2$-stability condition is satisfied.
The proof is similar to those of
Proposition \ref{prop;06.6.10.2},
Proposition \ref{prop;06.5.23.40}
and Proposition \ref{prop;06.5.23.41}.

\begin{prop}
\label{prop;06.6.10.1}
\mbox{{}}
Under the setting in the subsubsection 
{\rm \ref{subsubsection;06.5.17.100}},
the following claims hold:
\begin{itemize}
\item
 We have the equivariant obstruction theory
 of the master space $\Mhat$.
\item
 The induced obstruction theory of the fixed point sets
  $\Mhat_i$ and $\Mhat^{G_m}(\gbigi)$
 are equivalent to 
 the moduli theoretic obstruction theory.
\item
 The virtual normal bundle $\gbign(\Mhat_i)$
 of $\Mhat_i$ in $\Mhat$
 is given by $\nbigo_{\rel}\bigl((-1)^{i-1}\bigr)$ 
 with the $G_m$-action of the weight $(-1)^i$.
\item
 The virtual normal bundle $\gbign(\Mhat^{G_m}(\gbigi))$
 of $\Mhat^{G_m}(\gbigi)$ in $\Mhat$
 is given by the following:
\begin{multline}
 \gbign(E^{\Mhat}_1,E^{\Mhat}_2)\otimes I_{1+r_1/r_2}
+\gbign(E^{\Mhat}_2,E^{\Mhat}_1)\otimes I_{-1-r_1/r_2}
+\gbign(L,E^{\Mhat}_2)\otimes I_{1+r_1/r_2} \\
+\gbign_D(E^{\Mhat}_{1\,\ast},
 E^{\Mhat}_{2\,\ast})\otimes I_{1+r_1/r_2}
+\gbign_D(E^{\Mhat}_{2\,\ast},
 E^{\Mhat}_{1\,\ast})\otimes I_{-1-r_1/r_2}
\end{multline}
Here, 
$I_{w}$ denotes the trivial line bundle
with the $G_m$-action of weight $w$,
and the terms are as in 
{\rm(\ref{eq;06.5.23.15})},
{\rm(\ref{eq;06.5.23.16})}
and {\rm(\ref{eq;06.5.23.17})}.
(See also {\rm(\ref{eq;06.5.23.25})},
{\rm(\ref{eq;06.5.23.26})}
and {\rm(\ref{eq;06.5.23.27})}.)
\hfill\qed
\end{itemize}
\end{prop}

\subsubsection{The case of oriented reduced $\vecL$-Bradlow pairs}
\label{subsubsection;06.7.3.37.2}

We give only the statement about the $G_m$-equivariant
obstruction theory of the master space
in the case of the oriented reduced $\vecL$-Bradlow pairs,
under the setting of the subsubsection \ref{subsubsection;06.5.21.555}.
We prepare some notation.
For any decomposition type $\gbigi$,
the elements
 $\gbign(E_i^{\Mhat},E_j^{\Mhat})$,
$\gbign_D(E_{i\,\ast}^{\Mhat},E_{j\,\ast}^{\Mhat})$
and $\gbign(L_1,E_2^{\Mhat})$
of $K^{G_m}(\Mhat(\gbigi))$
are given as in 
{\rm(\ref{eq;06.5.23.15})},
{\rm(\ref{eq;06.5.23.16})}
and {\rm(\ref{eq;06.5.23.17})}.
(See also {\rm(\ref{eq;06.5.23.25})},
{\rm(\ref{eq;06.5.23.26})}
and {\rm(\ref{eq;06.5.23.27})}.)
Let $\varphi:\Mhat\lrarr\nbigm(m,\vecyhat,[\vecL])$
denote the naturally defined morphism.
We put $\nbigi^{(2)}:=\varphi^{\ast}\nbigo_{\rel}^{(2)}(-1)$,
which is naturally provided with the $G_m$-action.
We  have the following element
of $K^{G_m}(\Mhat^{G_m}(\gbigi))$:
\[
 \gbign(L_2\otimes\nbigi^{(2)},E^{\Mhat}_1):=
 \sum_{i=0,1,2}(-1)^i R^ip_{X\,\ast}
 \nhom\bigl(L_2\otimes\nbigi^{(2)},
 E^{\Mhat}_1\bigr)
\]

\begin{lem}
\mbox{{}}
Under the setting in the subsubsection 
{\rm \ref{subsubsection;06.5.17.100}},
the following holds:
\begin{itemize}
\item
 We have the equivariant obstruction theory
 of the master space $\Mhat$.
\item
 The induced obstruction theory of the fixed point sets
  $\Mhat_i$ and $\Mhat^{G_m}(\gbigi)$
 are equivalent to 
 the moduli theoretic obstruction theory.
\item
 The virtual normal bundle $\gbign(\Mhat_i)$
 of $\Mhat_i$ in $\Mhat$
 is given by $\nbigo^{(1)}_{\rel}\bigl((-1)^{i-1}\bigr)$ 
 with the $G_m$-action of the weight $(-1)^i$.
\item
 The virtual normal bundle
 $\gbign\bigl(\Mhat^{G_m}(\gbigi)\bigr)$
 of $\Mhat^{G_m}(\gbigi)$ in $\Mhat$
 is given by the following:
\begin{multline}
 \gbign(E^{\Mhat}_1,E^{\Mhat}_2)\otimes I_{1+r_1/r_2}
+\gbign(E^{\Mhat}_2,E^{\Mhat}_1)\otimes I_{-1-r_1/r_2}
+\gbign(L_1,E^{\Mhat}_2)\otimes I_{1+r_1/r_2}
+\gbign(L_2\otimes\nbigi^{(2)},E^{\Mhat}_1)\otimes I_{-1-r_1/r_2}\\
+\gbign_D(E^{\Mhat}_{1\,\ast},
 E^{\Mhat}_{2\,\ast})\otimes I_{1+r_1/r_2}
+\gbign_D(E^{\Mhat}_{2\,\ast},
 E^{\Mhat}_{1\,\ast})\otimes I_{-1-r_1/r_2}
\end{multline}
\item
We have the following equality
on $\nbigs$:
\[
 F^{\ast}\gbign(L_2\otimes\nbigi^{(2)},E^{\Mhat}_1)
=G^{\ast}\bigl(
 \gbign(L_2,E_1^u)
 \otimes \nbigo_{1,\rel}(1+r_1/r_2)
 \otimes \nbigo_{2,\rel}(1)
\bigr)
\]
Here, we put $\gbign(L_2,E_1^u):=Rp_{X\,\ast}\nhom(L_2,E_1^u)$.
\hfill\qed
\end{itemize}
\end{lem}

\section{Virtual fundamental Classes}
\label{section;06.7.3.40}

\subsection{Perfectness of the Obstruction Theories
for some Stacks}
\label{subsection;06.7.3.41}
\subsubsection{The moduli stacks}
\label{subsubsection;06.6.13.36}

Let $\vecy$ be an element of $\Type$,
and let $\alpha_{\ast}$  be a system of weights.
Let $L$ be a line bundle on $X$.
We use the notation in the subsection
\ref{subsection;06.6.13.5}.
The proof of the following propositions
will be given in the subsubsection 
\ref{subsubsection;06.6.13.121}
after the preparation in the subsubsection
\ref{subsubsection;06.6.13.55}.
The expected dimensions can be calculated
formally. 
We give the results in the subsubsection 
\ref{subsubsection;06.6.13.125}.

\begin{prop}
 \label{prop;06.5.5.40}
Let $m$ be a sufficiently large integer.
\begin{itemize}
\item
 The obstruction theory
 $\Ob(m,\vecyhat)$ of $\nbigm^{ss}(\vecyhat,\alpha_{\ast})$
 is perfect in the sense of Definition
 {\rm\ref{df;06.6.13.12}}.
\item
 Let $\delta$ be an element of $\nbigp^{\br}$.
 The obstruction theory
 $\Ob(m,\vecy,L)$ of $\nbigm^{ss}(\vecy,L,\alpha_{\ast},\delta)$
 is perfect.
\end{itemize}
\end{prop}

\begin{prop}
\label{prop;06.6.13.21}
Assume $\rank(\vecy)>1$.
Let $\delta$ be an element of $\nbigp^{\br}$.
Let $m$ be a sufficiently large integer.
\begin{itemize}
\item
 The obstruction theory 
 $\Ob(m,\vecyhat,L)$ of
 $\nbigm^{ss}(\vecyhat,L,\alpha_{\ast},\delta)$
 is perfect.
\item
 The obstruction theory 
 $\Ob(m,\vecyhat,[L])$ of
 $\nbigm^{ss}(\vecyhat,[L],\alpha_{\ast},\delta)$
 is perfect.
\end{itemize}
\end{prop}

\begin{prop}
\label{prop;06.5.11.550}
Assume $\rank(\vecy)=1$.
Let $m$ be a sufficiently large integer.
\begin{itemize}
\item
 We have the vanishing
 $\nbigh^i\bigl(\Ob(m,\vecyhat,[L])\bigr)$
 unless $i=0$.
 In particular,
 the moduli $\nbigm(\vecyhat)$ is smooth.
\item
 If the $2$-vanishing condition is satisfied for $(\vecy,L)$,
then the obstruction theory
$\Ob(m,\vecyhat,L)$ 
of $\nbigm(\vecyhat,L)$ is perfect,
and the obstruction theory
$\Ob(m,\vecyhat,[L])$ 
of $\nbigm(\vecyhat,[L])$ is perfect.
\end{itemize}
We remark $\Ob(m,\vecy,L)$ is always perfect
as in Proposition {\rm\ref{prop;06.5.5.40}}.
\end{prop}

\begin{prop}
 \label{prop;06.6.13.17}
Let $\vecL=(L_1,L_2)$ be a pair of line bundles on $X$.
Let $\vecdelta=(\delta_1,\delta_2)$ be a pair of
sufficiently small parameters $\delta_i$
as in Lemma {\rm\ref{lem;06.6.13.11}}.
Let $\alpha_{\ast}$ be a system of weights.
Assume that the $2$-vanishing condition
holds for $(\vecy,L_2,\alpha_{\ast})$.
Then, 
the obstruction theory
$\Ob(m,\vecyhat,[\vecL])$ is perfect
on $\nbigm^{ss}(\vecyhat,[\vecL],\alpha_{\ast},\vecdelta)$.
\end{prop}

\begin{notation}
\index{$[\nbigm]$}
 Due to Proposition
 {\rm\ref{prop;06.5.5.40}--\ref{prop;06.6.13.17}},
 we obtain the perfect obstruction theories
 of the moduli stacks $\nbigm$ of the corresponding stable objects,
 which induces the virtual fundamental classes
 due to Proposition {\rm\ref{prop;06.6.13.20}}.
 They are denoted by $[\nbigm]$.
 We use the notation
 $\int_{\nbigm}\Phi$ for the evaluation of
 a cohomology class via $[\nbigm]$.
 (See the subsection {\rm\ref{subsection;06.6.19.110}}.)
\hfill\qed
\end{notation}

\subsubsection{The master space and the related stacks}

We also obtain the following propositions.
\begin{prop}
The obstruction theory $\Ob(\Mhat)$ of $\Mhat$
is perfect.
(See the subsubsections {\rm\ref{subsubsection;06.5.8.5}},
 {\rm\ref{subsubsection;06.5.21.600}},
 {\rm\ref{subsubsection;06.6.13.16}}
 for $\Ob(\Mhat)$.)
\end{prop}
\pf
We consider the obstruction theory for 
the enhanced master space $\Mhat$
given in the subsubsection \ref{subsubsection;06.5.8.5}.
We have the naturally defined smooth morphism
$\gminip:\Mhat\lrarr\nbigm(m,\vecyhat,[L])$.
We remark that the image of $\gminip$ is contained
in the open substack
$\nbigm:=\nbigm^{ss}(\vecyhat,[L],\alpha_{\ast},\delta)$.
Then, the claim immediately follows from
the diagram (\ref{eq;06.6.13.1})
and Proposition \ref{prop;06.6.13.21}.

We obtain the perfectness of
the obstruction theories $\Ob(\Mhat)$
(subsubsection \ref{subsubsection;06.5.21.600})
by the the same argument.
We obtain the perfectness of
the obstruction theories $\Ob(\Mhat)$
in the subsubsection \ref{subsubsection;06.6.13.16}
by using the same argument and Proposition \ref{prop;06.6.13.17}.
\hfill\qed

\vspace{.1in}
The following claims can be shown by a similar argument.

\begin{prop}
\label{prop;06.6.13.52}
\mbox{{}}
\begin{itemize}
\item
The obstruction theory
$\obtilde(m,\vecyhat)$ of
$\nbigmtilde^s(\vecyhat,\alpha_{\ast},+)$
is perfect.
\item
The obstruction theory 
$\obtilde\bigl(m,\vecy,L\bigr)$
of the moduli stack
$\nbigmtilde^s\bigl(\vecyhat,L,\alpha_{\ast},(\delta,\ell)\bigr)$
is perfect.
\item
Assume one of the following:
\begin{itemize}
\item$\rank(\vecy)>1$.
\item $\rank(\vecy)=1$ and
the $2$-vanishing condition holds for $(\vecy,L,\alpha_{\ast},\delta)$.
\end{itemize}
The obstruction theory 
$\ob\bigl(m,\vecyhat,[L]\bigr)$
of the moduli stack
$\nbigmtilde^s\bigl(\vecyhat,[L],\alpha_{\ast},(\delta,\ell)\bigr)$
is perfect.
\hfill\qed
\end{itemize}
\end{prop}

\vspace{.1in}

Recall that we obtained obstruction theory 
$\obtilde(\Mhat^{G_m}(\gbigi))$
of $\Mhat^{G_m}(\gbigi)\times A^1$ over $A^1$
in Proposition \ref{prop;06.5.22.100}.
The specialization at $t=a$ is denoted by
$\obtilde_a(\Mhat^{G_m}(\gbigi))$.
\begin{prop}
The obstruction theories
$\Obtilde_a(\Mhat^{G_m}(\gbigi))$ are perfect
for any $a\in k$.
\end{prop}
\pf
Due to the distinguished triangle (\ref{eq;06.5.22.165}),
we have only to show that
$\Obtilde(m,\vecy_1,L)$ and 
$\Obtilde(m,\vecyhat_2)$ gives the obstruction
theories of
$\nbigmtilde^s\bigl(\vecy_1,L,\alpha_{\ast},(\delta,k_0)\bigr)$
and $\nbigmtilde^s(\vecyhat_2,\alpha_{\ast},+)$
as in Proposition \ref{prop;06.6.13.52}.
\hfill\qed

\begin{notation}
We obtain the virtual fundamental
classes of $M^{G_m}(\gbigi)$ with respect to
the obstruction theories
$\Ob_a(M^{G_m}(\gbigi))$ for each $a\in k$.
It is independent of a choice of $a$
(See Proposition {\rm 7.2} of {\rm\cite{bf}}.)
Therefore, we denote it by $[M^{G_m}(\gbigi)]$.
\hfill\qed
\end{notation}

We use the notation in the subsubsection
\ref{subsubsection;06.5.17.20}.
We put
$\nbigmtilde_{\spl}:=
 \nbigmtilde^{ss}\bigl(\vecy_1,L,\alpha_{\ast},(\delta,k_0)\bigr)
\times
 \nbigmtilde^{ss}(\vecyhat_2,\alpha_{\ast},+)$.
\begin{prop}
\label{prop;06.5.11.300}
In the diagram {\rm(\ref{eq;06.5.17.21})},
we have the relation
$F^{\ast}\bigl([\Mhat^{G_m}(\gbigi)]\bigr)=
 G^{\prime\ast}\bigl([\nbigmtilde_{\spl}]\bigr)$.
\end{prop}
\pf
It follows from the diagram (\ref{eq;06.5.22.160}).
\hfill\qed

\vspace{.1in}
We can show the following propositions
by the same argument.
\begin{prop}
\label{prop;06.6.10.20}
Assume that the $2$-stability condition holds
for $(\vecy,L,\alpha_{\ast},\delta)$.
\begin{itemize}
\item
 We have the perfectness of
 the obstruction theories
 $\obtilde_a(\Mhat^{G_m}(\gbigi))$
 in Proposition {\rm\ref{prop;06.6.10.11}}.
\item
 They give the virtual fundamental class
 $[\Mhat^{G_m}(\gbigi)]$.
\item
 We have the relation 
 $F^{\ast}\bigl([\Mhat^{G_m}(\gbigi)]\bigr)
=G^{\prime\,\ast}\bigl([\nbigm^{s}(\vecy_1,L,\alpha_{\ast},\delta)
 \times \nbigm^s(\vecyhat_2,\alpha_{\ast})]\bigr)$
in the diagram {\rm(\ref{eq;06.5.22.300})}.
\hfill\qed
\end{itemize}
\end{prop}

\begin{prop}
\label{prop;06.6.10.25}
Under the situation of the subsubsection
{\rm\ref{subsubsection;06.6.12.1}},
the following claims hold:
\begin{itemize}
\item
 We have the perfectness of
 the obstruction theories
 $\obtilde_a(\Mhat^{G_m}(\gbigi))$
 in Proposition {\rm\ref{prop;06.6.10.15}}.
\item
 They give the virtual fundamental class
 $[\Mhat^{G_m}(\gbigi)]$.
\item
 We have the relation 
 $F^{\ast}\bigl([\Mhat^{G_m}(\gbigi)]\bigr)
=G^{\prime\,\ast}\bigl([\nbigm^{s}(\vecy_1,L,\alpha_{\ast},\delta)
 \times \nbigm^s(\vecyhat_2,\alpha_{\ast})]\bigr)$
in the diagram {\rm(\ref{eq;06.6.22.10})}.
\hfill\qed
\end{itemize}
\end{prop}

\subsubsection{Vanishing of some cohomology groups}
\label{subsubsection;06.6.13.55}

We use the notation in the section \ref{section;06.6.4.250}.
Let $E$ be a torsion-free sheaf on $X$.
Let $V_{\cdot}=(V_{-1}\rarr V_0)$  be
a locally free resolution of $E$.
Then we put 
$C(E):=\gminig^{\lor}(V_{\cdot})$,
where $\gminig^{\lor}(V_{\cdot})$ denotes
the dual of $\gminig^{\lor}(V_{\cdot})$
as $\nbigo_X$-complexes.
We also put
$C^{\circ}(E):=\gminig^{\circ\,\lor}(V_{\cdot})$
and $C^{d}(E):=\gminig^{d\,\lor}(V_{\cdot})$.

Let $F_{\ast}$ be a quasi-parabolic structure of $E$ at $D$.
We have $\gminig(V_{\cdot|D})$
and $\gminig_D(V_{\cdot},F_{\ast})$ 
with the natural morphism
$\gminig_D(V_{\cdot},F_{\ast})\lrarr
 \gminig(V_{\cdot|D})$ on $D$.
(See the subsubsection \ref{subsubsection;06.5.5.5}.)
We have the dual complexes
$\gminig^{\lor}(V_{\cdot|D})$ 
and $\gminig^{\lor}_D(V_{\cdot},F_{\ast})$
as $\nbigo_D$-complexes.
We have the natural morphism
$C(E)\lrarr
 \gminig^{\lor}(V_{\cdot})_{|D}
=\gminig^{\lor}(V_{\cdot|D})$.
Thus we have the morphism
$\alpha:
 C(E)\oplus\gminig_D^{\lor}(V_{\cdot},F_{\ast})
\lrarr \gminig^{\lor}(V_{\cdot|D})$.
We put
$C(E,F_{\ast}):=\cone(\alpha)[-1]$.

Recall we have the decompositions
$\gminig(V_{\cdot|D})=\gminig^{\circ}(V_{\cdot|D})
\oplus\gminig^{d}(V_{\cdot|D})$ and
$\gminig_D(V_{\cdot},F_{\ast})=
 \gminig_D^{\circ}(V_{\cdot},F_{\ast})\oplus
 \gminig_D^d(V_{\cdot},F_{\ast})$.
(See the subsubsection \ref{subsubsection;06.5.11.200}.)
Thus, we obtain the complexes
$C^{\circ}(E,F_{\ast})$ and 
$C^{d}(E,F_{\ast})$.
It is easy to see $C^d(E,F_{\ast})=C^d(E)$.

\begin{lem}
 \label{lem;06.5.5.20}
Let $(E,F_{\ast})$ be as above.
The hyper-cohomology group
$\hyperh^i\bigl(X,C^{\circ}(E,F_{\ast})\bigr)$ vanishes
unless $i=-1,0,1$.
If $(E,F_{\ast})$ is stable
with respect to some system of weights $\alpha_{\ast}$,
we also have $\hyperh^{-1}(X,C^{\circ}(E,F_{\ast}))=0$.
\end{lem}
\pf
The $i$-th cohomology sheaf
$\nbigh^i\bigl(C^{\circ}(V_{\cdot},F_{\ast})\bigr)$ vanishes
unless $-2\leq i\leq 1$ by construction.
It is easy to see that
the morphisms
$C^{\circ}(E)^{-2}\lrarr C^{\circ}(E)^{-1}$ and
$\gminig_D^{\circ\lor}(V_{\cdot},F_{\ast})^{-2}
 \lrarr \gminig_D^{\circ\lor}(V_{\cdot},F_{\ast})^{-1}$
are injective.
Therefore, 
we have $\nbigh^{-2}\bigl(C^{\circ}(E,F_{\ast})\bigr)=0$.
Since the morphism
$C^{\circ}(E)^1\lrarr \gminig^{\circ\lor}(V_{\cdot})^{1}$
is surjective,
we have $\nbigh^1\bigl(C^{\circ}(E,F_{\ast})\bigr)=0$.
Let $P\not\in D$ be a point where $E$ is locally free.
Then we have $\nbigh^0\bigl(C^{\circ}(E,F_{\ast})\bigr)=0$
around $P$.
Thus the support of the sheaf
$\nbigh^0\bigl(C^{\circ}(E,F_{\ast})\bigr)$
is $1$-dimensional.

Then, we obtain the vanishing
$\hyperh^i\bigl(X,C^{\circ}(E,F_{\ast})\bigr)$
unless $i=-1,0,1$ by using the spectral sequence.
It is easy to see that
$\hyperh^{-1}\bigl(X,C^{\circ}(E,F_{\ast})\bigr)$
is the set of endomorphisms of $(E,F_{\ast})$
whose trace is $0$.
If $(E,F_{\ast})$ is assumed to be stable with respect to
some weight $\alpha_{\ast}$,
we obtain $\hyperh^{-1}\bigl(X,C^{\circ}(E,F_{\ast})\bigr)=0$.
\hfill\qed

\vspace{.1in}
In the rank one case,
we have the following result.
\begin{lem}
\label{lem;06.5.5.55}
Let $(E,F_{\ast})$ be as above,
and we assume $\rank(E)=1$.
Then, we have 
$\hyperh^i\bigl(X,C^{\circ}(E,F_{\ast})\bigr)=0$
unless $i=0$.
\end{lem}
\pf
Let $x$ be a point of $X$ with one of the following:
\begin{itemize}
\item
 $x$ is contained in $X-D$, and  $E$ is locally free around $x$.
\item
 $x$ is contained in $D$, 
 $E$ is locally free around $x$ as $\nbigo_X$-module,
 and $\Cok_i(E)$ are locally free around $x$
 as $\nbigo_D$-module.
\end{itemize}
Around such a point, we can compute
the cohomology sheaves of $C^{\circ}(E,F_{\ast})$
in the case $V_0=E$ and $V_{-1}=0$.
Hence it is easy to check $C^{\circ}(E,F_{\ast})\simeq 0$
around such a point.

We know
$\nbigh^i\bigl(C^{\circ}(E,F_{\ast})\bigr)=0$
unless $i=-1,0$, in general.
By the above consideration,
we know that the support of
$\nbigh^i\bigl(C^{\circ}(E,F_{\ast})\bigr)$
is $0$-dimensional.
Therefore,
we obtain that
$\hyperh^i\bigl(X,C^{\circ}(E,F_{\ast})\bigr)=0$
unless $i=-1,0$.
Since $(E,F_{\ast})$ is always stable
in the case $\rank(E)=1$.
we also obtain the vanishing
$\hyperh^{-1}\bigl(X,C^{\circ}(E,F_{\ast})\bigr)=0$.
Thus we are done.
\hfill\qed

\vspace{.1in}

Let $L$ be a line bundle on $X$.
Let $\phi$ be an $L$-section of $E$.
We take a locally free resolution
$P_{\cdot}=(P_{-1}\lrarr P_0)$
of the line bundle $L$ so that we have a lift
$\phitilde:P_{\cdot}\lrarr V_{\cdot}$ of $\phi$.
We have the complex
$\gminig_{\rel}(V_{\cdot},\phitilde)
=\nhom(P_{\cdot},V_{\cdot})^{\lor}$
and the dual complex
$\gminig_{\rel}^{\lor}(V_{\cdot},\phitilde)$.
We have the natural morphism
$\gminig_{\rel}(V_{\cdot},\phitilde)[-1]
\lrarr \gminig(V)$.
Thus we have the morphism
$\gamma_L:
 C(E)\lrarr \gminig^{\lor}_{\rel}(V_{\cdot},\phitilde)[1]$.
It induces the morphisms
$\alpha_L:
 C(E,F_{\ast})\lrarr \gminig^{\lor}_{\rel}(V_{\cdot},\phitilde)[1]$
and 
$\alpha_L^{\circ}:C^{\circ}(E,F_{\ast})
 \lrarr \gminig^{\lor}_{\rel}(V_{\cdot},\phitilde)[1]$.
We put
$C(E,F_{\ast},\phi):=\cone(\alpha_L)[-1]$
and $C^{\circ}(E,F_{\ast},\phi):=\cone(\alpha_L^{\circ})[-1]$.
It is easy to see that $C(E,F_{\ast},\phi)$ and 
$C^{\circ}(E,F_{\ast},\phi)$ are well defined in $D(X)$.

\begin{lem}
 \label{lem;06.5.5.25}
Let $(E,F_{\ast},\phi)$ be as above.
Assume $\phi\neq 0$.
Then the hyper-cohomology groups
$\hyperh^i\bigl(X,C(E,F_{\ast},\phi)\bigr)$
vanish unless $i=-1,0,1$.
If $(E,F_{\ast},\phi)$ is $(\alpha_{\ast},\delta)$-stable
for some $\delta\in\nbigp^{\br}$
and some weight $\alpha_{\ast}$,
we have
$\hyperh^{-1}\bigl(X,C(E,F_{\ast},\phi)\bigr)=0$.
\end{lem}
\pf
We have $\nbigh^i\bigl(C(E,F_{\ast},\phi)\bigr)=0$
unless $-2\leq i\leq 1$ by construction.
We also have $\nbigh^{-2}\bigl(C(E,F_{\ast},\phi)\bigr)=0$
as in the proof of Lemma \ref{lem;06.5.5.20}.
Since the morphism
$\nhom(P_0,V_0)\lrarr \nhom(P_{-1},V_0)$ is
surjective,
we obtain the vanishing
$\nbigh^1\bigl(C(E,F_{\ast},\phi)\bigr)=0$.
Let $x$ be a point of $X-D$
with $\phi(x)\neq 0$
such that $E$ is locally free around $x$.
Then, it is easy to show the vanishing of
$\nbigh^0\bigl(C(E,F_{\ast},\phi)\bigr)$
around $x$.
Therefore, the support of
$\nbigh^0\bigl(C(E,F_{\ast},\phi)\bigr)$ is
one dimensional.
Then we obtain 
$\hyperh^i\bigl(X,C(E,F_{\ast},\phi)\bigr)=0$
unless $i=-1,0,1$ by using the spectral sequence.
Since
 $\hyperh^{-1}\bigl(X,C(E,F_{\ast},\phi)\bigr)$
is the set of endomorphisms of $(E,F_{\ast},\phi)$,
we obtain
$\hyperh^{-1}\bigl(X,C(E,F_{\ast},\phi)\bigr)=0$
from the stability assumption of $(E,F_{\ast},\phi)$.
\hfill\qed

\begin{lem}
 \label{lem;06.5.5.30}
Let $(E,F_{\ast},\phi)$ be as above.
We assume $\phi\neq 0$ and $\rank(E)>1$.
Then, the hyper-cohomology groups
$\hyperh^i\bigl(X,C^{\circ}(E,F_{\ast},\phi)\bigr)$ vanish
unless $i=-1,0,1$.
\end{lem}
\pf
The $i$-th cohomology sheaf 
$\nbigh^i\bigl(C^{\circ}(E,F_{\ast},\phi)\bigr)$ vanish
unless $-2\leq i\leq 1$ by construction.
We also have the vanishings
of $\nbigh^{j}\bigl(C^{\circ}(E,F_{\ast},\phi)\bigr)$
$j=-2,1$, as in the proof of Lemma \ref{lem;06.5.5.25}.
Let $x$ be a point of $X-D$ with $\phi(x)\neq 0$
such that $E$ is locally free around $x$.
Under our assumption $\rank(E)>1$,
we can easily check the vanishing
$\nbigh^0\bigl(C^{\circ}(E,F_{\ast},\phi)\bigr)=0$
around such a point $x$.
Thus the claim can be shown
by using the spectral sequence.
\hfill\qed

\vspace{.1in}
In the rank one case,
we have the following:
\begin{lem}
\label{lem;06.5.11.150}
Let $(E,F_{\ast},\phi)$ be as above.
We assume $\phi\neq 0$ and $\rank(E)=1$.
Moreover, we assume 
$H^2\bigl(X,L^{-1}\otimes E\bigr)=0$.
Then, the hyper-cohomology groups
$\hyperh^i(X,C^{\circ}(E,F_{\ast},\phi))$ vanish
unless $i=0,1$.
\end{lem}
\pf
We have the distinguished triangle
$\gminig_{\rel}^{\lor}(V,\phitilde)
\lrarr C^{\circ}(E,F_{\ast},\phi)
\lrarr C^{\circ}(E,F_{\ast})
\lrarr \gminig_{\rel}^{\lor}(V,\phitilde)[1]$.
Hence we obtain the long exact sequence:
\[
\cdots\lrarr
 H^{i}\bigl(X,L^{-1}\otimes E\bigr)
\lrarr
 \hyperh^i\bigl(X,C^{\circ}(E,F_{\ast},\phi)\bigr)
\lrarr
 \hyperh^i\bigl(X,C(E,F_{\ast})\bigr)
\lrarr
 H^{i+1}\bigl(X,L^{-1}\otimes E\bigr)
\lrarr\cdots
\]
Then the claim follows from the assumption
and Lemma \ref{lem;06.5.5.55}.
\hfill\qed

\vspace{.1in}

Let $\vecL=(L_1,L_2)$ be a pair of line bundles on $X$.
Let $\phi_i$ be $L_i$-sections of $E$.
By taking appropriate resolutions $P_{\cdot}^{(i)}$ of $L_i$
and lifts $\phitilde_i:P_{\cdot}^{(i)}\lrarr V_{\cdot}$
of $\phi_i$,
we obtain the complexes
$C(E,F_{\ast},\phitilde_1,\phitilde_2)$
and $C^{\circ}(E,F_{\ast},\phitilde_1,\phitilde_2)$
as above.
They are well defined in $D(X)$.

\begin{lem}
\label{lem;06.5.5.35}
Assume 
$H^2\bigl(X,L_2^{-1}\otimes E\bigr)=0$
and $\phi_j\neq 0$ for $j=1,2$.
We also assume $\rank(E)>1$.
Then, we have
$\hyperh^i\bigl(X,C^{\circ}(E,F_{\ast},\phi_1,\phi_2)\bigr)=0$
unless $i=-1,0,1$.
If $(E,F_{\ast},\phi_1,\phi_2)$ is
$(\alpha_{\ast},\delta_1,\delta_2)$-stable,
then we also have
$\hyperh^{-1}\bigl(X,C^{\circ}(E,F_{\ast},\phi_1,\phi_2)\bigr)=0$.
\end{lem}
\pf
We have the exact sequence:
$ 0\lrarr
 \nhom\bigl( P^{(2)},V_{\cdot} \bigr)
\lrarr
 C(E,F_{\ast},\phi_1,\phi_2)
\lrarr
 C(E,F_{\ast},\phi_1)
\lrarr 0$.
Then the claims can be reduced
to Lemma \ref{lem;06.5.5.25}
and Lemma \ref{lem;06.5.5.30}.
\hfill\qed

\subsubsection{Proof of the propositions in the subsubsection 
 \ref{subsubsection;06.6.13.36} }
\label{subsubsection;06.6.13.121}

Proposition \ref{prop;06.5.5.40}
immediately follows from the following lemma.
\begin{lem}
\label{lem;06.6.13.30}
The obstruction theory
$\Ob(m,\vecyhat)$ of $\nbigm(m,\vecyhat)$
is perfect in the sense of Definition
{\rm\ref{df;06.6.13.12}}.
The obstruction theory 
$\Ob(m,\vecy,L)$ of $\nbigm(m,\vecy,L)$
is perfect.
\end{lem}
\pf
Let us discuss $\Ob(m,\vecyhat)$.
We would like to show
$\Ob(m,\vecyhat)$ is quasi-isomorphic to
a complex $E^{-1}\rarr E^{0}\rarr E^1$ 
of locally free sheaves on $\nbigm(m,\vecyhat)$.
Since $\Ob(m,\vecyhat)$ is obtained 
from the push-forward of
perfect complexes on
$\nbigm(m,\vecyhat)\times X$
and $\nbigm(m,\vecyhat)\times D$,
it is easy to show that 
$\Ob(m,\vecyhat)$ is quasi isomorphic to
a bounded complex of locally free sheaves
by using the projectivity of $X$ and $D$.

Therefore, we have only to check
$H^i\bigl(i_z^{\ast}\Ob(m,\vecyhat)\bigr)=0$
unless $i=-1,0,1$
for any point $z\in\nbigm(m,\vecyhat)$,
where $i_z$ denotes the inclusion of $z$
to $\nbigm(m,\vecyhat)$.
Let $(E,F_{\ast},\rho)$ be the parabolic
oriented torsion-free sheaf corresponding to $z$.
Then the dual of
$\nbigh^i\bigl(i_z^{\ast}\Ob(m,\vecyhat)\bigr)$
is isomorphic to
$\hyperh^{-i}\bigl( X,C^{\ori}(E,F_{\ast}) \bigr)
\oplus H^1(X,\nbigo)[0]$.
Therefore, the claim follows from
Lemma \ref{lem;06.5.5.20}.

The perfectness of $\Ob(m,\vecy,L)$
can be shown similarly,
by using Lemma \ref{lem;06.5.5.25}.
\hfill\qed

\vspace{.1in}

\begin{lem}
\label{lem;06.6.13.39}
Assume $\rank(\vecy)>1$.
 The obstruction theory
 $\Ob(m,\vecyhat,L)$ of $\nbigm(m,\vecyhat,L)$
 is perfect.
The obstruction theory $\Ob(m,\vecyhat,[L])$
 of $\nbigm(m,\vecyhat,[L])$ is perfect.
\end{lem}
\pf
We have the following commutative diagram:
\[
 \begin{CD}
 \Ob_{\rel}(m,\yhat)[-1]
 @>>>
 \Ob_{\rel}(m,y,L)\oplus\Ob_{\rel}(m,\yhat)[-1]
 @>>>
 \Ob_{\rel}(m,y,L)[-1]\\
 @VVV @VVV @VVV \\
 \Ob^d(m,\vecy) @>>>
 \Ob(m,\vecy) @>>>
 \Ob^{\circ}(m,\vecy)
 \end{CD}
\]
We put 
$C_1:=\Cone\bigl(
 \Ob_{\rel}(m,\yhat)[-1]\lrarr
 \Ob^d(m,\vecy) \bigr)$
and $C_2:=\Cone\bigl(
 \Ob_{\rel}(m,y,L)[-1]\lrarr
 \Ob^{\circ}(m,\vecy) \bigr)$.
We have only to show that
$C_1$ and $C_2$ are perfect of amplitude
in $[-1,1]$.
It is easy to see that
$C_1$ is isomorphic to
$H^1(X,\nbigo)^{\lor}\otimes\nbigo[0]$.
To check the claim for $C_2$,
we have only to see $i_z^{\ast}C_2$
as in the proof of Lemma \ref{lem;06.6.13.30}.
Then the dual of
$\nbigh^i\bigl(i_z^{\ast}C_2\bigr)$ are
isomorphic to 
$\hyperh^{-i}\bigl(X,C^{\circ}(E,F_{\ast},\phi)\bigr)$.
Thus the claim for $C_2$ follows from
Lemma \ref{lem;06.5.5.30}.
Therefore, we obtain the first claim of the lemma.

Let us show the second claim.
We have only to show the vanishing 
of $\nbigh^i\bigl(\Ob(m,\vecyhat,[L])\bigr)=0$
for $i<-1$.
Let $\pi:\nbigm^s(\vecyhat,L,\alpha_{\ast},\delta)
\lrarr \nbigm^s(\vecyhat,[L],\alpha_{\ast},\delta)$
be the natural morphism,
which is smooth.
Hence,
we have only to show
$\nbigh^i\bigl(\pi^{\ast}\Ob(m,\vecyhat,[L])\bigr)=0$
for $i<-1$.
Since we have
$\nbigh^i\bigl(\pi^{\ast}\Ob(m,\vecyhat,[L])\bigr)
\simeq
 \nbigh^i\bigl(\Ob(m,\vecyhat,L)\bigr)$
for $i<-1$,
the claim follows from Lemma \ref{lem;06.6.13.30}.
\hfill\qed

\vspace{.1in}

Let us see Proposition \ref{prop;06.5.11.550}.
The first claim can be shown 
by using Lemma \ref{lem;06.5.5.55}
and the argument in the proof of 
Lemma \ref{lem;06.6.13.30}.
The second and third claims can be shown
by using Lemma \ref{lem;06.5.11.150}
and the argument in the proof of
Lemma \ref{lem;06.6.13.39}.

Proposition \ref{prop;06.6.13.17}
can be shown by using Lemma \ref{lem;06.5.5.35}
and the argument in the proof of Lemma \ref{lem;06.6.13.39}.
\hfill\qed

\subsubsection{Expected dimension}
\label{subsubsection;06.6.13.125}

We put $p_g:=\dim H^2(X,\nbigo)$
and $\chi(\nbigo_X):=1-\dim H^1(X,\nbigo_X)+\dim H^2(X,\nbigo_X)$.
It is easy to obtain the formulas of the expected dimension
in the non-parabolic case.
Let $y$ be an element of $H^{\ev}(X)$.
The degree $2i$-part of $y$ is denoted by $y_i$.
We put as follows:
\[
 d(y):=\bigl(y_1^2-2y_0\cdot y_2)\cap[X]
-y_0^2\cdot \chi(\nbigo_X)
\]
For a line bundle $L$,
we put as follows:
\[
 d_{\rel}(y,L):=
\bigl(\ch(L^{-1})\cdot y\cdot\Td(X)\bigr)\cap[X]
\]

\begin{prop}
For the moduli of the non-parabolic objects,
the expected dimensions are as follows:
\begin{itemize}
\item
$\dim^f\nbigm(\yhat)= d(y)+1+p_g$.
\item
$\dim^f\nbigm(y,L,\delta)=d(y)+d_{\rel}(y,L)$.
\item
$\dim^f\nbigm(\yhat,L,\delta)=
 d(y)+1+p_g+d_{\rel}(y,L)$.
\item
$\dim^f\nbigm(\yhat,[L],\delta)
=d(y)+p_g+d_{\rel}(y,L)$.
\end{itemize}
\end{prop}
\pf
Let us consider the first case.
For any $(E,\rho,F_{\ast})\in \nbigm(\yhat)$,
the virtual tangent space is $K$-theoretically
given by
$\sum_{i}(-1)^i\hyperh^{i}\bigl(
 X,C^{\circ}(E) \bigr)
+H^1(X,\nbigo_X)$.
The Euler number 
can be easily calculated,
and it is given by 
$d(y)+\chi(\nbigo_X)+\dim H^1(X,\nbigo_X)
=d(y)+1+p_g$.
The second and third cases can be discussed similarly.

Let $\gminip$ denote the natural morphism
$\nbigm(\yhat,L,\delta)\lrarr\nbigm(\yhat,[L],\delta)$.
We have the distinguished triangle
on $\nbigm(m,y,L)$.
\[
 L_{\nbigm(m,\yhat,L)/\nbigm(m,\yhat,[L])}[-1]
\lrarr \gminip^{\ast}\Ob(m,\yhat,[L])\lrarr 
 \Ob(m,\yhat,L)\lrarr L_{\nbigm(m,\yhat,L)/\nbigm(m,\yhat,[L])}
\]
Then the fourth claim is obtained.
\hfill\qed

\vspace{.1in}

The contribution of the parabolic structure
can also be calculated formally.
Let $\vecy=(y,y_1,y_2,\ldots,y_l)$ be an element of $\Type$.
In this case,
$y_i$ can be regarded as elements
of $H^{\ev}(D)$.
We put $s_i+w_i:=\sum_{j\leq i}y_j$,
where $s_i$ (resp. $w_i$) denotes 
the element of $H^{0}(D)$ (resp. $H^2(D)$).
Let $g$ denote the genus of $D$.
Let $d_{\rel}(y,\vecy)$ denote
the Euler number of the complex
$\gminig_{\rel}(V_{\cdot},F_{\ast})$
for $(E,F_{\ast})$ of type $\vecy$.
The result is as follows:
\[
 d_{\rel}(y,\vecy):=
 2g-2
 \sum_{i=1}^l s_i(s_i-s_{i+1})
+\sum_{i=1}^{l-1}
 \int_D\bigl( s_{i+1}\cdot w_i-s_i\cdot w_{i+1}\bigr).
\]
\begin{prop}
\mbox{{}}
\begin{itemize}
\item
$\dim^f\nbigm(\vecyhat,\alpha_{\ast})
=d(y)+1+p_g+d_{\rel}(y,\vecy)$
\item
$\dim^f\nbigm(\vecyhat,[L],\alpha_{\ast},\delta)
=d(y)+p_g+d_{\rel}(y,\vecy)+d_{\rel}(y,L)$.
\hfill\qed
\end{itemize}
\end{prop}

\subsection{Comparison of the Oriented Reduced Case
and the Unoriented Unreduced Case}
\label{subsection;06.7.3.42}
\subsubsection{Statements}

We have the natural morphism
$\kappa:\nbigm(m,\vecyhat,[L])
\lrarr \nbigm(m,\vecy,L)$
which is etale and finite of degree $\rank(\vecy)^{-1}$.
But the obstruction theories are not same,
in the case $p_g:=\dim H^2(X,\nbigo)>0$.
We would like to compare them.

\begin{prop}
\label{prop;06.5.10.50}
We have the following commutative diagram
in the derived category $D(\nbigm(m,\vecyhat,[L]))$:
\[
 \begin{CD}
 \kappa^{\ast}\Ob(m,\vecy,L)
@>>>
\Ob(m,\vecyhat,[L]) \\
@VVV @VVV \\
 \kappa^{\ast}L_{\nbigm(m,\vecy,L)}
@>{\simeq}>>
 L_{\nbigm(m,\vecyhat,[L])}
 \end{CD}
\]
We also have the following distinguished triangle:
\[
\begin{CD}
 \kappa^{\ast}\Ob(m,\vecy,L)
@>>>
 \Ob(m,\vecyhat,[L])
@>>>
 H^2(X,\nbigo)^{\lor}
\otimes\nbigo_{\nbigm(m,\vecyhat,[L])}[1]
 @>>>
 \kappa^{\ast}\Ob(m,\vecy,L)[1]
\end{CD}
\]
\end{prop}

A proof will be given in the next subsubsection.
Before going into the proof,
we give some consequences.

\vspace{.1in}

We have the natural morphism
$\nbigm^s(\vecyhat,[L],\alpha_{\ast},\delta)
\lrarr
 \nbigm^s(\vecy,L,\alpha_{\ast},\delta)$
which is etale and finite of degree $\rank(\vecy)^{-1}$.
It is also denoted by $\kappa$.

\begin{prop}
\label{prop;06.5.31.1}
Assume $\rank(\vecy)>1$.
In the case $p_g>0$,
we have the vanishing
$\bigl[\nbigm^s(\vecy,L,\alpha_{\ast},\delta)\bigr]=0$.
In the case $p_g=0$,
we have the following relation:
\[
 \kappa^{\ast}\bigl(
\bigl[\nbigm^s(\vecy,L,\alpha_{\ast},\delta)\bigr]
\bigr)=
\bigl[\nbigm^s(\vecyhat,[L],\alpha_{\ast},\delta)\bigr]
\]
\end{prop}
\pf
Under the assumption $\rank(\vecy)>1$,
the obstruction theory $\ob(m,\vecyhat,[L])$ 
of $\nbigm^s(\vecyhat,[L],\alpha_{\ast},\delta)$
is perfect.
From Proposition \ref{prop;06.5.10.50},
we have the relation:
\[
 \Eu\bigl(
 H^2(X,\nbigo)\otimes\nbigo
 \bigr)\cap
 \bigl[\nbigm^s(\vecyhat,[L],\alpha_{\ast},\delta)\bigr]
=\bigl[\nbigm^s(\vecy,L,\alpha_{\ast},\delta)\bigr]
\]
Then the claim is clear.
\hfill\qed

\vspace{.1in}
When $\rank(\vecy)=1$,
the obstruction theory of
$\nbigm(\vecyhat,[L])$ is not perfect,
and hence, a similar vanishing result does not hold,
in general.
But, we obtain the following proposition
by the same argument.
Since the stability condition is trivial,
we omit to denote ``$s$'',
$\alpha_{\ast}$ and $\delta$.

\begin{prop} 
\label{prop;06.5.11.500}
Let $\vecy$ be an element of $\Type$
such that $\rank(\vecy)=1$.
\begin{itemize}
\item
In the case $p_g>0$,
we assume 
$H^2\bigl(X,L^{-1}\otimes E\bigr)=0$
for any $L$-Bradlow pair $(E_{\ast},\phi)\in\nbigm(\vecy,L)$.
Then we have
$\bigl[\nbigm(\vecy,L)\bigr]=0$.
\item
In the case $p_g=0$,
we have $\kappa^{\ast}[\nbigm(\vecy,L)]
=[\nbigm\bigl(\vecyhat,[L]\bigr)]$.
\end{itemize}
We remark that the assumption
in the first claim always holds 
in the case $p_g=0$.
\hfill\qed
\end{prop}

Similarly,
we have the natural morphism
$\kappa:\nbigmtilde^s(\vecyhat,[L],\alpha_{\ast},\delta,\ell)
\lrarr
 \nbigmtilde^s(\vecy,L,\alpha_{\ast},\delta,\ell)$,
which is etale and finite of degree $1/\rank(\vecy)$.
By the same argument,
we obtain the following  propositions.
\begin{prop}
\label{prop;06.5.11.301}
Assume $\rank(\vecy)>1$.
In the case $p_g>0$,
we have the vanishing
$[\nbigmtilde^s(\vecy,L,\alpha_{\ast},\delta,\ell)]=0$.
In the case $p_g=0$,
we have the relation
$\kappa^{\ast}\bigl[
 \nbigmtilde^s(\vecy,L,\alpha_{\ast},\delta,\ell)
 \bigr]
=\bigl[\nbigmtilde^s(\vecyhat,[L],\alpha_{\ast},\delta,\ell)\bigr]$.
\hfill\qed
\end{prop}

\begin{prop}
\label{prop;06.5.11.302}
Assume $\rank(\vecy)=1$.
\begin{itemize}
\item
In the case $p_g>0$,
we assume
$H^2\bigl(X,L^{-1}\otimes E\bigr)=0$
for any $(E_{\ast},\phi)\in \nbigm(\vecy,L)$.
Then, we have the vanishing
$[\nbigmtilde^s(\vecy,L,\alpha_{\ast},\delta,\ell)]=0$.
\item
If $p_g=0$,
we have the relation
$\kappa^{\ast}\bigl[
 \nbigmtilde^s(\vecy,L,\alpha_{\ast},\delta,\ell)
 \bigr]
=\bigl[\nbigmtilde^s(\vecyhat,[L],\alpha_{\ast},\delta,\ell)\bigr]$.
\hfill\qed
\end{itemize}
\end{prop}

We obtain the following vanishing result
for the virtual fundamental class
of the fixed point set of the master space.
\begin{prop}
Let $\gbigi=(\vecy_1,\vecy_2,I_1,I_2)$ be a decomposition type
(Definition {\rm\ref{df;06.5.22.1}}).
Assume $p_g>0$.
\begin{itemize}
\item
In the case $\rank(\vecy_1)>1$,
we have $[M^{G_m}(\gbigi)]=0$.
\item
In the case $\rank(\vecy_1)=1$,
we assume $H^2\bigl(X, L^{-1}\otimes E\bigr)=0$
for any torsion-free sheaf $E$ of type $y$,
moreover.
Then we have $[M^{G_m}(\gbigi)]=0$.
\end{itemize}
\end{prop}
\pf
It follows from the propositions \ref{prop;06.5.11.300},
\ref{prop;06.5.11.301} and \ref{prop;06.5.11.302}.
\hfill\qed

\subsubsection{Proof of Proposition \ref{prop;06.5.10.50}}

The canonical resolutions of the universal sheaves on
$\nbigm(m,\vecy,[L])\times X$,
$\nbigm(m,\vecy,L)\times X$
and $\nbigm(m,\vecyhat,[L])\times X$
are denoted by
$\nbigv[L]$, $\nbigv(L)$ and $\nbigvhat[L]$,
respectively.

Let us see $\Ob(m,\vecy,L)$.
From (\ref{eq;06.5.10.20}),
we obtain the following diagram:
\begin{equation}
 \label{eq;06.5.10.35}
 \begin{CD} 
 \Ob_{\rel}(\nbigv(L),\phitilde)[-1]
 @>>>
 \Ob(m,\vecy)\\
 @V{\gminif_1}VV @VVV \\
 \nbigo[-1]@>>>\nbigo[-1]
 \end{CD}
\end{equation}
We put  as follows:
\[
 \Obbar_{\rel}(\nbigv(L),\phitilde):=
 \Cone\bigl(\Ob_{\rel}(\nbigv(L),\phitilde)
 \lrarr \nbigo\bigr)[-1],
\quad
\Obbar(m,\vecy):=
 \Cone\bigl(
 \Ob(m,\vecy)\lrarr\nbigo[-1] \bigr)[-1]
\]
Then we obtain the following commutative diagram:
\[
 \begin{CD}
 \Obbar_{\rel}(\nbigv(L),\phitilde)[-1]
 @>>>
 \Ob_{\rel}(\nbigv(L),\phitilde)[-1]
 @>>>
 L_{\nbigm(m,\vecy,L)/\nbigm(m,\vecy)}[-1]\\
 @VVV @VVV @VVV\\
 \Obbar(m,\vecy)
 @>>>
 \Ob(m,\vecy)
 @>>>
 L_{\nbigm(m,\vecy)}
 \end{CD}
\]
Thus we put as follows:
\[
 \Obbar(m,\vecy,L):=
 \Cone\bigl(
 \Obbar_{\rel}(\nbigv(L),\phi)[-1]
\lrarr
 \Obbar(m,\vecy)
 \bigr)
\]
Then we obtain the morphisms
$\Obbar(m,\vecy,L)
\lrarr
 \Ob(m,\vecy,L)
\lrarr L_{\nbigm(m,\vecy,L)}$.
Since the first morphism is quasi isomorphic,
the composite gives the equivalent obstruction theory.

\vspace{.1in}

On the other hand,
let $\pi_2:\nbigm(m,\vecy,L)\lrarr\nbigm(m,\vecy,[L])$
denote the natural projection.
Due to Lemma \ref{lem;06.5.10.30},
we have the following commutative diagram
on $\nbigm(m,\vecy,L)$
from the diagram (\ref{eq;06.5.10.35}):
\[
 \begin{CD}
 \pi_2^{\ast}\Ob_{\rel}'(\nbigv[L],[\phitilde])[-1]
 @>{\simeq}>>
 \Ob_{\rel}(\nbigv(L),\phitilde)[-1]  
 @>>> \Ob(m,\vecy)\\
 @VVV @VVV @VVV\\ 
 L_{\nbigm(m,\vecy)_{G_m}/\nbigm(m,\vecy)}
 @>{\simeq}>> \nbigo[-1]
 @>>>\nbigo[-1]
 \end{CD}
\]
In particular, we have the isomorphism
\begin{equation}
\label{eq;06.5.10.40}
\pi_2^{\ast}\Ob_{\rel}(\nbigv[L],[\phitilde])
\simeq
\Obbar(\nbigv(L),\phitilde)
\end{equation}
On $\nbigm(m,\vecyhat,[L])$,
we have
$\Ob_{\rel}(\nbigvhat[L],\rho^u):=
 \Cone\Bigl(
 \Ob^d(\nbigvhat[L])\lrarr L_{\Pic}
 \Bigr)$.
The trace map induces the following commutative
diagram:
\[
 \begin{CD}
\Ob_{\rel}(\nbigvhat[L],\rho^u)[-1]
 @>>>
\Ob(m,\vecy)\\
@VVV @V{\tr}VV\\
\nbigo[-1]@>>>\nbigo[-1]
 \end{CD}
\]
We put 
$\Obbar_{\rel}(\nbigvhat[L],\rho^u):=
 \Cone\Bigl(
 \Ob_{\rel}(\nbigvhat[L],\rho^u)
\lrarr\nbigo
 \Bigr)[-1]$.
Then we have the following commutative diagram:
\[
\begin{CD}
 \kappa^{\ast}\pi^{\ast}\Obbar_{\rel}(\nbigv[L],[\phitilde])
 \oplus
 \Obbar_{\rel}(\nbigvhat[L],\rho^u)
 @>>>
 \Obbar(m,\vecy)\\
 @VVV @VVV\\
 \Ob_{\rel}(\nbigvhat[L],[\phitilde])
\oplus
 \Ob_{\rel}(\nbigvhat[L],\rho^u)
 @>>>
 \Ob(m,\vecy)
\end{CD}
\]
We put as follows:
\[
 \Obbar(m,\vecyhat,[L]):=
 \Cone\Bigl(
 \kappa^{\ast}\pi^{\ast}\Obbar_{\rel}(\nbigv[L],[\phitilde])
 \oplus
 \Obbar_{\rel}(\nbigvhat[L],\rho^u)
\lrarr
 \Obbar(m,\vecy)
 \Bigr)
\]
Then we have the morphisms
$\Obbar(m,\vecyhat,[L])
\lrarr \Ob(m,\vecyhat,[L])
\lrarr L_{\nbigm(m,\vecyhat,[L])}$.
Since the first morphism is quasi isomorphic,
the composite gives the equivalent obstruction theory.

By the construction,
we have the following isomorphism:
\[
\Obbar_{\rel}(\nbigvhat[L],\rho^u)
\simeq
 H^2(X,\nbigo)^{\lor}\otimes\nbigo[2]
\]
Then the claim of the proposition immediately follows.
\hfill\qed

\subsection{Rank One Case}
\label{subsection;06.6.1.10}
\subsubsection{The moduli of $L$-abelian pairs}
\label{subsubsection;06.6.30.3}

Let $L$ be a line bundle on $X$.
An $L$-Bradlow pair $(E,\phi)$
is called an $L$-abelian pair, if $E$ is a line bundle.
Similarly,
a reduced $L$-Bradlow pair $(E,[\phi])$
is called a reduced $L$-abelian pair,
if $E$ is a line bundle.

Let $c$ be an element of $H^2(X)$.
We denote by $M(c,L)$
the moduli of $L$-abelian pairs
$(E,\phi)$ such that $c_1(E)=c$.
We denote by $M(\widehat{c},[L])$
the moduli of oriented reduced $L$-abelian pairs
$(E,[\phi],\rho)$ such that
$c_1(E)=c$.
We have the isomorphism 
$M(\widehat{c},[L])\simeq M(c,L)$ of schemes.
We have the projection
$M(c,L)\lrarr \Pic(c)$,
which is a projectivization of a cone
over $\Pic(c)$.

The universal object on
$M(\widehat{c},[L])\times X$ is denoted by
$(\nbiglhat^u,[\phi^u])$.
The line bundle $\nbiglhat^u$ is the pull back
of the Poincar\'e bundle on $\Pic(c)\times X$.
On the other hand,
the universal object on $M(c,L)\times X$ is denoted by
$(\nbigl^u,\phi^u)$.
We have the relation
$\nbigl^u=\nbiglhat^u\otimes\nbigo_{\rel}(1)$.

The obstruction theories of $M(c,L)$
and $M(\widehat{c},[L])$ are denoted by
$\Ob\bigl(M(c,L)\bigr)$ and 
$\Ob\bigl(M(\widehat{c},[L])\bigr)$
respectively.
They are not equivalent in general.
In fact, $\Ob\bigl(M(\widehat{c},[L])\bigr)$
is not perfect in the case $p_g>0$ in general,
unless $H^2(X,\nbigl)=0$
for any $\nbigl$ such that $c_1(\nbigl)=c$.

\vspace{.1in}
In the case $H^1(X,\nbigo_X)=0$,
we have a simple description of the moduli
of abelian pairs.
We work on the complex number fields.
Let $c$ be an element of $H^2(X,\seisuu)$.
Let $\nbigl$ be a line bundle on $X$
such that $c_1(\nbigl)=c$.
Due to the assumption $H^1(X,\nbigo_X)=0$,
any line bundle $\nbigl'$ with $c_1(\nbigl')$
is isomorphic to $\nbigl$.
Hence, the moduli $M(c,L)$ is isomorphic
to the projective space
$\proj\bigl(H^0(X,L^{-1}\otimes\nbigl)^{\lor}\bigr)$.

Since the moduli is smooth,
the obstruction theory of $M(c,L)$
is given by the obstruction bundle $O(c)$,
and the virtual fundamental class is the Euler class of $O(c)$.
The vector bundle $O(c)$ is obtained as follows:
We have the universal sheaf
$\nbigl^u=\nbigl\otimes\nbigo_{\rel}(1)$
and the universal $L$-section 
$\phi:p_{M(c,L)}^{\ast}L\lrarr \nbigl^u$
over $M(c,L)\times X$.
We have the cokernel sheaf $\Cok:=\Cok(\phi)$.
Then the sheaves $p_{X\ast}(\Cok)$
and $R^1p_{X\ast}(\Cok)$ are locally free
$\nbigo_{M(c,L)}$-modules.
The vector bundle
$p_{X\,\ast}(\Cok)$ gives the tangent bundle of the moduli,
 and $R^1\pi_{X\,\ast}(\Cok)$
gives the obstruction bundle $O(c)$.
The virtual fundamental class is given by
the Euler class $Eu(O(c))\cap [M(c,L)]_0$,
where $[M(c,L)]_0$ is the ordinary fundamental class
of the smooth variety $M(c,L)$.

Let $d(c,L)$ denote the dimension of
the smooth variety $M(c,L)$,
which is same as
$\dim(H^0(X,L^{-1}\otimes\nbigl))-1$.
We put
$\chi(L^{-1}\otimes\nbigl):=
 \dim H^0(X,L^{-1}\otimes\nbigl)
-\dim H^1(X,L^{-1}\otimes\nbigl)
+\dim H^2(X,L^{-1}\otimes\nbigl)$.

\begin{prop}
\label{prop;06.5.13.70}
Assume $H^1(X,\nbigo_X)=0$
and $p_g=\dim H^2(X,\nbigo_X)> 0$.
Moreover,
we assume that the virtual fundamental class
of $M(c,L)$ is not $0$.
Let $\nbigl$ be  a line bundle such that $c_1(\nbigl)=c$.
Then, we have
$\chi(L^{-1}\otimes\nbigl)=\chi(\nbigo_X)$,
and the expected dimension of $M(c,L)$ is $0$.

The total Chern class of the obstruction bundle $O(c)$
is $\bigl(1+c_1(\nbigo_{\rel}(1))\bigr)^{d(c,L)-p_g}$,
and the virtual fundamental class is as follows:
\[
 \frac{\prod_{i=1}^{d(c,L)}(i-p_g)}{d(c,L)!}\cdot[p]
\]
Here $[p]$ denotes the cohomology class of a point $M(c,L)$.
We also have the inequalities $d(c,L)< p_g$
and $\dim H^1(X,L^{-1}\otimes\nbigl)<\dim H^2(X,L^{-1}\otimes\nbigl)$.
\end{prop}
\pf
Since we have assumed $H^1(X,\nbigo_X)=0$,
the obstruction bundle $O(c)$
can be decomposed on $M(c,L)$
as follows:
\[
 0\lrarr 
 H^1(X,L^{-1}\otimes\nbigl)\otimes\nbigo_{\rel}(1)
 \lrarr
 O(c)
 \lrarr
 H^2(X,\nbigo)\otimes\nbigo_{M(c,L)}
 \lrarr
 H^2(X,L^{-1}\otimes\nbigl)\otimes\nbigo_{\rel}(1)
 \lrarr 0.
\]
Thus the total Chern class of the vector bundle $O(c)$ is 
$\bigl(
 1+c_1(\nbigo_{\rel}(1))
 \bigr)^{\dim(H^1(X,L^{-1}\otimes\nbigl))-\dim(H^2(X,L^{-1}\otimes\nbigl))}$.
On the other hand, the rank of $O(c)$ is
$\dim\bigl(H^1(X,L^{-1}\otimes\nbigl)\bigr)
-\dim\bigl(H^2(X,L^{-1}\otimes\nbigl)\bigr)+p_g$.
Thus, if $\dim(H^1(X,L^{-1}\otimes\nbigl))-\dim(H^2(X,L^{-1}\otimes\nbigl))\geq 0$,
the Euler class of $O(c)$ is $0$ due to our assumption $p_g>0$.
It contradicts with our assumption
that the virtual fundamental class is not $0$.
Therefore,
we have $-e=\dim (H^1(X,L^{-1}\otimes\nbigl))-\dim(H^2(X,L^{-1}\otimes\nbigl))<0$.
Then, the coefficient of $c_1(\nbigo_{\rel}(1))^{d(c,L)}$
in the polynomial
$\bigl(1+c_1(\nbigo)_{\rel}\bigr)^{-e}$ is not $0$.
Thus, the rank of the $O\bigl(c\bigr)$ must be
$d(c,L)=\dim(H^0(X,L^{-1}\otimes\nbigl))-1$.
Therefore,
we obtain $\chi(\nbigo_X)=\chi(L^{-1}\otimes\nbigl)$,
and the expected dimension of the moduli is $0$.
By a formal calculation,
we obtain the virtual fundamental class.
\hfill\qed

\begin{rem}
In particular, we have the vanishing of
the virtual fundamental class
in the case $H^2(X,L^{-1}\otimes\nbigl)=0$.
We can derive it also from Proposition
{\rm\ref{prop;06.5.11.500}}.
\hfill\qed
\end{rem}

Let us consider the virtual fundamental class
of $M(\widehat{c},[L])$
when $H^2(X,L^{-1}\otimes\nbigl)=0$ holds
for any line bundle $\nbigl$ such that $c_1(\nbigl)=c$.
We have $M(\widehat{c},[L])\simeq M(c,L)\simeq 
\proj\bigl(H^0(X,L^{-1}\otimes\nbigl)^{\lor}\bigr)$
for such a line bundle $\nbigl$.
The obstruction theory of
$M(\widehat{c},[L])$ is perfect
(Proposition \ref{prop;06.5.11.550}).
Since the moduli is smooth,
we have the obstruction bundle
$\Ohat(c)$,
and the virtual fundamental class
is the Euler class of $\Ohat(c)$.
The following proposition can be checked
directly from the construction
of the obstruction theory.
\begin{prop}
We assume $H^1(X,\nbigo)=0$
and 
$H^2(X,L^{-1}\otimes\nbigl)=0$
for any line bundle $\nbigl$ such that
$c_1(\nbigl)=c$.
Let $\nbigl$ be a line bundle such that $c_1(\nbigl)=c$.
The expected dimension of $M(\widehat{c},[L])$ is
$\chi(L^{-1}\otimes\nbigl)$.
The obstruction bundle $\Ohat(c)$
is given by $H^1(X,\nbigl)\otimes\nbigo_{\rel}(1)$.
The total Chern class of $\Ohat(c)$
is $\bigl(1+c_1(\nbigo_{\rel}(1))\bigr)^{\dim H^1(X,\nbigl)}$.
In particular,
the virtual fundamental class is given by
$c_1(\nbigo_{\rel}(1))^{\dim H^1(X,\nbigl)}
\cap [M(\widehat{c},[L])]_0$,
where $[M(\widehat{c},[L])]$ denotes
the ordinary fundamental class.
\hfill\qed
\end{prop}

\subsubsection{Parabolic Hilbert scheme}
\label{subsubsection;06.7.3.431}

Let $\vecy$ be an element of $\Type$
such that $\rank(\vecy)=1$.
We assume that the first Chern class of $y$
is trivial.
Then, let $X^{[\vecy]}$ denote the moduli
of the oriented  parabolic sheaves $E_{\ast}$ 
of rank one such that $\det(E)=\nbigo_X$.
In other words,
$X^{[\vecy]}$ denotes the moduli space
of ideal sheaves of $0$-dimensional schemes
with parabolic structure of an appropriate type.
We call $X^{[\vecy]}$ the parabolic Hilbert scheme
of type $\vecy$.
The universal sheaf $\nbigi^u$ over
$X^{[\vecy]}\times X$ can be regarded
as the ideal sheaf of the relatively
$0$-dimensional scheme $\nbigz(\vecy)$.
The relative length is given by
$-y_2$, where $y_2$ denotes the $H^4(X)$-component
of $y$.

\begin{prop}
\label{prop;06.6.13.80}
When $D$ is smooth,
the parabolic Hilbert scheme is smooth.
The expected dimension is same as 
the ordinary dimension.
\end{prop}
\pf
$X^{[\vecy]}$ is the fiber of the smooth morphism
$\nbigm(\vecyhat)\lrarr \Pic$.
Hence the claim follows from the first claim of
Proposition \ref{prop;06.5.11.550}.
\hfill\qed

\vspace{.1in}

Due to Proposition \ref{prop;06.6.13.80},
the obstruction theory of $X^{[\vecy]}$
is obvious.
But, we give another expression of the obstruction theory
of $X^{[\vecy]}$ for later use.

Let $V_{\cdot}=(V_{-1}\rarr V_0)$ be a locally free resolution
of the universal sheaf $E$ over $X^{[\vecy]}\times X$.
We take vector spaces $W_i$
such that $\rank W_i=\rank V_i$.
We put $\SGL(W_{\cdot}):=
 \bigl\{(g_{-1},g_0)\in \GL(W_{\cdot}),\big|\,
 \det(g_{-1})\det (g_0)=1\bigr\}$.
We put $\Ybar(W):=
 N(W_{-1\,X},W_{0\,X})_{\SGL(W_{\cdot})}$.
Then, we have the classifying map
$\Phibar(E):X^{[\vecy]}\times X\lrarr \Ybar(W_{\cdot})$.
Thus, we obtain the morphism
$\Phibar(E)^{\ast}L_{\Ybar(W_{\cdot})/X}
\lrarr L_{X^{[\vecy]}\times X/X}$.
We can show that 
$\Phibar(E)^{\ast}L_{\Ybar(W_{\cdot})/X}$
is expressed by
$\gminig^{\circ}(V_{\cdot})_{\leq 1}$.
Therefore, we obtain the morphism
$\gminig^{\circ}(V_{\cdot})
\lrarr L_{X^{[\vecy]}\times X/X}$,
and it induces 
$\Ob^{\circ}(V_{\cdot}) \lrarr L_{X^{[\vecy]}}$.
Since it is uniquely determined in the derived category
$D(X^{[\vecy]})$,
we denote it by $\Ob^{\circ}(y)$.

We put 
$\Ybar_D(W_{\cdot}):=N(W_{-1\,D},W_{0\,D})_{\SGL(W_{\cdot})}$.
We have the naturally defined right action of
$\prod_{i=2}^l\GL(W^{(i)})\times
 \SGL(W_{\cdot})$ on 
$\prod_{i=1}^l N\bigl(W^{(i+1)},W^{(i)}\bigr)$.
The quotient stack is denoted by
$\Ybar_D(W_{\cdot},W^{\ast})$.
We have the naturally defined morphism
$\pi:\Ybar_D(W_{\cdot},W^{\ast})\lrarr \Ybar_D(W_{\cdot})$.

We have the morphisms
$\Phibar(V_{\cdot|D}):
 X^{[\vecy]}\times D\lrarr \Ybar_D(W_{\cdot})$
and $\Phibar_D(V_{\cdot|D},V_D^{\ast}):
 X^{[\vecy]}\times D\lrarr \Ybar(W_{\cdot},W^{\ast})$.
We have $\pi\circ\Phibar_D(V_{\cdot|D},V^{\ast})=\Phibar(V_{\cdot|D})$.
Thus, we obtain the morphisms:
\[
\begin{CD}
 L_{X^{[\vecy]}\times D/D}
 @<<<
 \Phibar_D(V_{\cdot|D},V_D^{\ast})^{\ast}
 L_{\Ybar_D(W_{\cdot},W^{\ast})/D}
@<<<
 \Phibar(V_{\cdot|D})^{\ast}
 L_{\Ybar_D(W_{\cdot})/D}
\end{CD}
\]
It is easy to see that
$\Phibar_D(V_{\cdot|D},V_D^{\ast})^{\ast}
 L_{\Ybar_D(W_{\cdot},W^{\ast})/D}$
and 
$\Phibar(V_{\cdot|D})^{\ast}
 L_{\Ybar_D(W_{\cdot})/D}$
are expressed by
$\gminig_D^{\circ}(V_{\cdot},F_{\ast})_{\leq 1}$
and $\gminig^{\circ}(V_{\cdot|D})_{\leq 1}$.
Thus, we obtain the morphisms:
\[
\begin{CD}
 L_{X^{[\vecy]}\times D/D}
@<<<
 \gminig_D^{\circ}(V_{\cdot},F_{\ast})
@<<<
 \gminig^{\circ}(V_{\cdot|D})
\end{CD}
\]
We put 
$\Ob_{D}^{\circ}(\vecy)
:=Rp_{D\,\ast}\bigl(
\gminig_D^{\circ}(V_{\cdot},F_{\ast})\otimes\omega_D
\bigr)$.
We also put
$\Ob^{\circ}_D(y):=
 Rp_{D\,\ast}\gminig^{\circ}(V_{\cdot|D})\otimes\omega_D$.
They are independent of the choice of $V_{\cdot}$
in the derived category.
Then, we obtain the following commutative diagram:
\[
 \begin{CD}
 \Ob^{\circ}_D(y) @>{i_1}>>
 \Ob^{\circ}(y) \\
 @V{i_2}VV @VVV \\
 \Ob_D^{\circ}(\vecy)
 @>>> L_{X^{[\vecy]}}
 \end{CD}
\]
The cone of
$(i_1,-i_2):\Ob^{\circ}_D(y)\lrarr
 \Ob(y)\oplus \Ob_D^{\circ}(\vecy)$
is denoted by $\Ob^{\circ}(\vecy)$.
We obtain the morphism
$\ob^{\circ}(\vecy):\Ob^{\circ}(\vecy)
 \lrarr L_{X^{[\vecy]}}$.

\begin{lem}
The morphism $\ob^{\circ}(\vecy)$ is quasi-isomorphic.
\end{lem}
\pf
We have the inclusion
$X^{[\vecy]}\lrarr \nbigm(\vecyhat)$.
We take a locally free resolution of the universal sheaf
on $\nbigm(\vecyhat)\times X$.
The restriction to $X^{[\vecy]}\times X$
gives the locally free resolution of the universal sheaf
on $X^{[\vecy]}\times X$.
Then, we obtain the following naturally defined commutative diagram:
\[
 \begin{CD}
 X^{[\vecy]}\times X
 @>{j_1}>>
 \nbigm(\vecyhat)\times X
 @>{j_2}>> \Pic\times X \\
 @VVV @VVV @VVV \\
 \Ybar(W_{\cdot}) @>>>
 Y(W_{\cdot}) @>>>
 X_{G_m}
 \end{CD}
\]
Then, we can obtain the following commutative diagram
on $X^{[\vecy]}$:
\[
 \begin{CD}
 L_{X^{[\vecy]}} @<<<
 j_1^{\ast}L_{\nbigm(\vecyhat)} @<<<
 j_1^{\ast}j_2^{\ast}L_{\Pic}\\
 @AAA @AAA @AAA \\
 \Ob^{\circ}(y) @<<<
 \Ob(y) @<<<
 \Ob^d(y)
 \end{CD}
\]
By the construction of $\Ob(\yhat)$,
we can obtain the following:
\[
 \begin{CD}
 L_{X^{[\vecy]}} @<<<
 j_1^{\ast}L_{\nbigm(\vecyhat)} @<<<
 j_1^{\ast}j_2^{\ast}L_{\Pic}\\
 @AAA @AAA @AAA \\
 \Ob^{\circ}(y) @<<<
 \Ob(\yhat) @<<<
j_1^{\ast}j_2^{\ast}L_{\Pic}
 \end{CD}
\]
It is easy to observe that the both of 
the horizontal rows are distinguished.

We also have the following diagram:
\[
 \begin{CD}
 X^{[\vecy]}\times X @>>>
 \Ybar_D(W_{\cdot},W^{\ast})@>>>
 \Ybar_D(W_{\cdot})\\
 @VVV @VVV @VVV \\
 \nbigm(\vecyhat)  @>>>
 Y_D(W_{\cdot},W^{\ast}) @>>>
 Y_D(W_{\cdot})
 \end{CD}
\]
We remark 
$\Ob(\vecyhat)\simeq \Cone\bigl(\Ob_D(y)\lrarr \Ob_D(\vecy)\oplus
 \Ob(\yhat)\bigr)$.
Then, we can derive the following commutative diagram
by construction of the complexes:
\[
\begin{CD}
 L_{X^{[\vecy]}} @<<< 
 \Ob^{\circ}_D(\vecy) @<<< 
 \Ob^{\circ}_D(y) \\
 @AAA @AAA @AAA \\
 j_1^{\ast}L_{\nbigm(\vecyhat)}
 @<<<
 \Ob_D(\vecy)
 @<<<
 \Ob_D(y)
\end{CD}
\]
Here, $\Ob_D(\vecy):=Rp_{D\,\ast}\gminig_D(V,F_{\ast})\otimes\omega_D$
and $\Ob_D(y):=Rp_{D\,\ast}\gminig(V_{\cdot|D})\otimes\omega_D$
on $\nbigm(\vecyhat)$.
Then, we obtain the following morphism of the distinguished
triangles:
\[
 \begin{CD}
 L_{X^{[\vecy]}} @<<<
 j_1^{\ast}L_{\nbigm(\vecyhat)} @<<<
 j_1^{\ast}j_2^{\ast}L_{\Pic}\\
 @AAA @A{a_1}AA @A{a_2}AA \\
 \Ob^{\circ}(\vecy) @<<<
 \Ob(\vecyhat) @<<<
j_1^{\ast}j_2^{\ast}L_{\Pic}
 \end{CD}
\]
The morphisms $a_i$ $(i=1,2)$ are isomorphic.
Therefore, the claim of the lemma follows.
\hfill\qed

\subsubsection{Splitting into the moduli of the abelian pairs
and the parabolic Hilbert schemes}
\label{subsubsection;06.5.13.40}

Let $\vecy$ be an element of $\Type$
such that $\rank(\vecy)=1$.
We have the following description of
$\nbigm\bigl(\vecy,L\bigr)$.
Let $c$ denote the $H^2(X)$-part of $y$,
and $\vecy(-c):=\vecy\cdot \exp(-c)$.
We have the universal line bundle $\nbigl_c$
on $M(c,L)\times X$.
We also have the structure sheaf
$\nbigo_{\nbigz(\vecy(-c))}$
of the subscheme
$\nbigz(\vecy(-c))\subset X^{[\vecy(-c)]}\times X$.
The pull back of them
via the projection
$M(c,L)\times X^{[\vecy(-c)]}\times X$
onto $M(c,L)\times X$ and $X^{[\vecy(-c)]}\times X$
are denoted by the same notation.
We put as follows:
\[
 \nbigk:=\nbigl_c^u\otimes 
 L^{-1}\otimes \nbigo_{\nbigz(\vecy(-c))}.
\]
It is easy to see that
$\gbigv:=p_{X\,\ast}\nbigk$ is a locally free sheaf
on $M(c,L)\times X^{[\vecy(-c)]}$.
We have the natural morphism
$p_{X\,\ast}\nbigl^u_c\otimes L^{-1}
\lrarr \gbigv$.
We also have the natural section
of $p_{X\,\ast}(L^{-1} \otimes\nbigl^u_c)$
induced by the universal section $\phi^u$
for $M(c,L)$.
Therefore,
we obtain the section 
$\psi:M(c,L)\times X^{[\vecy(-c)]}\lrarr \gbigv$.
It is easy to observe that
$\psi^{-1}(0)$ is isomorphic to the moduli $\nbigm(\vecy,L)$.

We have the following Cartesian diagram:
\begin{equation}
\label{eq;06.5.25.1}
 \begin{CD}
 \nbigm(\vecy,L) @>{F}>> M(c,L)\times X^{[\vecy(-c)]} \\
  @V{j}VV @VV{\psi}V \\
 M(c,L)\times X^{[\vecy(-c)]} @>{i}>>
 \gbigv
 \end{CD}
\end{equation}
Here $i$ denotes the $0$-section.
Since $i$ is a regular embedding,
we can define the Gysin map $i^{!}$.

\begin{prop}
\label{prop;06.5.13.35}
\mbox{{}}
\begin{itemize}
\item
We have the relation
$i^!\bigl([M(c,L)]\times [X^{[\vecy(-c)]}]\bigr)
=[\nbigm(\vecy,L)]$
among the virtual fundamental classes.
In particular,
we have
$i_{\ast}[\nbigm(\vecy,L)]
=\Eu(\gbigv)\cap \bigl([M(c,L)]\times [X^{[\vecy(-c)]}]\bigr)$,
where $\Eu(\gbigv)$ denotes the Euler class of $\gbigv$.

\item
Assume $H^2\bigl(X,\nbigl\bigr)=0$
for any $\nbigl\in \Pic(c)$ .
Then, we have the relation
$i^!\bigl( [M(\widehat{c},L)]\times [X^{[\vecy(-c)]}]
\bigr)
=\bigl[\nbigm(\vecyhat,[L])\bigr]$
and
$i_{\ast}[\nbigm(\vecyhat,[L])]
=\Eu(\gbigv)\cap \bigl([M(\widehat{c},[L])]\times [X^{[\vecy(-c)]}]\bigr)$.
\end{itemize}
\end{prop}
The proof will be given in the next subsubsections
\ref{subsubsection;06.5.13.61}--\ref{subsubsection;06.5.13.60}.

\vspace{.1in}
Before going into the proof,
we remark that 
the study of $\bigl[\nbigm(\vecy,L)\bigr]$
is reduced to the study of $X^{[\vecy(-c)]}$
by Proposition \ref{prop;06.5.13.70}
and Proposition \ref{prop;06.5.13.35},
when we assume $H^1(X,\nbigo)=0$ and $p_g>0$.
We are mainly interested in the cap product
of some cohomology class $\Phi$ 
and the fundamental class $\bigl[\nbigm(\vecy,L)\bigr]$.
In the interesting cases,
the cohomology class $\Phi$ is defined on
$M(c,L)\times X^{[\vecy(-c)]}$,
and thus we have only to consider the product of
$\Phi$ and $i_{\ast}\bigl[\nbigm(\vecy,L)\bigr]$.
Let $\nbigl$ be any line bundle on $X$
such that $c_1(\nbigl)=c$.
We put 
$\nbigktilde:=
 \nbigo_{\nbigz(\vecyhat(-c))}\otimes \nbigl\otimes L^{-1}$.
We put $\widetilde{\gbigv}:=p_{X\,\ast}\nbigktilde$
which is the vector bundle on $X^{[\vecy(-c)]}$.

\begin{prop}
Assume $H^1(X,\nbigo)=0$ and $p_g>0$.
We have the following formula in the cohomology ring
$H^{\ast}\bigl(M(c,L)\times X^{[\vecy(-c)]}\bigr)
\simeq
 H^{\ast}\bigl(M(c,L)\bigr)
\otimes H^{\ast}(X^{[\vecy(-c)]})$:
\begin{equation}
 \label{eq;06.5.13.30}
 i_{\ast}\bigl[\nbigm(\vecy,L)\bigr]
=\frac{\prod_{i=1}^{d(c,L)}(i-p_g)}{d(c,L)!}\cdot
 [p]\times
\bigl(
 \Eu(\widetilde{\gbigv})\cap[X^{[\vecy(-c)]}]
\bigr)
\end{equation}
Here $[p]$ denotes the cohomology class
of a point of $M(c,L)$.
\end{prop}
\pf
When $[\nbigm(\vecy,L)]\neq 0$,
the expected dimension of $M(c,L)$ is $0$,
and $[M(c,L)]$ is same as
$\bigl(d(c,L)!\bigr)^{-1}\cdot
 \prod_{i=1}^{d(c,L)}(i-p_g)$
points,
due to Proposition \ref{prop;06.5.13.70}.
Therefore, we have only to consider 
$\Eu(\gbigv)\cap\bigl([p]\times [X^{[\vecy(-c)]}]\bigr)$.
Moreover,
we can replace $\gbigv$
with $\gbigv_{|\{p\}\times X^{[\vecy(-c)]}}
=\gbigvtilde_{|\{p\}\times X^{[\vecy(-c)]}}$.
Then we obtain the formula (\ref{eq;06.5.13.30}).
\hfill\qed

\subsubsection{The morphism to the moduli of
 abelian pairs and the obstruction theory}
\label{subsubsection;06.5.13.61}

We use the notation in the subsubsection
\ref{subsubsection;06.5.4.1}.
Recall the diagram (\ref{eq;06.5.13.1}).
We would like to obtain a similar diagram
from an $L$-Bradlow pair
in the case $\rank(E)=1$.
Since $E$ is contained in $\det(E)$,
we have the naturally induced section
$\phi':L\lrarr \det(E)$.
Assume 
that $U$ is connected, for simplicity.
Hence $c=c_1\bigl(\det(E)_{|\{u\}\times X}\bigr)$
is independent of the choice of $u\in U$.
Therefore, we obtain the morphism
$\det_{E,\phi}:U_2\lrarr M(c,L)$.
\begin{prop}
 \label{prop;06.5.13.2}
We have the following commutative diagram:
\begin{equation}
\label{eq;06.5.13.2}
 \begin{CD}
 L_{U_2} @<<< \Ob(V_{\cdot},\phitilde)\\
 @AAA @AAA \\
 \det_{E,\phi}^{\ast}L_{M(c,L)} @<<<
 \det_{E,\phi}^{\ast}\Ob\bigl(M(c,L)\bigr)
 \end{CD}
\end{equation}
Here, we put
$\Ob(V_{\cdot},\phitilde):=
 Rp_{X\,\ast}\bigl(\gminig(V_{\cdot},\phitilde)\otimes\omega_X\bigr)$,
and $\Ob\bigl(M(c,L)\bigr)$ denotes the obstruction theory of
 $M(c,L)$.
\end{prop}
\pf
Let us begin with a general non-sense.
Let $\gbigv_0$ and $\gbigv_{-1}$ be locally free sheaves
on a stack $Z$
with a morphism $\gminif:\gbigv_{-1}\lrarr \gbigv_0$.
Assume $\rank \gbigv_0-\rank \gbigv_{-1}=1$.
Let $v_1,\ldots,v_l$ be a local frame of $\gbigv_{-1}$,
and let $u$ be a local section of $\gbigv_0$.
We put $\omega:=v_1\wedge\cdots\wedge v_l$.
Then we put 
$\Lambda_{\gminif}(u):=\bigl(
u\wedge \gminif(\omega)\bigr)\otimes \omega^{-1}$.
It is easy to see $\Lambda_{\gminif}(u)$ is independent of
a choice of the frame.
Therefore, we have the morphism 
$\Lambda_{\gminif}:\gbigv_0\lrarr \det(\gbigv_{\cdot})$.
It is easy to see that $\Lambda_{\gminif}\circ f=0$.
Hence,
we have the morphism of the complexes
$\cone(f)\lrarr \det(\gbigv_{\cdot})$.

Applying the above construction,
we obtain the morphism
$\Lambda_f:V_{\cdot}\lrarr \det(V_{\cdot})$.
We also obtain the morphism
$\Lambda_f\circ\phitilde_0:P_0\lrarr\det(V_{\cdot})\simeq\det(E)$.
Since we have
$\Lambda_f\circ \phitilde_0(P_{-1})=0$,
the morphisms
$\phi'$ and $\Lambda_f\circ\phitilde_0$ are essentially same.
We put
$\gminig_{\rel}:=\nhom\bigl(P_{\cdot},\det(V_{\cdot})\bigr)^{\lor}$.
The morphism $\Lambda_f\circ\phitilde_0$
naturally induces the morphism
$\gminig_{\rel}[-1]\lrarr \nbigo[-1]$.
The cone is denoted by
$\gminig\bigl(\det(V_{\cdot}),\phi'\bigr)$.

It is convenient to make a minor change in the construction
of $\Ob(V_{\cdot},\phitilde)$.
We have the natural right $\GL(W_{\cdot})$-action
on $N(W_{-1\,X},W_{0\,X})\times N(P_{-1},W_{0\,X})$.
The quotient stack is denoted by
$\Ytilde_0(W_{\cdot},P_{\cdot})$.
We put $\Ytilde_1(W_{\cdot},P_{\cdot}):=Y(W_{\cdot})$
and $\Ytilde_2(W_{\cdot},P_{\cdot})=Y_2(W_{\cdot},P_{\cdot})$.
We put
$\Ytilde(W_{\cdot},P_{\cdot}):=
 \Ytilde_1(W_{\cdot},P_{\cdot})\times_{\Ytilde_0(W_{\cdot},P_{\cdot})}
 \Ytilde_2(W_{\cdot},P_{\cdot})$.
We have the natural morphisms
$\Ytilde_i(W_{\cdot},P_{\cdot})\lrarr Y_i(W_{\cdot},P_{\cdot})$.
The induced morphism 
$\Ytilde(W_{\cdot},P_{\cdot})\lrarr 
 Y(W_{\cdot},P_{\cdot})$ is isomorphic.

We have the classifying map
$\Phitilde:U_2\times X\lrarr \Ytilde(W_{\cdot},P_{\cdot})$
and the induced maps
$\Phitilde_i:
 U_2\times X\lrarr \Ytilde_i(W_{\cdot},P_{\cdot})$.
In the construction of the subsubsection
\ref{subsubsection;06.5.4.1},
we can replace 
$\Phi(V_{\cdot},\phitilde)^{\ast}L_{Y(W_{\cdot},P_{\cdot})/X}$
with $\Phitilde^{\ast}L_{\Ytilde(W_{\cdot},P_{\cdot})/X}$.

We have the weight $1$-action of $G_m$ on $\det(W_{\cdot})$.
It induces the $G_m$-action on
$N\bigl(P_i,\det(W_{\cdot\,X})\bigr)$ $(i=0,1)$.
We put as follows:
\[
 Z_0(W_{\cdot},P_{\cdot}):=N\bigl(P_{-1},\det(W_{\cdot\,X})\bigr)_{G_m},
\quad
 Z_1(W_{\cdot},P_{\cdot}):=X_{G_m},
\quad
 Z_2(W_{\cdot},P_{\cdot}):=N\bigl(P_0,\det(W_{\cdot\,X})\bigr)_{G_m}
\]
The fiber product
$Z_1(W_{\cdot},P_{\cdot})\times_{Z_0(W_{\cdot},P_{\cdot})}
 Z_2(W_{\cdot},P_{\cdot})$ is 
isomorphic to $N(\nbigo_X,\det(W_{\cdot\,X}))_{G_m}$,
and it is denoted by $Z(W_{\cdot},P_{\cdot})$.
From the section
$\Lambda_f\circ\phitilde_0:P_0\lrarr \det(V_{\cdot})$,
we obtain the classifying map
$\Psi:U_2\times X\lrarr Z(W_{\cdot},P_{\cdot})$
and the induced maps
$\Psi_i:U_2\times X\lrarr Z_i(W_{\cdot},P_{\cdot})$
$(i=0,1,2)$.

We have the morphism
$\Gamma_2':
 N(W_{-1},W_0)\times N(P_{-1},W_{-1})\times N(P_0,W_0)
\lrarr N(P_0,\det(W_{\cdot\,X}))$
given by
$\Gamma_2'(e,a_{-1},a_0):=\Lambda_e\circ a_0$.
We have the morphism
$\GL(W_{\cdot})\lrarr G_m$
given by $(g_{-1},g_0)\longmapsto \det(g_{-1})^{-1}\cdot \det(g_0)$.
Then $\Gamma_2'$ is equivariant
with respect to the actions of $\GL(W_{\cdot})$ and $G_m$.
Therefore, we obtain the morphism
$\Gamma_2:\Ytilde_2(W_{\cdot},P_{\cdot})\lrarr
 Z_2(W_{\cdot},P_{\cdot})$.
Similarly, we have the equivariant morphism
$\Gamma_0':N(P_{-1},W_0)\lrarr 
 N\bigl(P_{-1},\det(W_{\cdot\,X})\bigr)$
given by $\Gamma_0'(e,a)=\Lambda_e\circ a$,
which induces
$\Gamma_0:\Ytilde_0(W_{\cdot},P_{\cdot})
\lrarr Z_0(W_{\cdot},P_{\cdot})$.
We also have the obvious map
$\Gamma_1:\Ytilde_1(W_{\cdot},P_{\cdot})
\lrarr Z_1(W_{\cdot},P_{\cdot})$.
We have $\Gamma_i\circ\Phitilde_i=\Psi_i$
for $i=0,1,2$.
We have the induced map
$\Gamma:\Ytilde(W_{\cdot},P_{\cdot})\lrarr Z(W_{\cdot},P_{\cdot})$
and the relation $\Gamma\circ \Phitilde=\Psi$.

We have the universal $L$-abelian pair
$\bigl(\nbigl^u,\phi^u\bigr)$ over $M(c,L)\times X$.
We have the classifying map
$\Phi_i(\nbigl^u,\phi^u):
M(c,L)\times X\lrarr Z_i(W_{\cdot},P_{\cdot})$.
We obtain the following commutative diagram:
\begin{equation}
\label{eq;06.5.23.100}
 \begin{CD}
 U_2\times X 
@>{\Phitilde_i}>> 
 \Ytilde_i(W_{\cdot},P_{\cdot})\\
 @V{\det_{E,\phi,X}}VV @V{\Gamma_i}VV \\
 M(c,L)\times X
@>{\Phi_i(\nbigl^u,\phi^u)}>> 
 Z_i(W_{\cdot},P_{\cdot})
 \end{CD}
\end{equation}
We obtain the following morphism on $U_2\times X$:
\[
\varphi:
 \cone\Bigl(
 \Psi_0^{\ast}L_{Z_0(W_{\cdot},P_{\cdot})/X}
\rarr
\bigoplus_{i=1,2}
\Psi_i^{\ast}L_{Z_i(W_{\cdot},P_{\cdot})/X}
 \Bigr)
\lrarr
 \cone\Bigl(
 \Phitilde_0^{\ast}L_{\Ytilde_0(W_{\cdot},P_{\cdot})/X}
\rarr\bigoplus_{i=1,2}
 \Phitilde_i^{\ast}L_{\Ytilde_i(W_{\cdot},P_{\cdot})/X}
\Bigr)
\]
By the argument in the subsubsection
\ref{subsubsection;06.4.29.15},
we can show that
$\varphi$ is expressed by the morphism of the complexes
$\gminig\bigl(\det(V_{\cdot}),\phi'\bigr)
\lrarr \cone\bigl(\gamma(\phitilde)_{\leq 1}\bigr)$.

\begin{lem}
\label{lem;06.5.23.101}
The morphism $\varphi$ naturally factors through
$\gminig(\det(V_{\cdot}),\phi')
\lrarr \gminig(V_{\cdot},\phitilde)
\lrarr \cone\bigl(\gamma(\phitilde)_{\leq 1}\bigr)$.
\end{lem}
\pf
We give only an indication.
The following diagram is commutative:
\[
\begin{CD}
 \nbigo @>>>
 \nhom(V_0,V_0)\oplus\nhom(V_{-1},V_{-1})\\
 @VVV @VVV \\
 0 @>>>\nhom(V_{-1},V_{0})
\end{CD} 
\]
In fact, it is a part of the morphism of the complexes
$\nbigo\lrarr \nhom(V_{\cdot},V_{\cdot})^{\lor}$.
The following diagram is commutative:
\[
\begin{CD}
 \nhom\bigl(\det(V_{\cdot}),P_0\bigr)
 @>>>\nhom(V_{-1},P_{-1})\oplus\nhom(V_{0},P_0)\\
 @VVV @VVV\\
 0 @>>> \nhom(V_{-1},P_0)
\end{CD}
\]
In fact, it is a part of the morphism of the complexes
$\nhom\bigl(P_{\cdot},\det(V_{\cdot})\bigr)^{\lor}
\lrarr
 \nhom\bigl(P_{\cdot},V_{\cdot}\bigr)^{\lor}$.
Then the claim of the lemma can be checked easily.
\hfill\qed

\vspace{.1in}

Let us finish the proof of Proposition \ref{prop;06.5.13.2}.
We obtain the following commutative diagram
from (\ref{eq;06.5.23.100}) and Lemma \ref{lem;06.5.23.101}:
\[
 \begin{CD}
 \gminig(V_{\cdot},\phitilde)
 @<<<
 \gminig\bigl(\det(V_{\cdot}),\phi'\bigr) \\
 @VVV @VVV \\
 L_{U_2\times X/X}
 @<<<
 \det_{E,\phi,X}^{\ast}L_{M(c,L)\times X/X}
 \end{CD}
\]
It is easy to observe
$\Ob\bigl(M(c,L)\bigr)=
 Rp_{X\,\ast}\bigl(
\gminig\bigl(\det(V_{\cdot}),\phi'\bigr)\otimes\omega_X\bigr)$.
Thus, we obtain the desired commutative diagram
(\ref{eq;06.5.13.2}).
\hfill\qed

\vspace{.1in}
By a similar argument,
we obtain a similar commutative diagram
in the reduced case.
We use the notation in the subsubsection
{\rm\ref{subsubsection;06.5.4.110}}.
Assume $\rank(E)=1$.
Then the reduced $L$-pair
$\bigl(\det(E),[\phi']\bigr)$ is induced.
Therefore, we have the morphism
$\det_{E,[\phi]}:U_3\lrarr M(c,[L])$.

\begin{lem}
We have the following commutative diagram:
\begin{equation}
 \label{eq;06.5.13.6}
 \begin{CD}
 L_{U_3}@<<< \Ob(V_{\cdot},[\phitilde])\\
 @AAA @AAA \\
 \det_{E,[\phi]}^{\ast}L_{M(c,[L])} @<<<
 \Ob\bigl(M(c,[L])\bigr)
 \end{CD}
\end{equation}
Here, we put $\Ob(V_{\cdot},[\phitilde]):=
 \Cone\bigl(\Ob_{\rel}(V_{\cdot},[\phitilde])
 \lrarr \Ob(V_{\cdot})
 \bigr)$,
and $\Ob\bigl(M(c,[L])\bigr)$ denotes
the obstruction theory of $M(c,[L])$.
\end{lem}
\pf
We indicate only an outline.
From the weight $(-1)$-actions of $G_m$
on $P_{\cdot}$ and $\nbigo_X$,
we obtain the $G_m$-actions
on $\Ytilde_i(W_{\cdot},P_{\cdot})$
and $Z_i(W_{\cdot},P_{\cdot})$.
The quotient stacks are denoted by
$\Ytilde_i(W_{\cdot},[P_{\cdot}])$
and $Z_i(W_{\cdot},[P_{\cdot}])$.
The natural morphisms
$\Gamma_i:\Ytilde_i(W_{\cdot},[P_{\cdot}])\lrarr 
 Z_i(W_{\cdot},[P_{\cdot}])$ are induced.

We have the induced map
$\Phitilde'_i:U_3\times X\lrarr \Ytilde_i(W_{\cdot},[P_{\cdot}])$.
We put $\Psi_i':=\Gamma_i\circ \Phitilde_i$.
We also have the classifying maps
$\Phi_i(\nbigl^u,[\phi^u]):M(c,[L])\times X
\lrarr Z_i(W_{\cdot},[P_{\cdot}])$
for the universal objects on $M(c,[L])\times X$.
We have the following commutative diagram:
\[
 \begin{CD}
 U_3\times X @>{\Phitilde'_i}>>
 Y_i(W_{\cdot},[P_{\cdot}]) @>>>
 X_{G_m}\\
 @V{\det_{E,[\phi],X}}VV @V{\Gamma_i}VV @V=VV \\
 M(c,[L])\times X @>{\Phi(\nbigl^u,[\phi^u])}>> 
 Z_i(W_{\cdot},[P_{\cdot}])  @>>>
 X_{G_m}
 \end{CD}
\]
We have the following induced morphism on $U_3\times X$:
\begin{multline}
 \varphi':
 \cone\Bigl(
 \Psi_0^{\prime\,\ast}
 L_{Z_0(W_{\cdot},[P_{\cdot}])/X_{G_m}}
\rarr
 \bigoplus_{i=1,2}
 \Psi_i^{\prime\,\ast}
 L_{Z_i(W_{\cdot},[P_{\cdot}])/X_{G_m}}
 \Bigr)
\lrarr \\
 \cone\Bigl(
 \Phitilde_0^{\prime\,\ast}
 L_{Y_0(W_{\cdot},[P_{\cdot}])/X_{G_m}}
\rarr
 \bigoplus_{i=1,2}
 \Phitilde_i^{\prime\,\ast}
 L_{Y_i(W_{\cdot},[P_{\cdot}])/X_{G_m}}
 \Bigr)
\end{multline}
We can show that $\varphi'$ is expressed by 
the morphism of the complexes
$\gminig(\det(V_{\cdot}),[\phi'])\lrarr
 \Cone(\gamma[\phitilde]_{\leq 0}[-1])$,
where $\gminig\bigl(\det(V_{\cdot}),[\phi']\bigr)$
is given by $\nhom\bigl(\det(V_{\cdot}),[P_{-1}]\bigr)
\lrarr \nhom\bigl(\det(V_{\cdot}),[P_0]\bigr)
\lrarr\nbigo$,
as in the proof of Proposition \ref{prop;06.5.13.2}.
It factors through
$\gminig'(V_{\cdot},[\phitilde])$.
Therefore, we obtain the following diagram:
\[
 \begin{CD}
 \gminig'(V_{\cdot},[\phitilde]) @<<<
 \gminig\bigl(\det(V_{\cdot}),[\phi']\bigr)\\
 @VVV @VVV \\
 L_{U_3\times X/X_{G_m}} @<<<
 \det_{E,[\phi],X}^{\ast}L_{M(c,[L])\times X/X_{G_m}}
 \end{CD}
\]
Then, we obtain the following:
\[
 \begin{CD}
 Rp_{X\,\ast}\bigl(
 \gminig'(V_{\cdot},[\phitilde])\otimes\omega_X\bigr)
@>{\varphi_1}>>
 L_{U_3/k_{G_m}} 
@>{\varphi_2}>>
 L_{k_{G_m}/k}\\
@AAA @AAA @AAA\\
 Rp_{X\,\ast}\bigl(
 \gminig\bigl(\det(V_{\cdot}),[\phi']\bigr)\otimes\omega_X
 \bigr)
@>{\psi_1}>>
 \det_{E,[\phi]}^{\ast}L_{M(c,[L])/k_{G_m}}
@>{\psi_2}>>
 L_{k_{G_m}/k}[1]
 \end{CD}
\]
It is easy to observe 
$\Ob(V_{\cdot},[\phitilde])\simeq
 \cone\bigl( \varphi_2\circ\varphi_1 \bigr)$
and
$\det_{E,[\phi]}^{\ast}\Ob\bigl(M(c,[L])\bigr) \simeq
 \cone\bigl(\psi_2\circ\psi_1 \bigr)$.
Thus, we obtain the desired diagram (\ref{eq;06.5.13.6}).
\hfill\qed

\vspace{.1in}
When $E$ has an orientation $\rho$,
we have the morphisms
$\Ob_{\rel}(V_{\cdot},\rho)[-1]\lrarr \Ob(V_{\cdot})
\lrarr \Ob(V_{\cdot},[\phitilde])$.
Let $\Ob(V_{\cdot},\rho,[\phitilde])$ denote
the cone of the composite of the morphisms.
On the other hand,
we have the morphism
$\det_{E,[\phi],\rho}:U_3\lrarr M(\widehat{c},[L])$.

\begin{prop}
\label{prop;06.5.23.105}
We have the following commutative diagram:
\[
 \begin{CD}
 L_{U_3} @<<< \Ob(V_{\cdot},\rho,[\phitilde])\\
 @AAA @AAA \\
 \det_{E,[\phi],\rho}^{\ast}L_{M(\widehat{c},[L])} @<<<
 \det_{E,[\phi],\rho}^{\ast}\Ob\bigl(M(\widehat{c},[L])\bigr)
 \end{CD}
\]
Here, $\Ob\bigl(M(\widehat{c},[L])\bigr)$
denotes the obstruction theory of $M\bigl(\widehat{c},[L]\bigr)$.
\end{prop}
\pf
The orientation induces the morphisms
$U_3\lrarr M(\widehat{c},[L])\lrarr \Pic(c)$,
and we have the following diagram:
\[
 \begin{CD}
 \Ob^d(V_{\cdot}) @>>>
 \Ob^{d}(V_{\cdot}) @>>>
 \Ob(V_{\cdot})\\
 @VVV @VVV @VVV \\
 \det_{E,\rho}^{\ast}L_{\Pic} @>>> 
 \det_{E,[\phi],\rho}^{\ast}L_{M(\widehat{c},[L])}@>>> 
 L_{U_3}
 \end{CD}
\]
Then the claim follows from the construction
of the relative obstruction complex for orientations.
\hfill\qed

\subsubsection{The morphism to the parabolic Hilbert scheme
 and the obstruction theories}

Let $U$ be a $k$-scheme.
Let $(E,F_{\ast})$ be a quasi-parabolic 
torsion free sheaf of type $\vecy$.
Assume $\rank(\vecy)=1$.
Let $c$ be the $H^2(X)$-part of $y$.
Then we put $\vecy(-c):=\vecy\cdot \exp(-c)$.
We put $I(E):=\det(E)^{-1}\otimes E$,
which has the induced quasi-parabolic structure $F_{\ast}$.
The type of $\bigl(I(E),F_{\ast}\bigr)$ is $\vecy(-c)$.
Thus, we have the morphism
$\Xi(E):U\lrarr X^{[\vecy(-c)]}$.

\begin{lem}
\label{lem;06.5.13.11}
We have the following commutative diagram:
\begin{equation}
 \label{eq;06.5.13.7}
 \begin{CD}
 L_{U} @<<<  \Ob(V_{\cdot},F_{\ast}) \\
 @AAA @AAA \\
 \Xi(E)^{\ast}L_{X^{[\vecy(-c)]}} @<<< \Xi(E)^{\ast}\Ob(X^{[\vecy(-c)]})
 \end{CD}
\end{equation}
Here, we put
$\Ob(V_{\cdot},F_{\ast}):=
 \Cone\bigl(\Ob_{\rel}(V,F_{\ast})[-1]\rarr \Ob(V_{\cdot})\bigr)$,
and
$\Ob(X^{[\vecy(-c)]})$ denotes
the obstruction theory of $X^{[\vecy(-c)]}$.
Moreover,
$\varphi$ factors through $\Ob^{\circ}(V_{\cdot},F_{\ast})$.
\end{lem}
\pf
Let $\nbigi^u$ denote the universal sheaf
on $X^{[\vecy(-c)]}\times X$.
We take a locally free resolution
$\nbigv_{\cdot}$ of $\nbigi^u$.
It is easy to observe that
$V_{\cdot}:=\Xi(E)_X^{\ast}\nbigv_{\cdot}\otimes\det(E)$
is a locally free resolution of $E$.
We have
$\gminig(V_{\cdot})=
 \Xi(E)_X^{\ast}\gminig(\nbigv_{\cdot})$,
$\gminig(V_{\cdot|D})=
 \Xi(E)_D^{\ast}\gminig(\nbigv_{\cdot|D})$
and
$\gminig(V_{\cdot},F_{\ast})=
 \Xi(E)_D^{\ast}\gminig(\nbigv_{\cdot},F_{\ast})$.

We take vector spaces $W_i$ $(i=-1,0)$
such that $\rank W_i=\rank V_i$.
We have $\rank W_0-\rank W_{-1}=1$.
In that case,
we have the homomorphism
$\GL(W_{\cdot})\lrarr \SGL(W_{\cdot})$,
given as follows:
\[
 \bigl(g_{-1},g_{0}\bigr)
\longmapsto
 \bigl(
 \det(g_{\cdot})^{-1}\!\cdot\! g_{-1},\,\,
 \det(g_{\cdot})^{-1}\!\cdot\! g_0
\bigr)
\]
Here, $\det(g_{\cdot})$ denotes
$\det(g_0)\cdot \det(g_{-1})^{-1}$.
Therefore, we have the isomorphism
$Y(W_{\cdot})\simeq \Ybar(W_{\cdot})\times_X X_{G_m}$.
In particular, we have the morphism:
\begin{equation}
 \gminiw_1:Y(W_{\cdot})\lrarr \Ybar(W_{\cdot})
\end{equation}
Then, we obtain the following commutative diagram:
\[
 \begin{CD}
 U\times X @>>> Y(W_{\cdot}) \\
 @VVV @VVV \\
 X^{[\vecy(-c)]}\times X @>>> \Ybar(W_{\cdot})
 \end{CD}
\]
Therefore, we obtain the following commutative diagram
on $U\times X$:
\[
 \begin{CD}
 L_{U\times X/X}@<<<  
 \Phi(V_{\cdot})^{\ast}L_{Y(W_{\cdot})}
 @<<< \gminig(V_{\cdot})\\
 @AAA @AAA @AAA\\
 \Xi(E)_X^{\ast}L_{X^{[\vecy(-c)]}\times X/X} @<<<
 \Xi(E)_X^{\ast}\Phi(\nbigv_{\cdot})^{\ast}
 L_{\Ybar(W_{\cdot})}
 @<<< \Xi(E)^{\ast}\gminig^{\circ}(\nbigv_{\cdot})
 \end{CD}
\]
Hence we obtain the following commutative diagram:
\[
\begin{CD}
 L_{U}@<<< \Ob(V_{\cdot})\\
 @AAA @AAA \\
 \Xi(E)^{\ast}L_{X^{[\vecy(-c)]}} @<<<
 \Xi(E)^{\ast}\Ob^{\circ}(\nbigv_{\cdot})
\end{CD}
\]
We obtain a similar diagram from $V_{\cdot|D}$
and $\nbigv_{\cdot|D}$.
Moreover, we obtain the following diagram
by the argument of Lemma \ref{lem;06.5.13.10}:
\[
 \begin{CD}
 L_U @<<< \Ob(V_{\cdot}) @<<< \Ob(V_{\cdot|D})\\
 @AAA @AAA @AAA \\
 \Xi(E)^{\ast}L_{X^{[\vecy(-c)]}} @<<<
 \Xi(E)^{\ast}\Ob^{\circ}(\nbigv_{\cdot}) @<<<
 \Xi(E)^{\ast}\Ob^{\circ}(\nbigv_{\cdot|D})
 \end{CD}
\]
We also have the following commutative diagram:
\[
 \begin{CD}
 L_U @<<< 
 Rp_{D\,\ast}\bigl(
 \gminig_D(V_{\cdot},F_{\ast})\otimes\omega_D\bigr)
 @<<< 
 \Ob(V_{\cdot|D})\\
 @AAA @AAA @AAA \\
 \Xi(E)^{\ast}L_{X^{[\vecy(-c)]}}
 @<<< 
 \Xi(E)^{\ast} Rp_{D\,\ast}\bigl(
 \gminig_D^{\circ}(\nbigv_{\cdot},F_{\ast})\otimes\omega_D\bigr)
 @<<< \Xi(E)^{\ast}\Ob^{\circ}(\nbigv_{\cdot|D})
 \end{CD}
\]
We remark the following:
\[
 \Ob(V,F_{\ast})\simeq
 \Cone\Bigl(
 \Ob(V_{\cdot|D})\lrarr 
 \Ob(V_{\cdot})\oplus
 Rp_{D\,\ast}\bigl(
 \gminig_D(V_{\cdot},F_{\ast})\otimes\omega_D\bigr)
 \Bigr)
\]
\[
  \Ob(X^{[\vecy(-c)]})\simeq
 \Cone\Bigl(
 \Ob^{\circ}(\nbigv_{\cdot|D})\lrarr
 \Ob^{\circ}(\nbigv_{\cdot})\oplus
 Rp_{D\,\ast}\bigl(
 \gminig_D^{\circ}(\nbigv_{\cdot},F_{\ast})
 \otimes\omega_D
\bigr)
 \Bigr)
\]
Hence, we obtain the desired commutative diagram
(\ref{eq;06.5.13.7}).
\hfill\qed

\subsubsection{The mixed case}
\label{subsubsection;06.6.13.100}

Let $(E,F_{\ast},\phi)$ be a quasi parabolic $L$-Bradlow pair
of type $\vecy$ over $U\times X$.
Assume $\rank(\vecy)=1$.
Let $c$ denote the $H^2(X)$-part of $y$,
and we put $\vecy(-c):=\vecy\cdot \exp(-c)$.
Then, we have the morphisms
$\det_{E,\phi}:U\lrarr M(c,L)$ and
$\Xi(E):U\lrarr X^{[\vecy(-c)]}$.

Assume we have a locally free resolution $V_{\cdot}$ of $E$,
a locally free resolution $P_{\cdot}$ of $L$,
and a lift $\phitilde:P_{\cdot}\lrarr V_{\cdot}$ of $\phi$.
We have the natural morphism
$i_1:\Ob(V_{\cdot})\lrarr \Ob(V_{\cdot},\phitilde)$
and $i_2:\Ob(V_{\cdot})\lrarr \Ob(V_{\cdot},F_{\ast})$.
We put as follows:
\[
 \Ob(V_{\cdot},F_{\ast},\phi):=
 \Cone\Bigl(
 \Ob(V_{\cdot})\stackrel{(i_1,-i_2)}{\lrarr}
 \Ob(V_{\cdot},F_{\ast})
\oplus
 \Ob(V_{\cdot},\phi)
 \Bigr)
\]
We have the induced map
$\Ob(V_{\cdot},F_{\ast},\phi)\lrarr
 L_{U}$.

\begin{lem}
\label{lem;06.6.13.101}
We have the following commutative diagram:
\[
\begin{CD}
 \det_{E,\phi}^{\ast}\Ob\bigl(M(c,L)\bigr)
\oplus
 \Xi(E)^{\ast}\Ob(X^{[\vecy(-c)]})
@>>>
\Ob(V_{\cdot},F_{\ast},\phi)\\
@VVV @VVV \\
 \det_{E,\phi}^{\ast}L_{M(c,L)}
\oplus
 \Xi(E)^{\ast} L_{X^{[\vecy(-c)]}}
@>>>
 L_{U}
\end{CD}
\]
\end{lem}
\pf
It follows from Proposition \ref{prop;06.5.13.2}
and Lemma \ref{lem;06.5.13.11}.
\hfill\qed

\vspace{.1in}

Let $(E,F_{\ast},\rho,[\phi])$ be
an oriented quasi-parabolic reduced $L$-Bradlow pair
of type $\vecy$ over $U\times X$.
We have the morphism
$\det_{E,\rho,[\phi]}:U\lrarr M(\widehat{c},[L])$
and $\Xi(E):U\lrarr X^{[\vecy(-c)]}$.

Assume we have a locally free resolution $V_{\cdot}$ of $E$,
a locally free resolution $P_{\cdot}$ of $L$,
and a lift $[\phitilde]$ of $[\phi]$.
We have the natural morphisms
$i_1:\Ob(V_{\cdot})\lrarr \Ob(V_{\cdot},\rho,[\phitilde])$
and $i_2:\Ob(V_{\cdot})\lrarr \Ob(V_{\cdot},F_{\ast})$.
We put as follows:
\[
 \Ob(V_{\cdot},F_{\ast},[\phitilde],\rho):=
 \Cone\Bigl(
 \Ob(V_{\cdot})\stackrel{(i_1,-i_2)}{\lrarr}
 \Ob(V_{\cdot},[\phitilde],\rho)
\oplus
 \Ob(V_{\cdot},F_{\ast})
 \Bigr)
\]

\begin{lem}
\label{lem;06.5.13.50}
We have the following commutative diagram:
\[
 \begin{CD}
 \det_{E,\rho,[\phi]}^{\ast}\Ob\bigl(M(\widehat{c},[L])\bigr)
 \oplus
 \Xi(E)^{\ast}\Ob(X^{[\vecy(-c)]}) 
 @>>>
 \Ob(V_{\cdot},F_{\ast},[\phi],\rho)
 \\
 @VVV @VVV \\
 \det_{E,\rho,[\phi]}^{\ast}L_{M(\widehat{c},[L])}
\oplus
 \Xi(E)^{\ast}L_{X^{[\vecy(-c)]}}
@>>>
 L_U
 \end{CD}
\]
\end{lem}
\pf
It follows from Proposition \ref{prop;06.5.23.105}
and Lemma \ref{lem;06.5.13.11}.
\hfill\qed

\subsubsection{Proof of Proposition \ref{prop;06.5.13.35}}
\label{subsubsection;06.5.13.60}

Applying the construction in the subsubsection 
\ref{subsubsection;06.6.13.100}
to the universal object
$(\nbige^u,\phi)$ on $\nbigm(\vecy,L)$,
we obtain the morphism
$F:\nbigm(\vecy,L)\lrarr
 M(c,L)\times X^{[\vecy(-c)]}$.
It is same as the inclusion given in the subsubsection 
\ref{subsubsection;06.5.13.40}.
Let us denote the obstruction theory of 
$\nbigm(\vecy,L)$ by $\Ob(\vecy,L)$.
We have the following commutative diagram,
due to Lemma \ref{lem;06.6.13.101}:
\begin{equation}
 \label{eq;06.5.23.111}
 \begin{CD}
 F^{\ast}\Bigl(
 \Ob\bigl(M(c,L)\bigr)\oplus \Ob(X^{[\vecy(-c)]})
\Bigr)
 @>>>
 \Ob\bigl(\vecy,L\bigr)\\
 @VVV @VVV \\
 F^{\ast}\Bigl( L_{M(c,L)\times X^{[\vecy(-c)]}} \Bigr)
 @>>>
L_{\nbigm(\vecy,L)}
 \end{CD}
\end{equation}

Let $\nbigl_c$ denote the universal line bundle
on $M(c,L)\times X$.
We have the natural inclusion
$\nbige^u\otimes L^{-1}\lrarr
 F^{\ast}\nbigl_c^u\otimes L^{-1}$,
and the quotient is isomorphic to $\nbigk_{|\nbigm(\vecy,L)\times X}$.
Thus, we have the following:
\begin{multline}
 \label{eq;06.5.13.55}
 \cone\Bigl(
F^{\ast}\bigl[
 \Ob\bigl(X^{[\vecy(-c)]}\bigr)
\oplus
 \Ob\bigl(M(c,L)\bigr)\bigr]
\lrarr
 \Ob(\vecy,L)
 \Bigr) \\
\simeq
 Rp_{X\,\ast}\Bigl(
 \Cone\bigl[
  \nrhom\bigl(F^{\ast}\nbigl_c^u, L\bigr)
\rarr
 \nrhom\bigl(\nbige^u,L\bigr)
 \bigr]\otimes\omega_X
\Bigr) \\
\simeq
 Rp_{X\,\ast}\bigl(\nbigk\bigr)^{\lor}_{|\nbigm(\vecy,L)}[1]
\simeq
 \gbigv^{\lor}[1]_{|\nbigm(\vecy,L)}
\simeq
 j^{\ast}L_{M(c,L)\times X^{[\vecy(-c)]}/\gbigv}
\end{multline}

From the diagram (\ref{eq;06.5.23.111}),
we obtain the morphism
$\nu:
 j^{\ast}L_{M(c,L)\times X^{[\vecy(-c)]}/\gbigv}
\lrarr
 L_{\nbigm(\vecy,L)/M(c,L)}$

\begin{lem}
The morphism $\nu$ is same as 
the morphism
$j^{\ast}L_{M(c,L)\times X^{[\vecy(-c)]}/\gbigv}
\lrarr
 L_{\nbigm(\vecy,L)/M(c,L)}$
obtained from the diagram {\rm(\ref{eq;06.5.25.1})}.
In particular,
the obstruction theories of $M(c,L)\times X^{[\vecy(-c)]}$
and $\nbigm(\vecy,L)$
are compatible over the morphism $i$.
\end{lem}
\pf
We take a locally free resolution $V_{\cdot}$ of $\nbige^u$.
We take vector spaces $W_i$ such that $\rank W_i=\rank V_i$.
We obtain the following commutative diagram:
\[
 \begin{CD}
 \nbigm(\vecy,L)\times X @>{\Phi(V_{\cdot},\phi)}>>
 Y(W_{\cdot},P_{\cdot})\\
 @VVV @V{\gminiw_1\times \Gamma}VV \\
 X^{[\vecy(-c)]}\times M(c,L)
 @>{\Phi(\det(E)^{-1}\otimes V_{\cdot})\times \Phi(\nbigl,\phi)}>>
 \Ybar(W_{\cdot})\times Z(W_{\cdot},P_{\cdot})
 \end{CD}
\]
It induces the morphism:
\[
\begin{CD}
 L_{\nbigm(m,\vecy)\times X/M(c,L)\times X^{[\vecy(-c)]}\times X}
 @<<<
 \Phi(V_{\cdot},\phi)^{\ast}
 L_{Y(W_{\cdot},P_{\cdot})/
 \Ybar(W_{\cdot})\times Z(W_{\cdot},P_{\cdot})}.
\end{CD}
\]
We use the notation in the subsubsection
\ref{subsubsection;06.5.25.2}.
We put 
$\Vtilde_{i}=V_i$ $(i=0,-1)$
and $\Vtilde_1=\nbigl_c^u$.
By the argument in the subsubsection 
\ref{subsubsection;06.4.29.15},
we can show that
$\Phi(V_{\cdot},\phi)^{\ast}
 L_{Y(W_{\cdot},P_{\cdot})/
 \Ybar(W_{\cdot})\times Z(W_{\cdot},P_{\cdot})}$
is expressed by
$\gminik(V_{\cdot},P_{\cdot},\phi,\phitilde)_{\leq 0}$.
Thus we obtain the following morphism:
\[
\begin{CD}
 \gminik(\Vtilde_{\cdot},P_{\cdot},\phi,\phitilde)
@>>> L_{\nbigm(\vecy,L)\times X/M(c,L)\times X^{[\vecy(-c)]}\times X}
\end{CD}
\]
We put $\Ob^H(V_{\cdot},P_{\cdot},\phi,\phitilde):=
 Rp_{X\,\ast}\bigl(
\gminik(\Vtilde_{\cdot},P_{\cdot},\phi,\phitilde)
 \otimes\omega_X\bigr)$.
Then, 
$\Ob^H(V_{\cdot},P_{\cdot},\phi,\phitilde)$
is isomorphic to
the cone of the morphism
$\Ob(M(c,L))\oplus\Ob(X^{[\vecy(-c)]})
\lrarr \Ob(\nbigm(\vecy,L))$.
The induced morphism
$\Ob^H(V_{\cdot},P_{\cdot},\phi,\phitilde)
\lrarr L_{\nbigm(\vecy,L)/M(c,L)\times X^{[\vecy(-c)]}}$
is same as $\nu$.

We have the following factorization:
\[
 \begin{CD}
 \nbigm(\vecy,L)\times X @>{\Phi_1}>>
 Z_1(\Vtilde,P_{\cdot})@>>>
 Y(W_{\cdot},P_{\cdot})\\
 @VVV @VVV @VVV \\
 M(c,L)\times X^{[\vecy(-c)]}\times X
 @>>>
 Z_2(\Vtilde,P_{\cdot})@>>>
 \Ybar(W_{\cdot})
\times Z(W_{\cdot},P_{\cdot})
 \end{CD}
\]
Thus, we have the following morphisms:
\[
\begin{CD}
 L_{\nbigm(\vecy,L)\times X/M(c,L)\times X^{[\vecy(-c)]}\times X}
@<<<
 \Phi_1^{\ast}L_{Z_1(\Vtilde,P_{\cdot})/Z_2(\Vtilde,P_{\cdot})}
@<{\simeq}<<
 \Phi(V_{\cdot},\phi)^{\ast}
 L_{Y(W_{\cdot},P_{\cdot})/
 \Ybar(W_{\cdot})\times Z(W_{\cdot},P_{\cdot})}
\end{CD}
\]
Therefore,
the morphism $\nu$ is same as
$\gminir(\Vtilde_{\cdot},P_{\cdot},\phi,\phitilde)$
in the subsubsection \ref{subsubsection;06.5.25.2}.
Then, the claim follows from Proposition 
\ref{prop;06.5.24.5}.
\hfill\qed

\vspace{.1in}

Recall that $X^{[\vecy]}$ is smooth.
Let $m$ be a sufficiently large integer
such that $H^i\bigl(X,\nbigl(m)\bigr)=0$ $(i=1,2)$
for any line bundle $\nbigl$ such that $c_1(\nbigl)=c$.
Then the moduli stacks
$\nbigm\bigl(\vecy,\nbigo(-m)\bigr)$
and $M\bigl(c,\nbigo(-m)\bigr)$ are smooth.
Let $\iota:\nbigo(-m)\lrarr L$ be an inclusion.
Then we have the inclusions
$\nbigm(\vecy,L)\lrarr\nbigm\bigl(\vecy,\nbigo(-m)\bigr)$
and $M(c,L)\lrarr M\bigl(c,\nbigo(-m)\bigr)$.
On $M\bigl(c,\nbigo(-m)\bigr)\times X^{[\vecy(-c)]}\times X$,
we put $\nbigk':=\nbigl_c^u\otimes L^{-1}\otimes\nbigo_{\nbigz(\vecy(-c))}$.
Then $\gbigv':=p_{X\,\ast}\nbigk'$ gives the vector bundle
over $M\bigl(c,\nbigo(-m)\bigr)\times X^{[\vecy(-c)]}$
such that
$\gbigv'_{|M(c,L)\times X^{[\vecy(-c)]}}=\gbigv$.
Therefore, we obtain the following Cartesian diagrams:
\[
\begin{CD}
 \nbigm(\vecy,L) @>>>  M(c,L)\times X^{[\vecy(-c)]}
 @>>> M\bigl(c,\nbigo(-m)\bigr) \\
 @VVV @V{i}VV @V{i_1}VV\\
 M(c,L)\times X^{[\vecy(-c)]}
 @>>>
 \gbigv
 @>>>
 \gbigv'
\end{CD}
\]
We have $i^!=i_1^!$.
Therefore,
we obtain the relation 
$i^!\bigl([M(c,L)]\times [X^{[\vecy]}]\bigr)
=[\nbigm(\vecy,L)]$ from (\ref{eq;06.5.13.55}), 
due to Proposition \ref{prop;06.6.11.150}.
The relation $i^!\bigl([M(\widehat{c},[L])]\times [X^{[\vecy]}]\bigr)
=[\nbigm(\vecyhat,[L])]$ can be shown by a similar argument.
Thus we finish the proof of Proposition \ref{prop;06.5.13.35}.
\hfill\qed

\subsection{Bradlow Perturbation}
\label{subsection;06.7.3.44}
\subsubsection{Statements}
\label{subsubsection;06.6.4.1}

Let $L$ be a line bundle on $X$.
If $\delta\in\nbigp^{\br}$ is sufficiently small,
we have the projective morphism:
\begin{equation}
 \label{eq;06.6.23.2}
 \gbigf:\nbigm^s(\vecyhat,[L],\alpha_{\ast},\delta)
\lrarr
 \nbigm^{ss}(\vecyhat,\alpha_{\ast})
\end{equation}
To discuss $\gbigf$,
let us consider the following condition
for $(\vecy,L,\alpha_{\ast})$:
\begin{description}
\item[($i$-vanishing condition)]
\index{$i$-vanishing condition}
We have $H^j\bigl(X,E\otimes L^{-1}\bigr)=0$
for any $j\geq i$ and
for any $E_{\ast}\in\nbigm^{ss}(\vecy,\alpha_{\ast})$.
\end{description}
The $1$-vanishing condition obviously implies
the $2$-vanishing condition.

\begin{prop}
\label{prop;06.6.11.200}
Assume the $1$-vanishing condition holds for
$(\vecy,L,\alpha_{\ast})$.
\begin{itemize}
\item
The morphism
$\gbigf:\nbigm^s(\vecyhat,[L],\alpha_{\ast},\delta)
\lrarr \nbigm^{ss}(\vecyhat,\alpha_{\ast})$ is smooth.
\item
Assume, moreover,
that the $1$-stability condition holds for 
$(\vecy,\alpha_{\ast})$.
Then, we have the following relation:
\[
\gbigf^{\ast}\bigl([\nbigm^s(\vecyhat,\alpha_{\ast})]\bigr)
=\bigl[\nbigm^s(\vecyhat,[L],\alpha_{\ast},\delta)\bigr]
\]
\end{itemize}
\end{prop}
\pf
We give only a remark.
The smoothness of $\gbigf$ is clear.
We put
$\nbigm_1:=\nbigm^{ss}(\vecyhat,\alpha_{\ast})$
and 
$\nbigm_2:= \nbigm^s(\vecyhat,[L],\alpha_{\ast},\delta)$.
It is easy to obtain the following morphism
of distinguished triangles:
\[
 \begin{CD}
 L_{\nbigm_2/\nbigm_1}[-1]
 @>>>
 \gbigf^{\ast}L_{\nbigm_1/k }
 @>>>
 L_{\nbigm_2/k }
 @>>>
 L_{\nbigm_2/\nbigm_1}\\
 @A{\simeq}AA @AAA @AAA @A{\simeq}AA \\
 \Ob_{\rel}(m,\vecy,[L])[-1]
@>>>
 \gbigf^{\ast}\bigl(\Ob(m,\vecy)\bigr)
 @>>>
 \Ob(m,\vecy,[L])
 @>>>
 \Ob_{\rel}(m,\vecy,[L])
 \end{CD}
\]
Due to Proposition \ref{prop;06.6.11.150},
we obtain the equality
$\Psi^{\ast}[\nbigm_1]=[\nbigm_2]$.
\hfill\qed

\vspace{.1in}

Let $\vecL=(L_1,L_2)$ be a pair of line bundles on $X$.
We take $\delta_i\in\nbigp^{\br}$
such that both of $\delta_i$ are sufficiently small.
If $\delta_2$ is sufficiently smaller than $\delta_1$,
we have the projective morphism:
\[
 \gbigf_1:\nbigm^{s}(\vecyhat,[\vecL],\alpha_{\ast},\vecdelta)
\lrarr \nbigm^{s}(\vecyhat,[L_1],\alpha_{\ast},\delta_1)
\]

\begin{prop}
\label{prop;06.6.11.605}
Assume that
the $1$-vanishing condition
holds for $(\vecy,\alpha_{\ast},L_2)$.
The morphism $\gbigf_1$ is smooth,
and we have the following relation:
\[
 \gbigf_1^{\ast}\bigl(\bigl[
 \nbigm^{s}(\vecyhat,[L_1],\alpha_{\ast},\delta_1)
 \bigr]
 \bigr)
=\bigl[\nbigm^{s}(\vecyhat,[\vecL],\alpha_{\ast},\vecdelta)
 \bigr]
\]
\end{prop}
\pf
It can be shown by an argument
similar to the proof of Proposition \ref{prop;06.6.11.200}.
\hfill\qed

\vspace{.1in}

If the $1$-vanishing condition does not hold,
$\gbigf_1$ is not smooth, in general.
The following proposition will be proved in 
the next subsubsections
\ref{subsubsection;06.6.11.600}--\ref{subsubsection;06.6.11.601}.

\begin{prop}
\label{prop;06.6.11.250}
Assume that the $2$-vanishing condition holds
for $(\vecy,\alpha_{\ast},L_2)$.
Then, there exists a Deligne-Mumford stack
$\gbigb$ over $\nbigm^s(\vecyhat,[L_1],\alpha_{\ast},\delta_1)$
with the vector bundle $\gbigv$ and the section $\psi$
such that the following holds:
\begin{itemize}
\item
 The morphism
 $\gbigg:\gbigb\lrarr \nbigm^s(\vecyhat,[L_1],\alpha_{\ast},\delta_1)$
 is smooth.
\item
 $\nbigm^{s}(\vecyhat,[\vecL],\alpha_{\ast},\vecdelta)$
 is $\psi^{-1}(0)$.
\item
 We have the following relation:
\[
 \psi^!\bigl(
 \gbigg^{\ast}[\nbigm^s(\vecyhat,[L_1],\alpha_{\ast},\delta_1)]
\bigr)
=[\nbigm^s(\vecyhat,[\vecL],\alpha_{\ast},\vecdelta)]
\]
Here, $\psi^!$ denotes the Gysin map
for the inclusion
$\nbigm^s(\vecyhat,[\vecL],\alpha_{\ast},\vecdelta)
\lrarr \gbigv$.
\end{itemize}
\end{prop}

Before going into the proof of Proposition 
\ref{prop;06.6.11.250},
we give a corollary.
Let $P$ be any point of
 $\nbigm^{s}(\vecyhat,[L_1],\alpha_{\ast},\delta_1)$.
The fiber $\gbigg^{-1}(P)$ is smooth,
and $\gbigf_1^{-1}(P)$ is the $0$-set of
the section $\psi_{|\gbigg^{-1}(P)}$
of $\gbigv_{|\gbigg^{-1}(P)}$.
Therefore, we obtain the following.

\begin{cor}
\label{cor;06.6.11.701}
For any $k$-valued point $P\in\nbigm^s(\vecyhat,[L_1],\alpha_{\ast},\delta)$,
the fiber $\gbigf_1^{-1}(P)$ is provided with
the perfect obstruction theory.
We also have the following formula:
\[
\int_{\nbigm^s(\vecyhat,[\vecL],\alpha_{\ast},\vecdelta)}
 \Phi\cdot\Psi
=\int_{\gbigf_1^{-1}(P)}\Psi
\cdot
 \int_{\nbigm^s(\vecyhat,[L_1],\alpha_{\ast},\delta_1)}
 \Phi
\]
Here, $\Phi$ and $\Psi$ are cohomology classes
on $\nbigm^s(\vecyhat,[L_1],\alpha_{\ast},\delta_1)$
and $\nbigm^s(\vecyhat,[\vecL],\alpha_{\ast},\vecdelta)$,
respectively,
such that $\deg(\Phi)=2\dim^f\nbigm^s(\vecy,\alpha_{\ast})$.
(See the subsection {\rm\ref{subsection;06.6.19.110}}
for the cohomology classes and the evaluation
 considered in this paper.)
\hfill\qed
\end{cor}

\begin{rem}
Similar claims also hold
for the morphism $\gbigf$ in {\rm(\ref{eq;06.6.23.2})},
if the $2$-vanishing condition holds
for $(\vecy,L,\alpha_{\ast})$.
The proof is similar.
\hfill\qed
\end{rem}

\subsubsection{Construction of $\gbigb$}
\label{subsubsection;06.6.11.600}

Let $m$ be a sufficiently large 
such that the $1$-vanishing condition
holds for $(\vecy,\nbigo(-m),\alpha_{\ast})$.
We put $L_2':=\nbigo(-m)$ and $\vecL_2':=(L_1,L_2')$.
We put $\gbigb:=\nbigm^s(\vecyhat,[\vecL'],\alpha_{\ast},\vecdelta)$.
Then, we have the natural smooth morphism
$\gbigg:\gbigb\lrarr \nbigm^s(\vecyhat,[L_1],\alpha_{\ast},\delta_1)$.

We take an inclusion $\iota:L_2'\lrarr L_2$
such that the cokernel $L_2/L_2'$
is a line bundle on some smooth divisor of $X$.
It naturally induces the morphism
$L_2^{-1}\lrarr L_2^{\prime\,-1}$.
The cokernel is denoted by $\Cok$.
We also obtain the inclusion
$\nbigm^s(\vecyhat,[\vecL],\alpha_{\ast},\vecdelta)
\lrarr \gbigb$.

We have the morphisms
$\nbigm^s(\vecyhat,[\vecL],\alpha_{\ast},\vecdelta)
\lrarr \nbigm(\vecyhat,[L_i])$ $(i=1,2)$.
The pull back of the relative tautological 
line bundles are denoted by $\nbigo_{\rel}^{(i)}(1)$.
Similarly, we have the morphisms of
$\gbigb$
to $\nbigm(\vecyhat,[L_1])$
and $\nbigm(\vecyhat,[L_2'])$.
The pull back of the relative tautological bundles
are also denoted by $\nbigo_{\rel}^{(1)}(1)$
and $\nbigo_{\rel}^{(2)}(1)$,
respectively.

Let $\Ehat^{\prime\,u}$ denote the universal sheaf
over $\gbigb\times X$.
We have the universal reduced $[L_2']$-section
$[\phi_2]:\nbigo_{\rel}^{(2)}(-1)\otimes L_2'
  \lrarr \Ehat^{\prime\,u}$.
We put
$E^{\prime\,u}:=
 \Ehat^{\prime\,u}\otimes \nbigo^{(2)}_{\rel}(1)$.
We put as follows:
\[
 \gbigv:=p_{X\,\ast}\bigl(
 E^{\prime\,u}\otimes \Cok
 \bigr).
\]
It is easy to see 
$R^ip_{X\,\ast}\bigl(E^{\prime\,u}\otimes\Cok\bigr)=0$
for $i=1,2$,
and hence $\gbigv$ is a locally free sheaf on $\gbigb$.
The universal reduced $L_2'$-section
$[\phi^u_2]$ induces the section $\psi$ of $\gbigv$
over $\gbigb$.
It is easy to observe that 
$\psi^{-1}(0)$ is isomorphic to
$\nbigm^s(\vecyhat,[\vecL],\alpha_{\ast},\vecdelta)$.

Thus, we obtain the following Cartesian diagram:
\begin{equation}
\label{eq;06.5.25.10}
 \begin{CD}
 \nbigm^s(\vecyhat,[\vecL],\alpha_{\ast},\vecdelta) 
 @>{i_1}>>\gbigb \\
 @V{j}VV @VV{\psi}V \\
 \gbigb  @>{i}>> \gbigv
 \end{CD}
\end{equation}
Here $i$ denotes the $0$-section.
For the proof of Proposition \ref{prop;06.6.11.250},
we have only to show
$i^![\gbigb]=[\nbigm^s(\vecyhat,[\vecL],\alpha_{\ast},\vecdelta)]$.

\subsubsection{Compatibility of the obstruction theories}

\begin{lem}
\label{lem;06.6.11.500}
The obstruction theories of
$\nbigm^s(\vecyhat,[\vecL],\alpha_{\ast},\vecdelta)$
and $\gbigb$ are compatible over $i$
in the diagram {\rm(\ref{eq;06.5.25.10})}.
\end{lem}
\pf
We discuss in some more general situation.
Let $L$ be a line bundle on $X$.
We take a sufficiently large integer $m$.
We take an inclusion $\iota:\nbigo(-m)\lrarr L$.
It naturally induces the morphism
$L^{-1}\lrarr \nbigo(m)$.
We assume that the cokernel $\Cok$
is a line bundle on some smooth divisor of $X$.

We use the notation $\nbigm_1$ and $\nbigm_2$
to denote the moduli stacks
$\nbigm(m,\vecyhat,[L])$
and $\nbigm(m,\vecyhat,[\nbigo(-m)])$.
The inclusion $\nbigm(m,\vecyhat,[L])
\lrarr \nbigm(m,\vecyhat,[\nbigo(-m)])$
is induced by $\iota$.
Let $\Ehat^u_2$ denote the universal sheaf 
over $\nbigm_2\times X$.
We put $E^u_2:=\Ehat^u_2\otimes\nbigo_{\rel}(1)$.
We put $\gbigv:=p_{X\,\ast}\bigl(\Ehat_2^u\otimes\Cok\bigr)$.
The universal reduced $\nbigo(-m)$-section of $\Ehat^u_2$
induces the section $\psi$ of $\gbigv$.
We have $\nbigm_1\simeq \psi^{-1}(0)$
and the following commutative diagram:
\begin{equation}
\label{eq;06.6.11.300}
 \begin{CD}
 \nbigm_1 @>{i_1}>> \nbigm_2\\
 @V{j}VV @V{\psi}VV \\
 \nbigm_2 @>{i}>>\gbigv
 \end{CD}
\end{equation}
To prove the claim of the lemma,
we have only to show that 
the obstruction theories of $\nbigm_i$
are compatible over $i$ 
in the diagram (\ref{eq;06.6.11.300}).

Let $\Ehat_1^u$ denote the universal sheaf
on $\nbigm_1$, 
and let  $[\phi_1^u]$
denote the universal reduced $L$-section.
We put
$\nbigv_{1,0}:=p_X^{\ast}p_{X\,\ast}\Ehat_1^u(m)$,
and the kernel of the surjection
$\nbigv_{1,0}\lrarr\Ehat_1^u(m)$ is denoted by
$\nbigv_{1,-1}$.
We take a locally free resolution $P_{\cdot}$
of $L(m)$ such that $P_0$ is a direct sum of
$\nbigo_X$.
We have the canonical lift
$[\phitilde_1]:
 p_{\nbigm_1}^{\ast}P_{\cdot}\otimes
 p_X^{\ast}\nbigo_{\rel}(-1)
 \lrarr \nbigv_{1,\cdot}$.
Recall 
$\gminig'_{\rel}(\nbigv_{1,\cdot},[\phitilde_1]):=
 \nhom\bigl(
 p_{\nbigm_1}^{\ast}P_{\cdot}\otimes
 p_X^{\ast}\nbigo_{\rel}(-1),\,
 \nbigv_{1,\cdot}
\bigr)^{\lor}$
and 
$\gminig(\nbigv_{1,\cdot}):=
 \nhom\bigl(\nbigv_{1,\cdot},
 \nbigv_{1,\cdot}\bigr)^{\lor}[-1]$.

Similarly,
we take a locally free resolution
$\nbigv_{2\,\cdot}$ of 
the universal sheaf $\Ehat_2$ on $\nbigm_2\times X$.
Let $[\phi_2]$ denote the universal reduced
$\nbigo(-m)$-section of $\Ehat_2$.
We have the canonical lift
$[\phitilde_2]:p_{\nbigm_2}^{\ast}\nbigo_X\otimes
 p_X^{\ast}\nbigo_{\rel}(-1)
\lrarr \nbigv_{2,0}$ of $[\phi_2]$.
In this case, we have 
$\gminig'_{\rel}(\nbigv_{2\cdot},[\phitilde_2])=
 \nhom\bigl(  p_{\nbigm_2}^{\ast}\nbigo
 \otimes p_{X}^{\ast}\nbigo_{\rel}(1),\,\,
\nbigv_{2,\cdot}
 \bigr)^{\lor}$,
and $\gminig(\nbigv_{2\,\cdot})
=\nhom\bigl(\nbigv_{2\,\cdot},\nbigv_{2\,\cdot}\bigr)^{\lor}[-1]$.

The inclusion $\iota:\nbigo\lrarr L(m)$
has the canonical lift $\nbigo_X\lrarr P_{\cdot}$.
Therefore, we have the following commutative diagram
on $\nbigm_1$:
\begin{equation}
\label{eq;06.5.14.15}
 \begin{CD}
 i_1^{\ast}\gminig'_{\rel}(\nbigv_{2\,\cdot},[\phitilde_2])[-1]
 @>>>
 \gminig'_{\rel}(\nbigv_{1\,\cdot},[\phitilde_1])[-1]\\
 @VVV @VVV \\
 i_1^{\ast}\gminig(\nbigv_{2\,\cdot})
 @>>>
 \gminig(\nbigv_{1\,\cdot})
 \end{CD}
\end{equation}

\begin{lem}
\label{lem;06.6.11.400}
The diagram {\rm(\ref{eq;06.5.14.15})}
is compatible with the morphisms to 
the cotangent complexes:
\[
\begin{CD}
 i_1^{\ast}L_{\nbigm_2\times X/
 (\nbigm(m,\vecy)\times X)_{G_m}}[-1]
 @>>>
 L_{\nbigm_1\times X/
 (\nbigm(m,\vecy)\times X)_{G_m}}[-1]\\
 @VVV @VVV \\
 i_1^{\ast}\pi_2^{\ast}L_{\nbigm(m,\vecy)\times X/X_{G_m}}
 @>>>
 \pi_1^{\ast}L_{\nbigm(m,\vecy)\times X/X_{G_m}}
\end{CD}
\]
Here, $\pi_i$ denote the natural morphism
of $\nbigm_i$ to $\nbigm(m,\vecy)$.
(See the subsubsection {\rm\ref{subsubsection;06.5.22.15}}
for the compatibility of diagrams.)
\end{lem}
\pf
We give only an indication.
We use the notation in the subsubsection
\ref{subsubsection;06.5.4.110}.
We construct the stack $Y(W_{\cdot},[\nbigo_X])$
by replacing $P_{\cdot}$ with the complex
$(0\rarr \nbigo_X)$.
Then we have the following commutative diagram:
\begin{equation}
 \label{eq;06.5.25.15}
 \begin{CD}
 \nbigm_1\times X@>>> 
 Y(W_{\cdot},[P_{\cdot}]) @>>> 
 X_{G_m}\\
 @VVV @VVV @VVV \\
 \nbigm_2\times X @>>> 
 Y(W_{\cdot},[\nbigo_X]) @>>>
 X_{G_m}\\
 @VVV @VVV @VVV \\
 (\nbigm(m,\vecy)\times X)_{G_m} @>>>
 Y(W_{\cdot})_{G_m} @>>>
 X_{G_m}
 \end{CD}
\end{equation}
Then the desired compatibility follows from
the construction of the complexes.
\hfill\qed

\begin{lem}
We have the following commutative diagram:
\begin{equation}
 \label{eq;06.5.25.20}
 \begin{CD}
 i_1^{\ast}\Ob\bigl(m,\vecy,[\nbigo(-m)]\bigr)
@>>>
 \Ob(m,\vecy,[L])  \\
 @VVV @VVV \\
 i_1^{\ast}L_{\nbigm_2}
 @>>>
 L_{\nbigm_1}
 \end{CD}
\end{equation}
\end{lem}
\pf
We obtain the following commutative diagram
from (\ref{eq;06.5.14.15}):
\begin{equation}
 \label{eq;06.6.11.444}
 \begin{CD}
 i_1^{\ast}\Ob'_{\rel}\bigl(\nbigv_{2\,\cdot},[\phitilde_2]\bigr)[-1]
 @>>>
 \Ob'_{\rel}\bigl(\nbigv_{1\cdot},[\phitilde_1]\bigr)[-1]\\
 @VVV @VVV\\
 i_1^{\ast}\Ob\bigl(\nbigv_{2,\cdot}\bigr)
 @>>>
 \Ob\bigl(\nbigv_{1\,\cdot}\bigr)
 \end{CD} 
\end{equation}
It is compatible with the morphisms
to the cotangent complexes
due to Lemma \ref{lem;06.6.11.400}:
\[
 \begin{CD}
 i_1^{\ast}L_{\nbigm_2/\nbigm(m,\vecy)_{G_m}}[-1]
 @>>>
 L_{\nbigm_1/\nbigm(m,\vecy)_{G_m}}[-1]\\
 @VVV @VVV\\
 i_1^{\ast}\pi_2^{\ast}L_{\nbigm(m,\vecy)_{G_m}}
 @>>>
 \pi_1^{\ast}L_{\nbigm(m,\vecy)_{G_m}}
 \end{CD}
\]
By a modification as in the subsubsection
\ref{subsubsection;06.5.4.110},
we obtain the commutative diagram
from (\ref{eq;06.6.11.444}):
\[
  \begin{CD}
 i_1^{\ast}\Ob_{\rel}\bigl(\nbigv_{2\,\cdot},[\phitilde_2]\bigr)[-1]
 @>>>
 \Ob_{\rel}\bigl(\nbigv_{1\cdot},[\phitilde_1]\bigr)[-1]\\
 @VVV @VVV\\
 i_1^{\ast}\Ob\bigl(\nbigv_{2,\cdot}\bigr)
 @>>>
 \Ob\bigl(\nbigv_{1\,\cdot}\bigr)
 \end{CD} 
\]
It is compatible with the following commutative diagram:
\[
   \begin{CD}
 i_1^{\ast}L_{\nbigm_2/\nbigm(m,\vecy)}[-1]
 @>>>
 L_{\nbigm_1/\nbigm(m,\vecy)}[-1]\\
 @VVV @VVV\\
 i_1^{\ast}\pi_2^{\ast}L_{\nbigm(m,\vecy)}
 @>>>
 \pi_2^{\ast}L_{\nbigm(m,\vecy)}
 \end{CD}
\]
Then we obtain the desired commutative diagram
(\ref{eq;06.5.25.20})
by construction.
\hfill\qed

\vspace{.1in}

We put
$\nbigvtilde_{1,\cdot}:=
 \nbigv_{1,\cdot}\otimes \nbigo_{\rel}(1)$.
We put as follows:
\[
 \gminig_{\rel}:=
 \cone\Bigl(
 i_1^{\ast}\nhom\bigl(
p_{\nbigm_2}^{\ast}\nbigo_X,\,
 \nbigvtilde_{2,\cdot}
 \bigr)^{\lor}
 \lrarr
 \nhom\bigl(
 p_{\nbigm_1}^{\ast}P_{\cdot},\,
 \nbigvtilde_{1,\cdot}
 \bigr)^{\lor}
 \Bigr),
\quad\quad
\Ob_{\rel}:=Rp_{X\,\ast}\bigl(
  \gminig_{\rel}\otimes\omega_X\bigr)
\]
It is easy to observe that
$\Ob_{\rel}$ naturally isomorphic to
$\cone\bigl(
\Ob(m,\vecyhat,[\nbigo(-m)])
\lrarr
 \Ob(m,\vecyhat,[L])
 \bigr)$.
From the diagram (\ref{eq;06.5.25.15}),
we obtain the following morphism:
\begin{equation}
 \label{eq;06.6.11.445}
\begin{CD}
 L_{\nbigm_1\times X/\nbigm_2\times X}
@<<<
 \Phi(\nbigv_{\cdot},[\phitilde])^{\ast}
 L_{Y(W_{\cdot},[P_{\cdot}])/Y(W_{\cdot},[\nbigo_X])}
@<<<
 \gminig_{\rel}
\end{CD}
\end{equation}
It is easy to observe that
the composite of (\ref{eq;06.6.11.445})
induces the morphism
$\beta_1:\Ob_{\rel}\lrarr L_{\nbigm_1/\nbigm_2}$,
which is same as the morphism induced from
(\ref{eq;06.5.25.20}).

We have the naturally defined morphism:
\[
 h:N\bigl(\nbigo_{\rel}(-1)\otimes P_{-1}, \nbigv_{1,-1}\bigr)
 \times 
 N\bigl(\nbigo_{\rel}(-1)\otimes P_{0},\nbigv_{1,0}\bigr)
\lrarr 
 N\bigl(\nbigo_{\rel}(-1)\otimes P_{-1},\nbigv_{1,0}\bigr)
\]
We put $Z_1:=h^{-1}(0)$.
We also put $Z_2:=N\bigl(\nbigo_{\rel}(-1),\nbigv_{2,0}\bigr)$
which is the vector bundle over $\nbigm_2\times X$.
Then, the lifts $[\phitilde_{i}]$ induce
the following factorization:
\[
 \begin{CD}
 \nbigm_1\times X @>{\Phi_1}>>
 Z_1 @>>> Y(W_{\cdot},[P_{\cdot}])\\
 @VVV @VVV @VVV\\
 \nbigm_2\times X @>>>
 Z_2 @>>> Y(W_{\cdot},[\nbigo_X])
 \end{CD}
\]
We obtain the following factorization:
\[
\begin{CD}
 \Phi(\nbigv_{\cdot},[\phitilde])^{\ast}
 L_{Y(W_{\cdot},[P_{\cdot}])/Y(W_{\cdot},[\nbigo_X])}
@>{\simeq}>>
 \Phi_1^{\ast}L_{Z_1/Z_2}
@>>>
 L_{\nbigm_1\times X/\nbigm_2\times X}
\end{CD} 
\]
We put $Z_3:=N(p_{\nbigm_2}^{\ast}\nbigo_X,\,p_X^{\ast}\gbigv)$,
which is the vector bundle over $\nbigm_2\times X$.
We have the naturally defined morphism
$\nbigv_{2,0}\otimes\nbigo_{\rel}(1)\lrarr p_X^{\ast}\gbigv$.
Then, we obtain the following diagram: 
\[
 \begin{CD}
 \nbigm_1\times X
 @>{\Phi_1}>>
 Z_1 @>{\Phi_2}>> 
 \nbigm_2\times X\\
 @VVV @VVV @VVV\\
\nbigm_2\times X @>>>
Z_2 @>>> Z_3
\end{CD}
\]
We naturally have
$\Phi_1^{\ast}\Phi_2^{\ast}L_{\nbigm_2\times X/Z_3}
\simeq j_X^{\ast}L_{\nbigm_2\times X/p_X^{\ast}\gbigv}$.
Thus, we obtain the morphism:
$\alpha:j_X^{\ast}L_{\nbigm_2\times X/p_X^{\ast}\gbigv}
\lrarr
 \Phi_1^{\ast}L_{Z_1/Z_2}$.
It can be checked that 
$\alpha$ factors through $\gminig_{\rel}$.
Namely, we have the following morphisms:
\[
 j_X^{\ast}L_{\nbigm_2\times X/p_X^{\ast}\gbigv}
\lrarr
 \gminig_{\rel}
\lrarr 
 L_{\nbigm_1\times X/\nbigm_2\times X}.
\]
It induces the morphism:
\[
\begin{CD}
 j^{\ast}L_{\nbigm_2/\gbigv}\otimes Rp_{X\,\ast}\omega_X
@>{\beta_0}>>
 \Ob_{\rel}
@>{\beta_1}>>
 L_{\nbigm_1/\nbigm_2}
\end{CD}
\]
It is easy to observe that $\beta_0$
induces the isomorphism
$\beta_2:j^{\ast}L_{\nbigm_2/\gbigv}\lrarr \Ob_{\rel}$,
and the composite
$\beta_1\circ\beta_2$ is same as
the morphism obtained from the diagram (\ref{eq;06.6.11.300}).
Thus the proof of Lemma \ref{lem;06.6.11.500}
is finished.
\hfill\qed

\subsubsection{Ambient smooth stack}
\label{subsubsection;06.6.11.601}

Since $\gbigb$ is not smooth in general,
we construct an ambient smooth stack
to use Proposition \ref{prop;06.6.11.150}.
Let $C$ denote the support of $\Cok$.
Recall that $C$ is smooth 
and that $\Cok$ is isomorphic to
a line bundle $\nbigl_C$ on $C$,
due to our choice of $\iota$.

For a $k$-scheme $T$,
let $F(T)$ denote the set of  the quotients
$q:p_T^{\ast}V_{m,C}\lrarr \nbige$
over $T\times C$,
satisfying the following condition:
\begin{itemize}
\item
 $\nbige$ is flat over $T$,
 and the type of $\nbige$ is
 $y_{|C}\cdot\ch(\nbigl_C)$.
\item
 For any point $u\in T$,
 $H^1\bigl(C,\nbige_{|\{u\}\times C}\bigr)=0$.
\end{itemize}
Then, we obtain the functor $F$
of the category of $k$-schemes to the category of sets.
The functor $F$ is representable
by a scheme $Q^{\circ}_1$.

\begin{lem}
The scheme $Q^{\circ}_1$ is smooth.
\end{lem}
\pf
We have the perfect obstruction theory $\Ob\bigl(Q^{\circ}_1\bigr)$
of $Q^{\circ}_1$
(Proposition \ref{prop;06.5.2.1}).
Let $z=(q,\nbige)$ be a point of $Q^{\circ}_1$,
and let $i_z$ denote the inclusion $\{z\}\lrarr Q^{\circ}_1$.
Let $\nbigk$ denote the kernel of $V_{m,C}\lrarr \nbige$.
We have only to show
$\nbigh^{-1}\bigl(i_z^{\ast}\Ob(Q^{\circ}_1)\bigr)=0$.

The dual
$\nbigh^i\bigl(i_z^{\ast}\Ob(Q^{\circ}_1)\bigr)^{\lor}$
is isomorphic to 
$\Ext^{-i}\bigl(\nbigk,\nbige\bigr)$.
We have the exact sequence
$\Ext^1\bigl(V_{m,C},\nbige\bigr)
\lrarr
 \Ext^1\bigl(\nbigk,\nbige\bigr)
\lrarr
 \Ext^2\bigl(\nbige,\nbige\bigr)$.
We have the vanishing of the first term
by definition of $Q^{\circ}_1$.
Since $C$ is a smooth curve,
we also have the vanishing of the third term.
Therefore, we obtain the desired vanishing,
and hence $Q^{\circ}_1$ is smooth.
\hfill\qed

\vspace{.1in}

We have the universal quotient sheaf
$p_{Q^{\circ}_1}V_{m,C}\lrarr \nbigc$
on $Q^{\circ}_1\times C$.
The push-forward $p_{C\,\ast}\nbigc$
gives the vector bundle on $Q^{\circ}_1$.
We denote it by $\gbigv_1$.

We use the notation in the subsubsections
\ref{subsubsection;06.5.14.30}.
We put 
$Q^{ss}(m,\vecyhat,\alpha_{\ast}):=
 Q^{ss}(m,\vecy,\alpha_{\ast})\times_{Q(m,y)}Q(m,\yhat)$.
Let $\Zhat_m$ be as in (\ref{eq;06.5.21.30}).
We have the $\GL(V_m)$-closed immersion
$Q^{ss}(m,\vecyhat,\alpha_{\ast})
\lrarr \Zhat_m\times \prod_i G_{m,i}$.
(See the subsubsection \ref{subsubsection;06.5.14.1}
for $G_{m,i}$.)

We also have the naturally defined
$\GL(V_m)$-equivariant morphism
$Q^{ss}(m,\vecyhat,\alpha_{\ast})
\lrarr Q^{\circ}_1$
by the correspondence
$(q,\nbige,F_{\ast},\rho)\longmapsto
 (q',\nbige(-m)\otimes\Cok)$,
where $q'$ denotes the naturally induced map
$V_{m,C}\lrarr \nbige(-m)\otimes\Cok$.

We have the natural right $\GL(V_m)$-action
on $\Zhat_m\times \prod_{i}G_{m,i}\times Q^{\circ}_1
 \times \proj_m^{(1)}\times\proj_m^{(2)}$.
The quotient stack is denoted by $\gbigbtilde'$.
The bundle $\gbigv_1$ induce the $\GL(V_m)$-vector bundle
on $\Zhat_m\times\prod_iG_{m,i}\times\proj_m\times Q^{\circ}_1$,
and hence the vector bundle $\gbigvtilde$ on $\gbigbtilde'$.

We have the $\GL(V_m)$-equivariant immersion
$Q^{\circ}(m,\vecyhat,[L_1])\times \proj_m^{(2)}
\lrarr \Zhat_m\times\prod_iG_{m,i}\times 
 Q^{\circ}_1\times\proj_m^{(1)}\times \proj_m^{(2)}$,
which induces the immersion
$\gbigb\lrarr \gbigbtilde'$.
Since $\gbigb$ is Deligne-Mumford,
we can take an open neighbourhood $\gbigbtilde$
of $\gbigb$ in $\gbigbtilde'$,
which is Deligne-Mumford and smooth.
The restriction of $\gbigvtilde$ to $\gbigbtilde$ 
is denoted by $\gbigvtilde$.
By the construction,
the restriction of $\gbigvtilde$
to $\gbigb$ is same as $\gbigv$.
We put $\nbigm=\nbigm^{ss}(\vecyhat,[\vecL],\alpha_{\ast},\vecdelta)$.
Then, we obtain the following diagram:
\[
 \begin{CD}
 \nbigm @>>> \gbigb @>>> \gbigbtilde\\
 @V{i_1}VV @V{i}VV @V{i_2}VV \\
 \gbigb @>{\psi}>> \gbigv @>>>\gbigvtilde
 \end{CD}
\]
Here $i$ and $i_2$ denote the $0$-section.
We have $i^!=i_2^!$.
The obstruction theories of $\nbigm$ and $\gbigb$
are compatible (Lemma \ref{lem;06.6.11.500}).
Then, we obtain
$i^!\bigl([\nbigm_2]\bigr)=[\nbigm_1]$,
due to Proposition \ref{prop;06.6.11.150}.
Thus, the proof of Proposition \ref{prop;06.6.11.200}
is finished.
\hfill\qed

\subsection{Comparison with Full Flag Bundles}

Let $\vecy$ be an element of $\Type$.
Let $\alpha_{\ast}$ be a system of weights.
Let $L$ be a line bundle on $X$.
We take a sufficiently large integer $m$.
Let $\nbigmtilde^{s}(\vecyhat,[L],\alpha_{\ast},\delta)$
denote the full flag bundle
associated to the bundle $p_{X\,\ast}\Ehat^u(m)$.
We have the natural smooth morphism:
\[
 \gbigf_1:\nbigmtilde^s(\vecyhat,[L],\alpha_{\ast},\delta)
\lrarr \nbigm^s(\vecyhat,[L],\alpha_{\ast},\delta)
\]

\begin{lem}
\label{lem;06.6.13.150}
We have the following relation:
\[
\gbigf_1^{\ast}\bigl(
[\nbigm^s(\vecyhat,[L],\alpha_{\ast},\delta)]\bigr)
=[\nbigmtilde^s(\vecyhat,[L],\alpha_{\ast},\delta)]
\]
\end{lem}
\pf
Let $\pi$ denote the projection
$\nbigmtilde(m,\vecyhat,[L])\lrarr \nbigm(m,\vecyhat,[L])$.
By construction in the subsubsection
\ref{subsubsection;06.6.13.40},
we have the following morphism of the distinguished triangles
on $\nbigmtilde(m,\vecyhat,[L])$:
\[
 \begin{CD}
 \pi^{\ast}\Ob(m,\vecyhat,[L])
 @>>>
 \Obtilde(m,\vecyhat,[L])
 @>>>
 L_{\nbigmtilde(m,\vecyhat,[L])/\nbigm(m,\vecyhat,[L])}
 @>>>
 \pi^{\ast}\Ob(m,\vecyhat,[L])[1] \\
 @VVV @VVV @V{\simeq}VV @VVV \\
 \pi^{\ast}L_{\nbigm(m,\vecyhat,[L])}
 @>>>
 L_{\nbigmtilde(m,\vecyhat,[L])}
 @>>>
 L_{\nbigmtilde(m,\vecyhat,[L])/\nbigm(m,\vecyhat,[L])}
 @>>>
  \pi^{\ast}L_{\nbigm(m,\vecyhat,[L])}[1]
 \end{CD}
\]
Hence we obtain the compatibility of
the obstruction theories $\Ob(m,\vecyhat,[L])$
and $\Obtilde(m,\vecyhat,[L])$
of $\nbigm^s(\vecyhat,[L],\alpha_{\ast},\delta)$
and $\nbigmtilde^s(\vecyhat,[L],\alpha_{\ast},\delta)$.
Then, the claim follows from
Proposition \ref{prop;06.6.11.150}.
\hfill\qed

\vspace{.2in}

Let $\vecy$ and $\alpha_{\ast}$ be as above.
We take a sufficiently large integer $m$.
Let $\nbigmtilde^s(\vecyhat,\alpha_{\ast},+)$
be as in the subsubsection \ref{subsubsection;06.6.13.51}.
Let $\epsilon$ be sufficiently small positive number.
Then, we have the naturally defined morphism:
\[
 \gbigf_2:\nbigmtilde^s(\vecyhat,\alpha_{\ast},+)
\lrarr \nbigm^s(\vecyhat,[\nbigo(-m)],\alpha_{\ast},\epsilon)
\]

\begin{lem}
\label{lem;06.6.13.201}
We have the following relation:
\begin{equation}
 \label{eq;06.6.13.200}
 \gbigf_2^{\ast}\bigl(
[\nbigm^s(\vecyhat,[\nbigo(-m)],\alpha_{\ast},\epsilon)]
\bigr)
=[\nbigmtilde^s(\vecyhat,\alpha_{\ast},+)]
\end{equation}
\end{lem}
\pf
In this case,
we have the other way of construction for
the obstruction theory
$\Ob(m,\vecyhat,[\nbigo(-m)])$
on $\nbigm(m,\vecyhat,[\nbigo(-m)])$
by the method in the subsubsection
\ref{subsubsection;06.5.8.5}.
Let $V_m$ be an $H_y(m)$-dimensional vector space.
We put $B(W):=k_{\GL(V_m)}$
and $B(W,[P]):=\proj(V_m^{\lor})_{\GL(V_m)}$
as in the subsubsections \ref{subsubsection;06.5.4.100}
and \ref{subsubsection;06.5.4.150}.
Let $\Ehat^u$ denote the universal sheaf
over $\nbigm(m,\vecyhat)\times X$.
Then, $\nbigm(m,\vecyhat,[\nbigo(-m)])$
is the projectivization of
$p_{X\,\ast}(\Ehat^u(m))$.
Hence, we have the following naturally defined
Cartesian diagram:
\[
 \begin{CD}
 \nbigm(m,\vecyhat,[\nbigo(-m)]) @>{\Psi}>> B(W,[P])\\
 @V{\pi}VV @VVV \\
 \nbigm(m,\vecyhat) @>{\varphi}>> B(W)
 \end{CD}
\]
Since $\varphi$ factors through
$\Ob(m,\vecyhat)$,
we obtain the following morphism:
\[
\begin{CD}
 \Psi^{\ast}L_{B(W,[P])/B(W)}[-1]
 @>>>
 \pi^{\ast}L_{B(W)}
 @>>> 
 \pi^{\ast}\Ob(m,\vecyhat)
\end{CD}
\]
The cone of the composite is denoted by
$\Ob_2(m,\vecyhat,[\nbigo(-m)])$.
Then, we obtain the morphism:
\[
 \ob_2(m,\vecyhat,[\nbigo(-m)]):
 \Ob_2(m,\vecyhat,[\nbigo(-m)])
\lrarr
 L_{\nbigm(m,\vecyhat,[\nbigo(-m)])}
\]
We use the following lemma.
\begin{lem}
\label{lem;06.6.23.1}
We have $\ob_2(m,\vecyhat,[\nbigo(-m)])=\ob(m,\vecyhat,[\nbigo(-m)])$
in $D(\nbigm(m,\vecyhat,[\nbigo(-m)]))$.
\end{lem}
\pf
We use the result in the subsubsection 
\ref{subsubsection;06.5.4.150}.
We take the canonical locally free resolution
of $\Ehat^u(m)$.
Namely,
we put $\nbigv_0:=p_X^{\ast}\bigl(p_{X\,\ast}\Ehat^u(m)\bigr)$,
and the $\nbigv_{-1}$ denotes the kernel
of the canonical morphism
$\nbigv_0\lrarr \Ehat^u(m)$.
The reduced $\nbigo_X$-section $[\phi]$ of $\Ehat^u(m)$
is canonically lifted to the reduced $\nbigo_X$-section
$[\phitilde]$ of $\nbigv_0$.
In this case, the diagram (\ref{eq;06.5.4.300})
is as follows:
\begin{equation}
 \begin{CD}
 L_{\nbigm(m,\vecyhat,[L])/\nbigm(m,\vecyhat)}[-1] @<<<
 \Ob_{\rel}(\nbigv_{\cdot},[\phitilde])[-1] @<<{a_1}< 
 \Ob^G_{\rel}(\nbigv_{\cdot},[\phitilde])[-1]\\
 @VVV @VVV @VVV \\
 \pi^{\ast}L_{\nbigm(m,\vecyhat)}
 @<<<
 \Ob(\nbigv_{\cdot}) @<<< 
 \Ob^G(\nbigv_{\cdot})
 \end{CD}
\end{equation}
The diagram (\ref{eq;06.5.4.250}) is as follows:
\begin{equation}
 \begin{CD}
 \pi^{\ast}L_{\nbigm(m,\vecyhat)} @<<<
 \Phi(\nbigv_{\cdot})^{\ast}L_{B(W)/k} @<{\tau_1}<<
 \Ob^G(V_{\cdot}) \\
 @AAA @AAA @AAA \\
 L_{\nbigm(m,\vecyhat,[L])/\nbigm(m,\vecyhat)}[-1] @<<<
 \Psi^{\ast}
 L_{B(W,[P])/B(W)}[-1]
 @<{\tau_3}<<
 \Ob^{G}_{\rel}(\nbigv_{\cdot},[\phitilde])[-1]
 \end{CD}
\end{equation}
The morphisms $\tau_1$ and $\tau_3$ are isomorphic,
in general (Lemma \ref{lem;06.5.4.201}
and Lemma \ref{lem;06.5.4.200}).
The morphism $a_1$ is also isomorphic in this case.
Then, the claim of Lemma \ref{lem;06.6.23.1}
immediately follows.
\hfill\qed

\vspace{.1in}

It is easy to see the compatibility
of the obstruction theories
$\obtilde(m,\vecyhat)$
and $\obtilde_2(m,\vecyhat,[\nbigo(-m)])$
on the moduli stacks $\nbigmtilde^s(\vecyhat,\alpha_{\ast},+)$
and $\nbigm^s(\vecyhat,[\nbigo(-m)],\alpha_{\ast},\epsilon)$.
Hence, we obtain (\ref{eq;06.6.13.200})
due to Proposition \ref{prop;06.6.11.150}.
Thus the proof of Lemma \ref{lem;06.6.13.201}
is finished.
\hfill\qed

\subsection{Parabolic Perturbation}
\label{subsection;06.6.1.11}
Since we do not use the result in this subsection later,
the reader can skip here.

\subsubsection{Statement}

Let $y$ be an element of $H^{\ast}(X)$.
Let $\vecy$ be an element of $\Type$
whose $H^{\ast}(X)$-component is $y$,
and let $\alpha_{\ast}$ denote a system of weights.
Let us discuss the relation of the obstruction theories
and the virtual fundamental classes of
the moduli stacks $\nbigm^s(\yhat)$ and
$\nbigm^{s}(\vecyhat,\alpha_{\ast})$.

We assume that
any semistable torsion-free sheaf of type $y$
is also $\mu$-stable,
and that
$\alpha_i$ are sufficiently close to $1$.
Then we have the morphism
$\gbigf:\nbigm^s(\vecyhat,\alpha_{\ast})
\lrarr \nbigm^s(\yhat)$.

\begin{prop}
\label{prop;06.5.13.100}
There exists a Deligne-Mumford stack
$\gbigb(\vecyhat,\alpha_{\ast})$ 
over $\nbigm^s(\yhat)$
with a vector bundle $\gbigv$ and the section $\psi$,
such that the following holds:
\begin{itemize}
\item
 The morphism $\gbigg:\gbigb(\vecyhat,\alpha_{\ast})
 \lrarr \nbigm^s(\yhat)$ is smooth.
\item
 $\nbigm^s(\vecyhat,\alpha_{\ast})$
 is isomorphic to $\psi^{-1}(0)$.
\item
 We have the following relation:
\begin{equation}
\label{eq;06.5.14.5}
 \psi^!\bigl(\gbigg^{\ast}[\nbigm^s(\yhat)]\bigr)
=\bigl[\nbigm^s(\vecyhat,\alpha_{\ast})\bigr]
\end{equation}
Here $\psi^{!}$ denote the Gysin map for
the inclusion $\nbigm^s(\vecyhat,\alpha_{\ast})\lrarr \gbigv$.
\end{itemize}
\end{prop}

The proof will be given in the next subsubsections.
Before going into the proof,
we give some remarks.

\begin{cor}
Let $P$ be any point of $\nbigm^s(\yhat)$.
The fiber $\gbigg^{-1}(P)$ is smooth,
and $\gbigf^{-1}(P)$ is the $0$-set of
the section $\psi_{|\gbigg^{-1}(P)}$
of $\gbigv_{|\gbigg^{-1}(P)}$
in the situation of Proposition {\rm\ref{prop;06.5.13.100}}.
Therefore, we obtain the perfect obstruction
theory of $\gbigf^{-1}(P)$.
\hfill\qed
\end{cor}

We are mainly interested in the cap product
of some cohomology classes and
the virtual fundamental classes.
\begin{cor}
Let $\Phi$ be a cohomology class on $\nbigm^s(\yhat)$,
and let $\Psi$ be a cohomology class on
$\nbigm^s(\vecyhat,\alpha_{\ast})$.
Assume $\deg(\Phi)=2\dim^f\nbigm^s(\yhat)$.
Proposition {\rm\ref{prop;06.5.13.100}}
implies the following relation
for any $k$-valued point $P$ of $\nbigm^s(\yhat)$:
\begin{equation}
 \label{eq;06.5.11.101}
\int_{\nbigm^s(\vecyhat,\alpha_{\ast})}
 \Phi\cdot\Psi=
\int_{\gbigf^{-1}(P)} \Psi
\times
\int_{\nbigm^s(\yhat)} \Phi
\end{equation}
(See the subsection {\rm\ref{subsection;06.6.19.110}}
for the cohomology classes and the evaluation
considered in this paper.)
\hfill\qed
\end{cor}

Let $L$ be a line bundle on $X$,
and let $\delta$ be an element of $\nbigp^{\br}$.
We can discuss a similar relation
for $\nbigm^s(\yhat,[L],\delta)$ and
$\nbigm^s(\vecyhat,[L],\alpha_{\ast},\delta)$.
We assume that
any $\delta$-semistable $L$-Bradlow pair
is $\mu$-$\delta$-stable,
and that $\alpha_i$ are sufficiently close to $1$.
Then we have the morphism
$\gbigf_L:\nbigm^s(\vecyhat,[L],\alpha_{\ast},\delta)
\lrarr \nbigm^s(\yhat,[L],\delta)$.

\begin{prop}
 \label{prop;06.5.13.102}
There exists a $\nbigm^s(\yhat,[L],\delta)$-scheme
$\gbigb(\vecyhat,[L],\alpha_{\ast},\delta)$
with a vector bundle $\gbigv$ and its section $\psi$,
such that the following holds:
\begin{itemize}
\item
 The morphism
 $\gbigg_L:\gbigb(\vecyhat,[L],\alpha_{\ast},\delta)
 \lrarr \nbigm^s(\yhat,[L],\delta)$ is smooth.
\item
 $\nbigm^s(\vecyhat,[L],\alpha_{\ast},\delta)$
 is isomorphic to $\psi^{-1}(0)$.
\item
 We have the following relation:
\[
 \psi^!\gbigg_L^{\ast}\bigl[\nbigm^s(\yhat,[L],\delta)\bigr]
=\bigl[\nbigm^s(\vecyhat,[L],\alpha_{\ast},\delta)\bigr]
\]
Here $\psi$ denotes the Gysin map
for $\nbigm^s(\vecyhat,[L],\alpha_{\ast},\delta)\lrarr \gbigv$.
\end{itemize}

As a result, the following formula holds,
for any $k$-valued point $P$ 
of $\nbigm^s(\yhat,[L],\delta)$:
\[
 \int_{\nbigm^s(\vecyhat,[L],\alpha_{\ast},\delta)}\Phi\cdot \Psi
=\int_{\gbigg_L^{-1}(P)} \Psi
\cdot
\int_{\nbigm^s(\yhat,[L],\delta)} \Phi
\]
Here $\Phi$ and $\Psi$ denote
cohomology classes on $\nbigm^s(\yhat,L,\delta)$
and $\nbigm^s(\vecyhat,L,\alpha_{\ast},\delta)$
respectively,
and we assume $\deg(\Phi)=2\dim^f\nbigm^s(\yhat,[L],\delta)$.
\hfill\qed
\end{prop}

We will give the proof of Proposition \ref{prop;06.5.13.100}
in the next subsubsections.
The proof of Proposition \ref{prop;06.5.13.102} is similar,
and hence we omit to give it.

\subsubsection{The construction of a stack $\gbigb$
 and the obstruction theory}
\label{subsubsection;06.5.14.3}

Let $m$ be a sufficiently large integer.
Let $\nbige^u$ denote the universal sheaf
on $\nbigm(m,\yhat)$.
We put $\nbigv_0:=p_X^{\ast}p_{X\,\ast}\nbige^u(m)$,
and the kernel of the natural morphism
$\nbigv_0\lrarr \nbige^u(m)$ is denoted by
$\nbigv_{-1}$.
We obtain the vector bundle
$\nbigv_{0\,|\,D}$ on $\nbigm(m,\yhat)\times D$.

Let $g:T\lrarr \nbigm(m,\yhat)$ be a morphism.
Let $F(T)$ denote the set of the sequence of the quotients
$g_{D}^{\ast}\nbigv_{0\,|\,D}=\nbigc_{l+1}\rarr
 \nbigc_{l}\rarr \nbigc_{l-1}\rarr \cdots \rarr \nbigc_2\rarr \nbigc_1$
satisfying the following conditions:
\begin{itemize}
\item
 $\nbigc_i$ are flat over $T$.
\item
 For any point $u\in T$,
 the induced morphisms
 $H^0\bigl(D,\nbigc_{l+1\,|\,\{u\}\times D}\bigr)
\lrarr H^0\bigl(D,\nbigc_{i\,|\,\{u\}\times D}\bigr)$
 are surjective $(i=1,\ldots,l)$.
\item
 $H^1\bigl(D,\nbigc_{i\,|\,\{u\}\times D}\bigr)=0$
 for any $u\in T$ and for any $i=1,\ldots,l$.
\item
 The type of $C_i$ is same as $\sum_{j\leq i}y_j(m)$.
\end{itemize}
Then, we obtain the functor $F$
of the category of $\nbigm(m,\yhat)$-schemes
to the category of sets.
The functor is representable by
the $\nbigm(m,\yhat)$-scheme,
which we denote by $\gbigb$.
Let $\pi:\gbigb\lrarr \nbigm(m,\yhat)$ denote 
the natural projection.

Let us discuss the obstruction theory of $\gbigb$.
We put $\nbigv_D^{(1)}:=\pi_D^{\ast}\nbigv_{0|D}$
on $\gbigb\times D$.
We have the universal quotients
$\nbigv_D^{(1)}\lrarr \nbigc_{i}^u$ $(i=1,\ldots,l)$.
We put
$\nbigv_D^{(i)}:=\Ker\bigl(\nbigv_D^{(1)}\lrarr \nbigc^u_{i-1}\bigr)$
for $i=2,\ldots,l+1$.
We also have the locally free sheaf
$\nbigv_{-1\,|D}$.
By changing slightly the construction
in the subsubsection \ref{subsubsection;06.5.23.10},
we consider the following complex 
$\Cbar(\nbigv_{\cdot|D},\nbigv_D^{\ast})$
\[
 \nhom(\nbigv_{0|D},\nbigv_{-1|D})
 \stackrel{d^{-1}}{\lrarr}
 \bigoplus_{i=1}^{l+1}
 \nhom(\nbigv_D^{(i)},\nbigv_D^{(i)})
\oplus
 \nhom(\nbigv_{-1|D},\nbigv_{-1|D})
\stackrel{d^0}{\lrarr}
 \bigoplus_{i=1}^l\nhom(\nbigv_D^{(i+1)},\nbigv_D^{(i)})
\oplus
 \nhom(V_{-1|D},V_{0|D})
\]
The first term stands in the degree $0$.
The morphism $d^{-1}$ is the composite
of the following morphisms:
\[
 \nhom(\nbigv_{0|D},\nbigv_{-1|D})
 \stackrel{a_1}{\lrarr}
 \nhom(\nbigv_{0|D},\nbigv_{0|D})
\oplus
 \nhom(\nbigv_{-1|D},\nbigv_{-1|D})
 \stackrel{a_2}{\lrarr}
 \bigoplus_{i=1}^{l+1}
 \nhom(\nbigv_D^{(i)},\nbigv_D^{(i)})
\oplus
 \nhom(\nbigv_{-1|D},\nbigv_{-1|D})
\]
Here $a_1$ is the differential of the complex
$\nhom(\nbigv_{\cdot|D},\nbigv_{\cdot|D})$,
and $a_2$ is the inclusion
via $\nbigv_{0|D}=\nbigv_D^{(1)}$.
The morphism $d^0$ is made
of the following maps $b_i$ $(i=1,2)$:
\[
 b_1:
 \bigoplus_{i=1}^{l+1}
 \nhom(\nbigv_D^{(i)},\nbigv_D^{(i)})
\lrarr \bigoplus_{i=1}^l
 \nhom(\nbigv_D^{(i+1)},\nbigv_D^{(i)})
\]
\[
 b_2:
 \nhom(\nbigv_D^{(1)},\nbigv_D^{(1)})
\oplus \nhom(\nbigv_{-1|D},\nbigv_{-1|D})
\lrarr
 \nhom(\nbigv_{-1|D},\nbigv_{0|D})
\]
Here $b_1$ is given as in (\ref{eq;06.6.14.2}),
and $b_2$ is the differential of the complex
$\nhom(\nbigv_{\cdot|D},\nbigv_{\cdot|D})$.
We put 
$\gminigbar_D(\nbigv_{\cdot},\nbigv_D^{\ast})
 :=\Cbar(\nbigv_{\cdot|D},\nbigv_D^{\ast})^{\lor}[-1]$.
We have the naturally defined morphism
$\Cbar(\nbigv_{\cdot|D},\nbigv_D^{\ast})
\lrarr \nhom(\nbigv_{\cdot|D},\nbigv_{\cdot|D})$.
It induces the morphism
$\gminig(\nbigv_{\cdot|D})\lrarr 
 \gminigbar_D(\nbigv_{\cdot},\nbigv^{\ast}_D)$

We take vector spaces $W^{(i)}$ over $k$
such that $\rank W^{(i)}=\rank \nbigv_D^{(i)}$.
We also take a vector space $W_{-1}$ over $k$
such that $\rank W_{-1}=\rank \nbigv_{-1}$,
and we put $W_0:=W^{(1)}$.
We put $W_D^{(i)}:=W^{(i)}\otimes\nbigo_D$
and $W_{i\,D}:=W_i\otimes\nbigo_D$.
We have the naturally defined right 
$\prod_{i=1}^{l+1}\GL(W^{(i)})$-action
on $\prod_{i=1}^{l}N\bigl(W_D^{(i+1)},W_D^{(i)}\bigr)$.
We also have natural right action
of $\GL(W^{(1)})\times\GL(W_{-1})$
on $N(W_{-1\,D},W_{0\,D})$
by the identification $W_{0}=W^{(1)}$.
Therefore,
we have the naturally defined right action of
$\prod_{i=1}^{l+1}\GL(W^{(i)})\times\GL(W_{-1})$
on $\prod_{i=1}^{l}N\bigl(W_D^{(i+1)},W_D^{(i)}\bigr)
\times N\bigl(W_{-1\,D},W_{0\,D}\bigr)$,
where the latter fiber product is taken over $D$.
The quotient stack is denoted by $\Ybar_D(W_{\cdot},W^{\ast})$.
(We remark that we used the notation $\Ybar_D(W_{\cdot},W^{\ast})$
in a different meaning in the subsection \ref{subsection;06.6.1.10}.)
On the other hand,
we use the stack $Y_D(W_{\cdot})$
introduced in the subsubsection
\ref{subsubsection;06.5.5.5}.
The morphism
$\Ybar_D(W_{\cdot},W^{\ast})\lrarr
 Y_D(W_{\cdot})$
is induced by
the natural projections
$\prod_{i=1}^{l}N(W_D^{(i+1)},W_D^{(i)})
\times N(W_{-1\,D},W_{0\,D})$
and $\prod_{i=1}^{l+1}\GL(W^{(i)})\times \GL(W_{-1})$
onto $N(W_{-1\,D},W_{0\,D})$
and $\GL(W_0)\times \GL(W_{-1})$,
respectively.

From $\nbigv_D^{(i+1)}$ and $\nbigv_{-1|D}$,
we have the classifying map
$\Phi(\nbigv_{\cdot},\nbigv_D^{\ast}):
 \gbigb\times D\lrarr \Ybar_D(W_{\cdot},W^{\ast})$.
We also have the classifying map
$\Phi(\nbigv_{\cdot|D}):\nbigm(m,\yhat)\times D\lrarr Y_D(W_{\cdot})$.
They give the following commutative diagram:
\begin{equation}
\label{eq;06.5.25.100}
 \begin{CD}
 \gbigb\times D @>>> 
 \Ybar_D(W_{\cdot},W^{\ast})\\
 @VVV @VVV \\
 \nbigm(m,\yhat)\times D @>>>
 Y_D(W_{\cdot})
 \end{CD}
\end{equation}
It can be shown that
$\Phibar(\nbigv_{\cdot},\nbigv_D^{\ast})^{\ast}
 L_{\Ybar_D(W_{\cdot},W^{\ast})}$
is expressed by $\gminigbar_D(\nbigv_{\cdot},\nbigv_D^{\ast})_{\leq 1}$.
Moreover,
the diagram (\ref{eq;06.5.25.100}) 
induces the following commutative diagram:
\[
\begin{CD}
 L_{\gbigb\times D/D} @<<< 
 \Phi(\nbigv_{\cdot},
 \nbigv_D^{\ast})^{\ast}L_{\Ybar_D(W_{\cdot},W^{\ast})}
 @<<< \gminigbar_D(\nbigv,\nbigv_D^{\ast})\\
 @AAA @AAA @AAA\\
 \pi_D^{\ast}
 L_{\nbigm(m,\yhat)\times D/D} @<<<
 \Phi(\nbigv_{\cdot|D})^{\ast}L_{Y_D(W_{\cdot})}
 @<<<\gminig(V_{\cdot})
\end{CD}
\]

We put 
$\gminig_{\rel}(\nbigv_{\cdot},\nbigv_{D}^{\ast}):=
 \cone\bigl(\gminig(\nbigv_{\cdot|D})
 \lrarr\gminigbar_D(\nbigv_{\cdot},\nbigv_D^{\ast})\bigr)$.
Then, we put as follows:
\[
\Ob(\nbigv_{\cdot|D}):=
 Rp_{D\,\ast}\bigl(
\gminig(\nbigv_{\cdot|D})\otimes\omega_D\bigr),
\quad
 \Ob_{\rel}(\gbigb):=
 Rp_{D\,\ast}\bigl(
 \gminig_{\rel}(\nbigv_{\cdot},\nbigv_D^{\ast})\otimes\omega_D
\bigr)
\]
Then, we obtain the following commutative diagram:
\[
 \begin{CD}
 \pi^{\ast}L_{\nbigm(m,\yhat)}@<<<
 \Ob(\nbigv_{\cdot|D})\\
 @AAA @AAA \\
 L_{\gbigb/\nbigm(m,\yhat)}[-1]  @<<<
 \Ob_{\rel}(\gbigb)[-1]
 \end{CD}
\]
Therefore,
we obtain the following commutative diagram:
\[
 \begin{CD}
 \Ob_{\rel}(\gbigb)[-1] @>>> 
 \pi^{\ast}\Ob(m,\yhat)\\
 @VVV @VVV \\
 L_{\gbigb/\nbigm(m,\yhat)}[-1]
 @>>> \pi^{\ast}L_{\nbigm(m,\yhat)}
 \end{CD}
\]
We put
$\Ob(\gbigb):=\Cone\bigl(\Ob_{\rel}(\gbigb)\lrarr 
 \pi^{\ast}\Ob(m,\yhat)\bigr)$.
Then we obtain the morphism:
$\Ob(\gbigb)\lrarr L_{\gbigb}$.

\begin{prop}
\label{prop;06.5.13.120}
The morphism $\Ob_{\rel}(\gbigb)\lrarr L_{\gbigb/\nbigm(m,\yhat)}$
gives a relative obstruction theory of $\gbigb$
over $\nbigm(m,\yhat)$.
The complex $\Ob_{\rel}(\gbigb)$ is quasi-isomorphic to
the $0$-th cohomology sheaf.
\end{prop}
\pf
The first claim follows from Lemma \ref{lem;06.5.3.15}.
Let us show the second claim.
Due to an argument in the proof of Proposition
\ref{prop;06.5.5.40},
we have only to check
$H^i\bigl(i_z^{\ast}\Ob_{\rel}(\gbigb)\bigr)=0$
for $i\neq 0$,
where $i_z$ denotes the inclusion
of any point $z$ into $\gbigb$.
Let $(E,V_D^{\ast})$ denote the tuple corresponding to $z$.
Then, $H^i\bigl(i_z^{\ast}\Ob_{\rel}(\gbigb)\bigr)$
is the dual of the hyper-cohomology group
$\hyperh^{-i}(D,Q)$,
where $Q$ is the complex:
\[
 \bigoplus_{i=1}^l\nhom(V_D^{(i+1)},V_D^{(i+1)})
\lrarr
 \bigoplus_{i=1}^{l}\nhom(V_D^{(i+1)},V_D^{(i)})
\]
Here, the first term stands in the degree $-1$.
We use the following lemma.

\begin{lem}
\label{lem;06.5.13.121}
We have the vanishing
$H^1\Bigl(D,
 \nhom\bigl(V_D^{(i+1)},V_D^{(i)}/V_D^{(i+1)}\bigr)
 \Bigr)=0$.
\end{lem}
\pf
By definition of $\gbigb$,
we have the vanishing
$H^1\bigl(D,V_D^{(1)}/V_D^{(i)}\bigr)=0$
for any $i=1,\ldots,l+1$.
We have the following exact sequence:
\[
 H^0\bigl(D,V_D^{(1)}/V_D^{(i+1)}\bigr)
\stackrel{\nu_1}
{\lrarr}
 H^0\bigl(D,V_D^{(1)}/V_D^{(i)}\bigr)
\lrarr
 H^1\bigl(D,V_D^{(i)}/V_D^{(i+1)}\bigr)
\lrarr
 H^1\bigl(D,V_D^{(1)}/V_D^{(i+1)}\bigr)=0
\]
By definition of $\gbigd$,
we have the surjectivity of $\nu_1$.
Then we obtain $H^1\bigl(D,V_D^{(i)}/V_D^{(i+1)}\bigr)=0$.
From the exact sequence,
$V_D^{(i+1)}\lrarr V_D^{(1)}\lrarr V_D^{(1)}/V_D^{(i+1)}$,
we have the following exact sequence:
\begin{equation}
 \label{eq;06.6.23.3}
 \Ext^1\bigl(V_D^{(1)}, V_D^{(i)}/V_D^{(i+1)}\bigr)
\lrarr
 \Ext^1\bigl(V_D^{(i+1)}, V_D^{(i)}/V_D^{(i+1)}\bigr)
\lrarr
 \Ext^2\bigl(V_D^{(1)}/V_D^{(i+1)},V_D^{(i)}/V_D^{(i+1)}\bigr)
\end{equation}
Recall that $V_D^{(1)}$ is a direct sum of
$\nbigo_D$,
and hence we have the vanishing of the first term
in (\ref{eq;06.6.23.3}).
Since the divisor $D$ is smooth,
we have the vanishing of the third term
in (\ref{eq;06.6.23.3}).
Therefore,
we obtain the desired vanishing.
\hfill\qed

\vspace{.1in}
From Lemma \ref{lem;06.5.13.121},
we can easily obtain the vanishing
of $\hyperh^i\bigl(D,Q\bigr)$
unless $i=0$.
Thus the proof of Proposition \ref{prop;06.5.13.120}
is finished.
\hfill\qed

\vspace{.1in}

As a result,
$\Ob(\gbigb)\lrarr L_{\gbigb}$ is an obstruction
theory of $\gbigb$,
and the morphism $\gbigb\lrarr \nbigm(m,\yhat)$ is smooth.

\subsubsection{Compatibility of the obstruction theories
 of $\gbigb$ and $\nbigm(m,\vecyhat)$}
\label{subsubsection;06.5.14.4}

On $\gbigb\times D$,
we have the filtration
$\nbigv^{(l+1)}\subset\nbigv^{(l)}
\subset\cdots\subset\nbigv^{(1)}$.
We put as follows:
\[
 \gbigv:=p_{D\,\ast}\nhom\bigl(
 \nbigv^{(l+1)},\pi_D^{\ast}\nbige^u(m)_{|D}\bigr)
\]
We have the canonical section $\psi$,
which is given by the composite
$\nbigv^{(l+1)}\subset\nbigv^{(1)}\lrarr\pi_D^{\ast}\nbige^u(m)_{|D}$.
\begin{lem}
\label{lem;06.5.14.2}
$\gbigv$ is a locally free sheaf on $\gbigb$.
\end{lem}
\pf
Let $z=(E,V_D^{\ast})$ be any point of $\gbigb$.
We have only to check
$\Ext^1\bigl(V_D^{(l+1)},E(m)_{|D}\bigr)=0$.
We have the following exact sequence:
\[
 \Ext^1\bigl( V_D^{(1)},E(m)_{|D} \bigr)
\lrarr
 \Ext^1\bigl(V_D^{(l+1)},E(m)_{|D}\bigr)
\lrarr
 \Ext^2\bigl(V_D^{(1)}/V_D^{(l+1)},E(m)_{|D}\bigr)
\]
Since $V_D^{(1)}$ is a direct sum of $\nbigo_D$,
the first term vanishes.
Since $D$ is a smooth curve, the last term vanishes.
Therefore, we obtain the desired vanishing.
\hfill\qed

\vspace{.1in}

It is easy to observe 
$\psi^{-1}(0)=\nbigm(m,\vecyhat)$.
Therefore, we have the following Cartesian diagram:
\begin{equation}
\label{eq;06.5.25.40}
 \begin{CD}
 \nbigm(m,\vecyhat) @>{i_1}>> \gbigb\\
 @V{j}VV @V{\psi}VV\\
 \gbigb @>{i}>> \gbigv
 \end{CD}
\end{equation}
Here $i$ is the $0$-section.

Let us compare the obstruction theories of 
$\gbigb$ and $\nbigm(m,\vecyhat)$.
We take an isomorphism
$\gbigi:W^{(l+1)}\simeq W_{-1}$.
It induces the morphism
$Y_D(W_{\cdot},W^{\ast})\lrarr
 \Ybar_D(W_{\cdot},W^{\ast})$.
Then, we obtain the following commutative diagram:
\[
 \begin{CD}
 \nbigm(m,\vecy) @>>>
 \gbigb @>>>
 \nbigm(m,y)\\
 @VVV @VVV @VVV \\
 Y_D(W_{\cdot},W^{\ast})
 @>>>
 \Ybar_D(W_{\cdot},W^{\ast})
 @>>>
 Y_D(W_{\cdot})
 \end{CD}
\]
It induces the following diagram
on $\nbigm(m,\vecy)$:
\[
 \begin{CD}
 L_{\nbigm(m,\vecyhat)} @<<<
 L_{\gbigb}@<<<
 L_{\nbigm(m,\yhat)}\\
 @AAA @AAA @AAA \\
 L_{Y_D(W_{\cdot},W^{\ast})}
 @<<<
 L_{\Ybar_D(W_{\cdot},W^{\ast})}
 @<<<
 L_{Y_D(W_{\cdot})}\\
 @AAA @AAA @AAA \\
 \gminig(\nbigv_{\cdot|D},\nbigv_D^{\ast})
 @<<<
 \gminigbar(\nbigv_{\cdot|D},\nbigv_D^{\ast})
 @<<<
 \gminig(\nbigv_{\cdot|D})
 \end{CD}
\]
Therefore, we obtain the following commutative diagram
on $\nbigm(m,\vecy)$:
\begin{equation}
\label{eq;06.5.25.111}
\begin{CD}
 i_1^{\ast}\Ob(\gbigb)@>{\beta}>>\Ob(m,\vecyhat)\\
 @VVV @VVV \\
 i_1^{\ast}L_{\gbigb}@>>> L_{\nbigm(m,\vecyhat)}
\end{CD}
\end{equation}

\begin{lem}
\label{lem;06.5.25.201}
The cone of $\beta$ is isomorphic to
$j^{\ast}L_{\gbigb/\gbigv}$,
and the morphism
$\cone(\beta)\lrarr L_{\nbigm(m,\vecyhat)/\gbigb}$
obtained from {\rm(\ref{eq;06.5.25.111})}
is same as the morphism
obtained from the diagram {\rm(\ref{eq;06.5.25.40})}.
In particular,
the obstruction theories of
$\nbigm(m,\vecy)$ and $\gbigb$
are compatible over $i$.
\end{lem}
\pf
We put 
$\gminig_{\rel}:=
 \nhom\bigl(\nbigv_{-1\,|\,D}[1],
\nbigv_{\cdot|D}\bigr)^{\lor}$.
Then it is easy to see the following:
\[
 \cone\bigl(\gminigbar(\nbigv_{\cdot|D},\nbigv^{\ast})\lrarr
 \gminig(\nbigv_{\cdot|D},\nbigv^{\ast})\bigr)
\simeq
 L_{Y_D(W_{\cdot},W^{\ast})/\Ybar_D(W_{\cdot},W^{\ast})}
\simeq
\gminig_{\rel}.
\]
We put
$\Ob_{\rel}:=Rp_{D\,\ast}\bigl(
\gminig_{\rel}\otimes\omega_D\bigr)$.
We have the induced morphism
$a:\Ob_{\rel}\lrarr L_{\nbigm(m,\vecyhat)/\gbigb}$.
We have 
$\Ob_{\rel}\simeq
 \cone\bigl(i_1^{\ast}\Ob(\gbigb)\lrarr \Ob(m,\vecyhat)\bigr)$,
and hence we have the morphism
$b:\Ob_{\rel}\lrarr L_{\nbigm(m,\vecyhat)/\gbigb}$
obtained from the diagram (\ref{eq;06.5.25.111}).
It is easy to observe $a=b$.

We have the natural
$\GL(W_{-1})\times\GL(W^{(1)})\times
 \GL(W^{(l+1)})$-action
on $N(W_{D}^{(l+1)},W^{(1)}_D)
 \times_D N(W_{-1D},W^{(1)}_D)$.
The quotient stack is denoted by $\gbigy_2$.
The isomorphisms
$W_{-1}\simeq W^{(l+1)}$
and $W_0\simeq W^{(1)}$ induce
$Y_D(W_{\cdot})\lrarr \gbigy_2$.
Then, we have the following diagram:
\[
 \begin{CD}
 \nbigm(m,\vecy)\times D @>>>
 Y_D(W_{\cdot},W^{\ast})@>{\gamma_1}>>
 Y_D(W_{\cdot})\\
 @VVV @VVV @VVV \\
 \gbigb\times D @>>>
 \Ybar_D(W_{\cdot},W^{\ast})@>>>
 \gbigy_2
 \end{CD}
\]
The induced morphism
$\Phi(\nbigv_{\cdot|D},\nbigv_D^{\ast})^{\ast}
 \gamma_1^{\ast}
L_{Y_D(W_{\cdot})/\gbigy_2}
\lrarr
 \Phi(\nbigv_{\cdot|D},\nbigv_D^{\ast})^{\ast}
 L_{Y_D(W_{\cdot},W^{\ast})/\Ybar_D(W_{\cdot},W^{\ast})}$
is isomorphic.

We have the natural
$\GL(W_{-1})\times\GL(W_0)\times\GL(W^{(l+1)})$-action
on $N(W_{-1},W_0)\times N(W^{(l+1)},W_{-1})$.
The quotient stack is denoted by
$\gbigy_1$.
The isomorphism $\gbigi$ induces
the following map:
\begin{equation}
\label{eq;06.5.25.200}
 N(W_{-1},W_0)\lrarr N(W_{-1},W_0)\times N(W^{(l+1)},W_{-1}),
\quad
 f\longmapsto (f,\gbigi)
\end{equation}
We also have the homomorphism
$\GL(W_{-1})\lrarr \GL(W_{0})\times\GL(W^{(l+1)})$
induced by $\gbigi$.
The morphism (\ref{eq;06.5.25.200})
is equivariant with respect to the actions.
Therefore, we obtain the morphism
$Y_D(W_{\cdot})\lrarr \gbigy_1$.
It is easy to observe that the morphism
is an open immersion.
Hence, we have the following diagram:
\[
 \begin{CD}
 \nbigm(m,\vecy)\times D@>>>
 Y_D(W_{\cdot})
 @>{\rm open}>> 
 \gbigy_1 \\
 @VVV @VVV @VVV \\
 \gbigb\times D @>>> 
 \gbigy_2 @>{=}>> \gbigy_2
 \end{CD}
\]

\vspace{.1in}
We put $\nbigw_{i}:=N(W^{(l+1)},W_{i})$ $(i=-1,0)$.
We have the natural right 
$\GL(\nbigw_{-1})\times\GL(\nbigw_{0})$-action
on $N(\nbigw_{-1},\nbigw_0)\times N(k,\nbigw_{-1})$
and $N(\nbigw_{-1},\nbigw_0)\times N(k,\nbigw_0)$.
The quotient stacks are denoted by 
$\Zbar_1$ and $\Zbar_2$ respectively.
We have the naturally induced map $\gbigy_i\lrarr \Zbar_i$,
and we obtain the following commutative diagram:
\begin{equation}
 \label{eq;06.5.25.50}
\begin{CD}
 \nbigm(m,\vecy)\times D@>{a_1}>>
 \gbigy_1
 @>{a_2}>>
 \Zbar_1\\
 @VVV @VVV @VVV \\
 \gbigb @>>>
 \gbigy_2 @>>>
 \Zbar_2
\end{CD}
\end{equation}
Thus, we obtain the following:
\[
\begin{CD}
 a_1^{\ast}a_2^{\ast}
 L_{\Zbar_1/\Zbar_2}
 @>{\simeq}>>
 a_1^{\ast}L_{\gbigy_1/\gbigy_2}
 @>{\simeq}>>
 L_{Y_D(W_{\cdot},W^{\ast})/\Ybar_D(W_{\cdot},W^{\ast})}
\lrarr
 L_{\nbigm(m,\vecy)/\gbigb}
\end{CD}
\]

We would like to use the result in 
the latter part of the subsubsection
\ref{subsubsection;06.5.25.2}.
We have the resolution
of $\nhom\bigl(\nbigv^{(l+1)},\nbige^u(m)\bigr)$
given by
$\nhom\bigl(\nbigv^{(l+1)},\nbigv_{-1}\bigr)
\lrarr \nhom\bigl(\nbigv^{(l+1)},\nbigv_0\bigr)$
on $\gbigb\times D$.
We put as follows:
\[
V_{0}=\nhom\bigl(\nbigv^{(l+1)},\nbigv_{-1}\bigr),
\quad
V_1=\nhom\bigl(\nbigv^{(l+1)},\nbigv_0\bigr)
\]
We have the naturally defined map
$\phi:\nbigo\lrarr V_1$ on $\gbigb\times D$,
and the lift
$\phitilde:\nbigo\lrarr V_0$ on $\nbigm(m,\vecy)\times D$.
The section $\varphi$ of $\nhom\bigl(\nbigv^{(l+1)},\nbige^u(m)\bigr)$
is naturally induced by $\phi$.
We put $Z_1=N(\nbigo,V_0)$
and $Z_2=N(\nbigo,V_1)$.
Then, we have the following commutative diagram:
\begin{equation}
 \label{eq;06.5.25.51}
 \begin{CD}
 \nbigm(m,\vecy)\times D @>>>
 Z_1 @>>> \Zbar_1\\
 @VVV @VVV @VVV \\
 \gbigb\times D@>>>
 Z_2 @>>> \Zbar_2
 \end{CD}
\end{equation}
We have the coincidence of
the composite of the horizontal arrows in the diagrams
(\ref{eq;06.5.25.50}) and (\ref{eq;06.5.25.51}).
Therefore,
we obtain 
$\Ob_{\rel}=\Ob^H(V_{\cdot},\varphi)$,
and the induced morphism
$\Ob_{\rel}\lrarr L_{\nbigm(m,\vecy)/\gbigb}$
is same as the morphism obtained from 
the diagram {\rm(\ref{eq;06.5.25.40})},
due to Proposition \ref{prop;06.5.25.200}.
Thus the proof of Lemma \ref{lem;06.5.25.201}
is finished.
\hfill\qed

\subsubsection{Smooth ambient stack}

Let $V_m$ be an $H(m)$-dimensional 
vector space,
where $H$ denotes the Hilbert polynomial
associated to $y$.
For any $k$-scheme $T$,
let $F_1(T)$ denote the set of
the sequences of quotients
$p_T^{\ast}V_{m,D}=\nbigc_{l+1}\rarr
 \nbigc_{l}\rarr \nbigc_{l-1}\rarr \cdots \rarr \nbigc_2\rarr \nbigc_1$
satisfying the following conditions:
\begin{itemize}
\item
$\nbigc_i$ and $\nbigf$
are flat over $T$.
\item
 For any point $u\in T$,
 the induced morphisms
 $H^0\bigl(D,\nbigc_{l+1\,|\,\{u\}\times D}\bigr)
\lrarr H^0\bigl(D,\nbigc_{i\,|\,\{u\}\times D}\bigr)$
 are surjective $(i=1,\ldots,l)$.
 The induced morphism
 $H^0\bigl(D,\nbigc_{1\,|\,u}\bigr)\lrarr H^0\bigl(D,\nbigf_{|u}\bigr)$
 is also surjective.
\item
 We have
 $H^1\bigl(D,\nbigf_{1\,|\,u}\bigr)=0$
and 
 $H^1\bigl(D,\nbigc_{i\,|\,\{u\}\times D}\bigr)=0$
 for any $u\in T$ and for any $i=1,\ldots,l$.
\item
 The type of $\nbigc_i$ is same as $\sum_{j\leq i}y_j(m)$.
 The type of $\nbigf$ is same as $\sum_j y_j(m)$.
\end{itemize}
Let $F_2(T)$ denote the set of the quotients
$p_T^{\ast}V_{m,D}\lrarr \nbigf$
satisfying the following conditions:
\begin{itemize}
\item
 $\nbigf$ is flat over $T$.
\item
 For any point $u\in T$,
 the induced morphism
 $V_m\lrarr H^0(D,\nbigf_{|u})$ is surjective.
\item
 $H^1\bigl(D,\nbigf_{|u}\bigr)=0$ 
 for any $u\in T$.
\item
 The type of $\nbigf$ is $\sum_{j} y_j(m)$.
\end{itemize}

Then, we obtain the functors $F_i$ $(i=1,2)$
of the category of $k$-schemes
to the category of sets.
The functors $F_i$ are representable by the $k$-schemes,
which we denote by $\gbigb_i$.

We have the natural right $\GL(V_m)$-action
on $Q^{\circ}(m,\yhat)\times\gbigb_1$.
The quotient stack is isomorphic to $\gbigb$.
By considering the restriction to $D$,
we obtain the natural morphism
$Q^{\circ}(m,\yhat)\lrarr \gbigb_2$,
which is $\GL(V_m)$-equivariant.

Let $\Zhat_m$ be as in (\ref{eq;06.5.21.30}).
Then, we have the natural right $\GL(V_m)$-action
on $\Zhat_m\times \gbigb_1\times\gbigb_2$.
The quotient stack is denoted by $\gbigbtilde$.
We have the $\GL(V_m)$-equivariant morphism
$Q^{\circ}(m,\yhat)\times \gbigb_1
\lrarr \Zhat_m\times \gbigb_1\times\gbigb_2$,
which is immersion.
Therefore, we obtain the immersion
$\gbigb\lrarr \gbigbtilde$.

We have the universal filtration on
$\Zhat_m\times \gbigb_1\times\gbigb_2$:
\[
 V_D^{(l+1)}\subset V_D^{(l)}\subset
\cdots \subset V_D^{(2)}\subset V_D^{(1)}
=V_m\otimes\nbigo_{\Zhat_m\times\gbigb_1\times\gbigb_2}
\]
We also have the universal subsheaf
on $\Zhat\times\gbigb_1\times\gbigb_2$:
\[
 V_{-1\,D}\subset V_{0\,D}=
 V_{m}\otimes\nbigo_{\Zhat\times\gbigb_1\times\gbigb_2}
\]
The $\GL(V_m)$-action on $\Zhat_m\times\gbigb_1\times\gbigb_2$
is naturally lifted to the action on them.
The descents are denoted by $\nbigv_D^{(i)}$ and $\nbigv_{-1\,D}$.
We put as follows:
\[
 \gbigvtilde:=p_{D\,\ast}\Bigl(
\nhom\bigl(
\nbigv_D^{(l+1)},\,\,\nbigv_D^{(1)}/\nbigv_{-1\,D}
 \bigr)\Bigr)
\]
By the same argument
as the proof of Lemma \ref{lem;06.5.14.2},
it can be shown that
$\gbigvtilde$ is locally free.
We also have 
$\gbigvtilde_{|\gbigb}\simeq\gbigv$.
Therefore, we obtain the following diagram:
\[
 \begin{CD}
 \nbigm(m,\vecyhat) @>>>
 \gbigb @>>>
 \gbigbtilde\\
 @V{i_1}VV @V{i}VV @V{i_2}VV\\
 \gbigb @>{\psi}>> \gbigv @>>>\gbigvtilde
 \end{CD}
\]
Here $i$ and $i_2$ denote the $0$-section.

\subsubsection{Proof of Proposition \ref{prop;06.5.13.100}}

Let us finish the proof of Proposition \ref{prop;06.5.13.100}.
We take a sufficiently large integer $m$
such that the condition $O_m$ holds
for any $(E,F_{\ast},\rho)\in\nbigm^s(\vecyhat,\alpha_{\ast})$.
Then we have the open immersion
$\nbigm^s(\vecyhat,\alpha_{\ast})\lrarr
 \nbigm(m,\vecyhat)$
and $\nbigm^s(\yhat)\lrarr \nbigm(m,y)$.
We take the stack $\gbigb$ as in the subsubsection
\ref{subsubsection;06.5.14.3}.
We put
$\gbigb(\vecyhat,\alpha_{\ast}):=
 \gbigb\times_{\nbigm(m,\yhat)}\nbigm^s(\yhat)$.
Due to Proposition \ref{prop;06.5.13.120},
it is smooth over $\nbigm^s(\yhat)$.
The restrictions of $\gbigv$ and $\psi$
to $\gbigb(\vecyhat,\alpha_{\ast})$
are denoted by the same notation.
It is clear $\psi^{-1}(0)=\nbigm^s(\vecyhat,\alpha_{\ast})$.

We have the immersion 
$\gbigb(\vecyhat,\alpha_{\ast})
\lrarr \gbigbtilde$.
Since $\gbigb(\vecyhat,\alpha_{\ast})$ is Deligne-Mumford,
there exists an open neighbourhood
$\gbigbtilde(\vecyhat,\alpha_{\ast})$
of $\gbigb(\vecyhat,\alpha_{\ast})$
in $\gbigbtilde$,
which is Deligne-Mumford and smooth.
The restriction of $\gbigvtilde$ to
$\gbigbtilde(\vecyhat,\alpha_{\ast})$ are denoted by the same notation.
Then, we obtain the following diagram:
\[
 \begin{CD}
 \nbigm^s(\vecyhat,\alpha_{\ast}) @>>>
 \gbigb(\vecyhat,\alpha_{\ast}) @>>>
 \gbigbtilde(\vecyhat,\alpha_{\ast})\\
 @VVV @VV{i}V @VV{i_2}V \\
 \gbigb(\vecyhat,\alpha_{\ast}) @>>>
 \gbigv @>>> \gbigvtilde
 \end{CD}
\]
Due to Lemma \ref{lem;06.5.25.201},
the obstruction theories of 
$\nbigm^s(\vecyhat,\alpha_{\ast})$
and $\gbigb(\vecyhat,\alpha_{\ast})$
are compatible.
Therefore, we obtain the relation
$ i^!\bigl([\gbigb(\vecyhat,\alpha_{\ast})]\bigr)
=\bigl[\nbigm^s(\vecyhat,\alpha_{\ast})\bigr]$
due to Proposition \ref{prop;06.6.11.150}.
We also have the relation
$\gbigg^{!}[\nbigm^s(\yhat)]=[\gbigb(\vecyhat,\alpha_{\ast})]$.
We also remark that $i^!=\psi^!$.
Thus we obtain the relation (\ref{eq;06.5.14.5}).
\hfill\qed

\section{Invariants}
\label{section;06.6.30.5}
For simplicity, we assume that the ground field $k$
is the complex number field $\cnum$ in this section.
Let $H^{\ast}(A)$ and $H_{\ast}(A)$
denote the singular cohomology and homology 
groups of a topological space $A$ with $\rnum$-coefficient.
Let $X$ be a smooth projective surface
over $\cnum$,
and let $D$ be a smooth divisor of $X$.
We denote the Picard variety of $X$ by $\Pic$.

\subsection{Preliminary}
\label{subsection;06.6.19.110}

\subsubsection{The ring $\nbigr$}

Let 
$\Map_{f}\bigl(\seisuu_{\geq\,0}^2,H^{\ast}(X)\bigr)$
denote the set of the maps
$\varphi:\seisuu_{\geq\,0}^2\lrarr H^{\ast}(X)$
such that
$\bigl\{(n_1,n_2)\,\big|\,\varphi(n_1,n_2)\neq 0\bigr\}$
is finite.
We use the notation $\Map_f\bigl(\seisuu_{\geq\,0}^3,H^{\ast}(D)\bigr)$
in a similar meaning.
The sets $\Map_f\bigl(\seisuu_{\geq\,0}^2,H^{\ast}(X)\bigr)$
and $\Map_f\bigl(\seisuu_{\geq\,0}^3,H^{\ast}(D)\bigr)$
are naturally vector spaces over $\rnum$.
We use the notation $\Sym(V)$ to denote the symmetric product
of a vector space $V$.
Then, we put as follows:
\index{$\nbigr$}
\[
 \nbigr:=H^{\ast}(\Pic)\otimes \nbigr',
 \quad\quad
\nbigr':=
 \Sym \Bigl(\Map_f\bigl(\seisuu_{\geq\,0}^2,H^{\ast}(X)
 \bigr)\Bigr)
\otimes
 \Sym\Bigl(\Map_f\bigl(\seisuu_{\geq\,0}^3,H^{\ast}(D)\bigr)\Bigr)
\]
An element of $\nbigr$ is described as 
a sum of the elements of the following form:
\begin{equation}
\label{eq;06.6.17.5}
P=\gminic\cdot\prod_{i=1}^{m_1}(a_i,\vecv_i)\cdot 
 \prod_{j=1}^{m_2}(b_j,\vecu_j)
\end{equation}
\begin{itemize}
\item
$\gminic$ is an element of $H^{\ast}(\Pic)$.
\item
 $a_i\in H^{\ast}(X)$
and $\vecv_i=\bigl(v_{i}(1),v_i(2)\bigr)\in\seisuu_{\geq\,0}^2$.
We identify $(a_i,\vecv_i)$ with the map
$\varphi_i:\seisuu_{\geq\,0}^2\lrarr H^{\ast}(X)$:
\[
 \varphi_i(\vecv)=\left\{
 \begin{array}{ll}
 a_i & (\vecv=\vecv_i)\\
 0 & (\vecv\neq \vecv_i)
 \end{array}
 \right.
\]
\item
 $b_j\in H^{\ast}(D)$
and $\vecu_j=\bigl(u_j(1),u_j(2),u_j(3)\bigr)\in\seisuu_{\geq\,0}^3$.
 We identify $(b_j,\vecu_j)$ with the map
 $\psi_j:\seisuu_{\geq\,0}^3\lrarr H^{\ast}(D)$:
\[
 \psi_j(\vecu)=\left\{
 \begin{array}{ll}
 b_j & (\vecu=\vecu_j)\\
 0 & (\vecu\neq \vecu_j)
 \end{array}
 \right.
\]
\end{itemize}
We put
$d_1(P):=\sum_{i=1}^{m_1}v_i(1)\cdot v_i(2)
   +\sum_{j=1}^{m_2}u_j(1)\cdot u_j(2)-2m_1-m_2$.
When $a_i$, $b_j$ and $\gminic$ are homogeneous,
we put $d_2(P):=\sum_i\deg(a_i)/2+\sum_j \deg(b_j)/2+\deg(\gminic)/2$
and $d(P):=d_1(P)+d_2(P)$.

\subsubsection{The ring $\nbigr_l$}

More generally, we put as follows
for any $l\geq 1$:
\[
 \nbigr_l:=H^{\ast}(\Pic)\otimes\nbigr_l',\quad\quad
\nbigr_l':=
 \Sym \Bigl(\Map_f\bigl(
 \seisuu_{\geq\,0}^{2l},\,
 H^{\ast}(X)
 \bigr)\Bigr)
\otimes
 \Sym\Bigl(
 \Map_f\bigl(
 \seisuu_{\geq\,0}^{3l},\,
 H^{\ast}(D)\bigr)\Bigr)
\]
\index{$\nbigr_l$}

An element of $\nbigr_l$ is described as 
a sum of the elements of the following form
\begin{equation}
P=\gminic\cdot\prod_{i=1}^{m_1}(a_i,\vecV_i)\cdot 
 \prod_{j=1}^{m_2}(b_j,\vecU_j)
\end{equation}
\begin{itemize}
\item
$\gminic$ is an element of $H^{\ast}(\Pic)$.
\item
 $a_i\in H^{\ast}(X)$
and
  $\vecV_i=
 \bigl(v^{(h)}_{i}(1),v^{(h)}_i(2)\,\big|\,h=1,\ldots,l\bigr)
 \in \seisuu_{\geq\,0}^{2l}$.
We identify $(a_i,\vecV_i)$ with the map
$\varphi_i:\seisuu_{\geq\,0}^{2l}\lrarr H^{\ast}(X)$:
\[
 \varphi_i(\vecV)=\left\{
 \begin{array}{ll}
 a_i & (\vecV=\vecV_i)\\
 0 & (\vecV\neq \vecV_i)
 \end{array}
 \right.
\]
\item
 $b_j\in H^{\ast}(D)$
and $\vecU_j=
 \bigl(u^{(h)}_j(1),u^{(h)}_j(2),u^{(h)}_j(3)\bigr)
 \in\seisuu_{\geq\,0}^{3l}$.
 We identify $(b_j,\vecU_j)$ with the map
 $\psi_j:\seisuu_{\geq\,0}^{3l}\lrarr H^{\ast}(D)$:
\[
 \psi_j(\vecU)=\left\{
 \begin{array}{ll}
 b_j & (\vecU=\vecU_j)\\
 0 & (\vecU\neq \vecU_j)
 \end{array}
 \right.
\]
\end{itemize}
We put
$d_1(P):=\sum_h\sum_{i=1}^{m_1}v^{(h)}_i(1)\cdot v^{(h)}_i(2)
   +\sum_h\sum_{j=1}^{m_2}u^{(h)}_j(1)\cdot u^{(h)}_j(2)-2m_1-m_2$.
If $a_i$, $b_j$ and $\gminic$ are homogeneous,
we put $d_2(P):=\sum_i\deg(a_i)/2+\sum_j \deg(b_j)/2+\deg(\gminic)/2$
and $d(P):=d_1(P)+d_2(P)$.

\vspace{.1in}
Let $(a,\vecV)\in H^{\ast}(X)\times \seisuu_{\geq\,0}^{2l}$
be as above.
We regard $\vecV$ as a tuple
$(\vecv^{(1)},\ldots,\vecv^{(l)})\in (\seisuu_{\geq\,0}^2)^l$.
Let $\Delta_X^l$ denote the diagonal map
$X\lrarr X^l$,
and $\Delta_{X\,\ast}^l$ denotes the Gysin map
$H^{\ast}(X)\lrarr H^{\ast}(X^l)=H^{\ast}(X)^{\otimes\,l}$.
We have the expression as follows:
\[
 \Delta_{X\,\ast}^l (a)=
 \sum_h \prod_{i=1}^l\alpha_{i,h}
\]
Then, we put as follows:
\begin{equation}
 \gminiq'_l\bigl((a,\vecV)\bigr)
=\sum_h\bigl(\alpha_{1,h},\vecv^{(1)}\bigr)
 \otimes\cdots\otimes \bigl(\alpha_{l,h},\vecv^{(l)}\bigr)
\end{equation}
Let $(b,\vecU)\in H^{\ast}(D)\times \seisuu_{\geq\,0}^{3l}$
be as above.
We regard $\vecU$ as a tuple
$(\vecu^{(1)},\ldots,\vecu^{(l)})
 \in (\seisuu_{\geq\,0}^3)^l$.
Let $\Delta_D^l$ denote the diagonal map
$D\lrarr D^l$.
We have the expression as follows:
\begin{equation}
 \Delta_{D\,\ast}^l(b)=
 \sum_h \prod_{i=1}^l\beta_{i,h}
\end{equation}
Then, we put as follows:
\begin{equation}
 \gminiq'_l\bigl((b,\vecU)\bigr)
=\sum_h \bigl(\beta_{1,h},\vecu^{(1)}\bigr)
 \otimes\cdots\otimes \bigl(\beta_{l,h},\vecu^{(l)}\bigr)
\end{equation}
They induce the algebra homomorphism
$\gminiq'_l:\nbigr_l'\lrarr \nbigr_1^{\prime\,\otimes\,l}$.
We also have the morphism
$\gminiq_l:H^{\ast}(\Pic)\lrarr H^{\ast}(\Pic)^{\otimes\,l}$
induced by the multiplication of $\Pic$.
Then, we obtain the algebra homomorphism:
\begin{equation}
\label{eq;06.6.19.1}
 \gminiq_l:\nbigr_l\lrarr\nbigr^{\otimes\,l}
\end{equation}

We have the naturally defined algebra homomorphism
$\nbigr^{\prime\,\otimes\,l}\lrarr \nbigr'$.
Hence, $\gminiq'_l$ induces the algebra homomorphism
$\gminir_l:\nbigr_l\lrarr \nbigr$.

\subsubsection{Homomorphisms}
\label{subsubsection;06.6.19.25}

Let $\nbigy$ be an algebraic stack over $\Pic$.
When we are given a tuple of
parabolic sheaves
$\vecE_{\ast}=(E_{1\,\ast},\ldots,E_{l\,\ast})$
over $\nbigy\times (X,D)$,
we put $\nbigr(\vecE_{\ast}):=\nbigr_l$.
\index{$\nbigr(\vecE_{\ast})$, $\nbigr(\vecE)$}
In that case, $(a,\vecV)$ and $(b,\vecU)$
are symbolically denoted as follows:
\begin{equation}
 \label{eq;06.6.19.5}
 \Bigl(
 \prod_{h=1}^l
 \ch_{v^{(h)}(1)}^{v^{(h)}(2)}(E_h)
 \Bigr)\Big/a,
\quad
 \Bigl(
 \prod_{h=1}^l\ch_{u^{(h)}(1)}^{u^{(h)}(2)}
 \bigl(\Gr_{u^{(h)}(3)}(E_h)\bigr)
 \Bigr)\Big/b
\end{equation}
In particular,
we use the notation $\nbigr(E_{\ast})$
in the case $l=1$. \index{$\nbigr(E_{\ast})$, $\nbigr(E)$}
We will often omit to denote the parabolic structure
if there are no risk of confusion,
i.e.,
$\nbigr(\vecE)$ and $\nbigr(E)$ are used
instead of $\nbigr(\vecE_{\ast})$ and $\nbigr(E_{\ast})$.

When we are given a direct sum
$E_{\ast}=E_{1\,\ast}\oplus E_{2\,\ast}$,
we have the algebra homomorphism
$\varphi_{E_{\ast}}^{E_{1\ast},E_{2\ast}}:
\nbigr(E_{\ast})\lrarr \nbigr(E_{1\,\ast},E_{2\,\ast})$
induced by the following correspondence:
\[
 \varphi_{E_{\ast}}^{E_{1\ast},E_{2\ast}}
 \bigl(\ch_{i}^j(E)/a\bigr)
=\sum_{h=0}^j\frac{j!}{h!(j-h)!}
  \bigl(\ch_{i}^h(E_1)\cdot \ch_{i}^{j-h}(E_2)\bigr)\big/a,
\]
\[
  \varphi_{E_{\ast}}^{E_{1\ast},E_{2\ast}}
 \Bigl(\ch_{i}^j\bigl(\Gr_k(E)\bigr)\big/b\Bigr)
=\sum_{h=0}^j\frac{j!}{h!(j-h)!}
 \Bigl(\ch_{i}^h\bigl(\Gr_k(E_1)\bigr)\cdot
 \ch_i^{j-h}\bigl(\Gr_k(E_2)\bigr)\Bigr)\Big/b
\]
As the composite of
$\varphi_{E_{\ast}}^{E_{1\ast},E_{2\ast}}$
and $\gminiq_2$,
we obtain the algebra morphism
$\nbigr(E_{\ast})\lrarr \nbigr(E_{1\ast})\otimes\nbigr(E_{2\,\ast})$.
For an element $P(E)\in\nbigr(E_{\ast})$,
we denote the image by
$P(E_1\oplus E_2)\in \nbigr(E_{1\,\ast})\otimes\nbigr(E_{2\,\ast})$.

\vspace{.1in}

The algebra homomorphism $\gminir_l$
gives $\nbigr(\overbrace{E_{\ast},E_{\ast},\ldots,E_{\ast}}^l)
\lrarr \nbigr(E_{\ast})$,
which is also denoted by $\gminir_l$.
Let $t$ be a formal variable.
We have the element
$\ch_i^j(E\otimes e^t)/a\in \nbigr(E)[t]$
given as follows:
\begin{equation}
\label{eq;06.6.20.1}
\ch_i^j(E\otimes e^t)/a:=
 \sum_{\sum_{k=0}^i j_k=j}
 \frac{j!}{\prod_{k=0}^ij_k!}
 \gminir_{i+1}\left(
 \Bigl(
 \prod_{k=0}^i
 \ch_{k}^{j_k}(E)
\Bigr)\Big/a
 \right)\cdot \frac{t^{ji-\sum j_k\cdot k}}
 {\prod_{k=0}^i\bigl( (i-k)!\bigr)^{j_k}}
\end{equation}
Similarly, we have the element
$\ch_i^j(\Gr_h(E)\otimes e^t)/b\in \nbigr(E)[t]$
given as follows:
\[
\ch_i^j\bigl(\Gr_h(E)\otimes e^t\bigr)\big/b:=
 \sum_{\sum_{k=0}^i j_k=j}
 \frac{j!}{\prod_{k=0}^ij_k!}
 \gminir_{i+1}\left(
 \Bigl(
 \prod_{k=0}^i
 \ch_{k}^{j_k}\bigl(\Gr_h(E)\bigr)
\Bigr)\Big/b
 \right)\cdot \frac{t^{ji-\sum j_k\cdot k}}
 {\prod_{k=0}^i\bigl( (i-k)!\bigr)^{j_k}}
\]
By the correspondences
$\ch_i^j(E)/a\longmapsto \ch_i^j(E\cdot e^t)/a$
and 
$\ch_i^j(\Gr_h(E))/b\longmapsto
 \ch_i^j(\Gr_h(E)\cdot e^t)/b$,
we obtain the algebra isomorphism
$\nbigr(E)[t]\lrarr \nbigr(E)[t]$.
The image of $P(E)\in\nbigr(E)$ is denoted by
$P(E\cdot e^t)$.

\begin{rem}
The formula {\rm(\ref{eq;06.6.20.1})}
is just a formal development
of $\bigl(\sum_{h=0}^i\ch_{h}(E)\cdot (i-h)!^{-1}t^{i-h}\bigr)^j$.
\hfill\qed
\end{rem}

\subsubsection{Twist by line bundle}
\label{subsubsection;06.6.19.27}

Let $\nbigl$ be a line bundle on $\nbigy$.
We put $\omega:=c_1(\nbigl)$.
Formally, we often use the notation $e^{\omega}$
to denote $\nbigl$, 
if there are no risk of confusion.
Let $E_{\ast}$ be a parabolic sheaf over $\nbigy\times (X,D)$.
Then, we have the natural isomorphism
$\nbigr(E_{\ast},\nbigy)
\simeq
 \nbigr\bigl(E_{\ast}\otimes e^{\omega},\,\nbigy\bigr)$
(see Notation \ref{notation;06.6.23.10})
given by the following correspondence:
\begin{equation}
\label{eq;06.6.19.30}
\ch_i^j(E\otimes e^{\omega})/a:=
 \sum_{\sum_{k=0}^i j_k=j}
 \frac{j!}{\prod_{k=0}^ij_k!}
 \gminir_{i+1}\left(
 \Bigl(
 \prod_{k=0}^i
 \ch_{k}^{j_k}(E)
\Bigr)\Big/a
 \right)\cdot \frac{\omega^{ji-\sum j_k\cdot k}}
 {\prod_{k=0}^i\bigl( (i-k)!\bigr)^{j_k}}
\end{equation}
\begin{equation}
\label{eq;06.6.19.31}
\ch_i^j\bigl(\Gr_h(E)\otimes e^{\omega}\bigr)\big/b:=
 \sum_{\sum_{k=0}^i j_k=j}
 \frac{j!}{\prod_{k=0}^ij_k!}
 \gminir_{i+1}\left(
 \Bigl(
 \prod_{k=0}^i
 \ch_{k}^{j_k}\bigl(\Gr_h(E)\bigr)
\Bigr)\Big/b
 \right)\cdot \frac{\omega^{ji-\sum j_k\cdot k}}
 {\prod_{k=0}^i\bigl( (i-k)!\bigr)^{j_k}}
\end{equation}
Thus we can naturally regard
$P(E\otimes e^{\omega})\in\nbigr(E\otimes e^{\omega},\nbigy)$
as an element of $\nbigr(E,\nbigy)$.

\begin{rem}
We will often use ``$\cdot$'' instead of ``$\otimes$''
to save the space.
\hfill\qed
\end{rem}

Let us consider the case
$\nbigy=\nbigy_1\times\nbigy_2$.
Assume that $L$ comes from the line bundle on $\nbigy_1$,
and that $E$ comes from the parabolic sheaf on $\nbigy_2\times (X,D)$.
The formulas (\ref{eq;06.6.19.30})
and (\ref{eq;06.6.19.31}) determines the element 
$P(E\cdot e^{\omega})\in A^{\ast}(\nbigy_1)\otimes \nbigr(E)$.

\subsubsection{Evaluation}
\label{subsubsection;06.6.19.2}

Let $\nbigy$ be a proper Deligne-Mumford stack over $\Pic$.
Assume that we are given a parabolic sheaf $E_{\ast}$
over $\nbigy\times (X,D)$ with a parabolic structure.
Let $A_{\ast}(\nbigy)$ denote the rational Chow group.
Let us take an element $P\in\nbigr(E_{\ast})$
of the following form:
\[
 P=\gminic\cdot\prod_{i=1}^{m_1}\bigl(
 \ch_{v_i(1)}^{v_i(2)}(E)/a_i\bigr)
\cdot\prod_{j=1}^{m_2}\Bigl(
 \ch_{u_j(1)}^{u_j(2)}
 \bigl(\Gr_{u_j(3)}(E) \bigr)\big/b_j\Bigr)
\]
We assume that $\gminic\in H^{\ast}(\Pic)$, 
$a_i\in H^{\ast}(X)$ and $b_j\in H^{\ast}(D)$
are homogeneous for simplicity.
We would like to construct the linear morphism of
$A_{\ast}(\nbigy)$ to $\rnum$.
Let $\pi_{X,i}$ denote the projection of
$\nbigy\times X^{m_1}\times D^{m_2}$
onto the product of $\nbigy$ and the $i$-th $X$.
Let $\pi_{D,j}$ denote the projection of
$\nbigy\times X^{m_1}\times D^{m_2}$
onto the product of $\nbigy$ and $j$-th $D$.
Let $\gminip$ denote the natural morphism
$\nbigy\times X^{m_1}\times D^{m_2}$
to $\Pic\times X^{m_1}\times D^{m_2}$.
Let $\nbigz$ be a $d$-dimensional algebraic cycle on $\nbigy$.
Then, we obtain the following element of
the Chow group 
$A_{\ast}\bigl(\Pic\times X^{m_1}\times D^{m_2}\bigr)$
of $\Pic\times X^{m_1}\times D^{m_2}$
with rational coefficient:
\[
 \Lambda_P(E_{\ast},\nbigz):=
 \gminip_{\ast}\left(
 \prod_{i=1}^{m_1}
 \ch_{v_i(1)}\bigl(
 \pi_{X,i}^{\ast}E
 \bigr)^{v_i(2)}\cdot
 \prod_{j=1}^{m_2}
 \ch_{u_j(1)}\bigl(
 \pi_{D,j}^{\ast}\Gr_{u_j(3)}(E)
 \bigr)^{u_j(2)}
\cap [\nbigz\times X^{m_1}\times D^{m_2}]
 \right)
\]
Thus we obtain the linear map
$\Lambda_P(E_{\ast},\cdot):A_{d}(\nbigy)\lrarr
 A_{d-d_1(P)}(\Pic\times X^{m_1}\times D^{m_2})$.
The cycle $\Lambda_P(E_{\ast},\nbigz)$ determines
the homology class
$\cycl(\Lambda_P(E_{\ast},\nbigz))$
of $H_{2(d-d_1(P))}(\Pic\times X^{m_1}\times D^{m_2})$.
Let $\pi_{\ast}$ denote 
the push-forward for $\Pic\times X^{m_1}\times D^{m_2}$
to a point $\pt$.
Then, we obtain the following:
\[
  \deg\bigl(P(E_{\ast})\cap [\nbigz]\bigr):=
\pi_{\ast}\Bigl(
 \gminic\cdot
 \prod_{i=1}^{m_1}a_i\cdot
 \prod_{j=1}^{m_2}b_j\cap
 \cycl\bigl(\Lambda_P(E_{\ast},\nbigz)\bigr)
\Bigr)
 \in H_{2(d-d(P))}(\pt)
\]
It is trivial in the case $d\neq d(P)$.
We identify $H_{0}(\pt)\simeq \rnum$.
Thus, we obtain the linear map
$\deg\bigl(P(E_{\ast})\cap \cdot\bigr):
 A_{\ast}(\nbigy)\lrarr \rnum$.

Let $A^{\ast}(\nbigy)$ denote the operational Chow ring
of $\nbigy$.
Let $\nbigz$ be an algebraic cycle of $\nbigy$.
Let $F$ be an element of $A^{\ast}(\nbigy)$.
Then, we obtain the number
$\int_{\nbigz}P(E_{\ast})\cdot F:=
 \deg\Bigl(
 P(E_{\ast})\cap 
F\bigl([\nbigz]\bigr)
 \Bigr)$.

\begin{notation}
\label{notation;06.6.23.10}
\index{$\nbigr(E_{\ast},\nbigy)$}
Let $\nbigy$ and $E_{\ast}$ be as above.
We put
$\nbigr(E_{\ast},\nbigy):=
 \nbigr(E_{\ast})\otimes A^{\ast}(\nbigy)$.
We can naturally regard $\nbigr(E_{\ast},\nbigy)$ as
an $\bigl(\nbigr(E_{\ast}),\,A^{\ast}(\nbigy)\bigr)$-bimodule.
\hfill\qed
\end{notation}

We have the linear morphism
$\nbigr(E_{\ast},\nbigy)\otimes A_{\ast}(\nbigy)\lrarr 
 Hom(A_{\ast},\rnum)$  by the above construction.

\begin{rem}
Formally, $\deg\bigl(P(E_{\ast})\cap [\nbigz]\bigr)$
is the following number:
\begin{equation}
 \label{eq;06.6.17.1}
\int_{\nbigz}
 \gminic\cdot
 \prod_{i=1}^{m_1}
 \bigl(
 \ch_{v_i(1)}\bigl(E\bigr)^{v_i(2)}/a_i
 \bigr)
\cdot\prod_{j=1}^{m_2}
 \bigl(
 \ch_{u_j(1)}\bigl(\Gr_{u_j(3)}(E)\bigr)^{u_j(2)}\big/b_j
 \bigr)
\end{equation}
The author does now know an appropriate reference
for the cohomology and homology theories
of Deligne-Mumford stacks
with the good cycle maps from the Chow groups,
$G_m$-localization theory
and any other expected properties.
That is the reason to avoid {\rm(\ref{eq;06.6.17.1})}
as the definition.
However, it is easy to observe that
the formal argument using {\rm(\ref{eq;06.6.17.1})}
is valid.
\hfill\qed
\end{rem}

We are especially interested in the following examples.
\begin{example}{\rm
Let $\vecy$ be an element of $\Type$,
and let $\alpha_{\ast}$ be a system of weights.
We have the universal sheaf $\Ehat$ over
$\nbigm^{ss}(\vecyhat,\alpha_{\ast})\times X$
with the parabolic structure at
$\nbigm^{ss}(\vecyhat,\alpha_{\ast})\times D$.
Let $P$ be an element of  $\nbigr$.
By the identification $\nbigr(\Ehat_{\ast})=\nbigr$,
we have $\Phi=P(\Ehat^u)\in\nbigr(\Ehat_{\ast})$.
Assume that the $1$-stability condition holds
for $(\vecy,\alpha_{\ast})$.
We put as follows:
\[
 \int_{\nbigm^{ss}(\vecyhat,\alpha_{\ast})}
 \Phi
:=\deg\Bigl(
 \Phi\cap
 [\nbigm^{ss}(\vecyhat,\alpha_{\ast})]
 \Bigr)
\]
In other words,
we obtain the linear map
$\int_{\nbigm^{ss}(\vecyhat,\alpha_{\ast})}:
 \nbigr\lrarr \rnum$,
under the assumption that 
the $1$-stability condition holds
for $(\vecy,\alpha_{\ast})$.
We will later discuss how to obtain such a morphism
in general.
\hfill\qed
}
\end{example}

\begin{example}
Let $\vecy$ and $\alpha_{\ast}$ be as above.
Let $L$ be a line bundle on $X$,
and let $\delta$ be an element of $\nbigp^{\br}$
such that the $1$-stability condition holds for
$(\vecy,L,\alpha_{\ast},\delta)$.
We denote by $\omega$ the first Chern class of
the tautological line bundle $\nbigo_{\rel}(1)$
on $\nbigm^s(\vecyhat,[L],\alpha_{\ast},\delta)$.
For any $P\in\nbigr$,
we have the element
$\Phi=P(\Ehat^u)\cdot \omega^k\in
 \nbigr\bigl(\Ehat^u,\nbigm^s(\vecyhat,[L],\alpha_{\ast},\delta)\bigr)$.
Thus, we obtain the following number:
\[
 \int_{\nbigm^s(\vecyhat,[L],\alpha_{\ast},\delta)}
\Phi
 :=
 \deg\Bigl(
 \Phi\cap\bigl[\nbigm^s(\vecyhat,[L],\alpha_{\ast},\delta)\bigr]
 \Bigr)
\]
If the $1$-vanishing condition holds for
$(\vecy,L,\alpha_{\ast},\delta)$, moreover,
then we have the relative tangent bundle $T_{\rel}$
of the smooth map
$\nbigm^s(\vecyhat,[L],\alpha_{\ast},\delta)
 \lrarr \nbigm(\vecyhat,[L])$.
Let $\Eu(T_{\rel})$ denote the Euler class of $T_{\rel}$.
For any $P\in\nbigr$,
we have
$\Phi=P(\Ehat^u)\cdot\Eu(T_{\rel})\in
 A^{\ast}\bigl(\nbigm^s(\vecyhat,[L],\alpha_{\ast},\delta)\bigr)$.
Thus, we obtain the integral
$\int_{\nbigm^s(\vecyhat,[L],\alpha_{\ast},\delta)}\Phi$.
\hfill\qed
\end{example}

Let us consider the case 
where we are given a tuple of parabolic sheaves
$\vecE_{\ast}=(E_{1\,\ast},\ldots,E_{l\,\ast})$
on $\nbigy\times (X,D)$.
Let us take an element $P\in\nbigr_l(\vecE_{\ast})$
of the following form:
\[
 P=\gminic\cdot
 \prod_{i=1}^{m_1}\left(
 \prod_{h=1}^l
 \ch_{v^{(h)}_i(1)}^{v^{(h)}_i(2)}(E)/a_i
\right)
\cdot\prod_{j=1}^{m_2}\left(
 \prod_{h=1}^l
 \ch_{u_j(1)}^{u_j(2)}
 \bigl(\Gr_{u_j(3)}(E)\bigr)\big/b_j
\right)
\]
We assume that $\gminic\in H^{\ast}(\Pic)$, 
$a_i\in H^{\ast}(X)$ and $b_j\in H^{\ast}(D)$
are homogeneous for simplicity.
Let $\nbigz$ be a $d$-dimensional algebraic cycle on $\nbigy$.
Then, we obtain the following element of
the Chow group 
$A_{\ast}\bigl(\Pic\times X^{m_1}\times D^{m_2}\bigr)$
of $\Pic\times X^{m_1}\times D^{m_2}$
with rational coefficient:
\begin{multline}
 \Lambda_P(\vecE_{\ast},Z):= \\
 \gminip_{\ast}\left(
\prod_{h=1}^l\Bigl(
 \prod_{i=1}^{m_1}
 \ch_{v^{(h)}_i(1)}\bigl(
 \pi_{X,i}^{\ast}E_h
 \bigr)^{v^{(h)}_i(2)}\cdot
 \prod_{j=1}^{m_2}
 \ch_{u^{(h)}_j(1)}\bigl(
 \pi_{D,j}^{\ast}\Gr_{u^{(h)}_j(3)}(E_h)
 \bigr)^{u^{(h)}_j(2)}
 \Bigr)
\cap [\nbigz\times X^{m_1}\times D^{m_2}]
 \right)
\end{multline}
Thus we obtain the linear map
$\Lambda_P(E_{\ast},\cdot):A_{d}(\nbigy)\lrarr
 A_{d-d_1(P)}(\Pic\times X^{m_1}\times D^{m_2})$.
The cycle $\Lambda_P(E_{\ast},\nbigz)$ determines
the homology class
$\cycl(\Lambda_P(E_{\ast},\nbigz))$
of $H_{2(d-d_1(P))}(\Pic\times X^{m_1}\times D^{m_2})$.
Let $\pi_{\ast}$ denote 
the push-forward for $\Pic\times X^{m_1}\times D^{m_2}$
to a point $\pt$.
Then, we obtain the following:
\[
  \deg\bigl(P(\vecE_{\ast})\cap [\nbigz]\bigr):=
\pi_{\ast}\Bigl(
 \gminic\cdot
 \prod_{i=1}^{m_1}a_i\cdot
 \prod_{j=1}^{m_2}b_j\cap
 \cycl\bigl(\Lambda_P(E_{\ast},\nbigz)\bigr)
\Bigr)
 \in H_{2(d-d(P))}(\pt)
\]
We identify $H_{0}(\pt)=\rnum$.
Thus, we obtain the linear map
$\deg\bigl(P(\vecE_{\ast})\cap \cdot\bigr):
 A_{\ast}(\nbigy)\lrarr \rnum$.

\begin{notation}
\index{$\nbigr(\vecE_{\ast},\nbigy)$}
We put $\nbigr(\vecE_{\ast},\nbigy):=
 \nbigr(\vecE_{\ast})\otimes A^{\ast}(\nbigy)$.
We obtain the linear map
$\nbigr(\vecE_{\ast},\nbigy)\otimes A_{\ast}(\nbigy)
\lrarr \rnum$
by the above construction.
\hfill\qed
\end{notation}

We have the following commutative diagram,
which we will use implicitly.
\[
 \begin{CD}
 \nbigr(\overbrace{E_{\ast},\ldots,E_{\ast}}^l,\nbigy)
 @>>>
 \Hom(A_{\ast}(\nbigy),\rnum)\\
 @VV{\gminir_l}V @VV{=}V \\
 \nbigr(E_{\ast},\nbigy) @>>>
 \Hom(A_{\ast}(\nbigy),\rnum)
 \end{CD}
\]

\vspace{.1in}
Assume that we have the decomposition
$\nbigy=\prod_{i=1}^l \nbigy_i$
such that $E_{i\,\ast}$ are pull back of 
the parabolic sheaves over $\nbigy_i\times (X,D)$
via the natural projection,
where $\nbigy_i$ are the stacks over $\Pic$
and the map $\nbigy\lrarr \Pic$ is the composite of
$\prod \nbigy_i\lrarr \Pic^l$
and the multiplication of $\Pic$.
We have the naturally defined morphism:
\[
 \Gamma_1:
 \bigotimes_{i=1}^l\nbigr(E_{i\,\ast},\nbigy_i)
\lrarr
 \Hom_{\rnum}\left(
 \bigotimes_{i=1}^lA_{\ast}(\nbigy_i),\,
 \rnum\right)
\]
We also have the following morphism:
\[
 \Gamma_2:\nbigr(\vecE_{\ast})
 \otimes\bigotimes_{i=1}^l A^{\ast}(\nbigy_i)
 \lrarr 
  \Hom_{\rnum}\left(
 \bigotimes_{i=1}^lA_{\ast}(\nbigy_i),\,
 \rnum\right)
\]
The algebra homomorphism $\gminiq_l$
in (\ref{eq;06.6.19.1}) induces 
$\nbigr(\vecE_{\ast})\otimes\bigotimes_{i=1}^lA^{\ast}(\nbigy_i)
 \lrarr
\bigotimes_{i=1}^l\nbigr(E_{i\ast},\nbigy_i)$,
which is also denoted by $\gminiq_l$.
We will use the following lemma without mention,
which can be checked by a formal calculation.
\begin{lem}
\label{lem;06.6.18.110}
We have
$\Gamma_1\circ\gminiq_l=\Gamma_2$.
\hfill\qed
\end{lem}

\subsubsection{Equivariant case}

We continue to use the setting in the subsubsection
\ref{subsubsection;06.6.19.2}.
We recall the equivariant Chow groups
for the torus action.
See \cite{edidin-graham} and \cite{gp}
for more detail.
Let $T$ denote an $l$-dimensional torus $G_m^l$.
Let $R(T)=\rnum[t_1,\ldots,t_l]$ denote 
the representation ring of $T$.
Let us consider the case where
$\nbigy$ is provided with a $T$-action.
Let $A^{m+1}$ denote the $m+1$-dimensional
linear space,
on which $G_m$ acts by component-wise multiplication.
We have the naturally defined $T$-action on
$\bigl(A^{m+1}-\{0\}\bigr)^l\times \nbigy$.
The quotient stack is denoted by $\nbigy^{(m)}$.
The $T$-equivariant Chow group of $A^T_d(\nbigy)$
is defined to be $A_{d+lm}(\nbigy^{(m)})$
for sufficiently large $m$ in this case.
We take a linear inclusion $\iota:A^{m+1}\lrarr A^{m+2}$,
which induces the regular embedding
$\nbigy^{(m)}\lrarr \nbigy^{(m+1)}$.
Thus we obtain the morphism
$A_{d+l(m+1)}(\nbigy^{(m+1)})\lrarr A_{d+lm}(\nbigy^{(m)})$,
which is independent of the choice of $\iota$.
It is isomorphism when $m$ is sufficiently large.
Thus $A^T_d(\nbigy)$ is well defined.

We have the naturally defined morphism
$\nbigy^{(m)}\lrarr (\proj^m)^l$.
Let $\nbigo^{(i)}(1)$ denote the tautological line bundle
of $i$-th $\proj^{m}$.
We have the action of the first Chern class
$c_1\bigl(\nbigo^{(i)}(1)\bigr):
 A_{d+lm}(\nbigy^{(m)})\lrarr A_{d-1+lm}(\nbigy^{(m)})$,
which induces the action
$t_i:A_{d}^T(\nbigy)\lrarr A_{d-1}(\nbigy)$.
Thus, we can naturally regard 
$A^T_{\ast}(\nbigy)$ as the $R(T)$-module.

Assume that $E_{\ast}$ is provided with a $T$-action.
Let $P$ be an element of $\nbigr(E_{\ast})$.
We naturally obtain 
the parabolic sheaf $E^{(m)}_{\ast}$
on $\nbigy^{(m)}$.
Let $\nbigz$ be an element of $A_d^{T}(\nbigy)$,
and let $\nbigz^{(m)}$ be the corresponding element of
$A_{d+lm}(\nbigy^{(m)})$.
We obtain the following element of
$A^{T}_{d-d_1(P)+lm}\bigl(
 (\proj^m)^l\times\Pic\times X^{m_1}\times D^{m_2}\bigr)$:
\[
 \Lambda_P(E^{(m)}_{\ast},\nbigz^{(m)})\!\!:=\!\!
 \gminip_{\ast}\left(
 \prod_{i=1}^{m_1}
 \ch_{v_i(1)}\bigl(
 \pi_{X,i}^{\ast}E^{(m)}
 \bigr)^{v_i(2)}
 \prod_{j=1}^{m_2}
 \ch_{u_j(1)}\bigl(
 \pi_{D,j}^{\ast}\Gr_{u_j(3)}(E^{(m)})
 \bigr)^{u_j(2)}
\cap [\nbigz^{(m)}\times X^{m_1}\times D^{m_2}]
 \right)
\]
Then, it is easy to observe that
$\Lambda_P(E^{(m)}_{\ast},\nbigz^{(m)})$
determines the element
$\Lambda_P^T(E_{\ast},\nbigz)$
of $A_{d-d_1(P)}^T(\Pic\times X^{m_1}\times D^{m_2})$.
Thus, we obtain the $R(T)$-morphism
$\Lambda^T_P(E_{\ast},\cdot):
A^T_{\ast}(\nbigy)\lrarr 
A^T_{\ast-d_1(P)}(\Pic\times X^{m_1}\times D^{m_2})$.
Then, we obtain the equivariant homology class
$\cycl(\Lambda^T_P(E_{\ast},\cdot))
\in H^T_{2(\ast-d_1(P))}(\Pic\times X^{m_1}\times D^{m_2})$.
Since the $T$-action on $\Pic\times X^{m_1}\times D^{m_2}$ are trivial,
$\gminic$, $a_i$ and $b_j$ naturally give the equivariant cohomology class.
Let $\pi^T_{\ast}$ denote the equivariant Gysin map
for $\Pic\times X^{m_1}\times D^{m_2}\lrarr \pt$.
Then, we obtain the following:
\[
 \deg^T\bigl(P(E_{\ast})\cap [\nbigz]\bigr):=
\pi_{\ast}\Bigl(
 \gminic\cdot
 \prod_{i=1}^{m_1}a_i\cdot
 \prod_{j=1}^{m_2}b_j\cap
 \cycl\bigl(\Lambda^T_P(E,\nbigz)\bigr)
\Bigr)
 \in H_{2(d-d(P))}(\pt)
\]
Since $H^{T}_{\ast}(\pt)$ is the free $R(T)$-module
of rank one with the special base $1\in H_0^T(\pt)$,
we identify it with $R(T)$.
The above procedure is compatible with
the actions of $t_i$.
Hence we obtain the morphism
$\deg^T\bigl(P(E_{\ast})\cap \cdot \bigr):
 A^T_{\ast}(\nbigy)\lrarr R(T)$
of $R(T)$-modules.

\begin{notation}
\index{$\nbigr_T(E_{\ast},\nbigy)$}
Let $A_T^{\ast}(\nbigy)$ denote the $T$-equivariant
operational Chow group of $\nbigy$.
We put $\nbigr_T(E_{\ast},\nbigy):=
 \nbigr(E_{\ast})\otimes_{\rnum} A^{\ast}_T(\nbigy)$.
\hfill\qed
\end{notation}
Due to the above construction,
we have the morphism
$\nbigr_T(E_{\ast},\nbigy)\lrarr
\Hom_{R(T)}\bigl(A^T_{\ast}(\nbigy),R(T)\bigr)$.

\vspace{.1in}

\vspace{.1in}

We are especially interested in the following example:
\begin{example}{\rm
\label{example;06.6.17.5}
Let $\Mhat$ be the master space as in
the subsubsections {\rm\ref{subsubsection;06.5.17.10}},
{\rm\ref{subsubsection;06.5.17.100}}
or {\rm\ref{subsubsection;06.5.21.555}}.
We have the sheaf $\Ehat^{\Mhat}$
on $\Mhat\times X$ with the parabolic structure 
at $\Mhat\times D$.
We have the $G_m$-action $\rhobar$
on $\Mhat$, which is naturally lifted 
to the action on $\Ehat^{\Mhat}$.
Let $\Phi(\Ehat^{\Mhat}_{\ast})$
be an element of $\nbigr_T(\Ehat,\Mhat)$.
As explained in the subsection {\rm\ref{subsection;06.6.17.2}},
the perfect obstruction theory of $\Mhat$
is lifted to the $G_m$-equivariant obstruction theory.
Thus, we obtain $[\Mhat]\in A^{G_m}_{\ast}(\Mhat)$.
Then, we obtain the following number:
\[
 \int_{\Mhat}\Phi(\Ehat_{\ast}^{\Mhat})
:=\deg^{G_m}\Bigl(
 \Phi(\Ehat_{\ast}^{\Mhat})\cap [\Mhat]
 \Bigr)
\]
We have the isomorphism:
\[
 A^{G_m}_{\ast}(\Mhat)\otimes_{\rnum[t]} \rnum[t,t^{-1}]
\simeq
\bigoplus_{i=1,2}
 A_{\ast}(\Mhat_i)\otimes_{\rnum}\rnum[t,t^{-1}]
\oplus
 \bigoplus_{\gbigi}
 A_{\ast}\bigl(\Mhat^{G_m}(\gbigi)\bigr)
 \otimes_{\rnum}\rnum[t,t^{-1}]
\]
Under the isomorphism,
we obtain the decomposition in
$A^{G_m}_{\ast}(\Mhat)\otimes_{\rnum[t]} \rnum[t,t^{-1}]$
due to {\rm\cite{gp}}:
\[
 [\Mhat]=
 \sum_{i=1,2}
 \Eu\bigl(\gbign(\Mhat_i)\bigr)^{-1}\cap[\Mhat_i]
+\sum_{\gbigi}
\Eu\bigl(\gbign(\Mhat^{G_m}(\gbigi))\bigr)^{-1}\cap
[\Mhat^{G_m}(\gbigi)]
\]
Then, we obtain the following equality
in $A_{\ast}(\Pic\times X^{m_1}\times D^{m_2})[t,t^{-1}]$.
\[
  \Lambda_P^{G_m}
 \bigl(\Ehat^{\Mhat}_{\ast},[\Mhat]\bigr)
=\sum_{i=1,2}
 \Lambda_P^{G_m}\left(
\bigl(
 \Ehat_{\ast}^{\Mhat}
\bigr)_{|\Mhat_i},
 \frac{[\Mhat_i]}{\Eu(\gbign\bigl(\Mhat_i)\bigr)}
 \right)
+\sum_{\gbigi}
 \Lambda_P^{G_m}\left(
 \bigl(
 \Ehat_{\ast}^{\Mhat}
\bigr)_{|\Mhat^{G_m}(\gbigi)},
 \frac{[\Mhat^{G_m}(\gbigi)]}
 {\Eu\bigl(\gbign(\Mhat^{G_m}(\gbigi))\bigr)}
 \right)
\]
The equality holds in
$A_{\ast}(\Pic\times X^{m_1}\times D^{m_2})[t]$.
Then, we obtain the following equality in $\rnum[t]$:
\[
 \int_{\Mhat}\Phi(\Ehat_{\ast}^{\Mhat})
=\sum_{i=1,2}
 \int_{\Mhat_i}
 \frac{\Phi\bigl((\Ehat_{\ast}^{\Mhat})_{|\Mhat_i}\bigr)}
 {\Eu(\gbign(\Mhat_i))}
+\sum_{\gbigi}
 \int_{\Mhat^{G_m}(\gbigi)}
 \frac{\Phi\bigl((\Ehat_{\ast}^{\Mhat})_{|\Mhat^{G_m}(\gbigi)}\bigr)}
 {\Eu(\gbign(\Mhat^{G_m}(\gbigi)))}
\]
\hfill\qed }
\end{example}

\subsubsection{The ring $\nbigr_{\CH}$}

\index{$\nbigr_{\CH}$}
We put as follows:
\[
 \nbigr_{\CH}:=
 \Sym \Bigl(\Map_f\bigl(\seisuu_{\geq\,0}^2,A^{\ast}(X)
 \bigr)\Bigr)
\otimes
 \Sym\Bigl(\Map_f\bigl(\seisuu_{\geq\,0}^3,A^{\ast}(D)\bigr)\Bigr)
\otimes A^{\ast}(\Pic)
\]
An element of $\nbigr_{\CH}$ is described as 
in the case of $\nbigr$.
When a parabolic sheaf $E_{\ast}$ is given,
we put $\nbigr_{\CH}(E_{\ast}):=\nbigr$.
Similarly, $\nbigr_{l,\CH}$ are defined for each $l\geq 1$.
When a tuple of parabolic sheaves
$\vecE_{\ast}=(E_{1\,\ast},\ldots,E_{l,\,\ast})$ is given,
we put $\nbigr_{\CH}(\vecE_{\ast}):=\nbigr_{l,\CH}$.
We use a convention (\ref{eq;06.6.19.5})
to denote elements of $\nbigr_{\CH}(\vecE_{\ast})$.
The maps $A^{\ast}(X)\lrarr H^{2\,\ast}(X)$
and $A^{\ast}(D)\lrarr H^{2\,\ast}(D)$
induce the algebra homomorphisms
$\gbigt:\nbigr_{\CH}(\vecE_{\ast})\lrarr \nbigr(\vecE_{\ast})$.

Let $\nbigz$ be an element of $A_{\ast}(\nbigy)$.
For $\prod_{h=1}^l\ch_{v^{(h)}(1)}^{v^{(h)}(2)}(E_h)\big/a$,
we put as follows:
\[
 \gbigq\Bigl(
 \prod_{h=1}^l\ch_{v^{(h)}(1)}^{v^{(h)}(2)}(E_h)\big/a
 \Bigr)\cap [\nbigz]:=
 p_{X\,\ast}\Bigl(
 \prod_{h=1}^l\ch_{v^{(h)}(1)}(E_h)^{v^{(h)}(2)}
 \cdot p_{\nbigy}^{\ast}(a)\cap [\nbigz\times X]
 \Bigr)
\]
For $\prod_{h=1}^l\ch_{u^{(h)}(1)}^{u^{(h)}(2)}
 \bigl(\Gr_{u^{(h)}(3)}(E_h)\bigr)\big/b$,
we put as follows:
\[
 \gbigq\Bigl(
 \prod_{h=1}^l\ch_{u^{(h)}(1)}^{u^{(h)}(2)}\bigl(
 \Gr_{u^{(h)}(3)}(E_h)\bigr)\big/b
 \Bigr)\cap [\nbigz]:=
 p_{D\,\ast}\Bigl(
 \prod_{h=1}^l\ch_{u^{(h)}(1)}\bigl(
 \Gr_{u^{(h)}(3)}(E_h)\bigr) ^{u^{(h)}(2)}
 \cdot p_{\nbigy}^{\ast}(b)\cap [\nbigz\times D]
 \Bigr)
\]
They induce the algebra homomorphism
$\gbigq:\nbigr_{\CH}(\vecE_{\ast})\lrarr A^{\ast}(\nbigy)$.
The following lemma can be checked by a formal calculation.
We will use it without mention.

\begin{lem}
Let $Q$ be an element of $\nbigr_{\CH}(\vecE_{\ast})$.
Let $P$ be an element of $\nbigr(\vecE_{\ast},\nbigy)$.
Let $\nbigz$ be an element of $A_{\ast}(\nbigy)$.
Then, we have the following equality:
\[
\deg\bigl(\gbigt(Q)\cdot P\cap[\nbigz]\bigr)
=
\deg\bigl(P\cap \gbigq(\nbigz)\bigr)
\]
\hfill\qed
\end{lem}

\subsubsection{The equivariant Euler class}

\label{subsubsection;06.6.19.20}
\label{subsubsection;06.6.19.50}
\index{equivariant Euler class}

For an algebra $R$,
let $R[[t^{-1},t]$ denote the algebra
of series $\sum a_j\cdot t^j$
such that the set $\bigl\{j>0|a_j\neq 0\bigr\}$ is finite.
\index{$R[[t^{-1},t]$}
We put $\gbigr(t):=\rnum[[t^{-1},t]$.
\index{$\gbigr(t)$, $\gbigr(t_1,\ldots,t_k)$}
Inductively, we put as follows:
\[
 \gbigr(t_1,t_2,\ldots,t_k):=
 \gbigr(t_2,\ldots,t_k)[[t_1^{-1},t_1]
\]

Let $E_{i\,\ast}$ $(i=1,2)$ be parabolic sheaves
over $\nbigy\times (X,D)$.
Assume that we are given $G_m$-actions on $E_{i\,\ast}$.
Then, we obtain the following complexes of
$G_m$-equivariant sheaves on $\nbigy$:
\begin{equation}
 \label{eq;06.6.19.10}
 \nbigf_1:=
 Rp_{X\,\ast}\Bigl(
 \nrhom\bigl(E_{1},E_2\bigr)
 \Bigr),
\quad
 \nbigf_2:=
 Rp_{D\,\ast}\Bigl(
 \nrhom'_2\bigl(E_{1|D\,\ast},E_{2|D\,\ast}\bigr)
 \Bigr)
\end{equation}
(See the subsubsection \ref{subsubsection;06.5.23.10}
for the notation $\nrhom'_2$.)
We have the $G_m$-equivariant Euler class 
$\Eu(\nbigf_a)\in A^{\ast}_{G_m}(\nbigy)
 \otimes_{\rnum[t]}
 \rnum[t,t^{-1}]$
of  $\nbigf_a$ $(a=1,2)$.
They are formally
$\sum_{i\leq d}c_{d(a)-i}(\nbigf_a)\cdot \bigl(w\cdot t\bigr)^i
 \in A^{\ast}(\nbigy)\otimes\gbigr(t)$,
where $d(a)$ denotes the expected rank
of $\nbigf_a$.
Due to the Grothendieck-Riemann-Roch theorem,
we have the element
$\Eu'(\nbigf_a)\in \nbigr_{\CH}(E_{1\,\ast},E_{2\,\ast})\otimes\gbigr(t)$
such that $\gbigq\bigl(\Eu'(\nbigf_a)\bigr)=\Eu(\nbigf_a)$.
For abbreviation,
we use the notation $\Eu(\nbigf_a)$
to denote $\Eu'(\nbigf_a)$ and the image of $\Eu'(\nbigf_a)$
via the composite of the morphisms:
\begin{equation}
\label{eq;06.6.19.15}
\begin{CD}
 \nbigr_{\CH}(E_{1\,\ast},E_{2\,\ast})\otimes\gbigr(t)
@>{\gbigt}>>
\nbigr(E_{1\,\ast},E_{2\,\ast})\otimes\gbigr(t)
@>{\gminiq_2}>>
\nbigr(E_{1\,\ast})\otimes 
 \nbigr(E_{2\,\ast})\otimes\gbigr(t)
\end{CD}
\end{equation}

Let $T$ denote the $k$-dimensional torus.
If $E_{i\,\ast}$ are provided with $T$-actions,
then we have the $T$-action on $\nbigf_a$.
Therefore, we obtain the $T$-equivariant
Euler classes $\Eu(\nbigf_a)$.
Due to the Grothendieck-Riemann-Roch theorem,
we have the element $\Eu'(\nbigf_a)$ of
$\nbigr_{\CH}(E_{1\,\ast},E_{2\,\ast})\otimes
 \gbigr(t_1,\ldots,t_k)$
such that
$\gbigq\bigl(\Eu(\nbigf_a)'\bigr)=
 \Eu(\nbigf_a)$.
For abbreviation,
we use the notation $\Eu(\nbigf_a)$
to denote $\Eu(\nbigf_a)'$
and the image of it via the composite
of the morphisms in (\ref{eq;06.6.19.15}).

\subsection{Transition Formulas in the Simple Cases}
\label{subsection;06.7.3.52}
\subsubsection{Basic case}
\label{subsubsection;06.6.6.1}

Let $\vecy$ be an element of $\Type$,
and let $\alpha_{\ast}$ be a system of weights.
Let $L$ be a line bundle on $X$.
Let $\delta$ denote an element of $\nbigp^{\br}$.
Let $P$ be an element of $\nbigr$,
and let $k$ be a non-negative integer.
Let $\omega$ denote the first Chern class
of $\nbigo_{\rel}(1)$
on $\nbigm^{s}(\vecyhat,[L],\alpha_{\ast},\delta)$.
We obtain the element
$\Phi:=P(\Ehat^u)\cdot\omega^k$
of $\nbigr\bigl(\Ehat^u,
 \nbigm^s(\vecyhat,[L],\alpha_{\ast},\delta)\bigr)$.
In the case where $\delta$ is not critical,
we obtain the following number
by the procedure explained in the subsubsection
\ref{subsubsection;06.6.19.2}:
\[
 \Phi(\vecyhat,[L],\alpha_{\ast},\delta):=
\int_{\nbigm(\vecyhat,[L],\alpha_{\ast},\delta)}
 \Phi
\]

Let $\delta$ be a critical parameter.
We take parameters $\delta_-<\delta<\delta_+$
such that
$\delta_{\kappa}$ $(\kappa=\pm)$
are sufficiently close to $\delta$.
We would like to describe
$\Phi(\vecyhat,[L],\alpha_{\ast},\delta_+)
-\Phi(\vecyhat,[L],\alpha_{\ast},\delta_-)$
as the sum of the integrals over the products
of the moduli stacks of the objects with lower ranks.
Such a description is called the transition formula.

For that purpose, we prepare some notation.
Let $S(\vecy,\alpha_{\ast},\delta)$
denote the set of the decomposition types:
\[
 S(\vecy,\alpha_{\ast},\delta):=
 \bigl\{\gbigi=(\vecy_1,\vecy_2)\in\Type^2\,\big|\,
 \vecy_1+\vecy_2=\vecy,\,\,\,
  P^{\alpha_{\ast},\delta}_{\vecy_1}
=P^{\alpha_{\ast}}_{\vecy_2}
=P^{\alpha_{\ast},\delta}_{\vecy}
 \bigr\}
\]
For a given $(\vecy_1,\vecy_2)\in S(\vecy,\alpha_{\ast},\delta)$,
we put $r_i=\rank \vecy_i$.
We also put as follows:
\[
 \nbigm(\vecy_1,\vecyhat_2,L,\alpha_{\ast},\delta):=
 \nbigm^{ss}(\vecy_1,L,\alpha_{\ast},\delta)
\times
 \nbigm^{ss}(\vecyhat_2,\alpha_{\ast}).
\]
We remark that the $2$-stability condition for
$(\vecy,L,\alpha_{\ast},\delta)$
implies the $1$-stability conditions for
$(\vecy_1,L,\alpha_{\ast},\delta)$
and $(\vecy_2,\alpha_{\ast})$.
Let $E_1^u$ denote the sheaf on
$\nbigm(\vecy_1,\vecyhat_2,L,\alpha_{\ast},\delta)\times X$
which is the pull back of the universal sheaf
over $\nbigm^s(\vecy_1,L,\alpha_{\ast},\delta)\times X$
via the natural morphism.
We use the notation $\Ehat_2^u$ in a similar meaning.
We put $\omega_1:=c_1(\Or(E_1^u))/\rank r_1$,
and $e^{w\cdot\omega_1}$ denotes $\Or(E_1^u)^{w/r_1}$
formally.

Let $G_m$ be the one dimensional torus.
Let $e^{w\cdot t}$ denote
the trivial bundle on
$\nbigm(\vecy_1,\vecyhat_2,L,\alpha_{\ast},\delta)$
with $G_m$-action of weight $w$.
We have the following element
of the $K$-group of the $G_m$-equivariant
coherent sheaves
on $\nbigm(\vecy_1,\vecyhat_2,L,\alpha_{\ast},\delta)$:
\begin{multline}
 \gbign_0(\vecy_1,\vecy_2)=
-Rp_{X\,\ast}\Bigl(
 \nrhom\bigl(E_1^u\!\cdot\! e^{-t},\,
 \Ehat_2^u\!\cdot\! e^{r_1(t-\omega_1)/r_2}\bigr)
\Bigr)
-Rp_{X\,\ast}\Bigl(
 \nrhom\bigl(\Ehat_2^u\!\cdot\! e^{r_1(t-\omega_1)/r_2},\,
 E^u_1\!\cdot\! e^{-t}\bigr)\Bigr) \\
-Rp_{D\,\ast}\Bigl(
 \nrhom'_2\bigl(E_{1|D\,\ast}^u\!\cdot\! e^{-t},
 \Ehat_{2|D\,\ast}^u\!\cdot\! e^{r_1(t-\omega_1)/r_2}\bigr)
\Bigr)
-Rp_{D\,\ast}\Bigl(
 \nrhom'_2\bigl(E_{2|D\,\ast}^u\!\cdot\! e^{r_1(t-\omega_1)/r_2},\,
 \Ehat_{1|D\,\ast}^u\!\cdot\! e^{-t}\bigr)
 \Bigr)\\
+Rp_{X\,\ast}\Bigl(
 \nhom\bigl(L\!\cdot\! e^{-t},\,
 \Ehat_2^u\!\cdot\! e^{r_1(t-\omega_1)/r_2}\bigr)
\Bigr)
\end{multline}
(See the subsubsection \ref{subsubsection;06.5.23.10}
for the notation $\nrhom'_2$.)
We have the equivariant Euler class
$\Eu\bigl(\gbign_0(\vecy_1,\vecy_2)\bigr)
 \in \nbigr\bigl(E_1^u,\nbigm^s(\vecy_1,L,\alpha_{\ast},\delta)\bigr)
 \otimes\nbigr(\Ehat_2)[[t^{-1},t]$.
(See the subsubsection \ref{subsubsection;06.6.19.20}
for our convention on the equivariant Euler class.)

By the homomorphisms in the subsubsection
\ref{subsubsection;06.6.19.25}
and the twist in the subsubsection
\ref{subsubsection;06.6.19.27},
we have the element
$P\bigl(E_1^u\cdot e^{-t}\oplus
 \Ehat_2^u\cdot e^{r_1(t-\omega_1)/r_2}\bigr)$
of $\nbigr\bigl(E_1^u,\nbigm^s(\vecy_1,L,\alpha_{\ast},\delta)\bigr)
\otimes \nbigr(\Ehat_2^u)[t]$.
Thus, we obtain the following element of
$\nbigr\bigl(E_1^u,\nbigm^s(\vecy_1,L,\alpha_{\ast},\delta)\bigr)
\otimes\nbigr(\Ehat_2^u)[[t^{-1},t]$:
\[
 \frac{
 P\bigl(E_1^u\cdot e^{-t}\oplus
 \Ehat_2^u\cdot e^{r_1(t-\omega_1)/r_2}\bigr)\cdot t^k}
 {\Eu\bigl(\gbign_0(\vecy_1,\vecy_2)\bigr)}
\]
By taking the residue with respect to $t$,
we obtain the following element of
$\nbigr\bigl(E_1^u,\nbigm^s(\vecy_1,L,\alpha_{\ast},\delta)\bigr)
\otimes
 \nbigr(\Ehat_2^u)$:
\begin{equation}
 \label{eq;06.6.19.26}
 \Psi(\vecy_1,\vecy_2)
=
\underset{t=0}{\res}
\left(
 \frac{P\bigl(E_1^u\cdot e^{-t}\oplus
 \Ehat_2^u\cdot e^{r_1(t-\omega_1)/r_2}\bigr)
 \cdot t^k}
 {\Eu\bigl(\gbign_0(\vecy_1,\vecy_2)\bigr)}
\right)
\end{equation}

\begin{thm}
\label{thm;06.5.31.100}
Assume that the $2$-stability condition holds
for $(\vecy,L,\alpha_{\ast},\delta)$.
Then, we have the following equality:
\begin{equation}
\label{eq;06.5.31.8}
 \Phi(\vecyhat,[L],\alpha_{\ast},\delta_+)
-\Phi(\vecyhat,[L],\alpha_{\ast},\delta_-)
=\sum_{(\vecy_1,\vecy_2)\in S(\vecy,\alpha_{\ast},\delta)} 
 \int_{\nbigm(\vecy_1,\vecyhat_2,L,\alpha_{\ast},\delta)}
 \Psi(\vecy_1,\vecy_2)
\end{equation}
The elements
$\Psi(\vecy_1,\vecy_2)\in
 \nbigr\bigl(E_1^u,\nbigm^s(\vecy_1,L,\alpha_{\ast},\delta)\bigr)
\otimes
 \nbigr(\Ehat_2^u)$
are given as in {\rm(\ref{eq;06.6.19.26})}.

The contribution of $(\vecy_1,\vecy_2)\in
S(\vecy,\alpha_{\ast},\delta)$
vanishes in the case $p_g>0$
and $\rank(\vecy_1)\neq 1$,
\end{thm}
\pf
Let $\Mhat$ denote the master space 
connecting $\nbigm^s(\vecyhat,[L],\alpha_{\ast},\delta_+)$
and $\nbigm^s(\vecyhat,[L],\alpha_{\ast},\delta_-)$
as in the subsubsection \ref{subsubsection;06.5.17.100}.
Let $\varphi:\Mhat\lrarr\nbigm(m,\vecyhat,[L])$
be the naturally defined morphism.
Let $\nbigt(1)$ denote the trivial line bundle
on $\nbigm(m,\vecyhat,[L])$ with the $G_m$-action of weight $1$.
We have the natural $G_m$-action on
$\Ehat^{\Mhat}$ and $\varphi^{\ast}\nbigo_{\rel}(1)$.
We consider the following elements
of $\nbigr_{G_m}\bigl(\Ehat^{\Mhat},\Mhat\bigr)$:
\[
 \Phi_t:=P(\Ehat^{\Mhat})\cdot 
  c_1\bigl(\varphi^{\ast}\nbigo_{\rel}(1)\bigr)^k,
\quad
 \Phitilde_t:=\Phi_t
\cdot c_1\bigl(\varphi^{\ast}\nbigt(1)\bigr).
\]
We  use Proposition \ref{prop;06.6.10.1},
then we obtain the polynomial 
$\int_{\Mhat}\Phitilde_t$ of $t$,
as explained in Example \ref{example;06.6.17.5}.
When we forget the $G_m$-action,
we have 
$c_1\bigl(\varphi^{\ast}\nbigt(1)\bigr)=0$.
Hence we have  
$\int_{\Mhat}\Phitilde_{t|t=0}=0$.
On the other hand,
we have the following equality in $\rnum[t^{-1},t]$,
due to the localization of the virtual fundamental classes
(\cite{gp}):
\[
\int_{\Mhat}
 \Phitilde_t=
 \sum_{i=1,2}
\int_{\Mhat_i}
\frac{\Phitilde_t}
{\Eu\bigl(\gbign(\Mhat_i)\bigr)}
+\sum_{\gbigi\in S(\vecy,\alpha_{\ast},\delta)}
 \int_{\Mhat^{G_m}(\gbigi)}
 \frac{\Phitilde_t}
 {\Eu\bigl(\gbign(\Mhat^{G_m}(\gbigi))\bigr)}
\]
Here, $\gbign(\Mhat_i)$ and $\gbign(\Mhat^{G_m}(\gbigi))$
denote the virtual normal bundles with the $G_m$-action
given in Proposition \ref{prop;06.6.10.1}.
We have
$c_1\bigl(\varphi^{\ast}\nbigt(1)\bigr)_{|\Mhat_i}=t$
and $c_1\bigl(\varphi^{\ast}\nbigt(1)\bigr)_{|\Mhat^{G_m}(\gbigi)}=t$.
Therefore, we obtain the following equality
in $\rnum[t,t^{-1}]$
\begin{equation}
\label{eq;06.5.31.5}
\sum_{i=1,2}
\int_{\Mhat_i}
\underset{t=0}{\Res}
\left(
 \frac{\Phi_t}
 {\Eu\bigl(\gbign(\Mhat_i)\bigr)}
\right)
+\sum_{\gbigi}
\int_{\Mhat^{G_m}(\gbigi)}
 \underset{t=0}{\Res}
\left(
 \frac{\Phi_{t}}
 {\Eu\bigl(\gbign(\Mhat^{G_m}(\gbigi))\bigr)}
\right)
=0.
\end{equation}

Let us see the contributions
from the components $\Mhat_i$.
We have 
$\iota_i^{\ast}\Ehat^{\Mhat}=\Ehat^u$
and $\iota_i^{\ast}\varphi^{\ast}\nbigo_{\rel}(1)=\nbigo_{\rel}(1)$
with the trivial $G_m$-action.
Due to Proposition \ref{prop;06.6.10.1},
we have the following equality:
\[
 \Eu\bigl(\gbign(\Mhat_i)\bigr)^{-1}
=(-1)^i\cdot (t-\omega)^{-1}
=(-1)^i\cdot\frac{1}{t}\cdot\sum_{j=0}^{\infty}
 \left(\frac{\omega}{t}\right)^j
\]
Therefore, we obtain the following:
\begin{equation}
\label{eq;06.5.31.6}
\sum_{i=1,2}
 \int_{\Mhat_i}
 \underset{t=0}\Res
\left(
 \frac{\Phi_t}
 {\Eu(\gbign(\Mhat_i))}
\right)
=\sum_{i=1,2}(-1)^i\cdot \int_{\Mhat_i}\Phi
=-\Phi(\vecyhat,[L],\alpha_{\ast},\delta_+)
  +\Phi(\vecyhat,[L],\alpha_{\ast},\delta_-)
\end{equation}

Let us calculate the contribution
from the component $\Mhat^{G_m}(\gbigi)$.
We remark that
$\varphi^{\ast}\nbigo_{\rel}(1)_{|\Mhat^{\ast}}$
and $\varphi^{\ast}\nbigt(1)_{|\Mhat^{\ast}}$
are naturally isomorphic as $G_m$-equivariant line bundles.
We use the relation of $E_i^{\Mhat}$,
$E_1^u$ and $\Ehat_2^u$
in Proposition \ref{prop;06.6.10.5}.
Then, we obtain the following equality
in $\nbigr\bigl(
 G^{\prime\ast}E_1^u,G^{\prime\,\ast}\Ehat_2^u,
 \nbigs\bigr)[t]$:
\[
 F^{\ast}\Phi_t
=G^{\prime\,\ast}
 P\bigl(E_1^u\cdot e^{-t}\oplus
 \Ehat_2^u\cdot e^{r_1(t-\omega_1)/r_2}\bigr)\cdot t^k
\]
We also have the following equality
in $\nbigr\bigl(G^{\prime\ast}E_1^u\bigr)
 \otimes\nbigr\bigl(G^{\prime\ast}\Ehat_2\bigr)[[t^{-1},t]$:
\[
 F^{\ast}\Eu\bigl(\gbign(\Mhat^{G_m}(\gbigi))\bigr)
=G^{\prime\,\ast}\Eu\bigl(\gbign_0(\vecy_1,\vecy_2)\bigr)
\]
We have the equality of the virtual fundamental classes
in Proposition \ref{prop;06.6.10.20}.
Thus, we obtain the following equality:
\begin{equation}
\label{eq;06.5.31.7}
\int_{\Mhat^{G_m}(\gbigi)}
 \underset{t=0}\Res
\left(
 \frac{\Phi_t}{\Eu\bigl(\gbign_0(\vecy_1,\vecy_2)\bigr)}
 \right)=\\
\int_{\nbigm(\vecy_1,\vecyhat_2,L,\alpha_{\ast},\delta)}
\underset{t=0}\Res
 \left(
 \frac{P\bigl(E_1^u\cdot e^{-t}\oplus 
 \Ehat_2^u\cdot e^{r_1(t-\omega_1)/r_2}\bigr)\cdot t^k}
 {\Eu\bigl(\gbign_0(\vecy_1,\vecy_2)\bigr)} \right)
\end{equation}
The desired equality (\ref{eq;06.5.31.8})
follows from (\ref{eq;06.5.31.5}),
(\ref{eq;06.5.31.6}) and (\ref{eq;06.5.31.7}).

The second claim of the theorem
immediately follows from Proposition \ref{prop;06.5.31.1}
and the first claim.
\hfill\qed

\begin{cor}
Assume $p_g>0$.
Assume that the $2$-stability condition
holds for $(\vecy,L,\alpha_{\ast},\delta)$,
and that the $2$-vanishing condition
holds for $(\vecy,L,\alpha_{\ast},\delta)$.
Then we have 
$\Phi(\vecyhat,[L],\alpha_{\ast},\delta_+)
=\Phi(\vecyhat,[L],\alpha_{\ast},\delta_-)$.
\end{cor}
\pf
It immediately follows from 
Theorem \ref{thm;06.5.31.100}
and Proposition \ref{prop;06.5.11.550}.
\hfill\qed

\subsubsection{The case of the Euler class of the relative tangent bundle}

Assume that the $1$-vanishing condition holds
for $(\vecy,L,\alpha_{\ast},\delta)$.
Let $T_{\rel}$ denote the relative tangent bundle
of the smooth morphism
$\nbigm^s(\vecyhat,[L],\alpha_{\ast},\delta)\lrarr\nbigm(\vecyhat)$.
We put as follows:
\[
 N_L(y):=\int_X\Td(X)\cdot y\cdot \ch(L).
\]
We will be interested in the integral
of the following element of
$\nbigr\bigl(\Ehat^u,
 \nbigm^s(\vecyhat,[L],\alpha_{\ast},\delta)\bigr)$:
\[
\Phi=
\frac{\Eu(T_{\rel})} {N_L(y)}\cdot P(\Ehat^u)
\]
We take parameters $\delta_-<\delta<\delta_+$
such that $\delta_{\kappa}$ $(\kappa=\pm)$
are sufficiently close to $\delta$.
The transition formula for
$\Phi(\vecyhat,[L],\alpha_{\ast},\delta)$
is rather simple,
if the $1$-vanishing condition and the $2$-stability condition
hold for $(\vecyhat,L,\alpha_{\ast},\delta)$.
Strictly speaking,
we do not need it,
because we prove a more general formula
later.
However, we give it as an explanation of the argument.
In the case $p_g>0$,
 the problem is easy.
\begin{prop}
Assume $p_g>0$.
Assume that the $1$-vanishing condition and the $2$-stability 
condition hold for $(\vecy,L,\alpha_{\ast},\delta)$.
Then, we have the equality
$\Phi(\vecyhat,[L],\alpha,\delta_-)
=\Phi(\vecyhat,[L],\alpha,\delta_+)$.
\end{prop}
\pf
By the same argument as the proof of
Theorem \ref{thm;06.5.31.100},
we can express $\Phi(\vecyhat,[L],\alpha,\delta_+)
-\Phi(\vecyhat,[L],\alpha,\delta_-)$
as the sum of the integrals
over the fixed point set $\Mhat^{G_m}(\gbigi)$.
Under the assumption of the proposition,
we have $[\Mhat^{G_m}(\gbigi)]=0$
due to Proposition \ref{prop;06.6.10.20},
Proposition \ref{prop;06.5.31.1}
and Proposition \ref{prop;06.5.11.500}.
Thus we are done.
\hfill\qed

\vspace{.1in}

Let us discuss the case $p_g=0$.
For a decomposition type
$(\vecy_1,\vecy_2)\in S(\vecy,L,\alpha_{\ast})$,
we put as follows:
\[
 \nbigm(\vecyhat_1,\vecyhat_2,[L],\alpha_{\ast},\delta):=
 \nbigm^s(\vecyhat_1,[L],\alpha_{\ast},\delta)
\times
 \nbigm^s(\vecyhat_2,\alpha_{\ast}).
\]
Let $e^{w\cdot s}$ denote the trivial
line bundle on $\nbigm(\vecyhat_1,\vecyhat_2,[L],\alpha_{\ast},\delta)$
with the $G_m$-action of weight $w$.
We have the following element of
the $K$-group 
$K^{G_m}(\nbigm(\vecyhat_1,\vecyhat_2,[L],\alpha_{\ast},\delta))$
of the $G_m$-equivariant coherent sheaves
on $\nbigm(\vecyhat_1,\vecyhat_2,[L],\alpha_{\ast},\delta)$
(see the subsubsection \ref{subsubsection;06.5.31.11}
for the notation):
\begin{multline}
 \label{eq;06.6.4.2}
 \gbign(\Ehat_1^u,\Ehat_2^u)
\cdot e^{s/r_1+s/r_2}
+\gbign(\Ehat_2^u,\Ehat_1^u)
\cdot e^{-s/r_1-s/r_2}
+\gbign_D(\Ehat_{1\,\ast}^u,\Ehat_{2\,\ast}^u)
 \cdot e^{s/r_1+s/r_2}
+\gbign_D(\Ehat_{2\,\ast}^u,\Ehat_{1\,\ast}^u)
\cdot e^{-s/r_1-s/r_2}
\end{multline}
Let $Q(\Ehat_1\cdot e^{-s/r_1},\Ehat_2\cdot e^{s/r_2})
 \in\nbigr(\Ehat_1)\otimes\nbigr(\Ehat_2)[[s^{-1},s]$
denote the equivariant Euler class of (\ref{eq;06.6.4.2}).

\begin{thm}
\label{thm;06.6.4.10}
Assume that the $2$-stability condition 
and the $1$-vanishing condition hold
for $(\vecy,L,\alpha_{\ast},\delta)$.
In the case $p_g=0$,
we have the following equality:
\begin{equation}
 \label{eq;06.6.10.50}
 \Phi(\vecyhat,[L],\alpha_{\ast},\delta_+)
-\Phi(\vecyhat,[L],\alpha_{\ast},\delta_-)
=\sum_{(\vecy_1,\vecy_2)\in S(\vecy,\alpha_{\ast},\delta)} 
\int_{\nbigm(\vecyhat_1,\vecyhat_2,[L],\alpha_{\ast},\delta)}
 \Psi(\vecy_1,\vecy_2)
\end{equation}
The elements $\Psi(\vecy_1,\vecy_2)\in
 \nbigr\bigl(\Ehat_1^u,\nbigm^s(\vecyhat_1,[L],\alpha_{\ast},\delta)\bigr)
 \otimes\nbigr\bigl(\Ehat_2^u\bigr)$
 are given as follows:
\begin{equation}
\label{eq;06.5.31.20}
 \Psi(\vecy_1,\vecy_2)=
 \frac{N_{L}(y_1)}{N_L(y)}\cdot
\underset{s=0}\Res
\left(
 \frac{P\bigl(\Ehat_1^u\cdot e^{-s/r_1}\oplus \Ehat_2^u\cdot e^{s/r_2}\bigr)}
 {Q\bigl(\Ehat_1^u\cdot e^{-s/r_1},\,\Ehat_2^u\cdot e^{s/r_2} \bigr)}
\right)
\cdot \frac{\Eu(T_{1,\rel})}{N_L(y_1)}
\end{equation}
Here $T_{1\,\rel}$ denote the vector bundle
on $\nbigm(\vecyhat_1,\vecyhat_2,[L],\alpha_{\ast},\delta)$
induced by the relative tangent bundle
of the smooth morphism
$\nbigm^s(\vecyhat_1,[L],\alpha_{\ast},\delta)
\lrarr \nbigm(\vecyhat_1)$.
\end{thm}
\pf
By using the same argument
as the proof of Theorem \ref{thm;06.5.31.100},
we can reduce the problem
to the calculation of the contributions
from $\Mhat^{G_m}(\gbigi)$
for $\gbigi=(\vecy_1,\vecy_2)\in S(\vecy,\alpha_{\ast},\delta)$.
We use the notation in the subsubsection
\ref{subsubsection;06.5.17.100}.
Let $\varphi_{\gbigi}$ denote the inclusion
$\Mhat^{G_m}(\gbigi)\lrarr \Mhat$.
We remark the relation of the virtual fundamental
classes in the case $p_g=0$ given in Proposition \ref{prop;06.5.31.1}.
Thus, we put
$[\nbigs]:=F^{\ast}\bigl([\Mhat^{G_m}(\gbigi)]\bigr)
=G^{\ast}\bigl[\nbigm(\vecyhat_1,\vecyhat_2,[L],\alpha_{\ast},\delta)\bigr]$.
Let $I_w$ denote the trivial line bundle 
with the action of $G_m$ of weight $w$.
On $\nbigs$,
we have the following decomposition
of the equivariant vector bundles:
\[
 F^{\ast}\varphi_{\gbigi}^{\ast}
 T_{\rel}=G^{\ast}T_{1\,\rel}\oplus
 G^{\ast}p_{X\,\ast}\nhom\bigl(
 L,\Ehat_2^{u}\bigl)
\otimes\nbigo_{1,\rel}(-r_1/r_2)
\otimes I_{1+r_1/r_2}
\]
Therefore, we have the following equality
in $A_{G_m}^{\ast}(\nbigs)$:
\[
 F^{\ast}\varphi_{\gbigi}^{\ast}
 \Eu\bigl(T_{\rel}\bigr)
=G^{\ast}\Eu(T_{1\,\rel})\cdot
 G^{\ast}\Eu\Bigl(
 p_{X\,\ast}\nhom\bigl(
 L,\Ehat_2^{u}
 \bigr)\cdot e^{(1+r_1/r_2)\cdot t-r_1\omega_1/r_2}
 \Bigr)
\]
The second term in the right hand side
also appears in $Eu(\gbign_0(\vecy_1,\vecy_2))$,
and hence they are cancelled out
in the evaluation of $[\nbigs]$.
Then, we obtain the following equality
in $\rnum[t^{-1},t]$:
\begin{equation}
 \label{eq;06.6.4.5}
\int_{\nbigs}
 F^{\ast} \left(
 \frac{\varphi_{\gbigi}^{\ast}\bigl(
 P\bigl(\Ehat^{\Mhat}\bigr)\cdot
\Eu(T_{\rel})\bigr)
 }{\Eu\bigl(\gbign(\Mhat^{G_m}(\gbigi))\bigr)}
\right)
=\int_{\nbigs}G^{\ast}\left(
\frac{
 P\bigl(\Ehat_1^u\cdot e^{-t+\omega_1}
\oplus \Ehat_2^u\cdot e^{r_1(t-\omega_1)/r_2}\bigr)\cdot
 \Eu(T_{1\,\rel})}
 {Q\bigl(\Ehat_1^u\cdot e^{-t+\omega_1},
 \Ehat_2^u\cdot e^{r_1(t-\omega_1)/r_2} \bigr)}
\right)
\end{equation}
Here, 
$Q\bigl(\Ehat_1^u\cdot e^{-t+\omega_1},\,
 \Ehat_2^u\cdot e^{r_1(t-\omega_1)/r_2}\bigr)$
is the Euler class of the following element
of $K^{G_m}\bigl(\nbigm(
 \vecyhat_1,\vecyhat_2,[L],\alpha_{\ast},\delta)\bigr)$,
for $A=1+r_1/r_2$:
\[
 \gbign(\Ehat_1^u,\Ehat_2^u)
\cdot e^{A(t-\omega_1)}
+\gbign(\Ehat_2^u,\Ehat_1^u)
\cdot e^{-A(t-\omega_1)}
+\gbign_D(\Ehat_{1\,\ast}^u,\Ehat_{2\,\ast}^u)
 \cdot e^{A(t-\omega_1)}
+\gbign_D(\Ehat_{2\,\ast}^u,\Ehat_{1\,\ast}^u)
\cdot e^{-A(t-\omega_1)}
\]
We remark that the integrand of the right hand side 
of (\ref{eq;06.6.4.5}) is of the form
$\sum_j A_j\cdot (t-\omega_1)^j$, where
$A_j\in\nbigr\bigl(\Ehat_1^u,
 \nbigm(\vecyhat_1,[L],\alpha_{\ast},\delta)\bigr)
 \otimes\nbigr\bigl(\Ehat_2^u\bigr)$.
By a direct calculation, we can check the following:
\[
 \underset{t=0}\Res\,
 (t-\omega_1)^j=
\left\{
\begin{array}{ll}
 1 & (j=-1)\\
 0 & (j\neq -1)
\end{array}
 \right.
\]
Hence,
we have $\Res_{t=0}(t-\omega_1)^j=\Res_{t=0}t^j$
for any $j$.
In particular,
we have the following equality:
\[
 \underset{t=0}\Res
\Bigl(
 \sum_j A_j\cdot (t-\omega_1)^j\Bigr)
=\underset{t=0}\Res \Bigl(
\sum_j A_j\cdot t^j\Bigr)
\]
Thus, we obtain the following equality:
\begin{multline}
\int_{\nbigs}
\underset{t}\Res\,
 G^{\ast}\left(
\frac{
 P\bigl(\Ehat_1^u\cdot e^{-t+\omega_1}
\oplus \Ehat_2^u\cdot e^{r_1(t-\omega_1)/r_2}\bigr)\cdot
 \Eu(T_{1\,\rel})}
 {Q\bigl(\Ehat_1^u\cdot e^{-t+\omega_1},
 \Ehat_2^u\cdot e^{r_1(t-\omega_1)/r_2} \bigr)}
\right) \\
=\int_{\nbigs}
\underset{t}\Res\,
 G^{\ast}\left(
\frac{
 P\bigl(\Ehat_1^u\cdot e^{-t}
\oplus \Ehat_2^u\cdot e^{r_1t/r_2}\bigr)\cdot
 \Eu(T_{1\,\rel})}
 {Q\bigl(\Ehat_1^u\cdot e^{-t},
 \Ehat_2^u\cdot e^{r_1t/r_2} \bigr)}
\right)
\end{multline}
Therefore, we obtain the following:
\begin{equation}
\int_{\Mhat^{G_m}(\gbigi)}
 \underset{t=0}{\Res}\left(
 \frac{P(\Ehat^{\Mhat})\cdot \Eu(T_{\rel})}
 {\Eu(\gbign\bigl(\Mhat^{G_m}(\gbigi)\bigr)}
\right)
=\\
r_1\cdot 
\int_{\nbigm(\vecyhat_1,\vecyhat_2,[L],\alpha_{\ast},\delta)}
 \underset{t=0}{\Res}\left(
 \frac{ P\bigl(\Ehat_1^u\cdot e^{-t}\oplus 
 \Ehat_2^u\cdot e^{r_1t/r_2}\bigr)
 \cdot \Eu(T_{1,\rel})}
{Q\bigl(\Ehat_1^u\cdot e^{-t},\,\,
 \Ehat_2^u\cdot e^{r_1t/r_2}\bigr)}
\right)
\end{equation}
By putting $r_1 t=s$,
we obtain the desired formula (\ref{eq;06.5.31.20}).
\hfill\qed

\subsubsection{The $\vecL$ case}
\label{subsubsection;06.6.18.1}

Let $\vecL=(L_1,L_2)$ be a pair of line bundles on $X$.
We assume that the $1$-vanishing condition holds
for $(\vecy,L_1,\alpha_{\ast})$,
and that the $2$-vanishing condition holds
for $(\vecy,L_2,\alpha_{\ast})$,
in the sense of the subsubsection
\ref{subsubsection;06.6.4.1}.
We discuss the transition formula
under the situation of the subsubsection
\ref{subsubsection;06.5.21.555}.

Let $\vecdelta=(\delta_1,\delta_2)$ be an element of
$\nbigp^{\br\,2}$.
Assume that both of $\delta_i$ are sufficiently
small as in the subsubsection \ref{subsubsection;06.5.21.555}.
Recall that the $1$-stability condition does not hold
for $(\vecy,\vecL,\alpha_{\ast},\vecdelta)$,
if and only if the following conditions hold:
\begin{itemize}
\item
$\delta_1/r_1=\delta_2/r_2$ holds for some
pair of positive integers $(r_1,r_2)$ such that $r_1+r_2$
\item
There exists a decomposition
$\vecy_1+\vecy_2=\vecy$ such that
$\rank\vecy_i=r_i$
and $P^{\alpha_{\ast}}_{\vecy_1}=P^{\alpha_{\ast}}_{\vecy_2}$.
\end{itemize}

Assume that the $1$-stability condition
does not hold for $(\vecy,\vecL,\alpha_{\ast},\vecdelta)$.
We take elements $\delta_-,\delta_+\in\nbigp^{\br}$
such that $\delta_-<\delta_1<\delta_+$
and that $|\delta_{\kappa}-\delta_1|$ $(\kappa=\pm)$
are sufficiently small.
We put $\vecdelta_{\kappa}=(\delta_{\kappa},\delta_2)$
for $\kappa=\pm$.
Let $T_{\rel}^{(1)}$ denote the relative tangent 
bundle of the smooth morphism
$\nbigm^{s}(\vecyhat,[\vecL],\alpha_{\ast},\vecdelta_{\kappa})
\lrarr \nbigm(\vecyhat,[L_2])$.
Let $\nbigo^{(i)}_{\rel}(1)$ denote the pull back
of the tautological line bundle on 
$\nbigm(\vecyhat,[L_i])$ via the morphism
$\nbigm^{ss}(\vecyhat,[\vecL],\alpha_{\ast},\delta)
\lrarr \nbigm(\vecyhat,[L_i])$.
We put $\omega^{(i)}:=c_1\bigl(\nbigo^{(i)}_{\rel}(1)\bigr)$.
We consider the following element of
$\nbigr\bigl(\Ehat,
 \nbigm^s(\vecyhat,[\vecL],\alpha_{\ast},\vecdelta_{\kappa})\bigr)$:
\[
 \Phi=\frac{\Eu(T_{\rel}^{(1)})}{N_{L_1}(y)}
\cdot P(\Ehat^u)\cdot\omega^{(2)\,k}
\]
We put as follows, for $\kappa=\pm$:
\[
 \Phi(\vecyhat,[\vecL],\alpha_{\ast},\vecdelta_{\kappa}):=
\int_{\nbigm^s(\vecyhat,[\vecL],\alpha_{\ast},\vecdelta_{\kappa})}
 \Phi
\]
We would like to discuss the transition formula
between $\Phi(\vecyhat,[\vecL],\alpha_{\ast},\vecdelta_+)$
and $\Phi(\vecyhat,[\vecL],\alpha_{\ast},\vecdelta_-)$.

\begin{prop}
\label{prop;06.6.4.15}
In the case $p_g>0$,
we have the equality:
\[
 \Phi(\vecyhat,[\vecL],\alpha_{\ast},\vecdelta_+)
=\Phi(\vecyhat,[\vecL],\alpha_{\ast},\vecdelta_-)
\]
\end{prop}
\pf
Let $\varphi:\Mhat\lrarr \nbigm(m,\vecyhat,[\vecL])$
denote the naturally defined morphism.
Let $\nbigt(1)$ denote the trivial line bundle
on $\nbigm(m,\vecyhat,[\vecL])$
with the $G_m$-action of weight $1$.
We put $\nbigi_i:=\varphi^{\ast}\nbigo^{(i)}_{\rel}(-1)$.
Let us consider the $G_m$-equivariant cohomology class:
\[
 \Phi_t:=P(\varphi^{\ast}\Ehat^u)\cdot
 \frac{\Eu\bigl(\varphi^{\ast}T^{(1)}_{\rel}\bigr)}{N_{L_1}(y)}\cdot
 c_1(\nbigi_2^{-1})^k,
\quad
 \Phitilde_t:=\Phi_t\cdot
 c_1\bigl(\varphi^{\ast}\nbigt(1)\bigr)
\]
By applying the argument as in the proof of
Theorem \ref{thm;06.5.31.100}
to $\int_{\Mhat}\Phitilde_t$,
we can obtain the description to express
$\Phi(\vecyhat,[\vecL],\alpha_{\ast},\vecdelta_+)
-\Phi(\vecyhat,[\vecL],\alpha_{\ast},\vecdelta_-)$
as the sum of the integrals
over the fixed point sets $\Mhat^{G_m}(\gbigi)$.
Under the assumption of the proposition,
we have $[\Mhat^{G_m}(\gbigi)]=0$,
due to Proposition \ref{prop;06.6.10.25},
Proposition \ref{prop;06.5.31.1}
and Proposition \ref{prop;06.5.11.500}.
Thus we are done.
\hfill\qed

\vspace{.1in}
To discuss the case $p_g=0$,
we prepare some notation.
We put as follows:
\[
 S(\vecy,\alpha_{\ast},\vecdelta):=
 \bigl\{
 (\vecy_1,\vecy_2)\in\Type^2\,\big|\,
 P_{\vecy_1}^{\alpha_{\ast}}=P^{\alpha_{\ast}}_{\vecy_2},
 \,\,\,
 \delta_1/r_1=\delta_2/r_2
 \bigr\}
\]
For any $(\vecy_1,\vecy_2)\in S(\vecy,\alpha_{\ast},\vecdelta)$,
we put as follows:
\[
 \nbigm(\vecyhat_1,\vecyhat_2,[\vecL],\alpha_{\ast},\vecdelta)
:=\nbigm^{s}(\vecyhat_1,[L_1],\alpha_{\ast},\delta_1)
\times
 \nbigm^{s}(\vecyhat_2,[L_2],\alpha_{\ast},\delta_2)
\]

Let $e^{w\cdot s}$ denote the trivial line bundle
on $\nbigm(\vecyhat_1,\vecyhat_2,
 [\vecL],\alpha_{\ast},\vecdelta)$
with the $G_m$-action of weight $w$.
We have the following element of the $K$-group
$K^{G_m}\bigl(\nbigm(\vecyhat_1,\vecyhat_2,[\vecL],\alpha_{\ast},
 \vecdelta)\bigr)$
of the $G_m$-equivariant coherent sheaves
on $\nbigm(\vecyhat_1,\vecyhat_2,[\vecL],\alpha_{\ast},\vecdelta)$:
\[
  \gbign(\Ehat_1^u,\Ehat_2^u)
\cdot e^{s/r_1+s/r_2}
+\gbign(\Ehat_2^u,\Ehat_1^u)
\cdot e^{-s/r_1-s/r_2}
+\gbign_D(\Ehat_{1\,\ast}^u,\Ehat_{2\,\ast}^u)
 \cdot e^{s/r_1+s/r_2}
+\gbign_D(\Ehat_{2\,\ast}^u,\Ehat_{1\,\ast}^u)
\cdot e^{-s/r_1-s/r_2}
\]
The equivariant Euler class is denoted by
$Q(\Ehat_1^u\cdot e^{s/r_1},\Ehat_2^u\cdot e^{s/r_2})
\in \nbigr(\Ehat_1)\otimes\nbigr(\Ehat_2)[[s^{-1},s]$.
Let $\nbigo_{2,\rel}(1)$ denote the tautological line bundle
on $\nbigm(\vecyhat_2,[L_2])$.
The pull back is also denoted by the same notation.
We put $\omega_2:=c_1\bigl(\nbigo_{2,\rel}(1)\bigr)$,
and we use the notation $e^{w\cdot \omega_2}$
to denote $\nbigo_{2,\rel}(w)$ formally.
We also have the following element
of $K^{G_m}\bigl(\nbigm(\vecyhat_1,
 \vecyhat_2,[\vecL],\alpha_{\ast},\delta)\bigr)$:
\[
 Rp_{X\,\ast}\nhom(L_2,E_1^u)\cdot e^{-s/r_1-s/r_2+\omega_2}
\]
The equivariant Euler class
is denoted by $R(L_2\cdot e^{-\omega_2+s/r_2},\,\Ehat_1^u\cdot
e^{-s/r_1})\in\nbigr(\Ehat_1)\otimes 
 A^{\ast}\bigl(\nbigm^s(\vecyhat_2,[L_2],\alpha_{\ast},\delta)\bigr)
 [[s^{-1},s]$.

\begin{prop}
\label{prop;06.5.31.30}
The following equality holds
in the case $p_g=0$:
\[
 \Phi(\vecyhat,[\vecL],\alpha_{\ast},\vecdelta_+)
-\Phi(\vecyhat,[\vecL],\alpha_{\ast},\vecdelta_-)
=\sum_{(\vecy_1,\vecy_2)\in S(\vecy,\alpha_{\ast},\vecdelta)}
 \int_{\nbigm(\vecyhat_1,\vecyhat_2,[\vecL],\alpha_{\ast},\delta)}
 \Psi(\vecy_1,\vecy_2)
\]
The elements
 $\Psi(\vecy_1,\vecy_2)\in
 \nbigr\bigl(\Ehat_1,\nbigm^s(\vecyhat_1,[L_1],\alpha_{\ast},\delta_1)
 \bigr)\otimes
 \nbigr\bigl(\Ehat_2,
 \nbigm^s(\vecyhat_2,[L_2],\alpha_{\ast},\delta_2)\bigr)$
are given as follows:
\begin{equation}
\label{eq;06.5.31.25}
 \Psi(\vecy_1,\vecy_2):=
 \frac{N_{L_1}(y_1)}{N_{L_1}(y)}\cdot
 \underset{s=0}\Res\left(
 \frac{P\bigl(\Ehat_1^u\!\cdot\! e^{-s/r_1}\oplus 
 \Ehat_2^u\!\cdot\! e^{s/r_2}\bigr)\cdot (\omega_2-s/r_2)^k}
 {Q\bigl(\Ehat_1^u\cdot e^{-s/r_1},\Ehat_2^u\cdot e^{s/r_2}\bigr)
\cdot R\bigl(L_2\!\cdot\! e^{-\omega_2+s/r_2},
 \Ehat_1^u\!\cdot\! e^{-s/r_1}\bigr)}
\right)
 \cdot \frac{\Eu(T_{1,\rel})}{N_{L_1}(y_1)} 
\end{equation}
Here, $T_{1,\rel}$ denotes the bundle
on $\nbigm(\vecyhat_1,\vecyhat_2,[\vecL],\alpha_{\ast},\delta)$
induced by
the relative tangent bundle
of the smooth morphism
$\nbigm^s(\vecyhat_1,[L_1],\alpha_{\ast},\delta)
\lrarr \nbigm(m,\vecyhat_1)$.
\end{prop}
\pf
The argument is essentially same as the proof of
Theorem \ref{thm;06.6.4.10}.
Applying the same argument as the proof of Theorem
\ref{thm;06.5.31.100}
to $\int_{\Mhat}\Phitilde_t$ in the proof of Proposition
\ref{prop;06.6.4.15},
we obtain the following expression of
$\Phi(\vecyhat,[\vecL],\alpha_{\ast},\vecdelta_+)
-\Phi(\vecyhat,[\vecL],\alpha_{\ast},\vecdelta_-)$:
\[
 \Phi(\vecyhat,[\vecL],\alpha_{\ast},\vecdelta_+)
-\Phi(\vecyhat,[\vecL],\alpha_{\ast},\vecdelta_-)
=\sum_{\gbigi\in S(\vecy,\alpha_{\ast},\vecdelta)}
 \int_{\Mhat^{G_m}(\gbigi)}
 \underset{t=0}\Res
 \left(\frac{\Phi_t}{\Eu\bigl(\gbign(\Mhat^{G_m}(\gbigi))\bigr)}\right)
\]

We have the $G_m$-equivariant isomorphism
$ F^{\ast}\varphi_{\gbigi}^{\ast}\nbigi_2^{-1}
\simeq G^{\ast}\bigl(
\nbigo_{2,\rel}(1)\bigr)
\otimes\nbigo_{1,\rel}(r_1/r_2)
 \otimes e^{-r_1\cdot t/r_2}$.
Therefore, we have the equality
$ F^{\ast}c_1\bigl(\nbigi_2^{-1}\bigr)
=G^{\ast}\bigl(
 \omega_2+r_1\cdot (\omega_1-t)/r_2
 \bigr)$.
Hence we have the following
equality in 
$\nbigr\bigl(G^{\ast}\Ehat_1,G^{\ast}\Ehat_2,\nbigs\bigr)[t]$:
\begin{multline}
 F^{\ast}\varphi_{\gbigi}^{\ast}(\Phi_t)=
\frac{1}{N_{L_1}(y)}\cdot
G^{\ast}
 P\bigl(\Ehat_1^u\cdot e^{\omega_1-t}\oplus
 \Ehat_2^u\cdot e^{-r_1(\omega_1-t)/r_2}\bigr)\cdot
 G^{\ast}\bigl(\omega_2+r_1(\omega_1-t)/r_2 \bigr)^k \\
\times G^{\ast}\Eu(T_{1,\rel})
\cdot
 G^{\ast}
 R\bigl(L_1\cdot e^{-t},\Ehat_2^u\cdot e^{-r_1(\omega_1-t)/r_2}\bigr)
\end{multline}
Here, $R\bigl(L_1\cdot e^{-t},\Ehat_2^u\cdot
e^{-r_1(\omega_1-t)/r_2}\bigr)$
denotes the equivariant Euler class of the following 
$G_m$-vector bundle:
\[
 p_{X\ast}
 \nhom(L_1,\Ehat_2^u)\cdot e^{-r_1\omega_1/r_2+(1+r_2/r_1)t}
\]
We also have the following equality
in $\nbigr(G^{\ast}\Ehat_1^u,G^{\ast}\Ehat_2^u,\nbigs)[[t^{-1},t]$:
\begin{multline}
\label{eq;06.6.29.1}
 F^{\ast}\Bigl(
\Eu\bigl(\gbign(\Mhat^{G_m}(\gbigi))\bigr)
\Bigr)
=G^{\ast}Q\bigl(\Ehat_1^u\cdot e^{\omega_1-t},\,
 \Ehat_2^u\cdot e^{-r_1(\omega_1-t)/r_2}\bigr) \\
\times 
 G^{\ast}
R\bigl(L_1\otimes e^{-t},
\Ehat_2^u\cdot e^{-r_1(\omega_1-t)/r_2}\bigr)
\cdot 
G^{\ast}R\bigl(L_2\otimes e^{-\omega_2-r_1(\omega_1-t)/r_2},\,\,
 \Ehat_1^u\cdot e^{-r_1(\omega_1-t)/r_2}
 \bigr)
\end{multline}
We put $t-\omega_1=\tbar$.
Therefore, we obtain the following:
\begin{multline}
\label{eq;06.6.4.11}
 \int_{\Mhat^{G_m}(\gbigi)}
 \underset{t=0} \Res\left(
 \frac{\Phi_t}{\Eu\bigl(\gbign(\Mhat^{G_m}(\gbigi))\bigr)}
\right)
= \\
 \frac{r_1}{N_{L_1}(y)}
 \int_{\nbigm(\vecyhat_1,\vecyhat_2,[\vecL],\alpha_{\ast},\vecdelta)}
\underset{t=0}\Res \left(
 \frac{P\bigl(\Ehat_1^u\cdot e^{-\tbar}\oplus
 \Ehat_2^u\cdot e^{r_1\tbar/r_2}\bigr)\cdot
 \bigl(\omega_2-r_1\tbar/r_2 \bigr)^k}
 {Q\bigl(E_1^u\cdot e^{-\tbar},\Ehat_2^u\cdot
 e^{r_1\tbar/r_2}\bigr)
\cdot 
 R\bigl(L_2\cdot e^{-\omega_2+r_1\tbar/r_2},
 \Ehat_1^u\cdot e^{-\tbar}\bigr) }
\right)
 \cdot
\Eu(T^{(1)}_{1,\rel}) 
\end{multline}
By an argument in the proof of Theorem \ref{thm;06.6.4.10},
we can replace $\tbar$ with $t$
in (\ref{eq;06.6.4.11}).
We put $s=r_1\cdot t$,
and then
we can show that the right hand side of (\ref{eq;06.6.4.11})
is same as the integral of $\Psi(\vecy_1,\vecy_2)$
over $\nbigm(\vecyhat_1,\vecyhat_2,[\vecL],\alpha_{\ast},\vecdelta)$.
\hfill\qed

\vspace{.1in}
Let us consider the integral of
the following element of
$\nbigr\bigl(\Ehat,\nbigm^s(\vecyhat,[\vecL],\alpha_{\ast},\vecdelta)\bigr)$,
assuming that the $1$-vanishing condition holds
for $(\vecy,L_2,\alpha_{\ast})$:
\begin{equation}
\label{eq;06.6.4.30}
 \Phi=\frac{\Eu(T^{(1)}_{\rel})}{N_{L_1}(y)}
 \cdot \frac{\Eu(T^{(2)}_{\rel})}{N_{L_2}(y)}
 \cdot P(\Ehat^u)
\end{equation}

\begin{lem}
\label{lem;06.6.4.35}
Assume that the $1$-vanishing condition holds
for $(\vecy,L_2,\alpha_{\ast})$.
Let $\Phi$ be as in {\rm(\ref{eq;06.6.4.30})}.
Then, the element
$\Psi(\vecy_1,\vecy_2)\in
 \nbigr\bigl(\Ehat_1^u,
 \nbigm^s(\vecyhat_1,[L_1],\alpha_{\ast},\delta_1)\bigr)
\otimes
 \nbigr
 \bigl(\Ehat_2^u,
\nbigm^s(\vecyhat_2,[L_2],\alpha_{\ast},\delta_2)
 \bigr)$
is given as follows:
\[
\frac{N_{L_1}(y_1)\cdot N_{L_2}(y_2)} {N_L(y_1)\cdot N_L(y_2)}
 \cdot
 \underset{s=0}\Res
 \left(
  \frac{P\bigl(\Ehat_1^u\cdot e^{-s/r_1}\oplus \Ehat_2^u\cdot e^{s/r_2}\bigr)}
 {Q\bigl(\Ehat_1^u\cdot e^{-s/r_1},\Ehat_2^u\cdot e^{s/r_2}\bigr)}
\right)
 \cdot \frac{\Eu(T_{1,\rel})}{N_{L_1}(y_1)}
 \cdot \frac{\Eu(T_{2,\rel})}{N_{L_2}(y_2)}
\]
Here $T_{2,\rel}$ denote vector bundle
on $\nbigm(\vecyhat_1,\vecyhat_2,[\vecL],\alpha_{\ast},\delta)$
obtained from the relative tangent bundle
of the smooth morphism
$\nbigm^s(\vecyhat_2,[L_2],\alpha_{\ast},\delta_2)
\lrarr\nbigm(\vecyhat_2)$.
\end{lem}
\pf
We put $\tbar=t-\omega_1$.
Then, we formally have the following decomposition:
\[
 F^{\ast}\varphi_{\gbigi}^{\ast}
 \Eu(T^{(2)}_{\rel})
=G^{\ast}\Eu(T_{2,\rel})
\cdot G^{\ast}R\bigl(
 L_2\cdot e^{-\omega_2+r_1\cdot \tbar/r_2},
 \Ehat_1^u\cdot e^{-\tbar}\bigr)
\]
We also have $G^{\ast}R\bigl(
 L_2\cdot e^{-\omega_2+r_1\cdot \tbar/r_2},
 \Ehat_1^u\cdot e^{-\tbar}\bigr)$
in the decomposition (\ref{eq;06.6.29.1}) of 
$ F^{\ast}\bigl(
\Eu\bigl(\gbign(\Mhat^{G_m}(\gbigi))\bigr)
\bigr)$,
which are cancelled out in the evaluation of $\nbigs$.
Then, the claim can be easily obtained.
\hfill\qed

\subsection{Invariants}
\label{subsection;06.7.3.53}
\subsubsection{Construction}
\label{subsubsection;06.6.18.3}

Let $\vecy$ be an element of $\Type$,
and let $\alpha_{\ast}$ be a system of weights.
Let $P$ be an element of $\nbigr$.
When the $1$-stability condition holds
for $(\vecy,\alpha_{\ast})$,
we obtain the number
$\int_{\nbigm^s(\vecyhat,\alpha_{\ast})}P(\Ehat^u)$.
We would like to obtain such a number even in the case
where the $1$-stability condition does not hold.

We take a line bundle $L$ on $X$
such that the $1$-vanishing condition
holds for $(\vecy,L,\alpha_{\ast})$
in the sense of the subsubsection
\ref{subsubsection;06.6.4.1}.
Let $\delta$ be a sufficiently small element of
$\nbigp^{\br}$
such that there are no critical value
smaller than $\delta$.
Then the $1$-stability condition holds
for $(\vecy,L,\alpha_{\ast},\delta)$.
Let $T_{\rel}$ denote the relative tangent bundle of
$\nbigm^s(\vecyhat,[L],\alpha_{\ast},\delta)\lrarr
 \nbigm(\vecyhat)$.
Then, we obtain the following number:
\[
 \Phi(L):=
\int_{\nbigm^s(\vecyhat,[L],\alpha_{\ast},\delta)}
P(\Ehat^u)\cdot \frac{\Eu(T_{\rel})}{N_L(y)}
\]

\begin{prop}
\label{prop;06.5.31.45}
\mbox{{}}
\begin{itemize}
\item
In the case $p_g>0$,
the number $\Phi(L)$ is independent of
the choice of $L$.
\item
In general,
let $L'$ and $L$ be line bundles on $X$.
We assume that $L^{-1}$ is ample.
Then there exists the limit
$\lim_{m\to\infty}\Phi(L'\otimes L^m)$,
and it is independent of the choice of 
$L'$ and $L$.
\item
Assume that the equality $P_{E_1}=P_{E_2}$ holds
for any $E_{1\,\ast}\oplus E_{2\,\ast}\in \nbigm^{ss}(\vecy,\alpha_{\ast})$.
Then, 
$\Phi\bigl(\nbigo(-m)\bigr)$ is independent of the choice of 
any sufficiently large $m$.
\end{itemize}
\end{prop}
\pf
Let $L_1'$, $L_1$ and $L_2$ be line bundles on $X$.
Assume that $L_1^{-1}$ is ample,
and that the $1$-vanishing condition holds
for $(\vecy,L_2,\alpha_{\ast})$.
We put $L_{1}^{(m)}:=L_1'\otimes L_1^m$.
If we take a sufficiently large integer $m$,
the $1$-vanishing condition holds
for $(\vecy,L_1^{(m)},\alpha_{\ast})$.
We put $\vecL^{(m)}:=(L_1^{(m)},L_2)$.
We consider the following number:
\[
 g(L_1^{(m)},L_2,\delta_1,\delta_2):=
\int_{\nbigm^s(\vecyhat,[\vecL^{(m)}],\alpha_{\ast},\vecdelta)}
P(\Ehat^u)\cdot
\frac{\Eu(T^{(1)}_{\rel})}{N_{L_1^{(m)}}(y)}
\frac{\Eu(T^{(2)}_{\rel})}{N_{L_2}(y)}
\]
Here, $T_{\rel}^{(1)}$ denotes the relative tangent bundle
of the smooth map
$\nbigm^s(\vecyhat,[\vecL^{(m)}],\alpha_{\ast},\vecdelta)
\lrarr \nbigm(\vecyhat,[L_2])$.
We use the notation $T_{\rel}^{(2)}$ in a similar meaning.
We assume that both of $\delta_i$ are sufficiently small.
When $\delta_1$ is sufficiently smaller than $\delta_2$,
we have the following equality
due to Proposition \ref{prop;06.6.11.605}:
\begin{equation}
\label{eq;06.6.4.25}
 g(L_1^{(m)},L_2,\delta_1,\delta_2)=\Phi(L_2)
\end{equation}
Similarly,
when $\delta_2$ is sufficiently smaller than $\delta_1$,
we have the following equality:
\begin{equation}
\label{eq;06.6.4.26}
 g(L_1^{(m)},L_2,\delta_1,\delta_2)=\Phi(L_1^{(m)})
\end{equation}

We fix $\delta_2$, and we move $\delta_1$.
The transition of the values $g(L_1^{(m)},L_2,\delta_1,\delta_2)$
occurs if the following holds:
\begin{itemize}
\item
$\delta_1=\delta_2\cdot r_1/r_2$ holds for some
pair of positive integers $(r_1,r_2)$.
\item
There exists a decomposition
$\vecy_1+\vecy_2=\vecy$ such that
$\rank\vecy_i=r_i$
and $P^{\alpha_{\ast}}_{\vecy_1}=P^{\alpha_{\ast}}_{\vecy_2}$.
\end{itemize}
We put as follows:
\[
 S(\vecy,\alpha_{\ast}):=
 \bigl\{(\vecy_1,\vecy_2)\,\big|\,
 P^{\alpha_{\ast}}_{\vecy_1}
=P^{\alpha_{\ast}}_{\vecy_2},
\quad\vecy_1+\vecy_2=\vecy
 \bigr\}
\]
Therefore, we have the expression of
$\Phi(L_1^{(m)})-\Phi(L_2)$
from (\ref{eq;06.6.4.25})
and (\ref{eq;06.6.4.26}):
\begin{equation}
\label{eq;06.6.4.40}
 \Phi(L_1^{(m)})-\Phi(L_2)
=\sum_{S(\vecy,\alpha_{\ast})} \gbigg(\vecy_1,\vecy_2)
\end{equation}

Due to Proposition \ref{prop;06.6.4.15},
the contributions $\gbigg(\vecy_1,\vecy_2)$
are trivial in the case $p_g>0$.
Thus we obtain the first claim.
Let us show the second claim.
We use an induction on $\rank(\vecy)$.
Due to Proposition \ref{prop;06.5.31.30}
and Lemma \ref{lem;06.6.4.35}:
\begin{multline}
\gbigg(\vecy_1,\vecy_2)=
\frac{N_{L_1^{(m)}}(y_1)\cdot N_{L_2}(y_2)}
{N_{L_1^{(m)}}(y)\cdot N_{L_2}(y)}\times\\
\int_{\nbigm(\vecyhat_1,\vecyhat_2,[\vecL^{(m)}],\alpha_{\ast}, \vecdelta')}
\underset{s=0}\Res
 \left(
 \frac{P(\Ehat_1^u\cdot e^{-s/r_1}\oplus \Ehat_2^u\cdot e^{s/r_2})}
 {Q(\Ehat_1^u\cdot e^{-s/r_1},\Ehat_2^u\cdot e^{s/r_2})}
\right)
 \cdot \frac{\Eu(T_{1,\rel})}{N_{L_1^{(m)}}(y_1)}
 \cdot \frac{\Eu(T_{2,\rel})}{N_{L_2}(y_2)}
\end{multline}
Here $\vecdelta'=(\delta_1',\delta_2')$ is
any element of $\nbigp^{\br\,2}$
such that $\delta_i'$ are sufficiently small.
We have $\rank(\vecy_1)<\rank(\vecy)$.
Recall $\nbigm(\vecyhat_1,\vecyhat_2,[\vecL^{(m)}],\alpha_{\ast},\vecdelta)
=\nbigm^s(\vecyhat_1,[L_1^{(m)}],\alpha_{\ast},\delta_1)
\times\nbigm^s(\vecyhat_2,[L_2],\alpha_{\ast},\delta_2)$.
By applying the hypothesis of the induction,
we have the following limit:
\[
  \lim_{m\to\infty}
\int_{\nbigm(\vecyhat_1,\vecyhat_2,[\vecL^{(m)}],\alpha_{\ast},\vecdelta')}
 \underset{s=0}\Res
\left(
 \frac{P(\Ehat_1^u\cdot e^{-s/r_1}\oplus \Ehat_2^u\cdot e^{s/r_2})}
 {Q(\Ehat_1^u\cdot e^{-s/r_1},\Ehat_2^u\cdot e^{s/r_2})}
\right)
 \cdot \frac{\Eu(T_{1,\rel})}{N_{L_1^{(m)}}(y_1)}
 \cdot \frac{\Eu(T_{2,\rel})}{N_{L_2}(y_2)}
\]
Moreover, it is independent of the choices of $L_1$
and $L_1'$.
We obviously have the limit:
\[
 \lim_{m\to\infty}
\frac{N_{L^{(m)}_1}(y_1)\cdot N_{L_2}(y_2)}
 {N_{L_1^{(m)}}(y)\cdot N_{L_2}(y)}
=\frac{r_1\cdot N_{L_2}(y_2)}
 {r\cdot N_{L_2}(y)}
\]
Therefore, we obtain the existence 
of the limit of the sequence 
$\bigl\{\Phi(L_1^{(m)})\bigr\}$,
and it is independent of the choice of $L_1$.
Hence, the second claim is shown.

Let us show the third claim.
We use an induction on $\rank\vecy$.
We put $L_1'=\nbigo_X$ and $L_1=\nbigo_X(-1)$.
By the hypothesis of the induction,
the following terms are independent of the choice of $m$:
\[
 \int_{\nbigm(\vecyhat_1,\vecyhat_2,[\vecL^{(m)}],\alpha_{\ast},
 \vecdelta')}
 \underset{s=0}\Res
\left(
 \frac{P(\Ehat_1^u\cdot e^{-s/r_1}\oplus \Ehat_2^u\cdot e^{s/r_2})}
 {Q(\Ehat_1^u\cdot e^{-s/r_1},\Ehat_2^u\cdot e^{s/r_2})}
\right)
 \cdot \frac{\Eu(T_{1,\rel})}{H_{y_1}(m)}
 \cdot \frac{\Eu(T_{2,\rel})}{N_{L_2}(y_2)}
\]
We also have the following equality,
due to the assumption $P_{y_1}=P_{y}$:
\[
 \frac{H_{y_1}(m)\cdot N_{L_2}(y_2)}
{H_y(m)\cdot N_{L_2}(y)}
=\frac{r_1\cdot N_{L_2}(y_2)}{r\cdot N_{L_2}(y)}
\]
Therefore, we obtain the desired independence.
\hfill\qed

\begin{df}
\label{df;06.6.4.56}
Let $P$ be an element of $\nbigr$.
We take a line bundle $L$
such that $L^{-1}$ is ample,
and we take a sufficiently small $\delta\in\nbigp^{\br}$.
Then, we put as follows:
\begin{equation}
\label{eq;06.6.4.55}
\int_{\nbigm^{ss}(\vecyhat,\alpha_{\ast})}
P(\Ehat^u)
:=
\lim_{m\to\infty}
\int_{\nbigm^s(\vecyhat,[L^m],\alpha_{\ast},\delta)}
P(\Ehat^u)\cdot\frac{\Eu(T_{\rel})}{N_{L^m}(y)}
\end{equation}
 It is well defined due to Proposition
{\rm\ref{prop;06.5.31.45}}.
\hfill\qed
\end{df}
Thus, we obtain the linear map
$\int_{\nbigm^{ss}(\vecyhat,\alpha_{\ast})}:
 \nbigr\lrarr \rnum$.

\subsubsection{Easy properties}
\label{subsubsection;06.6.18.15}

\begin{lem}
\label{lem;06.6.18.10}
Assume $p_g>0$.
We take a line bundle $L$ such that
the $1$-vanishing condition holds
for $(\vecy,L,\alpha_{\ast})$.
Then, the following equality holds:
\begin{equation}
\int_{\nbigm^{ss}(\vecyhat,\alpha_{\ast})}
P(\Ehat^u)
=\int_{\nbigm^s(\vecyhat,[L],\alpha_{\ast},\delta)}
P(\Ehat^u)\cdot\frac{\Eu(T_{\rel})}{N_{L}(y)}
\end{equation}
\end{lem}
\pf
It follows from Proposition \ref{prop;06.5.31.45}.
\hfill\qed

\vspace{.1in}

The following lemma is clear from the construction
and Proposition \ref{prop;06.6.11.200}.
\begin{lem}
 When the $1$-stability condition 
 holds for $(\vecy,\alpha_{\ast})$,
 Definition {\rm\ref{df;06.6.4.56}}
 is compatible with the ordinary definition.
\hfill\qed
\end{lem}

\begin{prop}
\label{prop;06.6.4.60}
Assume that the equality $P_{E_1}=P_{E_2}$ holds
for any $E_{1\,\ast}\oplus E_{2\,\ast}\in \nbigm^{ss}(\vecy,\alpha_{\ast})$.
Then, we have the following equality,
for any sufficiently large $m$:
\[
 \int_{\nbigm^{ss}(\vecyhat,\alpha_{\ast})}P(\Ehat^u)
=\int_{\nbigm^s(\vecyhat,[\nbigo(-m)],\alpha_{\ast},\delta)}
 P(\Ehat^u)\cdot \frac{\Eu(T_{\rel})}{H_y(m)}
\]
\end{prop}
\pf
It immediately follows from the third claim
of Proposition \ref{prop;06.5.31.45}.
\hfill\qed

\vspace{.1in}

We say a system of weights $\alpha_{\ast}$
is not critical,
if $\nbigm^{ss}(\vecy,\alpha_{\ast})=\nbigm^{ss}(\vecy,\alpha'_{\ast})$
for any $\alpha'_{\ast}$
such that $|\alpha_i-\alpha_i'|$ are sufficiently small.

\begin{cor}
\label{cor;06.6.4.100}
Assume one of the following:
\begin{itemize}
\item
$\alpha_{\ast}$ is not critical.
\item
The parabolic part of $\vecy$ is trivial.
\end{itemize}
Then, we have the following equality
for any sufficiently large $m$:
\[
\int_{ \nbigm^{ss}(\vecyhat,\alpha_{\ast})}
P(\Ehat^u)
=
\int_{\nbigm^s(\vecyhat,[\nbigo(-m)],\alpha_{\ast},\delta)}
P(\Ehat^u)\cdot\frac{\Eu(T_{\rel})}{H_y(m)}
\]
\end{cor}
\pf
We have only to check the condition
in Proposition \ref{prop;06.6.4.60}.
In the case where the parabolic part of $\vecy$ is trivial,
the condition is trivial.
Let us show the second claim.
Assume that $\alpha_{\ast}$ is not critical,
and take 
$E_{1\,\ast}\oplus E_{2\,\ast}\in
 \nbigm^{ss}(\vecy,\alpha_{\ast})$.
Let $\vecy_i$ be the types of $E_{i\,\ast}$.
Then we have 
$P_{\vecy_1}^{\alpha'_{\ast}}
=P_{\vecy_2}^{\alpha'_{\ast}}$
for any $\alpha_{\ast}'$ which are sufficiently close to 
$\alpha_{\ast}$.
It implies $P_{y_1}=P_{y_2}$,
i.e.,
$P_{E_1}=P_{E_2}$.
Thus we are done.
\hfill\qed

\subsubsection{The integral over 
 $\nbigmtilde^{ss}(\vecyhat,\alpha_{\ast},+)$}

Let $\nbigmtilde^{ss}(\vecyhat,\alpha_{\ast},+)$
be as in the subsubsection \ref{subsubsection;06.5.17.20}.
Let $\Ttilde_{\rel}$ denote the relative tangent bundle
of the smooth morphism
$\nbigmtilde^{ss}(\vecyhat,\alpha_{\ast},+)
\lrarr \nbigm(\vecyhat)$.
Recall the description of $\nbigmtilde^{ss}(\vecyhat,\alpha_{\ast},+)$
as the full flag bundle over
$\nbigm^{s}(\vecyhat,[\nbigo(-m)],\alpha_{\ast},\epsilon)$
for any sufficiently small positive number $\epsilon$.
We also remark the equality of the virtual fundamental
classes in Lemma \ref{lem;06.6.13.201}.
Then, we can easily derive the following equality:
\begin{equation}
 \label{eq;06.6.18.105}
\int_{\nbigmtilde^{ss}(\vecyhat,\alpha_{\ast},+)}
P(\Ehat^u)\cdot 
\frac{\Eu(\Ttilde_{\rel})}{H_y(m)!}
=\int_{\nbigm^{s}(\vecyhat,[\nbigo(-m)],\alpha_{\ast},\epsilon)}
P(\Ehat^u)\cdot 
 \frac{\Eu(T_{\rel})}{H_y(m)}
\end{equation}
Here $T_{\rel}$ denotes the relative tangent bundle
of the smooth morphism
$\nbigm^{s}(\vecyhat,[\nbigo(-m)],\alpha_{\ast},\epsilon)
\lrarr \nbigm(m,\vecyhat)$.

\begin{lem}
\label{lem;06.5.31.105}
Assume one of the following:
\begin{itemize}
\item $p_g>0$
\item
The condition in Proposition {\rm\ref{prop;06.6.4.60}}
is satisfied.
\end{itemize}
Then, we have the following equality:
\begin{equation}
 \label{eq;06.6.4.101}
 \int_{\nbigmtilde^{ss}(\vecyhat,\alpha_{\ast},+)}
P(\Ehat^u)\cdot 
\frac{\Eu(\Ttilde_{\rel})}{H_y(m)!}
=\int_{\nbigm^{ss}(\vecyhat,\alpha_{\ast})}
P(\Ehat^u)
\end{equation}
In particular,
if one of the conditions in Corollary
{\rm\ref{cor;06.6.4.100}} is satisfied,
the equality {\rm(\ref{eq;06.6.4.101})} holds.
\hfill\qed
\end{lem}

\subsubsection{Another expression}

In the case $p_g>0$,
we have another way to express
the integral (\ref{eq;06.6.4.55}).

\begin{lem}
\label{lem;06.6.11.700}
Assume that the $2$-vanishing condition
holds for $(\vecy,L,\alpha_{\ast})$.
We also assume $p_g>0$ and
$d:=\chi\bigl(y\cdot \ch(L)^{-1}\bigr)-1\geq 0$.
Let $P$ be an element of $\nbigr$.
The following equality holds,
for any sufficiently small $\delta\in\nbigp^{\br}$:
\[
\int_{\nbigm^{s}(\vecyhat,[L],\alpha_{\ast},\delta)}
P(\Ehat^u)\cdot \omega^d
=\int_{\nbigm^{ss}(\vecyhat,\alpha_{\ast})}
P(\Ehat^u)
\]
\end{lem}
\pf
We use the argument in the proof of
Proposition \ref{prop;06.5.31.45}.
We also use the notation in the subsubsection
\ref{subsubsection;06.6.18.1}.
We take a line bundle $L_1$ on $X$
such that the $1$-vanishing condition
holds for $(\vecy,L_1,\alpha_{\ast})$.
The pair $(L_1,L)$ is denoted by $\vecL$.
We put as follows:
\begin{equation}
 \label{eq;06.5.31.50}
 \overline{g}(L_1,L,\delta_1,\delta):=
\int_{\nbigm^{s}(\vecyhat,[\vecL],\alpha_{\ast},\vecdelta)}
P(\Ehat^u)\cdot 
 \frac{\Eu(T^{(1)}_{\rel})}{N_{L_1}(y)}
\cdot
 \omega^{(2)\,d}
\end{equation}
When $\delta_1$ is sufficiently smaller than $\delta$,
we have the following:
\[
 \overline{g}(L_1,L,\delta_1,\delta)=
\int_{\nbigm^{s}(\vecyhat,[L],\alpha_{\ast},\delta)}
P(\Ehat^u)\cdot \omega^{(2)\,d}
\]
When $\delta$ is sufficiently smaller than $\delta_1$,
we have the following equality,
due to Corollary \ref{cor;06.6.11.701}:
\[
 \overline{g}(L_1,L,\delta_1,\delta)=
 \int_{\nbigm^{ss}(\vecyhat,\alpha_{\ast})}
P(\Ehat^u)
\]
We move $\delta_1$,
and then the transitions are trivial
due to Proposition \ref{prop;06.6.4.15}.
Thus, we obtain the desired equality (\ref{eq;06.5.31.50}).
\hfill\qed

\vspace{.1in}
We recall that the $2$-vanishing condition
can be controlled numerically,
in general.
We give it for later use.

\begin{lem}
Let $y$ be an element of $\Type^{\circ}$.
Assume 
$P_y(t)>P_K(t)$
for any sufficiently large $t$,
where $K$ denotes the canonical line bundle of $X$.
Then, the $2$-vanishing condition holds for $(y,\nbigo)$.
\end{lem}
\pf
Take $E_{\ast}\in\nbigm^{ss}(y)$.
If $H^2(X,E)\neq 0$,
we have a non-trivial morphism $E\lrarr K_X$.
It implies $P_{y}(t)\leq P_K(t)$,
which contradicts the assumption.
\hfill\qed

\subsection{Rank 2 Case}
\label{subsection;06.7.3.54}
\subsubsection{Reduction to 
  the integrals over the products of Hilbert schemes}

\label{subsubsection;06.6.4.150}
In this subsubsection,
we assume $H^1(X,\nbigo)=0$ for simplicity.
Let $y$ be an element of $\Type^{\circ}_2$.
In the following,
the $H^2(X)$-part of $y$ is denoted by $a$,
and the $H^4(X)$-part of $y$ is denoted by $b$.
The second Chern class corresponding to $y$
is denoted by $n$.
We have the relation $b=a^2/2-n$.
We assume $P_y>P_K$
and $\chi(y)-1\geq 0$.
We would like to give the expression
of $\int_{\nbigm(\yhat)}P(\Ehat^u)$
as the sum of the integrals
over the products of Hilbert schemes
for any $P\in\nbigr$.

Let $NS^1(X)$ denote the subgroup of $H^2(X,\seisuu)$
generated by the $1$-cycles on $X$.
For any element $a_1\in NS^1(X)$,
we put $a_2:=a-a_1$.
Let $e^{a_i}$ denote the holomorphic line bundle
whose first Chern class is $a_i$.
Since we have assumed $H^1(X,\nbigo)=0$,
it is uniquely determined up to isomorphisms.
Let $\nbigi_i^u$ denote the universal
ideal sheaves over $X^{[n_i]}\times X$.
Let $\nbigz_i$ denote the universal $0$-scheme
over $X^{[n_i]}\times X$.
Let $\Xi_i$ denote $p_{X\,\ast}\bigl(
\nbigo_{\nbigz_i}\otimes e^{a_i}\bigr)$.
We use the same notation to denote
the pull back of them via  appropriate morphisms.
Let $G_m$ denote the one dimensional torus.
Let $e^{w\cdot s}$ denote the trivial line bundle
with the $G_m$-action of weight $w$.
Then, we have the following element of
the $K$-groups of $G_m$-equivariant sheaves
on $X^{[n_1]}\times X^{[n_2]}$
\[
-Rp_{X\,\ast}\nrhom\bigl(\nbigi_1^u\cdot e^{a_1-s},\,
 \nbigi_2^u\cdot e^{a_2+s} \bigr)
-Rp_{X\,\ast}\nrhom\bigl(\nbigi_2^u\cdot e^{a_2+s},\,
 \nbigi_1^u\cdot e^{a_1-s} \bigr)
\]
The equivariant Euler class is denoted by
$Q\bigl(\nbigi_1^u\cdot e^{a_1-s},
 \nbigi_2^u\cdot e^{a_2+s}
 \bigr)\in
 \nbigr(\nbigi_1^u\cdot e^{a_1})\otimes
 \nbigr(\nbigi_2^u\cdot e^{a_2})[[s^{-1},s]$.
We have the element
$P\bigl(\nbigi_1^u\cdot e^{a_1-s},\,\nbigi_2^u\cdot e^{a_2+s}\bigr)\in 
 \nbigr(\nbigi_1^u\cdot e^{a_1})\otimes
 \nbigr(\nbigi_2^u\cdot e^{a_2})[s]$,
which is induced as explained
in the subsubsection \ref{subsubsection;06.6.19.25}.
We also have the equivariant Euler class
$\Eu(\Xi_2\cdot e^{2s})\in 
A_{G_m}^{\ast}(X^{[n_1]}\times X^{[n_2]})$.
Thus, we obtain the following element of
$\nbigr(\nbigi_1^u\cdot e^{a_1})
 \otimes\nbigr(\nbigi_2^u\cdot e^{a_2})$:
\[
 \underset{s=0}\Res
\left(
 \frac{P\bigl(\nbigi_1^u\cdot e^{a_1-s}\oplus 
 \nbigi_2^u\cdot e^{a_2+s}\bigr)}
 {Q\bigl(\nbigi_1^u\cdot e^{a_1-s},
 \nbigi_2^u\cdot e^{a_2+s}\bigr)}
 \cdot \frac{\Eu(\Xi_1)\cdot \Eu(\Xi_2\cdot e^{2s})}
 {(2s)^{n_1+n_2-p_g}}
\right)
\]

In the case 
$\bigl(c_1(\nbigo(1)),a_1\bigr)<
\bigl(c_1(\nbigo(1)),a_2\bigr)$,
we put as follows:
\[
 \nbiga(a_1,y):=\sum_{n_1+n_2=n-a_1\cdot a_2}
\int_{X^{[n_1]}\times X^{[n_2]}}
 \underset{s=0}\Res
\left(
 \frac{P\bigl(\nbigi_1^u\cdot e^{a_1-s}\oplus 
 \nbigi_2^u\cdot e^{a_2+s}\bigr)}
 {Q\bigl(\nbigi_1^u\cdot e^{a_1-s},
 \nbigi_2^u\cdot e^{a_2+s}\bigr)}
 \cdot \frac{\Eu(\Xi_1)\cdot \Eu(\Xi_2\cdot e^{2s})}
 {(2s)^{n_1+n_2-p_g}}
\right)
\]
In the case
$\bigl(c_1(\nbigo(1)),a_1\bigr)
=\bigl(c_1(\nbigo(1)),a_2\bigr)$,
we put as follows:
\[
 \nbiga(a_1,y):=\sum_{\substack{n_1+n_2=n-a_1\cdot a_2\\ n_1>n_2}}
\int_{X^{[n_1]}\times X^{[n_2]}}
 \underset{s=0}\Res
\left(
 \frac{P\bigl(\nbigi_1^u\cdot e^{a_1-s}\oplus 
 \nbigi_2^u\cdot e^{a_2+s}\bigr)}
 {Q\bigl(\nbigi_1^u\cdot e^{a_1-s},
 \nbigi_2^u\cdot e^{a_2+s}\bigr)}
 \cdot \frac{\Eu(\Xi_1)\cdot \Eu(\Xi_2\cdot e^{2s})}
 {(2s)^{n_1+n_2-p_g}}
\right)
\]

We put as follows:
\[
 \nSW(X,y):=\bigl\{
 a_1\in NS^1(X)\,\big|\,
 \bigl[\nbigm(e^{a_1},\nbigo)\bigr]\neq 0,\,\,
 \bigl(a_1,c_1(\nbigo_X(1))\bigr)\leq
 \bigl(a,c_1(\nbigo_X(1))\bigr)/2
 \bigr\}
\]
Recall that the expected dimension of
$\nbigm(e^{a_1},\nbigo)$ is $0$,
if $[\nbigm(e^{a_1},\nbigo)]\neq 0$
(Proposition \ref{prop;06.5.13.70}).
Therefore, we can regard 
$[\nbigm(e^{a_1},\nbigo)]$ as the number,
and we denote it by $\SW(a_1)$.

\begin{thm}
\label{thm;06.6.6.10}
Assume $p_g>0$ and $H^1(X,\nbigo)=0$.
Assume $P_y>P_K$ and $\chi(y)-1\geq 0$.
We have the following equality:
\[
 \int_{\nbigm^{ss}(\yhat)}P(\Ehat^u)
+\sum_{a_1\in \nSW(X,y)}
\SW(a_1)\cdot 2^{1-\chi(y)}\cdot \nbiga(a_1,y)=0.
\]
\end{thm}
\pf
Let $\delta$ be an element of $\nbigp^{\br}$
which is not critical.
Let $\omega$ denote the first Chern class
of the relative tautological line bundle
on $\nbigm^s(\yhat,[\nbigo],\delta)$.
Then we put as follows:
\[
 \Phi(\delta):=
\int_{\nbigm^s(\yhat,[\nbigo],\delta)}
 P(\Ehat^u)\cdot \omega^{d}
\]

When $\delta\in\nbigp^{\br}$ is critical,
we put as follows:
\[
 S(y,\delta):=\bigl\{
(y_1,y_2)\in\Type\,\big|\,
 y_1+y_2=y,\,\,
 H_{y_1}+\delta=H_{y_2}
 \bigr\}
\]
We take parameters $\delta_-<\delta<\delta_+$
which are sufficiently close to $\delta$.
Let $E^u_1$ denote the universal sheaf over
$\nbigm(y_1,\nbigo)\times X$,
and $\Ehat^u_1$ denote the universal sheaf
over $\nbigm(\yhat_1,[\nbigo])\times X$.
Recall we have the isomorphism
$\nbigm(y_1,\nbigo)\simeq\nbigm(\yhat_1,[\nbigo])$,
although it does not preserve the virtual fundamental class.
We identify the moduli spaces via the isomorphism.
Then, we have the relation $E^u_1=\Ehat^u_1\cdot e^{\omega_1}$,
where  $\omega_1$ denotes the first Chern class of
the relative tautological line bundle
of $\nbigm(y_1,\nbigo)\simeq\nbigm(\yhat_1,[\nbigo])$.
Due to Theorem \ref{thm;06.5.31.100},
we have the following equality:
\begin{equation}
  \Phi(\delta_+)-\Phi(\delta_-)
=\sum_{(y_1,y_2)\in S(y,\delta)}
\int_{\nbigm(y_1,\nbigo)\times\nbigm(\yhat_2)}
 \underset{s=0}\Res
\left(
 \frac{
 P\bigl(\Ehat_1^u\cdot e^{\omega_1-s}\oplus \Ehat_2^u\cdot e^{s-\omega_1}\bigr)
 \cdot s^d}
 {Q\bigl(\Ehat_1^u\cdot e^{\omega_1-s},\,\Ehat_2^u\cdot e^{s-\omega_1}\bigr)
\cdot \Eu\bigl(Rp_{X\,\ast}\Ehat_2^u\cdot e^{2s-\omega_1}\bigr)}
\right)
\end{equation}
We put as follows:
\[
 S(y)=\bigl\{(y_1,y_2)\in\Type^2\,\big|\,
 y_1+y_2=y,\,\,P_{y_1}<P_{y_2}\bigr\}
\]
Recall $\nbigm^{ss}(\yhat,[\nbigo],\delta)=\emptyset$
for any sufficiently large $\delta$
(Proposition \ref{prop;06.4.11.1}).
On the other hand,
we have $\Phi(\delta)=\int_{\nbigm(\yhat)}P(\Ehat^u)$
for any sufficiently small $\delta$
(Lemma \ref{lem;06.6.11.700}).
Therefore, we obtain the following equality:
\begin{equation}
\label{eq;06.6.4.121}
 \int_{\nbigm(\yhat)}P(\Ehat^u)
+\sum_{(y_1,y_2)\in S(y)}
 \int_{\nbigm(y_1,\nbigo)\times\nbigm(\yhat_2)}
 \underset{s=0}\Res
\left(
 \frac{
 P\bigl(\Ehat_1^u\cdot e^{\omega_1-s}\oplus \Ehat_2^u\cdot e^{s-\omega_1}\bigr)
 \cdot s^d}
 {Q\bigl(\Ehat_1^u\cdot e^{\omega_1-s},\,\Ehat_2^u\cdot e^{s-\omega_1}\bigr)
\cdot \Eu\bigl(Rp_{X\,\ast}\Ehat_2^u\cdot e^{2s-\omega_1}\bigr)}
\right)=0
\end{equation}

If $[\nbigm(y_1,\nbigo)]\neq 0$,
the expected dimension of $\nbigm(e^{a_1},\nbigo)$ is $0$
(Proposition \ref{prop;06.5.13.70}).
Hence, we can omit $\omega_1$
in the right hand side of (\ref{eq;06.6.4.121}).
By using Proposition \ref{prop;06.5.13.35} and
$\Ehat_i^u=\nbigi_i^u\cdot e^{a_i}$,
we obtain the following:
\begin{multline}
\int_{\nbigm(y_1,\nbigo)\times\nbigm(\yhat_2)}
 \underset{s=0}\Res
\left(
 \frac{P\bigl(\Ehat_1^u\cdot e^{\omega_1-s}\oplus 
 \Ehat_2^u\cdot e^{s-\omega_1}\bigr)\cdot s^d}
 {Q\bigl(\Ehat_1^u\cdot e^{\omega_1-s},\,\Ehat_2^u\cdot e^{s-\omega_1}\bigr)
\cdot \Eu\bigl(Rp_{X\,\ast}\Ehat_2^u\cdot e^{2s-\omega_1}\bigr)}
\right)
=\\
 \SW(a_1)\cdot
\int_{X^{[n_1]}\times X^{[n_2]}}
\underset{s=0}\Res
\left(
 \frac{P\bigl(\nbigi_1^u\cdot e^{a_1-s}\oplus 
 \nbigi_2^u\cdot e^{a_2+s}\bigr)
 \cdot s^d}
 {Q\bigl(\nbigi_1^u\cdot e^{a_1-s},\,
 \nbigi_2\cdot e^{a_2+s}\bigr)
 \cdot \Eu\bigl(Rp_{X\,\ast}(\nbigi_2^u\cdot e^{a_2+2s})\bigr)
 }\right)
\cdot\Eu(\Xi_1)
\end{multline}
By a formal calculation, we obtain the following:
\[
\frac{1}{\Eu\bigl(Rp_{X\,\ast}
 (\nbigi_2^u\cdot e^{a_2+2s})\bigr)}
=\frac{\Eu(\Xi_2\cdot e^{2s})}{(2s)^{\chi(a_2)}}
\]
Since $[\nbigm(e^{a_1},\nbigo)]\neq 0$,
we have $\chi(a_1)=1+p_g$
(Proposition \ref{prop;06.5.13.70}).
Hence, we have the following:
\[
 \chi(a_2)=\chi(y_2)+n_2
=\chi(y)-\chi(y_1)+n_2=\chi(y)-1-p_g+n_1+n_2.
\]
Then, the desired equality can be obtained 
by a direct calculation.
\hfill\qed

\subsubsection{Dependence on the polarization}

Let $y$ be an element of $\Type_2^{\circ}$.
We use a similar convention as in the subsubsection
\ref{subsubsection;06.6.4.150}.
To distinguish the dependence on the polarization $H$,
we use the notation
$\nbigm_H(\yhat)$ to denote the moduli stack
of torsion-free sheaves of type $y$
which are semistable with respect to $H$.
Let $\xi$ be an element of $NS(X)$
such that $\xi+a$ is divisible by $2$ in $NS(X)$.
We also assume $a^2-4n\leq \xi^2<0$.
We put $W^{\xi}:=\bigl\{c\in
NS(X)\otimes\real\,\big|\,(c,\xi)=0\bigr\}$,
which is called the wall determined by $\xi$.
It is well known that
the ample cone is divided into the chambers
by such walls,
and the moduli $\nbigm_H(\yhat)$ depends
only on the chambers to which $H$ belongs.
For $\Phi:=P(\Ehat^u)\in\nbigf\bigl(\nbigm_H(\yhat)\bigr)$,
we put as follows:
\[
 \Phi_H(\yhat):=\int_{\nbigm_H(\yhat)}\Phi
\]
We would like to discuss how $\Phi_H(\yhat)$ changes
when the polarizations vary across the wall $W^{\xi}$.
We put as follows:
\[
 S(y,\xi)=\bigl\{(y_0,y_1)\in \Type(X)^2\,\big|\,
 y_0+y_1=y,\,\,
 a_0-a_1=m\cdot \xi\,\,(m>0)\,\,
 \mbox{\rm in }H^{2}(X)
 \bigr\}
\]
For each $(y_0,y_1)\in S(y,\xi)$,
we put 
$\nbigm(\yhat_0,\yhat_1):=
 \nbigm(\yhat_0)\times\nbigm(\yhat_1)$.
Let $\Ehat_i$ denote the sheaf over $\nbigm(\yhat_0,\yhat_1)\times X$
which is the pull back of the universal sheaves
over $\nbigm(\yhat_i)\times X$ via the appropriate projection.
Let $e^{w\cdot s}$ denote the trivial line bundle
with the $G_m$-action of weight $w$.
We have the following element of the $K$-group
of the $G_m$-equivariant sheaves on $\nbigm(\yhat_0,\yhat_1)$:
\[
-Rp_{X\,\ast}
 \nrhom\bigl(\Ehat_0\cdot e^{-s},\,\Ehat_1\cdot e^s\bigr)
-Rp_{X\,\ast}
 \nrhom\bigl(\Ehat_1\cdot e^s,\,\Ehat_0\cdot e^{-s}\bigr)
\]
The equivariant Euler class is denoted by
$Q\bigl(\Ehat_0\cdot e^{-s},\Ehat_1\cdot e^{s}\bigr)
 \in\nbigr(\Ehat_0)\otimes\nbigr(\Ehat_1)[[s^{-1},s]$.
By the homomorphisms in the subsubsection
\ref{subsubsection;06.6.19.25},
we have $P\bigl(\Ehat_0\cdot e^{-s}\oplus \Ehat_1\cdot e^s\bigr)
 \in \nbigr(\Ehat_0)\otimes\nbigr(\Ehat_1)[s]$.

\begin{thm}
\label{thm;06.6.4.201}
Let $C_+$ and $C_-$ be chambers
which are divided by the wall $W^{\xi}$.
Let $H_+$ and $H_-$ be ample line bundles
contained in $C_+$ and $C_-$, respectively.
We assume $(H_-,\xi)<0<(H_+,\xi)$.
\begin{itemize}
\item
 In the case $p_g>0$,
 we have $\Phi_{H_+}(\yhat)=\Phi_{H_-}(\yhat)$.
 Namely, the invariant does not depend
 on the choice of generic polarization.
\item
 In the case $p_g=0$,
we have the following equality:
\begin{equation}
\label{eq;06.6.4.200}
 \Phi_{H_+}(\yhat)-\Phi_{H_-}(\yhat)
=\sum_{(y_0,y_1)\in S(y,\xi)}
 \int_{\nbigm(\yhat_0,\yhat_1)}
 \underset{s=0}\Res\left(
 \frac{P\bigl(\Ehat_0^u\cdot e^{-s}\oplus \Ehat_1^u\cdot e^s\bigr)}
 {Q\bigl(\Ehat_0\cdot e^{-s},\,\Ehat_1\cdot e^s\bigr)}
\right)
\end{equation}
\end{itemize}
\end{thm}

We give two arguments to prove Theorem \ref{thm;06.6.4.201}.
Both of them are based on the following observation.
\begin{lem}
[\cite{eg}, \cite{mw}]
\label{lem;06.6.4.202}
Let $H$ be an ample line bundle contained in $W^{\xi}$.
We can take $H_{\kappa}\in C_{\kappa}$ $(\kappa=\pm)$
satisfying the following:
\begin{itemize}
\item
There exists a very ample curve $C$
such that $H_+=H\otimes\nbigo(C)$
and $H_-=H\otimes\nbigo(-C)$.
\item
The following holds for a torsion-free sheaf $E$ of type $y$:
\begin{itemize}
\item
$E$ is $H_+$-semistable
if and only if $E(C)$ is $H$-semistable.
\item
$E$ is $H_-$ -semistable
if and only if $E(-C)$ is $H$-semistable.
\hfill\qed
\end{itemize}
\end{itemize}
\end{lem}

Ellingsrud and G\"ottsche proved the formula
(\ref{eq;06.6.4.200})
for the Donaldson invariant
under the assumption that the wall $W^{\xi}$ is {\em good}.
They used the parabolic structure
$E(-C)\subset E(C)$ with weight $\alpha$
for torsion-free sheaves $E$.
Let $\nbigm^{ss}(\yhat,\alpha)$ 
denote the moduli stack
of torsion-free sheaves
with the parabolic structure as above,
which are semistable with respect to the polarization $H$.
Due to Lemma \ref{lem;06.6.4.202},
we have 
$\nbigm^{ss}(\yhat,1)=\nbigm_{H_+}(\yhat)$
and $\nbigm^{ss}(\yhat,0)=\nbigm_{H_-}(\yhat)$.
We say that $\alpha$ is critical,
if $\nbigm^{ss}(\yhat,\alpha)\neq \nbigm^{ss}(\yhat,\alpha')$
for any sufficiently close $\alpha'\neq \alpha$.
We can easily show that
there are only finitely many critical values,
by using some boundedness result.
By investigating the transition of the invariants
at critical parabolic weights,
they obtained the formula (\ref{eq;06.6.4.200}).
Using their framework and our transition formula 
(\ref{eq;06.6.10.50}),
we will show the formula (\ref{eq;06.6.4.200})
without the assumption that the wall is good,
in the subsubsection \ref{subsubsection;06.6.4.210}.

We will give another argument
for the proof of Theorem \ref{thm;06.6.4.201}
in the subsubsection \ref{subsubsection;06.6.4.220}.
Perhaps, it may be a little more suitable
when we discuss a similar problem
in the higher rank case.

\vspace{.1in}
We give some preliminary for the argument.
In the rest of this subsection,
we use the polarization $H$.
Namely, the $\mu$-semistability condition
and the semistability condition are considered
with respect to $H$.
Let $\nbigs$ denote the family of
$\mu$-semistable torsion-free sheaves of type $y$.
Let $\nbigsbar$ denote the family
of torsion-free sheaves $E'$ of rank one
with the following property:
\begin{itemize}
\item
 $\mu(E')=\mu(y)$.
\item
 There is a member $E$ of $\nbigs$,
 such that $E'$ is a saturated subsheaf of $E$.
\end{itemize}
The families $\nbigs$ and $\nbigsbar$ are bounded.
We take a sufficiently large integer $m$
such that the family $\nbigs$ satisfies the condition $O_m$.
In the rest of this subsection,
$\delta$ denotes a polynomial of degree $0$.

\subsubsection{The proof of Theorem \ref{thm;06.6.4.201} (I)}
\label{subsubsection;06.6.4.210}

For a torsion-free sheaf $E$,
we denote by $F^{(\alpha)}$
the parabolic structure
$\bigl(E(-C)\subset E(C),\alpha\bigr)$.
The following lemma is clear.
\begin{lem}
If $(E,F^{(\alpha)},\phi)$ is
a $\delta$-semistable parabolic $L$-Bradlow pair,
$E$ is $\mu$-semistable.
(Recall $\delta$ is assumed to be a polynomial of degree $0$.)
\hfill\qed
\end{lem}

Let $\nbigm^{ss}\bigl(\yhat,[L],\alpha,\delta\bigr)$
denote the moduli stack of the oriented 
parabolic reduced $L$-Bradlow pairs
of type $y$ with weight $\alpha$,
where the parabolic structure is given as above.

Let $\alpha$ be any real number,
and let $E'$ be any member of $\nbigsbar$.
We take $E\in\nbigs$ such that $E'$ is a saturated subsheaf of $E$.
We put $E'':=E/E'$.
Then, we put $F_{\alpha}(E'):=P^{\alpha}_{E'}-P^{\alpha}_{E''}$.
The number is determined by $\ch(E')$,
and hence it is independent of the choice of $E$.
When we fix $\alpha$,
the function $F_{\alpha}:\nbigsbar\lrarr\real$
has finitely many values,
due to the boundedness of $\nbigsbar$.
It is clear that $\alpha$ is critical
if and only if there exists a member $E'\in\nbigsbar$
such that $F_{\alpha}(E')=0$.

Let $\alpha_0$ be critical,
and let $\epsilon$ be any sufficiently small positive number.
We can take a small positive number $\eta>0$
such that the following holds
for any $\alpha'$ with $|\alpha'-\alpha_0|<\eta$:
\begin{itemize}
\item
 $F_{\alpha_0}(E')>0$
 $\Longleftrightarrow$
 $F_{\alpha'}(E')>\epsilon$
\item
 $F_{\alpha_0}(E')<0$
 $\Longleftrightarrow$
 $F_{\alpha'}(E')<-\epsilon$
\item
$F_{\alpha_0}(E')=0$
 $\Longleftrightarrow$
 $|F_{\alpha'}(E)|<\epsilon$
\end{itemize}

\begin{lem}
\label{lem;06.6.5.1}
We have 
$\nbigm^{ss}\bigl(\yhat,[L],\alpha_0,\epsilon\bigr)
\simeq \nbigm^{ss}\bigl(\yhat,[L],\alpha',\epsilon\bigr)$
for any $\alpha'$ such that $|\alpha_0-\alpha'|<\eta$.
Moreover,
we have
$\nbigm^{ss}(\yhat,[L],\alpha_0,\epsilon)
=\nbigm^{s}(\yhat,[L],\alpha_0,\epsilon)$.
\end{lem}
\pf
Let $(E,\rho,F^{(\alpha')},[\phi])$ be 
an oriented parabolic reduced $L$-Bradlow pair,
such that $E$ is $\mu$-semistable.
Let $E'\subset E$ be a member of $\nbigsbar$.
We put $E'':=E/E'$.
We have $P_{E'}^{\alpha_0}+\epsilon<P^{\alpha_0}_{E''}$
if and only if $P_{E'}^{\alpha'}+\epsilon <P^{\alpha'}_{E''}$.
We have $P_{E'}^{\alpha_0}<P^{\alpha_0}_{E''}+\epsilon$
if and only if $P_{E'}^{\alpha'}<P^{\alpha'}_{E''}+\epsilon$.
Then, the first claim of the lemma is clear.
The second claim is also easy to see.
\hfill\qed

\vspace{.1in}

To compare $\Phi_{H_+}(\yhat)$ and $\Phi_{H_-}(\yhat)$,
we would like to consider the invariants
obtained from the moduli stacks
$\nbigm^{s}\bigl(\yhat,[\nbigo(-m)],\alpha,\delta\bigr)$.
We remark that the divisor of the parabolic structure
is not reduced in this case,
contrast to that we assumed the smoothness of the divisor
in the section \ref{section;06.6.4.250}.
We can deal with the point 
by the following two arguments:
\begin{enumerate}
\item
 We do not have to think
 the contribution of the parabolic structure
 to the obstruction theory,
 because the filtration is canonically determined
 by the sheaf.
 In fact, $\nbigm^{s}(\yhat,[L],\alpha,\delta)$
 is an open substack of $\nbigm(m,\yhat,[L])$
 for a sufficiently large $m$, in this case.
 Thus we obtain the perfect obstruction theory
 $\Ob(m,\yhat,[L])$  and the virtual fundamental class.
\item
 In the case $\alpha>1/2$,
 we consider the parabolic structure
 $F^{(\alpha)}_1$ given by
 $(E\subset E(C),2\alpha-1)$
  for a torsion-free sheaf $E$.
 Then,
 $(E,F^{(\alpha)},[\phi])$
 is $\delta$-semistable
 if and only if
 $(E,F^{(\alpha)}_1,[\phi])$ is $\delta$-semistable.
 In the case $\alpha<1/2$,
 we consider the parabolic structure
 $F^{(\alpha)}_2$ given by
 $(E(-C)\subset E,2\alpha)$
 for a torsion-free sheaf $E$.
 Then, $(E,F^{(\alpha)},[\phi])$ 
 is $\delta$-semistable if and only if
 $(E,F^{(\alpha)}_1,[\phi])$ is $\delta$-semistable.
\end{enumerate}
Therefore, we can freely apply our previous results.

Let $T_{\rel}$ denote the relative tangent bundle
of the smooth morphism
$\nbigm(\yhat,[\nbigo(-m)],\alpha,\delta)
\lrarr \nbigm(m,\yhat)$.
If $\delta$ is not critical
with respect to $(y,\nbigo(-m),\alpha)$,
we obtain the following number:
\[
 \Phitilde(\alpha,\delta):=
 \int_{\nbigm^s(\yhat,[\nbigo(-m)],\alpha,\delta)}
 \Phitilde.
\]
Here, we put as follows:
\[
 \Phitilde:=\Phi\cdot \frac{\Eu(T_{\rel})}{H_y(m)}
 \in \nbigr\bigl(\Ehat,\nbigm(\yhat,[\nbigo(-m)],\alpha,\delta)\bigr)
\]
The following lemma is the special case
of Corollary \ref{cor;06.6.4.100}.
\begin{lem}
Assume that $\alpha$ is not critical.
Then, $\Phitilde(\alpha,\delta)$ is independent of
the choice of $m$,
if $\delta$ is sufficiently small.
The number is denoted by $\Phi(\alpha)$.
\hfill\qed
\end{lem}

We have $\Phi(1)=\Phi_{H_+}(y)$
and $\Phi(0)=\Phi_{H_-}(y)$.
Therefore, we have only to see 
the transition of the invariants
at critical weights.
For a critical $\alpha$,
we put as follows:
\[
 S(y,\xi,\alpha):=\bigl\{
 (y_0,y_1)\in S(y,\xi)\,\big|\,
 P^{\alpha}_{y_0}=P^{\alpha}_{y_1}
 \bigr\}
\]

\begin{prop}
\label{prop;06.6.5.3}
We take real numbers $\alpha_+>\alpha>\alpha_-$
such that $|\alpha_{\kappa}-\alpha|<\eta$.
In the case $p_g>0$,
we have $\Phi(\alpha_+)-\Phi(\alpha_-)=0$.
In the case $p_g=0$,
we have the following formula:
\begin{equation}
\label{eq;06.6.5.4}
  \Phi(\alpha_+)-\Phi(\alpha_-)
=\sum_{(y_0,y_1)\in S(y,\xi,\alpha)}
 \int_{\nbigm(\yhat_0,\yhat_1)}
 \underset{s=0}\Res\left(
 \frac{P\bigl(\Ehat_0^u\cdot e^{-s}\oplus \Ehat_1^u\cdot e^{s}\bigr)}
 {Q(\Ehat_0^u\cdot e^{-s},\Ehat_1^u\cdot e^s)}
\right)
\end{equation}
\end{prop}
\pf
We have 
$\nbigm^{s}(\yhat,[\nbigo(-m)],\alpha_{+},\epsilon)
\simeq
 \nbigm^{s}(\yhat,[\nbigo(-m)],\alpha_-,\epsilon)$
due to Lemma \ref{lem;06.6.5.1}.
Therefore,
we obtain the following:
\begin{equation}
 \label{eq;06.6.5.2}
 \Phi(\alpha_+)-\Phi(\alpha_-)
=\bigl(\Phi(\alpha_-,\epsilon)-\Phi(\alpha_-)\bigr)
-\bigl(\Phi(\alpha_+,\epsilon)-\Phi(\alpha_+)\bigr)
\end{equation}

Let us see the first term in the right hand side
of (\ref{eq;06.6.5.2}).
We see the transition of $\Phi(\alpha_+,\delta)$
when we move $\delta$ from $0$ to $\epsilon$.
We use Theorem \ref{thm;06.6.4.10}.
The transition occurs 
when $P^{\alpha_-}_{y_0}+\delta=P^{\alpha_-}_{y_1}$
holds for some $(y_0,y_1)\in S(y,\xi,\alpha)$.
In the case $p_g>0$,
the transitions are trivial.
Hence we obtain $\Phi(\alpha_-,\epsilon)-\Phi(\alpha_-)=0$.
For any $(y_0,y_1)\in S(y,\xi,\alpha)$,
we put 
$\nbigm(\yhat_0,\yhat_1,[\nbigo(-m)])
:=\nbigm\bigl(\yhat_0,[\nbigo(-m)]\bigr)
\times\nbigm\bigl(\yhat_1\bigr)$.
In the case $p_g=0$, we obtain the following equality:
\begin{multline}
\label{eq;06.6.5.25}
 \Phi(\alpha_-,\epsilon)-\Phi(\alpha_-)
=\!\!\!\!\!\sum_{(y_0,y_1)\in S(y,\xi,\alpha)}
 \frac{H_{y_0}(m)}{H_{y}(m)}\cdot
\int_{\nbigm(\yhat_0,\yhat_1,[\nbigo(-m)])}
 \underset{s=0}\Res\left(
 \frac{P\bigl(\Ehat_0^u\cdot e^{-s}\oplus \Ehat_1^u\cdot e^{s}\bigr)}
 {Q(\Ehat_0^u\cdot e^{-s},\Ehat_1^u\cdot e^s)}
\right)
 \cdot \frac{\Eu(T_{0,\rel})}{H_{y_0}(m)} \\
 =\sum_{(y_0,y_1)\in S(y,\xi,\alpha)}
 \frac{H_{y_0}(m)}{H_{y}(m)}\cdot
\int_{\nbigm(\yhat_0,\yhat_1)}
 \underset{s=0}\Res\left(
 \frac{P\bigl(\Ehat_0^u\cdot e^{-s}\oplus \Ehat_1^u\cdot e^{s}\bigr)}
 {Q(\Ehat_0^u\cdot e^{-s},\Ehat_1^u\cdot e^s)}
\right)
\end{multline}
Similarly, we have
$\Phi(\alpha_+,\epsilon)-\Phi(\alpha_+)=0$
in the case $p_g>0$,
and we have the following equality in the case $p_g=0$:
\begin{multline}
\label{eq;06.6.5.30}
\Phi(\alpha_+,\epsilon)-\Phi(\alpha_+) 
=\!\!\!\!\!\sum_{(y_0,y_1)\in S(y,\xi,\alpha)}
 \frac{H_{y_1}(m)}{H_{y}(m)}\cdot
\int_{\nbigm(\yhat_1,\yhat_0)}
 \underset{s=0}\Res
\left(
 \frac{P\bigl(\Ehat_1^u\cdot e^{-s}\oplus \Ehat_0^u\cdot e^{s}\bigr)}
 {Q(\Ehat_1^u\cdot e^{-s},\Ehat_0^u\cdot e^s)}
\right)
 \\
=\sum_{(y_0,y_1)\in S(y,\xi,\alpha)}
- \frac{H_{y_1}(m)}{H_{y}(m)}\cdot
\int_{\nbigm(\yhat_0,\yhat_1)}
 \underset{s=0}\Res\left(
 \frac{P\bigl(\Ehat_0^u\cdot e^{-s}\oplus \Ehat_1^u\cdot e^{s}\bigr)}
 {Q(\Ehat_0^u\cdot e^{-s},\Ehat_1^u\cdot e^s)}
\right)
\end{multline}
Therefore, we obtain the following:
\begin{multline}
 \Phi(\alpha_+)-\Phi(\alpha_-)
=\!\!\!\!\!\sum_{(y_0,y_1)\in S(y,\xi,\alpha)}
 \left(
 \frac{H_{y_0}(m)}{H_y(m)}
+\frac{H_{y_1}(m)}{H_y(m)}
 \right)
\cdot 
 \int_{\nbigm(\yhat_0,\yhat_1)}
 \underset{s=0}\Res\left(
 \frac{P\bigl(\Ehat_0^u\cdot e^{-s}\oplus \Ehat_1^u\cdot e^{s}\bigr)}
 {Q(\Ehat_0^u\cdot e^{-s},\Ehat_1^u\cdot e^s)}
\right)
 \\
=\sum_{(y_0,y_1)\in S(y,\xi,\alpha)}
 \int_{\nbigm(\yhat_0,\yhat_1)}
 \underset{s=0}\Res\left(
 \frac{P\bigl(\Ehat_0^u\cdot e^{-s}\oplus \Ehat_1^u\cdot e^{s}\bigr)}
 {Q(\Ehat_0^u\cdot e^{-s},\Ehat_1^u\cdot e^s)}
\right)
\end{multline}
Thus the proof of Proposition \ref{prop;06.6.5.3}
is finished.
\hfill\qed

\vspace{.1in}
The first claim of Theorem \ref{thm;06.6.4.201}
obviously follows from the first claim of
Proposition \ref{prop;06.6.5.3}.
We have $S(y,\xi)=\bigcup_{0<\alpha<1}S(y,\xi,\alpha)$
since the intersection pairing $(C,\xi)$ is sufficiently large.
Then, the formula (\ref{eq;06.6.4.200})
immediately follows from (\ref{eq;06.6.5.4}).
Thus, the first proof of Theorem \ref{thm;06.6.4.201}
is finished.
\hfill\qed

\subsubsection{The proof of Theorem \ref{thm;06.6.4.201} (II)}

\label{subsubsection;06.6.4.220}

We put $y(C):=y\cdot \ch(\nbigo(C))$.
We use the notation $y(-C)$ in a similar meaning.
We regard them as the element of $\Type$
whose parabolic parts are trivial.
We put $\nbigo(-m,C):=\nbigo(-m)\otimes\nbigo(C)$
and $\nbigo(-m,-C):=\nbigo(-m)\otimes\nbigo(-C)$.
Let $(E,[\phi])$ be a reduced $\nbigo(-m)$-Bradlow pair
such that $\phi\neq 0$.
We naturally obtain the reduced 
$\nbigo(-m,C)$-Bradlow pair $(E(C),[\phi_C])$.
Similarly, we obtain 
the reduced $\nbigo(-m,-C)$-Bradlow pair $(E(-C),[\phi_{-C}])$.
Let $T_{\rel}$ denote the relative tangent bundle
of the smooth map
$\nbigm^s\bigl(\yhat(C),[\nbigo(-m,C)],\delta\bigr)
\lrarr \nbigm^s(m,\yhat(C))$.
When $\delta$ is not critical,
we put as follows:
\begin{equation}
\label{eq;06.6.6.16}
 \Phitilde_C(\delta):=
\int_{\nbigm^s(\yhat(C),[\nbigo(-m,C)],\delta)}
 \Phi\cdot \frac{\Eu(T_{\rel})}{H_y(m)},
\quad\quad
  \Phitilde_{-C}(\delta):=
\int_{\nbigm^s(\yhat(-C),[\nbigo(-m,-C)],\delta)}
 \Phi\cdot \frac{\Eu(T_{\rel})}{H_y(m)}.
\end{equation}
When $\delta$ is sufficiently small,
we have $\Phitilde_C(\delta)=\Phi_{H_+}(y)$
and $\Phitilde_C(\delta)=\Phi_{H_-}(y)$.

Let $\Type(\nbigsbar)$ denote the set
$\bigl\{\ch(E')\in H^{\ast}(X)\,\big|\,E'\in\nbigsbar\bigr\}$.
For any $y_0\in\nbigsbar$,
we put $y_1:=y-y_0$,
and then $y_1$ is also an element of $\Type(\nbigsbar)$.
We put $y_i(C):=y_i\cdot \ch\bigl(\nbigo(C)\bigr)$.
We use the notation
$y_i(-C)$ in a similar meaning.
We remark that
$P_{y_0(C)}-P_{y_1(C)}$ and
$P_{y_0(-C)}-P_{y_1(-C)}$
are the polynomials of degree $0$.

Let $a_i$, $b_i$ and $n_i$ denote 
the first Chern class, the second Chern character
and the second Chern class corresponding to $y_i$.
We have $(a_0,H)=(a_1,H)$,
and $H$ is a generic element of $W^{\xi}$.
Therefore, we have $a_0-a_1=A\cdot \xi$
for some $A\in\rnum$ in $H^2(X)$.
Since $(C,\xi)$ is assumed to be sufficiently large,
we have $P_{y_0(C)}-P_{y_1(C)}\neq 0$
unless $y_0=y_1$.
We also have $P_{y_0(-C)}-P_{y_1(-C)}\neq 0$
unless $y_0=y_1$.

We put as follows:
\[
 S(y,C):=
 \bigl\{(y_0,y_1)\,\big|\,
 y_0+y_1=y,\,\,y_i\in\Type(\nbigsbar),\,\,
 P_{y_0(C)}<P_{y_1(C)}
 \bigr\}
\]
\[
   S(y,-C):=
 \bigl\{(y_0,y_1)\,\big|\,
 y_0+y_1=y,\,\,y_i\in\Type(\nbigsbar),\,\,
 P_{y_0(-C)}<P_{y_1(-C)}
 \bigr\}.
\]

We take a positive constant $\delta_0$
satisfying the following inequalities:
\[
 \delta_0>
 \max\bigl\{
 |P_{y_0(C)}-P_{y_1(C)}|\,
\big|\,\,
 (y_0,y_1)\in S(y,C)
 \bigr\},
\quad
 \delta_0>
\max\bigl\{
 |P_{y_0(-C)}-P_{y_1(-C)}|\,
\big|\,\,
 (y_0,y_1)\in S(y,-C)
 \bigr\}.
\]

\begin{lem}
\label{lem;06.6.5.10}
$(E(C),\phi_C)$ is $\delta_0$-semistable,
if and only if 
the following conditions hold:
\begin{itemize}
\item
 $E$ is $\mu$-semistable.
\item
 For any subobject 
 $(E',\phi')\subset (E,\phi)$ such that $\phi'\neq 0$,
 we have $\mu(E')<\mu(E)$.
\end{itemize}
Moreover,
$(E(C),\phi_C)$ is $\delta_0$-stable,
if it is $\delta_0$-semistable.
\end{lem}
\pf
Assume that $(E(C),\phi_C)$ is $\delta_0$-semistable.
Since $\delta_0$ is a polynomial of degree $0$,
it is easy to see that $E(C)$ is $\mu$-semistable.
Hence, the first condition holds.
Let $(E',\phi')\subset (E,\phi)$ be a subobject
such that $\phi'\neq 0$.
We put $E''=E/E'$.
Assume $\mu(E')=\mu(E)$.
Then $E'$ is a member of $\nbigsbar$,
and hence we have
$|P_{E'(C)}(t)-P_{E''(C)}(t)|<\delta_0$,
due to our choice of $\delta_0$.
Therefore, we obtain
$P_{E'(C)}^{\delta_0}(t)>P_{E''(C)}(t)$,
which contradicts the $\delta_0$-semistability of
$\bigl(E(C),\phi_C\bigr)$.
Hence the second condition holds.

Let us assume that the two conditions are satisfied.
Let $(E',\phi')\subset (E,\phi)$ be a subobject
such that $\phi'\neq 0$.
Since we have $\mu(E')<\mu(E)$,
the inequality
$P_{(E',\phi')}^{\delta_0}(t)<P_{(E,\phi)}^{\delta}$
holds for any sufficiently large $t$.
Take a subobject $E'\subset (E,\phi)$.
We have the inequality
$\mu(E')\leq \mu(E)$.
When the strict inequality holds,
we obviously have 
$P_{(E',\phi')}^{\delta_0}(t)<P_{(E,\phi)}^{\delta}$
for any sufficiently large $t$.
Assume $\mu(E')=\mu(E)$.
Then $E'$ is a member of $\nbigsbar$.
We put $E''=E/E'$.
Then, we have 
$P_{E'(C)}(t)<P_{E''(C)}^{\delta_0}$
due to our choice of $\delta_0$.
Thus we obtain the semistability of
$(E(C),\phi_C)$.

From the above argument,
we also obtain that
$P_{(E'(C),\phi'_C)}=P_{(E(C),\phi_C)}$
cannot hold.
Therefore, we obtain the second claim.
\hfill\qed

\begin{lem}
\label{lem;06.6.5.15}
$(E(C),\phi_C)$ is $\delta_0$-semistable,
if and only if 
$(E(-C),\phi_{-C})$ is $\delta_0$-semistable.
Moreover,
the $1$-stability condition holds
for $(y(C),\nbigo(-m,C),\delta_0)$
and $(y(-C),\nbigo(-m,-C),\delta_0)$.
\end{lem}
\pf
By the same argument as the proof of Lemma
\ref{lem;06.6.5.10},
we can show that 
$(E(-C),\phi_{-C})$ is $\delta_0$-semistable
if and only if
the two conditions in Lemma \ref{lem;06.6.5.10} hold.
Then the first claim immediately follows.
The second claim can be shown similarly.
\hfill\qed

\vspace{.1in}

Due to Lemma \ref{lem;06.6.5.15},
we obtain the equality
$\Phitilde_C(\delta_0)=\Phitilde_{-C}(\delta_0)$.
Therefore, we obtain the following equality:
\begin{equation}
\label{eq;06.6.5.20}
 \Phi_{H_+}(\yhat)-\Phi_{H_-}(\yhat)
=\bigl(\Phitilde_{-C}(\delta_0)-\Phi_{H_-}(\yhat)\bigr)
-\bigl(\Phitilde_C(\delta_0)-\Phi_{H_+}(\yhat)\bigr)
\end{equation}

Let us see the first term in the right hand side of
(\ref{eq;06.6.5.20}).
We see the transition of $\Phitilde_{-C}(\delta)$
when we move $\delta$ from $0$ to $\delta_0$.
The transition occurs
when $P_{y_0(-C)}+\delta=P_{y_1(-C)}$ holds
for some $(y_0,y_1)\in S(y,-C)$.
In the case $p_g>0$,
the transitions are trivial.
Hence we obtain
$\Phitilde_{-C}(\delta_0)-\Phi_{H_-}(\yhat)=0$.
In the case $p_g=0$,
we obtain the following equality,
as in the equality (\ref{eq;06.6.5.25}):
\begin{equation}
\Phitilde_{-C}(\delta_0)-\Phi_{H_-}(\yhat)=
 \sum_{(y_0,y_1)\in S(y,-C)}
 \frac{H_{y_0}(m)}{H_y(m)}\cdot
 \int_{\nbigm(\yhat_0,\yhat_1)}
 \underset{s=0}\Res\left(
 \frac{P(\Ehat_0^u\cdot e^{-s}\oplus \Ehat_1^u\cdot e^s)}
 {Q(\Ehat_0^u\cdot e^{-s},\Ehat_1^u\cdot e^s)}
\right)
\end{equation}
Similarly, we obtain
$\Phitilde_C(\delta_0)-\Phi_{H_+}(\yhat)=0$
in the case $p_g>0$,
and we have the following equality in the case $p_g=0$.
\begin{equation}
\Phitilde_{C}(\delta_0)-\Phi_{H_+}(\yhat)=
\sum_{(y_0,y_1)\in S(y,C)}
 \frac{H_{y_0}(m)}{H_{y}(m)}\cdot
 \int_{\nbigm(\yhat_0,\yhat_1)}
 \underset{s=0}\Res\left(
 \frac{P(\Ehat_0^u\cdot e^{-s}\oplus \Ehat_1^u\cdot e^s)}
 {Q(\Ehat_0^u\cdot e^{-s},\Ehat_1^u\cdot e^s)}
\right)
\end{equation}

Now, we have already obtained $\Phi_{H_+}(\yhat)-\Phi_{H_-}(\yhat)=0$
in the case $p_g>0$.
Namely the first claim of Theorem \ref{thm;06.6.4.201}
is proved.
To show the claim in the case $p_g=0$,
we see the sets $S(y,C)$ and $S(y,-C)$ more closely.
We put as follows:
\[
 S_1:=\bigl\{
 (y_0,y_1)\,\big|\,
 y_0+y_1=y,\,\,y_i\in\Type(\nbigsbar),\,\,
 r_0=r_1=1,\,\,a_0=a_1,\,\,
 b_0<b_1
 \bigr\}
\]
We also put 
$ S'(y,\xi):=
 \bigl\{(y_0,y_1)\,\big|\,(y_1,y_0)\in S(y,\xi)\bigr\}$.
Then, it is easy to observe
$S(y,-C)=S(y,\xi)\sqcup S_1$
and $S(y,C)=S'(y,\xi)\sqcup S_1$.
We remark the equality (\ref{eq;06.6.5.30}).
Therefore, we obtain the following equalities:
\begin{multline}
\Phi_{H_+}(\yhat)-\Phi_{H_-}(\yhat)=
 \sum_{(y_0,y_1)\in S(y,\xi)}
 \frac{H_{y_0}(m)}{H_{y}(m)}\cdot
 \int_{\nbigm(\yhat_0,\yhat_1)}
 \underset{s=0}\Res
\left(
 \frac{P(\Ehat_0^u\cdot e^{-s}\oplus \Ehat_1^u\cdot e^s)}
 {Q(\Ehat_0^u\cdot e^{-s},\Ehat_1^u\cdot e^s)}\right)\\
+\sum_{(y_0,y_1)\in S_1}
 \frac{H_{y_0}(m)}{H_{y}(m)}\cdot
 \int_{\nbigm(\yhat_0,\yhat_1)}
 \underset{s=0}\Res\left(
 \frac{P(\Ehat_0^u\cdot e^{-s}\oplus \Ehat_1^u\cdot e^s)}
 {Q(\Ehat_0^u\cdot e^{-s},\Ehat_1^u\cdot e^s)}
 \right)
 \\
+\sum_{(y_0,y_1)\in S(y,\xi)}
 \frac{H_{y_1}(m)}{H_{y}(m)}\cdot
 \int_{\nbigm(\yhat_0,\yhat_1)}
 \underset{s=0}\Res\left(
 \frac{P(\Ehat_0^u\cdot e^{-s}\oplus \Ehat_1^u\cdot e^s)}
 {Q(\Ehat_0^u\cdot e^{-s},\Ehat_1^u\cdot e^s)} 
 \right)\\
-\sum_{(y_0,y_1)\in S_1}
 \frac{H_{y_0}(m)}{H_{y}(m)}\cdot
 \int_{\nbigm(\yhat_0,\yhat_1)}
 \underset{s=0}\Res\left(
 \frac{P(\Ehat_0^u\cdot e^{-s}\oplus \Ehat_1^u\cdot e^s)}
 {Q(\Ehat_0^u\cdot e^{-s},\Ehat_1^u\cdot e^s)}
 \right)\\
=\sum_{(y_0,y_1)\in S(y,\xi)}
 \int_{\nbigm(\yhat_0,\yhat_1)}
 \underset{s=0}\Res\left(
 \frac{P(\Ehat_0^u\cdot e^{-s}\oplus \Ehat_1^u\cdot e^s)}
 {Q(\Ehat_0^u\cdot e^{-s},\Ehat_1^u\cdot e^s)}
 \right)
\end{multline}
Thus we are done.
\hfill\qed

\subsection{Higher Rank Case ($p_g>0$)}
\label{subsection;06.7.3.55}
\subsubsection{Transition formula in the case $p_g>0$}

Let $\vecy$ be an element of $\Type$,
and let $\alpha_{\ast}$ be a system of weights.
Let $L$ be a line bundle on $X$.
Let $\delta\in\nbigp^{\br}$ be a non-critical parameter.
We denote by $\omega$ the first Chern class of 
$\nbigo_{\rel}(1)$ on $\nbigm^s(\vecyhat,[L],\alpha_{\ast},\delta)$.
Let $P$ be an element of $\nbigr$.
When the $1$-stability condition holds
for $(\vecy,L,\alpha_{\ast},\delta)$,
we have $P(\Ehat^u)\cdot \omega^k\in
 \nbigr\bigl(\Ehat^u,
 \nbigm^{s}(\vecyhat,[L],\alpha_{\ast},\delta)\bigr)$,
and we put as follows:
\[
 \Phi(\vecyhat,[L],\alpha_{\ast},\delta):=
 \int_{\nbigm^s(\vecyhat,[L],\alpha_{\ast},\delta)}
 P(\Ehat^u)\cdot \omega^k
\]

Let us discuss the transition formula,
when the $2$-stability condition does not hold
for $(\vecy,L,\alpha_{\ast},\delta)$.
In the case $p_g>0$,
the problem is simpler.
Actually, 
it can be shown that the same formula in 
Theorem \ref{thm;06.5.31.100}
holds.
We use the notation 
in the subsubsection \ref{subsubsection;06.6.6.1}.
We put as follows:
\[
 S_1(\vecy,\alpha_{\ast},\delta):=
\bigl\{
(\vecy_1,\vecy_2)\in S(\vecy,\alpha_{\ast},\delta)\,\big|\,
 \rank\vecy_1=1
\bigr\}
\]
For $(\vecy_1,\vecy_2)\in S_1(\vecy,\alpha_{\ast},\delta)$,
we have $\nbigm(\vecy_1,L)=\nbigm^s(\vecy_1,L,\alpha_{\ast},\delta)$.

\begin{thm}
\label{thm;06.6.6.20}
The following equality holds:
\begin{equation}
\label{eq;06.6.6.11}
 \Phi(\vecyhat,[L],\alpha_{\ast},\delta_+)
-\Phi(\vecyhat,[L],\alpha_{\ast},\delta_-)
=\sum_{(\vecy_1,\vecy_2)\in S_1(\vecy,\alpha_{\ast},\delta)} 
 \int_{\nbigm(\vecy_1,\vecyhat_2,L,\alpha_{\ast},\delta)}
 \Psi(\vecy_1,\vecy_2)
\end{equation}
The elements
$\Psi(\vecy_1,\vecy_2)
\in\nbigr\bigl(E_1^u,
 \nbigm(\vecy_1,L)\bigr)
\otimes
 \nbigr(\Ehat_2^u) $ are given as in {\rm(\ref{eq;06.6.19.26})}:
\[
  \Psi(\vecy_1,\vecy_2)=
\underset{t=0}{\res}
\left(
 \frac{P\bigl(E_1^u\cdot e^{-t}\oplus
 \Ehat_2^u\cdot e^{(t-\omega_1)/(r-1)}\bigr)
 \cdot t^k}
 {\Eu\bigl(\gbign_0(\vecy_1,\vecy_2)\bigr)}
\right)
\]
Here, we put $\omega_1:=c_1\bigl(\Or(E_1^u)\bigr)$.
\end{thm}
\pf
We put $\nbigmtilde(\delta_{\kappa}):=
 \nbigmtilde^{ss}(\vecyhat,[L],\alpha_{\ast},\delta_{\kappa})$.
Let $\Ttilde_{\rel}$ denote the relative tangent
bundle of the smooth morphism
of $\nbigmtilde(m,\vecyhat,[L])$
to $\nbigm(m,\vecyhat,[L])$.
We have the open embedding
$\nbigmtilde(\delta,\kappa)\subset \nbigmtilde(m,\vecyhat,[L])$.
The restriction of $\Ttilde_{\rel}$
is denoted by the same notation.

Let $\Mhat$ be the master space
connecting $\nbigmtilde(\delta_+)$
and $\nbigmtilde(\delta_-)$
constructed in the subsubsection \ref{subsubsection;06.5.17.10}.
Let $\varphi:\Mhat\lrarr\nbigmtilde(m,\vecyhat,[L])$ denote
the naturally defined morphism.
Let $\nbigt(1)$ denote the trivial line bundle
on $\nbigmtilde(m,\vecyhat,[L])$ with the $G_m$-action of weight $1$.
We consider the following element of
$\nbigr_{G_m}\bigl(\Ehat^{\Mhat},\Mhat\bigr)$:
\[
 \Phitilde_t:=
 P\bigl(\varphi^{\ast}\Ehat^u\bigr)\cdot 
  c_1\bigl(\varphi^{\ast}\nbigo_{\rel}(1)\bigr)^k
\cdot \frac{\Eu(\Ttilde_{\rel})}{H_y(m)!},
\quad\quad
 \Phibar_t:=\Phitilde_t
\cdot c_1\bigl(\varphi^{\ast}\nbigt(1)\bigr).
\]
We use Proposition \ref{prop;06.6.10.2}.
Then, we obtain the polynomial 
$\int_{\Mhat}\Phibar_t$ of $t$.
When we forget the $G_m$-action,
we have 
$c_1\bigl(\varphi^{\ast}\nbigt(1)\bigr)=0$.
Hence we have  
$\int_{\Mhat}\Phibar_{t|t=0}=0$.
On the other hand,
we have the following equality in $\rnum[t^{-1},t]$,
due to the localization of the virtual fundamental classes
(\cite{gp}):
\[
\int_{\Mhat}
 \Phibar_t=
 \sum_{i=1,2}
\int_{\Mhat_i}
\frac{\Phibar_t}{\Eu\bigl(\gbign(\Mhat_i)\bigr)}
+\sum_{\gbigi\in \Dec(m,\vecy,\alpha_{\ast},\delta)}
 \int_{\Mhat^{G_m}(\gbigi)}
 \frac{\Phibar_t}{\Eu\bigl(\gbign(\Mhat^{G_m}(\gbigi))\bigr)}
\]
Here, $\Dec(m,\vecy,\alpha_{\ast},\delta)$ denotes 
the set of the decomposition types
(Definition \ref{df;06.5.22.1}),
and $\gbign(\Mhat_i)$ and $\gbign(\Mhat^{G_m}(\gbigi))$
denote the virtual normal bundles with the $G_m$-action
given in Proposition \ref{prop;06.5.23.40}
and Proposition \ref{prop;06.5.23.41}.
We have
$c_1\bigl(\varphi^{\ast}\nbigt(1)\bigr)_{|\Mhat_i}=t$
and $c_1\bigl(\varphi^{\ast}\nbigt(1)\bigr)_{|\Mhat^{G_m}(\gbigi)}=t$.
Therefore, we obtain the following equality:
\begin{equation}
\label{eq;06.6.10.100}
\sum_{i=1,2}
\int_{\Mhat_i}
\underset{t=0}{\Res}
\left(
 \frac{\Phitilde_t}{\Eu\bigl(\gbign(\Mhat_i)\bigr)}
\right)
+\sum_{\gbigi}
\int_{\Mhat^{G_m}(\gbigi)}
 \underset{t=0}{\Res}
\left(
 \frac{\Phitilde_{t}}{\Eu\bigl(\gbign(\Mhat^{G_m}(\gbigi))\bigr)}
\right)
=0.
\end{equation}
As in the proof of Theorem \ref{thm;06.5.31.100},
the first term of the left hand side of 
(\ref{eq;06.6.10.100}) can be rewritten as follows:
\begin{equation}
\label{eq;06.6.10.101}
 -\int_{\Mhat_1}\Phitilde
+\int_{\Mhat_2}\Phitilde
=-\int_{\nbigmtilde(\delta_+)}\Phitilde
 +\int_{\nbigmtilde(\delta_-)}\Phitilde
\end{equation}
Recall $\nbigmtilde(\delta_{\kappa})$ is
the full flag bundle over
$\nbigm^{ss}\bigl(\vecyhat,[L],\alpha_{\ast},\delta_{\kappa}\bigr)$,
and we have the equality of the virtual
fundamental classes as in Lemma \ref{lem;06.6.13.150}.
Hence, we obtain the following equality:
\begin{equation}
\sum_{i=1,2}
\int_{\Mhat_i}
\underset{t=0}{\Res}
\left(
 \frac{\Phitilde_t}{\Eu\bigl(\gbign(\Mhat_i)\bigr)}
\right)
= -\Phi(\vecyhat,[L],\alpha_{\ast},\delta_+)
+\Phi(\vecyhat,[L],\alpha_{\ast},\delta_-)
\end{equation}

The contributions from $\Mhat^{G_m}(\gbigi)$
can be calculated by the arguments
in the proof of Theorem \ref{thm;06.5.31.100}
and Theorem \ref{thm;06.6.4.10}.
For any decomposition type
$\gbigi=(\vecy_1,\vecy_2,I_1,I_2)\in
 \Dec(m,\vecy,\alpha_{\ast},\delta)$,
we put as follows:
\[
 \nbigmtilde(\gbigi):=
 \nbigmtilde^{ss}\bigl(\vecy_1,L,\alpha_{\ast},\delta,\gminik(\gbigi)\bigr)
\times
 \nbigmtilde^{ss}(\vecyhat_2,\alpha_{\ast},+)
\]
Let $\Ttilde_{1,\rel}$ denote the vector bundle
over $\nbigmtilde(\gbigi)$ 
induced by the relative tangent bundle of
the smooth map
$\nbigmtilde^{ss}(\vecy_1,L,\alpha_{\ast},\delta,\gminik(\gbigi))
\lrarr
 \nbigm(m,\vecy_1,L)$.
Let $\Ttilde_{2,\rel}$ denote the vector bundle
over $\nbigmtilde(\gbigi)$ 
induced by the relative tangent bundle of
the smooth map
$\nbigmtilde^{ss}(\vecyhat_2,\alpha_{\ast},+)
\lrarr \nbigm(m,\vecyhat_2)$.
Let $N_0$ be as in the subsubsection 
\ref{subsubsection;06.5.31.11}.
Then, we have the following decomposition of the vector bundles:
\[
 \varphi_{\gbigi}^{\ast}
\Ttilde_{\rel}
=\Ttilde_{1,\rel}\oplus \Ttilde_{2,\rel}\oplus N_0.
\]
We remark that
$\varphi^{\ast}\nbigo_{\rel}(1)_{|\Mhat^{\ast}}$
and $\varphi^{\ast}\nbigt(1)_{|\Mhat^{\ast}}$
are naturally isomorphic as $G_m$-equivariant line bundles.
We use the relation of $E_i^{\Mhat}$,
$E_1^u$ and $\Ehat_2^u$
in Corollary \ref{cor;06.5.22.50}.
Then, we have the following equality 
in $\nbigr(G^{\prime\,\ast}E_1^u,
 G^{\prime\,\ast}\Ehat_2^u,\nbigs)[t]$:
\[
 F^{\ast}\Phitilde_t=
 \frac{H_{y_1}(m)!H_{y_2}(m)!}{H_y(m)!}
\cdot
G^{\prime\,\ast}
 P\bigl(E_1^u\cdot e^{-t}\oplus
 \Ehat_2^u\cdot e^{r_1(t-\omega_1)/r_2}\bigr)\cdot t^k
 \cdot
 \frac{\Eu(\Ttilde_{1,\rel})}{H_{y_1}(m)!}
 \frac{\Eu(\Ttilde_{2,\rel})}{H_{y_2}(m)!}
 \cdot \Eu(N_0)
\]
We also have the following equality
of the equivariant Euler classes
in $A^{\ast}(\nbigs)[[t^{-1},t]$:
\[
 F^{\ast}\Eu\bigl(\gbign(\Mhat^{G_m}(\gbigi))\bigr)
=G^{\prime\,\ast}\Eu\bigl(\gbign_0(\vecy_1,\vecy_2)\bigr)
 \cdot \Eu(N_0)
\]
Recall that we have the equality of the virtual fundamental classes
in Proposition \ref{prop;06.5.11.300}.
Therefore, the contribution from $\Mhat^{G_m}(\gbigi)$
is as follows:
\begin{equation}
\label{eq;06.6.6.5}
 \frac{H_{y_1}(m)!\cdot H_{y_2}(m)!}{H_y(m)!}
 \int_{\nbigmtilde(\gbigi)}
 \underset{t=0}\Res
 \left(
 \frac{P\bigl(E_1^u\cdot e^{-t}\oplus
 \Ehat_2^u\cdot e^{r_1(t-\omega_1)/r_2}\bigr)\cdot t^k}
 {\Eu\bigl(\gbign_0(\vecy_1,\vecy_2)\bigr)}
 \right)
\frac{\Eu(\Ttilde_{1,\rel})}{H_{y_1}(m)!}
\frac{\Eu(\Ttilde_{2,\rel})}{H_{y_2}(m)!}
\end{equation}
We remark that
the virtual fundamental class of
$\nbigmtilde(\gbigi)$ is $0$,
and hence (\ref{eq;06.6.6.5}) vanishes,
if the conditions $\rank(\vecy_1)>1$ and $p_g>0$
are satisfied
(Proposition \ref{prop;06.5.31.1}).
In the case $\rank(\vecy_1)=1$,
the $(\delta,\ell)$-semistability condition is trivial.
We also remark that
the integrand of (\ref{eq;06.6.6.5}) is 
the element of
$\nbigr\bigl(E_1^u,
 \nbigmtilde^{ss}(\vecy_1,L,\alpha_{\ast},\delta,\gminik(\gbigi))\bigr)
\otimes
 \nbigr\bigl(\Ehat_2^u,\nbigmtilde^{ss}(\vecyhat_2,\alpha_{\ast},+)\bigr)$,
and hence we have only to consider the component-wise
integration.
By using Lemma \ref{lem;06.5.31.105},
we can rewrite (\ref{eq;06.6.6.5}) as follows:
\begin{equation}
 \frac{H_{y_1}(m)!\cdot H_{y_2}(m)!}{H_y(m)!}
 \int_{\nbigm(\vecy_1,\vecyhat_2,L,\alpha_{\ast},\delta)}
 \underset{t=0}\Res
 \left(
 \frac{P\bigl(E_1^u\cdot e^{-t}\oplus
 \Ehat_2^u\cdot e^{(t-\omega_1)/(r-1)}\bigr)\cdot t^k}
 {\Eu\bigl(\gbign_0(\vecy_1,\vecy_2)\bigr)}
 \right)
\end{equation}
The number of the decompositions
$(I_1,I_2)$ of $\{1,\ldots,H_y(m)\}$
satisfying $|I_i|=H_{y_i}(m)$
is $H_y(m)!\cdot H_{y_1}(m)!^{-1}\cdot H_{y_2}(m)!^{-1}$.
Therefore, the second term in the left hand side of
(\ref{eq;06.6.10.100})
is same as the left hand side of (\ref{eq;06.6.6.11}).
Thus, we obtain the desired formula.
\hfill\qed

\subsubsection{Reduction to the integrals over the products of Hilbert
   schemes}
\label{subsubsection;06.6.30.2}

We assume $p_g>0$ and $\dim H^1(X,\nbigo)=0$.
Let $y$ be an element of $\Type_r^{\circ}$.
Assume 
$P_y(t)>P_K(t)$ for any sufficiently large $t$,
where $K$ denotes the canonical line bundle of $X$.
We also assume $\chi(y)-1\geq (r-1)\cdot (1+p_g)$.
We give a straightforward generalization of 
Theorem \ref{thm;06.6.6.10}.
For a given element $y_i\in \Type$,
we use the notation $r_i$, $a_i$, $b_i,$ and $n_i$
to denote the rank, the first Chern class,
the second Chern character and the second Chern class,
in the following argument.

We put as follows:
\[
 T_i=\sum_{j<i}\frac{t_j}{r-j}-t_i,\,\,\,
(i=1,\ldots,r-1),\quad
 T_{r}=\sum_{j<r}\frac{t_j}{r-j}
\]
We put as follows:
\[
 S(y)=\Bigl\{
(y_1,y_2,\ldots,y_r)\in\Type_1^{\circ\,r}
 \,\Big|\,
 \sum y_i=y,\,\,
 H_{y_i}<(r-i)^{-1}\sum_{j>i} H_{y_{j}},\,\,
 \chi(a_i)=1+p_g\,\,(i<r)
 \Bigr\}
\]
Let $(y_1,\ldots,y_r)$ be an element of $S(y)$.
We put $X^{[\vecn]}:=\prod_{i=1}^{r} X^{[n_i]}$.
The universal ideal sheaves over $X^{[n_i]}\times X$
are denoted by $\nbigi_i^u$.
Let $\nbigz_i$ denote the universal scheme 
over $X^{[n_i]}\times X$
of length $n_i$.
We also use the same notation
to denote the pull back of $\nbigi_i^u$
and $\nbigz_i$ via the projection
$X^{[\vecn]}\times X\lrarr
 X^{[n_i]}\times X$.
Let $e^{a_i}$ denote the holomorphic line bundle
corresponding to $a_i$.
Since we have assumed $H^1(X,\nbigo_X)=0$
in this subsubsection, such a holomorphic line
bundle is uniquely determined up to isomorphisms.
Let $G$ denote the $(r-1)$-dimensional torus
$\Spec k[\tau_1,\tau_{1}^{-1},\ldots,\tau_{r-1},\tau_{r-1}^{-1}]$.
Let $e^{w\cdot t_i}$ denote the trivial line bundle
with the $G$-action,
which is induced by the $\Spec k[\tau_i,\tau_i^{-1}]$-action
of weight $w$.
Thus, we obtain the $G$-equivariant sheaf
$\nbigi_i^u\cdot e^{a_i+T_i}$.
We use the notation 
$Q(\nbigi_i^u\cdot e^{a_i+T_i},\nbigi_j^u\cdot e^{a_j+T_j})$
to denote the $G$-equivariant Euler class
of the following:
\[
 -Rp_{X\,\ast}\Bigl(
 \nrhom\bigl(
 \nbigi_i^u\cdot e^{a_i+T_i},\,\,
 \nbigi_j^u\cdot e^{a_j+T_j} \bigr)\Bigr)
-Rp_{X\,\ast}\Bigl(
 \nrhom\bigl(
 \nbigi_j^u\cdot e^{a_j+T_j},\,\,
 \nbigi_i^u\cdot e^{a_i+T_i}
 \bigr)\Bigr)
\]
We regard it as the element of 
$\bigotimes_{i=1}^{r}\nbigr\bigl(\nbigi_i^u\cdot e^{a_i}\bigr)
\otimes_{\rnum} \gbigr(t_1,\ldots,t_{r-1})$.
(See the subsubsection \ref{subsubsection;06.6.19.50}
 for the ring $\gbigr(t_1,\ldots,t_{r-1})$.)
We put as follows:
\[
 Q\bigl(
 \nbigi_1^u\cdot e^{a_1+T_1},
 \nbigi_2^u\cdot e^{a_2+T_2},
 \ldots
 \nbigi_{r}^u\cdot e^{a_r+T_r}
 \bigr)
=\prod_{i<j}
 Q\bigl(\nbigi_i^u\cdot e^{a_i+T_i},\,\,
 \nbigi_j^u\cdot e^{a_j+T_j}
 \bigr)
\]
We also have the element
$P\bigl(\bigoplus_{i=1}^{r}\nbigi_i^u\cdot e^{a_i+T_i}\bigr)$
of $\bigotimes_{i=1}^{r}
 \nbigr(\nbigi_i^u\cdot e^{a_i})[t_1,\ldots,t_{r-1}]$
by the homomorphisms in the subsubsection
 \ref{subsubsection;06.6.19.25}.
Let $\Xi_i$ denote the vector bundle
$p_{X\,\ast}\bigl(\nbigo_{\nbigz_i}\otimes e^{a_i}\bigr)$.

\begin{thm}
Assume $p_g>0$ and $H^1(X,\nbigo)=0$.
We also assume $\chi(y)-1\geq (r-1)\cdot (1+p_g)$
and $P_y>P_K$.
We have the following formula:
\[
\int_{\nbigm(\yhat)}P(\Ehat^u)
=(-1)^{r-1}
 \sum_{(y_1,\ldots,y_r)\in S(y)}
  \prod_{i=1}^{r-1}\SW(a_i)\cdot
\int_{X^{[\vecn]}}
 \Psi(y_1,\ldots,y_r)
\]
The elements
$\Psi(y_1,\ldots,y_r)\in\bigotimes_{i=1}^r
 \nbigr\bigl(\nbigi_i^u\cdot e^{a_i}\bigr)$
are given as follows:
\begin{multline}
\label{eq;06.5.31.121}
 \Psi(y_1,\ldots,y_r)=\\
  \underset{t_{r-1}}\Res
 \cdots
 \underset{t_{1}}\Res
 \left(
 \frac{P\bigl(\bigoplus_{i=1}^{r}\nbigi_i^u\cdot e^{a_i+T_i} \bigr)}
 {Q\bigl(\nbigi_1^u\cdot e^{a_1+T_1},\ldots,
 \nbigi_r^{u}\cdot e^{a_r+T_r}\bigr)}
 \prod_{i<j}\frac{\Eu(\Xi_j\cdot e^{T_j-T_i})}
 {(T_j-T_i)^{\chi(a_j)}}
\cdot \prod_{i=1}^{r-1}t_i^{\sum_{j\geq i} \chi(y_j)-1}
\right)
 \cdot \prod_{i=1}^{r-1}\Eu(\Xi_i)
\end{multline}
\end{thm}
\pf
We put as follows:
\[
 S_1(y):=\bigl\{
 (y_1,y_2)\in\Type_1\times\Type_{r-1}\,\big|\,
 y_1+y_2=y,\,\,
 P_{y_1}<P_{y_2}
 \bigr\}
\]
By using the transition formula (\ref{eq;06.6.6.11})
and the same argument as the proof of
Theorem \ref{thm;06.6.6.10},
we obtain the following equality:
\begin{equation}
\label{eq;06.6.6.17}
\int_{\nbigm(\yhat)}\!\! P(\Ehat^u)=
-\!\!\!\!\!\sum_{(y_1,y_2)\in S_1(y)}\!\!
 \SW(a_1)\cdot
 \int_{X^{[n_1]}\times\nbigm(\yhat_2)}
 \underset{t=0}\Res
 \left(
 \frac{P\bigl(\nbigi_1^u\cdot e^{a_1-t}\oplus
 \Ehat_2^u\cdot e^{t/(r-1)}\bigr)\cdot t^{\chi(y)-1}}
 {Q(\nbigi_1^u\cdot e^{-t},\,
 \Ehat_2^u\cdot e^{t/(r-1)})
 \cdot R\bigl(\nbigo_X\cdot e^{-t},\,
 \Ehat_2^u\cdot e^{t/(r-1)}\bigr)}
 \right)
\end{equation}
Here, $R\bigl(\nbigo_X\cdot e^{-t},\,\Ehat_2^u\cdot e^{t/(r-1)}\bigr)$
denotes the equivariant Euler class of
$p_{X\,\ast}\bigl(
\nhom(\nbigo_X\cdot e^{-t},\,\Ehat_2^u\cdot e^{t/(r-1)})\bigr)$.

We remark that $\chi(y_2)-1\geq 0$
and $P_{y_2}>P_y>P_{K}$ and 
$\chi(y_2)-1\geq (r-2)(1+p_g)$.
We also remark that the integrand of (\ref{eq;06.6.6.17})
is the element of 
$\nbigr\bigl(\nbigi_1\cdot e^{a_1}\bigr)\otimes
 \nbigr\bigl(\Ehat_2\bigr)$,
and hence we have only to consider the component-wise
integration.
Thus, we may use (\ref{eq;06.6.6.17})
inductively.
We put 
$\Tbar_i:=\sum_{j<i}t_j/(r-j)$.
Then,
we obtain the following equality:
\begin{equation}
 \int_{\nbigm(\yhat)}P(\Ehat^u)
+(-1)^r\sum_{(y_1,\ldots,y_r)\in S(y)}
 \prod_{i=1}^{r-1}\SW(a_i)\cdot
 \int_{X^{[\vecn]}}\Psi_1(y_1,\ldots,y_r)
=0
\end{equation}
The elements
$\Psi_1(y_1,\ldots,y_r)\in
 \bigotimes_{i=1}^{r}\nbigr(\nbigi_i^u\cdot e^{a_i})$ 
are given as follows:
\[
 \Psi_1(y_1,\ldots,y_r)=
 \underset{t_{r-1}}\Res\cdots
 \underset{t_1}\Res
 \left(
\frac{
 P\bigl(\bigoplus \nbigi_i^u\cdot e^{a_i+T_i}\bigr)
 \cdot \prod_{i=1}^{r-1}t_i^{\sum_{j\geq i}\chi(y_j)-1}
 \prod_{i=1}^{r-1}\Eu(\Xi_i)}
{\prod_{i<j}
 Q\bigl(\nbigi_i^u\cdot e^{a_i-t_i},
 \nbigi_j^u\cdot e^{a_j+T_j-\Tbar_i}\bigr)
\cdot
 R\bigl(\nbigo_X\cdot e^{-t_i},\,
 \nbigi_j^u\cdot e^{a_j+T_j-\Tbar_i}\bigr)}
 \right)
\]
We have
$Q\bigl(\nbigi_i^u\cdot e^{a_i-t_i},
 \,\nbigi_j^u\cdot e^{a_j+T_j-\Tbar_i}\bigr)
=Q\bigl(\nbigi_i^u\cdot e^{a_i+T_i},\,
 \nbigi_j^u\cdot e^{a_j+T_j}\bigr)$.
We also obtain the following by a formal calculation:
\[
 R\bigl(\nbigo_X\cdot e^{-t_i},\,\nbigi_j^u\cdot e^{T_j-\Tbar_i}\bigr)
=\frac{(T_j-T_i)^{\chi(a_j)}}
{\Eu\bigl(\Xi_j\cdot e^{T_j-T_i}\bigr)}
\]
Hence, we obtain $\Psi_1(y_1,\ldots,y_r)=\Psi(y_1,\ldots,y_r)$.
Thus we are done.
\hfill\qed

\subsubsection{Independence from the polarization
 in the case $p_g>0$}

\label{subsubsection;06.6.6.121}

Let $y$ be an element of $\Type_r^{\circ}$.
We use the notation $\nbigm_H(\yhat)$ 
to denote the moduli stack of torsion-free sheaves
of type $y$, which are semistable with respect to 
a polarization $H$.
For $\Phi=P(\Ehat^u)\in\nbigr\bigl(\Ehat,\nbigm_H(\yhat)\bigr)$,
we put as follows:
\[
 \Phi_H(\yhat):=\int_{\nbigm_H(\yhat)}\Phi
\]

\begin{thm}
\label{thm;06.6.8.150}
The invariant $\Phi_H(\yhat)$
is independent of the choice of a generic polarization
in the case $p_g>0$.
\end{thm}
\pf
Let $\xi$ be an element of $NS(X)$
such that $\xi^2< 0$.
Let $W^{\xi}$ be the wall
in the ample cone determined by $\xi$.
It is well known that the ample cone is
divided into the chambers
by such walls,
and $\nbigm_{H}(\yhat)$ depends only 
on the chamber to which $H$ belongs.
Moreover, the set of such walls are locally finite
(\cite{mw}).

Let $C_+$ and $C_-$ be chambers
which are divided by a wall $W^{\xi}$.
Let $H_+$ and $H_-$ be ample line bundles
contained in $H_+$ and $H_-$ respectively.
We assume $(\xi,H_+)>0>(\xi,H_-)$.
For the proof of the theorem,
we have only to show $\Phi_{H_+}(\yhat)=\Phi_{H_-}(\yhat)$.
We take an ample line bundle $H$
which is generic in $W^{\xi}$.
As in Lemma \ref{lem;06.6.4.202},
we may assume the following:
\begin{itemize}
\item
$H_+=H\otimes \nbigo(C)$ 
and $H_-=H\otimes \nbigo(-C)$
for some sufficiently ample divisor $C$.
\item
A torsion free sheaf $E$ of type $y$
is $H_+$-semistable
if and only if $E(C)$ is $H$-semistable.
\item
$E$ is $H_-$ -semistable
if and only if $E(-C)$ is $H$-semistable.
\end{itemize}

We use the notation and the argument
in the subsubsection \ref{subsubsection;06.6.4.220}.
Let $\nbigs$ be the family of
torsion-free sheaves of type $y$,
which are $\mu$-semistable.
Let $\nbigsbar$ be the family of
torsion-free sheaves $E'$
with the following properties:
\begin{itemize}
\item $\mu(E')=\mu(y)$ and $\rank(E')<r$.
\item
There exists a member $E$ of $\nbigs$
such that $E'$ is a saturated subsheaf of $E$.
\end{itemize}

We take a large integer $m$ 
such that the condition $O_m$ holds for the family $\nbigs$.
As in the subsubsection \ref{subsubsection;06.6.4.220},
we consider the integrals over
$\nbigm\bigl(\yhat(C),[\nbigo(-m,C)],\delta\bigr)$
and $\nbigm\bigl(\yhat(-C),\nbigo(-m,-C),\delta\bigr)$
given by (\ref{eq;06.6.6.16}).
When $\delta$ is sufficiently small,
we have 
$\Phitilde_C(\delta)=\Phi_{H_+}(\yhat)$
and $\Phitilde_{-C}(\delta)=\Phi_{H_-}(\yhat)$.

Let $\Type(\nbigsbar)$ denote 
the set of the types of members of $\nbigsbar$.
For each $y_0\in \Type(\nbigsbar)$,
we have $y_1:=y-y_0\in\Type(\nbigsbar)$.
We remark that
$P_{y_0(C)}-P_{y_1(C)}$
and $P_{y_0(-C)}-P_{y_1(-C)}$ are polynomials
of degree $0$.

Let $r_i$, $a_i$, $b_i$ and $n_i$ denote 
the rank, the first Chern class, the second Chern character
and the second Chern class corresponding to $y_i$.
We have $\mu(y_0)=\mu(y_1)$,
and $H$ is a generic element of $W^{\xi}$.
Therefore, we have $a_0/r_0-a_1/r_1=A\cdot \xi$
for some $A\in\rnum$ in $H^2(X)$.
Since the intersection number $(C,\xi)$ is
assumed to be sufficiently large,
we have $P_{y_0(C)}-P_{y_1(C)}\neq 0$
unless $y_0/r_0=y_1/r_1$.
We also have $P_{y_0(-C)}-P_{y_1(-C)}\neq 0$
unless $y_0/r_0=y_1/r_1$.

We put as follows:
\[
 S(y,C):=
 \bigl\{(y_0,y_1)\,\big|\,
 y_0+y_1=y,\,\,y_i\in\Type(\nbigsbar),\,\,
 P_{y_0(C)}<P_{y_1(C)}
 \bigr\}
\]
\[
   S(y,-C):=
 \bigl\{(y_0,y_1)\,\big|\,
 y_0+y_1=y,\,\,y_i\in\Type(\nbigsbar),\,\,
 P_{y_0(-C)}<P_{y_1(-C)}
 \bigr\}.
\]

We take a positive constant $\delta_0$
satisfying the following inequalities:
\[
 \delta_0>
 \max\bigl\{
 |P_{y_0(C)}-P_{y_1(C)}|\,
\big|\,\,
 (y_0,y_1)\in S(y,C)
 \bigr\},
\quad
 \delta_0>
\max\bigl\{
 |P_{y_0(-C)}-P_{y_1(-C)}|\,
\big|\,\,
 (y_0,y_1)\in S(y,-C)
 \bigr\}.
\]
The following lemma can be shown 
by the same argument as the proof of 
Lemma \ref{lem;06.6.5.10} and
Lemma \ref{lem;06.6.5.15}.
\begin{lem}
\label{lem;06.6.11.55}
$(E(C),\phi_C)$ is $\delta_0$-semistable,
if and only if 
$(E(-C),\phi_{-C})$ is $\delta_0$-semistable.
Moreover, the $\delta_0$-semistability
implies the $\delta_0$-stability in the both cases.
\hfill\qed
\end{lem}

Due to the lemma,
we obtain the following:
\begin{equation}
 \Phi_{H_+}(\yhat)-\Phi_{H_-}(\yhat)
=\bigl(\Phitilde_{-C}(\delta_0)-\Phi_{H_-}(\yhat)\bigr)
-\bigl(\Phitilde_{C}(\delta_0)-\Phi_{H_+}(\yhat)\bigr)
\end{equation}
Let us move the parameter $\delta$
from $0$ to $\delta_0$,
and we see 
the transition of the invariants
$\Phitilde_{-C}(\delta)$ and $\Phitilde_C(\delta)$.
We use Theorem \ref{thm;06.6.6.20}.
Since the condition $O_m$ holds for $\nbigs$,
it is easy to see that
the contributions from 
the decomposition types $(y_0,y_1)$
are trivial even in the case $\rank (y_0)=1$
(Proposition \ref{prop;06.5.11.500}).
Hence we obtain 
$\Phitilde_{-C}(\delta_0)-\Phi_{H_-}(\yhat)=0$
and $\Phitilde_C(\delta_0)-\Phi_{H_+}(\yhat)=0$.
Thus, the proof of Theorem \ref{thm;06.6.8.150}
is finished.
\hfill\qed

\subsection{Transition Formula in the Case $p_g=0$}
\label{subsection;06.6.21.30}
\subsubsection{Statement}
\label{subsubsection;06.6.23.15}

Let us discuss the transition formula
when the $2$-stability condition does not hold for 
$(\vecy,L,\alpha_{\ast},\delta)$
in the case $p_g=0$.
Let $P$ be an element of $\nbigr$.
We restrict ourselves to the case
where the $1$-vanishing condition holds
for $(\vecy,L,\alpha_{\ast},\delta)$,
and we discuss the integral of the element
$\nbigr\bigl(\Ehat^u,\nbigm^s(\vecyhat,[L],\alpha_{\ast},\delta)\bigr)$
of the following form:
\begin{equation}
\label{eq;06.6.11.1}
 \Phi=P(\Ehat^u)\cdot\frac{\Eu(T_{\rel})}{N_L(y)}
\end{equation}
We put as follows, for non critical parameter $\delta$:
\[
 \Phi(\delta):=\int_{\nbigm^{s}(\vecyhat,[L],\alpha_{\ast},\delta)}\Phi.
\]

Let $\delta$ be a critical parameter.
We take $\delta_-<\delta<\delta_+$ 
such that $\delta_{\kappa}$ $(\kappa=\pm)$
are sufficiently close to $\delta$.
We would like to obtain the formula
$\Phi(\delta_+)-\Phi(\delta_-)$
to express the sum as the integrals
over the products of the moduli stacks
of the objects with lower ranks.
We also impose the following condition
to $(\vecy,L,\alpha_{\ast},\delta)$:
\begin{condition}
\label{condition;06.6.6.100}
 For any $(E_{1\ast},\phi)\oplus E_{2\,\ast}\oplus E_{3\,\ast}\in
 \nbigm^{ss}(\vecy,L,\alpha_{\ast},\delta)$,
 the equality $P_{E_2}=P_{E_3}$ holds  
 for the reduced Hilbert polynomials of $E_2$ and $E_3$.
\hfill\qed
\end{condition}

For each positive integer $k$
we put as follows:
\[
 S_k(\vecy,\delta):=
 \bigl\{
 \vecY=(\vecy_1,\ldots,\vecy_k)\in\Type^k\,\big|\,
 P^{\alpha_{\ast}}_{\vecy_i}
=P^{\alpha_{\ast},\delta}_{\vecy}\,\,\,
 (i=1,\ldots,k)
 \bigr\}
\]
For each element
$\vecY=(\vecy_1,\ldots,\vecy_k)\in S_k(\vecy,\delta)$,
we put $|\vecY|=\sum_{i=1}^k\vecy_i$.
We also put as follows:
\[
 W(\vecY):=\prod_{i=1}^k
 \frac{\rank(\vecy_i)}
 {\sum_{1\leq j\leq i}\rank(\vecy_j)}
\]
We put as follows:
\[
 \Sbar_k(\vecy,\delta):=
 \bigl\{
 (\vecy_0,\vecY)\in\Type\times S_k(\vecy,\delta)
 \,\big|\,\vecy_0+|\vecY|=\vecy
 \bigr\},
\quad
 \Sbar(\vecy,\delta):=
 \coprod_{k}\Sbar_k(\vecy,\delta)
\]
For any $(\vecy_0,\vecY)\in \Sbar_k(\vecy,\delta)$,
we put as follows:
\[
 \nbigm(\vecyhat_0,\vecYhat,[L]):=
 \nbigm^s\bigl(\vecyhat_0,[L],\alpha_{\ast},\delta_-\bigr)
\times
 \prod_{i=1}^k\nbigm^{ss}(\vecyhat_i,\alpha_{\ast})
\]
Let $\Ehat_0^u$ denote the sheaf
over $\nbigm(\vecyhat_0,\vecYhat,[L])\times X$
which is obtained as the pull back of
the universal sheaf over
$\nbigm^s(\vecyhat_0,[L],\alpha_{\ast},\delta_-)\times X$
via the natural projection.
We use the notation $\Ehat_i^u$ in a similar meaning.

When $(\vecy_0,\vecY)\in \Sbar_k(\vecy,\delta)$ is given,
we put as follows:
\[
 T_0=-\sum_{j>0}\frac{t_j}{\sum_{0\leq h<j}\rank(\vecy_h)},
\quad
 T_i=-\sum_{j>i}\frac{t_j}{\sum_{0\leq h<j}\rank(\vecy_h)}
+\frac{t_i}{\rank(\vecy_i)}
\]
Here, $t_1,\ldots,t_k$ are variables.
Let $G$ denote the $k$-dimensional torus
$\Spec k[\tau_1,\tau_1^{-1},\ldots,\tau_k,\tau_k^{-1}]$.
Let $e^{w\cdot t_i}$ denote the trivial line bundle
with the $G$-action
which is induced by the $\Spec k[\tau_i,\tau_i^{-1}]$-action
of weight $w$.
We have the following element of 
the $K$-group of the $G$-equivariant
coherent sheaves on $\nbigm(\vecyhat_0,\vecYhat,[L])$:
\begin{multline}
-Rp_{X\,\ast}\nrhom\bigl(\Ehat^u_i\cdot e^{T_i},\,\Ehat^u_j\cdot e^{T_j}\bigr)
-Rp_{X\,\ast}\nrhom\bigl(\Ehat^u_j\cdot e^{T_j},\,\Ehat^u_i\cdot e^{T_i}\bigr) \\
-Rp_{D\,\ast}\nrhom'_2\bigl(\Ehat^u_{i|D\,\ast}\cdot e^{T_i},\,
 \Ehat^u_{j|D\,\ast}\cdot e^{T_j}\bigr)
-Rp_{D\,\ast}\nrhom'_2\bigl(\Ehat^u_{j|D\,\ast}\cdot e^{T_j},\,
 \Ehat^u_{i|D\,\ast}\cdot e^{T_i}\bigr)
\end{multline}
(See the subsubsection \ref{subsubsection;06.5.23.10}
for the notation $\nrhom'_2$.)
The equivariant Euler class is denoted by
$Q\bigl(\Ehat^u_i\cdot e^{T_i},\Ehat^u_j\cdot e^{T_j}\bigr)$,
which is regarded as the element of
$\bigotimes_{i=0}^k\nbigr(\Ehat^u_i)
 \otimes\gbigr(t_k,t_{k-1},\ldots,t_2,t_{1})$.
We put as follows:
\begin{equation}
 \label{eq;06.6.11.20}
 Q\bigl(\Ehat^u_0\cdot e^{T_0},\Ehat^u_1\cdot e^{T_1},\ldots,
 \Ehat^u_k\cdot e^{T_i}\bigr):=
 \prod_{i<j}
 Q\bigl(\Ehat^u_i\cdot e^{T_i},\Ehat^u_j\cdot e^{T_j}\bigr)
\end{equation}
Let $T_{0\,\rel}$ denote the vector bundle
over $\nbigm(\vecyhat_0,\vecYhat,[L])$
obtained from the relative tangent bundle
of the smooth map 
$\nbigm^s(\vecyhat_0,[L],\alpha_{\ast},\delta_-)\lrarr 
 \nbigm(m,\vecyhat_0,[L])$.
Then, we have the following element:
\[
 \Psi(\vecy_0,\vecY):=
 \underset{t_1=0}\Res\cdots
 \underset{t_k=0}\Res
\left(
\frac{P\bigl(\bigoplus_{i=0}^k\Ehat_i^u\cdot e^{T_i}\bigr)}
 {Q\bigl(\Ehat^u_0\cdot e^{T_0},\ldots,\Ehat^u_k\cdot e^{T_k}\bigr)}
\right)
 \cdot \frac{\Eu(T_{0,\rel})}{N_L(y_0)}
\in
\nbigr\bigl(\Ehat_0,\nbigm^s(\vecyhat_0,[L],\alpha_{\ast},\delta_-)\bigr)
\otimes
\bigotimes_{i=1}^k \nbigr\bigl(\Ehat_i\bigr)
\]

\begin{thm}
\label{thm;06.6.6.30}
Assume that the $1$-vanishing 
and Condition {\rm\ref{condition;06.6.6.100}} hold
for $(\vecy,L,\alpha_{\ast},\delta)$.
For $\Phi$  as in {\rm(\ref{eq;06.6.11.1})},
we have the following transition formula:
\begin{equation}
\label{eq;06.6.6.105}
 \Phi(\delta_+)-\Phi(\delta_-)
=\sum_{(\vecy_0,\vecY)\in \Sbar(\vecy,\delta)}
 \frac{N_L(y_0)}{N_L(y)}\cdot
 W(\vecY)\cdot
 \int_{\nbigm(\vecyhat_0,\vecYhat,[L])}
 \Psi(\vecy_0,\vecY)
\end{equation}
\end{thm}

The proof will be given in the next subsubsections
\ref{subsubsection;06.6.11.12}--\ref{subsubsection;06.6.11.36}.

\subsubsection{Step 1}
\label{subsubsection;06.6.11.12}

We put $\nbigmtilde(\delta,\ell):=
 \nbigmtilde^{ss}\bigl(\vecy,[L],\alpha_{\ast},(\delta,\ell)\bigr)$
for a positive integer $\ell$.
We put $\nbigmtilde(\delta,0):=
 \nbigmtilde^{ss}(\vecy,[L],\alpha_{\ast},\delta_+)$.
We put $\nbigmtilde(\delta_{-}):=
 \nbigmtilde^{ss}(\vecy,[L],\alpha_{\ast},\delta_{-})$.
Let $\Ttilde_{\rel}$ denote the relative tangent bundle
of the smooth morphism
$\nbigmtilde(m,\vecyhat,[L])\lrarr\nbigm(m,\vecyhat,[L])$.
We have the open immersion of
$\nbigmtilde(\delta,\ell)$
and $\nbigmtilde(\delta_-)$
into $\nbigmtilde(m,\vecyhat,[L])$.
We use the same notation to denote
the restriction of $\Ttilde_{\rel}$
to $\nbigmtilde(\delta,\ell)$
and $\nbigmtilde(\delta_-)$.
We put as follows:
\begin{equation}
\label{eq;06.6.11.11}
 \Phitilde:=\Phi\cdot \frac{\Eu(\Ttilde_{\rel})}{H_y(m)!},
\quad
 \Phitilde(\delta,\ell):=
 \int_{\nbigmtilde(\delta,\ell)}
 \Phitilde,
\quad
 \Phitilde(\delta_{\kappa}):=
 \int_{\nbigmtilde(\delta_{\kappa})}
 \Phitilde
\end{equation}
Then, we have $\Phitilde(\delta,0)=\Phi(\delta_+)$
and $\Phitilde(\delta_-)=\Phi(\delta_-)$.

Recall that $\Dec(m,\vecy,\alpha_{\ast},\delta)$
denotes the set of the decomposition types
(Definition {\rm\ref{df;06.5.22.1}}).
We put as follows:
\[
 S(\ell):=\bigl\{
 \gbigi=(\vecy_1,\vecy_2,I_1,I_2)\in
 \Dec(m,\vecy,\alpha_{\ast},\delta)
 \,\big|\,
 \ellbar\subset I_1
 \bigr\}
\]
For any decomposition type
$\gbigi=(\vecy_1,\vecy_2,I_1,I_2)\in 
 \Dec(m,\vecy,\alpha_{\ast},\delta)$,
we put as follows:
\[
 \nbigmtilde(\gbigihat):=
 \nbigmtilde^{ss}\bigl(\vecyhat_1,[L],\alpha_{\ast},\delta,\gminik(\gbigi)\bigr)
\times
 \nbigm^{ss}(\vecyhat_2,\alpha_{\ast})
\]
Let $\Ttilde_{1,\rel}$ denote the vector bundle
over $\nbigmtilde(\gbigihat)$ 
induced by the relative tangent bundle of
the smooth map
$\nbigmtilde^{ss}\bigl(
 \vecyhat_1,[L],\alpha_{\ast},\delta,\gminik(\gbigi)\bigr)
\lrarr
 \nbigm(m,\vecyhat_1,[L])$.
Let $T_{1,\rel}$ denote the vector bundle
over $\nbigmtilde(\gbigihat)$
obtained from the relative tangent bundle
of the smooth morphism
$\nbigm(m,\vecyhat_1,[L])
\lrarr \nbigm(m,\vecyhat_1)$.

\begin{prop}
\label{prop;06.6.6.55}
We have the following equality:
\begin{equation}
\label{eq;06.6.6.50}
\Phitilde(\delta,\ell)-\Phi(\delta_{-})
=\sum_{\gbigi\in S(\ell)} 
 \frac{N_L(y_1)}{N_L(y)}
 \frac{H_{y_1}(m)!\cdot H_{y_2}(m)!}{H_y(m)!} 
\int_{\nbigmtilde(\gbigihat)}\Psi(\gbigi).
\end{equation}
The elements 
$\Psi(\gbigi)\in
 \nbigr\bigl(\Ehat_1,
 \nbigmtilde^{ss}(\vecyhat_1,[L],\alpha_{\ast},\delta,\gminik(\gbigi))
 \bigr)
 \otimes
 \nbigr(\Ehat_2)$
are given as follows:
\begin{equation}
\label{eq;06.6.6.51}
\Psi(\gbigi)=
 \underset{s=0}\Res
 \left(
 \frac{P\bigl(\Ehat_1^u\cdot e^{-s/r_1}\oplus
 \Ehat_2^u\cdot e^{s/r_2}\bigr)}
 {Q(\Ehat_1^u\cdot e^{-s/r_1},\Ehat^u_2\cdot e^{s/r_2})}
 \right)
\cdot \frac{\Eu(T_{1,\rel})}{N_L(y_1)}
\frac{\Eu(\Ttilde_{1,\rel})}{H_{y_1}(m)!}
\end{equation}
\end{prop}
\pf
The argument is essentially same as
the proof of Theorem \ref{thm;06.6.6.20}
and Theorem \ref{thm;06.6.4.10}.
Hence, we give only an indication.
Let $\Mhat$ denote the master space
connecting $\nbigmtilde(\delta,\ell)$
and $\nbigmtilde(\delta_-)$
constructed in the subsubsection
\ref{subsubsection;06.5.17.10}.
By using $\Mhat$ and the argument in the proof of
Theorem \ref{thm;06.6.6.20},
we obtain the following expression:
\[
 \Phitilde(\delta,\ell)-\Phi(\delta_-)
=\sum_{\gbigi\in S(\ell)}
 \int_{\Mhat^{G_m}(\gbigi)}
 \frac{\Phitilde_t}{\Eu\bigl(\gbign(\Mhat^{G_m}(\gbigi))\bigr)}
\]

The contributions from $\Mhat^{G_m}(\gbigi)$
can be calculated by the arguments
in the proof of Theorem \ref{thm;06.6.6.20}
and Theorem \ref{thm;06.6.4.10}.
We remark that we can use Lemma \ref{lem;06.5.31.105},
due to Condition \ref{condition;06.6.6.100}.
Then, we arrive at the formulas
(\ref{eq;06.6.6.50}) and (\ref{eq;06.6.6.51})
\hfill\qed

\vspace{.1in}

For our later argument,
we reword Proposition \ref{prop;06.6.6.55}.
Let $I$ be a finite subset of $\seisuu_{>0}$
such that $|I|=H_{y}(m)$.
We naturally regard $I$ as the totally ordered set.
Let $i_0$ be any element of $I$.
Let $\nbigmtilde^{ss}\bigl(\vecyhat,[L],\alpha_{\ast},
 (\delta,i_0),I\bigr)$ be the moduli stack
of the objects $(E_{\ast},[\phi],\rho,\nbigf)$
as follows:
\begin{itemize}
\item
 $(E_{\ast},[\phi],\rho)$ is a $\delta$-semistable
 oriented reduced $L$-Bradlow pairs.
\item
 $\nbigf$ is a full flag of $H^0(X,E(m))$
 indexed by $I$.
\item
 The tuple $(E_{\ast},[\phi],\rho,\nbigf)$
 is $(\delta,i_0)$-semistable in following sense.
 We take a partial Jordan-H\"older filtration 
 of $(E_{\ast},[\phi])$:
\[
 E^{(1)}_{\ast}\subsetneq E^{(2)}_{\ast}
\subsetneq
\cdots\subsetneq E^{(j-1)}_{\ast}\subsetneq
(E^{(j)}_{\ast},[\phi])\subsetneq
(E^{(j+1)}_{\ast},[\phi])\subsetneq\cdots\subsetneq
(E^{(k-1)}_{\ast},[\phi])\subsetneq
(E^{(k)}_{\ast},[\phi])
\]
Then, $\nbigf_{i_0}\cap H^0\bigl(X,E^{(j-1)}(m)\bigr)=\{0\}$
and $\nbigf_{i_0}\not\subset H^0\bigl(X,E^{(k-1)}(m)\bigr)$.
\end{itemize}
We have the bijection
$\varphi:I\lrarr \{1,\ldots,H_y(m)\}$
preserving the order.
Then, we have the isomorphism
$\nbigmtilde^{ss}\bigl(\vecyhat,[L],\alpha_{\ast},(\delta,i_0),I\bigr)
\simeq
 \nbigmtilde^{ss}\bigl(\vecyhat,[L],\alpha_{\ast},(\delta,\varphi(i_0))\bigr)$.
We have the open immersion of
$\nbigmtilde^{ss}\bigl(\vecyhat,[L],\alpha_{\ast},(\delta,i_0),I\bigr)$
into $\nbigmtilde(m,\vecyhat,[L])$.
The pull back of $\Ttilde_{\rel}$ is denoted by
the same notation.
Let $\Phitilde$ be as in (\ref{eq;06.6.11.11}).
We put as follows:
\[
 \Phitilde(\delta,i_0,I):=
\int_{\nbigmtilde^{ss}(\vecyhat,[L],\alpha_{\ast},(\delta,i_0),I)}
 \Phitilde
\]

Let $\Dec(m,\vecy,\alpha_{\ast},\delta,I)$ denote the set 
of the tuples $\gbigi=(\vecy_1,\vecy_2,I_1,I_2)$
satisfying the following conditions:
\[
 \vecy_1+\vecy_2=\vecy,
\quad
P^{\alpha_{\ast},\delta}_{\vecy}
=P^{\alpha_{\ast},\delta}_{\vecy_1}
=P^{\alpha_{\ast}}_{\vecy_2},
\quad
 I_1\sqcup I_2=I,\quad
 |I_i|=H_{y_i}(m)
\]
For any $\gbigi=(\vecy_1,\vecy_2,I_1,I_2)\in
   \Dec(m,\vecy,\alpha_{\ast},\delta,I)$,
we put as follows:
\[
 \gminik(\gbigi):=
 \max\bigl\{i\in I_1\,\big|\,
 i<\min(I_2)
 \bigr\}
\]
We also put as follows:
\[
 \nbigmtilde(\gbigi):=
 \nbigmtilde^{ss}\bigl(\vecyhat_1,[L],\alpha_{\ast},
 (\delta,\gminik(\gbigi)),I_1 \bigr)
\times
 \nbigm^{ss}(\vecyhat_2,\alpha_{\ast}),
\]
We put
$S(i_0,I):=\bigl\{\gbigi\in \Dec(m,\vecy,\alpha_{\ast},\delta,I)
 \,\big|\,
 \gminik(\gbigi)\geq i_0\bigr\}$.

Then, Proposition \ref{prop;06.6.6.55} can be
reworded as follows.
\begin{prop}
\label{prop;06.6.11.30}
We have the following equality:
\begin{equation}
\Phitilde(\delta,i_0,I)-\Phi(\delta_-)
=\sum_{\gbigi\in S(i_0,I)} 
 \frac{N_L(y_1)}{N_L(y)}
 \frac{H_{y_1}(m)!\cdot H_{y_2}(m)!}{H_y(m)!} 
\int_{\nbigmtilde(\gbigihat)}\Psi(\gbigi).
\end{equation}
Here, $\Psi(\gbigi)$ are given 
as in {\rm(\ref{eq;06.6.6.51})}.
\hfill\qed
\end{prop}

\subsubsection{Step 2}

We define the set $\Dec^{(j)}(m,\vecy,\alpha_{\ast},\delta)$ inductively.
Put $\Dec^{(1)}(m,\vecy,\alpha_{\ast},\delta):=
 \Dec(m,\vecy,\alpha_{\ast},\delta)$.
 Assume that $\Dec^{(j-1)}(m,\vecy,\alpha_{\ast},\delta)$ is 
 already given.
 Let $\Dec^{(j)}(m,\vecy,\alpha_{\ast},\delta)$ be the set
 of the tuple 
$\gbigi^{(j)}=(\vecy_1^{(j)},\vecY_2^{(j)},I_1^{(j)},\vecI_2^{(j)})$
 as follows:
\begin{itemize}
\item
 $\vecY_2^{(j)}$ denotes an element 
 $(\vecy_2^{(j)},\vecy_2^{(j-1)},\ldots, \vecy_2^{(1)})$
 of $S_j(\vecy,\delta)$.
\item
 $\vecI_2^{(j)}$ denotes a tuple
 $(I_2^{(j)},I_2^{(j-1)},\ldots,I_2^{(1)})$
 of subsets of $\{1,\ldots,H_y(m)\}$.
 Assume $\min(I_2^{(i)})>\min(I_2^{(i-1)})$
 for $i=2,\ldots,j$
\item 
 We assume 
 $\{1,\ldots,H_y(m)\}=I_1^{(j)}\sqcup \coprod_{i=1}^j I_2^{(i)}$.
\item
 We put $\vecy_1^{(j-1)}:=\vecy_1^{(j)}+\vecy_2^{(j)}$
 and $I_1^{(j-1)}:=I_1^{(j)}\sqcup I_2^{(j)}$.
 Then,
 $(\vecy_1^{(j)},\vecy_2^{(j)},I_1^{(j)},I_2^{(j)})$
 is an element of $\Dec(m,\vecy_1^{(j-1)},\alpha_{\ast},I_1^{(j-1)})$,
 in the sense of the subsubsection 
 \ref{subsubsection;06.6.11.12}.
\item
 We put $\vecY_2^{(j-1)}:=(\vecy_2^{(j-1)},\ldots,\vecy_2^{(1)})$
 and $\vecI_2^{(j-1)}=(I_2^{(j-1)},\ldots,I_{2}^{(1)})$.
 Then, 
 $(\vecy_1^{(j-1)},\vecY_2^{(j-1)},I_1^{(j-1)},\vecI_2^{(j-1)})$
 is an element of $\Dec^{(j-1)}(m,\vecy,\alpha_{\ast},\delta)$.
\end{itemize}

Let $\gbigi^{(j)}=(\vecy_1^{(j)},\vecY_2^{(j)},I_1^{(j)},\vecI_2^{(j)})$
be an element of $\Dec^{(j)}(m,\vecy,\alpha_{\ast},\delta)$.
We put as follows:
\[
 \gminik(\gbigi^{(j)}):=
 \max\bigl\{i\in I_1^{(j)}\,\big|\,
 i<\min(I_2^{(j)})
 \bigr\}
\]
We also put as follows:
\[
 \nbigmtilde(\gbigi^{(j)}):=
 \nbigmtilde^{ss}\bigl(\vecyhat_1^{(j)},[L],\alpha_{\ast},
 (\delta,\gminik(\gbigi^{(j)})),I_1^{(j)}\bigr)
\times
 \prod_{i=1}^j
 \nbigm^{ss}(\vecyhat_2^{(i)},\alpha_{\ast})
\]
\[
 \nbigm_-(\gbigi^{(j)}):=
 \nbigm^{ss}\bigl(\vecyhat_1^{(j)},[L],\alpha_{\ast},\delta_-\bigr)
\times
 \prod_{i=1}^j
 \nbigm^{ss}(\vecyhat_2^{(i)},\alpha_{\ast})
\]

Let $\Ehat^{(j)}_1$ denote the universal sheaves on 
$\nbigmtilde^{ss}\bigl(\vecyhat_1^{(j)},[L],\alpha_{\ast},
 (\delta,\gminik(\gbigi^{(j)})),I_1^{(j)}\bigr)\times X$
and
$\nbigm^{ss}\bigl(\vecyhat_1^{(j)},
 [L],\alpha_{\ast},\delta_-\bigr)\times X$.
The appropriate pull backs are denoted by the same notation.
Let $\Ehat^{(i)}_2$ $(i=1,\ldots,j)$
denote the universal sheaves on 
$\nbigm^{ss}\bigl(\vecyhat_2^{(i)},\alpha_{\ast}\bigr)\times X$.
The appropriate pull backs are denoted by the same notation.

Let $T^{(j)}_{1,\rel}$ denote the relative tangent bundle
of the smooth morphism of
$\nbigm\bigl(m,\vecyhat^{(j)}_1,[L]\bigr)$
to $\nbigm\bigl(m,\vecyhat^{(j)}_1\bigr)$.
Let $\Ttilde^{(j)}_{1,\rel}$ denote the relative tangent bundle
of $\nbigmtilde(m,\vecyhat_1^{(j)},[L])
\lrarr \nbigm(m,\vecyhat_1^{(j)},[L])$.
The appropriate pull backs are denoted by the same notation.

We put as follows for variables $s^{(1)},\ldots,s^{(j)}$:
\[
 T_1^{(j)}:=-\sum_{h\leq j}\frac{s^{(h)}}{r_1^{(h)}},
\quad
 T_2^{(i)}:=-\sum_{h<i}\frac{s^{(h)}}{r_1^{(h)}},
 +\frac{s^{(i)}}{r_2^{(i)}},
\]
Here, we put as follows:
\[
r_1^{(j)}:=\rank \vecy_1^{(j)},
\quad
 r_1^{(h)}:=r_1^{(j)}+\sum_{h<p\leq j}r_2^{(p)}
\]
Let $G$ be the $j$-dimensional torus
$\Spec k[\sigma^{(1)},\ldots,\sigma^{(j)}]$.
Let $e^{w\cdot s^{(i)}}$ denote the trivial line bundle
with $G$-action
which is induced by the
$\Spec k[\sigma^{(i)},\sigma^{(i)\,-1}]$-action
of weight $w$.
We have the following element of 
the $K$-group of the $G$-equivariant
coherent sheaves on $\nbigm(\vecyhat_0,\vecYhat,[L])$:
\begin{multline}
-Rp_{X\,\ast}\nrhom\bigl(\Ehat_i^{(a)}\cdot e^{T_i^{(a)}},\,
 \Ehat^{(b)}_l\cdot e^{T_l^{(b)}}\bigr)
-Rp_{X\,\ast}\nrhom\bigl(\Ehat^{(b)}_l\cdot e^{T_l^{(b)}},\,
 \Ehat^{(a)}_i\cdot e^{T_i^{(a)}}\bigr) \\
-Rp_{D\,\ast}\nrhom'\bigl(\Ehat^{(a)}_{i|D\,\ast}\cdot e^{T_i^{(a)}},\,
 \Ehat^{(b)}_{l|D\,\ast}\cdot e^{T_l^{(b)}}\bigr)
-Rp_{D\,\ast}\nrhom'\bigl(\Ehat^{(b)}_{l|D\,\ast}\cdot e^{T_l^{(b)}},\,
 \Ehat^{(a)}_{i|D\,\ast}\cdot e^{T_i^{(a)}}\bigr)
\end{multline}
The equivariant Euler class is denoted by
$Q\bigl(\Ehat^{(a)}_i\cdot e^{T_i^{(a)}},
 \Ehat^{(b)}_l\cdot e^{T_l^{(b)}}\bigr)$,
which is regarded as the element of
$\nbigr(\Ehat^{(j)}_1)
\otimes
\bigotimes_{h=1}^j \nbigr(\Ehat^{(h)}_2)\otimes
\gbigr\bigl(s^{(1)},s^{(2)},\ldots,s^{(j-1)},s^{(j)}\bigr)$.
We put as follows:
\begin{equation}
 Q\bigl(\Ehat^{(j)}_1 e^{T_1^{(j)}},
 \Ehat^{(j)}_2e^{T_2^{(j)}},
\ldots,
 \Ehat^{(1)}_2e^{T_2^{(1)}}\bigr):=
 \prod_{h<i\leq j}
 Q\bigl(\Ehat^{(h)}_2e^{T_2^{(h)}},
 \Ehat^{(i)}_2e^{T^{(i)}_2}\bigr)
\times
 \prod_{h\leq j}
 Q\bigl(\Ehat^{(j)}_1e^{T_1^{(j)}},
 \Ehat^{(h)}_2e^{T^{(h)}_2}\bigr)
\end{equation}
Then, we put as follows:
\[
 \Psi^{(j)}(\vecy_1^{(j)},\vecY_2^{(j)}):=
 \underset{s^{(j)}}\Res\cdots
 \underset{s^{(1)}}\Res
\left(
\frac{P\bigl(\Ehat^{(j)}_1\cdot e^{T_1^{(j)}}
 \oplus
 \bigoplus_{i=0}^k\Ehat_2^{(i)}\cdot e^{T_2^{(i)}}\bigr)}
 {Q\bigl(\Ehat_1^{(j)}\cdot e^{T_1^{(j)}},
 \Ehat_2^{(j)}\cdot e^{T_2^{(j)}},
 \ldots,\Ehat_2^{(1)}\cdot e^{T_2^{(1)}}\bigr)}
\right)
 \cdot \frac{\Eu(T^{(j)}_{1,\rel})}{N_L(y_1^{(j)})}
\]
We also put as follows:
\[
 \Psitilde^{(j)}(\vecy_1^{(j)},\vecY_2^{(j)}):=
 \Psi^{(j)}(\vecy_1^{(j)},\vecY_2^{(j)})\cdot
 \frac{\Eu(\Ttilde_{1,\rel}^{(j)})}
 {H_{y_1^{(j)}}(m)!}
\]

\begin{lem}
For each $j$,
we have the following formula:
\begin{multline}
 \label{eq;06.6.11.25}
 \Phi(\delta_+)-\Phi(\delta_-)
=\sum_{i<j}\sum_{\gbigi^{(i)}\in \Dec^{(i)}(m,\vecy,\alpha_{\ast},\delta)}
 \frac{N_L(y_1^{(i)})}{N_L(y)}
 \frac{H_{y_1^{(i)}}(m)!\cdot \prod_{h=1}^i H_{y_2^{(h)}}(m)!}
 {H_{y}(m)!}
 \int_{\nbigm_-(\gbigi^{(i)})}
 \Psi^{(i)}(\vecy_1^{(i)},\vecY_2^{(i)})\\
+\sum_{\gbigi^{(j)}\in \Dec^{(j)}(m,\vecy,\alpha_{\ast},\delta)}
 \frac{N_L(y_1^{(j)})}{N_L(y)}
 \frac{H_{y_1^{(j)}}(m)!\cdot \prod_{h=1}^j H_{y_2^{(h)}}(m)!}
 {H_{y}(m)!}
 \int_{\nbigmtilde(\gbigi^{(j)})}
 \Psitilde^{(j)}(\vecy_1^{(j)},\vecY_2^{(j)})\\
\end{multline}
\end{lem}
\pf
We use an induction on $j$.
In the case $j=1$,
the claim is the proven in
Proposition \ref{prop;06.6.6.55} (the case $\ell=0$).
Assume that the formula (\ref{eq;06.6.11.25}) holds for $j$,
and we will prove it for $j+1$.
By definition, 
we have the naturally defined morphism
$\pi_j:\Dec^{(j)}\lrarr \Dec^{(j-1)}$.
Due to Proposition \ref{prop;06.6.11.30},
we obtain the following equality:
\begin{multline}
 \int_{\nbigmtilde(\gbigi^{(j)})}
 \Psitilde^{(j)}(\vecy_1^{(j)},\vecY_2^{(j)})
-\int_{\nbigm_-(\gbigi^{(j)})}
 \Psi^{(j)}(\vecy_1^{(j)},\vecY_2^{(j)}) =\\
\sum_{
\substack{\gbigi^{(j+1)}\in \Dec^{(j+1)}(m,\vecy,\alpha_{\ast},\delta)\\
 \pi_j(\gbigi^{(j+1)})=\gbigi^{(j)} }}
 \frac{N_L(y_1^{(j+1)})}{N_L(y_1^{(j)})}
 \cdot 
\frac{H_{y_1^{(j+1)}}(m)\cdot H_{y_2^{(j+1)}}(m)}
 {H_{y_1^{(j)}}(m)!}
 \int_{\nbigmtilde(\gbigi^{(j+1)})}
 \Psitilde^{\prime\,(j+1)}(\vecy_1^{(j+1)},\vecY_2^{(j+1)})
\end{multline}
Here, the elements
$\Psitilde^{\prime\,(j+1)}(\vecy_1^{(j+1)},\vecY_2^{(j+1)})
 \in
 \nbigr(\Ehat_1^{(j+1)})\otimes
 \bigotimes_{h=1}^{j+1}\nbigr(\Ehat_2^{(h)})$
are given as follows:
\begin{multline}
 \Psitilde^{\prime\,(j+1)}(\vecy_1^{(j+1)},\vecY_2^{(j+1)})=\\
\underset{s^{(j+1)}}\Res
 \underset{s^{(j)}}\Res
 \cdots
 \underset{s^{(1)}}\Res
\left(
 \frac{P\bigl(E_1^{(j+1)}e^{T_1^{(j+1)}}\oplus
 \bigoplus_{i=1}^j E_2^{(i)}\cdot e^{T_2^{(i)}}\bigr)}
 {\prod_{h<j+1}
 Q\bigl(E_1^{(j+1)}e^{T_1^{(j+1)}},E_2^{(h)}e^{T_2^{(h)}}\bigr)
\cdot
 \prod_{h<i\leq j+1}Q\bigl(E_2^{(h)}e^{T_2^{(h)}},
 E_2^{(i)}e^{T_2^{(i)}}\bigr) } \right.\\
 \left.
\times
 \frac{1}{Q\bigl(E_1^{(j+1)}e^{-s^{(j+1)}/r_1^{(j+1)}},
 E_2^{(j+1)}e^{-s^{(j+1)}/r_2^{(j+1)}} \bigr)}
 \right)
\end{multline}
It is easy to observe 
$\Psitilde^{\prime(j+1)}(\vecy_1^{(j+1)},\vecY_2^{(j+1)})
=\Psitilde^{(j+1)}(\vecy_1^{(j+1)},\vecY_2^{(j+1)})$.
Thus the formula (\ref{eq;06.6.11.25}) holds
for $j+1$.
\hfill\qed

\subsubsection{Step 3}
\label{subsubsection;06.6.11.36}

When $j$ is sufficiently large,
we have $\Dec^{(j)}(m,\vecy,\alpha_{\ast},\delta)=\emptyset$.
Therefore, we obtain the following:
\begin{equation}
\label{eq;06.6.11.30}
  \Phi(\delta_+)-\Phi(\delta_-)
=\sum_{k}\sum_{\gbigi^{(k)}\in \Dec^{(k)}(m,\vecy,\alpha_{\ast},\delta)}
 \frac{N_L(y_1^{(k)})}{N_L(y)}
 \frac{H_{y_1^{(k)}}(m)!\cdot \prod_{h=1}^k H_{y_2^{(h)}}(m)!}
 {H_{y}(m)!}
 \int_{\nbigm_-(\gbigi^{(k)})}
 \Psi^{(k)}(\vecy_1^{(k)},\vecY_2^{(k)})
\end{equation}

We have the map 
$\rho_k:\Dec^{(k)}(m,\vecy,\alpha_{\ast},\delta)
\lrarr \Sbar_k(\vecy,\delta)$
given as follows:
\[
 \rho_k(\gbigi^{(k)})
=\bigl(\vecy_1^{(k)},\vecy_2^{(k)},\vecy_2^{(k-1)},
\ldots,\vecy_2^{(1)}\bigr)
=(\vecy_0,\vecy_1,\vecy_2,\ldots,\vecy_k)
\]
Then, (\ref{eq;06.6.11.30}) can be rewritten as follows:
\begin{equation}
\label{eq;06.6.11.35}
  \Phi(\delta_+)-\Phi(\delta_-)
=\sum_k
 \sum_{(\vecy_0,\vecY)\in \Sbar_k(\vecy,\delta)}
 \frac{N_L(y_0)}{N_L(y)}
 \frac{\prod_{i=0}^kH_{y_i}(m)!}{H_y(m)!}
 \cdot 
 \#\rho_k^{-1}(\vecy_0,\vecY)
 \cdot
 \int_{\nbigm(\vecyhat_0,\vecYhat,[L])}
 \Psi(\vecy_0,\vecY)
\end{equation}

\begin{lem}
\label{lem;06.6.6.101}
Under Condition {\rm\ref{condition;06.6.6.100}},
we have the following equality:
\[
 \frac{\prod_{i=0}^kH_{y_i}(m)! }{H_{y}(m)!}
 \cdot \#\rho_k^{-1}(\vecy_0,\vecY)
=W(\vecY)
\]
\end{lem}
\pf
It is easy to observe that
$\rho_k^{-1}(\vecy_0,\vecY)$ is bijective
to the following set:
\[
 \Bigl\{
 (I_0,I_1,\ldots,I_k)\,\Big|\,
 \coprod_{i=0}^kI_i=\{1,\ldots,H_y(m)\},\,\,\,
 |I_i|=H_{y_i(m)},\,\,
 \min(I_k)<\min(I_{k-1})<\cdots<\min(I_1)
 \Bigr\}
\]
The correspondence is given by
$(I_1^{(k)},I_2^{(k)},\ldots,I_2^{(1)})
=(I_0,I_1,\ldots,I_k)$.

We put $N=H_y(m)$ and $N_i=H_{y_i}(m)$.
We have the following:
\begin{multline}
\#\rho_k^{-1}(\vecy_0,\vecY)=
 \frac{N!}{N_0!(N-N_0)!}
\cdot
 \frac{(N-N_0-1)!}{(N_k-1)!(N-N_0-N_k)!}
\cdot
 \frac{(N-N_0-N_k-1)!}{(N_{k-1}-1)!(N-N_0-N_k-N_{k-1})!}
\cdot \cdots\\
=
\frac{N!}{N_0!(N-N_0)!}
\cdot
\prod_{i=1}^k
 \frac{(N-N_0-\sum_{j>i}N_j-1)!}
 {(N_{i}-1)!(N-N_0-\sum_{j\geq i}N_j)!}
=\frac{N!}{\prod_{i=0}^kN_i!}
 \cdot
 \frac{\prod_{i=1}^k N_i}
 {\prod_{i=1}^k\sum_{1\leq j\leq i}N_j}
\end{multline}
Under the condition \ref{condition;06.6.6.100},
we have $N_i/r_i=N_j/r_j$ for $1\leq i,j\leq k$.
Therefore, we obtain the following:
\[
 \#\rho_k^{-1}(\vecy_0,\vecY)\cdot
 \frac{\prod_{i=1}^kN_i!}{N!}
= \frac{\prod_{i=1}^k N_i}
 {\prod_{i=1}^k\sum_{1\leq j\leq i}N_j}
=\frac{\prod_{i=1}^k r_i}
 {\prod_{i=1}^k\sum_{1\leq j\leq i}r_j}
\]
Thus the proof of Lemma \ref{lem;06.6.6.101}
is finished.
\hfill\qed

\vspace{.1in}
We immediately obtain (\ref{eq;06.6.6.105})
from (\ref{eq;06.6.11.35}) and Lemma \ref{lem;06.6.6.101}.
\hfill\qed

\subsection{Weak Wall Crossing Formula}
\label{subsection;06.7.3.57}
\subsubsection{Statement}

\label{subsubsection;06.6.7.30}

Let $y$ be an element of $\Type_r^{\circ}$.
We use the notation $\nbigm_H(\yhat)$ 
to denote the moduli stack of torsion-free sheaves
of type $y$, which are semistable with respect to 
a polarization $H$.
For $\Phi=P(\Ehat^u)\in\nbigr(\Ehat^u)$,
we put as follows:
\begin{equation}
 \label{eq;06.6.19.100}
 \Phi_H(\yhat):=\int_{\nbigm_H(\yhat)}\Phi
\end{equation}

Let $\xi$ be an element of $NS(X)$
such that $\xi^2< 0$.
Let $W^{\xi}$ be the wall
in the ample cone determined by $\xi$.
It is well known that the ample cone is
divided into the chambers by such walls,
and $\nbigm_{H}(\yhat)$ depends only 
on the chamber to which $H$ belongs.
Moreover, the set of such walls are locally finite
(\cite{mw}).

Let $C_+$ and $C_-$ be chambers
which are divided by a wall $W^{\xi}$.
Let $H_+$ and $H_-$ be ample line bundles
contained in $H_+$ and $H_-$ respectively.
We assume $(\xi,H_+)>(\xi,H_-)$.
We would like to obtain the description
$\Phi_{H_+}-\Phi_{H_-}$
as the sum of the integrals
over the products of the moduli stacks
of the objects with lower ranks.
In the following,
$r_i$, $a_i$ and $b_i$
 denote the rank, the first Chern class
 and the second Chern character of 
 a given $y_i\in \Type^{\circ}$,
 respectively.
We take a generic $H\in W^{\xi}$.

Let $S_{+,k}$ denote the set of
$\gminiy=(y_0,y_1,\ldots,y_k)\in \Type^{\circ\,k+1}$
satisfying the following conditions:
\begin{itemize}
\item 
 There exist $A_i\geq 0$ $(i=1,\ldots,k)$
 such that $a_i/r_i-a_{i-1}/r_{i-1}=A_i\cdot \xi$ 
 in $H^2(X)$.
\item
 In the case $A_1=0$,
 the inequality $b_0/r_0<b_1/r_1$ holds.
\item
 In the case $A_i=0$ for some $i\geq 2$,
 we have $b_i/r_i\leq b_{i+1}/r_{i+1}$.
\end{itemize}
Let $S_{-,k}$ denote the set of
$\gminiy=(y_0,y_1,\ldots,y_k)\in \Type^{\circ\,k+1}$
satisfying the following conditions:
\begin{itemize}
\item 
 There exist $A_i\leq 0$ $(i=1,\ldots,k)$
 such that $a_i/r_i-a_{i-1}/r_{i-1}=A_i\cdot \xi$.
\item
 In the case $A_1=0$,
 the inequality $b_0/r_0<b_1/r_1$ holds.
\item
 In the case $A_i=0$ for some $i\geq 2$,
 we have $b_i/r_i\leq b_{i+1}/r_{i+1}$.
\end{itemize}
Then, we put as follows:
\[
 S_+:=\bigcup_{k}S_{+,k},
\quad
 S_-:=\bigcup_{k}S_{-,k}
\]

Let $\gminiy=(y_0,\ldots,y_k)$ be an element of $S_{+,k}$
or $S_{-,k}$.
Then the integers $1=i(1)<i(2)<\cdots<i(s)\leq k$
are determined by the following condition:
\[
 P_{y_{i(j)-1}}\neq P_{y_{i(j)}}=\cdots =P_{y_{i(j+1)-1}}
 \neq
 P_{y_{i(j+1)}}
\]
We put $i(s+1)=k+1$ formally.
Then, we put as follows:
\begin{equation}
\label{eq;06.6.7.151}
\Bbar(\gminiy):=\frac{r_0}{r}B(\gminiy),
\quad
 B(\gminiy):=\frac{\prod_{i=1}^kr_i}
 {\prod_{j=1}^s\prod_{h=i(j)}^{i(j+1)-1}
 \sum_{i(j)\leq l\leq h}r_l}
\end{equation}
We also put as follows,
for $\kappa=\pm$
\[
 \nbigm_{H_{\kappa}}(\gminiyhat)
=\prod_{i=0}^k
 \nbigm_{H_{\kappa}}(\yhat_i)
\]
Let $\Ehat_i^u$ denote the sheaf
which is the pull back of the universal sheaf
over $\nbigm(\yhat_i)\times X$.

When $\gminiy$ is given,
we put as follows for variables $t_1,\ldots,t_k$:
\[
 T_0=-\sum_{j>0}\frac{t_j}{\sum_{0\leq h< j}\rank(\vecy_h)},
\quad
 T_i=-\sum_{j>i}\frac{t_j}{\sum_{0\leq h<j}\rank(\vecy_h)}
+\frac{t_i}{\rank(\vecy_i)}
\]
Let $G$ denote the $k$-dimensional torus
$\Spec k[\tau_1,\tau_1^{-1},\ldots,\tau_k,\tau_k^{-1}]$.
Let $e^{w\cdot t_i}$ denote the trivial line bundle
with $G$-action which is induced by
the $\Spec k[\tau_i,\tau_i^{-1}]$-action of weight $w$.
We have the following element of the $K$-group
of the $G$-equivariant coherent sheaves over
$\nbigm(\gminiyhat)$
as in the subsubsection \ref{subsubsection;06.6.23.15}:
\[
 -Rp_{X\,\ast}\Bigl(
 \nrhom\bigl(\Ehat_i\cdot e^{T_i},\,\Ehat_j\cdot e^{T_j}\bigr)
 \Bigr)
-Rp_{X\,\ast}\Bigl(
\nrhom\bigl(\Ehat_j\cdot e^{T_j},\,\Ehat_i\cdot e^{T_i}\bigr)
 \Bigr)
\]
The equivariant Euler class is denoted by
$Q\bigl(\Ehat_i\cdot e^{T_i},\Ehat_j\cdot e^{T_j}\bigr)$.
We regarded it as the element of
$\bigotimes_{i=0}^k\nbigr(\Ehat_i)\otimes\gbigr(t_k,\ldots,t_1)$.
We put as follows:
\[
 Q\bigl(\Ehat_0\cdot e^{T_0},\Ehat_1\cdot e^{T_1},\ldots,
 \Ehat_k\cdot e^{T_k}\bigr):=
 \prod_{i<j}
 Q\bigl(\Ehat_i\cdot e^{T_i},\Ehat_j\cdot e^{T_j}\bigr)
\]
Then, we put as follows:
\begin{equation}
 \label{eq;06.6.7.250}
 \Psi(\gminiy):=
 \underset{t_1}\Res\cdots
 \underset{t_k}\Res
\left(
 \frac{P\bigl(\bigoplus_{i=0}^k\Ehat_i^u\cdot e^{T_i}\bigr) }
 {Q\bigl(\Ehat^u_0\cdot e^{T_0},\ldots,\Ehat_k^u\cdot e^{T_k}\bigr)}
\right)
\end{equation}

\begin{thm}
\label{thm;06.6.7.100}
We have the following formula:
\begin{equation}
\label{eq;06.6.6.200}
 \Phi_{H_+}(\yhat)
-\Phi_{H_-}(\yhat)
=
-\sum_{\gminiy\in S_+}
 \Bbar(\gminiy)\cdot
 \int_{\nbigm_{H_+}(\gminiyhat)}
 \Psi(\gminiy)
+\sum_{\gminiy\in S_-}
 \Bbar(\gminiy)\cdot
 \int_{\nbigm_{H_-}(\gminiyhat)}
 \Psi(\gminiy)
\end{equation}
\end{thm}

We will write down the formula
for the rank $3$ case
in the subsubsection \ref{subsubsection;06.6.7.1}.

\subsubsection{Proof of the theorem}

We use the argument and the notation
in the proof of Theorem \ref{thm;06.6.8.150}.
\begin{lem}
\label{lem;06.6.11.50}
Let $\gminiy=(y_0,\ldots,y_k)$
be an element of $\Type^{\circ k+1}$
such that $\sum y_i=y$ and $\mu(y_i)=\mu(y)$.
We also assume that there are $\mu$-semistable
torsion-free sheaves of type $y_i$ $(i=0,\ldots,k)$.
Then, it is contained in $S_{+,k}$
if and only if 
the following inequality holds for any sufficiently large $t$:
\begin{equation}
\label{eq;06.6.6.150}
 P_{y_0(C)}(t)<P_{y_1(C)}(t)\leq
 P_{y_2(C)}(t)\leq\cdots \leq P_{y_k(C)}(t)
\end{equation}
Similarly,
$\gminiy$ is contained in $S_{-,k}$
if and only if
the following inequality holds for any sufficiently large $t$:
\[
 P_{y_0(-C)}(t)<P_{y_1(-C)}(t)\leq
 P_{y_2(-C)}(t)\leq \cdots\leq P_{y_k(-C)}(t)
\]
\end{lem}
\pf
We remark that finiteness 
of $\gminiy$ as in the assumption of the lemma,
due to the boundedness of $\nbigsbar$.
The condition (\ref{eq;06.6.6.150})
is equivalent to the following:
\begin{equation}
\label{eq;06.6.6.151}
 \frac{b_0+(a_0,C)}{r_0}<
 \frac{b_1+(a_1,C)}{r_1}\leq
 \frac{b_2+(a_2,C)}{r_2}\leq
 \cdots\leq
 \frac{b_k+(a_k,C)}{r_k}
\end{equation}
We have $(a_{i+1},H)/r_{i+1}=(a_{i},H)/r_i$,
and $H$ is generic in $W^{\xi}$.
Hence we have $a_{i+1}/r_{i+1}-a_i/r_i=A_i\cdot \xi$
for some rational numbers $A_i$,
and we have $(a_{i+1},C)/r_{i+1}-(a_i,C)/r_i=A_i\cdot (\xi,C)$.
Since $(\xi,C)$ is sufficiently large,
the condition (\ref{eq;06.6.6.151})
is equivalent to  $\gminiy\in S_{+,k}$.
Thus the first claim is proved.
The second claim can be shown similarly.
\hfill\qed

\vspace{.1in}

For any $\gminiy=(y_0,\ldots,y_k)\in S_+$,
we put as follows:
\[
 \nbigm\bigl(\gminiyhat(C),[\nbigo(-m,C)]\bigr):=
 \nbigm^{s}\bigl(\yhat_0(C),[\nbigo(-m,C)],\epsilon\bigr)
\times
 \prod_{i=1}^k \nbigm^{ss}\bigl(\yhat_i(C)\bigr)
\]
Here, $\epsilon$ denotes any sufficiently small positive number.
Similarly, we put as follows for any $\gminiy\in S_-$:
\[
 \nbigm\bigl(\gminiyhat(-C),[\nbigo(-m,-C)]\bigr):=
 \nbigm^s\bigl(\yhat_0(-C),[\nbigo(-m,-C)],\epsilon\bigr)
\times
 \prod_{i=1}^k \nbigm^{ss}\bigl(\yhat_i(-C)\bigr)
\]
Let $T_{\rel}$ denote the vector bundle
over $\nbigm\bigl(\gminiyhat(C),[\nbigo(-m,C)]\bigr)$
obtained from the relative tangent bundle
of the smooth morphism
$\nbigm^s\bigl(\yhat(C),[\nbigo(-m,C)],\epsilon\bigr)
\lrarr \nbigm(m,\yhat(C))$.
We use the same notation to denote the vector bundle
over $\nbigm^s\bigl(\gminiyhat(-C),[\nbigo(-m,-C)],\epsilon\bigr)$
obtained from the relative tangent bundle
of the smooth morphism
$\nbigm\bigl(\yhat(-C),[\nbigo(-m,-C)],\epsilon\bigr)
\lrarr\nbigm(m,\yhat(-C))$.

Let $\delta_0$ be as in the proof of Theorem 
\ref{thm;06.6.8.150}.
To obtain the description 
for $\Phitilde_C(\delta_0)-\Phi(H_+)$
and $\Phitilde_{-C}(\delta_0)-\Phi(H_-)$,
we use the transition formula
(Theorem \ref{thm;06.6.6.30}) inductively.
Due to Lemma \ref{lem;06.6.11.50},
we obtain the following equality:
\begin{equation}
 \Phitilde_C(\delta_0)
=\Phi_{H_+}(\yhat)
+\sum_{\gminiy\in S_+}
 \frac{H_{y_0}(m)}{H_{y}(m)}\cdot
 B(\gminiy) \cdot 
 \int_{\nbigm(\gminiyhat(C),[\nbigo(-m,C)])}
 \Psi(\gminiy)\cdot\frac{\Eu(T_{\rel})}{H_{y_0}(m)}
\end{equation}
We also obtain the following equality:
\begin{equation}
 \Phitilde_{-C}(\delta_0)
=\Phi_{H_-}(\yhat)
+\sum_{\gminiy\in S_-}
 \frac{H_{y_0}(m)}{H_{y}(m)}\cdot
 B(\gminiy) \cdot 
 \int_{\nbigm(\gminiyhat(-C),[\nbigo(-m,-C)])}
 \Psi(\gminiy)\cdot\frac{\Eu(T_{\rel})}{H_{y_0}(m)}
\end{equation}
From the equality
$\Phitilde_{-C}(\delta_0)=\Phitilde_C(\delta_0)$
(Lemma \ref{lem;06.6.11.55}),
we obtain the following equality:
\begin{multline}
\Phi_{H_+}(\yhat)
+\sum_{\gminiy\in S_+}
 \frac{H_{y_0}(m)}{H_{y}(m)}\cdot
 B(\gminiy) \cdot 
 \int_{\nbigm(\gminiyhat(C),[\nbigo(-m,C)])}
 \Psi(\gminiy)\cdot\frac{\Eu(T_{\rel})}{H_{y_0}(m)} \\
=\Phi_{H_-}(\yhat)
+\sum_{\gminiy\in S_-}
 \frac{H_{y_0}(m)}{H_{y}(m)}\cdot
 B(\gminiy) \cdot 
 \int_{\nbigm(\gminiyhat(-C),[\nbigo(-m,-C)])}
 \Psi(\gminiy)\cdot\frac{\Eu(T_{\rel})}{H_{y_0}(m)}
\end{multline}
By taking the limit $m\to\infty$,
we obtain the desired equality 
(\ref{eq;06.6.6.200}),
due to Proposition \ref{prop;06.5.31.45}.
\hfill\qed

\subsubsection{Weak wall crossing formula in the rank 3 case}
\label{subsubsection;06.6.7.1}

As an example,
we write down the formula (\ref{eq;06.6.6.200})
in the rank $3$ case.
In the following,
$a_i$ and $b_i$ denote the first Chern class
and the second Chern character of a given $y_i$.
We put as follows:
\[
 U_1:=\bigl\{
 (y_1,y_2)\in\Type_1\times\Type_2\,\big|\,
 y_1+y_2=y
 \bigr\},
\quad
 U_2:=\bigl\{
 (y_1,y_2)\in\Type_2\times\Type_1\,\big|\,
 y_1+y_2=y
 \bigr\},
\]
\[
  U_3:=\bigl\{
 (y_1,y_2,y_3)\in\Type_1^3\,\big|\,
 y_1+y_2+y_3=y
 \bigr\}
\]
We put as follows:
\[
 \begin{array}{l}
 S^{(1)}_{1,+}:=\bigl\{
 (y_1,y_2)\in U_1\,\big|\,
 a_2-2a_1=A\cdot \xi\,\,(A>0)
 \bigr\},\,\,\,
 S^{(2)}_{1,+}:=\bigl\{
 (y_1,y_2)\in U_1\,\big|\,
 2a_1=a_2,\,\,2b_1<b_2
 \bigr\}\\
 S^{(1)}_{2,+}:=\bigl\{
 (y_1,y_2)\in U_2\,\big|\,
 2a_2-a_1=A\cdot \xi\,\,(A>0)
 \bigr\},\,\,\,
   S^{(2)}_{2,+}:=\bigl\{
 (y_1,y_2)\in U_2
\big|\,
 a_1=2a_2,\,
 b_1<2b_2
 \bigr\}\\
 S^{(1)}_{3,+}:=\bigl\{
 (y_1,y_2,y_3)\in U_3\,\big|\,
 a_2-a_1=A\cdot \xi,\,\,
 a_3-a_2=B\cdot\xi\,\,
 (A,B>0)
 \bigr\}\\
 S^{(2)}_{3,+}:=\bigl\{
 (y_1,y_2,y_3)\in U_3\,\big|\,
 a_1+A\cdot \xi=a_2=a_3\,\,
 (A>0),\,\,
  b_2<b_3
 \bigr\}
 \\
 S^{(3)}_{3,+}:=\bigl\{
 (y_1,y_2,y_3)\in U_3\, \big|\,
 a_1=a_2=a_3-A\cdot \xi\,\,(A>0),\,\,
 b_1<b_2
 \bigr\}\\
 S^{(4)}_{3+}:=\bigl\{
 (y_1,y_2,y_3)\in U_3\,\big|\,
 a_1=a_2=a_3,\,\,
 b_1<b_2<b_3
 \bigr\}\\
 S^{(1)}_{4,+}:=\bigl\{
 (y_1,y_2,y_3)\in U_3\,\big|\,
 a_1+A\cdot\xi=a_2=a_3\,\,(A>0),\,\,
 b_2=b_3
 \bigr\}\\
 S^{(2)}_{4,+}:=\bigl\{
 (y_1,y_2,y_3)\in U_3\,\big|\,
 a_1=a_2=a_3,\,\,
 b_1<b_2=b_3
 \bigr\}
 \end{array}
\]

We also put as follows:
\[
 \begin{array}{l}
 S^{(1)}_{1,-}:=\bigl\{
 (y_1,y_2)\in U_1\,\big|\,
 a_2-2a_1=A\cdot \xi\,\,(A<0)
 \bigr\},\,\,\,
 S^{(2)}_{1,-}:=S'_{1,+}=\bigl\{
 (y_1,y_2)\in U_1\,\big|\,
 2a_1=a_2,\,\,2b_1<b_2
 \bigr\}\\
 S^{(1)}_{2,-}:=\bigl\{
 (y_1,y_2)\in U_2\,\big|\,
 2a_2-a_1=A\cdot \xi\,\,(A<0)
 \bigr\},\,\,\,
   S^{(2)}_{2,-}:=S^{(2)}_{2,+}=\bigl\{
 (y_1,y_2)\in U_2
\big|\,
 a_1=2a_2,\,
 b_1<2b_2
 \bigr\}\\
 S^{(1)}_{3,-}:=\bigl\{
 (y_1,y_2,y_3)\in U_3\,\big|\,
 a_2-a_1=A\cdot \xi,\,\,
 a_3-a_2=B\cdot\xi\,\,
 (A,B<0)
 \bigr\}\\
 S^{(2)}_{3,-}:=\bigl\{
 (y_1,y_2,y_3)\in U_3\,\big|\,
 a_1+A\cdot \xi=a_2=a_3\,\,
 (A<0),\,\,
  b_2<b_3
 \bigr\}
 \\
 S^{(3)}_{3,-}:=\bigl\{
 (y_1,y_2,y_3)\in U_3\, \big|\,
 a_1=a_2=a_3-A\cdot \xi\,\,(A<0),\,\,
 b_1<b_2
 \bigr\}\\
 S^{(4)}_{3,-}:=S^{(4)}_{3,+}=\bigl\{
 (y_1,y_2,y_3)\in U_3\,\big|\,
 a_1=a_2=a_3,\,\,
 b_1<b_2<b_3
 \bigr\}\\
 S^{(1)}_{4,-}:=\bigl\{
 (y_1,y_2,y_3)\in U_3\,\big|\,
 a_1+A\cdot\xi=a_2=a_3\,\,(A<0),\,\,
 b_2=b_3
 \bigr\}\\
 S^{(2)}_{4,-}:=S^{(2)}_{4,+}=
 \bigl\{
 (y_1,y_2,y_3)\in U_3\,\big|\,
 a_1=a_2=a_3,\,\,
 b_1<b_2=b_3
 \bigr\}
 \end{array}
\]

We use the notation $\nbigm(\yhat,H_{\kappa})$
to denote the moduli stack of
oriented torsion-free sheaves
of type $y$
which are semistable with respect to $H_{\kappa}$.
In the case $\rank(y)=1$,
we omit to denote $H_{\kappa}$.
Let $(y_0,y_1)$ be an element of $\Type^2$.
Let $G_m$ be a one dimensional torus.
From $G_m$-equivariant coherent sheaf $\nbige_i$ $(i=1,2)$
on $\nbigm(\yhat_1,H_{\kappa})\times
 \nbigm(\yhat_2,H_{\kappa})\times X$,
we obtain the following element of
the $k$-group of the $G$-equivariant coherent sheaves
on $\nbigm(\yhat_1,H_{\kappa})
 \times\nbigm(\yhat_2,H_{\kappa})$:
\[
 -Rp_{X\,\ast}\bigl(
 \nrhom(\nbige_1,\nbige_2)\bigr)
 -Rp_{X\,\ast}\bigl(
 \nrhom(\nbige_2,\nbige_1)\bigr)
\]
The $G_m$-equivariant Euler class is denoted by
$Q(\nbige_1,\nbige_2)$.
Let $P$ be an element of $\nbigr$.
Then, we put as follows:
\[
 \nbigp(\nbige_1,\nbige_2):=
 \nbigp(\nbige_1,\nbige_2)
:=\frac{
 P(\nbige_1\oplus\nbige_2)}
 {Q(\nbige_1,\nbige_2)} 
 \in \nbigr(\nbige_1)\otimes\nbigr(\nbige_2)
\otimes\gbigr(t)
\]

Let $(y_1,y_2,y_3)$ be an element of $\Type^3$.
Let $G_m^2$ denote the two dimensional torus.
Similarly,
we put as follows
for given $G_m^2$-equivariant
coherent sheaves $\nbige_i$ $(i=1,2,3)$:
\[
 \nbigp(\nbige_1,\nbige_2,\nbige_3)
:=\frac{
 P(\nbige_1\oplus \nbige_2\oplus \nbige_3)}
 {\prod_{i<j}Q(\nbige_i,\nbige_j)}
\in \nbigr(\nbige_1)\otimes\nbigr(\nbige_2)\otimes\nbigr(\nbige_3)
 \otimes\gbigr((t_2,t_1)).
\]

Then, we put as follows for $\kappa=\pm$:
\begin{multline}
 \Lambda(H_{\kappa}):=\Phi_{H_{\kappa}}(\yhat)
+\sum_{i=1,2} \sum_{(y_1,y_2)\in S^{(i)}_{1,{\kappa}}}
 \frac{1}{3}\int_{\nbigm(\yhat_1)\times\nbigm(\yhat_2,H_{\kappa})}
  \underset{t}\Res\Bigl(
 \nbigp\bigl(\Ehat_1^u\!\cdot\! e^{-t},\,
  \Ehat_2^u\!\cdot\! e^{t/2}\bigr) \Bigr)\\
+\sum_{i=1,2}\sum_{(y_1,y_2)\in S^{(i)}_{2,\kappa}}
 \frac{2}{3}\int_{\nbigm(\yhat_1,H_{\kappa})\times\nbigm(\yhat_2)}
  \underset{t}\Res\Bigl(
 \nbigp\bigl(\Ehat_1^u\!\cdot\! e^{-t/2},\,
  \Ehat_2^u\!\cdot\! e^{t}\bigr) \Bigr)\\
+\sum_{i=1}^{4}\sum_{(y_1,y_2,y_3)\in S^{(i)}_{3,\kappa}}
 \frac{1}{3}\int_{\prod_i\nbigm(\yhat_i)}
 \underset{t_1}\Res\,
 \underset{t_2}\Res\Bigl(
 \nbigp\bigl(\Ehat_1^u\!\cdot\! e^{-t_1-t_2/2},\,
 \Ehat_2^u\!\cdot\! e^{t_1-t_2/2},\,
 \Ehat_3^u\!\cdot\! e^{t_2}\bigr)\Bigr)\\
+\sum_{i=1,2}\sum_{(y_1,y_2,y_3)\in S^{(i)}_{4,\kappa}}
 \frac{1}{6}\int_{\prod_i\nbigm(\yhat_i)}
 \underset{t_1}\Res\,
 \underset{t_2}\Res\Bigl(
 \nbigp\bigl(\Ehat_1^u\!\cdot\! e^{-t_1-t_2/2},\,
 \Ehat_2^u\!\cdot\! e^{t_1-t_2/2},\,
 \Ehat_3^u\!\cdot\! e^{t_2}\bigr)\Bigr)
\end{multline}

Then, the formula (\ref{eq;06.6.6.200})  says 
$\Lambda(H_+)=\Lambda(H_-)$.
The contributions from $S^{(1)}_{3,\kappa}$
and $S^{(2)}_{4,\kappa}$
are cancelled out.

\subsubsection{Weak intersection rounding formula in the rank $3$ case}
\label{subsubsection;06.6.8.300}

The weak wall crossing formula (\ref{eq;06.6.6.200})
is not so easy to deal with,
even in the rank $3$ case.
We would like to derive a more accessible
quantity from our invariants.
For that purpose, 
we consider the ``intersection rounding formula''.
We will show it in the general case, later
(the subsection \ref{subsection;06.6.7.400}).
However, the general statement is a little complicated.
So we give the formula in the rank $3$ case,
in this subsubsection.

We take an element $\vecxi=(\xi_1,\xi_2)\in NS(X)^2$
such that $\xi_1$ and $\xi_2$ are linearly independent,
and we put $W^{\vecxi}:=W^{\xi_1}\cap W^{\xi_2}$.
A connected component $T$ of
$ W^{\vecxi}\setminus
 \bigcup_{W\neq W^{\xi_i}}W$
is called a tile.
For each tile $T$,
there exists four chambers 
whose closures contain $T$.
Let $C$ be such a chamber,
and let $H_C$ be an ample line bundle
contained in $C$.
Then the map $\varphi_C:\{1,2\}\lrarr \{1,-1\}$
is determined by
$\varphi_C(i)=\sign(H_C,\xi_i)$.
The chambers and the maps correspond
bijectively.
We denote by $C(\varphi)$
the chamber corresponding to a map
$\varphi:\{1,2\}\lrarr \{1,-1\}$.
We put $k(\varphi):=\#\{i\,|\,\varphi(i)=-1\}$.
We take any ample line bundle 
$H_{\varphi}\in C(\varphi)$.
Then, we put as follows:
\[
 \nbigd_{\vecxi}^T\Phi(\yhat):=
 \sum_{\varphi} (-1)^{k(\varphi)}
\Phi_{H_{\varphi}}(\yhat)
\]
We would like to show that 
$\nbigd_{\vecxi}^T\Phi(\yhat)$
is independent of a choice of $T$,
and we would like to express
$\nbigd_{\vecxi}^T\Phi(\yhat)$
as the sum of the integrals
over the products of Hilbert schemes.

Let $S(2,1)$ be the set of
$\veca=(a_0,a_1,a_2)\in NS(X)^3$
with the following property:
\begin{itemize}
\item
 $(a_0+a_1)/2-a_2=A_1\cdot \xi_1$
 and $a_0-a_1=A_2\cdot \xi_2$
 for some $A_i>0$.
\end{itemize}
Let $S(1,2)$ be the set of
$\veca=(a_0,a_1,a_2)\in NS(X)^3$
with the following property:
\begin{itemize}
\item
 $a_0-(a_1+a_2)/2=A_1\cdot \xi_1$
 and $a_1-a_2=A_2\cdot \xi_2$
 for some $A_i>0$.
\end{itemize}
For each $\veca$,
we put as follows:
\[
 X(y,\veca):=
 \coprod_{n_0+n_1+n_2=N(y,\veca)}
 \prod_{i=0}^2 (X^{[n_i]}
\times \Pic(a_i))
\quad
 N(y,\veca)=n+\frac{a_0^2+a_1^2+a_2^2-a^2}{2}
\]
Here, $a$ and $n$ denote the first Chern class
and the second Chern class of $y$.
Let $\Ehat_i^u$ denote the sheaf
over $\prod (X^{[n_i]}\times\Pic(a_i))\times X$
which is the pull back of the universal ideal sheaf 
over $X^{[n_i]}\times\Pic(a_i)\times X$.

Let $G=G_m^2$ denote the two dimensional torus.
Let $\nbige_j$ $(j=1,2,3)$ be $G$-equivariant coherent sheaves
on $\prod(X^{[n_i]}\times\Pic(a_i)) \times X$.
Then, $Q(\nbige_j,\nbige_k)$ denotes the $G$-equivariant
Euler class of the following:
\[
 -Rp_{X\,\ast}\bigl(
 \nrhom(\nbige_j,\nbige_k)\bigr)
 -Rp_{X\,\ast}\bigl(
 \nrhom(\nbige_k,\nbige_j)\bigr)
\]
We use the following notation:
\[
 \nbigp(\nbige_1,\nbige_2,\nbige_3):=
 \frac{
 P\bigl(\nbige_1\oplus\nbige_2\oplus\nbige_3\bigr)}
 {\prod_{i<j}Q(\nbige_i,\nbige_j)}
\in
 \nbigr(\nbige_1)\otimes\nbigr(\nbige_2)
\otimes\nbigr(\nbige_3)\otimes \gbigr(t_1,t_2)
\]

\begin{prop}
\mbox{{}}
\begin{itemize}
\item
 $\nbigd_{\vecxi}^T\Phi(\yhat)$ is independent of the choice of $T$.
 Therefore, we can omit to denote $T$.
\item
 The following equality holds:
\begin{multline}
\label{eq;06.6.7.2}
 \nbigd_{\vecxi}\Phi(\yhat)=
\sum_{\veca\in S(1,2)}
 \int_{X(y,\veca)}
 \underset{t_2}\Res\,
 \underset{t_1}\Res\Bigl(
 \nbigp(\Ehat_0^u e^{-t_1},\,\,
 \Ehat_1e^{t_1/2-t_2},\,\,
 \Ehat_2^ue^{t_1/2+t_2})
\Bigr) \\
+\sum_{\veca\in S(2,1)}
 \int_{X(y,\veca)}
 \underset{t_2}\Res\,
 \underset{t_1}\Res\Bigl(
\nbigp\bigl(\Ehat_0^ue^{-t_1/2-t_2},\,\,
 \Ehat^u_1e^{-t_1/2+t_2},\,\,\,
 \Ehat^u_2e^{t_1}\bigr)\Bigr)
\end{multline}
\end{itemize}
\end{prop}
\pf
We use the notation in the subsubsection 
\ref{subsubsection;06.6.7.1}.
We identify $\varphi$ and $(\varphi(1),\varphi(2))$.
We have the equality:
\begin{equation}
 \Lambda(H_{+,+})-\Lambda(H_{+,-})
=\Lambda(H_{+,-})-\Lambda(H_{-,-})
\end{equation}
We put as follows:
\[
 \begin{array}{l}
 S^{(1)}_{1,+}:=\bigl\{
 (y_1,y_2)\in U_1\,\big|\,
 a_2-2a_1=A\cdot \xi_1\,\,(A>0)
 \bigr\},\,\,\,
 S^{(2)}_{1,+}:=\bigl\{
 (y_1,y_2)\in U_1\,\big|\,
 2a_1=a_2,\,\,2b_1<b_2
 \bigr\}\\
 S^{(1)}_{2,+}:=\bigl\{
 (y_1,y_2)\in U_2\,\big|\,
 2a_2-a_1=A\cdot \xi_1\,\,(A>0)
 \bigr\},\,\,\,
   S^{(2)}_{2,+}:=\bigl\{
 (y_1,y_2)\in U_2
\big|\,
 a_1=2a_2,\,
 b_1<2b_2
 \bigr\}
 \end{array}
\]
We also put as follows:
\[
 \begin{array}{l}
 S^{(1)}_{1,-}:=\bigl\{
 (y_1,y_2)\in U_1\,\big|\,
 a_2-2a_1=A\cdot \xi_1\,\,(A<0)
 \bigr\},\quad
 S^{(2)}_{1,-}:=S'_{1,+}\\
 S^{(1)}_{2,-}:=\bigl\{
 (y_1,y_2)\in U_2\,\big|\,
 2a_2-a_1=A\cdot \xi_1\,\,(A<0)
 \bigr\},\quad
   S^{(2)}_{2,-}:=S^{(2)}_{2,+}
 \end{array}
\]
Then, we have the following equality for $\kappa=\pm$:
\begin{multline}
 \Lambda(H_{\kappa,+})
-\Lambda(H_{\kappa,-})
=\Phi(H_{\kappa,+})-\Phi(H_{\kappa,-})\\
+\sum_{i=1,2}\sum_{(y_1,y_2)\in S_{1,\kappa}^{(i)}}
  \frac{1}{3}
 \int_{\nbigm(\yhat_1)\times\nbigm(\yhat_2,H_{\kappa,+})}
  \underset{t}\Res\Bigl(
 \nbigp\bigl(\Ehat_1^u\cdot e^{-t},\,
  \Ehat_2^u\cdot e^{t/2}\bigr) \Bigr)\\
-\sum_{i=1,2}\sum_{(y_1,y_2)\in S_{1,\kappa}^{(i)}}
  \frac{1}{3}
\int_{\nbigm(\yhat_1)\times\nbigm(\yhat_2,H_{\kappa,-})}
  \underset{t}\Res\Bigl(
 \nbigp\bigl(\Ehat_1^u\cdot e^{-t},\,
  \Ehat_2^u\cdot e^{t/2}\bigr) \Bigr) \\
+\sum_{i=1,2}\sum_{(y_1,y_2)\in S_{2,\kappa}^{(i)}}
 \frac{2}{3}\int_{\nbigm(\yhat_1,H_{\kappa,+})\times\nbigm(\yhat_2)}
  \underset{t}\Res\Bigl(
 \nbigp\bigl(\Ehat_1^u\cdot e^{-t/2},\,
  \Ehat_2^u\cdot e^{t}\bigr) \Bigr)\\
-\sum_{i=1,2}\sum_{(y_1,y_2)\in S_{2,\kappa}^{(i)}}
 \frac{2}{3}\int_{\nbigm(\yhat_1,H_{\kappa,-})\times\nbigm(\yhat_2)}
  \underset{t}\Res\Bigl(
 \nbigp\bigl(\Ehat_1^u\cdot e^{-t/2},\,
  \Ehat_2^u\cdot e^{t}\bigr) \Bigr)
\end{multline}
By a formal calculations using the weak wall crossing formula
in the rank $2$ (Theorem \ref{thm;06.6.4.201}),
we can show the following equalities:
\begin{multline}
 \sum_{(y_1,y_2)\in S_{1,-}^{(1)}}
 \left(
 \int_{\nbigm(\yhat_1)\times\nbigm(\yhat_2,H_{-+})}
 \underset{t}\Res\,
 \nbigp\bigl(\Ehat_1^u\cdot e^{-t},\,\Ehat_2^u\cdot e^{t/2} \bigr)
- \int_{\nbigm(\yhat_1)\times\nbigm(\yhat_2,H_{--})}
 \underset{t}\Res\,
 \nbigp\bigl(\Ehat_1^u\cdot e^{-t},\,\Ehat_2^u\cdot e^{t/2} \bigr)
 \right) \\
=\sum_{\veca\in S(1,2)}
 \underset{t_2}\Res\,
 \underset{t_1}\Res\,
\int_{X(y,\veca)}
 \nbigp\bigl(\Ehat_0^u\cdot e^{-t_1},
 \Ehat_1^u\cdot e^{t_1/2-t_2},
 \Ehat_2^u\cdot e^{t_1/2+t_2}
 \bigr)
\end{multline}
\begin{multline}
 \sum_{(y_1,y_2)\in S^{(1)}_{2,-}}
 \left(
 \int_{\nbigm(\yhat_1,H_{-+})\times\nbigm(\yhat_2)}
 \underset{t}\Res\,
 \nbigp\bigl(\Ehat_1\cdot e^{-t/2},\Ehat_2\cdot e^{t}\bigr)
-\int_{\nbigm(\yhat_1,H_{--})\times\nbigm(\yhat_2)}
 \underset{t}\Res\,
 \nbigp\bigl(\Ehat_1\cdot e^{-t/2},\Ehat_2\cdot e^{t}\bigr)
 \right)\\
=\sum_{\veca\in S(2,1)}
 \int_{X(y,\veca)}
 \underset{t_2}\Res\,
 \underset{t_1}\Res\,
 \nbigp\bigl(\Ehat_0\cdot e^{-t_1/2-t_2},
 \Ehat_1\cdot e^{-t_1/2+t_2},
 \Ehat_2\cdot e^{t_1}
 \bigr)
\end{multline}
We also have the following:
\begin{multline}
 \sum_{(y_1,y_2)\in S^{(1)}_{1,+}}
 \left(
 \int_{\nbigm(\yhat_1)\times\nbigm(\yhat_2,H_{++})}
 \underset{t}\Res\,
 \nbigp\bigl(\Ehat_1^u\cdot e^{-t},\Ehat_2^ue^{t/2}\bigr)
-\int_{\nbigm(\yhat_1)\times\nbigm(\yhat_2,H_{+-})}
 \underset{t}\Res\,
 \nbigp\bigl(\Ehat_1^u\cdot e^{-t},\Ehat_2^ue^{t/2}\bigr)
 \right)\\
= -\sum_{(y_1,y_2)\in S^{(1)}_{1,+}}
 \left(
 \int_{\nbigm(\yhat_2,H_{++})\times\nbigm(\yhat_1)}
 \underset{t}\Res\,
 \nbigp\bigl(\Ehat_2^ue^{-t/2},\Ehat_1^u\cdot e^{t}
 \bigr)
-\int_{\nbigm(\yhat_2,H_{+-})\times\nbigm(\yhat_1)}
 \underset{t}\Res\,
 \nbigp\bigl(\Ehat_2^ue^{t/2},\Ehat_1^u\cdot e^{-t} \bigr)
 \right)\\
=-\sum_{\veca\in S(2,1)}
 \int_{X(y,\veca)}
 \underset{t_2}\Res\,
 \underset{t_1}\Res\,
 \nbigp\bigl(\Ehat_0\cdot e^{-t_1/2-t_2},
 \Ehat_1\cdot e^{-t_1/2+t_2},
 \Ehat_2\cdot e^{t_1}
 \bigr)
\end{multline}
\begin{multline}
\sum_{(y_1,y_2)\in S^{(1)}_{2,+}}
\left(
 \int_{\nbigm(\yhat_1,H_{++})\times\nbigm(\yhat_2)}
 \underset{t}\Res\,
 \nbigp\bigl(\Ehat_1\cdot e^{-t/2},\Ehat_2^u\cdot e^t\bigr)
-\int_{\nbigm(\yhat_1,H_{+-})\times\nbigm(\yhat_2)}
 \underset{t}\Res\,
 \nbigp\bigl(\Ehat_1\cdot e^{-t/2},\Ehat_2^u\cdot e^t\bigr)
 \right)\\
=-\sum_{(y_1,y_2)\in S^{(1)}_{2,+}}
\left(
 \int_{\nbigm(\yhat_2)\times\nbigm(\yhat_1,H_{++}) }
 \underset{t}\Res\,
 \nbigp\bigl(\Ehat_2^u\cdot e^{-t},\Ehat_1\cdot e^{t/2} \bigr)
-\int_{\nbigm(\yhat_2)\times\nbigm(\yhat_1,H_{+-}) }
 \underset{t}\Res\,
 \nbigp\bigl(\Ehat_2^u\cdot e^t,\Ehat_1\cdot e^{-t/2} \bigr)
 \right)\\
=-\sum_{\veca\in S(1,2)}
 \underset{t_2}\Res\,
 \underset{t_1}\Res\,
\int_{X(y,\veca)}
 \nbigp\bigl(\Ehat_0^u\cdot e^{-t_1},
 \Ehat_1^u\cdot e^{t_1/2-t_2},
 \Ehat_2^u\cdot e^{t_1/2+t_2}
 \bigr)
\end{multline}

Recall $S_{1,+}^{(2)}=S_{1,-}^{(2)}$.
Then, it is also easy to observe the following equalities
by using the weak wall crossing formula
in the rank $2$ case:
\begin{multline}
\sum_{(y_1,y_2)\in S_{1,+}^{(2)}}
 \int_{\nbigm(\yhat_1)\times\nbigm(\yhat_2,H_{+,+})}
  \underset{t}\Res\Bigl(
 \nbigp\bigl(\Ehat_1^u\cdot e^{-t},\,
  \Ehat_2^u\cdot e^{t/2}\bigr) \Bigr)\\
-\sum_{(y_1,y_2)\in S_{1,+}^{(2)}}
\int_{\nbigm(\yhat_1)\times\nbigm(\yhat_2,H_{+,-})}
  \underset{t}\Res\Bigl(
 \nbigp\bigl(\Ehat_1^u\cdot e^{-t},\,
  \Ehat_2^u\cdot e^{t/2}\bigr) \Bigr) \\
=\sum_{(y_1,y_2)\in S_{1,-}^{(2)}}
 \int_{\nbigm(\yhat_1)\times\nbigm(\yhat_2,H_{-,+})}
  \underset{t}\Res\Bigl(
 \nbigp\bigl(\Ehat_1^u\cdot e^{-t},\,
  \Ehat_2^u\cdot e^{t/2}\bigr) \Bigr)\\
-\sum_{(y_1,y_2)\in S_{1,-}^{(2)}}
\int_{\nbigm(\yhat_1)\times\nbigm(\yhat_2,H_{-,-})}
  \underset{t}\Res\Bigl(
 \nbigp\bigl(\Ehat_1^u\cdot e^{-t},\,
  \Ehat_2^u\cdot e^{t/2}\bigr) \Bigr) 
\end{multline}

\begin{multline}
\sum_{(y_1,y_2)\in S_{2,+}^{(2)}}
  \int_{\nbigm(\yhat_1,H_{+,+})\times\nbigm(\yhat_2)}
  \underset{t}\Res\Bigl(
 \nbigp\bigl(\Ehat_1^u\cdot e^{-t/2},\,
  \Ehat_2^u\cdot e^{t}\bigr) \Bigr)\\
-\sum_{(y_1,y_2)\in S_{2,+}^{(2)}}
  \int_{\nbigm(\yhat_1,H_{+,-})\times\nbigm(\yhat_2)}
  \underset{t}\Res\Bigl(
 \nbigp\bigl(\Ehat_1^u\cdot e^{-t/2},\,
  \Ehat_2^u\cdot e^{t}\bigr) \Bigr)\\
=\sum_{(y_1,y_2)\in S_{2,-}^{(2)}}
  \int_{\nbigm(\yhat_1,H_{-,+})\times\nbigm(\yhat_2)}
  \underset{t}\Res\Bigl(
 \nbigp\bigl(\Ehat_1^u\cdot e^{-t/2},\,
  \Ehat_2^u\cdot e^{t}\bigr) \Bigr)\\
-\sum_{(y_1,y_2)\in S_{2,-}^{(2)}}
  \int_{\nbigm(\yhat_1,H_{-,-})\times\nbigm(\yhat_2)}
  \underset{t}\Res\Bigl(
 \nbigp\bigl(\Ehat_1^u\cdot e^{-t/2},\,
  \Ehat_2^u\cdot e^{t}\bigr) \Bigr)
\end{multline}
Then, we obtain (\ref{eq;06.6.7.2})
by a formal calculation.
\hfill\qed

\subsubsection{Transition for a critical parabolic weight}
\label{subsubsection;06.7.3.575}

We give a generalization of 
Proposition \ref{prop;06.6.5.3}
in the higher rank case.
The argument is essentially same as
the proof of Theorem \ref{thm;06.6.7.100}.
Hence we state the result without a proof.

Let $\alpha$ be a critical parabolic weight.
We take $\alpha_-<\alpha<\alpha_+$
sufficiently closely.
For $\Phi=P(\Ehat^u)\in\nbigr(\Ehat^u)$,
we put as follows:
\[
 \Phi(\yhat,\alpha_{\kappa}):=
 \int_{\nbigm(\yhat,\alpha_{\kappa})} \Phi
\]

Let $S_{+,k}(y,\xi,\alpha)$ be the set of
$\gminiy=(y_0,\ldots,y_k)\in NS(X)^{k+1}$
with the following properties:
\begin{itemize}
\item
 $\sum y_i=y$ and $P^{\alpha}_{y_i}=P^{\alpha}_y$.
\item
 We have $a_0/r_0-a_1/r_1=m_0\cdot \xi$ for some $m_0<0$.
\item
 We have $a_i/r_i-a_{i+1}/r_{i+1}=m_i\cdot \xi$
 for some $m_i\leq 0$.
 In the case $a_i/r_i=a_{i+1}/r_{i+1}$,
 we have $y_i/r_i=y_{i+1}/r_{i+1}$.
\end{itemize}
Here $r_i$ and $a_i$  denote the rank
and the first Chern class of $y_i$, as usual.
Let $S_{-,k}(y,\xi,\alpha)$ denote the set of
$\gminiy=(y_0,\ldots,y_k)\in NS(X)^{k+1}$
with the following properties:
\begin{itemize}
\item
 $\sum y_i=y$ and $P^{\alpha}_{y_i}=P^{\alpha}_y$.
\item
 We have $a_0/r_0-a_1/r_1=m_0\cdot \xi$ for some $m_0>0$.
\item
 We have $a_i/r_i-a_{i+1}/r_{i+1}=m_i\cdot \xi$
 for some $m_i\geq 0$.
 In the case $a_i/r_i=a_{i+1}/r_{i+1}$,
 we have $y_i/r_i=y_{i+1}/r_{i+1}$.
\end{itemize}
Then, we put as follows:
\[
 S_{\kappa}(y,\xi,\alpha):=
 \bigcup_{k}S_{\kappa,k}(y,\xi,\alpha)
\]
For each $\gminiy=(y_0,\ldots,y_k)\in S_{\kappa,k}(y,\xi,\alpha)$,
we put as follows:
\[
 \nbigm(\gminiyhat,\alpha_{\kappa}):=
 \prod_{i=0}^k\nbigm(\yhat_i,\alpha_{\kappa})
\]

\begin{prop}
\label{prop;06.6.7.350}
The following equality holds:
\begin{equation}
\label{eq;06.6.7.300}
 \Phi(\yhat,\alpha_+)
-\Phi(\yhat,\alpha_-)
=
-\sum_{\gminiy\in S_+(y,\xi,\alpha)}
 \Bbar(\gminiy)\cdot
 \int_{\nbigm(\gminiyhat,\alpha_+)}
 \Psi(\gminiy)
+\sum_{\gminiy\in S_-(y,\xi,\alpha)}
 \Bbar(\gminiy)\cdot
 \int_{\nbigm(\gminiyhat,\alpha_-)}
 \Psi(\gminiy)
\end{equation}
Here, $\Psi(\gminiy)$ denote the elements of
$\bigotimes_{i=0}^k\nbigr(\Ehat^u_i)$
given by {\rm(\ref{eq;06.6.7.250})},
and $\Bbar(\gminiy)$ denote the numbers
given by {\rm(\ref{eq;06.6.7.151})}.
\hfill\qed
\end{prop}

As an example,
we write down the formula (\ref{eq;06.6.7.300})
in the case of rank $3$.
We use the notation in the subsubsection 
\ref{subsubsection;06.6.7.1}.
We put as follows:
\[
 U_1(\alpha):=
 \bigl\{(y_1,y_2)\in U_i\,\big|\,
 P^{\alpha}_{y_i}=P^{\alpha}_y\bigr\}\,\,\,\,
 (i=1,2),
\quad
 U_3(\alpha):=
 \bigl\{(y_1,y_2,y_3)\in U_3\,\big|\,
 P_{y_i}^{\alpha}=P_y^{\alpha}
 \bigr\}
\]
We put as follows:
\[
 \begin{array}{l}
 S_{1+}:=\bigl\{
 (y_1,y_2)\in U_1(\alpha)\,\big|\,
 2a_1-a_2=m\cdot \xi\,\,(m<0)
 \bigr\}\\
 S_{2,+}:=\bigl\{
 (y_1,y_2)\in U_2(\alpha)\,\big|\,
 a_1-2a_2=m\cdot \xi\,\,(m<0)
 \bigr\}\\
 S_{3+}:=\bigl\{
 (y_1,y_2,y_3)\in U_3(\alpha)\,\big|\,
 a_1-a_2=m_1\cdot\xi,\,
 a_2-a_3=m_2\cdot \xi\,\,
 (m_i<0)
 \bigr\}\\
 S_{4+}:=\bigl\{
 (y_1,y_2,y_3)\,\big|\,
 a_1-a_2=m\cdot \xi\,\,(m<0),\,\,\,
 y_2=y_3
 \bigr\}
 \end{array}
\]
Here, $a_i$ denotes the first Chern class of $y_i$.
We also put as follows:
\[
  \begin{array}{l}
 S_{1-}:=\bigl\{
 (y_1,y_2)\in U_1(\alpha)\,\big|\,
 2a_1-a_2=m\cdot \xi\,\,(m>0)
 \bigr\}\\
 S_{2,-}:=\bigl\{
 (y_1,y_2)\in U_2(\alpha)\,\big|\,
 a_1-2a_2=m\cdot \xi\,\,(m>0)
 \bigr\}\\
 S_{3-}:=\bigl\{
 (y_1,y_2,y_3)\in U_3(\alpha)\,\big|\,
 a_1-a_2=m_1\cdot\xi,\,
 a_2-a_3=m_2\cdot \xi\,\,
 (m_i>0)
 \bigr\}\\
 S_{4-}:=\bigl\{
 (y_1,y_2,y_3)\,\big|\,
 a_1-a_2=m\cdot \xi\,\,(m>0),\,\,\,
 y_2=y_3
 \bigr\}
 \end{array}
\]

We put as follows for $\kappa=\pm$:
\begin{multline}
 \Lambda(\alpha_{\kappa}):= \Phi(\yhat,\alpha_{\kappa})
+\sum_{(y_1,y_2)\in S_{1,{\kappa}}}
 \frac{1}{3}\int_{\nbigm(\yhat_1)\times\nbigm(\yhat_2,\alpha_{\kappa})}
  \underset{t}\Res\Bigl(
 \nbigp\bigl(\Ehat_1^u\cdot e^{-t},\,
  \Ehat_2^u\cdot e^{t/2}\bigr) \Bigr)\\
+\sum_{(y_1,y_2)\in S_{2,\kappa}}
 \frac{2}{3}\int_{\nbigm(\yhat_1,\alpha_{\kappa})\times\nbigm(\yhat_2)}
  \underset{t}\Res\Bigl(
 \nbigp\bigl(\Ehat_1^u\cdot e^{-t/2},\,
  \Ehat_2^u\cdot e^{t}\bigr) \Bigr)\\
+\sum_{(y_1,y_2,y_3)\in S_{3,\kappa}}
 \frac{1}{3}\int_{\prod_i\nbigm(\yhat_i)}
 \underset{t_1}\Res\,
 \underset{t_2}\Res\Bigl(
 \nbigp\bigl(\Ehat_1^u\cdot e^{-t_1-t_2/2},\,
 \Ehat_2^u\cdot e^{t_1-t_2/2},\,
 \Ehat_3^u\cdot e^{t_2}\bigr)\Bigr)\\
+\sum_{(y_1,y_2,y_3)\in S_{4,\kappa}}
 \frac{1}{6}\int_{\prod_i\nbigm(\yhat_i)}
 \underset{t_1}\Res\,
 \underset{t_2}\Res\Bigl(
 \nbigp\bigl(\Ehat_1^u\cdot e^{-t_1-t_2/2},\,
 \Ehat_2^u\cdot e^{t_1-t_2/2},\,
 \Ehat_3^u\cdot e^{t_2}\bigr)\Bigr)
\end{multline}
Then, Proposition \ref{prop;06.6.7.350}
says $\Lambda(\alpha_+)=\Lambda(\alpha_-)$.

\subsection{Weak Intersection Rounding Formula}
\label{subsection;06.6.7.400}
\subsubsection{Preliminary}

We prepare some terminology to state the theorem.
An oriented tree is a tree provided with an orientation
for each edge.
For an oriented tree $R$,
let $V(R)$ denote the set of the vertices of $R$.
We have the natural order $\leq_R$ on $V(R)$.
Let $V_{\max}(R)$ denote the set of the maximal
vertices of $R$,
and we put $V^{\circ}(R):=V(R)-V_{\max}(R)$.
An oriented tree $R$ is called 
an indexed rooted oriented binary plane tree of rank $r$,
if the following conditions are satisfied:
\begin{itemize}
\item 
 $R$ is embedded in $\real\times\real_{\geq\,0}$,
 and it intersects with $\real\times\{0\}$
 transversally.
 The intersection $R\cap(\real\times\{0\})$
 consists of the maximal points.
\item
 There exists the unique minimal vertex $o$
 with respect to the order induced by the orientation.
 The vertex is called the root.
\item
 The maximal vertices are
 $\{(0,0),(1,0),\ldots,(r-1,0)\}$.
\item
 For any $v\in V(R)-\bigl(\{o\}\cup V_{\max}(R)\bigr)$,
 there are three edges which contain $v$.
\item
 A bijective map $\varphi:V^{\circ}(R)\lrarr \{1,\ldots,r-1\}$
 is provided, which preserves the orders
 on the sets $V^{\circ}(R)$ and $\{1,\ldots,r-1\}$.
\end{itemize}
The notion of isotopy is naturally given.
Let $\gbigt(r)$ denote the set of the isotropy classes
of an indexed rooted oriented binary plane trees
of rank $r$.
It is easy to see
$\#\gbigt(2)=1$, $\#\gbigt(3)=2$
and $\#\gbigt(4)=4$, for example.

\begin{rem}
$\gbigt(r)$ parameterizes the composition rules.
\hfill\qed
\end{rem}

Let $R$ be an indexed rooted oriented binary plane tree
of rank $r$. 
For any $v\in V(R)$,
we put as follows:
\[
 r(v):=\#\bigl\{u\in V_{\max}(R)\,\big|\,
 u\geq_R v\bigr\}
\]
We put $j(v):=\max\bigl\{j\,\big|\,(j,0)\geq_R v\bigr\}$.

Take $v\in V^{\circ}(R)$.
It is easy to see that there are two edges
 $e_1$ and $e_2$
going out from $v$.
Let $v_i$ denote the vertex of $e_i$
which is not $v$.
We may assume
$j(v_1)<j(v_2)$,
and we put $v^l=v_1$
and $v^r=v_2$.

If $v$ is not the root,
there is the unique vertex $v^b$
of $R$ such that there exists the edge
from $v^b$ to $v$.
We put as follows:
\[
 \sign(v):=\left\{
 \begin{array}{ll}
 1 & (v=(v^b)^r)\\
 -1 & (v=(v^b)^l)
 \end{array}
 \right.
\]
We take $t_{1},\ldots, t_{r-1}$ be formal variables,
and we put as follows for $i=0,\ldots,r-1$:
\[
 T_R^i:=\sum_{\substack{v\leq_R (i,0)\\v\neq o}}
 \sign(v) \frac{t_{\varphi(v^b)}}{r(v)}
\]

We identify $(i,0)\in V_{\max}(R)$ and $i$.
For any $\veca=(a_0,\ldots,a_{r-1})\in NS^1(X)$
and each $v\in V(R)$,
we put as follows:
\[
 a_{v}:=\sum_{
 \substack{u\in V_{\max}(R)\\ u\geq_Rv }} a_{u}
\]
Let $a$ be an element of $NS^1(X)$.
Let $\vecxi=(\xi_1,\ldots,\xi_{r-1})$ be an element of
$NS^1(X)^{r-1}$ such that $\xi_i$ are linearly independent.
We denote by $S(r,a,\vecxi,R)$
the set of $\veca=(a_0,\ldots,a_{r-1})\in NS^1(X)$
satisfying  $\sum a_i=a$
and the following condition:
\begin{itemize}
\item
 There exists a positive rational number $A$
 such that the following equality holds  in $NS^1(X)\otimes\rnum$
 for any $v\in V^{\circ}(R)$:
\[
 \frac{a_{v^l}}{r(v^l)}
-\frac{a_{v^r}}{r(v^r)}
=A_v\cdot \xi_{\varphi(v)}
\]
\end{itemize}

Let $\veca$ be any element of $S(r,a,\vecxi,R)$.
Let $y$ be an element of $\Type$
such that $\rank(y)=r$ and $\det(y)=a$.
Let $n$ denote the second Chern class of $y$.
We put as follows:
\[
 \nbigm(\veca,y):=
 \coprod_{\sum_{i=0}^{r-1}n_i=N(\veca,y)}
 \prod_{i=0}^{r-1}\bigl(
 X^{[n_i]}\times\Pic(a_i)\bigr),
\quad
 N(\veca,y):=n-\frac{a^2}{2}+\sum_{i=1}^{r-1}\frac{a_i^2}{2}
\]

On $\nbigm(\veca,y)\times X$,
we have the sheaf $\Ehat_i^u$
which is the pull back of the universal
sheaf on $X^{[n_i]}\times\Pic(a_i)\times X$
via the naturally defined projection.
Let $G$ denote the $(r-1)$-dimensional torus
$\Spec k[\tau_1,\tau_1^{-1},\ldots,\tau_{r-1},\tau_{r-1}^{-1}]$.
Let $e^{w\cdot t_i}$ denote the trivial line bundle
with the $G$-action
induced by the action of $\Spec k[\tau_i,\tau_i^{-1}]$
of weight $w$.
We have the following element of
the $K$-group of the $G$-equivariant
coherent sheaves on $\nbigm(\veca,y)$:
\[
-Rp_{X\,\ast}\nrhom\bigl(\Ehat_i^u\cdot e^{T_i},\,\,
 \Ehat_j^u\cdot e^{T_j}
 \bigr)
-Rp_{X\,\ast}\nrhom\bigl(\Ehat_j^u\cdot e^{T_j},\,\,
 \Ehat_i^u\cdot e^{T_i}
 \bigr)
\]
The equivariant Euler class is denoted by
$Q\bigl(\Ehat_i^u\cdot e^{T_i},
 \Ehat_j^u\cdot e^{T_j}\bigr)$.
We regard it as the element of
$\bigotimes_{i=0}^{r-1}\nbigr(\Ehat_i^u)\otimes
\gbigr(t_1,\ldots,t_{r-1})$.
We put as follows:
\[
 Q\bigl(
 \Ehat_0^u\cdot e^{T_0},
 \Ehat_1^u\cdot e^{T_1},\ldots,
 \Ehat_{r-1}^u\cdot e^{T_{r-1}} \bigr)
:=\prod_{i<j}
 Q\bigl(\Ehat_i^u\cdot e^{T_i},\,
 \Ehat_j\cdot e^{T_j}\bigr)
\]
Let $P$ be an element of $\nbigr$.
We obtain
$P\bigl(\bigoplus_{i=0}^{r-1} \Ehat_i^u\cdot e^{T_i}\bigr)\in 
 \bigotimes_{i=0}^{r-1}\nbigr(\Ehat_i^u)
 [t_1,\ldots,t_{r-1}]$
by the homomorphisms in the subsubsection
\ref{subsubsection;06.6.19.25}.
Then, we put as follows:
\begin{equation}
\label{eq;06.6.19.101}
 \Psi(\veca,R):=
 \underset{t_{r-1}}\Res\,
 \underset{t_{r-2}}\Res\cdots
 \underset{t_{1}}\Res
\left(
\frac{
 P\bigl(\bigoplus_{i=0}^{r-1}
 \Ehat_i^ue^{T_i} \bigr)}
{Q\bigl(
 \Ehat_0^u e^{T_0},
 \Ehat_1^u e^{T_1},\ldots,
 \Ehat_{r-1}^u e^{T_{r-1}} \bigr)}
\right)
\in\bigotimes_{i=0}^{r-1}\nbigr(\Ehat_i^u)
\end{equation}

\subsubsection{Statement}

We take an element 
$\vecxi=(\xi_1,\xi_2,\ldots,\xi_l)\in NS^1(X)^l$
such that 
$\xi_1,\ldots,\xi_l$ are linearly independent.
We put $W^{\vecxi}:=\bigcap_{i=1}^lW^{\xi_i}$.
A connected component of
$ W^{\vecxi}\setminus
 \bigcup_{W\neq W^{\xi_i}}W$ is called a tile.

Let $\varphi:\{1,\ldots,l\}\lrarr \{1,-1\}$ be a map.
We have the chamber $C(\varphi)$ satisfying 
the following conditions:
\begin{itemize}
\item The closure of $C(\varphi)$ contains $T$.
\item Let $H_{\varphi}$ be any ample line bundle
 contained in $C(\varphi)$.
 Then $\sign\bigl(\xi_i,H\bigr)=\varphi(i)$.
\end{itemize}
We take line bundles $H_{\varphi}\in C(\varphi)$.
We put $k(\varphi):=\#\{i\,|\,\varphi(i)=-1\}$.
For $\Phi=P(\Ehat^u)\in\nbigr(\Ehat^u)$,
we put as follows:
\[
 \nbigd_{\vecxi}^T\Phi(\yhat):=
 \sum_{\varphi}
 (-1)^{k(\varphi)}\cdot
 \Phi_{H_{\varphi}}(\yhat)
\]
Here, $\Phi_{H_{\varphi}}(\yhat)$ is given as in 
(\ref{eq;06.6.19.100}).

\begin{thm}
\label{thm;06.6.7.200}
Assume $p_g=0$.
Let $y$ be an element of $\Type$,
and let $\vecxi=(\xi_1,\ldots,\xi_l)$ and $T$ be as above.
\begin{itemize}
\item 
 In the case $l\geq \rank(y)$,
 we have $\nbigd_{\vecxi}^T\Phi(y)=0$.
\item
 In the case $l=\rank(y)-1$,
 the number $\nbigd^{T}_{\vecxi}\Phi(y)$ is independent of 
 the choice of $T$,
 and the following equality holds:
\begin{equation}
\label{eq;06.6.7.10}
 \nbigd_{\vecxi}\Phi(y)
=\sum_{R\in \gbigt(r)}
 \sum_{\veca\in S(r,a,\vecxi,R)}
 \int_{\nbigm(\veca,y)} \Psi(\veca,R)
\end{equation}
Here, $\Psi(\veca,R)$ are given in 
{\rm(\ref{eq;06.6.19.101})}.
\end{itemize}
\end{thm}

The theorem will be proved in the subsubsection 
\ref{subsubsection;06.6.11.1000}
after the combinatorial preparation.
We have already written down the formula (\ref{eq;06.6.7.10})
for the rank $3$ case
in the subsubsection \ref{subsubsection;06.6.8.300}.

\subsubsection{Preparation from combinatorics}

Let $S$ be a finite set.
Let $\Map(S,\{\pm\})$ denote the set
of maps of $S$ to $\{\pm\}$.
For any $\varphi\in\Map(S,\{\pm\})$,
let $k(\varphi)$ denote the number
of $\bigl\{i\,\big|\,\varphi(i)=-\bigr\}$.
Let $I$ be a subset of $S$.
For any $\varphi_I\in\Map(I,\{\pm\})$,
let $\Ind_S^-(\varphi_I)$ and 
 $\Ind_S^+\varphi_I$
be elements of $\Map(S,\{\pm\})$
determined as follows:
\[
 \Ind_S^-\bigl(
 \varphi_I\bigr)(i):=\left\{
 \begin{array}{ll}
 \varphi_I(i) & (i\in I)\\
 - & (i\not\in I)
 \end{array}
 \right.
\quad
 \Ind_S^+\bigl( \varphi_I\bigr)
 (i):=\left\{
 \begin{array}{ll}
 \varphi_I(i) & (i\in I)\\
 + & (i\not\in I)
 \end{array}
 \right.
\]

\begin{lem}
\label{lem;06.6.7.15}
For any functions
$F_i:\Map(S,\{\pm\})\lrarr \cnum$ $(i=1,2)$,
we have the following equality:
\begin{multline}
\label{eq;06.6.7.5}
\sum_{\varphi\in\Map(S,\{\pm\})}(-1)^{k(\varphi)}
 F_1(\varphi)\cdot F_2(\varphi)
=\\
 \sum_{I\sqcup J=S}
\left(
 \sum_{\varphi_I\in \Map(I,\{\pm\})}
 (-1)^{k(\varphi_I)}
 F_1\bigl(\Ind_S^-(\varphi_I)\bigr)
\right)
\cdot
\left(
\sum_{\varphi_J\in\Map(J,\{\pm\})}
 (-1)^{k(\varphi_J)}
 F_2\bigl(\Ind_S^+(\varphi_J)\bigr)
 \right)
\end{multline}
\end{lem}
\pf
We use an induction on $|S|$.
In the case $|S|=1$,
we can identify $\Map(S,\{\pm\})$
and $\{\pm\}$.
We have the following:
\[
 F_1(+)\cdot F_2(+)-F_1(-)\cdot F_2(-)
=\bigl(F(+)-F(-)\bigr)\cdot G(+)
+F(-)\cdot \bigl(G(+)-G(-)\bigr)
\]
Thus, the claim holds in this case.

Assume that the claim holds in the case $|S|=k-1$,
and we will show that the claim in the case $|S|=k$.
We take any element $s\in S$,
and we put $S':=S-\{s\}$.
The left hand side of (\ref{eq;06.6.7.5})
can be rewritten as follows:
\begin{multline}
\label{eq;06.6.7.6}
 \sum_{\substack{\varphi\in\Map(S,\{\pm\})\\
\varphi(s)=+ }}
 (-1)^{k(\varphi)}\cdot F_1(\varphi)\cdot F_2(\varphi)
+\sum_{\substack{\varphi\in\Map(S,\{\pm\})\\
 \varphi(s)=-}}
 (-1)^{k(\varphi)}\cdot F_1(\varphi)\cdot F_2(\varphi) \\
=\sum_{\varphi'\in\Map(S',\{\pm\})}
 (-1)^{k(\varphi')}
 F_1\bigl(\Ind_S^+(\varphi^{\prime})\bigr)\cdot 
 F_2\bigl(\Ind_S^+(\varphi^{\prime})\bigr)
-\sum_{\varphi'\in\Map(S',\{\pm\})}
 (-1)^{k(\varphi')}
 F_1\bigl(\Ind_S^-(\varphi^{\prime})\bigr)\cdot 
 F_2\bigl(\Ind_S^+(\varphi^{\prime})\bigr)
\end{multline}
By using the hypothesis of the induction,
the right hand side of (\ref{eq;06.6.7.6})
can be rewritten as follows:
\begin{multline}
\label{eq;06.6.7.7}
 \sum_{I\sqcup J=S'}
 \sum_{\substack{\varphi_I\in \Map(I,\{\pm\})\\ 
 \varphi_J\in \Map(J,\{\pm\})}}
 (-1)^{k(\varphi_I)+k(\varphi_J)}
 F_1\bigl(
 \Ind_S^+(\Ind_{S'}^-(\varphi_I)) \bigr)
\cdot
 F_2\bigl(\Ind_S^+\Ind_{S'}^+(\varphi_J)
 \bigr)\\
-\sum_{I\sqcup J=S'}
 \sum_{\substack{\varphi_I\in \Map(I,\{\pm\})\\
 \varphi_J\in\Map(J,\{\pm\}) }}
 (-1)^{k(\varphi_I)+k(\varphi_J)}
 F_1\bigl(
 \Ind_S^{-}(\Ind_{S'}^-(\varphi_I)) \bigr)
\cdot F_2\bigl(
 \Ind_S^{-}(\Ind_{S'}^{+}(\varphi_J)) \bigr)  \\
=\sum_{I\sqcup J=S'}
 \sum_{\substack{\varphi_I\in \Map(I,\{\pm\})\\ 
 \varphi_J\in \Map(J,\{\pm\})}}
 (-1)^{k(\varphi_I)+k(\varphi_J)}
\Bigl(
 F_1\bigl(
 \Ind_S^+(\Ind_{S'}^-(\varphi_I)) \bigr)
-F_1\bigl(\Ind_S^-(\Ind_S'^-(\varphi_I))\bigr)
\Bigr)
\cdot
 F_2\bigl(\Ind_S^+\Ind_{S'}^+(\varphi_J)
 \bigr)\\
+\sum_{I\sqcup J=S'}
 \sum_{\substack{\varphi_I\in \Map(I,\{\pm\})\\
 \varphi_J\in\Map(J,\{\pm\}) }}
 (-1)^{k(\varphi_I)+k(\varphi_J)}
 F_1\bigl(
 \Ind_S^{-}(\Ind_{S'}^-(\varphi_I)) \bigr)
\cdot 
\Bigl(
 F_2\bigl(\Ind_S^+(\Ind_{S'}^+(\varphi_J))\bigr)
-F_2\bigl(\Ind_S^{-}(\Ind_{S'}^{+}(\varphi_J)) \bigr)
\Bigr)
\end{multline}
It is easy to observe the the right hand side of
(\ref{eq;06.6.7.7})
can be rewritten as follows:
\begin{multline}
 \sum_{\substack{I\sqcup J=S\\ s\in I}}
 (-1)^{k(\varphi_I)+k(\varphi_J)}
 F_1\bigl(\Ind_S^-(\varphi_I)\bigr)\cdot
 F_2\bigl(\Ind_S^+(\varphi_J)\bigr)
+\sum_{\substack{I\sqcup J=S\\ s\in J}}
 (-1)^{k(\varphi_I)+k(\varphi_J)}
 F_1\bigl(\Ind_S^-(\varphi_I)\bigr)\cdot
 F_2\bigl(\Ind_S^+(\varphi_J)\bigr)\\
=\sum_{I\sqcup J=S}
 (-1)^{k(\varphi_I)+k(\varphi_J)}
 F_1\bigl(\Ind_S^-(\varphi_I)\bigr)\cdot
 F_2\bigl(\Ind_S^+(\varphi_J)\bigr)
\end{multline}
Thus we are done.
\hfill\qed

\subsubsection{Proof of Theorem \ref{thm;06.6.7.200}}
\label{subsubsection;06.6.11.1000}

We consider the following statements:
\begin{description}
\item[$P(r)$:]
 Assume $\rank(\vecy)\leq r$.
 For any 
 $\vecxi=(\xi_1,\ldots,\xi_l)$ such that $l\geq \rank(\vecy)$,
 we have $\nbigd^T_{\vecxi}\Phi(\yhat)=0$.
\item[$Q(r)$:]
 Assume $\rank(\vecy)\leq r$.
 For any
 $\vecxi=(\xi_1,\ldots,\xi_l)$ such that $l= \rank(\vecy)-1$,
 the number $\nbigd^T_{\vecxi}\Phi(\yhat)$ is independent
 of the choice of the tile $T$.
 Moreover, the equality (\ref{eq;06.6.7.10}) holds.
\end{description}

It is easy to see $Q(r)$ implies $P(r)$.
We have already known that $Q(1)$ and $Q(2)$ hold.
We will prove $P(r)$, and hence $Q(r)$, by an induction on $r$.

Let $y_i$ $(i=1,\ldots,s)$ be elements of $\Type$,
and let $\Phi_i$ be elements of $\nbigr(\Ehat_i)$.
We put as follows:
\[
 \nbigd^T_{\vecxi}
\left(
 \prod_{i=1}^s\Phi_i(y_i)\right):=
 \sum_{\varphi}
 (-1)^{k(\varphi)}
 \prod_{i=1}^s\Phi_{i,H_{\varphi}}(y_i)
\]

For a given $\vecxi=(\xi_1,\ldots,\xi_l)$
and a subset $I\subset\{1,\ldots,l\}$,
we put $\vecxi_I:=\bigl(\xi_i\,|\,i\in I\bigr)$.
For a tile $T$ in $W^{\vecxi}$,
let $T^I_-$ be the tile of $W^{\vecxi_I}$
such that 
 $(H,\xi_j)<0$ for any $j\in \{1,\ldots,l\}\setminus I$
and for any element $H\in T^I_-$.
Similarly,
let $T^I_+$ be the tile of $W^{\vecxi_I}$
such that
$(H,\xi_j)>0$ for any $j\in \{1,,\ldots,l\}\setminus I$
and for any element $H\in T^I_+$.
We put $S=\{1,\ldots,l\}$.
Then, we obtain the following equality
from Lemma \ref{lem;06.6.7.15}:
\begin{equation}
\label{eq;06.6.7.20}
 \nbigd_{\vecxi}^T\left(
\prod_{i=1}^s\Phi_i(\yhat_i)\right)
=\sum_{I\sqcup J=S}
  \nbigd^{T^I_-}_{\vecxi_I}\bigl(
 \Phi_1(\yhat_1)\bigr)
\cdot
 \nbigd^{T^J_+}_{\vecxi_J}\left(
 \prod_{i=2}^s\Phi_i(\yhat_i)\right)
\end{equation}

\begin{lem}
\label{lem;06.6.7.25}
Assume that $P(r-1)$ holds.
Let $y_1,\ldots,y_s$ $(s\geq 2)$ be  elements of $\Type$
such that $\rank y_i\leq r-1$.
For $\vecxi=(\xi_1,\ldots,\xi_l)$ such that
$l\geq \sum_{i=1}^s\rank(y_i)-s$,
we have
$\nbigd_{\vecxi}^T\prod_{i=1}^s\Phi_i(y_i)=0$.
\end{lem}
\pf
We use an induction on $s$.
In the case $s=1$,
the claim follows from the assumption $P(r-1)$.
Assume that the claim holds in the case $s-1$,
and we will show the claim in the case $s$.
Due to the assumption $P(r-1)$,
the contribution from $(I,J)$ is $0$
in (\ref{eq;06.6.7.20}),
if the inequality $|I|\geq \rank(y_1)$ holds.
If $|I|<\rank(y_1)$ holds,
the inequality $|J|\geq \sum_{i=2}^s\rank y_i-s+1$ holds.
Therefore, we obtain the vanishing of the contribution
due to the hypothesis of the induction.
Thus we are done.
\hfill\qed

\begin{lem}
\label{lem;06.6.7.101}
Assume that $Q(r-1)$ holds.
Let $y_i$ $(i=1,2)$ be elements of $\Type$
such that $\rank y_i\leq r-1$.
For $\vecxi=(\xi_1,\ldots,\xi_{r-2})$,
the number $\nbigd_{\vecxi}^T(\Phi_1(\yhat_1)\cdot \Phi_2(\yhat_2))$
is independent of the choice of a tile $T$,
and the following equality holds:
\[
 \nbigd_{\vecxi}\bigl(\Phi_1(\yhat_1)\cdot\Phi_2(\yhat_2)\bigr)
=\sum_{\substack{I_1\sqcup I_2=S\\ |I_i|=\rank y_i}}
 \nbigd_{\vecxi_{I_1}}\bigl(
\Phi_1(\yhat_1)\bigr)
\cdot
 \nbigd_{\vecxi_{I_2}}\bigl(
 \Phi_2(\yhat_2)\bigr)
\]
\end{lem}
\pf
It can be shown by an argument
similar to the proof of Lemma \ref{lem;06.6.7.25}.
\hfill\qed

\vspace{.1in}

Let us show the claim $Q(r)$ assuming 
that $Q(r-1)$ holds.
We put $S':=\{2,\ldots,l\}$.
For a map $\lambda:S'\lrarr \{\pm\}$,
let $T(\lambda)$ denote the tile of
$W^{\xi_1}$ determined by the following conditions:
\begin{itemize}
\item The closure of $T(\lambda)$ contains $T$.
\item 
 Let $H$ be any ample line bundle contained in $T(\lambda)$.
 Then we have
 $\sign(H,\xi_i)=\lambda(i)$.
\end{itemize}
From the definition of $\nbigd_{\vecxi}^T\Phi(\yhat)$,
we obtain the following equality:
\begin{equation}
\label{eq;06.6.11.1010}
 \nbigd_{\vecxi}^T\bigl(\Phi(\yhat)\bigr)
=\sum_{\lambda} (-1)^{k(\lambda)}\cdot
 \nbigd_{\xi_1}^{T(\lambda)}\bigl(
\Phi(\yhat)\bigr)
\end{equation}
Let $H_{\kappa,\lambda}$ denote
an ample line bundle contained in
the chamber corresponding to
$\Ind_S^{\kappa}(\lambda)$.
We can rewrite (\ref{eq;06.6.11.1010})
as follows,
by using Theorem \ref{thm;06.6.7.100}:
\begin{equation}
\label{eq;06.6.11.1111}
 \nbigd_{\vecxi}^T\Phi(\yhat)
+\sum_{\lambda}(-1)^{k(\lambda)}
 \sum_{\gminiy\in S_+}\Bbar(\gminiy)
 \cdot \int_{\nbigm_{H_{+\lambda}}(\gminiyhat)}\Psi(\gminiy)
-\sum_{\lambda}(-1)^{k(\lambda)}
 \sum_{\gminiy\in S_-}\Bbar(\gminiy)
 \cdot \int_{\nbigm_{H_-\lambda}(\gminiyhat)}\Psi(\gminiy)
=0
\end{equation}

We put $\vecxi':=(\xi_2,\ldots,\xi_l)$.
Let $T_{\kappa}^{\vecxi'}$ $(\kappa=\pm)$
be the tiles of $W^{\vecxi'}$
determined by the following conditions:
\begin{itemize}
\item
 The closure of $T_{\kappa}^{\vecxi'}$
 contain $T$.
\item
 Let $H_{\kappa}$ be any element of $T_{\kappa}^{\vecxi'}$.
 Then, we have
 $(H_{-},\xi_1)<0<(H_+,\xi_1)$.
\end{itemize}

For $\gminiy=(y_0,\ldots,y_k)$,
we put as follows:
\[
 \nbigd_{\vecxi'}^{T^{\vecxi'}_{\kappa}}
 \bigl(\Psi(\gminiyhat)\bigr):=
 \nbigd_{\vecxi'}^{T^{\vecxi'}_{\kappa}}
 \Psi(\gminiy)(\yhat_0,\ldots,\yhat_k).
\]

Then, we can rewrite (\ref{eq;06.6.11.1111})
as follows:
\begin{equation}
\label{eq;06.6.7.102}
\nbigd_{\vecxi}^T\bigl(\Phi(\yhat)\bigr)
+\sum_{\gminiy\in S_+}\Bbar(\gminiy)\cdot
 \nbigd_{\vecxi'}^{T_+^{\vecxi'}}
 \bigl(\Psi(\gminiyhat)\bigr)
-\sum_{\gminiy\in S_-}\Bbar(\gminiy)\cdot
 \nbigd_{\vecxi'}^{T_-^{\vecxi'}}
 \bigl(\Psi(\gminiyhat)\bigr)
=0
\end{equation}
Due to Lemma \ref{lem;06.6.7.25},
the contribution from $\gminiy=(y_0,\ldots,y_k)$
is $0$, if $k\geq 2$ holds.
In the case $k=1$,
the contributions do not depend
on $T_{\kappa}^{\vecxi'}$.
So, we can omit to denote them in the following argument.

We put as follows:
\[
 S_1:=\bigl\{
 (y_0,y_1)\in\Type^2\,\big|\,
 y_0+y_1=y,\,\,a_0/r_0-a_1/r_1=A\cdot \xi_1\,\,(A<0)
 \bigr\},
\]
\[
 S_2:=\bigl\{
 (y_0,y_1)\in\Type^2\,\big|\,
 y_0+y_1=y,\,\,a_0/r_0-a_1/r_1=A\cdot \xi_1\,\,(A>0)
 \bigr\}
\]
\[
 S_3:=\bigl\{
 (y_0,y_1)\in\Type^2\,\big|\,
 y_0+y_1=y,\,\,a_0/r_0=a_1/r_1,\,\,
 b_0<b_1
 \bigr\}
\]
Then, we obtain the following equality
from  (\ref{eq;06.6.7.102})
and Lemma \ref{lem;06.6.7.101}:
\begin{multline}
0=\nbigd_{\vecxi}^T\bigl(
 \Phi(\yhat)\bigr)
+\sum_{(y_0,y_1)\in S_1}
 \frac{r_0}{r}
 \nbigd_{\vecxi'}\bigl(
\Psi(\yhat_0,\yhat_1)\bigr)
+\sum_{(y_0,y_1)\in S_3}
 \frac{r_0}{r}
 \nbigd_{\vecxi'}\bigl(
 \Psi(\yhat_0,\yhat_1)\bigr)\\
-\sum_{(y_0,y_1)\in S_2}
 \frac{r_0}{r}
 \nbigd_{\vecxi'}\bigl(
\Psi(\yhat_0,\yhat_1)\bigr)
-\sum_{(y_0,y_1)\in S_3}
 \frac{r_0}{r}
 \nbigd_{\vecxi'}\bigl(
\Psi(\yhat_0,\yhat_1)\bigr)
\end{multline}
The contributions from $S_3$ are cancelled out.
We have the bijection $S_1\lrarr S_2$
given by $(y_0,y_1)\longmapsto (y_1,y_0)$.
We also have 
$ \nbigd_{\vecxi'}\Psi(\yhat_0,\yhat_1)
=-\nbigd_{\vecxi'}\Psi(\yhat_1,\yhat_0)$
for $(y_0,y_1)\in S_1$ and $(y_1,y_0)\in S_2$.
Therefore, we obtain the following equality:
\begin{equation}
\label{eq;06.6.7.150}
 \nbigd_{\vecxi}^T\Phi(\yhat)
=\sum_{(y_0,y_1)\in S_2} \nbigd_{\vecxi'}\Psi(\yhat_0,\yhat_1)
\end{equation}
By using (\ref{eq;06.6.7.150})
inductively,
we can obtain the formula (\ref{eq;06.6.7.10}).
Then, it is clear that 
$\nbigd^T_{\vecxi}\Phi(y)$ is independent of
the choice of $T$.
Applying the hypothesis of the induction
to (\ref{eq;06.6.7.150})
we can easily derive the equality (\ref{eq;06.6.7.10}).
Thus, the claim $Q(r)$ is proved,
and the proof of Theorem \ref{thm;06.6.7.200}
is finished.
\hfill\qed

\noindent
{\it Address\\
Department of Mathematics,
Kyoto University,
Kyoto 606-8502, Japan.\\
takuro@math.kyoto-u.ac.jp
}\\

\noindent
{\it Current address\\
 Max-Planck Institute for Mathematics,
Vivatsgasse 7, D-53111 Bonn, Germany,\\
 takuro@mpim-bonn.mpg.de
}

\printindex


\begin{thebibliography}{99}

\bibitem{altman-kleiman}
A. Altman and S. Kleiman, 
{\em Introduction to Grothendieck duality theory},
 Lecture Notes in Mathematics {\bf 146}, Springer-Verlag, 
Berlin-New York (1970)

\bibitem{aoki}
M. Aoki, {\em Hom stacks},
Manuscripta Math. {\bf 119} (2006), 37--56.
 
\bibitem{a}
M. Artin, {\em Versal deformation and algebraic stacks},
 Invent. math. {\bf 27}, (1974), 165--189.

\bibitem{b} 
  K. Behrend,
  {\em Gromov-Witten invariants in Algebraic Geometry},
  Invent.  Math.  {\bf 127}, (1997), 601--617. 

\bibitem{bf}
 K. Behrend and B. Fantechi,
 {\em The intrinsic normal cone},
 Invent. Math. {\bf 128}, (1997), 45--88.

\bibitem{bhosle}
 U. N. Bhosle, {\em Parabolic vector bundles
 on curves},
 Ark. Math. {\bf 27}, (1989) 15--22.

\bibitem{d-m}
 P. Deligne and D. Mumford,
 {\em The irreducibility of the space of curve of given genus},
 IHES  Publ.  Math.  {\bf 36} (1969), 75--109.

\bibitem{d} 
 S. K. Donaldson,
 {\em The Seiberg-Witten equations 
 and $4$-manifold topology},
 Bull. Amer. Math. Soc. {\bf 33}, (1996), 45--70.

\bibitem{edidin-graham}
D. Edidin and W. Graham,
{\em Equivariant intersection theory},
Invent. Math. {\bf 131} (1998), 
595--634. 

\bibitem{eg}
 G. Ellingsrud and L. G\"ottsche, % \"
 {\em Variation of Moduli spaces and Donaldson invariants
   under change of polarization},
   J. Reine Angew. Math. {\bf 467}, (1995), 1--49.

\bibitem{egl}
 G. Ellingsrud, L. G\"ottsche and  M. Lehn, %\"
 {\em On the Cobordism Class of the Hilbert Scheme},
  J. Algebraic Geom. {\bf 10} (2001), 81--100. 

\bibitem{fl1}
 P. M. N. Feehan and T. G. Leness,
 {\em ${\rm PU(2)}$ monopoles, I:
   Regularity, Uhlenbeck compactness, and transversality},
   J. Diff. Geom. {\bf 49}, (1998), 265--410.

\bibitem{fl2}
 P. M. N. Feehan and T. G. Leness,
  {\em ${\rm PU(2)}$ monopoles. II. Top-level Seiberg-Witten
    moduli spaces and Witten's conjecture in low degrees},
    J. Reine Angew. Math. {\bf 538}, (2001), 135--212.

\bibitem{fl3}
 P. M. N. Feehan and T. G. Lennes,
 {\em Homotopy equivalence and Donaldson
 invariants when $b^+ = 1$. I: 
  Cobordisms of moduli spaces and continuity of gluing maps}. 
               dg-ga/9812060.

\bibitem{fq}
 R. Friedman and Z. Qin,
 {\em Flips of moduli spaces and transition formulas
   for Donaldson polynomial invariants of rational surfaces},
  Communications in Analysis and Geometry, {\bf 3},
 (1995),  11--83.

\bibitem{f-o}
 K. Fukaya and K. Ono,
{\em Arnold Conjecture and Gromov-Witten invariant},
 Topology, {\bf 38} (1999), 933--1048.

\bibitem{ful}
 W. Fulton, {\em Intersection theory},
 Springer, (1984)

\bibitem{ga}
 O. Garc\'\i a-Prada,
 {\em Seiberg-Witten invariants and vortex equations,
  Sym\'etries quantiques} 
 (Les Houches, 1995),
 North-Holland, Amsterdam, (1998),  885--934.

\bibitem{gil}
 H. Gillet,
 {\em Intersection theory on algebraic stacks and Q-varieties},
  J. Pure and Applied algebra {\bf 34} (1984), 193--240.

\bibitem{gottsche}
L. G\"{o}ttsche,
{\em Modular forms and Donaldson invariants for
$4$-manifolds with $b_{+}=1$} 
J. Amer. Math. Soc. {\bf 9} (1996), 827--843.

\bibitem{gny}
 L. G\"{o}ttsche, H. Nakajima and K. Yoshioka,
 {\em Instanton counting and Donaldson invariants},
 AG/0606180.

\bibitem{gm1} 
 W. M. Goldman and J. J. Millson,
{\em The deformation theory of representations
  of fundamental groups of compact K\"ahler manifolds}, 
Publ. Math. IHES, {\bf 67} (1988), 43--96.

\bibitem{gm2} 
 W. M. Goldman and J. J. Millson,
 {\em The homotopy invariance of the Kuranishi space},
  Illinois Journal of Mathematics, {\bf 34}, (1990), 337--367.

\bibitem{gp}
 T. Graber and R. Pandharipande,
 {\em Localization of virtual classes},
 Invent. math. {\bf 135}, (1999), 487--518.

\bibitem{gro}
A. Grothendieck,
 {\it Techniques de construction et th\'eor\`ems
 d'existence en g\'eom\'etrie alg\'ebrique IV:
 Les sch\'emas de Hilbert}, S\'eminaire Bourbaki expos\'e {\bf 221},
 (1961), IHP, Paris.

\bibitem{hartshorne}
 R. Hartshorne,
 {\em Residues and Duality},
 Lecture notes in Mathematics No. {\bf 20}. Springer. (1996).

\bibitem{hl1}
 D. Huybrechts and M. Lehn,
 {\em Stable pairs on curves and surfaces},
  J. Alg. Geom. {\bf 4}, (1995), 67--104.

\bibitem{hl2}
 D. Huybrechts and M. Lehn,
 {\em Framed modules and their moduli},
  Internat. J. Math. {\bf 6}, No.2 (1995), 297-324.

\bibitem{ill}
 L. Illusie,
 {\em Complexe cotangent et deformations I},
 Lecture notes in Math. {\bf 239} (Springer, 1971)

\bibitem{ks}
  M. Kashiwara and P. Schapira,
  {\em Sheaves on Manifolds},
  Springer (1990)

\bibitem{ks2}
 M. Kashiwara and P. Schapira,
 {\em Categories and Sheaves},
 Springer (2006)

\bibitem{kont}
 M. Kontsevich,
 {\em Enumeration of rationnal curves via torus actions},
 in  {\em The moduli space of curves}, 
 Progr. Math., {\bf 129}, Birkh\"auser (1995), 335--368. %\"

\bibitem{kont2}
 M. Kontsevich,
 {\em Deformation Quantization of Poisson manifolds I},
  q-alg/9709040.

\bibitem{km} 
 D. Kotschick and J. W. Morgan,
  {\em $SO(3)$-invariants for $4$-manifolds with $b_2^+=1$.II,}
   J. Diff. Geom. {\bf 39} (1994), 433--456.

\bibitem{kresch}
  A. Kresch, {\em Cycle groups for Artin stacks},
  Invent. Math. {\bf 138} (1999), 495--536.

\bibitem{Kronheimer}
 P. B. Kronheimer,
 {\em Four-manifold invariants from higher-rank bundles},
  J. Diff. Geom. {\bf 70}, (2005),
 59--112

\bibitem{lau}
 G. Laumon and L. Moret-Bailly,
  {\em Champs alg\`ebriques},
  Springer, (2000)

\bibitem{li}
 J. Li,
 {\em Algebraic geometric interpretation of
  Donaldson's polynomial invariants},
  J. Diff. Geom. {\bf 37} (1993), 417-466.

\bibitem{l-t}
  J. Li and G. Tian,
 {\em Virtual Moduli cycles and GW invariant
 of algebraic varieties},
 J.  Amer.  Math.  Soc.  {\bf 11} (1998), 119--174.

\bibitem{mar}
 M. Maruyama,
 {\em Moduli of stable sheaves, I},
  J. Math. Kyoto Univ. {\bf 17} (1977), 91--126.

\bibitem{my}
 M. Maruyama and K. Yokogawa,
 {\em Moduli of parabolic stable sheaves},
  Math. Ann. {\bf 293}, (1992), 77--99.

\bibitem{mw}
 K. Matsuki and R. Wentworth,
 {\em Munford-Thaddeus princile
   on the moduli space of vector bundles
  on an algebraic surface,}
  Internat. J. Math. {\bf 8} (1997), 97--148.

\bibitem{matsu}
 H. Matsumura,
 {\it Commutative ring theory},
  Cambridge University Press,
  Cambridge-New York, (1989).

\bibitem{mochi} 
 T. Mochizuki,
 {\em Gromov-Witten class and a Perturmation theory in
    Algebraic Geometry},
   Amer. J. Math. {\bf 123}, (2001), 343--381.

\bibitem{m3}
 T. Mochizuki,
 {\em On the geometry of the parabolic Hilbert schemes},
 AG/0210210

\bibitem{mochi2}
T. Mochizuki,
 {\em Asymptotic behaviour of tame harmonic bundles
 and an application to pure twistor $D$-modules},
 to appear in Mem. AMS.

\bibitem{mor}
  J. W. Morgan,
 {\em Comparison of the Donaldson polynomial invariants
          with their algebro geometric analogues},
  Topology {\bf 32} (1993), 449-488.

\bibitem{GIT}
 D. Mumford, {\em Geometric invariant theory},
 Berlin Heidelberg New York: Springer (1965).

\bibitem{ost}
 CH. Okonek, A. Schmitt and A. Teleman,
 {\em Master spaces for stable pairs},
  Topology {\bf 38} (1999), 117--139.

\bibitem{olsson}
  M. Olsson,
 {\em Sheaves on Artin stacks},
 preprint, avariable from his web page.

\bibitem{pt}
 V. Y. Pidstrigach and A. N. Tyurin,
 {\em Localisation of Donaldson invariants
   along the Seiberg-Witten classes},
   dg-ga/9507004.

\bibitem{sch}
 M. Schlessinger,
 {\em Functors of Artin rings},
 Trans.  Amer. Math. Soc.  {\bf 130}, (1968), 208--222.

\bibitem{si}
 C. Simpson,
 {\em Moduli of representations of the fundamental group
 of a smooth projective variety I},
   IHES Publ. Math.. {\bf 79} (1994), 47--129.

\bibitem{th1} 
 M. Thaddeus,
 {\em Stable pairs, linear systems and Verlinde formula},
  Invent. Math. {\bf 117}, (1994), 317--353.

\bibitem{th2} 
 M. Thaddeus,
 {\em Geometric invariant theory and flips},
  J. Amer. Math. Soc. {\bf 9}, (1996), 691--723.

\bibitem{toen}
 B. To\"{e}n,
 {\em Higher and derived stacks: a global overview},
 math.AG/0604504

\bibitem{toen-vaquie}
B. To\"{e}n and M. Vaqui\'{e},
{\em Moduli of objects in dg-categories},
math.AG/0503269

\bibitem{vistoli}
 A. Vistoli,
 {\em Intersection theory on algebraic stacks and
 on their moduli spaces},
 Invent.  math.  {\bf 97}, (1989), 613--670. 

\bibitem{yamada}
K. Yamada,
{\em Polarization change of moduli of vector bundles
on surfaces with $p\sb g>0$},
J. Math. Kyoto Univ. {\bf 42} (2002), 607--616.

\bibitem{y3} 
 K. Yokogawa,
 {\em Moduli of stable pairs}, J. Math. Kyoto Univ., {\bf 31}, (1991),
 311-327.

\bibitem{y}
 K. Yokogawa,
 {\em Compactification of moduli of parabolic sheaves
          and moduli of parabolic Higgs sheaves},
  J. Math. Kyoto Univ. {\bf 33}, (1993), 451--504.

\end{thebibliography}
\end{document}